\documentclass[a4paper]{article}

\makeatletter
\renewcommand{\@oddhead}{\hfill\Large\thepage}
\makeatother
\usepackage{indentfirst}     
\usepackage{tabularx}        
\usepackage{ifthen}          
\usepackage{amsmath,amssymb} 

\usepackage{theorem}
\usepackage[mathscr]{eucal}
\usepackage[all]{xy}
\usepackage{graphics}
\usepackage{epic}

\usepackage{makeidx}
\usepackage{multind}
\makeindex{mindex}
\makeindex{mpoint}

\usepackage[cp1251]{inputenc}
\usepackage[russian]{babel}

\fontsize {14.4}{21.6}\selectfont
\makeatletter 
\renewcommand{\@biblabel}[1]{#1.}
\def\@startsection#1#2#3#4#5#6{%
  \if@noskipsec \leavevmode \fi
  \par
  \@tempskipa #4\relax
  \@afterindenttrue
  \ifdim \@tempskipa <\z@
    \@tempskipa -\@tempskipa \@afterindentfalse
  \fi
  \if@nobreak
    \everypar{}%
  \else
    \addpenalty\@secpenalty\addvspace\@tempskipa
  \fi
  \@ifstar
    {\@ssect{#3}{#4}{#5}{#6}}%
    {\@dblarg{\@sect{#1}{#2}{#3}{#4}{#5}{#6}}}}
\renewcommand{\@seccntformat}[1]{\csname the#1\endcsname.\hspace{1ex}}%
\def\@sect#1#2#3#4#5#6[#7]#8{%
  \ifnum #2>\c@secnumdepth
    \let\@svsec\@empty
  \else
    \refstepcounter{#1}%
    \protected@edef\@svsec{\@seccntformat{#1}\relax}%
  \fi
  \@tempskipa #5\relax
  \ifdim \@tempskipa>\z@
    \begingroup
      #6{%
        \@hangfrom{\hskip #3\relax\@svsec}%
          \interlinepenalty \@M #8\@@par}%
    \endgroup
    \csname #1mark\endcsname{#7}%
    \addcontentsline{toc}{#1}{
      \ifnum #2>\c@secnumdepth \else
        \protect\numberline{\csname the#1\endcsname.}%
      \fi
      #7}%
  \else
    \def\@svsechd{%
      #6{\hskip #3\relax
      \@svsec #8}%
      \csname #1mark\endcsname{#7}%
      \addcontentsline{toc}{#1}{%
        \ifnum #2>\c@secnumdepth \else
          \protect\numberline{\csname the#1\endcsname.}%
        \fi
        #7}}%
  \fi
  \@xsect{#5}}
\renewcommand{\@seccntformat}[1]{{\LARGE\csname the#1\endcsname}\hspace{1ex}}%
\def\@sect#1#2#3#4#5#6[#7]#8{%
  \ifnum #2>\c@secnumdepth
    \let\@svsec\@empty
  \else
    \refstepcounter{#1}%
    \protected@edef\@svsec{\@seccntformat{#1}\relax}%
  \fi
  \@tempskipa #5\relax
  \ifdim \@tempskipa>\z@
    \begingroup
      #6{%
        \@hangfrom{\hskip #3\relax\@svsec}%
          \interlinepenalty \@M #8\@@par}%
    \endgroup
    \csname #1mark\endcsname{#7}%
    \addcontentsline{toc}{#1}{
      \ifnum #2>\c@secnumdepth \else
        \protect\numberline{\csname the#1\endcsname.}%
      \fi
      #7}%
  \else
    \def\@svsechd{%
      #6{\hskip #3\relax
      \@svsec #8}%
      \csname #1mark\endcsname{#7}%
      \addcontentsline{toc}{#1}{%
        \ifnum #2>\c@secnumdepth \else
          \protect\numberline{\csname the#1\endcsname.}%
        \fi
        #7}}%
  \fi
  \@xsect{#5}}
\renewcommand*\l@section{\@dottedtocline{1}{0em}{1.1em}}
\renewcommand*\l@subsection{\@dottedtocline{2}{1.1em}{2.3em}}
\renewcommand*\l@subsubsection{\@dottedtocline{3}{3.4em}{3.0em}}
\renewcommand\@pnumwidth{1.5em} 
\renewcommand\@tocrmarg{2.2em}  
\makeatother  
\makeatletter
\renewcommand{\section}{\setcounter{table}{0}%
     \@startsection%
     {section}
     {1}
     {0pt}
     {0.5\baselineskip}
     {1\baselineskip}
     {\centering\normalfont \LARGE \bf }}
\makeatother
\makeatletter
\renewcommand{\subsection}{\@startsection
     {subsection}
     {2}
     {\parindent}
     {0.7\baselineskip}
     {0.7\baselineskip}
     {\normalfont\LARGE\bf}}
\makeatother
\makeatletter
\renewcommand{\subsubsection}{\@startsection
     {subsubsection}
     {3}
     {\parindent}
     {0.5\baselineskip}
     {-2ex
       plus -\fontdimen3\font
       minus -\fontdimen4\font}
     {\normalfont\Large\bf }}
\makeatother

\newtheorem{theorem}{\sc Теорема}[subsection]
\newtheorem{lemma}[theorem]{\sc Лемма}
\newtheorem{proposition}[theorem]{\sc Предложение}
\newtheorem{corollary}[theorem]{\sc Следствие}

{\theorembodyfont{\rmfamily}
\newtheorem{note}[theorem]{\sc Замечание}
\newtheorem{example}[theorem]{\sc Пример}
\newtheorem{remark}[theorem]{\sc Замечание}}

\makeatletter \@addtoreset{equation}{subsection} \makeatother

\newcommand{\IND}[1]{\index{mindex}{#1}}
\newcommand{\PNT}[1]{\index{mpoint}{#1}}

\def\itemii{параграф}
\def\itemiii{пункт}

\def\newoplus{\mathop{\textstyle\bigoplus}}
\def\unit{\mathtt{1}\mspace{-6.5mu}\mathbb{I}}

\begin{document}
\renewcommand{\baselinestretch}{1.230}
\Large

\pagenumbering{roman}

Below is Leonid Vaksman's monograph ``Quantum bounded symmetric domains''
(in Russian), which remains unfinished due to his sudden death. The Russian
text is preceded with an English translation of the table of contents
(with our brief comments and references to preprints where applicable) and
an English translation of the (initial part of the) introduction. It is we
who are responsible for the translation.

The PS or PDF file that you are reading was automatically compiled from a
source file which was written in a Russian version of TEX (rather than a
standard one), and we cannot guarantee that no errors occurred during its
automatic compilation by the electronic archive program. Therefore we would
like to mention that a copy of the Russian text of the book (without errors
of this type) is available at
  http://www.ilt.kharkov.ua/vaksman/

Yevgen Kolisnyk and Sergey Sinel'shchikov

\newpage

\begin{center}\bf\LARGE
QUANTUM BOUNDED SYMMETRIC DOMAINS
\\ \Large\fboxrule=1pt\framebox{Leonid L. Vaksman}
\\ a brief description of contents
\end{center}

\bigskip\bigskip

{\bf 1. INTRODUCTION}

An outline of the book. A history of research in quantum groups.

\bigskip

{\bf 2. QUANTUM DISC

\medskip

2.1. A warning for the pedants}

A review of some logical problems that may be caused by various
inconsistencies in the standard mathematical texts.

\medskip

{\bf 2.2. Topology

\smallskip

2.2.1. Commutative $C^*$-algebras}

The definition and basic properties of commutative $C^*$-algebras.

{\bf 2.2.2. An elementary review of the general theory of $C^*$-algebras.}

Basic facts on general $C^*$-algebras are expounded. Examples are included.

{\bf 2.2.3. The Fock representation}

Irreducible representations of the $*$-algebra of polynomials on the
quantum disc are described.

{\bf 2.2.4. A proof of Lemma 2.2.11

2.2.5. Continuous functions on the quantum disc}

Some properties of the algebra of continuous functions on the quantum disc
are discussed.

{\bf 2.2.6. Holomorphic functions on the quantum disc}

Some properties of the algebra of holomorphic functions on the quantum disc
are discussed.

{\bf 2.2.7. A supplement on spectral sets}

Spectral sets and unitary dilations are discussed.

\medskip

{\bf\boldmath 2.3. Symmetry

\smallskip

2.3.1. The Hopf algebras $U\mathfrak{sl}_2$ and $U_q\mathfrak{sl}_2$}

The general notion of a universal enveloping algebra of a Lie algebra is
discussed, along with the simplest examples and the associated Hopf algebra
structures.

{\bf\boldmath 2.3.2. The Hopf $*$-algebras $U\mathfrak{su}_{1,1}$,
$U\mathfrak{su}_2$, $U_q\mathfrak{su}_{1,1}$, $U_q\mathfrak{su}_2$}

The general notion of a Hopf $*$-algebra is discussed, along with some
involutions in the simplest examples.

{\bf\boldmath 2.3.3. $U_q\mathfrak{sl}_2$-module algebras
$\mathbb{C}[z]_q$, $\mathbb{C}[z^*]_q$}

$U_q\mathfrak{sl}_2$-module algebra structures on polynomial algebras are
discussed.

{\bf 2.3.4. The universal R-matrix}

Some properties of the universal R-matrix are discussed.

{\bf\boldmath 2.3.5. The $U_q\mathfrak{su}_{1,1}$-module $*$-algebra
$\operatorname{Pol}(\mathbb{C})_q$}

An R-matrix is used to introduce the necessary structures. See
math.QA/0410530.

{\bf\boldmath 2.3.6. The quantum $SL_2$}

The Hopf algebra $\mathbb{C}[SL_2]_q$ is investigated.

{\bf 2.3.7. A supplement on rings and modules of quotients}

Rings and modules of quotients are introduced; the Ore condition is
discussed.

{\bf 2.3.8. A supplement on algebras and modules in tensor categories}

Some properties of algebras and modules in braided tensor categories are
discussed.

\medskip

{\bf 2.4. The invariant integral

\smallskip

2.4.1. Invariant integrals and $*$-representations}

Some properties of invariant integrals on $A$-module algebras are
considered.

{\bf\boldmath 2.4.2. Does there exist a nonzero nonnegative
$U_q\mathfrak{sl}_2$-invariant integral?}

The question posed is answered in the negative in the case of
$\operatorname{Pol}(\mathbb{C})_q$.

{\bf 2.4.3. Distributions and finite functions}

The passage from $\operatorname{Pol}(\mathbb{C})_q$ to finite functions and
distributions. See math.QA/9808015.

{\bf\boldmath 2.4.4. An invariant integral on $\mathscr{D}(\mathbb{D})_q$}

The existence and uniqueness of an invariant integral on
$\mathscr{D}(\mathbb{D})_q$ is established. See math.QA/9808037.

\smallskip

{\bf 2.5. Differential calculi

\smallskip

2.5.1. Differential forms with polynomial coefficients}

The general notion of a differential calculus over an algebra is discussed.
See math.QA/0606499.

{\bf 2.5.2. Differential forms with finite coefficients}

The passage from $\Omega(\mathbb{C})_q$ to $\Omega(\mathbb{D})_q$.

{\bf 2.5.3. Invariant Hermitian metrics and the invariant Laplacian}

Hermitian metrics on spaces of differential forms and the invariant
Laplacian are discussed. See math.QA/9808015, math.QA/9808037.

{\bf 2.5.4. The radial part of the invariant Laplacian}

The radial part of the invariant Laplacian is introduced, and related
special functions are considered. See math.QA/9808037, math.QA/9808047.

{\bf 2.5.5. Boundedness and invertibility of the invariant Laplacian}

The aforementioned properties are established. See math.QA/9808037.

{\bf 2.5.6. From functions to sections of linear bundles}

Some examples of $q$-analogs of smooth sections of vector bundles.

{\bf 2.5.7. A supplement on $q$-special functions}

Some $q$-special functions and $q$-differential operators are discussed.

\medskip

{\bf 2.6. Integral representations.

\smallskip

2.6.1. Algebras of kernels of integral operators}

Some properties of kernels and the associated integral operators are
discussed. See math.QA/9511007.

{\bf 2.6.2. The Green function for the operator $\square_q$}

The Green function for the quantum invariant Laplacian is constructed. See
math.QA/9808015.

{\bf 2.6.3. Quantum cones}

A $q$-analog of the function algebra on a cone is introduced. See
math.QA/9808047.

{\bf 2.6.4. The eigenfunctions of the operator $\square_q$}

The eigenfunctions and eigenvalues of $\square_q$ are described. \\ See
math.QA/9809002.

{\bf 2.6.5. Harmonic analysis}

A $q$-analog of the non-Euclidean Fourier transform is produced. See
math.QA/9808015.

{\bf 2.6.6. The Toeplitz operators in Bergman spaces}

Toeplitz and Hankel operators in $q$-analogs of weighted Bergman spaces are
discussed. See math.QA/9809018.

{\bf 2.6.7. The covariant symbols of linear operators}

Some properties of covariant symbols are discussed. See math.QA/9809018.

{\bf 2.6.8. The Berezin transform}

Some properties of the quantum Berezin transform are discussed. See
math.QA/9809018.

{\bf 2.6.9. The Berezin quantization}

The Berezin deformation of $\operatorname{Pol}(\mathbb{C})_q$ is
considered. See math.QA/9904173.

{\bf 2.6.10. A supplement on von Neumann-Shatten classes}

Some standard facts about von Neumann-Shatten classes are expounded.

\medskip

{\bf\boldmath 2.7. On kernels of intertwining integral operators

\smallskip

2.7.1. The embedding $\operatorname{Pol}(\mathbb{C})_q\hookrightarrow
\mathbb{C}[w_0SU_{1,1}]_{q,x}$}

The canonical embedding of $\operatorname{Pol}(\mathbb{C})_q$ is
constructed and its properties are described. See math.QA/9808047.

{\bf 2.7.2. Invariant integrals}

Explicit forms for positive invariant integrals on the quantum principal
homogeneous space $\mathbb{C}[w_0SU_{1,1}]_q$ and on the quantum cone are
found. See math.QA/9808047.

{\bf 2.7.3. A proof of Lemma 2.6.6}

{\bf 2.7.4. The standard intertwining operators and the c-function of
Harish-Chandra}

A $q$-analog of Harish-Chandra's c-function is introduced. \\ See
math.QA/9808015.

\bigskip

{\bf 3. FOUNDATIONS OF THE QUANTUM THEORY OF BOUNDED SYMMETRIC DOMAINS

\medskip

3.1. A survey of known results on quantum universal enveloping algebras

\smallskip

3.1.1. The Drinfeld-Jimbo algebras}

Classical and quantum universal enveloping algebras are described.

{\bf\boldmath 3.1.2. The weight $U_q\mathfrak{g}$-modules}

The properties of weight $U_q\mathfrak{g}$-modules are discussed.

{\bf 3.1.3. Monomial bases}

The Weyl group, root systems, and Lusztig automorphisms are considered. PBW
bases in $U_q\mathfrak{g}$ are described.

{\bf 3.1.4. The universal R-matrix}

The universal R-matrix and the associated linear maps between tensor
products of weight $U_q\mathfrak{g}$-modules are considered.

{\bf 3.1.5. Verma modules}

Verma modules are introduced, some properties are discussed.

{\bf\boldmath 3.1.6. Central elements of the algebra
$U_q^\mathrm{ext}\mathfrak{g}$}

The extended algebra $U_q^\mathrm{ext}\mathfrak{g}$ is introduced. The
center of this algebra and its action on Verma modules are described.

{\bf\boldmath 3.1.7. Reduction to the case $q=1$}

A method of reducing problems that involve $U_q\mathfrak{g}$-modules to
those related to $U\mathfrak{g}$-modules.

{\bf 3.1.8. Locally finite dimensional $U_q\mathfrak{g}$-modules}

Some useful facts on finite dimensional and locally finite dimensional
$U_q\mathfrak{g}$-modules are recalled. The category of weight finite
dimensional $U_q\mathfrak{g}$-modules is introduced.

\medskip

{\bf 3.2. A survey of known results on algebras of functions on compact
quantum groups

\smallskip

3.2.1. Algebraic quantum groups}

A $q$-analog of the algebra of regular functions on an affine algebraic
group is discussed.

{\bf 3.2.2. Commutation relations}

The commutation relations between matrix elements of finite dimensional
representations of $U_q\mathfrak{g}$ are discussed.

{\bf 3.2.3. Compact quantum groups}

An involution on $U_q\mathfrak{g}$ is introduced, which provides a quantum
analog of a compact form of a complex reductive group.

{\bf\boldmath 3.2.4. An example: the algebra of functions on the quantum
$SU_2$

3.2.5. Irreducible $*$-representations of the algebra $\mathbb{C}[K]_q$ and
an invariant integral}

Irreducible $*$-representations of $\mathbb{C}[K]_q$ are described,
together with an invariant integral.

\medskip

{\bf 3.3. Quantum vector spaces and Harish-Chandra modules

\smallskip

3.3.1. Generalized Verma modules}

Some standard reductive and parabolic subalgebras of $U_q\mathfrak{g}$ are
described. Generalized Verma modules are introduced, together with standard
morphisms. See math.QA/9703005, math.QA/0606499.

{\bf\boldmath 3.3.2. The quantum vector space $\mathfrak{u}^-$}

The construction of a $q$-analog of the polynomial algebra
$\mathbb{C}[\mathfrak{u}^-]$ as the graded dual of a generalized Verma
module. See math.QA/9703005.

{\bf\boldmath 3.3.3. The $(U_q\mathfrak{g},*)$-module $*$-algebra
$\operatorname{Pol}(\mathfrak{u}^-)_q$}

An R-matrix construction of the $(U_q\mathfrak{g},*)$-module $*$-algebra of
polynomials on the quantum space $\mathfrak{u}^-$. See math.QA/9703005.

{\bf\boldmath 3.3.4. The category $C(\mathfrak{g},\mathfrak{l})_q$ and
Hermitian-symmetric pairs}

The category $C(\mathfrak{g},\mathfrak{l})_q$ is introduced, and the
Zuckerman functor is considered. See math.QA/0410463.

{\bf\boldmath 3.3.5. Quantum analogs of the subalgebras
$U\mathfrak{u}^\pm\subset U\mathfrak{g}$}

$q$-analogs of the embeddings $\mathbb{C}[\mathfrak{u}^\pm]\hookrightarrow
U\mathfrak{g}$ are produced.

{\bf\boldmath 3.3.6. A quantum analog of the decomposition
$U\mathfrak{g}=U\mathfrak{u}^-\cdot U\mathfrak{l}\cdot U\mathfrak{u}^+$}

{\bf\boldmath 3.3.7. The quadratic algebra $\mathbb{C}[\mathfrak{p}^-]_q$}

The ideal of relations in $\mathbb{C}[\mathfrak{p}^-]_q$ is described. See
math.QA/0606499.

{\bf 3.3.8. A supplement on prehomogeneous vector spaces}

The list of prehomogeneous vector spaces of commutative parabolic type is
presented.

{\bf 3.3.9. A supplement on non-commutative Gr\"obner bases}

An algorithm that reduces words to a normal form is discussed.

\medskip

{\bf 3.4. Spaces of functions on quantum bounded symmetric domains

\smallskip

3.4.1. The Fock representation and finite functions}

An extension of the Fock representation to finite functions is discussed.

{\bf 3.4.2. The Fock representation and distributions}

An extension of the Fock representation to distributions is discussed.

{\bf\boldmath 3.4.3. The action of $U_q\mathfrak{g}$ in the spaces of test
functions and distributions}

The structures of $U_q\mathfrak{g}$-module
$\operatorname{Pol}(\mathfrak{p}^-)_q$-bimodules are introduced on some
spaces of test functions and distributions.

{\bf\boldmath 3.4.4. Existence, uniqueness, and an explicit form of the
invariant integral on $\mathscr{D}(\mathbb{D})_q$}

See math.QA/0410530, math.QA/9803110.

{\bf\boldmath 3.4.5. An isomorphism of $U_q\mathfrak{g}$-modules between
$\mathscr{D}(\mathbb{D})_q$ and
$U_q\mathfrak{g}\otimes_{U_q\mathfrak{k}}\mathbb{C}$}

{\bf\boldmath 3.4.6. Boundedness of the operators of the Fock
representation of the $*$-algebra $\operatorname{Pol}(\mathfrak{p}^-)_q$}

See math.QA/0410605.

\medskip

{\bf 3.5. The canonical embedding

\smallskip

3.5.1. Introduction}

Some auxiliary facts are expounded.

{\bf 3.5.2. A quantum analog of the element $w_0$ of the Weyl group}

The notion of a quantum Weyl group is expounded.

{\bf 3.5.3. Fundamental representations and special bases}

Special bases of finite dimensional $U_q\mathfrak{g}$-modules are
described.

{\bf\boldmath 3.5.4. The canonical embedding $\mathbb{C}[\mathfrak{p}^-]_q
\hookrightarrow\mathbb{C}[X_\mathbb{S}^-]_{q,t}$}

The canonical embedding of holomorphic polynomials is constructed.

{\bf\boldmath 3.5.5. The $*$-algebra
$(\mathbb{C}[\mathbb{X}_\mathbb{S}]_q,*)$}

{\bf\boldmath 3.5.6. The canonical embedding of the $*$-algebra
$\operatorname{Pol}(\mathfrak{p}^-)_q$}

The canonical embedding of the entire polynomial algebra is constructed.

{\bf\boldmath 3.5.7. $U_q\mathfrak{k}$-invariant polynomials}

The standard generators of the algebra of $U_q\mathfrak{k}$-invariant
polynomials are selected.

{\bf\boldmath 3.5.8. Representations of the $*$-algebra
$\operatorname{Pol}(\mathfrak{p}^-)_q$}

The existence and uniqueness of the faithful irreducible representation of
$\operatorname{Pol}(\mathfrak{p}^-)_q$ by bounded operators in a Hilbert
space is proved.

{\bf\boldmath 3.5.9. A supplement on the cone over $\partial\mathbb{D}$}

{\bf 3.5.10. A supplement on spherical
$(\mathfrak{g},\mathfrak{k})$-modules and the Hua-Schmid theorem}

A sequence of strictly orthogonal roots is described. The notions of the
rank and the genus of a bounded irreducible symmetric domain are
introduced.

\medskip

{\bf 3.6. Covariant differential calculi and invariant differential
operators

\smallskip

3.6.1. First order differential calculi}

A construction of the first order differential calculus on
$\mathbb{C}[\mathfrak{p}^-]_q$. See math.QA/0606499, math.QA/9905035.

{\bf 3.6.2. Universal enveloping differential calculi}

The definition and examples of a universal enveloping differential
calculus. See math.QA/0606499.

{\bf 3.6.3. Differential forms with polynomial coefficients}

A covariant differential calculus over
$\operatorname{Pol}(\mathfrak{p}^-)_q$ is produced. See math.QA/0606499.

{\bf 3.6.4. Differential forms with generalized coefficients}

See math.QA/0606499.

{\bf 3.6.5. Invariant differential operators on $K\backslash G$}

The notion of an invariant differential operator is introduced.

{\bf 3.6.6. Invariant differential operators on $K\backslash G$, continued}

{\bf 3.6.7. Invariant differential operators on the big cell of the
generalized flag manifold}

\bigskip

{\bf 4. CONCLUDING NOTES

\medskip

4.1. Boundaries and spherical principal series}

A list of some `typical problems' is presented.

{\bf 4.2. Analytic continuation of the holomorphic discrete series and the
Penrose transform}

An outline for further research.

\bigskip

{\bf REFERENCES

\bigskip

Subject index

\bigskip

Notation index}

\newpage

\centerline{\LARGE\bf Noncommutative complex analysis}

\centerline{\Large\it Leonid Vaksman}

{\Large\bf Introduction}

Beginning with the works of Murray and von Neumann on operator algebras
\cite{Neu1}, the construction of noncommutative analogues of the most
fundamental mathematical theories has been invariably attracting the
attention of experts. Among the widely known examples are the theory of
$C^*$-algebras, operator $K$-theory, noncommutative differential geometry,
and quantum group theory. These are the noncommutative analogues of
important branches of general topology, algebraic topology, differential
geometry, and group theory, respectively.

\bigbreak

We will be interested in noncommutative complex analysis.

\medbreak

In 1934, Lavrentiev obtained an important result in the theory of
approximation in the complex domain. He proved that every continuous
function on a compact subset with empty interior, $K$, of the complex plane
$\mathbb{C}$, can be uniformly approximated by polynomials. New approaches
to this and similar theorems were found in the 50s, which lead to the birth
of the theory of uniform algebras, the latter being closely related to
complex analysis \cite{Gamelin}.

\medbreak

The first substantial results of noncommutative complex analysis were
obtained by Arveson in 1969. In his foundational work \cite{Arv}, he began
the study of noncommutative analogues of uniform algebras and, in
particular, introduced the notion of the Shilov boundary of a closed
subalgebra of a $C^*$-algebra.

\medbreak

The next important step was made in the mid 90s. Noncommutative analogues
of bounded symmetric domains were found \cite{SV} in the framework of the
theory of quantum groups \cite{Jant}, which later lead to a noncommutative
analogue of function theory in such domains.

\medbreak

Let us note that bounded symmetric domains \cite[p. 342]{Helg}, invariably
attract the attention of geometers, algebraists and analysts, because they
serve as a source of exactly solvable problems in complex analysis,
noncommutative harmonic analysis and classical mathematical physics.

\bigbreak

The simplest bounded symmetric domain is the unit disc
$\mathbf{D}=\bigl\{z\in\mathbb{C}\,\bigl\lvert\,\lvert z\rvert<1\bigr\}$.
Its quantum analogue was introduced by Klimek and Lesniewski
\cite{KlimLesn}. For pedagogical reasons, one should begin the study of the
quantum theory of bounded symmetric domains with the quantum disc, cf. the
second chapter of the book.

\medbreak

The third chapter of the book lays the foundations of the theory of quantum
bounded symmetric domains, while the draft of the expanded version of the
book, available on the same website, contains a long list of open problems
and explains the applications of the theory of quantum bounded symmetric
domains to noncommutative function theory (in particular, results are
obtained on the $q$-analogues of the Shilov boundary, weighted Bergman
spaces and Toeplitz operators in such spaces; $q$-analogues of the Bergman
and Cauchy-Szeg\"o kernels are found).

\bigbreak\bigbreak

A wide array of literature is available on quantum bounded symmetric
domains, for instance, the set of lecture notes \cite{Vaks_lectures} and
the papers of Bershtein, Kolesnik, Proskurin, Shklyarov, Sinel'schikov,
Stolin, Turowska, Vaksman, Zhang \cite{Bershtein1, BSV, BKV, Turow,
ProskTurow, ShklZhang, ShZhang_diff, VakShkl, SSV-BER, Vak01, SSV01,
SStV_KV1}. The work \cite{Vaks} is especially notable. It uses methods of
the theory of unitary dilatations of Sz\"okefalvi-Nagy and Foia\v{s}
\cite{NagyFoias} to describe the Shilov-Arveson boundary of a quantum ball,
introduced by Pusz and Woronowicz \cite{PW}. These articles are simpler
than the third chapter of the book, since they either omit the proofs, or
discuss only one of Cartan's series of irreducible bounded symmetric
domains \cite[p. 387]{Helg}. However, the results and proofs of the third
chapter of the book are applicable to all quantum bounded symmetric domains
at once.

\bigbreak\bigbreak

I am deeply grateful to my students Bershtein, Kolesnik, Korogodsky,
Shklyarov; coauthors Soibelman and Stolin; and colleagues Jakobsen,
Koelink, Klimyk, Kolb, Samoilenko, Turowska, Schmyudgen, Zhang, for
numerous helpful discussions.

\medbreak

A special role in my life was played by Vladimir Drinfeld. In the mid 80s
he taught me the basics of the theory of quantum groups and helped me
return to mathematics after an involuntary break I had to take from it for
many years. Now, after decades, Volodya's help and support are still very
important to me.

\newpage

\pagenumbering{arabic}
\thispagestyle{empty}

\vspace*{10ex}

\begin{center}\huge\bf
КВАНТОВЫЕ ОГРАНИЧЕННЫЕ СИММЕТРИЧЕСКИЕ ОБЛАСТИ
\end{center}
\vspace{5ex}

\centerline{\LARGE\bf Л. Л. Ваксман}

\bigskip \centerline{\LARGE\bf при участии Е. К. Колесника и С. Д.
Синельщикова}

\vfill

\centerline{\LARGE\bf Харьков,\ \ \ 2007}

\newpage
\noindent
\renewcommand{\contentsname}{\centerline{\large\rm ОГЛАВЛЕНИЕ}}
\par\vspace{-2.0\baselineskip}\par
\tableofcontents

\newpage

\section{ ВВЕДЕНИЕ}\label{Introduction}

{\bf 1.} В предлагаемой Вашему вниманию книге изучаются квантовые аналоги
ограниченных симметрических областей. Такие области привлекают неизменное
внимание геометров, алгебраистов и аналитиков, поскольку служат источником
точно решаемых задач комплексного анализа, некоммутативного гармонического
анализа и классической математической физики.

Простейшей ограниченной симметрической областью является единичный круг
$\mathbb{D}=\{z\in\mathbb{C}|\:|z|<1\}$. Его квантовый аналог введен
Климеком и Лесневским \cite{KlimLesn}. Из методических соображений мы
начинаем с главы \ref{quantum_disc_part} о квантовом круге, чтобы, не
отвлекая читателя алгебраическими подробностями, дать общее представление
об интересующих нас задачах некоммутативного комплексного анализа и
некоммутативного гармонического анализа. Материал этой главы может стать
основой семинара для студентов 3-4 курсов механико-математических
факультетов университетов.

Задачи, обсуждаемые в главе \ref{quantum_disc_part}, можно поставить, а
многие из них и решить, в существенно более общем контексте -- в рамках
квантовой теории ограниченных симметрических областей \cite{SV,
SV_Cont_Math}, основы которой изложены в главе \ref{basic_theory}. Тем
самым перебрасывается мост между теорией квантовых групп и некоммутативным
комплексным анализом.

Более полное представление о главах \ref{quantum_disc_part},
\ref{basic_theory} дает оглавление к книге.

В следующих главах планировалось изложить результаты Берштейн, Ваксмана,
Колесника, Проскурина, Синельщикова, Столина, Туровской, Цанга, Шклярова
\cite{Bershtein1, BSV, BKV, Turow, ProskTurow, ShklZhang, VakShkl,SSV-BER,
Vak01, SSV01, SStV_KV1}, а также неопубликованные результаты о квантовых
ограниченных симметрических областях, некоторые из которых имеются в
библиотеке xxx.lanl.gov или появятся там в недалеком будущем, см.,
например, \cite{Vaks_lectures}.

К сожалению, этот план не удается осуществить по причинам, не связанным с
математикой. В сети будет размещена страница, посвященная квантовым
ограниченным симметрическим областям, и на ней -- черновик полного варианта
несостоявшейся книги, примерно вдвое превышающий ее по объему и включающий
обсуждение нерешенных задач.

Автор глубоко признателен своим ученикам Берштейн, Колеснику, Корогодскому,
Шклярову, соавторам Сойбельману, Столину и коллегам Кёлинку, Климыку,
Колбу, Самойленко, Туровской, Шмюдгену, Цангу, Якобсену за многочисленные
полезные обсуждения результатов книги.

Особую роль в моей судьбе сыграл Владимир Дринфельд. В середине
восьмидесятых годов он обучил меня основам теории квантовых групп и тем
самым помог вернуться в математику после многолетнего вынужденного
перерыва.

\bigskip \bigskip

{\bf 2.} Немного истории. В конце семидесятых годов изучение точно решаемых
задач статистической механики и квантовой теории поля привело Фаддеева и
его сотрудников к созданию квантового метода обратной задачи рассеяния
\cite{SF,STF,TF}. Ими было введено квантовое уравнение Янга-Бакстера и
показано, что его решениям отвечают серии точно решаемых задач
математической физики. В работе \cite{KS} приводятся известные к
восьмидесятому году решения этого уравнения и отмечается, что
''вырисовываются его глубокие связи с такими разделами математики, как
теория групп и алгебраическая геометрия''. В начале восьмидесятых годов
выявлению этих связей и изучению решений квантового уравнения Янга-Бакстера
была посвящена обширная литература. Отметим работы Склянина, Кулиша,
Дринфельда и Решетихина \cite{SklFAN, KulResh, Drinf83, SklUMN,
KulReshSkl}.

В 1984 году Дринфельд ввел квантовые аналоги универсальных обертывающих
алгебр (см. \itemiii\ \ref{DrJ}), и важным событием стал его доклад о
квантовых группах на семинаре Гельфанда. В 1985 году к квантовым аналогам
универсальных обертывающих алгебр пришел Джимбо. Большую роль в становлении
теории квантовых групп сыграли статьи \cite{Drinf85,Jimbo85} и обзоры
\cite{Drinf1, JimboLMP, RTF}.

Другой подход к теории квантовых групп принадлежит Вороновичу
\cite{Woron80, Woron1}. В дальнейшем его подход использовался в теории
компактных квантовых групп \cite{Woron91, DijkhKoorn, Van_Daele95},
простейшей из которых является квантовая группа $SU(2)$ \cite{Woron86,
Podles1, VakSoib88, Masuda, Noumi_Mim1, Noumi_mim2,MasudaJ,
KoornNederl,KoornSIAM, Kir_CG, Vak89, BiedLohe1, Klimyk_Kach1,
Klimyk_Kach2, Koelink1, Koelink2, Koelink3, Koelink4, Podles2, Podles3,
LesnRin, AhmedovSU}.

В последующие годы были найдены приложения теории квантовых групп в
малоразмерной топологии и в теории категорий , а также в конформной
квантовой теории поля \cite{Reshetikhin, Turaev, Kassel, Joyal, Lyubash1,
Lyubash2, Alvares, Pasquier}. Не касаясь этих и многих других приложений,
ограничимся работами, посвященными квантовым аналогам однородных
пространств некомпактных вещественных групп Ли.

Первые результаты в этом направлении получены в конце восьмидесятых годов в
заметке \cite{KorVak89}. В этой работе введена квантовая группа движений
евклидовой плоскости, что позволило решить ряд задач некоммутативного
гармонического анализа и теории специальных функций \cite{Celeghini1,
WoronE2_LMP,WoronE2_CMP, Chakrabarty, Lukiersky, KoornSw, Koelink55,
Koelink6, KoelinkSw1, KoelinkSw2, KoelinkAssche, Swarttouw, Bonechi,
Ahmedov}. Квантовую группу движений плоскости можно получить из квантовой
группы $SU(2)$ с помощью контракции Иненю-Вигнера \cite[стр.
234]{VinbGorbOn}; этот прием применим и к некоторым другим ''неоднородным''
квантовым группам \cite{Celeghini2, Celeghini3, Dobrev}.

Другим методом, применимым к овеществленным комплексным полупростым группам
Ли, Вороновичу и Подлещу удалось получить квантовый аналог группы Лоренца
\cite{PodlWor}. Их исследования были продолжены в работах
\cite{Pusz_Lorentz, Podl-RIMS, OlshRogov1, OlshRogov2, Buff1, Buff2}.

Большие трудности ожидали на пути к квантовому аналогу группы движений
плоскости Лобачевского (точнее, к квантовому аналогу некоторой локально
изоморфной группы). Дело дошло до того, что в статье \cite{WoronAff}
Воронович отказал этой квантовой группе в праве на существование. Впрочем,
позднее, по-видимому, познакомившись с работой Корогодского \cite{Korogod},
он сменил точку зрения. Спустя еще несколько лет построение требуемой
квантовой группы в рамках аксиоматики Кастерманса-Ваеса \cite{KusVaes1,
KusVaes2} было завершено Кёлинком и Кастермансом \cite{KoelinkKustermans}.

Приведенные выше ссылки на литературу не являются исчерпывающими, но дают
общее представление о развитии теории некомпактных квантовых групп в первые
годы ее существования.

\bigskip

В работе \cite{Korogod} бросалось в глаза резкое несоответствие между
простотой классического объекта исследования и сложностью его квантового
аналога. Возник план отказаться от использования алгебр функций на
квантовых группах при изучении алгебр функций на квантовых однородных
пространствах. Основными методами исследования при таком подходе могли бы
стать методы теории представлений вещественных редуктивных групп Ли
\cite{Wallach1, Wallach2}. Введенная Хариш-Чандрой в \cite{HCh6}
голоморфная дискретная серия представлений групп эрмитова типа реализуется
во взвешенных пространствах Бергмана голоморфных функций в ограниченных
симметрических областях \cite{Kor}, а соответствующие представления алгебр
Ли -- в алгебрах полиномов на предоднородных комплексных векторных
пространствах коммутативного параболического типа (см. \itemiii\
\ref{prehomogeneous}). Естественно возник вопрос о том, имеются ли
квантовые аналоги ограниченных симметрических областей, взвешенных
пространств Бергмана, голоморфных дискретных серий и полиномов
Сато-Бернштейна предоднородных векторных пространств коммутативного
параболического типа (о полиномах Сато-Бернштейна см. \cite{Muller, Sato}).

Положительный ответ на этот вопрос должен был привести к содержательным
задачам о квантовых предоднородных векторных пространствах и о модулях
Хариш-Чандры \cite{Schmid} над квантовыми универсальными обертывающими
алгебрами. Это могло привести к прорыву в некоммутативном комплексном
анализе -- области исследований, восходящей к работе У. Арвесона
\cite{Arv}. Надежды оправдались.

\bigskip

Во второй половине девяностых годов три группы математиков получили первые
результаты квантовой теории ограниченных симметрических областей.
Исследования проводились независимо, поскольку команды не знали друг о
друге и считали используемые ими методы самодостаточными.

Танисаки и его сотрудники в работах \cite{KMT, Kamita0, Morita, Kamita3}
ввели в рассмотрение $q$-аналоги предоднородных векторных пространств
коммутативного параболического типа и нашли явный вид их полиномов
Сато-Бернштейна.

Якобсен предложил более простой способ построения этих же квантовых
векторных пространств и осознал, что идет по направлению к квантовым
эрмитовым симметрическим пространствам некомпактного типа
\cite{Jak-Hermit}, \cite{JakZhang}.

Авторы перечисленных работ упустили из виду так называемую скрытую
симметрию рассматриваемых ими квантовых предоднородных векторных
пространств \cite{SV0, ShSinVak2001} и, по-видимому, вследствие этого, не
сделали решающий шаг на пути к квантовым ограниченным симметрическим
областям.

Основы квантовой теории ограниченных симметрических областей заложены в
статье \cite{SV}. Соответствие подходов работ \cite{Tanisaki, Jak-Hermit,
SV} установлено Шкляровым \cite{Shkl}.

\newpage

\section{ КВАНТОВЫЙ КРУГ }\label{quantum_disc_part}

\subsection{Предупреждение педантам}\label{to_pedants}

В некоторых случаях, не имея возможности включить в текст используемые
определения, мы заменяем их ссылками на литературу.

Полезно сознавать, что в учебниках встречаются близкие по смыслу понятия с
одинаковыми названиями. Иногда не совсем точно используются общепринятые
обозначения. Но всегда по контексту ясно, что имеет в виду автор. Приведем
несколько примеров.

Базисом векторного пространства $V$ обычно называют подмножество $B\subset
V$ линейно независимых векторов, порождающее $V$ \cite{Lang,GlLu}. Однако,
сопоставляя линейному оператору его матрицу, базис считают линейно
упорядоченным, тем самым подменяют определение (одно множество $B\subset V$
с двумя разными отношениями линейного порядка -- это теперь два базиса).

Пополнением метрического пространства $X$ называют полное метрическое
пространство $\bar{X}$ и изометрическое вложение $i:X\to\bar{X}$, для
которого $i(X)$ плотно в $\bar{X}$ \cite{KolmFomin}. Но некоторые авторы
\cite{Naimark} вместо любой пары $(\bar{X},i)$ с этими свойствами считают
пополнением лишь одну выделенную пару, построение которой приведено во всех
учебниках по функциональному анализу.

Ассоциативность декартова произведения множеств записывают следующим
образом: $(X_1\times X_2)\times X_3=X_1\times(X_2\times X_3)$. Но, в
действительности, речь здесь идет не о равенстве, а о каноническом
изоморфизме \cite{Bou_sets}:\\
$(X_1\times X_2)\times X_3\cong X_1\times(X_2\times X_3)$.

Аналогично, вместо равенства $(X_1\otimes X_2)\otimes
X_3=X_1\otimes(X_2\otimes X_3)$ следовало бы писать $ (X_1\otimes
X_2)\otimes X_3\cong X_1\otimes(X_2\otimes X_3)$.

Кстати, как и понятие пополнения метрического пространства, понятие
тензорного произведения векторных пространств трактуется в литературе
двояко. Некоторые авторы \cite{Kassel_QG} тензорным произведением векторных
пространств $V'$, $V''$ называют векторное пространство, обозначаемое
$V'\otimes V''$, и билинейное отображение $i:V'\times V''\to V'\otimes
V''$, для которых выполнено следующее требование универсальности. Если
\hbox{$f:V'\times V''\to W$} билинейное отображение, то существует и
единственно линейное отображение $\widetilde{f}:V'\otimes V''\to W$, для
которого $f=\widetilde{f}\;i$. Большинство авторов \cite{Bou_lin, Lang,
KostrikinManin} вместо любой такой пары со свойством универсальности
называют тензорным произведением одну, выделенную пару, построение которой
приведено в \cite{Bou_lin, Lang, KostrikinManin}.

\subsection{Топология}\label{topology}

\subsubsection{Коммутативные $C^*$-алгебры.}\label{Gelfand_l}

В конце шестидесятых годов в работе \cite{Arv} У.~Арвесон предложил
программу построения некоммутативного комплексного анализа. Основным
объектом исследования в \cite{Arv} считается пара алгебр, одна из которых
является некоммутативным аналогом алгебры $C(K)$ непрерывных функций на
компакте $K$, а другая -- некоммутативным аналогом подалгебры функций,
аналитических во внутренних точках компакта $K$.

Такой подход к обычному комплексному анализу хорошо зарекомендовал себя в
рамках теории (коммутативных) равномерных алгебр, тесно связанной с
задачами равномерной аппроксимации рациональными функциями в комплексной
области \cite{Gamelin}. В \cite{Arv} Арвесон существенно использует
основные понятия теории $C^*$-алгебр. Напомним некоторые из них.

Основным полем будет поле комплексных чисел $\mathbb{C}$.

Алгебру с единицей называют также унитальной алгеброй. \IND{алгебра !
унитальная} Рассмотрим категорию, вообще говоря, неунитальных алгебр. В
этой категории морфизм $F_1\to F_2$ не обязан переводить единицу алгебры
$F_1$ в единицу алгебры $F_2$, даже при наличии единиц.

Присоединением единицы \IND{присоединение единицы} называют функтор из
рассматриваемой категории в категорию унитальных алгебр, который каждой
алгебре $F$ сопоставляет алгебру с единицей $\widetilde{F}=\mathbb{C}\times
F$:
$$
(\alpha_1,f_1)(\alpha_2,f_2)=
(\alpha_1\alpha_2,\alpha_1f_2+\alpha_2f_1+f_1f_2),\qquad\alpha_1,\alpha_2
\in\mathbb{C},\;f_1,f_2\in F.
$$
Действие этого функтора на морфизмы определяется очевидным образом.

Напомним, что алгебра $F$ называется банаховой, \IND{алгебра ! банахова}
если она является банаховым пространством и
$$\|f_1\,f_2\|\le\|f_1\|\cdot\|f_2\|,\qquad f_1,f_2\in F.$$

\begin{example}\label{compact_l}
Алгебра $C(K)$ непрерывных функций на компакте $K$, наделенная нормой
$\|f\|=\max\limits_{x\in K}|f(x)|$, является банаховой алгеброй.
\end{example}

\begin{example}\label{loc_compact_l}
Пусть $X$ -- локально компактное топологическое пространство. Рассмотрим
несвязное объединение $\widetilde{X}=X\cup\{\infty\}$ множества $X$ и
одноточечного множества $\{\infty\}$. Наделим $\widetilde{X}$ слабейшей из
топологий, содержащих дополнения в $\widetilde{X}$ ко всем компактам
$K\subset X$ и все открытые подмножества $U\subset X$. Полученное
топологическое пространство компактно. Вложение
$i:X\hookrightarrow\widetilde{X}$ является непрерывным отображением и
называется тихоновской компактификацией $X$. Тихоновская компактификация
\IND{компактификация ! тихоновская} -- это функтор, сопоставляющий локально
компактному топологическому пространству $X$ компакт $\widetilde{X}$ с
отмеченной точкой $\infty$. Алгебра $C(\widetilde{X})$ всех непрерывных
функций на компакте $\widetilde{X}$ и ее подалгебра $C_0(X)=\{f\in
C(\widetilde{X})\;|\;f(\infty)=0\}$ с унаследованной нормой
$\|f\|=\max\limits_{x \in X} f(x)=\max\limits_{x\in\widetilde{X}}f(x)$
являются банаховыми алгебрами. Отметим, забывая о нормах, что
$C(\widetilde{X})$ может быть получена присоединением единицы к $C_0(X)$.
\end{example}

\medskip

Инволюция $*$ \IND{инволюция} -- это антилинейный инволютивный
антиавтоморфизм алгебры $F$:
$$
(f_1+f_2)^*=f_1^*+f_2^*,\quad(\lambda f)^*=
\overline{\lambda}f^*,\qquad\lambda\in\mathbb{C},\;f_1,f_2,f\in F,
$$
$$**=\operatorname{id},\qquad(f_1f_2)^*=f_2^*f_1^*,\qquad f_1,f_2\in F.$$
Из единственности единицы в $F$ следует, что $1^*=1$. Алгебру с инволюцией
для краткости называют $*$-алгеброй. \IND{$*$-алгебра}

Морфизмом $*$-алгебр \IND{морфизм ! $*$-алгебр} -- это гомоморфизм алгебр
$j:F_1\rightarrow F_2$, уважающий инволюции: $*j=j*$. Как обычно, обратимый
морфизм в категории $*$-алгебр называют изоморфизмом.

Банахову алгебру $F$ с инволюцией $*$ называют $C^*$-алгеброй,
\IND{$C^*$-алгебра} если для всех элементов $f\in F$ имеет место равенство
\begin{equation}\label{C-star_l}
\|f\|^2=\|f^*f\|.
\end{equation}
Морфизмом $C^*$-алгебр \IND{морфизм ! $C^*$-алгебр} называют непрерывное
линейное отображение банаховых пространств, являющееся морфизмом
$*$-алгебр. Отметим,что требование непрерывности можно не включать в это
определение, см. предложение \ref{c_star_l}.

Банахова алгебра $C_0(X)$ с инволюцией $*:f\mapsto\overline{f}$, очевидно,
является коммутативной $C^*$-алгеброй.

Напомним, что отображение топологических пространств называют собственным,
если прообраз каждого из компактов является компактом \cite[стр.
187]{Borisovich}. Рассмотрим категорию, объектами которой служат локально
компактные топологические пространства, а морфизмами -- непрерывные
собственные отображения. Каждому локально компактному топологическому
пространству $X$ сопоставлена $C^*$-алгебра $C_0(X)$. Пойдем дальше:
каждому непрерывному собственному отображению $\varphi:X_1\to X_2$
сопоставим морфизм $C^*$-алгебр:
$$C_0(X_2)\to C_0(X_1),\qquad f(x)\mapsto f(\varphi(x)).$$

Получен контравариантный функтор из категории локально компактных
топологических пространств в категорию коммутативных $C^*$-алгебр. Как
известно \cite[стр. 70]{BrRob}, эти категории антиэквивалентны. Напомним
построение квазиобратного функтора \cite[стр. 94]{GelfandManin} --
знаменитого преобразования Гельфанда. \IND{преобразование ! Гельфанда}

Характерами коммутативной $C^*$-алгебры называют ее гомоморфизмы в
$C^*$-алгебру $\mathbb{C}$.\footnote{Требования непрерывности и
инволютивности можно не включать в определение характера \cite[стр.
70]{BrRob}.} Множество $\sigma(F)$ характеров называют спектром
коммутативной $C^*$-алгебры $F$. \IND{спектр ! коммутативной $C^*$-алгебры}
Каждому элементу $f\in F$ отвечает функция
$\hat{f}(\chi)\stackrel{\operatorname{def}}{=}\chi(f)$ на спектре
$\sigma(F)$. Наделим $\sigma(F)$ слабейшей из топологий, в которых
непрерывны все функции $\{\hat{f}|\:f\in F\}$. Можно доказать, что спектр
$\sigma(F)$ локально компактен и преобразование Гельфанда
\begin{equation}\label{gelf_transform}
F \rightarrow C_0(\sigma(F)),\qquad f\mapsto\hat{f},
\end{equation}
является изоморфизмом $C^*$-алгебр. Спектр компактен, если и только если
алгебра $F$ унитальна \cite[стр. 70]{BrRob}. Действие преобразования
Гельфанда на морфизмы определяется естественным образом: каждому морфизму
$\varphi: F_1 \rightarrow F_2$ коммутативных $C^*$-алгебр отвечает
непрерывное собственное отображение спектров
$$\sigma(F_2)\rightarrow\sigma(F_1),\qquad\chi(f)\mapsto\chi(\varphi(f)).$$

\begin{example}\label{St-Ch}(компактификация Стоуна-Чеха)
\IND{компактификация ! Стоуна-Чеха} Существует ли такое непрерывное
вложение вещественной оси в компакт, что любая ограниченная непрерывная на
оси функция допускает единственное продолжение до непрерывной функции на
этом компакте? Рассмотрим коммутативную $C^*$-алгебру $C_b(\mathbb{R})$
ограниченных и непрерывных функций на оси с нормой
$\|f\|=\sup_{x\in\mathbb{R}}|f(x)|$. Каждой точке $x_0\in\mathbb{R}$
отвечает характер, сопоставляющий функции $f\in C_b(\mathbb{R})$ число
$f(x_0)$. Нетрудно показать, что полученное отображение вещественной оси в
компакт $\sigma(C_b(\mathbb{R}))$ обладает требуемыми свойствами \cite[стр.
23]{Dunford2}.
\end{example}

\begin{example}\label{quasi-period}(квазипериодические функции)
\IND{функция ! квазипериодическая} Пусть
$\gamma\in\mathbb{R}\setminus\mathbb{Q}$. Рассмотрим наименьшую
$C^*$-алгебру $F\subset C_b(\mathbb{R})$, содержащую как все
\hbox{$2\pi\,$-периодическиме} функции, так и все
\hbox{$2\pi\gamma$-периодические} функции. Опишем спектр $\sigma(F)$ этой
коммутативной $C^*$-алгебры. Рассмотрим пару ее элементов $u_1=e^{ix}$,
$u_2=e^{ix/\gamma}$. Они унитарны
$$u_1u_1^*=u_1^*u_1=1,\qquad u_2u_2^*=u_2u_2^*=1,$$
что доставляет непрерывное отображение компакта $\sigma(F)$ в двумерный тор
\begin{equation}\label{map_l}
\sigma(F)\to\mathbb{T}^2,\qquad\chi\mapsto(\chi(u_1),\chi(u_2)).
\end{equation}
Покажем, что это отображение является гомеоморфизмом. Достаточно установить
его биективность. Инъективность вытекает из того, что $F$ является
наименьшей $C^*$-подалгеброй, содержащей элементы $u_1$, $u_2$, как следует
из теоремы Стоуна-Вейерштрасса \cite[стр. 119]{ReedSimon1}. Сюръективность
вытекает из того, что образ прямой при естественном вложении в $\sigma(F)$
и последующем отображении \eqref{map_l} содержит плотную обмотку тора
$\{(e^{ix},e^{ix/\gamma})|\:x\in\mathbb{R}\}\subset\mathbb{T}^2$.
\end{example}

\begin{example}\label{Bohr}(почти периодические функции)
\IND{функция ! почти периодическая} Рассмотрим наименьшую $C^*$-подалгебру
$AP\subset C_b(\mathbb{R})$ и содержащую множество экспоненциальных функций
$\{e^{i\lambda x}|\:\lambda\in\mathbb{R}\}$. Эти функции являются общими
собственными векторами операторов сдвига
$$
T_t:C_b(\mathbb{R})\to C_b(\mathbb{R}),\qquad T_t:f(x)\mapsto f(x+t),\quad
t\in\mathbb{R}
$$
в $C_b(\mathbb{R})$. Известная теорема Г.~Бора утверждает, что функция
$f\in C_b(\mathbb{R})$ принадлежит подалгебре $AP$, если и только если
множество $\{T_tf|\:t\in\mathbb{R}\}$ относительно компактно в
$C_b(\mathbb{R})$. Такие функции называют почти периодическими, а связанную
с $C^*$-алгеброй $AP$ компактификацию прямой называют ее боровской
компактификацией.\IND{компактификация ! боровская}
\end{example}

\medskip

Один из элементарных результатов теории гомотопий состоит в том, что
единичная окружность $\mathbb{T}$ не является ретрактом замкнутого
единичного круга $\overline{\mathbb{D}}$ \cite[стр. 143]{Borisovich}.
Другими словами, не существует непрерывного отображения
$\overline{\mathbb{D}}$ в $\mathbb{T}$, тождественного на $\mathbb{T}$.
Переведем этот результат на язык теории банаховых алгебр, что позволит в
дальнейшем получить его аналог для квантового круга.

Рассмотрим точную последовательность
\begin{equation}\label{retract}
0\to C_0(\mathbb{D})\to
C(\overline{\mathbb{D}})\stackrel{j}{\to}C(\mathbb{T})\to 0
\end{equation}
в категории коммутативных $C^*$-алгебр, где $j$ -- сужение непрерывных в
замкнутом круге функций на его границу. Говорят, что эта точная
последовательность расщепляется, если существует такой гомоморфизм
коммутативных $C^*$-алгебр $i:C(\mathbb{T})\to C(\overline{\mathbb{D}})$,
что $ji=\operatorname{id}$. Воспользуемся антиэквивалентностью категорий.
Сформулированный выше топологический результат означает, что точная
последовательность \ref{retract} не расщепляется.

\subsubsection{Элементы общей теории $C^*$-алгебр.}\label{C-star_algebras}

$C^*$-алгебры являются некоммутативными аналогами алгебр непрерывных
функций. Приведем некоторые результаты теории $C^*$-алгебр, следуя
\cite{BrRob, Takesaki} и опуская большую часть доказательств.

Из равенства \eqref{C-star_l} для $f$ и $f^*$ следует равенство норм
$\|f^*\|=\|f\|$. Действительно, $\|f\|\le\|f^*\|$, поскольку
$\|f\|^2=\|f^*f\|\le\|f^*\|\cdot\|f\|$. Противоположное неравенство
получаем заменой элемента $f$ элементом $f^*$. Значит, $C^*$-алгебра $F$ --
это банахова $*$-алгебра, для которой
\begin{equation}\label{long_def}
\|f^*\|=\|f\|,\qquad\|f^*f\|=\|f^*\|\cdot\|f\|
\end{equation}
при всех $f \in F$. Некоторые авторы определяют $C^*$-алгебру равенствами
\eqref{long_def} вместо \eqref{C-star_l}.

Встречаются как унитальные, так и неунитальные $C^*$-алгебры. В последнем
случае унитальную $*$-алгебру $\widetilde{F}=\mathbb{C}\times F$ наделяют
структурой $C^*$-алгебры так, как описано в \cite[стр. 29]{BrRob}:
\IND{$C^*$-алгебра ! унитальная}
$$\|(\alpha,f)\|=\sup\{\|\alpha f_1+ff_1\|\;|\;f_1\in F,\;\|f_1\|=1\}.$$

Спектром элемента $f$ унитальной алгебры $F$ \IND{спектр ! элемента
алгебры} называют множество $\mathrm{spec}(f)$ тех $\lambda\in\mathbb{C}$,
для которых элемент $\lambda-f$ не имеет обратного в $F$. Спектром элемента
$f$ неунитальной алгебры $F$ называют его спектр $\mathrm{spec}(f)$ в
$\widetilde{F}$ \cite[стр. 10]{Bou_spectr}. Для любого элемента $f$
банаховой алгебры $F$ существует предел
$$
\rho(f)=\lim\limits_{n\to\infty}\|f^n\|^{\frac1{n}}=
\sup\{|\lambda|\,|\:\lambda\in\mathrm{spec}(f)\}
$$
\cite[стр. 33]{BrRob}, \cite[стр. 8]{Takesaki}, называемый спектральным
радиусом элемента $f$. \IND{спектральный радиус} Если $f$ является
элементом $C^*$-алгебры, то $\|f\|=\rho(f^*f)^{\frac12}$ \cite[стр.
36]{BrRob}, и мы приходим к следующему результату об автоматической
непрерывности \cite[стр. 21]{Takesaki}.

\begin{proposition}\label{c_star_l}
Если $\pi: F_1\to F_2$ -- морфизм $*$-алгебр и $F_1$, $F_2$ являются
$C^*$-алгебрами, то
$$\|\pi(f)\|\le\|f\|,\qquad f\in F.$$
\end{proposition}

Следующее утверждение достаточно доказать для коммутативных $C^*$-алгебр и
самосопряженного элемента $f$. В этом частном случае его легко получить с
помощью преобразования Гельфанда \cite[стр. 22]{Takesaki}.

\begin{proposition}\label{iso_l}
Если $\pi:F_1\to F_2$ -- мономорфизм $*$-алгебр и $F_1$, $F_2$ являются
$C^*$-алгебрами, то $\|\pi(f)\|=\|f\|$ при всех $f\in F$.
\end{proposition}

Получаем возможность исключить требование изометричности из определения
изоморфизма $C^*$-алгебр: на банаховой $*$-алгебре либо не существует
$C^*$-норм, либо существует единственная $C^*$-норма.

Приведем пример некоммутативной $C^*$-алгебры.

\begin{example}\label{B_H_l}
Рассмотрим алгебру $\mathcal{B}(H)$ \IND{алгебра ! $\mathcal{B}(H)$} всех
ограниченных линейных операторов в гильбертовом пространстве $H$. Наделим
ее операторной нормой $\|A\|=\sup\limits_{\|v\|=1}\|Av\|$. Легко
доказывается полнота полученной нормированной алгебры. Действительно, если
последовательность операторов фундаментальна, то она ограничена, и ее
сходимость на каждом векторе вытекает из полноты гильбертова пространства.
Наделим $\mathcal{B}(H)$ инволюцией, сопоставляющей оператору $A$
сопряженный оператор $A^*$, и покажем, что $\mathcal{B}(H)$ является
$C^*$-алгеброй:
\begin{equation}\label{C-star_B}
\|A\|^2=\|A^*A\|.
\end{equation}
Для этого докажем, что
\begin{equation}\label{ex_9a}
\|A\|=\sup\limits_{\|v\|=1}|(Av,v)|,
\end{equation}
если $A$ -- ограниченный самосопряженный линейный оператор. Неочевидно лишь
неравенство $\|A\|\le\sup\limits_{\|v\|=1}|(Av,v)|$. В частном случае
двумерного гильбертова пространства оно известно из курсов линейной
алгебры, а общий случай сводится к этому частному случаю с помощью
неравенства
\begin{equation}\label{2projector}
\|A\|\le\sup\limits_{\dim L=2}\|P_LA_{|L}\|,
\end{equation}
где $L$ пробегает множество двумерных подпространств, $A|_L$ -- сужение
линейного оператора $A$ на $L$, а $P_L$ -- ортогональный проектор на
$L$.\footnote{Для доказательства \eqref{2projector} следует рассмотреть
двумерные подпространства, содержащие пары векторов $\{v,Av\}$.} Проверку
\eqref{C-star_B} завершим с помощью \eqref{ex_9a}, как в \cite[стр.
210]{ReedSimon1},
$$
\|A^*A\|=\sup\limits_{\|v\|=1}|(A^*Av,v)|=\sup\limits_{\|v\|=1}(A^*Av,v)=
\sup\limits_{\|v\|=1}\|Av\|^2=\|A\|^2.
$$
\end{example}

\bigskip

Теорема Гельфанда-Наймарка \cite[стр. 32]{BrRob} утверждает, что каждая
$C^*$-алгебра изоморфна некоторой замкнутой самосопряженной подалгебре
$C^*$-алгебры $\mathcal{B}(H)$.

Пусть $J$ -- замкнутый двусторонний идеал $C^*$-алгебры $F$. Можно
показать, что $J=J^*$ и что факторпространство $F/J$ с нормой
$$\|f+J\|=\inf\limits_{x\in J}\|f+x\|,\qquad f\in F$$
является $C^*$-алгеброй \cite[стр. 48]{BrRob}.

\begin{example}\label{Calkin-algebra}
Пусть $H$ -- бесконечномерное сепарабельное гильбертово пространство.
Линейный оператор $A$ в $H$ называют конечномерным, если $\dim(AH)<\infty$.
Такие операторы образуют двусторонний идеал в $\mathcal{B}(H)$, и его
замыкание $\mathcal{K}(H)$ совпадает с замкнутым $*$-идеалом всех
компактных операторов \cite[стр. 223]{ReedSimon1}. Как известно \cite[стр.
90]{GKr}, $\mathcal{K}(H)$ -- единственный нетривиальный замкнутый
двусторонний идеал $C^*$-алгебры $\mathcal{B}(H)$. $C^*$-алгебра
$\mathcal{B}(H)/\mathcal{K}(H)$ называется алгеброй Калкина. \IND{алгебра !
Калкина}
\end{example}

\bigskip

Приведем результаты теории операторов, тесно связанные с алгеброй Калкина.
Существенным спектром $\mathrm{spec}_e(A)$ \IND{спектр ! существенный}
ограниченного линейного оператора $A\in\mathcal{B}(H)$ называют спектр
отвечающего ему элемента алгебры Калкина. Линейный оператор
$A\in\mathcal{B}(H)$ называют фредгольмовым, \IND{фредгольмов оператор}
если отвечающий ему элемент алгебры Калкина обратим. Согласно теореме
Аткинсона \cite[стр. 97]{Halmos_problems}, это означает, что ядро оператора
$A$ конечномерно, а образ замкнут и коразмерность образа конечна. Очевидно,
существенный спектр образуют те $\lambda\in\mathbb{C}$, для которых
оператор $\lambda-A$ не является фредгольмовым.

Индекс фредгольмова оператора \IND{фредгольмов оператор ! индекс}
$\operatorname{ind(A)}\stackrel{\operatorname{def}}{=}
\dim\operatorname{ker}A-\operatorname{codim}\operatorname{im}A$ не
изменяется при компактных возмущениях и, следовательно, является функцией
на группе обратимых элементов алгебры Калкина \cite[стр. 300]{Kato},
\cite[стр. 138]{Douglas}. Эта функция непрерывна, и ее множества уровня
линейно связны \cite[стр. 136]{Douglas}.

Простейшим примером фредгольмова оператора с ненулевым индексом служит
оператор одностороннего сдвига $S$, определяемый в ортонормированном базисе
$\{e_j\}_{j\in\mathbb{Z}_+}$ равенством $S:e_j\mapsto e_{j+1}$. Оператор
$S$ не допускает разложения в сумму нормального и компактного операторов,
поскольку индекс нормального оператора равен нулю. Это следует, из
спектральной теоремы \cite[стр. 50]{Dunford2}.

\subsubsection{Фоковское представление.}\label{Fock_l}

Напомним \cite[стр. 66]{KolmFomin}, что пополнением метрического
пространства $X$ называют его вложение $X\hookrightarrow\overline{X}$ со
всюду плотным образом в полное метрическое пространство $\overline{X}$.
Пополнение существует и в существенном единственно.

Пусть $\mathbb{D}=\{z\in\mathbb{C}|\:|z|<1\}$. Рассмотрим коммутативную
$C^*$-алгебру $C(\overline{\mathbb{D}})$ непрерывных функций в замкнутом
единичном круге $\overline{\mathbb{D}}$ и $*$-алгебру полиномиальных
функций $\operatorname{Pol}(\mathbb{C})$ на овеществленной комплексной
плоскости $\mathbb{C}\simeq \mathbb{R}^2$. Согласно теореме Вейерштрасса о
полиномиальной аппроксимации \cite[стр. 120]{ReedSimon1}, $*$-алгебра
$\operatorname{Pol}(\mathbb{C})$ плотна в $C(\overline{\mathbb{D}})$.
Другими словами, естественное вложение $\operatorname{Pol}(\mathbb{C})
\hookrightarrow C(\overline{\mathbb{D}})$ является пополнением $*$-алгебры
$\operatorname{Pol}(\mathbb{C})$ по соответствующей норме. Во многих
вопросах функционального анализа выбор пополнения несуществен, что
позволяет заменить конкретную $C^*$-алгебру $C(\overline{\mathbb{D}})$
любым пополнением $\operatorname{Pol}(\mathbb{C})$ по этой норме.
Воспользуемся этим подходом для построения некоммутативной $C^*$-алгебры --
квантового аналога коммутативной $C^*$-алгебры $C(\overline{\mathbb{D}})$.

Выберем число $q\in(0,1)$, которое будет параметром деформации: в пределе
$q\to 1$ формулы теории функций в квантовом круге будут переходить в
обычные формулы теории функций в круге. К обсуждению роли этого параметра и
процедуры квантования мы в дальнейшем не раз вернемся.

Начнем с квантовой алгебры полиномов. \IND{квантовая ! алгебра полиномов}
Рассмотрим свободную алгебру $\mathbb{C}\langle z,z^*\rangle$ с двумя
образующими $z$, $z^*$. Ее обычно отождествляют с алгеброй формальных
линейных комбинаций слов над алфавитом $\{z,z^*\}$. Очевидно, существует и
единственна инволюция $*$ в $\mathbb{C}\langle z,z^*\rangle$, отображающая
$z$ в $z^*$. Элемент $z^*z-q^2z\,z^*-1+q^2\in\mathbb{C}\langle
z,z^*\rangle$ самосопряжен. Значит, порожденный им двусторонний идеал также
является самосопряженным. Следовательно, алгебра
$\operatorname{Pol}(\mathbb{C})_q$ с образующими $z$, $z^*$ и определяющим
соотношением
\begin{equation}\label{Pol_l}
z^*z=q^2zz^*+1-q^2
\end{equation}
наследует инволюцию и становится $*$-алгеброй. Эта $*$-алгебра является
$q$-аналогом алгебры полиномов переменных $z$ и $\overline{z}$.

Для получения алгебры непрерывных функций в квантовом круге осталось
пополнить $\operatorname{Pol}(\mathbb{C})_q$. Как выбрать норму? Будем
следовать традициям квантовой механики. Подстановка $q=e^{-\frac{h}2}$ и
формальный предельный переход $h\to 0$
$$
\{z^*,z\}\stackrel{\rm{def}}{=}\lim\limits_{h\to
0}\frac{z^*z-zz^*}{ih}=\lim\limits_{h\to
0}\frac{(1-e^{-h})(1-zz^*)}{ih}=\frac1{i}\,(1-|z|^2),
$$
называемый квазиклассическим предельным переходом, наделяют плоскость
структурой пуассонова многообразия \cite[стр. 422]{Arnold}, а дополнение к
единичной окружности -- структурой симплектического многообразия \cite[стр.
175]{Arnold} с симплектической формой
\begin{equation}\label{omega_l}
\omega=\frac2{1-|z|^2}\;d\,\operatorname{Im}z\wedge d\operatorname{Re}z.
\end{equation}

Единичный круг является одной из двух компонент связности дополнения к
окружности -- его {\bf ограниченной} компонентой связности. Опыт квантовой
механики \cite{Ber_CMP} указывает, что нужно найти точное неприводимое
представление $*$-алгебры $\operatorname{Pol}(\mathbb{C})_q$ {\bf
ограниченными} операторами\footnote{Точность означает, что ненулевым
элементам алгебры отвечают ненулевые операторы представления, а
неприводимость -- отсутствие нетривиальных общих инвариантных
подпространств этих операторов.} и нормой элемента
$f\in\operatorname{Pol}(\mathbb{C})_q$ назвать норму отвечающего ему
оператора представления. Для осуществления этого плана важно знать, что
точное неприводимое $*$-представление существует и в существенном
единственно.

\begin{proposition}\label{unit_equiv}
Унитальная $*$-алгебра $\operatorname{Pol}(\mathbb{C})_q$ обладает
единственным с точностью до унитарной
эквивалентности\footnote{Представления $\pi_1$, $\pi_2$ алгебры $F$ в
гильбертовых пространствах $H_1$, $H_2$ называют унитарно эквивалентными,
если $ U\pi_1(f)=\pi_2(f)U$ для некоторого унитарного оператора $U: H_1 \to
H_2$ и всех $f\in F$.} точным неприводимым представлением ограниченными
операторами в гильбертовом пространстве.
\end{proposition}

{\bf Доказательство.} Рассмотрим предгильбертово пространство $\mathcal{H}$
с ортонормированным базисом $\{e_n\}_{n\in\mathbb{Z}_+}$ и представление
$T_F$ $*$-алгебры $\operatorname{Pol}(\mathbb{C})_q$ в $\mathcal{H}$,
определяемое равенствами
\begin{equation}\label{T_F}
T_F(z)e_n=\sqrt{1-q^{2(n+1)}}\,e_{n+1},\qquad T_F(z^*)e_n=
\begin{cases}
\sqrt{1-q^{2n}}\,e_{n-1}, & n\in\mathbb{N},
\\ 0, & n=0.
\end{cases}
\end{equation}

Пусть
\begin{equation}\label{y_l}
y=1-zz^*=q^{-2}(1-z^*z).
\end{equation}
Последнее равенство следует из $\eqref{Pol_l}$ и влечет коммутационные
соотношения
\begin{equation}\label{comm_disc}
z^*y=q^2yz^*,\qquad zy=q^{-2}yz.
\end{equation}
Значит, каждый элемент $f$ алгебры $\operatorname{Pol}(\mathbb{C})_q$
допускает разложение
\begin{equation}\label{expansion_l}
f=\sum_{j\in\mathbb{N}}z^j\psi_j(y)+\psi_0(y)+\sum_{j\in\mathbb{N}}
\psi_{-j}(y)z^{*j},
\end{equation}
где $\psi_j$ -- полиномы одной переменной. Единственность разложения
\eqref{expansion_l} и точность представления $T_F$ нетрудно доказать с
помощью равенства
$$T_F(p(y))e_n=p(q^{2n})e_n,\qquad n\in\mathbb{Z}_+.$$
Именно, при $j\ge 0$ полиномы $f_j$ легко восстанавливаются по матричным
элементам $(T_F(f)e_n,e_{n+j})$ оператора $T_F(f)$, а полиномы $f_{-j}$ --
по матричным элементам $(T_F(f)e_{n+j},e_n)$ этого оператора.

Операторы представления $T_F$ ограничены, поскольку $\|T_F(z)\|\le 1$,
$\|T_F(z^*)\|\le 1$. Следовательно, эти операторы допускают продолжение по
непрерывности с предгильбертова пространства $\mathcal{H}$ на его
пополнение $H$. Получаем точное представление $*$-алгебры
$\operatorname{Pol}(\mathbb{C})_q$ в гильбертовом пространстве, которое
будет называться фоковским и обозначаться $T_F$, как и его сужение на
$\mathcal{H}$. \IND{фоковское представление}

Докажем неприводимость представления $T_F$ в гильбертовом пространстве $H$.
Пусть $e_{jk}$ -- линейные операторы в $H$, определяемые равенством
$e_{jk}e_m=\delta_{k,m}e_j$, где $j,k,m\in\mathbb{Z}_+$, а $\delta_{k,m}$
-- символ Кронекера. Их называют матричными единицами. Рассмотрим
$C^*$-алгебру $\mathcal{B}(H)$ всех ограниченных линейных операторов в $H$
и замыкание $\mathcal{A}$ алгебры операторов представления $T_F$ в
$\mathcal{B}(H)$. Для доказательства неприводимости этого представления
достаточно установить, что $e_{jk}\in\mathcal{A}$ при всех
$j,k\in\mathbb{Z}_+$.

Воспользуемся функциональным исчислением для непрерывных функций \cite[стр.
247]{ReedSimon1} применительно к ограниченному самосопряженному линейному
оператору $T_F(y)$. Пусть $j\in\mathbb{Z}_+$ и $\psi_j(\lambda)$ --
непрерывная функция на оси, равная $1$ при $\lambda=q^{2j}$ и нулю во всех
остальных точках геометрической прогрессии $q^{2\mathbb{Z}_+}$. Тогда
$e_{jj}=\psi_j(T_F(y))$. Следовательно, диагональные матричные единицы
$e_{jj}$, $j\in\mathbb{Z}_+$, принадлежат $C^*$-алгебре $\mathcal{A}$.
Остальные матричные единицы принадлежат $\mathcal{A}$ потому, что, как
следует из \eqref{T_F},
$$
T_F(z)e_{jk}=\left(1-q^{2(j+1)}\right)^\frac12e_{(j+1)k},\qquad
e_{kj}T_F(z^*)=\left(1-q^{2(j+1)}\right)^\frac12 e_{k(j+1)}.
$$

Остается показать, что всякое точное неприводимое представление $T'$
$*$-алгебры $\operatorname{Pol}(\mathbb{C})_q$ ограниченными операторами в
гильбертовом пространстве $H'$ унитарно эквивалентно фоковскому
представлению $T_F$.

\begin{lemma}\label{spectr_y_l}
Пусть $T'$ -- неприводимое представление $*$-алгебры
$\operatorname{Pol}(\mathbb{C})_q$ ограниченными операторами в гильбертовом
пространстве $H'$ и $T'(y)\ne 0$. Тогда для некоторого вещественного числа
$a\ne 0$
\begin{equation}\label{progr_l}
\mathrm{spec}(T'(y))\subset\{0\}\cup aq^{2\mathbb{Z}}.
\end{equation}
\end{lemma}

\medskip

Доказательство приведено в следующем \itemiiiе.

\bigskip

Как вытекает из этой леммы, ненулевые точки спектра оператора $T'(y)$
принадлежат геометрической прогрессии $aq^{2\mathbb{Z}}$, $a\ne 0$, и,
следовательно, являются его собственными значениями. Пусть $\lambda_0$ --
наибольшее по модулю собственное значение и $v_0$ -- отвечающий ему
нормированный собственный вектор: $T'(y)v_0=\lambda_0v_0$, $\|v_0\|=1$. Из
\eqref{comm_disc} следует, что $T'(y)T'(z^*)v_0=q^{-2}\lambda_0T'(z^*)v_0$.
Но число $q^{-2}\lambda_0$ не является собственным значением оператора
$T'(y)$. Значит, $T'(z^*)v_0=0$. Следовательно,
\begin{equation}\label{vacuum_l}
T'(y)v_0=T'(1-zz^*)v_0=T'(1)v_0=v_0,
\end{equation}
и $\lambda_0=1$.

Аналогично, с использованием \eqref{comm_disc} доказывается, что векторы
$v_n=T'(z^n)v_0$ удовлетворяют равенствам $T'(y)v_n=q^{2n}v_n$ и
\begin{multline*}
\|v_{n+1}\|^2=\|T'(z)v_n\|^2=(T'(z^*z)v_n,v_n)=(T'(1-q^2y)v_n,v_n)=
\\ =(1-q^{2(n+1)})\|v_n\|^2.
\end{multline*}
Значит векторы $v_n\ne 0$ являются собственными векторами ограниченного
самосопряженного оператора $T'(y)$, отвечающими попарно различным
собственным значениям. Следовательно, они попарно ортогональны. Нетрудно
показать, используя равенства
$$
T'(z)v_n=v_{n+1},\qquad T'(z^*)v_{n+1}=(1-q^{2(n+1)})v_n,\qquad
T'(z^*)v_0=0,
$$
что отображение
$$H_F\to H',\qquad e_n\mapsto\frac{e'_n}{\|e'_n\|},\quad n\in\mathbb{Z}_+$$
изометрично и сплетает представления (является морфизмом). Оно сюръективно
и, следовательно, осуществляет унитарную эквивалентность рассматриваемых
представлений, так как представление $T'$ неприводимо. \hfill $\square$

\begin{remark}\label{formulas}
Фоковское представление было введено равенствами \eqref{T_F}, происхождение
которых не обсуждалось. На самом деле эти формулы можно получить, рассуждая
так, как в заключительной части доказательства предложения
\ref{unit_equiv}.
\end{remark}

\medskip

Рассмотрим одномерные представления $*$-алгебры
$\operatorname{Pol}(\mathbb{C})_q$, параметризуемые точками окружности:
\begin{equation}\label{one_l}
\rho_\varphi(z)=e^{i\varphi},\quad\rho_\varphi(z^*)=e^{-i\varphi},\quad
\qquad\varphi\in\mathbb{R}/(2\pi\mathbb{Z}).
\end{equation}
Эти представления попарно неэквивалентны.

\begin{corollary}\label{list_l}
Пусть $T$ -- неприводимое представление $*$-алгебры
$\operatorname{Pol}(\mathbb{C})_q$ ограниченными операторами в гильбертовом
пространстве. Представление $T$ унитарно эквивалентно одному из
представлений \eqref{one_l}, если $T(y)=0$, и унитарно эквивалентно
фоковскому представлению $T_F$, если $T(y)\ne 0$.
\end{corollary}

Таким образом, получен полный список неприводимых представлений $*$-алгебры
$\operatorname{Pol}(\mathbb{C})_q$, рассматриваемых с точностью до
унитарной эквивалентности.

\subsubsection{Доказательство леммы \ref{spectr_y_l}.}\label{spectr_proof}

Предположим, что утверждение леммы неверно. Так как $T'(y)\ne\{0\}$, то
спектр оператора $T'(y)$ пересекается с двумя геометрическими прогрессиями:
$$
\mathrm{spec}(T'(y))\cap
a_1q^{2\mathbb{Z}}\ne\varnothing,\qquad\mathrm{spec}(T'(y))\cap
a_2q^{2\mathbb{Z}}\ne\varnothing,
$$
где $a_1q^{2\mathbb{Z}}\ne a_2q^{2\mathbb{Z}}$.

Рассмотрим две непрерывные функции $\Psi_1(\lambda)$, $\Psi_2(\lambda)$,
инвариантные относительно гомотетии с коэффициентом $q^2$:
$$
\Psi_1(q^2\lambda)=\Psi_1(\lambda),\qquad\Psi_2(q^2\lambda)=\Psi_2(\lambda),
$$
которые принимают ненулевые значения на ''своих'' геометрических
прогрессиях
$$
\Psi_1(a_1 q^{2j})\ne 0,\quad\Psi_2(a_2q^{2j})\ne 0,\qquad j\in\mathbb{Z},
$$
и произведение которых равно $0$. Подразумевается, что $a_1, a_2 \ne 0$.
Так как
\begin{equation}\label{quasi_comm_l}
T'(z^*)T'(y)=T'(z^*)T'(1-zz^*)=T'(1-z^*z)T'(z^*)=q^2T'(y)T'(z^*),
\end{equation}
то для любого полинома $p$ одной переменной
$$T'(z^*)\,p(T'(y))=p(q^2 T'(y))\,T'(z^*).$$
Значит, для любой непрерывной функции $\Psi$ на вещественной оси
$$
T'(z^*)\,\Psi_j(T'(y))=\Psi_j(q^2T'(y))\,T'(z^*)=\Psi_j(T'(y))\,T'(z^*),
\qquad j=1,2.
$$
Линейные операторы $A_1=\Psi_1(T'(y))$, $A_2=\Psi_2(T'(y))$ коммутируют со
всеми операторами представления $T'$ и отличны от нуля, поскольку по
теореме об отображении спектров $\Psi_1(a_1)\in\mathrm{spec}(A_1)$,
$\Psi_2(a_2)\in\mathrm{spec}(A_2)$, и $\Psi_1(a_1)\ne 0$, $\Psi_2(a_2)\ne
0$. Значит, замыкания образов операторов $A_1$, $A_2$ являются ненулевыми
общими инвариантными подпространствами операторов представления $T'$. Но
представление $T'$ неприводимо. Следовательно, образы операторов $A_1$,
$A_2$ плотны в $H'$, что невозможно, поскольку
$A_1A_2=\Psi_1(T'(y))\Psi_2(T'(y))=\Psi_1\Psi_2\,(T'(y))=0$. \hfill
$\square$

\subsubsection{Непрерывные функции в квантовом круге.}\label{cont-disc}

Из точности фоковского представления $T_F$ следует, что функционал
$$
\|f\|\stackrel{\operatorname{def}}{=}\|T_F(f)\|,\qquad
f\in\operatorname{Pol}(\mathbb{C})_q,
$$
наделяет $\operatorname{Pol}(\mathbb{C})_q$ структурой нормированной
$*$-алгебры. Как объяснялось в \itemiiе \ref{Fock_l}, ее пополнение
\begin{equation}\label{compl_l}
\operatorname{Pol}(\mathbb{C})_q\hookrightarrow C(\overline{\mathbb{D}})_q
\end{equation}
является $q$-аналогом вложения
$\operatorname{Pol}(\mathbb{C})\hookrightarrow C(\overline{\mathbb{D}})$.
Следуя \cite{NagyNica, KlimLesn}, будем называть $C^*$-алгебру
$C(\overline{\mathbb{D}})_q$ алгеброй непрерывных функций в квантовом
круге. \IND{квантовый круг}

\bigskip

Перейдем к сравнению классического и квантового кругов.

Прежде всего, сравним спектр элемента $1-|z|^2\in C(\overline{\mathbb{D}})$
со спектром элемента $1-zz^* \in C(\overline{\mathbb{D}})_q$. Первый из
этих спектров является множеством значений непрерывной функции $1-|z|^2$ в
замкнутом единичном круге, то есть отрезком $[0,1]$. Второй спектр -- это
замыкание геометрической прогрессии $q^{2\mathbb{Z}_+}$, как следует из
определения фоковского представления и из изометричности вложения
$C(\overline{\mathbb{D}})_q \hookrightarrow\mathcal{B}(H)$. При возрастании
$q\in(0,1)$ геометрическая прогрессия становится все гуще и в пределе $q\to
1$ заполняет отрезок $[0,1]$, чем обеспечивается соответствие классической
и квантовой задач о спектре ''физической величины'' $y=1-zz^*$.

\bigskip

Спектры элементов $z$, $z^*$ как в классическом, так и в квантовом случае
совпадают с $\overline{\mathbb{D}}$. Действительно, $\|T_F(z)\|\le 1$ и все
числа $\lambda\in\mathbb{D}$ являются собственными значениями линейного
оператора $T_F(z^*)$.

\bigskip

Перейдем к некоммутативному аналогу обсуждавшегося в \itemiiе
\ref{Gelfand_l} результата о нерасщепляемости точной последовательности
\begin{equation}\label{new_l}
0\to C_0(\mathbb{D})\to
C(\overline{\mathbb{D}})\stackrel{j}{\to}C(\mathbb{T})\to 0,
\end{equation}
то есть о том, что окружность не является ретрактом круга.

Введем обозначение $C_0(\mathbb{D})_q$ для замкнутого двустороннего идеала
$C^*$-алгебры $C(\overline{\mathbb{D}})_q$, порожденного элементом $y$.

\begin{proposition}\label{to_standard}
Продолжение фоковского представления по непрерывности на всю алгебру
$C(\overline{\mathbb{D}})_q$ доставляет изоморфизм
\begin{equation}\label{iso_1_l}
C_0(\mathbb{D})_q\cong\mathcal{K}(H).
\end{equation}
\end{proposition}

{\bf Доказательство.} Как было установлено при доказательстве
неприводимости фоковского представления, образ $C^*$-алгебры
$C_0(\mathbb{D})_q$ в $\mathcal{B}(H)$ содержит все ''матричные единицы''
$e_{jk}$ и, в частности, содержит ортогональный проектор на вакуумное
подпространство $\mathbb{C}e_0$. Из неприводимости представления $T_F$
вытекает, что эта алгебра содержит все линейные операторы ранга 1. Значит,
образ $C^*$-алгебры $C_0(\mathbb{D})_q$ содержит все операторы конечного
ранга и, следовательно, все компактные операторы в $H$. С другой стороны,
оператор $T_F(y)$ компактен. \hfill $\square$

\bigskip

Сохраним обозначение $T_F$ для продолжения по непрерывности фоковского
представления на $C(\overline{\mathbb{D}})_q$.

\begin{proposition}\label{q-boundary}
$C(\overline{\mathbb{D}})_q/C_0(\mathbb{D})_q\cong C(\mathbb{T})$.
\end{proposition}

{\bf Доказательство.} Коммутативная $C^*$-алгебра
$C(\overline{\mathbb{D}})_q/C_0(\mathbb{D})_q$ естественно изоморфна
$C(\sigma_e(T_F(z)))$. Но $\sigma_e(T_F(z))=\mathbb{T}$. Действительно,
$\sigma_e(T_F(z)) \subset \mathbb{T}$ в силу унитарности элемента алгебры
Калкина, отвечающего оператору $T_F(z)$. Остается воспользоваться
непустотой множества $\sigma_e(T_F(z))$ и его инвариантностью относительно
вращений в $\mathbb{C}$. Инвариантность следует из того, что линейные
операторы $e^{i\varphi}T_F(z)$, где $\varphi \in \mathbb{R}$, попарно
унитарно эквивалентны. \hfill $\square$

\bigskip

Алгеброй Теплица \IND{алгебра ! Теплица} $\mathcal{T}$ называют
$C^*$-алгебру линейных операторов в $H$, порожденную идеалом
$\mathcal{K}(H)$ всех компактных операторов и изометрическим оператором
одностороннего сдвига $S$ (см. \itemiii\ \ref{C-star_algebras}).
Предложение \ref{to_standard} доставляет изоморфизм $C^*$-алгебр
$T_F:C(\overline{\mathbb{D}})_q\stackrel{\approx}{\to}\mathcal{T}$,
поскольку $T_F(z)-S\in\mathcal{K}$. Как следует из предложения
\ref{q-boundary}, $\mathcal{T}/\mathcal{K}(H)\cong C(\mathbb{T})$. Это
хорошо известный результат, см. \cite[стр. 68]{Wegge-Olsen}.

\bigskip

Мы получили точную последовательность $C^*$-алгебр
\begin{equation}\label{exact_seq_first}
0\to C_0(\mathbb{D})_q\to
C(\overline{\mathbb{D}})_q\stackrel{j_q}{\to}C(\mathbb{T}) \to 0,
\end{equation}
являющуюся $q$-аналогом точной последовательности \eqref{new_l} и
изоморфную точной последовательности
\begin{equation}\label{noncomm_exact}
0\to\mathcal{K}(H)\to\mathcal{T}\stackrel{j}{\to}C(\mathbb{T})\to 0.
\end{equation}

\bigskip

Нетрудно показать \cite[стр. 68]{Wegge-Olsen}, что точная
последовательность \eqref{noncomm_exact} не расщепляется. Действительно,
пусть $i: C(\mathbb{T})\to\mathcal{T}$ -- расщепляющий ее морфизм
$C^*$-алгебр: $ji=\operatorname{id}$. Образ алгебры $C(\mathbb{T})$ в
$\mathcal{T}$ состоит из нормальных операторов. Значит, любой оператор из
алгебры $\mathcal{T}$ допускает разложение в сумму нормального и
компактного операторов. Следовательно, индекс любого оператора из
$\mathcal{T}$ равен нулю. Приходим к противоречию, поскольку индекс
оператора $S$ отличен от нуля.

\bigskip

Следовательно, точная последовательность \eqref{exact_seq_first} также
нерасщепляема. Это и есть $q$-аналог результата о том, что окружность не
является ретрактом круга.

\begin{remark}
Гомоморфизм $j_q$ является $q$-аналогом линейного оператора сужения
непрерывных в круге $\overline{\mathbb{D}}$ функций на его границу
$\mathbb{T}$. Как видно из доказательства предложения \ref{q-boundary},
$j_q$ отображает элемент $z\in C(\overline{\mathbb{D}})_q$ в функцию
$e^{i\varphi}$ на окружности $\mathbb{T}\cong\mathbb{R}/(2\pi\mathbb{Z})$.
Из равенства $\|j_q\|=1$ вытекает неравенство
$$
|\rho_\varphi(f)|\le\|T_F(f)\|,\qquad f\in\operatorname{Pol}(\mathbb{C})_q.
$$
Участвующими в этом неравенстве представлениями исчерпывается множество
классов унитарно эквивалентных неприводимых представлений $*$-алгебры
$\operatorname{Pol}(\mathbb{C})_q$. Следовательно,
$C(\overline{\mathbb{D}})_q$ является универсальной обертывающей
$C^*$-алгеброй $*$-алгебры $\operatorname{Pol}(\mathbb{C})_q$ \cite[стр.
81]{Bou_spectr}.
\end{remark}

\subsubsection{Голоморфные функции в квантовом круге.}\label{A_q}
Начнем с полиномов. Пусть
$\mathbb{C}[z]_q\subset\operatorname{Pol}(\mathbb{C})_q$ -- подалгебра,
порожденная элементом $z$, а
$\mathbb{C}[\overline{z}]_q\subset\operatorname{Pol}(\mathbb{C})_q$ --
подалгебра, порожденная элементом $z^*$. Множество
$\{z^j\}_{j\in\mathbb{Z}_+}$ является базисом векторного пространства
$\mathbb{C}[z]_q$, а множество $\{z^{*j}\}_{j\in\mathbb{Z}_+}$ -- базисом
векторного пространства $\mathbb{C}[\overline{z}]_q$. Действительно,
линейная независимость мономов $z^j$ вытекает из линейной независимости
векторов $T_F(z^j)e_0$ и влечет линейную независимость мономов $z^{*j}$.

В дальнейшем предполагаются известными основные понятия полилинейной
алгебры, см., например, \cite{KostrikinManin}.
\begin{proposition}\label{hol-antihol}
Линейное отображение
$$
m:\mathbb{C}[z]_q\otimes\mathbb{C}[\overline{z}]_q\to
\operatorname{Pol}(\mathbb{C})_q,\qquad m:f_1\otimes f_2\mapsto f_1f_2
$$
биективно.
\end{proposition}

{\bf Доказательство.} Требуемое утверждение означает, что множество
$\{z^jz^{*k}\}_{j,k\in\mathbb{Z}_+}$ является базисом векторного
пространства $\operatorname{Pol}(\mathbb{C})_q$. Полнота этого множества
векторов следует из определения алгебры $\operatorname{Pol}(\mathbb{C})_q$,
и остается доказать линейную независимость.

Рассмотрим фоковское представление $T_F$ $*$-алгебры
$\operatorname{Pol}(\mathbb{C})_q$ в предгильбертовом пространстве
$\mathcal{H}$. предположим, что векторы $z^jz^{*k}$ линейно зависимы. Тогда
при некоторых $j_0,k_0\in\mathbb{Z}_+$
$$
\sum\limits_{j+k\ge j_0+k_0}c_{j,k}\,T_F(z^jz^{*k})=0,\qquad c_{j_0k_0}\ne
0,
$$
где подразумевается, что отлично от нуля лишь конечное число комплексных
коэффициентов $c_{j,k}$. Так как
$$
(T_F(z^jz^{*k})e_{k_0},e_{j_0})=(T_F(z^{*k})e_{k_0},T_F(z^{*j})e_{j_0})=
\begin{cases}
0, & k>k_0\;\text{или}\;j>j_0,
\\ \ne 0, & k=k_0\; \& \; j=j_0
\end{cases},
$$
то $c_{j_0,k_0}=0$. \hfill $\square$

\bigskip

Пусть $A(\overline{\mathbb{D}})$ -- банахова алгебра всех функций,
непрерывных в $\overline{\mathbb{D}}$ и голоморфных в $\mathbb{D}$,
изометрически вложенная в $C(\overline{\mathbb{D}})$. Рассмотрим замыкание
$A(\overline{\mathbb{D}})_q$ алгебры $\mathbb{C}[z]_q$ в
$C(\overline{\mathbb{D}})_q$. Пара
$(A(\overline{\mathbb{D}})_q,C(\overline{\mathbb{D}})_q)$ является
$q$-аналогом пары $(A(\overline{\mathbb{D}}),C(\overline{\mathbb{D}}))$.

Как известно из курсов комплексного анализа, максимум модуля непрерывной в
$\overline{\mathbb{D}}$ и аналитической в $\mathbb{D}$ функции достигается
на границе $\mathbb{T}$ круга $\mathbb{D}$. Это означает, что сужение
гомоморфизма
\begin{equation*}
j:\;C(\overline{\mathbb{D}})\to C(\mathbb{T}),\qquad j:\;f\mapsto
f|_\mathbb{T}
\end{equation*}
на подалгебру $A(\overline{\mathbb{D}})$ является изометрическим линейным
оператором. Докажем аналогичное утверждение для квантового круга.

Воспользуемся гомоморфизмом $j_q:C(\overline{\mathbb{D}})_q\to
C(\mathbb{T})$, введенным в \itemiiiе \ref{cont-disc}.

\begin{proposition}\label{max_disc}
Сужение гомоморфизма $j_q$ на подалгебру $A(\overline{\mathbb{D}})_q$
является изометрическим линейным оператором.
\end{proposition}

{\bf Доказательство.} Из неравенства $\|j_q\|\le 1$ следует, что
рассматриваемое ограничение гомоморфизма $j_q$ на подалгебру
$A(\overline{\mathbb{D}})_q$ является сжатием. Остается показать, что
$\|f\|\le\|j_q(f)\|$ для любого $f\in\mathbb{C}[z]_q$. Для этого достаточно
воспользоваться неравенством Неймана (предложением \ref{Neumann_ineq}) и
тем, что $\|T_F(z)\|\le 1$:
$$
\|f\|=\|T_F(f)\|\le\max\limits_{|z|=1}|f(z)|=
\max\limits_{\varphi\in\mathbb{R}/(2\pi\mathbb{Z})}\|\rho_\varphi(f)\|=
\|j_q(f)\|. \eqno \square$$

\begin{remark}
 Из доказательства предложения \ref{Neumann_ineq} следует, что приведенный
выше результат допускает усиление: ограничение гомоморфизма $j_q$ на
подалгебру $A(\overline{\mathbb{D}})_q$ является {\bf вполне}
изометрическим линейным оператором (определение см., например, в \cite[стр.
148]{Arv}).
\end{remark}

\subsubsection{Дополнение о спектральных множествах.}\label{Nagy-Foias}
Напомним, что спектром \IND{спектр ! ограниченного линейного оператора}
ограниченного линейного оператора $A$ называют множество тех
$\lambda\in\mathbb{C}$, для которых линейный оператор $\lambda I-A$ не
имеет ограниченного обратного. Спектр -- непустое компактное множество
\cite[стр. 213, 214]{ReedSimon1}.

Пусть $\mathfrak{M}$ -- непустое компактное подмножество комплексной
плоскости $\mathbb{C}$, дополнение которого связно. Говорят, что
$\mathfrak{M}$ является спектральным множеством ограниченного линейного
оператора $A$ \IND{спектральное множество ограниченного линейного
оператора} в гильбертовом пространстве $\mathcal{H}$, если для любого
полинома $f(\lambda)\in\mathbb{C}[\lambda]$
\begin{equation}\label{spectr_set_l}
\|f(A)\|\le\max_{\lambda\in\mathfrak{M}}|f(\lambda)|.
\end{equation}

Отметим, что спектр самосопряженного линейного оператора лежит на
вещественной оси \cite[стр. 217]{ReedSimon1} и, следовательно, его
дополнение в $\mathbb{C}$ связно.

\begin{proposition}\label{spectr_selfadj}
Cпектр $\mathfrak{M}$ ограниченного самосопряженного оператора в
гильбертовом пространстве является спектральным множеством этого линейного
оператора.
\end{proposition}

\medskip

Это утверждение доказано, например, в \cite[стр. 248]{ReedSimon1}, и, как
отмечают авторы, из него в конечном счете следуют все утверждения
''спектральной теоремы для непрерывных функций'' \cite[стр.
247]{ReedSimon1}.

Например, непрерывные функции ограниченного самосопряженного оператора $A$
в гильбертовом пространстве $\mathcal{H}$ вводятся следующим образом. Пусть
$\mathfrak{M}\subset\mathbb{R}$ -- спектр ограниченного самосопряженного
оператора $A$. Рассмотрим $*$-алгебру $\mathbb{C}[\mathfrak{M}]$
полиномиальных функций на $\mathfrak{M}$, то есть фактор-алгебру
$\mathbb{C}[\lambda]/J_{\mathfrak{M}}$, где $J_{\mathfrak{M}}$ -- идеал
полиномов, равных нулю на $\mathfrak{M}$. Если $f\in J_{\mathfrak{M}}$, то
$f(A)=0$, как следует из предложения \ref{spectr_selfadj}. Значит,
корректно определен гомоморфизм $f(\lambda)\mapsto f(A)$\ \ $*$-алгебры
$\mathbb{C}[\mathfrak{M}]$ в $*$-алгебру ограниченных линейных операторов в
$\mathcal{H}$. Как следует из \eqref{spectr_set_l}, этот гомоморфизм
продолжается по непрерывности на замыкание алгебры полиномиальных функций в
$C(\mathfrak{M})$. Теорема Вейерштрасса о полиномиальной аппроксимации
\cite[стр. 120]{ReedSimon1} утверждает, что это замыкание совпадает с
$C(\mathfrak{M})$. Мы получили функциональное исчисление для непрерывных
функций, то есть гомоморфизм коммутативной $C^*$-алгебры $C(\mathfrak{M})$
в $C^*$-алгебру ограниченных линейных операторов, при котором
$f(\lambda)|_{\mathfrak{M}}$ переходит в $f(A)$ для всех полиномов
$f(\lambda)$.

\medskip

Отметим важное следствие ''спектральной теоремы для непрерывных функций'':
если непрерывная функция $f(\lambda)$ определена на спектре ограниченного
самосопряженного оператора $A$ и отлична от нуля хотя бы в одной точке
спектра, то $f(A)\ne 0$.

\bigskip

Линейный оператор $T$ в гильбертовом пространстве $H$ называют сжатием,
\IND{сжатие} если $\|T\|\le 1$. Пусть $E\supset H$ -- гильбертово
пространство, $P$ -- ортогональный проектор в $E$ на подпространство $H$ и
$U$ -- унитарный оператор в $E$. Его называют унитарной дилатацией сжатия
$T$, \IND{дилатация сжатия ! унитарная} если $T^n\,=\,PU^n|_H$ при всех
$n\in\mathbb{Z}_+$ \cite{NagyFoias}. Отсюда следует равенство
$T^{*n}=PU^{-n}|_H$.

Аналогично вводится понятие изометрической дилатации сжатия.\IND{дилатация
сжатия ! изометрическая}

Приведем хорошо известный результат Б. Секефальви-Надя.

\begin{proposition}\label{dilation_l}
Каждое сжатие в гильбертовом пространстве обладает унитарной дилатацией.
\end{proposition}

Короткое доказательство см. в \cite[стр. 124]{Halmos_problems}. Ограничимся
построением унитарной дилатации в наиболее важным для нас частном случае
сжатия $T$, для которого
\begin{equation}\label{partial_dil}
\lim_{n\to\infty}\|T^{*\, n}v\|=0
\end{equation}
при всех $v\in H$. (Таким является, например, сжатие $T_F(z)$.)

Пусть $F$ -- замыкание образа линейного оператора $I-TT^*$ и
$(I-TT^*)^{\frac12}$ -- квадратный корень из ограниченного неотрицательного
линейного оператора $I-TT^*$. Рассмотрим гильбертово пространство
$l^2(\mathbb{Z}_+,F)$ последовательностей $\mathbf{v}=(v_0,v_1,v_2,\cdots)$
векторов из $F$, для которых
$\|\mathbf{v}\|^2=\sum\limits_{j\in\mathbb{Z}_+}|v_j|^2 < \infty$.
Используя \eqref{partial_dil} нетрудно показать, что линейное отображение
$$
v\mapsto((I-TT^*)^{\frac12}v,(I-TT^*)^{\frac12}T^*v,
(I-TT^*)^{\frac12}T^{*2}v,(I-TT^*)^{\frac12}T^{*3}v,\cdots)
$$
изометрически вкладывает $H$ в $l^2(\mathbb{Z}_+,F)$. Оператор сдвига
$$(v_0,v_1,v_2,v_3,\cdots)\mapsto(0,v_0,v_1,v_2,v_3,\cdots)$$
в пространстве $l^2(\mathbb{Z}_+, F)$ является изометрической дилатацией
оператора, унитарно эквивалентного $T$. Остается заметить, что двусторонний
сдвиг в $l^2(\mathbb{Z},F)$ является унитарной дилатацией изометрического
сдвига в $l^2(\mathbb{Z}_+,F)$. \hfill $\square$

\begin{remark}\label{a_c_spectr}
Построенный унитарный оператор $U$ кратен двустороннему сдвигу и,
следовательно, имеет абсолютно непрерывный спектр. (Это справедливо для
любого вполне неунитарного сжатия \cite{NagyFoias}.)
\end{remark}

\begin{proposition}\label{Neumann_ineq} (Неравенство Неймана)
\IND{неравенство Неймана} Замкнутый единичный круг $\overline{\mathbb{D}}$
является спектральным множеством всех сжатий в гильбертовых пространствах.
\end{proposition}

{\bf Доказательство} В частном случае унитарного оператора требуемое
утверждение доказывается так же, как предложение \ref{spectr_selfadj}. Общий
случай сводится к этому частному случаю с помощью предложения
\ref{dilation_l}. \hfill $\square$

\subsection{Симметрия}\label{symmetry_sl_2}

\subsubsection{Алгебры Хопфа $U\mathfrak{sl}_2$ и
$U_q\mathfrak{sl}_2$.}\label{q-sl_2}

Напомним, что алгеброй Ли \IND{алгебра ! Ли} называют неассоциативную,
вообще говоря, алгебру $\mathfrak{g}$ с антикоммутативным умножением
$[x,y]=-[y,x]$, удовлетворяющим тождеству Якоби
$[x,[y,z]]+[z,[x,y]]+[y,[z,x]]=0$. Каждая ассоциативная алгебра $A$
является алгеброй Ли с операцией $[x,y]=xy-yx$.

Рассмотрим алгебру Ли $\mathfrak{g}$ и ее линейное отображение $i$ в
ассоциативную унитальную алгебру $U\mathfrak{g}$, являющееся гомоморфизмом
алгебр Ли. Пару $(U\mathfrak{g},i)$ называют универсальной обертывающей
алгеброй, \IND{универсальная ! обертывающая алгебра} если она обладает
следующим свойством универсальности. Для любой ассоциативной унитальной
алгебры $A$ и любого гомоморфизма алгебр Ли $j:\mathfrak{g}\to A$
существует и единствен такой гомоморфизм ассоциативных унитальных алгебр
$U\mathfrak{g}\to A$, что коммутативна следующая диаграмма
$$\xymatrix{\mathfrak{g}\ar[dr]^j\ar[r]^i & U\mathfrak{g}\ar[d]\\ & A}$$
Это определение универсальной обертывающей алгебры приведено, например, в
\cite[стр. 186]{Bergman}. Универсальная обертывающая алгебра в существенном
единственна.

Для доказательства ее существования выберем базис
$\{e_\lambda\}_{\lambda\in\Lambda}$ векторного пространства $\mathfrak{g}$.
В этом базисе
$[e_\lambda,e_\mu]=\sum\limits_{\nu\in\Lambda}c^\nu_{\lambda\mu}e_\nu$, где
$c^\nu_{\lambda\mu}$ -- структурные константы алгебры Ли $\mathfrak{g}$.
Остается ввести в рассмотрение ассоциативную алгебру $U\mathfrak{g}$ с
одноименными образующими $\{e_\lambda\}_{\lambda\in\Lambda}$ и определяющими
соотношениями $ e_\lambda\cdot e_\mu-e_\mu\cdot
e_\lambda=\sum\limits_{\nu\in\Lambda}c^\nu_{\lambda\mu}e_\nu$, а также
гомоморфизм
$$i:\mathfrak{g}\to U\mathfrak{g},\qquad i:e_\lambda\mapsto e_\lambda.$$
Нетрудно показать \cite[стр. 186]{Bergman}, что этот гомоморфизм $i$
является вложением, а пара $(U\mathfrak{g},i)$ -- универсальной
обертывающей алгеброй.

\begin{example}\label{Ug_examplle}
Рассмотрим комплексную алгебру Ли $\mathfrak{sl}_2$ с образующими $H$, $E$,
$F$ и определяющими соотношениями
\begin{equation*}\label{sl_2_rel}
[H,E]=2E,\qquad[H,F]=-2F,\qquad[E,F]=H.
\end{equation*}
Она изоморфна алгебре Ли комплексных $2\times2$-матриц с нулевым следом:
\begin{equation}\label{2-dim_rep}
H\mapsto \begin{pmatrix} 1, & 0\\ 0, & -1\end{pmatrix},\qquad E\mapsto
\begin{pmatrix} 0, & 1\\ 0, & 0\end{pmatrix},\qquad
F\mapsto\begin{pmatrix} 0, & 0\\ 1, & 0\end{pmatrix}.
\end{equation}
Пусть $U\mathfrak{sl}_2$ -- унитальная ассоциативная алгебра с одноименными
образующими $H$, $E$, $F$ и определяющими соотношениями
\begin{equation*}\label{Usl_2_rel}
HE-EH=2E,\qquad HF-FH=-2F,\qquad EF-FE=H.
\end{equation*}
Гомоморфизм алгебр Ли
$$i:\mathfrak{sl}_2\to U\mathfrak{sl}_2,\qquad
i:H\mapsto H,\;i:E\mapsto E,\;i:F\mapsto F
$$
обладает свойством универсальности, и пара $(U\mathfrak{sl}_2,i)$ является
универсальной обертывающей алгеброй.
\end{example}

В дальнейшем мы будем отождествлять алгебру Ли $\mathfrak{g}$ с ее образом
при вложении в $U\mathfrak{g}$.

\medskip

Имеется другое, почти эквивалентное определение. Универсальной обертывающей
алгеброй называют не любую пару $(U\mathfrak{g},i)$ со свойством
универсальности, а одну, выделенную пару, построение которой приведено в
\cite[стр. 167]{Kir}, \cite[стр. 85]{Dix}, \cite[стр. 41]{GG}, \cite[стр.
20]{Bou1-3}.

\bigskip

В теории квантовых групп используется третье определение универсальной
обертывающей алгебры $U\mathfrak{g}$, согласно которому она является
алгеброй Хопфа. Обсудим подробно это понятие.

\medskip

Начнем издалека -- с мотивировок. Рассмотрим комплексную аффинную
алгебраическую группу $G$ и ее алгебру Ли $\mathfrak{g}$\cite{Hum}.
Например, группу
$$
SL_2=\left\{\left.\begin{pmatrix}t_{11}, & t_{12}\\ t_{21}, &
t_{22}\end{pmatrix}\;\right|\;t_{11}t_{22}\,-\,t_{12}t_{21}\,=\,1\right\}
$$
и ее алгебру Ли $\mathfrak{sl}_2$.

В абелевой категории рациональных представлений группы $G$, помимо операции
$\oplus$, имеется операция $\otimes$:
$$(T_1\otimes T_2)(g)=T_1(g)\otimes T_2(g),\qquad g\in G.$$
Кроме того, имеется тривиальное одномерное представление, и каждому
представлению $T$ отвечает двойственное представление $T^*$ в сопряженном
векторном пространстве:
$$T^*(g)=T(g^{-1})^*,\qquad g\in G. $$
Известно \cite{Kir}, что категория $U\mathfrak{g}$-модулей наследует эти
операции, и понятие алгебры Хопфа призвано выявить ответственные за них
алгебраические структуры в $U\mathfrak{g}$.

\bigskip

Перейдем от мотивировок к определениям. Векторное пространство $A$ называют
коалгеброй, \IND{коалгебра} если задано линейное отображение
${\triangle:A\to A\otimes A}$, \PNT{delta@$\delta$ -- коумножение в
биалгебре} обладающее свойством коассоциативности
$$
(\triangle\otimes\operatorname{id})\triangle(a)=
(\operatorname{id}\otimes\triangle)\triangle(a),\quad a\in A,
$$
и линейное отображение $\varepsilon:A\to \mathbb{C}$,\PNT{eps@$\varepsilon$
-- коединица в биалгебре} для которого
$$
(\varepsilon\otimes\operatorname{id})\triangle(a)=
(\operatorname{id}\otimes\varepsilon)\triangle(a)=a,\quad a\in A.
$$
Эти требования равносильны коммутативности следующих диаграмм:
\begin{equation*}
\xymatrix{A\ar[d]_\Delta\ar[r]^\Delta & A\otimes A
\ar[d]^{\Delta\otimes\mathrm{id}} \\
A\otimes A\ar[r]_{\mathrm{id}\otimes\Delta\quad} & A \otimes A \otimes A }
\qquad \xymatrix{A\ar[d]_\Delta\ar[dr]^{\mathrm{id}}\ar[r]^\Delta &
A\otimes A\ar[d]^{\mathrm{id}\otimes\varepsilon} \\
A\otimes A\ar[r]_{\varepsilon\otimes\mathrm{id}} & A}
\end{equation*}
Отображение $\triangle$ называют коумножением, \IND{коумножение} а
отображение $\varepsilon$ -- коединицей \IND{коединица}.

\begin{remark}
Два слова об этимологии. Изменяя направления стрелок, получаем диаграммы,
выражающие ассоциативность умножения
$$m:A\otimes A\to A,\qquad m:a_1\otimes a_2\mapsto a_1a_2$$
и наличие единицы
$$\unit:\mathbb{C}\to A,\qquad\unit:z\mapsto z\cdot 1$$
в алгебре $A$:
\begin{equation*}\label{Algebra}
\xymatrix{
A & A \otimes A\ar[l]_m
\\ A \otimes A\ar[u]^m & A\otimes A\otimes A\ar[l]^{\mathrm{id}\otimes m
\quad}\ar[u]_{m\otimes\mathrm{id}}}\qquad
\xymatrix{
A & A\otimes A\ar[l]_m
\\ A\otimes A\ar[u]^m & A\ar[l]^{\mathrm{id}\otimes\unit}
\ar[lu]_{\mathrm{id}}\ar[u]_{\unit\otimes\mathrm{id}}}
\end{equation*}
\end{remark}

\bigskip \medskip

Векторное пространство $A$ называют биалгеброй, \IND{биалгебра} если оно
наделено как структурой алгебры с единицей, так и структурой коалгебры с
коединицей, причем отображения $\triangle$ и $\varepsilon$ являются
гомоморфизмами алгебр. Подразумевается стандартная структура алгебры в
$A\otimes A$:
$$
(a'\otimes b')\cdot(a''\otimes b'')=a'a''\otimes b'b'',\qquad
a',a'',b',b''\in A.
$$

Отметим, что $\triangle(1)=1\otimes 1$, $\varepsilon(1)=1$. Сказанное
равносильно коммутативности следующих диаграмм:
\begin{equation*}\label{Large}
\xymatrix{
& A\otimes A\ar[r]^{\Delta\otimes\Delta\quad\quad}
& A\otimes A\otimes A\otimes A\ar[dd]^{\mathrm{id}\otimes\sigma\otimes
\mathrm{id}}
\\ A\ar[ur]^\Delta & &
\\ & A\otimes A\ar[lu]^m & A\otimes A\otimes A\otimes A\ar[l]^{m\otimes m
\quad\quad} }
\end{equation*}
\begin{equation*}\label{Respect-1}
\xymatrix{
A\ar[r]^{\varepsilon} & \mathbb{C}
\\ A\otimes A\ar[u]^m\ar[r]_{\varepsilon\otimes\varepsilon} & \mathbb{C}
\otimes\mathbb{C}\ar[u]_m}\qquad
\xymatrix{
A\ar[d]_{\Delta} & \mathbb{C}\ar[l]_{\unit}\ar[d]^{\Delta}
\\ A\otimes A & \mathbb{C}\otimes\mathbb{C}\ar[l]^{\unit\otimes\unit}}
\qquad \xymatrix{
A \ar@/^/[r]^{\varepsilon} & \mathbb{C}\ar@/^/[l]^{\unit}}
\end{equation*}

Здесь $\sigma$ -- оператор ''наивной'' перестановки тензорных сомножителей:
$$
\sigma:A\otimes A\to A\otimes A,\qquad\sigma:a_1\otimes a_2\mapsto
a_2\otimes a_1.
$$

\medskip\medskip

Коумножение $\Delta$ позволяет ввести операцию тензорного произведения
представлений алгебры $A$ ($\equiv$ {$A$-модулей}), \IND{тензорное
произведение представлений алгебры} обладающую привычными свойствами
операции тензорного произведения векторных пространств. Именно, пусть
$\pi_1$, $\pi_2$ -- представления биалгебры $A$ в векторных пространствах
$V_1,V_2$ соответственно. Их тензорным произведением $\pi_1\otimes\pi_2$
называют следующее представление биалгебры $A$ в пространстве $V_1\otimes
V_2$:
$$
A\overset{\triangle}{\to}A\otimes
A\overset{\pi_1\boxtimes\pi_2}{\longrightarrow}
\operatorname{End}(V_1)\otimes\operatorname{End}(V_2)\hookrightarrow
\operatorname{End}(V_1\otimes V_2).
$$

Тривиальным представлением биалгебры $A$ \IND{тривиальное представление
биалгебры} называют её одномерное представление, задаваемое коединицей
$\varepsilon:A\to\mathbb{C}\cong\operatorname{End}(\mathbb{C})$. Отвечающий
этому представлению $A$-модуль $\mathbb{C}$ называют тривиальным
$A$-модулем. \IND{модуль ! тривиальный}

Рассмотрим $A$-модуль $V$ и вектор $v\in V$. Этот вектор называют
$A$-инвариантным, \IND{инвариантный ! вектор} если линейное отображение
$\mathbb{C}\to A$, при котором $1\mapsto v$, является морфизмом
$A$-модулей. Другими словами, $A$-инвариантность вектора $v$ означает, что
$av=\varepsilon(a)$ при всех $a\in A$.

Из свойств коумножения и коединицы следует, что канонические изоморфизмы
векторных пространств
$$
V_1\otimes(V_2\otimes V_3)\cong(V_1\otimes V_2 )\otimes V_3,\qquad
V\otimes\mathbb{C}\cong\mathbb{C}\otimes V\cong V
$$
являются изоморфизмами $A$-модулей. Если $f_1:V_1\to W_1$, $f_1:V_1\to W_1$
являются морфизмами $A$-модулей, то и $f_1\otimes f_2:V_1\otimes W_1\to
V_2\otimes W_2$ -- морфизм $A$-модулей, то есть $A$-модули образуют
тензорную категорию, см. \itemiii\ \ref{braided_categories_sl_2}.

\bigskip

Перейдем от биалгебр к алгебрам Хопфа. Рассмотрим биалгебру $A$ с
умножением $m:A\otimes A\to A$, коумножением $\triangle:A\to A\otimes A$ и
коединицей $\varepsilon:A\to\mathbb{C}$. Линейное отображение $S:A\to A$
\PNT{S@$S$ -- антипод алгебры Хопфа} называют антиподом, \IND{антипод} если
$m(S\otimes\operatorname{id})\triangle(a)=m(\operatorname{id}\otimes
S)\triangle(a)=\varepsilon(a)1$ при всех $a\in A$, что равносильно
коммутативности следующей диаграммы
\begin{equation*}\label{S-1}
\xymatrix{& A\otimes A\ar[rr]^{S\otimes\mathrm{id}} & &
A\otimes A\ar[rd]^m & \\
A\ar[ur]^{\Delta}\ar[rr]^{\varepsilon}\ar[rd]_{\Delta} & & \mathbb{C}
\ar[rr]^{\unit} & & A \\
& A\otimes A\ar[rr]_{\mathrm{id}\otimes S} & & A\otimes A\ar[ru]_m}
\end{equation*}
Как вытекает из определений \cite[стр. 13]{Jo}, антипод единствен, является
антигомоморфизмом алгебры и антигомоморфизмом коалгебры, то есть
коммутативны следующие диаграммы
\begin{equation*}\label{S3_S2}
\xymatrix{A\otimes A\ar[d]_{S\otimes S}\ar[rr]^m & & A\ar[d]^S \\
A\otimes A\ar[r]_{\sigma} & A\otimes A\ar[r]_m & A}\qquad
\xymatrix{A\ar[d]_S\ar[rr]^{\Delta} & & A\otimes A\ar[d]^{S\otimes S}\\
A\ar[r]_{\Delta} & A\otimes A & A\otimes A\ar[l]^{\sigma}}
\end{equation*}

Всюду далее антипод $S$ будет предполагаться обратимым.

Биалгебру $A$, обладающую антиподом, называют алгеброй Хопфа. \IND{алгебра
! Хопфа}

\bigskip

Рассмотрим алгебру Хопфа $A$ и $A$-модуль $V$. Сопряженное векторное
пространство $V^*$ наделяется структурой $A$-модуля:
$$(al)(v)=l(S(a)v),\qquad a\in A,\;v\in V,\;l\in V^*.$$

Каноническое спаривание векторных пространств
$$V^*\otimes V\to\mathbb{C},\qquad l\otimes v\mapsto l(v)$$
является морфизмом $A$-модулей.

\bigskip

Универсальная обертывающая алгебра является алгеброй Хопфа. Действительно,
гомоморфизмы алгебр Ли
$$
\mathfrak{g}\to U\mathfrak{g}\otimes U\mathfrak{g},\quad x\mapsto
i(x)\otimes 1+1\otimes i(x),\text{\ \ и\ \ }\mathfrak{g}\to\mathbb{C},\quad
x\mapsto 0
$$
единственным образом продолжаются до гомоморфизмов алгебр
$$
\triangle:U\mathfrak{g}\to U\mathfrak{g}\otimes U\mathfrak{g}, \qquad
\varepsilon:U\mathfrak{g}\to\mathbb{C},
$$
а антиавтоморфизм алгебры Ли $\mathfrak{g}$, при котором $x$ переходит в
$-x$, единственным образом продолжается до антиавтоморфизма $S$ алгебры
$U\mathfrak{g}$.

\begin{example}\label{U-sl_2}
Коумножение $\Delta$, коединица $\varepsilon$ и антипод $S$ алгебры Хопфа
$U\mathfrak{sl}_2$ определяются равенствами
$$
\Delta(H)=H\otimes 1+1\otimes H,\quad\Delta(E)=E\otimes 1+1\otimes E,\quad
\Delta(F)=F\otimes 1+1\otimes F,
$$
$$\varepsilon(H)=\varepsilon(E)=\varepsilon(F)=0,$$
$$S(H)=-H,\qquad S(E)=-E,\qquad S(F)=-F.$$
\end{example}

\bigskip

Понятия тензорного произведения $U\mathfrak{g}$-модулей, тривиального
$U\mathfrak{g}$-модуля и двойственного $U\mathfrak{g}$-модуля согласуются с
соответствующими понятиями теории групп Ли \cite{Kir}.

\bigskip

В теории квантовых групп универсальной обертывающей алгеброй называют
описанную выше алгебру Хопфа $U\mathfrak{g}$.

Алгебра Хопфа $U\mathfrak{sl}_2$ кокоммутативна, то есть
$\Delta=\Delta^{\mathrm op}$, где $\Delta^{\mathrm op}=\sigma\Delta$. Кроме
того, $S^2={\mathrm id}$. Приведем пример некокоммутативной алгебры Хопфа с
антиподом бесконечного порядка. Напомним, что $q\in(0,1)$.

\begin{example}\label{U_q-sl_2}(\cite[глава 3]{Jant},
\cite[глава 7]{Kassel_QG}). Рассмотрим алгебру с образующими $K$,
$K^{-1}$,$E$, $F$ и следующими определяющими соотношениями:
\begin{equation}\label{U_q_first}
KK^{-1}=K^{-1}K=1,\quad K^{\pm 1}E=q^{\pm 2}EK^{\pm 1},\quad K^{\pm
1}F=q^{\mp 2}FK^{\pm 1},
\end{equation}
\begin{equation}\label{U_q_second}
EF-FE=\frac{K-K^{-1}}{q-q^{-1}}.
\end{equation}

Можно доказать существование гомоморфизмов
$$
\Delta:U_q\mathfrak{sl}_2\to U_q\mathfrak{sl}_2\otimes
U_q\mathfrak{sl}_2,\qquad\varepsilon:U_q\mathfrak{sl}_2\to\mathbb{C}
$$
и антигомоморфизма $S$ алгебры $U_q\mathfrak{sl}_2$, для которых
$$
\Delta(K^{\pm 1})=K^{\pm 1}\otimes K^{\pm 1},\quad\Delta(E)=E\otimes
1+K\otimes E,\quad\Delta(F)=F\otimes K^{-1}+1\otimes F,
$$
$$\varepsilon(E)=\varepsilon(F)=0,\qquad\varepsilon(K^{\pm 1})=1,$$
$$S(K^{\pm 1})=K^{\mp 1},\quad S(E)=-K^{-1}E,\quad S(F)=-FK.$$
Тем самым $U_q\mathfrak{sl}_2$ наделяется структурой алгебры Хопфа.
Определяющие ее соотношения связаны с определяющими соотношениями для
$U\mathfrak{sl}_2$ подстановкой
$$
q=\exp\left(-\frac{h}2\right),\qquad K^{\pm
1}=\exp\left(\mp\frac{hH}2\right)
$$
и формальным предельным переходом $h\to 0$.
\end{example}

\bigskip

О пути, который привел к этой алгебре Хопфа, рассказано в обзоре Дринфельда
\cite[стр. 807]{DrinfEng}. Она является простейшим представителем класса
квантовых универсальных обертывающих алгебр Дринфельда-Джимбо \cite[глава
4]{Jant}.

\subsubsection{$*$-Алгебры Хопфа $U\mathfrak{su}_{1,1}$, $U\mathfrak{su}_2$,
$U_q\mathfrak{su}_{1,1}$, $U_q\mathfrak{su}_2$.}\label{star-Hopf_algebras}

Пусть \hbox{$*:\mathfrak{g}\to\mathfrak{g}$} -- антилинейный инволютивный
антиавтоморфизм комплексной алгебры Ли $\mathfrak{g}$, \IND{антилинейный
инволютивный антиавтоморфизм комплексной алгебры Ли} то есть
$$
(\lambda a+\mu b)^*=\bar{\lambda}a^*+\bar{\mu}b^*,\qquad
a^{**}=a,\qquad[a,b]^*=[b^*,a^*]
$$
для любых $a,b\in\mathfrak{g}$, $\lambda,\mu\in\mathbb{C}$. Такая инволюция
$*$ единственным образом продолжается до антилинейного инволютивного
антиавтоморфизма $*:U\mathfrak{g}\to U\mathfrak{g}$ универсальной
обертывающей алгебры.

Множество антиэрмитовых элементов $\{a\in\mathfrak{g}|\:a^*=-a\}$ является
вещественной алгеброй Ли. В теории унитарных представлений групп Ли
ключевую роль играют представления вещественных алгебр Ли антиэрмитовыми
операторами, то есть $*$-представления алгебры $U\mathfrak{g}$. Это
позволяет в ряде случаев заменить выбор вещественной формы комплексной
алгебры Ли выбором инволюции.

Далее все инволюции предполагаются антилинейными антиавтоморфизмами, если
не оговорено противное.

\begin{example}\label{Usu_2}
Наделим алгебру Ли $\mathfrak{sl}_2$ инволюцией $\star$, полагая
\begin{equation}\label{classical_inv_star}
H^\star=H,\qquad E^\star=F,\qquad F^\star=E.
\end{equation}
Продолжим эту инволюцию на $U\mathfrak{sl}_2$ и введем обозначение
$U\mathfrak{su}_2=(U\mathfrak{sl}_2,\star)$. Очевидно, конечномерные
$*$-представления алгебры $U\mathfrak{su}_2$ находятся во взаимно
однозначном соответствии с конечномерными представлениями вещественной
алгебры Ли $\mathfrak{su}_2=\{a\in\mathfrak{sl}_2|\:a^\star=-a\}$
антиэрмитовыми операторами. Отметим, что вещественная алгебра Ли
$\mathfrak{su}_2$ изоморфна алгебре Ли группы движений двумерной сферы, ср.
\cite[стр. 78]{Vinb_reps}, \cite[стр. 109,110]{Vilenkin}.
\end{example}

\begin{example}\label{Usu_1_1} Наделим алгебру Ли $\mathfrak{sl}_2$
инволюцией $*$, полагая
\begin{equation}\label{classical_inv}
H^*=H,\qquad E^*=-F,\qquad F^*=-E.
\end{equation}
Продолжим эту инволюцию на $U\mathfrak{sl}_2$ и введем обозначение
$U\mathfrak{su}_{1,1}=(U\mathfrak{sl}_2,*)$. Вещественная алгебра Ли
$\mathfrak{su}_{1,1}\,=\,\{a\in\mathfrak{sl}_2\,|\,a^*=-a\}$ интересна нам
потому, что она изоморфна алгебре Ли группы дробно-линейных преобразований
единичного круга $\mathbb{D}$, ср. \cite[стр. 45]{Helg1}, \cite[стр. 294,
295]{Vilenkin}. \IND{алгебра ! $U\mathfrak{su}_{1,1}$}
\end{example}

Отметим, что введенные обозначения несущественно отличаются от
традиционных, согласно которым универсальные обертывающие алгебры
$U\mathfrak{su}_{1,1}$, $U\mathfrak{su}_2$ являются алгебрами над
$\mathbb{R}$, а не над $\mathbb{C}$.

\bigskip

$*$-Алгебра Хопфа $A$ -- это алгебра Хопфа с инволюцией $*$, являющейся
автоморфизмом коалгебры:
$$\Delta(a^*)=\Delta(a)^{*\otimes*},\qquad a\in A.$$
Нетрудно доказать, что
\begin{equation}\label{equal_star}
1^*=1,\quad *S*=S^{-1},\quad
\varepsilon(a^*)=\overline{\varepsilon(a)},\quad a\in A.
\end{equation}
Действительно, пусть $A^{\operatorname{cop}}$ -- алгебра Хопфа,
отличающаяся от $A$ заменой коумножения на противоположное. Первое и третье
равенства вытекают из единственности единицы и коединицы, а второе -- из
единственности и обратимости антипода, поскольку операторы $*S*$, $S^{-1}$
являются антиподами алгебры Хопфа $A^{\operatorname{cop}}$. Доказательство
единственности антипода см. в \cite[стр. 65]{Kassel_QG}, \cite[стр.
13]{Jo}.

\medskip

Рассмотрим два $*$-представления $\pi_1,\pi_2$ $*$-алгебры Хопфа $A$ в
предгильбертовых пространствах $\mathcal{H}_1$, $\mathcal{H}_2$. Из
определений следует, что их тензорное произведение является
$*$-представлением этой $*$-алгебры Хопфа в тензорном произведении
предгильбертовых пространств. Это утверждение сродни тому, что тензорное
произведение унитарных представлений группы является ее унитарным
представлением.

\medskip

Из определения алгебры Хопфа $A$ на языке коммутативных диаграмм видно, что
при условии $\dim A<\infty$ сопряженное векторное пространство $A^*$
является алгеброй Хопфа с умножением, единицей, коумножением, коединицей и
антиподом, сопряженными к коумножению, коединице, умножению, единице и
антиподу алгебры Хопфа $A$.

\begin{remark}\label{duality_of_Hopf}
Требование конечномерности алгебры Хопфа $A$ в предыдущих рассуждениях
можно ослабить. Именно, рассмотрим множество $A^\star$ линейных
функционалов на $A$, ядра которых содержат двусторонние идеалы конечной
коразмерности. Как известно \cite[стр. 27]{Jo}, $A^\star$ по двойственности
наделяется структурой алгебры Хопфа. Ее называют двойственной алгеброй
Хопфа. Если $A$ является $*$-алгеброй Хопфа, то равенство
\begin{equation}\label{invol1}
l^*(a)\stackrel{\operatorname{def}}{=}\overline{l(S(a)^*)},\qquad l\in
A^*,\quad a\in A
\end{equation}
наделяет $A^\star$ структурой $*$-алгебры Хопфа.
\end{remark}

\bigskip\medskip

Как следует из определений, $U\mathfrak{su}_{1,1}$ и $U\mathfrak{su}_2$
являются $*$-алгебрами Хопфа. Введем в рассмотрение их квантовые аналоги.

Рассмотрим инволюцию $*$ в $U_q\mathfrak{sl}_2$, для которой
\begin{equation}\label{U_q_su_1_1}
(K^{\pm 1})^*=K^{\pm 1},\qquad E^*=-KF,\qquad F^*=-EK^{-1}.
\end{equation}
Она существует, единственна и наделяет $U_q\mathfrak{sl}_2$ структурой
$*$-алгебры Хопфа, обозначаемой $U_q\mathfrak{su}_{1,1}$. \IND{алгебра !
$U_q\mathfrak{su}_{1,1}$}

Аналогично, инволюция $\star$ в $U_q\mathfrak{sl}_2$, для которой
\begin{equation}\label{U_q_su_2}
(K^{\pm 1})^\star=K^{\pm 1},\qquad E^\star=KF,\qquad F^\star=EK^{-1},
\end{equation}
существует, единственна и наделяет $U_q\mathfrak{sl}_2$ структурой
$*$-алгебры Хопфа, обозначаемой $U_q\mathfrak{su}_2$. \IND{алгебра !
$U_q\mathfrak{su}_2$}

$*$-Алгебра Хопфа $U_q\mathfrak{su}_2$ обладает двумерным
$*$-представлением $\pi$ в гильбертовом пространстве $\mathbb{C}^2$:
\begin{equation}\label{q2-dim_rep}
\pi(K^{\pm 1})=\begin{pmatrix} q^{\pm 1}, & 0\\
0, & q^{\mp 1}\end{pmatrix},\;\;\pi(E)=\begin{pmatrix} 0, &
q^{-\frac12}\\ 0, & 0\end{pmatrix},\;\;\pi(F)=\begin{pmatrix} 0, & 0\\
q^{\frac{1}{2}}, & 0\end{pmatrix}.
\end{equation}
Соответствующее представление алгебры $U_q\mathfrak{sl}_2$ в векторном
пространстве $\mathbb{C}^2$ называют ее векторным представлением.
\IND{векторное представление ! алгебры $U_q\mathfrak{sl}_2$}

\subsubsection{$U_q\mathfrak{sl}_2$-модульные алгебры
$\mathbb{C}[z]_q$, $\mathbb{C}[z^*]_q$.}\label{sL_2_mod_algebras}

Пусть $A$ -- алгебра Хопфа. Алгебру $F$ без единицы, наделенную структурой
$A$-модуля, называют $A$-модульной, если умножение
$$m:F\otimes F\to F,\qquad m:f_1\otimes f_2\mapsto f_1f_2$$
является морфизмом $A$-модулей. Это означает, что
$a(f_1f_2)=\sum\limits_ib_i(f_1)c_i(f_2)$, если $a\in A$ и
$\triangle(a)=\sum\limits_ib_i\otimes c_i$. В случае алгебры с единицей в
определение включается дополнительное требование, чтобы линейное
отображение
\begin{equation}\label{unit F}
\unit:\mathbb{C}\to F,\qquad\unit:1\mapsto 1
\end{equation}
было морфизмом $A$-модулей. Другими словами, единица алгебры $F$ должна
быть ее $A$-инвариантным элементом. \IND{ алгебра ! $A$-модульная}

В частном случае $A=U_q\mathfrak{sl}_2$ эти требования означают следующее:
\begin{equation}\label{K_Uq-algebra}
K^{\pm 1}(f_1f_2)=K^{\pm 1}(f_1)K^{\pm 1}(f_2),\qquad K^{\pm 1}(1)=1,
\end{equation}
\begin{equation}\label{E_Uq-algebra}
E(f_1f_2)=E(f_1)f_2+K(f_1)E(f_2),\qquad E(1)=0,
\end{equation}
\begin{equation}\label{F_Uq-algebra}
F(f_1f_2)=K^{-1}(f_1)F(f_2)+f_1F(f_2),\qquad F(1)=0.
\end{equation}

Гомоморфизмом $A$-модульных алгебр \IND{гомоморфизм $A$-модульных алгебр}
называется линейное отображение, являющееся как гомоморфизмом алгебр, так и
морфизмом $A$-модулей.

\bigskip

Аналогично понятию $A$-модульной алгебры вводятся понятия $A$-модульного
левого модуля, $A$-модульного правого модуля и $A$-модульного бимодуля $M$
над $A$-модульной алгеброй $F$: \IND{$A$-модульный $F$-модуль}
\begin{itemize}
\item[-]\ \ $M$ -- модуль над алгеброй Хопфа $A$,

\item[-]\ \ действия $m_\mathrm{left}:F\otimes M\to M$, и $m_\mathrm{
    right}:M\otimes F\to M$ являются морфизмами $A$-модулей,

\item[-]\ \ единице алгебры $F$ отвечает тождественный оператор в $M$.
\end{itemize}
Морфизмом $A$-модульных $F$-модулей \IND{морфизм ! $A$-модульных
$F$-модулей} называется линейное отображение, являющееся как морфизмом
$F$-модулей, так и морфизмом $A$-модулей.

\bigskip

Приведем примеры.

\begin{example}\label{EndV_sl_2}
Для любого модуля $V$ над алгеброй Хопфа $A$ алгебра
$\operatorname{End}(V)$ всех линейных операторов в $V$ является
$A$-модульной. Именно, если $a\in A$, $\triangle(a)=\sum\limits_ib_i\otimes
c_i$, то
\begin{equation}\label{End_mod}
(aT)v=\sum_ib_iT(S(c_i)v),\qquad T\in\operatorname{End}(V),\;v\in V.
\end{equation}
Естественное вложение
$$\mathbb{C}\hookrightarrow\operatorname{End}(V),\qquad 1\mapsto I,$$
является гомоморфизмом $A$-модульных алгебр.

Для любых $A$-модулей $V_1$ и $V_2$ векторное пространство
$\operatorname{Hom}(V_1,V_2)$ линейных операторов из $V_1$ в $V_2$
наделяется структурой $A$-модуля равенством, аналогичным \eqref{End_mod}.
\end{example}

\begin{example}\label{classic_plane}
Рассмотрим дифференциальные операторы первого порядка в пространстве
$\mathbb{C}[t_1^{\pm 1},t_2^{\pm 1}]$ полиномов Лорана двух переменных
\begin{equation}\label{sl_2-action}
H=t_1\frac{\partial}{\partial t_1}-t_2\frac{\partial}{\partial t_2},\quad
E=t_1\frac{\partial}{\partial t_2},\quad F=t_2\frac{\partial}{\partial t_1}.
\end{equation}
Они удовлетворяют определяющим соотношениям алгебры $U\mathfrak{sl}_2$,
являются дифференцированиями алгебры $\mathbb{C}[t_1^{\pm 1},t_2^{\pm 1}]$
и, следовательно, наделяют ее структурой $U\mathfrak{sl}_2$-модульной
алгебры. Отметим, что

\begin{equation}\label{H-t}
(Ht_1,Ht_2)=(t_1,t_2)\begin{pmatrix} 1, & 0\\ 0, & -1\end{pmatrix},
\end{equation}
\begin{equation}\label{EF-t_classic}
(Et_1,Et_2)=(t_1,t_2)\begin{pmatrix} 0, & 1\\ 0, & 0\end{pmatrix},
\qquad(Ft_1,Ft_2)=(t_1,t_2)\begin{pmatrix} 0, & 0\\ 1, & 0\end{pmatrix}.
\end{equation}
\end{example}

\medskip

\begin{example}\label{classic_z} Рассмотрим алгебру
$\mathbb{C}[z]$ полиномов одной переменной и вложение алгебр
$$
\mathbb{C}[z]\hookrightarrow\mathbb{C}[t_1^{\pm 1},t_2^{\pm 1}],\qquad
f(z)\mapsto f\left(\frac{t_1}{t_2}\right).
$$
Образ $\mathbb{C}[z]$ при этом вложении является подмодулем
$U\mathfrak{sl}_2$-модуля $\mathbb{C}[t_1^{\pm 1},t_2^{\pm 1}]$. Значит,
$\mathbb{C}[z]$ наследует структуру $U\mathfrak{sl}_2$-модульной алгебры:
\begin{equation}\label{new-sl_2-action}
H=2z\frac{d}{dz},\quad E=-z^2\frac{d}{dz},\quad F=\frac{d}{dz}.
\end{equation}
Очевидно,
\begin{equation}\label{1_on_z-action}
Hz=2z,\qquad Ez=-z^2,\qquad Fz=1.
\end{equation}
\end{example}

\bigskip

Построение $U_q\mathfrak{sl}_2$-модульных алгебр начнем с одного важного
частного случая. Рассмотрим векторное пространство $M$, его тензорный
квадрат $M\otimes M$ и подпространство $N\subset M\otimes M$. Квадратичной
алгеброй с пространством образующих $M$ и пространством соотношений $N$
называют фактор-алгебру $T(M)/(N)$ тензорной алгебры $T(M)$ по
двустороннему идеалу, порожденному множеством $N$ \cite[стр.
118]{Kassel_QG}, \cite[стр. 19]{Man1}. \IND{алгебра ! квадратичная} Если
$A$ -- алгебра Хопфа, $M$ -- модуль над $A$ и $N$ -- подмодуль его
тензорного квадрата, то квадратичная алгебра $T(M)/(N)$ является
$A$-модульной.

\begin{example}\label{manins-plane}
Рассмотрим алгебру $\mathbb{C}[t_1,t_2]_q=\mathbb{C}\langle t_1,t_2
\rangle/(t_1t_2-qt_2t_1)$ с образующими $t_1$, $t_2$ и определяющим
соотношением $t_2t_2=qt_2t_1$. Ее называют алгеброй регулярных функций на
квантовом векторном пространстве $\mathbb{C}^2$. Докажем, что равенства

\begin{equation}\label{K-t}
(K^{\pm 1}t_1,K^{\pm 1}t_2)=(t_1,t_2)\begin{pmatrix} q^{\pm 1}, & 0\\ 0, &
q^{\mp 1}\end{pmatrix},
\end{equation}
\begin{equation}\label{EF-t}
(Et_1,Et_2)=(t_1,t_2)\begin{pmatrix} 0, & q^{-\frac12}\\ 0, &
0\end{pmatrix}, \qquad(Ft_1,Ft_2)=(t_1,t_2)\begin{pmatrix} 0, & 0\\
q^{\frac12}, & 0\end{pmatrix}
\end{equation}
наделяют $\mathbb{C}[t_1,t_2]_q$ структурой $U_q\mathfrak{sl}_2$-модульной
алгебры. Для этого рассмотрим векторное представление алгебры Хопфа
$U_q\mathfrak{sl}_2$ в $M=\mathbb{C}^2$ и стандартный базис $\{t_1,t_2\}$.
Действие образующих алгебры $U_q\mathfrak{sl}_2$ на $t_1$, $t_2$
описывается равенствами \eqref{K-t}, \eqref{EF-t}. Как нетрудно показать,
вектор $t_1\otimes t_2-qt_2\otimes t_1$ является
$U_q\mathfrak{sl}_2$-инвариантным. Полагая $N=\mathbb{C}(t_1\otimes
t_2-qt_2\otimes t_1)$, приходим к $U_q\mathfrak{sl}_2$-модульной
квадратичной алгебре, естественно изоморфной $\mathbb{C}[t_1,t_2]_q$.

Отметим, что равенства \eqref{K-t}, \eqref{EF-t} являются $q$-аналогами
равенств \eqref{H-t},\eqref{EF-t_classic}.
\end{example}

\bigskip

$q$-Аналоги биномиальных коэффициентов были введены Гауссом \cite{Andrews}.
Используя алгебру $\mathbb{C}[t_1,t_2]_q$, легко объяснить появление
$q$-биномиальных коэффициентов в рассматриваемом круге вопросов. Множество
$\{t_2^{j_2}t_1^{j_1}\}_{j_1,j_2\in\mathbb{Z}_+}$ является базисом
векторного пространства $\mathbb{C}[t_1,t_2]_q$. (Действительно, линейная
независимость мономов следует из того, что алгебра $\mathbb{C}[t_1,t_2]_q$
является $\mathbb{Z}^2$-градуированной
$$\deg(t_1)=(1,0),\qquad\deg(t_2)=(0,1),$$
и бистепени рассматриваемых мономов попарно различны. Остается заметить,
что эти мономы отличны от нуля, поскольку отличны от нуля операторы
$T(t_2^{j_2}t_1^{j_1})$ представления $T$ алгебры $\mathbb{C}[t_1,t_2]_q$ в
векторном пространстве функций на множестве $q^{\mathbb{Z}}$, определяемого
равенствами: $T(t_1)f(u)=f(qu)$, $T(t_2)f(u)=uf(u)$.)

$q$-Аналоги $\binom{n}{k}_q$ биномиальных коэффициентов $\binom{n}{k}$
можно ввести равенством
$$(t_1+t_2)^n=\sum_{k=0}^n\binom{n}{k}_qt_2^kt_1^{n-k}.$$
Очевидно, $\lim\limits_{q\to 1}\binom{n}{k}_q=\binom{n}{k}$. Без труда
находятся явный вид $q$-биномиальных коэффициентов \cite[стр.
96]{Kassel_QG}: \IND{$q$-биномиальные коэффициенты}
\begin{equation}\label{q-binom_l}
\binom{n}{k}_q=\frac{(q;q)_n}{(q;q)_k(q;q)_{n-k}},
\end{equation}
где использовано стандартное обозначение
$$(a;q)_m=\prod\limits_{j=0}^{m-1}(1-aq^j).$$

\bigskip

Перейдем к описанию $q$-аналога $U\mathfrak{sl}_2$-модульной алгебры
$\mathbb{C}[t_1^{\pm 1}, t_2^{\pm 1}]$. Следующее утверждение очевидно.
\begin{lemma}\label{uniq_set}
Если полином Лорана $n$ переменных $u_1,u_2,\ldots,u_n$ обращается в нуль
на множестве $\left\{\left.(q^{j_1},q^{j_2},\ldots,q^{j_n})\;\right|\;
j_1,j_2,\ldots,j_n\in\mathbb{N}\right\}$, то он равен нулю тождественно.
\end{lemma}

\medskip

\begin{lemma}\label{Laurent_1} Существует и единственна такая тройка
полиномов Лорана $k(u_1,u_2),e(u_1,u_2),f(u_1,u_2)\in\mathbb{C}[u_1^{\pm
1},u_2^{\pm 1}]$, что при всех $j_1,j_2 \in \mathbb{N}$ в
$\mathbb{C}[t_1,t_2]_q$ имеют место равенства
\begin{equation}\label{k-action}
K^{\pm 1}(t_2^{j_2}t_1^{j_1})=k(q^{\pm j_1},q^{\pm
j_2})t_2^{j_2}t_1^{j_1},
\end{equation}
\begin{equation}\label{ef-action}
E(t_2^{j_2}t_1^{j_1})=e(q^{j_1},q^{j_2})t_2^{j_2+1}t_1^{j_1-1},\qquad
F(t_2^{j_2}t_1^{j_1})=f(q^{j_1},q^{j_2})t_2^{j_2-1}t_1^{j_1+1}.
\end{equation}
\end{lemma}

{\bf Доказательство.} Существование таких полиномов Лорана устанавливается
прямыми вычислениями с помощью формулы суммы членов геометрической
прогрессии. Единственность вытекает из леммы \ref{uniq_set}. \hfill
$\square$

\bigskip

\begin{example}\label{to_shilov_boundary}
Как легко показать, используя $\mathbb{Z}^2$-градуировку, алгебра
$\mathbb{C}[t_1,t_2]_q$ не имеет делителей нуля, и ее мультипликативное
подмножество $(t_1t_2)^{\mathbb{Z}_+}$ является множеством Оре (см.
\itemiii\ \ref{localization}). Локализацию обозначим $\mathbb{C}[t_1^{\pm
    1},t_2^{\pm 1}]_q$. Эту алгебру можно задать
образующими $t_1^{\pm 1}$, $t_2^{\pm 1}$ и определяющими соотношениями
$$
t_1t_2=qt_2t_1,\qquad t_1t_1^{-1}=t_1^{-1}t_1=1,\qquad
t_2t_2^{-1}=t_2^{-1}t_2=1.
$$
Как и прежде, доказывается, что множество
$\{t_2^{j_2}t_1^{j_1}\}_{j_1,j_2\in\mathbb{Z}}$ является базисом векторного
пространства $\mathbb{C}[t_1^{\pm 1},t_2^{\pm 1}]_q$.

Покажем, что структура $U_q\mathfrak{sl}_2$-модульной алгебры единственным
образом продолжается с $\mathbb{C}[t_1,t_2]_q$ на $\mathbb{C}[t_1^{\pm
1},t_2^{\pm 1}]_q$. Единственность продолжения легко следует из
определений. Остается предъявить продолжение.

Используя лемму \ref{Laurent_1}, введем линейные операторы $K^{\pm 1}$,
$E$, $F$ в пространстве $\mathbb{C}[t_1^{\pm 1},t_2^{\pm 1}]_q$, задав их
действие на элементы базиса
$\{t_2^{j_2}\,t_1^{j_1}\}_{j_1,j_2\in\,\mathbb{Z}}$ равенствами
\eqref{k-action}, \eqref{ef-action}. Нетрудно показать, что эти операторы
удовлетворяют определяющим соотношениям \eqref{U_q_first},
\eqref{U_q_second} алгебры $U_q\mathfrak{sl}_2$ и наделяют
$\mathbb{C}[t_1^{\pm 1},t_2^{\pm 1}]_q$ структурой
$U_q\mathfrak{sl}_2$-модульной алгебры. Ограничимся проверкой соотношения
\eqref{U_q_second}.

Так как сужения операторов $EF-FE$, $\frac{K-K^{-1}}{q-q^{-1}}$ на
подпространство $\mathbb{C}[t_1,t_2]_q$ равны, то
\begin{multline}\label{first_Laur}
e(qu_1,q^{-1}u_2)f(u_1,u_2)-f(q^{-1}u_1,qu_2)e(u_1,u_2)=
\\ =\frac{k(u_1,u_2)-k(u_1^{\pm 1},u_2^{\pm 1})}{q-q^{-1}}
\end{multline}
при всех $u_1,u_2\in q^\mathbb{N}$. Из леммы \ref{uniq_set} следует, что
\eqref{first_Laur} имеет место при всех $u_1,u_2\in q^\mathbb{Z}$. Получаем
требуемое равенство операторов в $\mathbb{C}[t_1^{\pm 1},t_2^{\pm 1}]_q$.

\medskip

Разумеется, $U_q\mathfrak{sl}_2$-модульная алгебра $\mathbb{C}[t_1^{\pm
1},t_2^{\pm 1}]_q$ является $q$-аналогом $U\mathfrak{sl}_2$-модульной
алгебры $\mathbb{C}[t_1^{\pm 1},t_2^{\pm 1}]$.
\end{example}

\begin{example}\label{q-vector_space}
Гомоморфизм
$$
\mathbb{C}[z]\to\mathbb{C}[t_1^{\pm 1},t_2^{\pm 1}]_q,\qquad f(z)\mapsto
f(t_2^{-1}\,t_1),
$$
инъективен, поскольку $\mathbb{C}[t_1^{\pm 1},t_2^{\pm 1}]_q$ не имеет
делителей нуля, образ этого гомоморфизма является
$U_q\mathfrak{sl}_2$-подмодулем, и $\mathbb{C}[z]$ наследует структуру
$U_q\mathfrak{sl}_2$-модульной алгебры:
\begin{equation}\label{q-C}
K^{\pm 1}z=q^{\pm 2}z,\qquad Ez=-q^{1/2}z^2,\qquad Fz=q^{1/2}.
\end{equation}
Введем для этой $U_q\mathfrak{sl}_2$-модульной алгебры обозначение
$\mathbb{C}[z]_q$. Равенства \eqref{q-C} являются $q$-аналогами равенств
\eqref{1_on_z-action}. Из \eqref{q-C}, \eqref{K_Uq-algebra},
\eqref{E_Uq-algebra}, \eqref{F_Uq-algebra} получаем: $K^{\pm
1}f(z)=f(q^{\pm 1}z)$,
\begin{equation}\label{EF-q-diff}
Ef(z)=-q^{1/2}z^2\frac{f(z)-f(q^2z)}{z-q^2z},\quad
Ff(z)=q^{1/2}\frac{f(q^{-2}z)-f(z)}{q^{-2}z-z}.
\end{equation}
Это $q$-аналоги равенств \eqref{new-sl_2-action}.
\end{example}

\bigskip

До сих пор мы изучали голоморфные полиномы Лорана в $\mathbb{C}^2$ и их
квантовые аналоги. Перейдем к антиголоморфным полиномам Лорана. Пусть
$\mathbb{C}[(t_1^*)^{\pm 1},(t_2^*)^{\pm 1}]_q$ -- алгебра с образующими
$t_1^*$, $(t_1^*)^{-1}$, $t_2^*$, $(t_2^*)^{-1}$ и определяющими
соотношениями
$$
t_2^*t_1^*=qt_1^*t_2^*,\qquad(t_1^*)^{-1}t_1^*=t_1^*(t_1^*)^{-1}=1,\qquad
(t_2^*)^{-1}t_2^*=t_2^*(t_2^*)^{-1}=1.
$$
Рассмотрим антилинейный антиизоморфизм
$$
*:\mathbb{C}[(t_1^*)^{\pm 1},(t_2^*)^{\pm 1}]_q\to\mathbb{C}[t_1^{\pm
1},t_2^{\pm 1}]_q,\qquad *:t_1^*\mapsto t_1,\;*:t_2^*\mapsto t_2.
$$
Для обратного антилинейного антиизоморфизма будем использовать то же самое
обозначение $*$. Наделим $\mathbb{C}[(t_1^*)^{\pm 1},(t_2^*)^{\pm 1}]_q$
структурой $U_q\mathfrak{sl}_2$-модульной алгебры, перенеся ее из
$\mathbb{C}[t_1^{\pm 1},t_2^{\pm 1}]_q$:
\begin{equation}\label{antihol_pol}
af\stackrel{\rm def}{=}(S(a)^*f^*)^*,\qquad a\in U_q\mathfrak{sl}_2,\;
f\in\mathbb{C}[(t_1^*)^{\pm 1},(t_2^*)^{\pm 1}]_q.
\end{equation}
(Равенство \eqref{antihol_pol} наделяет $\mathbb{C}[(t_1^*)^{\pm
1},(t_2^*)^{\pm 1}]_q$ структурой $U_q\mathfrak{sl}_2$-модуля потому, что
антилинейный оператор $*S$ является антиавтоморфизмом коалгебры
$U_q\mathfrak{sl}_2$. Равенства \eqref{K_Uq-algebra}-\eqref{F_Uq-algebra}
вытекают из аналогичных равенств для $\mathbb{C}[t_1^{\pm 1},t_2^{\pm
1}]_q$ и из того, что $*S$ является антиавтоморфизмом коалгебры
$U_q\mathfrak{sl}_2$.)

Так же, как прежде, с помощью вложения
$$
\mathbb{C}[z^*] \hookrightarrow\mathbb{C}[(t_1^*)^{\pm 1},(t_2^*)^{\pm 1}]_q,
\qquad f(z^*)\mapsto f(t_1^*(t_2^*)^{-1})
$$
вводится $U_q\mathfrak{sl}_2$-модульная алгебра $\mathbb{C}[\bar{z}]_q$,
для которой
\begin{equation}\label{U_q_zstar_sl_2}
K^{\pm 1}z^*=q^{\mp 2}z^*,\qquad Ez^*=q^{-3/2},\qquad Fz^*=-q^{5/2}z^{*2}.
\end{equation}

\subsubsection{Универсальная $R$-матрица.}\label{RM_sl_2}

Модуль $V$ над алгеброй $U_q\mathfrak{sl}_2$ будем называть весовым,
\IND{модуль ! весовой} если
$$
V=\newoplus\limits_{\lambda\in\mathbb{R}}V_\lambda,\qquad V_\lambda=\{v\in
V\,|\,K^{\pm 1}\,v\,=\,q^{\pm \lambda}\,v\}.
$$
Все рассматриваемые в этой работе $U_q\mathfrak{sl}_2$-модули являются
весовыми.\footnote{Конечномерные весовые $U_q\mathfrak{sl}_2$-модули
называют также $U_q\mathfrak{sl}_2$-модулями типа 1 \cite{Jant}.}

\begin{example}\label{weight_example}
$\mathbb{C}[z]_q$ -- бесконечномерный весовой $U_q\mathfrak{sl}_2$-модуль,
а простейший невесовой $U_q\mathfrak{sl}_2$-модуль одномерен: $K^{\pm
1}=-1$, $E=F=0$.
\end{example}

 Для весового  $U_q\mathfrak{sl}_2$-модуля $V$  корректно определен
  линейный оператор
$$
H:V\rightarrow V, \qquad H|_{V_{\mathbf{\lambda}}} = \lambda.
$$
Очевидно, $ K^{\pm1}v=q^{\pm H}v$ при всех $v\in V$.

Полная подкатегория весовых $U_q\mathfrak{sl}_2$-модулей замкнута
относительно перехода к подмодулям,\footnote{См., например, \cite{Kac}
предложение 1.5} фактор-модулям и тензорным произведениям, то есть является
абелевой тензорной категорией. Каждый $U_q\mathfrak{sl}_2$-модуль обладает
наибольшим весовым подмодулем, что позволяет ввести понятие двойственного
весового $U_q\mathfrak{sl}_2$-модуля. \IND{двойственный ! весовой модуль}
Именно,
$$
V^*\stackrel{\operatorname{def}}{=}
\newoplus\limits_{\lambda\in\mathbb{R}}\;(V^*)_\lambda,
\qquad (V^*)_\lambda=(V_{-\lambda})^*,
$$
если
$V=\newoplus\limits_{\lambda\in\mathbb{R}}V_{\lambda}$.

\medskip

Пусть $B$ -- подалгебра алгебры $A$. Говорят, что $A$-модуль $V$ является
локально $B$-конечномерным, \IND{модуль ! локально конечномерный} если
$\dim(Bv)<\infty$ для всех $v\in V$.

Рассмотрим подалгебру Хопфа $U_q\mathfrak{b}^+$, порожденную элементами
$\{K^{\pm1},E\}$, и подалгебру Хопфа $U_q\mathfrak{b}^-$, порожденную
элементами $\{K^{\pm1},F\}$. Полную подкатегорию весовых локально
\hbox{$U_q\mathfrak{b}^+$-конечномерных}
\hbox{$U_q\mathfrak{sl}_2$-модулей} обозначим $\mathcal{C}^+$, а полную
подкатегорию весовых локально \hbox{$U_q\mathfrak{b}^-$-конечномерных}
\hbox{$U_q\mathfrak{sl}_2$-модулей} -- $\mathcal{C}^-$. \IND{категории
$\mathcal{C}^+$ и $\mathcal{C}^-$}

\bigskip

Рассмотрим $U_q\mathfrak{sl}_2$-модули $V'$, $V''$ и их тензорные
произведения $V'\otimes V''$, $V''\otimes V'$. Наивная перестановка
тензорных сомножителей
$$
\sigma_{V',V''}:V'\otimes V''\to V''\otimes V',\quad
\sigma_{V',V''}:v'\otimes v''\mapsto v''\otimes v',
$$
не является, вообще говоря, морфизмом $U_q\mathfrak{sl}_2$-модулей. В.
Дринфельд ввел понятие универсальной $R$-матрицы и с его помощью определил
линейные операторы, играющие роль перестановок $\sigma_{V',V''}$ в теории
$U_q\mathfrak{sl}_2$-модулей. Опишем эти линейные операторы.

Пусть
$$
\exp_q(t)=\sum\limits_{i=0}^{\infty}
\left(\prod\limits_{j=1}^i\frac{1-q}{1-q^j}\right)t^i
$$
и $V',V''$ -- весовые $U_q\mathfrak{sl}_2$-модули. Введем в рассмотрение
линейный оператор в $V'\otimes V''$ равенством
\begin{equation}\label{Rmatrix_sl_2}
R_{V'V''}(v'\otimes v'')\;=\;\exp_{q^2}((q^{-1}-q)E\otimes
F)q^{-\frac{H\otimes H}2}\;(v'\otimes v'').
\end{equation}
Он корректно определен, если либо $V'\in\mathcal{C}^+$, либо
$V''\in\mathcal{C}^-$.

 Как показывают следующие предложения, линейные операторы
$$\check{R}_{V'V''}=\sigma_{V'V''}\;
R_{V',V''}$$ являются $q$-аналогами перестановок
 тензорных сомножителей $\sigma_{V'V''}$.

\begin{proposition}\label{braiding1_sl_2}

  1. Рассмотрим такие весовые $U_q\mathfrak{sl}_2$-модули $V'$, $V''$,
  что либо
$V'\in \mathcal{C}^+$, либо $V''\in \mathcal{C}^-$. Линейный оператор
$$\check{R}_{V'V''}: V'\otimes V''\rightarrow V''\otimes V'$$ обратим и
является морфизмом
$U_q\mathfrak{sl}_2$-модулей. \\
  2. Если $f':V'\rightarrow W'$, $f'':V''\rightarrow W''$ -- морфизмы весовых
$U_q\mathfrak{sl}_2$-модулей, и либо $V',W'\in \mathcal{C}^+$, либо
$V'',W''\in \mathcal{C}^-$, то
$$
(f''\otimes f') \cdot \check{R}_{V',V''}=\check{R}_{W',W''}\cdot(f'\otimes
f'').
$$
\end{proposition}

\medskip
\begin{proposition}\label{braiding2_sl_2}

  1. Если $V$, $V'$, $V''$ -- весовые $U_q\mathfrak{sl}_2$-модули
  и либо $V',V''\in \mathcal{C}^+$, либо $V\in \mathcal{C}^-$, то
$$
\check{R}_{V'\otimes V'',V}=(\check{R}_{V',V}\otimes
\operatorname{id}_{V''})(\operatorname{id}_{V'}\otimes \check{R}_{V'',V}).
$$
  2. Если $V$, $V'$, $V''$ -- весовые $U_q\mathfrak{sl}_2$-модули и
  либо $V\in \mathcal{C}^+$, либо $V',V''\in \mathcal{C}^-$, то
$$
\check{R}_{V,V'\otimes V''}=(\operatorname{id}_{V'}\otimes
\check{R}_{V,V''})(\check{R}_{V,V'}\otimes \operatorname{id}_{V''}).
$$
  3. Для любого весового $U_q\mathfrak{sl}_2$-модуля $V$
$$
\check{R}_{V,\mathbb{C}}=\check{R}_{\mathbb{C},V}=\operatorname{id}_V.
$$
\end{proposition}

\medskip Из приведенных свойств операторов $\check{R}_{V',V''}$
следует, что $\mathcal{C}^-$ и $\mathcal{C}^+$ являются сплетенными
тензорными категориями, см. \itemiii\ \ref{braided_categories_sl_2}. Обычно
предложения \ref{braiding1_sl_2}, \ref{braiding2_sl_2} формулируются в
предположении о конечномерности переставляемых весовых
$U_q\mathfrak{sl}_2$-модулей \cite[стр. 329]{ChP}), но из доказательств
видно, что требование конечномерности можно ослабить.

\medskip

\begin{example}
Из определений следует, что $U_q\mathfrak{sl}_2$-модуль
$\mathbb{C}[z]_q$ является весовым и локально
$U_q\mathfrak{b}^-$-конечномерным. Имеет место равенство
\begin{equation}\label{commutativity_sl_2}
m\,\check{R}_{\mathbb{C}[z]_q\mathbb{C}[z]_q}(f_1\otimes f_2)\,=\,
\,f_1f_2,\qquad f_1,f_2\in\mathbb{C}[z]_q,
\end{equation}
где $m:\mathbb{C}[z]_q\otimes\mathbb{C}[z]_q\to\mathbb{C}[z]_q$ --
умножение в $\mathbb{C}[z]_q$. Достаточно доказать это равенство при
$f_1,f_2\in\{1,z\}$. Неочевиден лишь случай $f_1=f_2=z$:
$$
m\check{R}_{\mathbb{C}[z]_q\mathbb{C}[z]_q}(z\otimes z)=
m\left(1\,+\,(q^{-1}-q)F\otimes E)q^{-\frac{H\otimes H}{2}}(z\otimes
z)\right)=
$$
$$=m\left(q^{-2}z\otimes z\,-\,(q^{-2}-1)(1\otimes z^2)\right)=z^2.$$
Равенство \eqref{commutativity_sl_2} означает, что
$U_q\mathfrak{sl}_2$-модульная алгебра $\mathbb{C}[z]_q$ коммутативна в
теоретико-категорном смысле, см. \itemiii\ \ref{braided_categories_sl_2}.
\end{example}

\bigskip

Следующие замечания не претендуют ни на полноту, ни на строгость. Они
призваны пояснить связь сформулированных выше предложений
\ref{braiding1_sl_2}, \ref{braiding2_sl_2} с введенным {В. Дринфельдом}
\cite{DrinfEng} понятием универсальной $R$-матрицы. Хорошее изложение см. в
\cite[главы 16, 17]{Kassel_QG}.

Рассмотрим топологическое кольцо формальных степенных рядов
$\mathbb{C}[[h]]$ с $h$-адической топологией -- топологией с базой
окрестностей нуля $\{h^{\,n}\,\mathbb{C}[[h]]\}_{n\in \mathbb{Z}_+}$
\cite[стр. 481]{Kassel_QG}. (Отметим, что $\mathbb{C}[[h]]$ не является
комплексным топологическим векторным пространством в $h$-адической
топологии, поскольку билинейное отображение
$$\mathbb{C}\times V\to V[[h]],\qquad (z,v)\mapsto zv,$$
не непрерывно.)

Следуя традициям теории деформаций \cite{EK1}, введем в рассмотрение
аддитивную категорию $\mathcal{A}$ \IND{категория ! $\mathcal{A}$}
топологически свободных $\mathbb{C}[[h]]$-модулей. Каждому векторному
пространству $V$ сопоставим топологический $\mathbb{C}[[h]]$-модуль
$V[[h]]$ формальных рядов с коэффициентами из $V$ и $h$-адической
топологией -- топологией с базой окрестностей нуля
$\{h^{\,n}\,V[[h]]\}_{n\in\mathbb{Z}_+}$.

Топологические $\mathbb{C}[[h]]$-модули, изоморфные $V[[h]]$ для некоторого
векторного пространства $V$, называют топологически свободными. Это объекты
категории $\mathcal{A}$, а ее морфизмами являются все морфизмы
топологических $\mathbb{C}[[h]]$-модулей. Прямые суммы и тензорные
произведения в $\mathcal{A}$ вводятся равенствами
$$
V'[[h]]\;\newoplus\,V''[[h]]\stackrel{\rm def}{=}(V'\newoplus
V'')[[h]],\qquad V'[[h]]\;\hat{\otimes}\,V''[[h]]\stackrel{\rm
def}{=}(V'\otimes V'')[[h]]
$$
на объектах \cite[стр. 487]{Kassel_QG} и естественным образом -- на
морфизмах. ''Наивная'' перестановка тензорных сомножителей в коэффициентах
формальных рядов из $(V'\otimes V'')[[h]]$ наделяет полученную тензорную
категорию структурой симметрической тензорной категории, см. \itemiii\
\ref{braided_categories_sl_2}.

Следуя В. Дринфельду \cite{Drinfeld_KZ}, построим алгебру Хопфа
$U_h\mathfrak{sl}_2$ \IND{алгебра ! Хопфа ! $U_h\mathfrak{sl}_2$} в
категории $\mathcal{A}$.

  Рассмотрим свободную некоммутативную алгебру
  $\mathbb{C}\langle H,E,F \rangle$  с образующими $H$, $E$, $F$
  и подпространство
  $V \subset \mathbb{C}\langle H,E,F \rangle$, порожденное
  множеством $\{H^iE^jF^k\;|\; i,j,k \in \mathbb{Z}_+\}$.

Топологический подмодуль $V[[h]]\subset\mathbb{C}\langle F,H,E
\rangle[[h]]$ не является подкольцом, так как в произведении
$(F^{i'}H^{j'}E^{k'})(F^{i''}H^{j''}E^{k''})$ нарушен порядок образующих.
Перенесем образующую $H$ влево, а образующую $F$ вправо с помощью
коммутационных соотношений
$$
HE-EH=2E,\qquad HF-FH=-2F,\qquad EF-FE=\frac{{\rm sh}(hH/2)}{{\rm sh}(h/2)},
$$
где ${\rm
sh}(t)=\sum\limits_{n=0}^\infty\frac{t^{2n+1}}{(2n+1)!}\in\mathbb{C}[[t]]$
-- ряд Тейлора гиперболического синуса. Получаем операцию умножения,
наделяющую $U_h\mathfrak{sl}_2\,\stackrel{\rm def}{=}\,V[[h]]$ структурой
алгебры в тензорной категории $\mathcal{A}$.

Равенствами
$$
\Delta(H)\,=\,H\otimes 1\,+\,1\otimes H,\qquad \varepsilon(H)=0, \qquad
S(H)\,=\,-H,
$$
$$
 \Delta(E)=E\otimes 1 \;+\; e^{-hH/2}\otimes E,\qquad
  \varepsilon(E)=0, \qquad S(E)\,=\, -e^{hH/2}E,
$$
$$
\Delta(F)=F\otimes e^{hH/2}\;+\; 1\otimes F,\qquad
  \varepsilon(F)=0, \qquad S(F)\,=\, -F e^{-hH/2}.
$$
вводятся
 $\Delta: U_h\mathfrak{sl}_2\to U_h\mathfrak{sl}_2 \hat{\otimes}
U_h\mathfrak{sl}_2$, \hbox{$\varepsilon: U_h\mathfrak{sl}_2\to
\mathbb{C}[[h]]$} и $S: U_h\mathfrak{sl}_2 \to U_h\mathfrak{sl}_2$--
коумножение, коединица и антипод, наделяющие $U_h\mathfrak{sl}_2$ структурой
алгебры Хопфа в симметрической тензорной категории $\mathcal{A}$.

\begin{remark}\label{Change-remark-sl_2}
Алгебра Хопфа $U_h\mathfrak{sl}_2$ связана с $U_q\mathfrak{sl}_2$ следующими
``заменами образующих и параметра деформации'':
\begin{equation}\label{change1_sl_2}
q=e^{-h/2},\quad K^{\pm 1}=e^{\mp hH/2}.
\end{equation}
\end{remark}

\bigskip \medskip

Элемент
\begin{equation}\label{Rmatrix_h}
R\;=\;\exp_{q^2}((q^{-1}-q)E\otimes F)q^{-\frac{H\otimes H}{2}}
\end{equation}
тензорного квадрата $U_h\mathfrak{sl}_2\hat{\otimes}U_h\mathfrak{sl}_2$
называется универсальной $R$-матрицей алгебры Хопфа $U_h\mathfrak{sl}_2$.
\IND{универсальная ! $R$-матрица алгебры Хопфа $U_h\mathfrak{sl}_2$}
Применяя стандартные вложения
$$
i_{1,2},i_{1,3},i_{2,3}:\quad U_h\mathfrak{g}\hat{\otimes}
U_h\mathfrak{g}\;\to\;U_h\mathfrak{g}\hat{\otimes}
U_h\mathfrak{g}\hat{\otimes}U_h\mathfrak{g},
$$
получаем формальные ряды $R^{1,2}$, $R^{1,3}$, $R^{2,3}$. В этом контексте
аналоги предложений \ref{braiding1_sl_2}, \ref{braiding2_sl_2} выглядят
следующим образом:
\begin{equation}\label{R_property1_sl_2}
\Delta^{op}(\xi)=R\cdot\Delta(\xi)\cdot R^{-1},\quad\xi\in
U_q\mathfrak{sl}_2,
\end{equation}
\begin{equation}\label{R_property2_sl_2}
(\Delta\otimes\operatorname{id})(R)=R^{1,3}\cdot R^{2,3},
\end{equation}
\begin{equation}\label{R_property3_sl_2}
(\operatorname{id}\otimes\Delta)(R)=R^{1,3}\cdot R^{1,2},
\end{equation}
где $\Delta^{\mathrm{op}}$ -- противоположное коумножение. Алгебры Хопфа с
выделенным элементом тензорного квадрата, обладающим такими свойствами,
ввел Дринфельд. Их называют сплетенными алгебрами Хопфа \IND{алгебра !
Хопфа ! сплетенная} \cite[стр. 219]{Kassel_QG}.

\bigskip

 Рассмотрим антиавтоморфизм $\tau$ топологически свободной алгебры
   $U_h\mathfrak{sl}_2$,  определяемый равенствами
 $$\tau(H)\,=\,H,\qquad \tau(E)\,=\, e^{-\frac{hH}{2}}F,\qquad
 \tau(F)\,=\,  Ee^{\frac{hH}{2}}.$$
 Следующие часто используемые равенства получены в
  \cite[стр. 34, 37]{Drinf2}:
\begin{equation}\label{S_theta_sl_2}
   (S \otimes S)R=R, \qquad (\tau \otimes \tau)R=R^{21},
\end{equation}
где $R^{21}$ отличается от $R$ перестановкой тензорных сомножителей.

\subsubsection{$U_q\mathfrak{su}_{1,1}$-модульная
$*$-алгебра $\rm{Pol}(\mathbb{C})_q$.}\label{Pol_C_q}

В \itemiiiе \ref{sL_2_mod_algebras} было дано определение $A$-модульной
алгебры $F$. Если $A$ является $*$-алгеброй Хопфа и $F$ -- $*$-алгебра, то
в это определение включается дополнительное требование согласованности
инволюций: \IND{$A$-модульная $*$-алгебра}
\begin{equation}\label{invol2}
(af)^*=S(a)^*f^*
\end{equation}
при всех $a\in A$,\ $f\in F$.

\begin{example}
Рассмотрим $*$-алгебру Хопфа $A=(U\mathfrak{g},*)$ и отвечающую ей
вещественную подалгебру Ли $\{a\in\mathfrak{g}\,|\,a^*=-a\}$, см. \itemiii\
\ref{star-Hopf_algebras}. Пусть
\begin{itemize}
\item[-] $X$ -- однородное пространство группы Ли с этой вещественной алгеброй
Ли,
\item[-] $a$ -- один из ее элементов,
\item[-] $F$ -- алгебра гладких комплексных функций на
$X$,
\item[-] $f\in F$ -- ее самосопряженный элемент.
\end{itemize}
 В этом случае равенство \eqref{invol2} принимает вид
$(af)^*=(-a)^*f^*=af$ и означает, что отвечающее элементу $a$ вещественное
векторное поле на $X$, рассматриваемое как дифференциальный оператор
первого порядка, отображает вещественные гладкие функции в вещественные
гладкие функции.
\end{example}

\bigskip Нам предстоит вложить $U_q\mathfrak{sl}_2$-модульную
алгебру $\mathbb{C}[z]_q$ ''голоморфных полиномов'' в
$U_q\mathfrak{su}_{1,1}$-модульную $*$-алгебру ''всех полиномов на
квантовой комплексной плоскости''. Наделим структурой
$U_q\mathfrak{sl}_2$-модульной алгебры \hbox{$U_q\mathfrak{sl}_2$-модуль}
$\mathbb{C}[z]_q \otimes \mathbb{C}[\bar{z}]_q$.

Очевидно, $U_q\mathfrak{sl}_2$-модули $\mathbb{C}[z]_q$,
$\mathbb{C}[\bar{z}]_q$ являются весовыми и $U_q\mathfrak{sl}_2$-модуль
$\mathbb{C}[z]_q$ локально $U_q\mathfrak{b}^-$-конечномерен, а
$U_q\mathfrak{g}$-модуль $\mathbb{C}[\bar{z}]_q$ --- локально
$U_q\mathfrak{b}^+$-конечномерен. Значит, определен морфизм
$U_q\mathfrak{sl}_2$-модулей
$$
\check{R}=\check{R}_{\mathbb{C}[\bar{z}]_q\mathbb{C}[z]_q}:
\mathbb{C}[\bar{z}]_q\otimes\mathbb{C}[z]_q\to\mathbb{C}[z]_q\otimes
\mathbb{C}[\bar{z}]_q.
$$
Наделим $\mathbb{C}[z]_q\otimes\mathbb{C}[\bar{z}]_q$ структурой алгебры,
задав умножение
$$
m:(\mathbb{C}[z]_q\otimes\mathbb{C}[\bar{z}]_q)\otimes(\mathbb{C}[z]_q
\otimes\mathbb{C}[\bar{z}]_q)\to\mathbb{C}[z]_q\otimes\mathbb{C}[\bar{z}]_q
$$
с помощью умножений
$$
m^+:\;\mathbb{C}[\bar{z}]_q\otimes\mathbb{C}[\bar{z}]_q\to
\mathbb{C}[\bar{z}]_q,\qquad m^-:\;\mathbb{C}[z]_q\otimes\mathbb{C}[z]_q\to
\mathbb{C}[z]_q
$$
в алгебрах $\mathbb{C}[\bar{z}]_q$, $\mathbb{C}[z]_q$. Именно, положим
\begin{equation}\label{mult_pm_sl_2}
m=(m^-\otimes m^+)(\mathrm{id}_{\mathbb{C}[z]_q}\otimes\check{R}\otimes
\mathrm{id}_{\mathbb{C}[\bar{z}]_q}).
\end{equation}
Более подробно,
$$
\xymatrix{(\mathbb{C}[z]_q\otimes\mathbb{C}[\bar{z}]_q)\otimes
(\mathbb{C}[z]_q\otimes\mathbb{C}[\bar{z}]_q)\ar[d]^\cong
\\ \mathbb{C}[z]_q\otimes(\mathbb{C}[\bar{z}]_q\otimes\mathbb{C}[z]_q)
\otimes\mathbb{C}[\bar{z}]_q\ar[d]^{\mathrm{id}\otimes\check{R}\otimes
\mathrm{id}}
\\ \mathbb{C}[z]_q\otimes(\mathbb{C}[z]_q\otimes\mathbb{C}[\bar{z}]_q)
\otimes\mathbb{C}[\bar{z}]_q\ar[d]^\cong
\\ (\mathbb{C}[z]_q\otimes\mathbb{C}[z]_q)\otimes(\mathbb{C}[\bar{z}]_q
\otimes\mathbb{C}[\bar{z}]_q)\ar[d]^{m^-\otimes m^+}
\\ \mathbb{C}[z]_q\otimes\mathbb{C}[\bar{z}]_q}
$$
Из предложений \ref{braiding1_sl_2}, \ref{braiding2_sl_2} вытекают следующие
равенства:
$$
(m^-\otimes\mathrm{id})\check{R}_{\mathbb{C}[\bar{z}]_q,\,\mathbb{C}[z]_q
\otimes\mathbb{C}[z]_q}=\check{R}_{\mathbb{C}[\bar{z}]_q\mathbb{C}[z]_q}
(\mathrm{id}\otimes m^-),
$$
$$
(\mathrm{id}\otimes
m^+)\check{R}_{\mathbb{C}[\bar{z}]_q\otimes\mathbb{C}[\bar{z}]_q,\,
\mathbb{C}[z]_q}=\check{R}_{\mathbb{C}[\bar{z}]_q\mathbb{C}[z]_q}(m^+\otimes
\mathrm{id}),
$$
$$
\check{R}_{\mathbb{C}[\bar{z}]_q,\,\mathbb{C}[z]_q\otimes\mathbb{C}[z]_q}=
(\mathrm{id}_{\mathbb{C}[z]_q}\otimes\check{R}_{\mathbb{C}[\bar{z}]_q
\mathbb{C}[z]_q})\cdot(\check{R}_{\mathbb{C}[\bar{z}]_q\mathbb{C}[z]_q}
\otimes\mathrm{id}_{\mathbb{C}[z]_q}),
$$
$$
\check{R}_{\mathbb{C}[\bar{z}]_q\otimes\mathbb{C}[\bar{z}]_q,\,
\mathbb{C}[z]_q}=(\check{R}_{\mathbb{C}[\bar{z}]_q\mathbb{C}[z]_q}\otimes
\mathrm{id}_{\mathbb{C}[\bar{z}]_q})\cdot
(\mathrm{id}_{\mathbb{C}[\bar{z}]_q}\otimes\check{R}_{\mathbb{C}[\bar{z}]_q
\mathbb{C}[z]_q}).
$$
C их помощью, используя ассоциативность умножений $m^\pm$ в алгебрах
$\mathbb{C}[\bar{z}]_q$, $\mathbb{C}[z]_q$, легко доказать ассоциативность
умножения $m$. (Это доказательство хорошо известно в теории сплетенных
тензорных категорий, см. \cite[стр. 439]{Majid_book}.) Остальные
утверждения следующего предложения вытекают из предложения
\ref{braiding1_sl_2}.

\begin{proposition}\label{assoc_sl_2}
Умножение $m$ наделяет $\mathbb{C}[z]_q \otimes \mathbb{C}[\bar{z}]_q$
 структурой $U_q\mathfrak{sl}_2$-модульной алгебры с единицей $1 \otimes 1$.
\end{proposition}

\bigskip
Из \eqref{Rmatrix_sl_2},\eqref{q-C}, \eqref{U_q_zstar_sl_2} вытекают
равенства
$$ \check{R}_{\mathbb{C}[\bar{z}]_q
\mathbb{C}[z]_q}(z^*\otimes z)\,=\,
 (1\, +\, (q^{-1}-q) F\otimes E)\,q^{-\frac{H\otimes H}{2}} (z\otimes
 z^*)\,=\,
$$
$$ =\,q^2z\otimes z^*\,+\, 1-q^2.$$

Значит,
\begin{equation}\label{Pol_revisited}
(1\otimes z^*)(z\otimes 1)\;=\; q^2\,(z\otimes 1)(1\otimes z^*)\,+\, 1-q^2.
\end{equation}

Следовательно, отображение
$$z\mapsto z\otimes 1,\qquad z^*\mapsto 1\otimes z^*$$
единственным образом продолжается до изоморфизма алгебры ${\rm
Pol}(\mathbb{C})_q$ и полученной алгебры
$(\mathbb{C}[z]_q\otimes\mathbb{C}[\bar{z}]_q,m)$. Наделяя эту алгебру
инволюцией
$$
*:\;z\otimes 1\mapsto 1\otimes z^*,\qquad*:1\otimes z^*\mapsto z\otimes 1,
$$
получаем изоморфизм $*$-алгебр. Равенство \eqref{Pol_revisited} можно
записать более лаконично: $z^*z=q^2zz^*+1-q^2$.

 Будем отождествлять ${\rm
Pol}(\mathbb{C})_q$ с построенной в этом
\itemiiiе $*$-алгеброй.

Докажем согласованность инволюций -- равенство \eqref{invol2}. Рассмотрим
представление
$$\pi(a)f=af,\qquad a\in U_q\mathfrak{sl}_2,\;f\in{\rm Pol}(\mathbb{C})_q,$$
алгебры $U_q\mathfrak{sl}_2$ в ${\rm Pol}(\mathbb{C})_q$ и ее представление
$$
\pi'(a)f=(\pi(S(a)^*)f^*)^*,\qquad a\in U_q\mathfrak{sl}_2,\;f\in{\rm
Pol}(\mathbb{C})_q,
$$
в том же пространстве. Мы воспользовались тем, что $*S$ является
антилинейным автоморфизмом алгебры $U_q\mathfrak{g}$. Представление $\pi'$,
как и представление $\pi$, наделяет ${\rm Pol}(\mathbb{C})_q$ структурой
$U_q\mathfrak{sl}_2$-модульной алгебры, поскольку антилинейный оператор
$*S$ является антилинейным антиавтоморфизмом коалгебры $U_q\mathfrak{g}$
(ср. с аналогичными рассуждениями в \itemiiiе \ref{sL_2_mod_algebras}).
Следовательно, равенство операторов
$$\pi(K^{\pm 1})=\pi'(K^{\pm 1}),\qquad\pi(E)=\pi'(E),\qquad\pi(F)=\pi'(F)$$
вытекает из их равенства на трехмерном подпространстве, порожденном
элементами $1$, $z$, $z^*$, а оно, в свою очередь, легко следует из
определений.

\bigskip

В \itemiiiе \ref{Fock_l} $*$-алгебра ${\rm Pol}(\mathbb{C})_q$ была выбрана
как-бы наугад, без серьезной мотивации. Описанный в этом \itemiiiе путь к
$*$-алгебре ${\rm Pol}(\mathbb{C})_q$ \IND{$*$-алгебра !
$\mathrm{Pol}(\mathbb{C})_q$} не только приводит к ней с необходимостью, но
и наделяет ее структурой $U_q\mathfrak{su}_{1,1}$-модульной алгебры.

\subsubsection{Квантовая группа $SL_2$.}\label{SL_2_sl_2}

Напомним, что $q\in(0,1)$. Впрочем, результаты настоящего \itemiiiа верны
для всех чисел $q$, не являющихся корнями из единицы \cite[главы 3,
4]{KlSch}, \cite[главы 6, 7]{Kassel_QG}.

Пусть $n\in\mathbb{Z}_+$ и $\bar{\omega}$ -- фундаментальный вес алгебры Ли
$\mathfrak{sl}_2$, то есть $\bar{\omega}(H)=1$. Введем обозначение
$L(n\bar{\omega})$ для $U_q\mathfrak{sl}_2$-модуля с образующей
$v(n\bar{\omega})$ и определяющими соотношениями
\begin{equation}\label{L_lambda_sl_2}
E\,v(n\bar{\omega})=0,\qquad K^{\pm 1}\,v(n\bar{\omega})=q^{\pm
n}v(n\bar{\omega}),\qquad F^{\,n+1}\,v(n\bar{\omega})=0.
\end{equation}

 Следующее утверждение вытекает из
 \eqref{U_q_first}, \eqref{U_q_second}.

\begin{lemma}\label{sl_2_center} Элемент
\begin{equation}\label{C_q}
C_q = FE + \frac{1}{(q^{-1}-q)^2} \left(q^{-1} K^{-1} + qK -
(q^{-1}+q)\right)
\end{equation}
принадлежит центру алгебры $U_q\mathfrak{sl}_2$.
\end{lemma}

Очевидно,
\begin{equation}\label{C_q_second_form}
C_q = EF + \frac{1}{(q^{-1}-q)^2} \left(q^{-1} K + qK^{-1} -
(q^{-1}+q)\right)
\end{equation}

 Из  леммы  \ref{sl_2_center} и из равенства
$$C_q\,v(n\bar{\omega}) = \frac{(q^{-\frac{n}{2}}-q^\frac{n}{2})
(q^{-(\frac{n}{2}+1)}-q^{\frac{n}{2}+1})} {(q^{-1}-q)^2}v(n\bar{\omega})
$$
получаем

\begin{corollary}
\begin{equation}\label{eigenvalues_C_q}
C_q|_{L(n\bar{\omega})}=\frac{(q^{-\frac{n}2}-q^\frac{n}2)
(q^{-(\frac{n}2+1)}-q^{\frac{n}2+1})}{(q^{-1}-q)^2}.
\end{equation}
\end{corollary}

Это $q$-аналог хорошо известной формулы для действия элемента Казимира
$C=FE+\frac{H(H+2)}4$ алгебры Ли $\mathfrak{sl}_2$ в конечномерных
$\mathfrak{sl}_2$-модулях \cite[стр. 149]{Zhelob}. \IND{$q$-аналог ! элемента
Казимира}

\begin{lemma}\label{sl_2_modules}
1. $L(n\bar{\omega})$ -- простой весовой $U_q\mathfrak{sl}_2$-модуль
размерности $n+1$.
\\ 2. Каждый конечномерный простой весовой $U_q\mathfrak{sl}_2$-модуль
изоморфен одному из $L(n\bar{\omega})$.
\\ 3. В категории $U_q\mathfrak{sl}_2$-модулей
\begin{equation}\label{CG_sl_2}
L(n_1\bar{\omega})\otimes
L(n_2\bar{\omega})\approx\newoplus\limits_{n=0}^{{\rm min}(n_1,n_2)}
L((n_1+n_2-2n)\bar{\omega}),\qquad n_1,n_2\in\mathbb{Z}_+.
\end{equation}
\end{lemma}

Доказательство аналогично доказательству в классическом случае $q=1$, см.
\cite[стр. 61-62, 72]{KlSch}.

\begin{corollary}\label{separation}
Если $L'$ и $L''$ -- неизоморфные конечномерные простые весовые
$U_q\mathfrak{sl}_2$-модули, то $C_q|_{L'}\ne C_q|_{L''}$.
\end{corollary}

Используя лемму \ref{sl_2_modules}, нетрудно показать \cite[стр.
65]{KlSch}, что каждый конечномерный весовой $U_q\mathfrak{g}$-модуль
полупрост, то есть, что каждый его подмодуль дополняем. Используя
полупростоту, легко получить (см. доказательство предложения
\ref{semisimple1}) следующий результат.

\begin{lemma}\label{loc_finite_sl_2} Если весовой
$U_q\mathfrak{sl}_2$-модуль $V$ локально конечномерен, то есть $ {\rm
dim}(U_q\mathfrak{sl}_2\cdot v) < \infty$ при всех $v \in V$, то он
изоморфен прямой сумме простых весовых конечномерных
$U_q\mathfrak{sl}_2$-модулей.
\end{lemma}


\bigskip
Матричными элементами представления $A\to{\rm End}(V)$ алгебры $A$ называют
элементы образа сопряженного линейного оператора ${\rm End}(V)^*\to A^*$.
\IND{матричный элемент представления} Прямая сумма и тензорное произведение
конечномерных весовых $U_q\mathfrak{sl}_2$-модулей являются конечномерными
весовыми $U_q\mathfrak{sl}_2$-модулями. Значит, матричные элементы таких
$U_q\mathfrak{sl}_2$-модулей образуют подалгебру Хопфа
$\mathbb{C}[SL_2]_q\subset(U_q\mathfrak{sl}_2)^\star$, см. замечание
\ref{duality_of_Hopf}. Ее называют алгеброй регулярных функций на квантовой
группе $SL_2$ \cite[стр. 96]{KorSoib}. \IND{алгебра ! регулярных функций !
на квантовой группе $SL_2$}

\medskip

Простота и конечномерность $L(n\bar{\omega})$ позволяют с помощью теоремы
Бернсайда \cite[стр. 177]{CurtisReiner}, \cite[стр. 100]{Zhelob} доказать,
что естественный гомоморфизм коалгебр ${\rm End}(L(n\bar{\omega}))^* \to
\mathbb{C}[SL_2]_q$ инъективен. Будем отождествлять векторное пространство
${\rm End}(L(n\bar{\omega}))^*$ с его образом при вложении в
$\mathbb{C}[SL_2]_q$. Получим аналог разложения Петера--Вейля.

\begin{proposition}\label{PeterWeyl_sl_2}
  $\mathbb{C}[SL_2]_q \cong
  \newoplus\limits_{n=0}^\infty
   \operatorname{End}(L(n\bar{\omega}))^*$.
\end{proposition}

{\bf Доказательство.} Из лемм \ref{sl_2_modules}, \ref{loc_finite_sl_2}
следует, что $\mathbb{C}[SL_2]_q$ -- сумма своих подпространств
$\operatorname{End}(L(n\bar{\omega}))^*$. Требуется доказать их линейную
независимость. Но они являются собственными подпространствами с попарно
различными собственными значениями линейного оператора
$$f(\xi)\mapsto f(\xi\,C_q),\qquad \xi \in U_q\mathfrak{sl}_2,$$
в векторном пространстве $\mathbb{C}[SL_2]_q$, как нетрудно показать с
помощью следствия \ref{separation}. \hfill $\square$

\begin{remark} \label{reg_actions_sl_2}
Пусть $U_q\mathfrak{sl}_2^\mathrm{op}$ -- алгебра Хопфа, отличающаяся от
$U_q\mathfrak{sl}_2$ заменой умножения на противоположное. Наделим
$\mathbb{C}[SL_2]_q$ структурой $U_q\mathfrak{sl}_2$-бимодуля, задав
действие $U_q\mathfrak{sl}_2$ и действие $U_q\mathfrak{sl}_2^{\mathrm{op}}$
следующим образом:
\begin{equation}\label{Rreg_Lreg_sl_2}
R_{\operatorname{reg}}(\xi):f(\eta)\mapsto f(\eta\xi),\qquad
L_{\operatorname{reg}}(\xi):f(\eta)\mapsto f(\xi\eta),
\end{equation}
где $f\in \mathbb{C}[SL_2]_q$, $\xi,\eta\in U_q\mathfrak{sl_2}$. Из
определений следует, что $\mathbb{C}[SL_2]_q$ является
  $U_q\mathfrak{sl}_2^{\mathrm{op}} \otimes U_q\mathfrak{sl}_2$-модульной алгеброй.
Предложение \ref{PeterWeyl_sl_2} доставляет разложение $\mathbb{C}[SL_2]_q$
в сумму простых попарно неизоморфных $U_q\mathfrak{sl}_2^{\rm op} \otimes
U_q\mathfrak{sl}_2$-модулей. (В дальнейшем мы будем отдавать предпочтение
представлению $R_{reg}$ перед представлением $L_{reg}$ и использовать
сокращенное обозначение $\xi f=R_{\operatorname{reg}}(\xi)f$.)
\end{remark}

\bigskip

Часто используется другое, более наглядное описание алгебры Хопфа
$\mathbb{C}[SL_2]_q$ -- ее описание в терминах образующих и соотношений.
Пусть $\{e_1,\, e_2\}$ -- стандартный базис эрмитова пространства
$\mathbb{C}^2$. Рассмотрим представление $\pi$ алгебры $U_q\mathfrak{sl}_2$
в $\mathbb{C}^2$, введенное равенствами \eqref{q2-dim_rep}, и его матричные
элементы $t_{ij}(\xi)\stackrel{\rm def}{=}(\pi(\xi)e_j,e_i)$. Нетрудно
доказать следующие утверждения, см. \cite{VakSoib88}.

\begin{proposition}\label{su_2_relations}
Множество $\{t_{11},t_{12},t_{21},t_{22}\}$ порождает алгебру
$\mathbb{C}[SL_2]_q$, и
\begin{align*}
t_{11}t_{12}=qt_{12}t_{11}, & \quad t_{11}t_{21}=qt_{21}t_{11},\quad
t_{12}t_{22}=qt_{22}t_{12},\quad t_{21}t_{22}=q\,t_{22}t_{21},
\\ \quad t_{12}t_{21}=t_{21}t_{12}, & \quad
t_{11}t_{22}-t_{22}t_{11}=(q-q^{-1})\,t_{12}t_{21},
\end{align*}
\begin{equation}\label{q-det_sl_2}
t_{11}t_{12}-qt_{12}t_{21}=1.
\end{equation}
Это множество соотношений между образующими $\{t_{ij}\}_{i,j=1,2}$ является
определяющим.\footnote{Другими словами, отвечающее ему множество элементов
свободной алгебры $\mathbb{C}<t_{11},t_{12},t_{21},t_{22}>$ порождает ядро
ее естественного гомоморфизма в $\mathbb{C}[SL_2]_q$.} Действие коумножения
$\Delta$, коединицы $\varepsilon$ и антипода $S$ описываются следующими
равенствами:
$$
\Delta(t_{ij})=t_{i1}\otimes t_{1j}+t_{i2}\otimes t_{2j},\quad
\varepsilon(t_{i\,j})=\delta_{i,j},\quad i,j=1,2,
$$
\begin{equation}\label{anti2}
\left(\begin{array}{cc}S(t_{11}), & S(t_{12})\\ S(t_{21}), & S(t_{22})
\end{array}\right)=
\left(\begin{array}{cc}t_{22}, & -q^{-1}t_{12}\\ -qt_{21}, &
t_{11}\end{array}\right).
\end{equation}
\end{proposition}

\bigskip
 Поясним  первое из этих утверждений.
 Как вытекает из \eqref{CG_sl_2}, простой $U_q\mathfrak{sl}_2$-модуль
$L(n\bar{\omega})$ допускает вложение в $L(\bar{\omega})^{\otimes\, n}$.
Значит, из предложения \ref{PeterWeyl_sl_2} следует, что подпространство
${\rm End}(L(\bar{\omega}))^*$ порождает алгебру $\mathbb{C}[SL_2]_q$. Но
это подпространство является линейной оболочкой множества $\{t_{11}, t_{12},
t_{21}, t_{22}\}$.

\begin{corollary}\label{*-reps_sl_2}
1. Равенства
\begin{align*}
\pi_+(t_{11})e_j&=\sqrt{q^{-2(j+1)}-1}\,e_{j+1}, &
\pi_+(t_{12})e_j&=q^{-j}e_j,
\\ \pi_+(t_{21})e_j&=-q^{-(j+1)}e_j, & \pi_+(t_{22})e_j&=
\begin{cases}
-\sqrt{q^{-2j}-1}\,e_{j-1}, & j\ne 0,
\\ 0, & j=0
\end{cases}
\end{align*}
определяют неприводимое бесконечномерное представление $\pi_+$ алгебры
$\mathbb{C}[SL_2]_q$ в векторном пространстве с базисом $\{e_j\}_{j \in
\mathbb{Z}_+}$.

2. Равенства
$$
 \rho(t_{11})f(z)=zf(z),\quad \rho(t_{12})f(z)=\rho(t_{21})f(z)=0,
\quad \rho(t_{22})f(z)=z^{-1}f(z)
$$
определяют бесконечномерное представление $\rho$ алгебры
$\mathbb{C}[SL_2]_q$ в векторном пространстве $\mathbb{C}[z^{\pm 1}]$
полиномов Лорана переменной $z$.
\end{corollary}

С помощью представления $\pi_+\otimes \rho$ легко доказать линейную
независимость элементов
\begin{equation}\label{PBW_sl_2}
\{t_{11}^at_{12}^bt_{21}^c,\; t_{12}^bt_{21}^ct_{22}^d\quad|\quad a,b,c,d
\in \mathbb{Z}_+\}
\end{equation}
алгебры $\mathbb{C}[SL_2]_q$. Значит, имеет место

\begin{corollary}\label{PBW_SL_2}
Множество \eqref{PBW_sl_2} является базисом векторного пространства
$\mathbb{C}[SL_2]_q$.
 \end{corollary}

\begin{remark}\label{integrity_sl_2} Используя приведенное выше
 описание алгебры $\mathbb{C}[SL_2]_q$ в
терминах образующих и соотношений, нетрудно доказать, что эта алгебра
является нетеровой областью целостности \cite[стр. 18]{BrGood}.
\end{remark}

\medskip

Рассмотрим коалгебру $A$ с коумножением $\triangle_A$ и коединицей
$\varepsilon_A$. Векторное пространство $F$ называют (правым)
$A$-комодулем, \IND{комодуль} если задано такое линейное отображение
$\triangle_F: F\to F\otimes A$, что
$$
(\triangle_F\otimes\operatorname{id}_A)\cdot\triangle_F=
(\operatorname{id}_F\otimes\triangle_A)\cdot\triangle_F,\qquad
(\operatorname{id}_F\otimes\varepsilon_A)\cdot\triangle_F=
\operatorname{id}_F.
$$

Каждому $U_q\mathfrak{sl}_2$-модулю $F$ отвечает линейное отображение
$${F\to\operatorname{Hom}(U_q\mathfrak{sl}_2,F)},$$ сопоставляющее вектору
$f\in F$ линейный оператор
 $$U_q\mathfrak{sl}_2 \to F,\qquad a\mapsto af.$$
Если весовой \hbox{$U_q\mathfrak{sl}_2$-модуль} $F$ локально конечномерен,
то образ этого отображения лежит в подпространстве ${F\otimes
\mathbb{C}[SL_2]_q\hookrightarrow\operatorname{Hom}(U_q\mathfrak{sl}_2,F)}$,
и возникающее линейное отображение $\triangle_F:F\to F\otimes
\mathbb{C}[SL_2]_q$ наделяет $F$ структурой $\mathbb{C}[SL_2]_q$-комодуля.
Это линейное отображение является гомоморфизмом унитальных алгебр, если $F$
-- унитальная $U_q\mathfrak{sl}_2$-модульная алгебра с весовым и локально
конечномерным действием $U_q\mathfrak{sl}_2$ в $F$. Тем самым $F$
наделяется структурой $\mathbb{C}[SL_2]_q$-комодульной алгебры.
\IND{алгебра ! комодульная}
\begin{example}\label{like_Kassel} $\mathbb{C}[t_1,t_2]_q$ является
$\mathbb{C}[SL_2]_q$ комодульной алгеброй:
$$\Delta:\; \mathbb{C}[t_1,t_2]_q \to \mathbb{C}[t_1,t_2]_q \otimes
\mathbb{C}[SL_2]_q, \qquad \Delta:\;t_j \mapsto t_1 \otimes t_{1j} +
t_2\otimes t_{2j}.$$
\end{example}

Важно отметить, что $U_q\mathfrak{sl}_2$-модули $\mathbb{C}[t_1^{\pm
1},t_2^{\pm 1}]_q$, $\mathbb{C}[z]_q$, играющие ключевую роль в
рассматриваемом круге вопросов, не являются локально конечномерными.
Поэтому мы в дальнейшем отдаем предпочтение языку
$U_q\mathfrak{sl}_2$-модульных алгебр даже в тех случаях, когда действие
алгебры $U_q\mathfrak{sl}_2$ в $F$ является весовым и локально
конечномерным, ср. \cite[глава 4]{Kassel_QG}.

\bigskip


\subsubsection{Дополнение о кольцах и модулях частных.}\label{localization}

Пусть $S$ -- мультипликативное подмножество кольца $R$.
\IND{мультипликативное подмножество кольца} Это означает, что $S\cdot
S\subset S$ и $1\in S$. Кольцо $R_S$ вместе с гомоморфизмом $\varphi:R\to
R_S$ называют кольцом частных \IND{кольцо частных} $R$ относительно $S$,
если выполнено следующее требование универсальности.

Для любого гомоморфизма колец $\psi:R\to Q$ из обратимости всех элементов
$\psi(s)$, $s\in S$, вытекает существование и единственность такого
гомоморфизма $\psi_S:R_S\to Q$, что $\psi=\psi_S\circ\varphi$.

Если $\varphi':R\to R_S'$, $\varphi'':R\to R_S''$ -- два кольца частных, то
существует и единствен такой изоморфизм колец $f:R_S'\to R_S''$, что
$\varphi''=f\, \varphi'$. Другими словами, кольцо частных в существенном
единственно.

При доказательстве его существования можно считать, что кольцо $R$ задано
своими образующими и определяющими соотношениями. Сопоставим каждому
элементу $s\in S$ дополнительную образующую $\overline{s}$ и два
дополнительных соотношения $s\overline{s}=1$, $\overline{s}s=1$. Пусть $R_S$
-- кольцо, определяемое расширенными таким образом множествами образующих и
соотношений. Естественное вложение множеств образующих приводит к
гомоморфизму $\varphi:R\to R_S$. Очевидно, это локализация $R$ относительно
$S$.

Если $R$ -- градуированное кольцо и $S$ состоит из однородных элементов, то
кольцом частных называют градуированное кольцо $R_S$ вместе с гомоморфизмом
$\varphi:R\to R_S$ градуированных колец, если выполнено аналогичное
приведенному выше требование универсальности. Существование и единственность
кольца частных доказывается в этом случае так же, как прежде.

Если $M$ -- правый $R$-модуль, то $M_S=M\otimes_RR_S$ является правым
$R_S$-модулем. Его называют модулем частных.

\bigskip
В приведенных выше определениях мультипликативное множество $S$ может быть
произвольным, и кольцо частных $R_S$ обладает неконтролируемыми свойствами.
Введем дополнительные определения и связанные с ними ограничения на $S$
\cite[стр. 299]{Lam}.

Кольцо частных $R_S$ вместе с гомоморфизмом $\varphi:R\to R_S$ называют
классическим {\bf правым} кольцом частных, если
\begin{itemize}
\item[1.] $\varphi(s)$ -- обратимый элемент кольца $R_S$ для всех $s\in S$,

\item[2.] каждый элемент $R_S$ имеет вид $\varphi(r)\varphi(s)^{-1}$, $r\in
R$, $s\in S$,

\item[3.] $\operatorname{Ker}\varphi=\{r\in R|\:rs=0\text{\ при некотором\
}s\in S\}$.
\end{itemize}

Следующий результат хорошо известен.

\begin{proposition}(\cite[стр. 110]{Bour_com}, \cite[стр. 300]{Lam})
\label{Ore_prop} Кольцо $R_S$ является классическим правым кольцом частных
в том и только в том случае, когда $S$ обладает следующими свойствами:
\begin{flalign}\label{Ore_right}
1.\; && rS\cap sR\ne 0\text{\ \ для всех\ \ }s\in S,\;r\in R, &&
\end{flalign}
$2.\, $ если $s'r=0$ при некоторых $s'\in S$, $r\in R$, то $rs=0$ для
некоторого $s \in S$.
\end{proposition}

Точно так же вводятся классические левые кольца частных, и для них имеет
место аналогичное утверждение. Роль \eqref{Ore_right} играет условие
\begin{flalign}\label{Ore_left}
&& Sr\cap Rs\ne\varnothing\text{\ \ для всех\ \ }s\in S,\;r\in R. &&
\end{flalign}
В большинстве приложений условия \eqref{Ore_left}, \eqref{Ore_right},
называемые условиями Оре, \IND{условия Оре} выполняются, и $S$ состоит из
неделителей нуля.


\subsubsection{Дополнение об алгебрах и модулях в тензорных категориях.}
\label{braided_categories_sl_2}

Многие из перечисленных в \itemiiiе \ref{RM_sl_2} свойств универсальной
$R$-матрицы допускают теоретико-категорную интерпретацию \cite{Joyal}.

 Например, как объясняется в заключительной части настоящего
\itemiiiа, введенная в \itemiiiе \ref{RM_sl_2} категория $\mathcal{C}^-$
является абелевой сплетенной тензорной категорией. Напомним, что
$U_q\mathfrak{sl}_2$-модуль принадлежит этой категории, если и только если
он является весовым и локально $U_q\mathfrak{b}^-$-конечномерным.

Алгеброй в категории $\mathcal{C}^-$ \IND{алгебра ! в категории
$\mathcal{C}^-$} называют унитальную $U_q\mathfrak{sl}_2$-модульную
алгебру, принадлежащую $\mathcal{C}^-$. Левым модулем над алгеброй $F$ в
категории $\mathcal{C}^-$ \IND{модуль ! над алгеброй ! в категории
$\mathcal{C}^-$} называют $U_q\mathfrak{g}$-модульный левый $F$-модуль,
принадлежащий $\mathcal{C}^-$. Аналогично вводятся понятия правого модуля и
бимодуля над алгеброй $F$ в категории $\mathcal{C}^-$.

\medskip

Алгебра $F$ в категории $\mathcal{C}^-$ называется коммутативной в этой
категории, \IND{алгебра ! коммутативная ! в категории $\mathcal{C}^-$} если
$m\;=\;m\,\check{R}_{FF}$, где $m:F\otimes F\to F$ -- умножение в $F$.

Рассмотрим бимодуль $E$ над коммутативной алгеброй $F$ в категории
$\mathcal{C}^-$:
$$m_{\rm left}:F\otimes E\to E,\qquad m_{\rm right}:E\otimes F\to E.$$
Следуя \cite{LychaginPrasolov}, будем называть его симметричным, если
$$m_{\rm right}=m_{\rm left}\,\check{R}_{EF}.$$
\IND{симметричный бимодуль ! над коммутативной алгеброй в категории
$\mathcal{C}^-$}

\begin{example}\label{subalgebra_in_category}
Алгебра $F$ в категории $\mathcal{C}^-$ является бимодулем над любой своей
подалгеброй $F'$. Если $F$ -- коммутативная алгебра в $\mathcal{C}^-$, то
этот бимодуль симметричен.
\end{example}

\begin{proposition}\label{tensor_in_category}\cite[Замечание 1.1]{Pareigis}
Если $E_1$, $E_2$ -- симметричные бимодули над алгеброй $F$ в категории
$\mathcal{C}^-$, то бимодуль $E_1\otimes_F E_2$ также симметричен.
\end{proposition}


\begin{proposition}\label{to_bimodule}\cite[Замечание 1.2]{Pareigis}
Пусть $E$ -- левый модуль \hbox{$m_{\rm left}:F\otimes E\to E$} над
коммутативной алгеброй $F$ в $\mathcal{C}^-$. Морфизм
$$m_{\rm right}=m_{\rm left}\check{R}_{EF}$$
наделяет $E$ структурой симметричного бимодуля над $F$ в $\mathcal{C}^-$.
\end{proposition}

\bigskip

Это означает, что каждый левый модуль $E$ над коммутативной алгеброй $F$ в
$\mathcal{C}^-$ является симметричным бимодулем над $F$ в $\mathcal{C}^-$.


\bigskip

Бимодуль $E$ над коммутативной алгеброй $F$ в категории $\mathcal{C}^-$
называют локальным, если \IND{локальный бимодуль ! над коммутативной
алгеброй в категории $\mathcal{C}^-$}
$$
m_{\rm right}=m_{\rm left}\,\check{R}_{EF},\qquad m_{\rm left}=m_{\rm
right}\,\check{R}_{FE}
$$
(этот вольный перевод введенного Парейгисом в \cite{Pareigis} термина
''dysleсtic'' используется в работах по конформной квантовой теории поля).

\begin{proposition}\label{dislectic}
\cite[предложение 2.4, теорема 2.5]{Pareigis} Пусть $E_1$, $E_2$--
локальные бимодули над коммутативной алгеброй $F$ в категории
$\mathcal{C}^-$ и $j$ -- канонический эпиморфизм $E_1\otimes E_2\to
E_1\otimes_F E_2$. Тогда
\begin{itemize}
\item[1.]\ $E_1\otimes_F E_2$ является локальным бимодулем над $F$;

\item[2.]\  существует и единствен такой морфизм
$$\bar{R}_{E_1E_2}:\;E_1\otimes_FE_2\to E_2\otimes_FE_1,$$
в категории $\mathcal{C}^-$, что
$\bar{R}_{E_1E_2}\,j=j\;\check{R}_{E_1E_2}$.
\end{itemize}
\end{proposition}

\medskip

Это означает, что категория локальных бимодулей над коммутативной алгеброй
в категории $\mathcal{C}^-$ является абелевой сплетенной тензорной
категорией (см. определения в заключительной части \itemiiiа).

\medskip Сформулированные утверждения
  получены  Парейгисом методами теории категорий в очень большой
  общности.

\bigskip
Обсудим некоторые понятия, не упоминаемые в общих курсах теории категорий,
например, в \cite{BucurDeleanu}.

Начнем с тензорных категорий. \IND{категория ! тензорная} Основным примером
служит категория $\mathcal{C}_A$ модулей над алгеброй Хопфа $A$.

Рассмотрим категорию $\mathcal{C}$ с выделенным объектом $I$ и функтор
\hbox{$\otimes: \mathcal{C}\times \mathcal{C} \to\mathcal{C}$}. Чтобы задать
такой функтор, нужно ввести тензорные произведения объектов и тензорные
произведения морфизмов так, что
\begin{equation*}\label{knots-2.6}
\rm{id}_{U\otimes V}\;=\; \rm{id}_U \otimes \rm{id}_V, \qquad (f''\otimes
g'')(f'\otimes g')\;=\; f''f' \otimes g''g',
\end{equation*}
 для любых двух объектов $U,V$ категории $\mathcal{C}$ и любых ее морфизмов
 $$f' \in \rm{Hom}(U,U'),\; f'' \in \rm{Hom}(U',U''),\; g'\in
\rm{Hom}(V,V'),\; g'' \in \rm{Hom}(V',V'').$$ Подробности см. в \cite[стр.
350]{Kassel_QG}, \cite[стр. 22]{Kassel_knots}.

Структура тензорной категории в $\mathcal{C}$ определяется семействами
изоморфизмов
\begin{equation}\label{ass_tensor}
\alpha_{X,Y,Z}:X\otimes(Y\otimes Z)\;\stackrel{\approx}{\to}\;(X\otimes
Y)\otimes Z,
\end{equation}
\begin{equation}\label{unit_tensor}
\lambda_X:\;I\otimes X\;\stackrel{\approx}{\to}\;X,\qquad\rho_X:\;X\otimes
I\;\stackrel{\approx}{\to}\;X,
\end{equation}
удовлетворяющими естественным требованиям, которые выполняются для
категории модулей над алгеброй Хопфа \cite[стр. 351-357]{Kassel_QG}. В этом
частном случае выделенным элементом $\mathbf{1}$ служит тривиальный
$A$-модуль $\mathbb{C}$, а изоморфизмами $\alpha_{X,Y,Z}$, $\lambda_X$,
$\rho_X$ -- канонические изоморфизмы векторных пространств, что позволяет
использовать знак равенства вместо знака $\approx$, не опасаясь ошибок в
рассуждениях. В общем случае принято заменять рассматриваемую тензорную
категорию эквивалентной ей тензорной категорией, для которой в
\eqref{ass_tensor}, \eqref{unit_tensor} выполняются равенства \cite[стр.
361]{Kassel_QG}, то есть {\bf строгой} тензорной категорией.
\IND{категория! тензорная ! строгая}


\medskip

Алгебра в тензорной категории $\mathcal{C}$ -- \IND{алгебра ! в тензорной
категории} это ее объект $F$ и пара морфизмов
$$m:F\otimes F\to F,\qquad\unit:I\to F,$$
удовлетворяющих следующим требованиям
\cite[стр.438]{Majid_book}:
$$
m(m\otimes\rm{id}_F)=m(\rm{id}_F\otimes m),\qquad m(\unit\otimes
F)=m(F\otimes\unit)=\rm{id}_F.
$$
Левый модуль над алгеброй $F$ в тензорной категории $\mathcal{C}$ --
\IND{модуль ! над алгеброй ! в тензорной категории} это ее объект $E$ и
морфизм $m_{\rm left}:F\otimes E\to E$, для которых
$$
m_{\rm left}(m\otimes\rm{id}_E)\,=\,m_{\rm left}(\rm{id}_F\otimes m_{\rm
left}),\qquad m_{\rm left}(\unit\otimes\rm{id}_E)\,=\,\rm{id}_E.
$$

\begin{example} Рассмотрим тензорную категорию $\mathcal{C}_A$ модулей над
алгеброй Хопфа $A$.
 Каждая унитальная $A$-модульная алгебра $F$ является алгеброй в
 $\mathcal{C}_A$, а каждый $A$-модульный левый $F$-модуль -- левым модулем
 над $F$ в  $\mathcal{C}_A$.
\end{example}

Аналогично вводится понятия правого модуля и бимодуля в тензорной категории.
В случае бимодуля требуется, чтобы левое и правое действия
$$
 m_{\rm left}:\, F\otimes E \to E, \qquad m_{\rm right}: \, E\otimes F \to E
$$
были перестановочны:\ \
 $m_{\rm right}(\rm{id}_F \otimes m_{\rm left})
 \,=\, m_{\rm left}(\rm{id}_F \otimes m_{\rm right})$.

\medskip

Если каждой паре $V',V''$ объектов строгой тензорной категории
$\mathcal{C}$ сопоставлен изоморфизм $\check{R}_{V'V''}:V'\otimes V''\to
V''\otimes V'$, обладающий свойствами, перечисленными в предложениях
\ref{braiding1_sl_2}, \ref{braiding2_sl_2}, то $\mathcal{C}$ называют
сплетенной тензорной категорией, \IND{категория ! тензорная ! сплетенная} а
изоморфизм $\check{R}_{V'V''}$ -- сплетением. Для нестрогих тензорных
категорий в определении участвуют изоморфизмы $\alpha_{X,Y,Z}$, см.
\cite[стр. 395]{Kassel_QG}.

\medskip

Сплетенная тензорная категория называется симметрической, \IND{категория !
тензорная ! сплетенная симметрическая} если
$\check{R}_{V''V'}\check{R}_{V'V''}={\rm id}_{v'\otimes V''}$ \cite[стр.
400]{Kassel_QG}.

\medskip

Алгебра $F$ в сплетенной тензорной категории называется коммутативной,
\IND{алгебра ! коммутативная ! в сплетенной тензорной категории} если
$m\;=\;m\,\check{R}_{FF}$, где $m:F\otimes F\to F$ -- умножение в $F$.
Модули над коммутативными алгебрами в сплетенных абелевых тензорных
категориях используются в алгебраической $K$-теории \cite{LychaginPrasolov}
и естественно возникают в конформной квантовой теории поля
\cite{KirillovOstrik,Fuchs}.\footnote{Абелевы тензорные категории
определяются очевидным образом \cite{Hai_embed}: функтор $\otimes$ должен
быть биаддитивным \cite[стр. 106, 98]{BucurDeleanu}.}

\subsection{Инвариантный интеграл}\label{inv_int_subsection}

\subsubsection{Инвариантные интегралы и $*$-представления.}\label{IntStar}

Рассмотрим \hbox{$A$-модульную} алгебру $F$. Линейный функционал $\nu:F\to
\mathbb{C}$ называют \hbox{$A$-инвариантным} интегралом, \IND{инвариантный
! интеграл} если $\nu$ является морфизмом \hbox{$A$-модулей}, то есть
$$\nu(af)=\varepsilon(a)f,\qquad a\in A,\;f\in F,$$ где
$\varepsilon$-коединица алгебры Хопфа $A$. Другими словами, инвариантный
интеграл -- это $A$-инвариантный элемент $A$-модуля $F^*$.

Рассмотрим модуль $V$ над алгеброй Хопфа $A$, $A$-модульную алгебру
$\operatorname{End}V$ всех линейных операторов в $V$ и ее $A$-модульную
подалгебру
$$\operatorname{End}_fV=\{T\in\operatorname{End}V|\:\dim(TV)<\infty\}.$$
Отметим, что $\operatorname{End}_fV$ -- алгебра без единицы, если $\dim
V=\infty$. Нашей ближайшей целью является построение ненулевого
$A$-инвариантного интеграла на $\operatorname{End}_fV$ в предположении, что
квадрат антипода $S$ алгебры Хопфа $A$ является внутренним автоморфизмом:
 при некотором обратимом $g\in A$
\begin{equation}\label{square_of_S}
S^2(a)=gag^{-1},\qquad a\in A.
\end{equation}
Пусть $\pi(a)$ -- представление алгебры $A$, отвечающее $A$-модулю $V$.
Тогда представление $\pi(S^2(a))$ эквивалентно представлению $\pi(a)$,
поскольку
$$\pi(S^2(a))=\pi(g)\pi(a) \pi(g)^{-1},\qquad a\in A.$$
Значит, линейное отображение
$$V\to V^{**},\qquad v\mapsto\pi(g)v,$$
является морфизмом $A$-модулей, и цепочка морфизмов $A$-модулей
$$
\operatorname{End}_fV\stackrel{\cong}{\to}V\otimes V^* \hookrightarrow
V^{**}\otimes V^*\to \mathbb{C}
$$
доставляет следующий хорошо известный результат.

\begin{proposition} \label{tr_g}
Линейный функционал
$$\nu:\; \operatorname{End}_f(V)\to\mathbb{C},\qquad \nu: \; T\mapsto\operatorname{tr}(T\pi(g))$$
является $A$-инвариантным интегралом.
\end{proposition}

В предположении, что элемент $g$ в \eqref{square_of_S} удовлетворяет
следующим дополнительным условиям
\begin{equation}\label{group_like}
\triangle(g)=g\otimes g,\qquad\varepsilon(g)=1,
\end{equation}
число $\operatorname{tr}_q(A)=\operatorname{tr}(A\pi(g))$ называют
$q$-следом \IND{$q$-след} линейного оператора $A$ в пространстве
представления $\pi$. Из \eqref{group_like} следует, что
$$
\operatorname{tr}_q(A\otimes B)=
\operatorname{tr}_qA\cdot\operatorname{tr}_qB.
$$
В случае конечномерного представления $\pi$ определен $q$-след единичного
оператора, называемый $q$-размерностью пространства представления:
\IND{$q$-размерность пространства представления}
$$\dim_qV\overset{\mathrm{def}}{=}\operatorname{tr}\pi(g).$$

\bigskip

Получим аналог формулы ''интегрирования по частям''.

\begin{proposition}\label{by_parts}
Если $\nu$ -- инвариантный интеграл на $A$-модульной алгебре $F$, то
\begin{equation}\label{int_by_parts}
\nu((af_1)\cdot f_2)=\nu(f_1\cdot S(a)f_2),\qquad f_1,f_2\in F,\;\;a\in A.
\end{equation}
\end{proposition}

{\bf Доказательство.} Покажем, что для любых $A$-модулей $U, V$ и любого
морфизма $A$-модулей $l: V \otimes U \rightarrow \mathbb{C}$ имеет место
равенство
\begin{equation}\label{l_new}
  l(av \otimes u) =  l(v \otimes S(u)),  \qquad   u \in U,\; v \in V,\; a \in A
\end{equation}
После этого останется положить $U=V=F$ и $l(f_1 \otimes
f_2)=\nu(f_1\,f_2)$.

Равенство \eqref{l_new} очевидно в частном случае естественного спаривания
$l_{\rm nat}:U^*\otimes U\rightarrow\mathbb{C}$. Общий случай сводится к
этому частному случаю с помощью морфизма $A$-модулей $V\rightarrow U^*$,
определяемого линейным функционалом $l$:
$$
V\cong V\otimes\mathbb{C}\rightarrow V\otimes U\otimes
U^*\stackrel{l\otimes\operatorname{id}}{\longrightarrow}\mathbb{C}\otimes
U^*\cong U^*.\eqno\square
$$

\bigskip

Опишем один из методов построения $*$-представлений $*$-алгебр Хопфа.
Напомним терминологию. Линейную по первой переменной и антилинейную по
второй переменной форму $(v',v'')$ на векторном пространстве $V$ называют
полуторалинейной формой. \IND{форма ! полуторалинейная} Ее называют
эрмитовой формой, \IND{форма ! эрмитова} если
$(v'',v')=\overline{(v',v'')}$ для всех $v',v''\in V$. Эрмитова форма
называется неотрицательной, \IND{форма ! эрмитова ! неотрицательная} если
$(v,v)\ge 0$ для всех $v\in V$, и положительной, \IND{форма ! эрмитова !
положительная} если $(v,v)\,>\,0$ для всех ненулевых $v\in V$.
Представление $\pi$\ $*$-алгебры $A$ в векторном пространстве $V$,
наделенном эрмитовой формой $(\cdot,\cdot)$, называется $*$-представлением,
если \IND{$*$-представление}
$$(\pi(a)v',v'')=(v',\pi(a^*)v''),\quad a\in A,\quad v',v''\in V.$$
На языке $A$-модулей это означает, что
\begin{equation}\label{star}
(av',v'')=(v',a^*v''),\quad a\in A,\quad v',v''\in V.
\end{equation}

\begin{proposition}\label{Star} Рассмотрим $*$-алгебру Хопфа
  $A$, $A$-модульную \hbox{$*$-алгебру} $F$ и инвариантный интеграл
  $\nu:F\rightarrow \mathbb{C}$,
удовлетворяющий условию вещественности
\begin{equation}\label{real}
  \nu(f^*)=\overline{\nu(f)},\quad f\in F.
\end{equation}
 Представление алгебры
  $A$ в векторном пространстве $F$ с эрмитовой
  формой
  \begin{equation}\label{sesq1}
    (f_1,f_2)=\nu(f_2^*f_1),\quad f_1,f_2\in F
  \end{equation}
  является $*$-представлением.
\end{proposition}

{\bf Доказательство.} Введем понятие $A$-инвариантной полуторалинейной
формы, докажем, что (\ref{star}) равносильно $A$-инвариантности, и, наконец,
установим $A$-инвариантность эрмитовой формы (\ref{sesq1}).

Как известно \cite[стр. 388]{Lang}, каждому векторному пространству $V$
отвечает антипространство \IND{антипространство} $\overline{V}$,
отличающееся от $V$ лишь действием поля $\mathbb{C}$:
$$
\lambda:v\mapsto\bar{\lambda}v,\qquad \lambda\in \mathbb{C},\; \; v\in V.
$$
 Структура модуля над $*$-алгеброй Хопфа $A$ переносится с $V$ на
$\overline{V}$:
$$
a:v\mapsto\,S(a)^*v,\qquad a\in A,\; \; v\in V,
$$
и векторное пространство полуторалинейных форм $V_1\times V_2\to
\mathbb{C}$ естественно изоморфно $(\overline{V_2}\otimes V_1)^*$. Наделим
векторное пространство полуторалинейных форм на $V$ структурой модуля над
алгеброй Хопфа $A$ с помощью этого изоморфизма.

 Рассмотрим  $A$-инвариантную полуторалинейную форму
$(\cdot,\cdot)$ в $V$. Ей отвечает линейный оператор
$$
\overline{V}\hookrightarrow V^*,\quad v\mapsto(\cdot,v),
$$
являющийся морфизмом $A$-модулей. Другими словами,
$$
(v',S(a)^*v'')=(S(a)v',v''),\quad a\in A,\quad v',v''\in V.
$$
Значит, в силу обратимости антипода $S$,
$$
(v',a^*v'')=(av',v''),\quad a\in A,\quad v',v''\in V.
$$
Из инвариантности полуторалинейной формы следует (\ref{star}).

Обратно, предположим, что полуторалинейная форма $(\cdot,\cdot)$ в $V$
удовлетворяет требованию (\ref{star}). Если
$\Delta(a)=\sum\limits_ib_i\otimes c_i$, то, как вытекает из (\ref{star}),
$$\sum_i(c_iv',S(b_i)^*v'')=(v',\sum_ic_i^*S(b_i)^*v'')=$$
$$
=(v',(\sum_iS(b_i)c_i)^*v'')=(v',{\varepsilon(a)}v'')=
\varepsilon(a)(v',v'').
$$
Это означает $A$-инвариантность полуторалинейной формы $(\cdot,\cdot)$.

Остается доказать $A$-инвариантность эрмитовой формы
 (\ref{sesq1}). Пусть, как и прежде,
 $\Delta(a)=\sum\limits_i b_i \otimes c_i$. Воспользуемся тем,  что
 $(*\,S)^2=\operatorname{id}$.
Используя это равенство и (\ref{invol2}), получаем:
$$
\sum_i(c_if_1,S(b_i)^*f_2)=\sum_i\nu(b_i(f_2^*)\cdot c_i(f_1)), \qquad f_1,
f_2 \in F.
$$
Но $F$ является $A$-модульной алгеброй, а $\nu$ -- $A$-инвариантным
интегралом. Значит,
$$
\sum_i\nu(b_i(f_2^*)\cdot
c_i(f_1))=\nu(a(f_2^*f_1))=\varepsilon(a)\nu(f_2^*f_1)=
\varepsilon(a)\cdot(f_1,f_2). \qquad \qquad \hfill\square
$$

\subsubsection{Существует ли ненулевой неотрицательный $U_q\mathfrak{sl}_2$-инвариантный
интеграл?}\label{Exist_IntStar}

Начнем с общеизвестных определений.

Векторные подпространства $\{F^{(j)}A\}_{j\in\mathbb{Z}}$ алгебры $A$
задают фильтрацию этой алгебры, \IND{фильтрация алгебры} если
\begin{itemize}
\item[1.] $F^{(j)}A\subset F^{(k)}A,\qquad j<k$,

\item[2.] $1\in F^{(0)}A$,

\item[3.] $F^{(j)}A\cdot F^{(k)}A \subset F^{(j+k)}A,\qquad j,k\in\mathbb{Z}$.
\end{itemize}
Алгебра с выделенной фильтрацией называется фильтрованной алгеброй.
\IND{алгебра ! фильтрованная}

\begin{example}\label{ideal_filtration}
Каждому двустороннему идеалу $I$ алгебры $A$ отвечает ее фильтрация
$
F^{(j)}A=\begin{cases} A, & j\in \mathbb{Z}_+,\\
                       I^{-j}, & j\in -\mathbb{N}
         \end{cases}.
$
Эта фильтрация называется $I$-адической. \IND{фильтрация алгебры !
$I$-адическая}
\end{example}

В дальнейшем будет предполагаться, что рассматриваемые фильтрации
  обладают следующими дополнительными свойствами:
 $\bigcup\limits_{j\in \mathbb{Z}} F^{(j)}A = A$,
 $\bigcap\limits_{j\in \mathbb{Z}} F^{(j)}A = 0$.
Степень ненулевого элемента $f\in A$ такой алгебры определяется равенством
$\deg(f)=\min\{j\in\mathbb{Z}\,|\,f\in F^{(j)}A\}$. \IND{степень элемента
фильтрованной алгебры}

\begin{example}\label{grad_filtr}
Если $A=\newoplus\limits_{j \in \mathbb{Z}} A_j$ -- градуированная алгебра,
то последовательность
$$
F^{(j)}A=\newoplus\limits_{k \leq j} A_k,\qquad j
\in \mathbb{Z}
$$
 является фильтрацией алгебры $A$.
\end{example}

   \medskip Векторное пространство
   $\operatorname{gr}(A)=
 \newoplus\limits_{j \in \mathbb{Z}_+} (F^{(j)}A/F^{(j-1)}A)$
   наделяется структурой алгебры следующим образом:
   $$ [x+F^{(j-1)}A] \cdot [y+F^{(k-1)}A] = xy + F^{(j+k-1)}A,$$
    где $x \in  F^{(j)}A/F^{(j-1)}A,\; y \in F^{(k)}A/F^{(k-1)}A$.
 Алгебра $\operatorname{gr}(A)$ является  $\mathbb{Z}$-градуированной с
однородными компонентами $\operatorname{gr}(A)_j=F^{(j)}A/F^{(j-1)}A$. Она
называется присоединенной градуированной алгеброй. \IND{алгебра !
присоединенная градуированная} О функторах, связанных с переходом от
фильтрованных алгебр к присоединенным градуированным алгебрам см.
\cite[стр. 278 -- 292]{NastOyst}.

\bigskip

 \begin{example}\label{Pol_filtres} Рассмотрим алгебру
$\operatorname{Pol}(\mathbb{C})_q$ и ее фильтрацию:
$$ \deg(z)=\deg(z^*)=1. $$ Присоединенной градуированной
алгеброй служит некоммутативная алгебра $\mathbb{C}[z,z^*]_{q^2}$ с
одноименными образующими $z, z^*$ и определяющим соотношением $z^*\, z = q^2
z\, z^*$.
\end{example}

Переход от $A$ к $\operatorname{gr}(A)$ является часто используемым
приемом. Например, как объяснялось в \itemiiiе \ref{sL_2_mod_algebras},
множество $\{z^jz^{*k}\,|\, j,k \in \mathbb{Z}_+\}$
 является базисом векторного пространства
 $\mathbb{C}[z,z^*]_{q^2}$. Значит, справедливо следующее утверждение.

\begin{proposition}\label{Wick_basis_l}
Множество $\{z^jz^{*k}\,|\, j,k \in \mathbb{Z}_+\}$
 является базисом векторного пространства
$\operatorname{Pol}(\mathbb{C})_q$.
\end{proposition}

  Аналогично доказывается, что множество $\{z^{*j}z^k\}_{j,k \in
  \mathbb{Z}_+}$ является базисом векторного пространства
  $\operatorname{Pol}(\mathbb{C})_q$.

Переходя к $\operatorname{gr}(\operatorname{Pol}(\mathbb{C})_q)$, нетрудно
показать, что алгебра $\operatorname{Pol}(\mathbb{C})_q$ не имеет делителей
нуля. Действительно, имеет место следующее несложное утверждение \cite[стр.
109]{BrGood}.
\begin{proposition}\label{integrity_domain}
Если присоединенная градуированная алгебра $\operatorname{gr}(A)$ является
областью целостности, то и $A$ -- область целостности.
\end{proposition}

\bigskip

Линейный функционал $\nu$ на $*$-алгебре $A$ называют неотрицательным,
\IND{линейный функционал на $*$-алгебре ! неотрицательный} если
$\nu(f^*f)\ge 0$ при всех $f\in A$. Неотрицательный линейный функционал
$\nu$ называют положительным, \IND{линейный функционал на $*$-алгебре !
положительный} если $\nu(f^*f)\ne 0$ при $f\ne 0$.

 Из
неотрицательности линейного функционала следует его вещественность
\eqref{real}, поскольку из неотрицательности матрицы
 $(\nu(f_i^*f_j))_{i,j=1,2}$ следует ее эрмитовость.

\bigskip
При ответе на вопрос, вынесенный в название \itemiiiа, будет использован
  элемент \eqref{C_q_second_form}  центра алгебры $U\mathfrak{sl}_2$ и равенство
\begin{equation}\label{zero_central_char}
C_q f=0,\qquad f\in \mathbb{C}[z]_q
\end{equation}
(его доказательство легко сводится к очевидным частным случаям $f=1$ и
$f=z$).

\begin{proposition} На ${\rm Pol}(\mathbb{C})_q$ не существует ненулевого
 неотрицательного  $U_q\mathfrak{sl}_2$-инвариантного интеграла.
\end{proposition}

{\bf Доказательство.} Пусть $\nu$ -- ненулевой вещественный
$U_q\mathfrak{sl}_2$-инва\-риантный интеграл на ${\rm Pol}(\mathbb{C})_q$.

Множество $\{z^{*m}z^n\}_{m,n \in \mathbb{Z}_+}$ является базисом векторного
пространства ${\rm Pol}(\mathbb{C})_q$. Из
 $\nu(K^{\pm 1}f)=\nu(f)$ следует, что  $\nu(z^{*m}z^n)=0$ при $m\neq n$,
 а из
 $\nu(Ff)=0$, $Fz=q^{\frac{1}{2}}$ вытекает равенство $\nu(1)=0$.

 Значит, $\nu(z^{*n}z^n)\neq 0$ при некотором
 $n \in \mathbb{N}$.  Но $z^{n+1}$ лишь числовым множителем
 отличается от $E^n z$.
 Следовательно, $\nu(z^*z)\neq 0$. Из неотрицательности $\nu$ вытекает
неравенство $\nu(z^*z)> 0$.

 Чтобы прийти к противоречию достаточно доказать наличие отрицательных чисел
 в последовательности $\{\nu((E^n z)^*(E^n z))\}_{n \in \mathbb{N}}$.
Как следует из $\eqref{invol2}$, $\eqref{int_by_parts}$, это равносильно
 наличию отрицательных чисел в последовательности
 $\{\nu(z^* (F^nE^n z))\}_{n \in \mathbb{N}}$.

Но знаки этих чисел чередуются, поскольку
$$
F^nE^nz= -F^{n-1}\frac{1}{(q^{-1}-q)^2} \left(q^{-1} K^{-1} + qK -
(q^{-1}+q)\right)E^{n-1}z.
$$
Последнее равенство вытекает из \eqref{C_q},
$$
F^nE^nz= F^{n-1}C_qE^{n-1}z-F^{n-1}\frac{1}{(q^{-1}-q)^2} \left(q^{-1}
K^{-1} + qK - (q^{-1}+q)\right)E^{n-1}z$$
 и из \eqref{zero_central_char}.\hfill $\square$

\bigskip
Полученный результат можно было ожидать. Действительно, $*$-алгебра
 ${\rm Pol}(\mathbb{C})_q$ тесно связана с квантовым аналогом единичного круга
$\mathbb{D}=\{z\in \mathbb{C}\,|\, |z|<1\}$, как показано в \itemiiе
\ref{topology}. В классическом случае $q=1$ инвариантный интеграл в
$\mathbb{D}$ известен \cite[стр. 46]{Helg1}
\begin{equation}\label{integral_l}
\nu(f)=\int\limits_\mathbb{D}f\cdot (1-|z|^2)^{-2} \; d\operatorname{Re}z
\;d\operatorname{Im}z,
\end{equation}
и $\nu(f^*f)=+\infty$ для любого ненулевого полинома $f$.

Чтобы избежать расходимости интегралов, достаточно заменить алгебру
полиномов алгеброй $\mathscr{D}(\mathbb{D})$ гладких функций с компактными
носителями в $\mathbb{D}$. Используемое в задачах гармонического анализа в
$\mathbb{D}$ гильбертово пространство $L^2(d\nu)$ не выпадает из поля
зрения, поскольку является пополнением $\mathscr{D}(\mathbb{D})$ по норме
$\|f\|=\nu(f^*f)^{\frac{1}{2}}$.

  В следующем \itemiiiе  вводится $U_q\mathfrak{su}_{1,1}$-модульная
 $*$-алгебра $\mathscr{D}(\mathbb{D})_q$, являющаяся $q$-аналогом
 алгебры $\mathscr{D}(\mathbb{D})$. Как будет показано, положительный
 $U\mathfrak{sl}_2$-инвариантный интеграл на $\mathscr{D}(\mathbb{D})_q$
 существует и единствен с точностью до числового множителя.

\subsubsection{Обобщенные функции и финитные функции.}\label{finite_sl_2}

Введем пространство $\mathscr{D}(\mathbb{D})'_q$, являющееся $q$-аналогом
пространства $\mathscr{D}(\mathbb{D})'$ обобщенных функций в квантовом
круге. С этой целью наделим векторное пространство ${\rm Pol}(\mathbb{C})_q$
топологией и пополним его.

Рассмотрим введенное в \itemiiiе \ref{Fock_l} точное $*$-представление
$T_F$ $*$-алгебры ${\rm Pol}(\mathbb{C})_q$ в предгильбертовом пространстве
$\mathcal{H}$ с выделенным ортонормированным базисом $\{e_n\}_{n\in
\mathbb{Z}_+}$. Наделим ${\rm Pol}(\mathbb{C})_q$ слабейшей из топологий
$\mathcal{T}$, в которых непрерывны все линейные функционалы
\begin{equation}\label{first-def-topology}
l'_{j,k}(f)\,\stackrel{\rm def}{=}\, (T_F(f)e_j,e_k), \qquad j,k \in
\mathbb{Z}_+,
\end{equation}
то есть все матричные элементы фоковского представления в выделенном базисе.

Имеется другое описание топологии $\mathcal{T}$. Как было показано в
\itemiiiе \ref{Fock_l}, каждый элемент $f\in {\rm Pol}(\mathbb{C})_q$
является конечной суммой вида \eqref{expansion_l}
\begin{equation}\label{expansion_ll}
 f=\sum_{m \in \mathbb{N}}z^m \psi_m(y)\,
+\, \psi_0(y)\, +\,
 \sum_{m \in \mathbb{N}} \psi_{-m}(y)\, z^{*m},
\end{equation}
и такое разложение единственно. Здесь $y=1-zz^*$ и $\psi_m$ -- полиномы
одной переменной.

\begin{proposition}\label{second-def-topology}
$\mathcal{T}$ -- слабейшая из топологий, в которых непрерывны все линейные
функционалы
\begin{equation}\label{l_two_prime}
l''_{m,n}(f)\,\stackrel{\rm def}{=}\,\psi_m(q^{2n}),\qquad
m\in\mathbb{Z},\;n\in\mathbb{Z}_+,
\end{equation}
то есть значения коэффициентов $\psi_m(y)$ разложения \eqref{expansion_ll} в
точках геометрической прогрессии $q^{2\mathbb{Z}_+}$.
\end{proposition}

{\bf Доказательство.} Пусть $v_m=T_F(z^m)e_0$ и $l_{j,k}(f)$ -- матричные
элементы операторов фоковского представления в базисе
$\{v_m\}_{m\in\mathbb{Z}_+}$. Из определений следует, что $\mathcal{T}$
является слабейшей из топологий, в которых непрерывны все линейные
функционалы $l_{j,k}$. Здесь $j,k\in\mathbb{Z}_+$. Остается заметить, что
$l_{j,k}$ лишь ненулевым числовым множителем отличается от линейного
функционала $l''_{k-j,\,{\rm min}(j,k)}$. \hfill $\square$

\bigskip

Наделим пространство всех функций на множестве $q^{2\mathbb{Z}_+}$
топологией поточечной сходимости. Оно является полным в этой топологии
\cite[стр. 50]{Edwards}, и сужения полиномов образуют плотное линейное
подмногообразие. Значит, пространство всех функций на $q^{2\mathbb{Z}_+}$,
наделенное топологией поточечной сходимости, является пополнением
пространства полиномов с индуцированной топологией \cite[стр. 51]{Edwards}.

Рассуждая аналогично, рассмотрим векторное пространство
$\mathscr{D}(\mathbb{D})'_q$ {\bf формальных рядов} \IND{пространство !
$\mathscr{D}(\mathbb{D})'_q$} вида $\eqref{expansion_ll}$ с коэффициентами
из пространства функций на множестве $q^{2\mathbb{Z}_+}$. Наделим
$\mathscr{D}(\mathbb{D})'_q$ топологией поточечной сходимости коэффициентов
$\psi_j$ на множестве $q^{2\mathbb{Z}_+}$. Очевидно,
$\mathscr{D}(\mathbb{D})'_q$ является полным топологическим векторным
пространством, и ${\rm Pol}(\mathbb{C})_q$ естественно изоморфно его
плотному линейному подмногообразию. Следовательно,
$\mathscr{D}(\mathbb{D})'_q$ -- пополнение топологического векторного
пространства ${\rm Pol}(\mathbb{C})_q$.

Инволюция $*$ допускает продолжение по непрерывности на
$\mathscr{D}(\mathbb{D})'_q$.

\medskip

Пусть $\overline{\mathcal{H}}$ -- прямое произведение одномерных векторных
пространств $\mathcal{H}_j=\mathbb{C}e_j$ с топологией прямого
произведения. Очевидно, $\mathcal{H}\subset\overline{\mathcal{H}}$ --
плотное линейное подмногообразие.

Линейные операторы из $\mathcal{H}$ в $\overline{\mathcal{H}}$ имеют вид
\begin{equation}\label{formal_operators}
Av=\sum\limits_{k,j=0}^\infty a_{kj}\,(v,e_j)\,e_k,\qquad v\in\mathcal{H},
\end{equation}
где $a_{kj}\in\mathbb{C}$. Будем называть их обобщенными линейными
операторами в $\mathcal{H}$ \IND{обобщенные ! линейные операторы} и
наделять векторное пространство таких операторов топологией
покоэффициентной сходимости формальных рядов \eqref{formal_operators}.
Приведем следствие предложения \ref{second-def-topology}.

\begin{corollary}\label{isom_general_sl_2}
  1. Фоковское представление алгебры $\mathrm{Pol}(\mathbb{C})_q$
допускает продолжение по непрерывности до линейного отображения
$\mathscr{D}(\mathbb{D})_q'$ в пространство обобщенных линейных операторов.
\\ 2. Это линейное отображение является изоморфизмом топологических
векторных пространств.
\end{corollary}

\medskip
Сохраним обозначение $T_F$ для полученного продолжения на векторное
пространство $\mathscr{D}(\mathbb{D})_q'$.

  Напомним, что ${\rm Pol}(\mathbb{C})_q$ является
$U_q\mathfrak{su}_{1,1}$-модульной алгеброй и, следовательно,
$U_q\mathfrak{sl}_2$-модульным ${\rm Pol}(\mathbb{C})_q$-бимодулем. Можно
показать, что операторы представлений алгебр ${\rm Pol}(\mathbb{C})_q$ и
$U_q\mathfrak{sl}_2$ продолжаются по непрерывности с ${\rm
Pol}(\mathbb{C})_q$ на пополнение $\mathscr{D}(\mathbb{D})'_q$, наделяя его
структурой $U_q\mathfrak{sl}_2$-модульного ${\rm
Pol}(\mathbb{C})_q$-бимодуля. Доказательство будет приведено позднее в
существенно большей общности, см. \itemii\ \ref{function_spaces}.

\begin{example}
 Равенства $K^{\pm 1}\psi(y)=\psi(y)$,
\begin{equation}\label{EF-psi}
E\psi(y)=-q^{\frac{1}{2}}zy\frac{\psi(y)-\psi(q^2y)}{y\,-\,q^2y},\quad
F\psi(y)=-q^{\frac{5}{2}}y\frac{\psi(y)-\psi(q^2y)}{y\,-\,q^2y}z^*
\end{equation}
нетрудно доказать по индукции для $\psi=y^n$. Значит, они верны для всех
полиномов $\psi(y)$, а, следовательно, и для всех функций $\psi(y)$ на
 $q^{2\mathbb{Z}_+}$.
\end{example}

\medskip

Линейный оператор $A$ в $\mathcal{H}$ назовем финитным, \IND{финитный
линейный оператор} если $Ae_j=0$ для всех $j\in\mathbb{Z}_+$, кроме
конечного их числа. Примером может служить одномерный ортогональный
проектор $P_0$ на вакуумное подпространство $\mathcal{H}_0=\mathbb{C}e_0$.
Стандартным образом вводится понятие сопряженного к финитному линейному
оператору.

С помощью предложения \ref{second-def-topology} легко доказывается

\begin{proposition}\label{def-finite} Следующие утверждения
для $f\in \mathscr{D}(\mathbb{D})'_q$ эквивалентны:
\begin{itemize}
 \item[1.]\ $T_F(f)$ является финитным линейным оператором в $\mathcal{H}$;
 \item[2.]\  конечны как число ненулевых
членов формального ряда \eqref{expansion_ll}, так и носители функций
$\psi_j(y)$ в \eqref{expansion_ll};
 \item[3.]\ $z^{*N}\,f\,=\,f\,z^N\,=\,0$ при некотором $N \in \mathbb{Z}_+$.
\end{itemize}
\end{proposition}

\medskip

Элементы $f\in\mathscr{D}(\mathbb{D})'_q$, удовлетворяющие условиям
предложения \ref{def-finite}, будем называть финитными функциями в
квантовом круге. \IND{функция ! финитная ! в квантовом круге} Они образуют
плотное в $\mathscr{D}(\mathbb{D})'_q$ линейного подмногообразие,
обозначаемое $\mathscr{D}(\mathbb{D})_q$.

Как следует из определений, отображение $T_F$ является изоморфизмом
векторного пространства $\mathscr{D}(\mathbb{D})_q$ и векторного
пространства финитных линейных операторов в $\mathcal{H}$. Это позволяет
наделить $\mathscr{D}(\mathbb{D})_q$ структурой \hbox{$*$-алгебры}.
Умножение в $\mathscr{D}(\mathbb{D})_q$ можно получить и другим способом --
по непрерывности, исходя из умножения в ${\rm Pol}(\mathbb{C})_q$. Так
можно
 доказать, что $\mathscr{D}(\mathbb{D})_q$ является
$U_q\mathfrak{su}_{1,1}$-модульной алгеброй. Более общий подход изложен в
\itemiiе \ref{function_spaces}.

\begin{example}\label{f_0_sl_2}
 Существует и единствен элемент $f_0\in \mathscr{D}(\mathbb{D})_q$, для
которого $T_F f_0=P_0$. Очевидно,
\begin{equation}\label{rel3_sl_2}
 z^*\,f_0\,=\,f_0\,z\,=\,0,
\end{equation}
\begin{equation}\label{rel1_sl_2} f_0^2=f_0,\qquad
 f_0^*=f_0.
\end{equation}
Коэффициенты разложения \eqref{expansion_ll} элемента $f_0$ равны
$$
\psi_0(y)=\begin{cases}1, &\; \text{при}\; y=1,\\
0, &\; \text{при}\; y\neq 1
\end{cases}, \qquad \psi_j(y)=0, \; \; \; \text{при}\; j \neq 0.
$$
Значит,
\begin{equation}\label{lemma_EFf0}
K^{\pm 1} f_0 = f_0,\qquad F f_0 = -\frac{q^{1/2}}{q^{-2}-1} f_0 z^*,\qquad
E f_0 = -\frac{q^{1/2}}{1-q^2} z f_0.
\end{equation}
\end{example}
Как следует из определений фоковского представления $T_F$ и элемента $f_0$,
\begin{equation}\label{f_0-generator_of_ideal}
\mathscr{D}(\mathbb{D})_q=\mathbb{C}[z]_q\cdot f_0\cdot
\mathbb{C}[\bar{z}]_q.
\end{equation}

\bigskip

\begin{remark}
Соображения, которые привели нас к $U_q\mathfrak{su}_{1,1}$-мо\-дуль\-ной
алгебре $\mathscr{D}(\mathbb{D})_q$, наглядны и достаточно хорошо
мотивированы. В случае квантовой ограниченной симметрической области общего
вида алгебра финитных функций будет построена иначе, более формально, и
элемент $f_0$ сыграет ключевую роль в построении. Вскоре после этого станет
очевидной эквивалентность обоих подходов.
\end{remark}

\subsubsection{Инвариантный интеграл на
$\mathscr{D}(\mathbb{D})_q$.}\label{inv_int_sl_2}

Начнем со следующего вспомогательного утверждения.

\begin{lemma}\label{for_inv_int}
Рассмотрим $A$-модульную алгебру $F$, $A$-модульный левый $F$-модуль $V$ и
отвечающее ему представление $\pi$ алгебры $F$ в векторном пространстве $V$
$$\pi(f)v=fv,\qquad f\in F,\; v \in V.$$
Гомоморфизм алгебр $\pi:F \to {\rm End}(V)$ является также морфизмом
$A$-модулей.
\end{lemma}

{\bf Доказательство.} Пусть $a \in A$,
$$
({\rm id}\otimes \Delta)\Delta (a)\,=\,(\Delta \otimes {\rm id})\Delta (a)=
\sum\limits_j\; b_j\otimes c_j \otimes d_j.
$$
Требуется доказать равенство $\pi(af)=a\pi(f)$ для всех $f \in F$, то есть
равенство
$$ (af)v\,=\, \sum\limits_j \, (b_jf)(c_jS(d_j)v),\qquad
f\in F,\, v \in V.
$$
Но $\sum\limits_j b_j \otimes c_jS(d_j)\,=\, a \otimes 1$ по определению
антипода $S$. \hfill $\square$

\bigskip
Из результатов предыдущего \itemiiiа следует, что
$\widetilde{\mathcal{H}}=\mathbb{C}[z]_q\,f_0$ является как
$U_q\mathfrak{b}^+$-модульным $\mathscr{D}(\mathbb{D})_q$-модулем, так и
$U_q\mathfrak{b}^+$-модульным ${\rm Pol}(\mathbb{C})_q$-модулем. Очевидно,
$$
\widetilde{\mathcal{H}}=\{f \in \mathscr{D}(\mathbb{D})_q|\;fz=0 \}= \{f \in
\mathscr{D}(\mathbb{D})_q|\; T_F(f)e_j=0 \; \text{при всех}\; j \neq 0\}.
$$

Линейное отображение
\begin{equation}\label{g_r_Fock} \widetilde{\mathcal{H}} \to
\mathcal{H},\qquad \psi(z)\,f_0 \mapsto T_F(\psi(z))e_0
\end{equation}
является изоморфизмом ${\rm Pol}(\mathbb{C})_q$-модулей.

С помощью \eqref{g_r_Fock} перенесем структуру $U_q\mathfrak{b}^+$-модуля
из $\widetilde{\mathcal{H}}$ в $\mathcal{H}$ и введем обозначение $\Gamma$
для полученного представления алгебры $U_q\mathfrak{b}^+$ в
предгильбертовом пространстве $\mathcal{H}$.

\begin{proposition}\label{explicit_int_sl_2}
 Линейный функционал
\begin{equation}\label{int_D_sl_2}
\int\limits_{\mathbb{D}_q} f d\nu\stackrel{\operatorname{def}}{=}(1-q^2)\,
\mathrm{tr}_q T_F(f)\,=\,(1-q^2)\,\mathrm{tr}(T_F(f)\Gamma(K^{-1}))
\end{equation}
на $\mathscr{D}(\mathbb{D})_q$ является положительным
$U_q\mathfrak{sl}_2$-инвариантным интегралом.
\end{proposition}

{\bf Доказательство.} При $f\in\mathscr{D}(\mathbb{D})_q$ операторы
$T_F(f)$ в $\mathcal{H}$ являются финитными. Значит, \eqref{int_D_sl_2}
определяет неотрицательный линейный функционал $\nu$ на
$\mathscr{D}(\mathbb{D})_q$: $\mathrm{tr}_qT_F(f^*f)=\mathrm{tr}(T_F(f)
\Gamma(K_{-2\rho})T_F(f)^*)\ge 0$. Его положительность вытекает из точности
представления $T_F$ алгебры $\mathscr{D}(\mathbb{D})_q$, см.
\itemiii\ \ref{finite_sl_2}.

Гомоморфизм алгебр $T_F:\mathscr{D}(\mathbb{D})_q\to{\rm End}(H)_f$
является морфизмом $U_q{\mathfrak b}^+$-модулей, как следует из леммы
\ref{for_inv_int}. Значит, $U_q\mathfrak{b}^+$-инвариантность линейного
функционала $\nu$ вытекает из предложения \ref{tr_g} и равенства
$$S^2(a)=K^{-1}aK,\qquad a\in U_q\mathfrak{sl}_2.$$
Из $U_q\mathfrak{b}^+$-инвариантности и равенства \eqref{invol2} следует
$U_q\mathfrak{b}^-$-ин\-ва\-ри\-ант\-ность, поскольку линейный функционал
$\nu$ является вещественным:
$$\nu(f^*)=\overline{\nu(f)},\qquad f\in\mathscr{D}(\mathbb{D})_q.$$
Остается заметить, что подалгебры Хопфа $U_q\mathfrak{b}^{\pm}$ порождают
$U_q\mathfrak{sl}_2$. \hfill $\square$

\medskip

\begin{remark} Напомним обозначение Джексона
$$ \int_0^1 f(t)d_{q^2}\,t\;=
\; (1-q^2)\,\sum\limits_{j=0}^\infty f(q^{2j})q^{2j}.
$$
Если $
 f=\sum\limits_{m \in \mathbb{N}}z^m \psi_m(y)\,
+\, \psi_0(y)\, +\,
 \sum\limits_{m \in \mathbb{N}} \psi_{-m}(y)\, z^{*m}\in
\mathscr{D}(\mathbb{D})_q $, то
\begin{equation}\label{inv_Jackson}
\int\limits_{\mathbb{D}_q} f d\nu \;=\; \int_0^1
\psi_0(y)\,y^{-2}\,d_{q^2}\,y.
\end{equation}
В формальном пределе $q \to 1$ получаем: $y=1-|z|^2$,
$$d\nu\,=\,\frac{d{\rm Re}\,z\; d{\rm Im}\,z}{\pi\, (1-|z|^2)^2}.$$
\end{remark}

\begin{proposition}\label{t2.3.9}
$U_q \mathfrak{sl}_2\cdot f_0=\mathscr{D}(\mathbb{D})_q$.
\end{proposition}

{\bf Доказательство.} $*$-Алгебра $\mathscr{D}(\mathbb{D})_q$ является
$U_q\mathfrak{su}_{1,1}$-модульной, и имеют место равенства \eqref{q-C},
\eqref{lemma_EFf0}, \eqref{U_q_zstar_sl_2}. Значит,
\begin{equation}\label{2.3.6}
 E(z^jf_0)\,=\,a_j\,z^{j+1}f_0,
\end{equation}
\begin{equation}\label{2.3.7}
\left \{
\begin{gathered}
F(f_0z^{*k})\,=\,b_k\,f_0z^{*(k+1)}
\\ F(z^{j+1}f_0z^{*k})\,=\,c_{jk}\,z^jf_0z^{*k}\;+
\;d_{jk}\,z^{j+1}f_0z^{*(k+1)}
\end{gathered}
\right.
\end{equation}
для некоторых ненулевых числовых множителей
  $a_j$,\  $b_k$,\  $c_{jk}$,\  $d_{jk}$.

Рассмотрим подпространства
$$
L_m=\left\{\left.f\,=\,\sum_{j=0}^m\sum_{i=0}^\infty\,a_{ij}\,z^if_0z^{*j},
\;\right|\;a_{ij}\in\mathbb{C}\right\}\subset\mathscr{D}(\mathbb{D})_q,
\qquad m\in\mathbb{Z}_+.
$$

Как следует из \eqref{f_0-generator_of_ideal}, достаточно доказать, что
 $U_q \mathfrak{sl}_2 \cdot f_0 \supset L_m$ при всех
 $m \in \mathbb{Z}_+$. Воспользуемся  индукцией по $m$. Во-первых,
из \eqref{2.3.6} следует, что $L_0=U_q \mathfrak{b}_+\cdot f_0$, поскольку
$a_j \ne 0$ при всех $j \in{\mathbb Z}_+$.
 Во-вторых, из \eqref{2.3.7} следует, что $L_{m+1}\subset
L_m\,+\, F\cdot L_m$, поскольку $b_k \ne 0$, $c_{jk}\ne 0$, $d_{jk}\ne 0$
при всех $j,k \in{\mathbb Z}_+$. \hfill $\square$

\begin{corollary}\label{uniq_int_sl_2}
 $U_q\mathfrak{sl}_2$-инвариантный интеграл
на $\mathscr{D}(\mathbb{D})_q$ единствен с точностью до числового множителя.
\end{corollary}

Мы доказали существование и единственность
$U_q\mathfrak{sl}_2$-инвариантного интеграла и нашли его явный вид
\eqref{int_D_sl_2}.

\bigskip

Предложение \ref{t2.3.9} можно уточнить.
\begin{proposition}\label{def-rel-D(D)}
Элемент $f_0$ является образующей $U_q\mathfrak{sl}_2$-модуля
$\mathscr{D}(\mathbb{D})_q$, а соотношение $K^{\pm 1}f_0=f_0$ -- его
определяющим соотношением.
\end{proposition}

{\bf Доказательство.} Рассмотрим $U_q\mathfrak{sl}_2$-модуль $V_0$ с
образующей $v_0$ и определяющим соотношением $K^{\pm 1}v_0=v_0$. Требуется
доказать биективность морфизма $U_q\mathfrak{sl}_2$-модулей $j:V_0\to
\mathscr{D}(\mathbb{D})_q$, отображающего $v_0$ в $f_0$.

Рассмотрим подпространства $V_0^{(i,k)}\subset V_0$,
$\mathscr{D}(\mathbb{D})_q^{(i,k)}\subset\mathscr{D}(\mathbb{D})_q$,
являющиеся линейными оболочками подмножеств
$$
\left\{\left.F^{k_1}E^{i_1}\,v_0\;\right|\; k_1\le k\;\&\;i_1\le i\right\},
\qquad\left\{\left.z^{i_1}\,f_0\,z^{*k_1}\;\right|\;k_1\le k\;\&\;i_1\le i
\right\}.
$$
Рассуждая так же, как при доказательстве предложения \ref{t2.3.9}, можно
получить включение $j\,V_0^{(i,k)}\supset\mathscr{D}(\mathbb{D})_q^{(i,k)}$
для всех $i,k\in\mathbb{Z}_+$. Действительно,
$$
\mathscr{D}(\mathbb{D})_q^{(i,k)}\subset\mathscr{D}(\mathbb{D})_q^{(i-1,k-1)}
+F\mathscr{D}(\mathbb{D})_q^{(i,k-1)},\qquad i,k\ne 0,
$$
$$
\mathscr{D}(\mathbb{D})_q^{(0,j)}\subset
F\mathscr{D}(\mathbb{D})_q^{(0,j-1)},\qquad
\mathscr{D}(\mathbb{D})_q^{(j,0)}\subset
E\mathscr{D}(\mathbb{D})_q^{(j-1,0)},\qquad j\ne 0,
$$
см. \eqref{2.3.6},\eqref{2.3.7}.

Но $ \dim V_0^{(i,k)}\leq\dim
\left(\mathscr{D}(\mathbb{D})_q^{(i,k)}\right)<\infty$. Значит, $j$
биективно отображает $V_0^{(i,k)}$ на $\mathscr{D}(\mathbb{D})_q^{(i,k)}$.
\hfill $\square$

\subsection{Дифференциальные исчисления}

\subsubsection{Дифференциальные формы с полиномиальными
коэффициентами.}\label{pol-forms_sl_2}

Введем необходимые определения, следуя \cite{KlSch}. Алгебры будут
предполагаться унитальными, если не оговорено противное.

Рассмотрим градуированную алгебру
$\Omega=\newoplus\limits_{i\in\mathbb{Z}_+}\Omega_i$ и линейный оператор
\hbox{$d:\Omega\to\Omega$} степени 1 в $\Omega$. Пару $(\Omega,d)$ называют
дифференциальной градуированной алгеброй, \IND{дифференциальная
градуированная алгебра} если $d^2=0$ и
\begin{equation}\label{Leibnitz_form}
d(\omega'\cdot\omega'')=d\omega'\cdot\omega''+(-1)^n\omega'\cdot
d\omega'',\qquad\omega'\in\Omega_n,\;\omega''\in\Omega.
\end{equation}
Гомоморфизм дифференциальных градуированных алгебр определяется очевидным
образом. Для умножения в $\Omega$ наряду с $\cdot$ используют символ
$\wedge$. Очевидно, $d1=0$.

Дифференциальное исчисление над алгеброй \IND{дифференциальное исчисление !
над алгеброй} $F$ -- это такая дифференциальная градуированная алгебра
$(\Omega,d)$, что $\Omega_0=F$ и $\Omega$ порождается элементами из
$F\oplus dF$. Рассмотрим алгебру Хопфа $A$ и $A$-модульную алгебру $F$.
Дифференциальное исчисление $(\Omega,d)$ называется ковариантным,
\IND{дифференциальное исчисление ! ковариантное} если алгебра $\Omega$
является $A$-модульной и $d$ -- эндоморфизм $A$-модуля $\Omega$.

\bigskip

Пусть $A=U_q\mathfrak{sl}_2$ и $(\Lambda(\mathbb{C})_q,d)$ -- ковариантное
дифференциальное исчисление над $\mathbb{C}[z]_q$. Из \eqref{q-C} следуют
равенства
\begin{equation}\label{KF-dz}
Fdz=0,\,\qquad K^{\pm 1}dz\,=\,q^{\pm 2}\,dz,
\end{equation}
\begin{equation}\label{E-dz}
Edz\,=\,-q^{1/2}(dz\,z\,+\,z\,dz).
\end{equation}

Напомним, что $\mathbb{C}[z]_q$ -- коммутативная алгебра в сплетенной
тензорной категории $\mathcal{C}^-$. Предположим, что
$\Lambda(\mathbb{C})_q=\oplus_{i\in\mathbb{Z}_+}\Lambda^i(\mathbb{C})_q$
является суперкоммутативной алгеброй в этой категории, то есть
\begin{equation}\label{supercommutativity_sl_2}
m\,\check{R}_{\Lambda(\mathbb{C})_q,\Lambda(\mathbb{C})_q}
(\omega_1\otimes\omega_2)\,=\,\,(-1)^{mn}\;\omega_1\cdot\omega_2
\end{equation}
для любых $\omega_1\in\Lambda^m(\mathbb{C})_q$,
$\omega_2\in\Lambda^n(\mathbb{C})_q$, ср. с \eqref{commutativity_sl_2}. Из
равенств \eqref{Rmatrix_sl_2}, \eqref{KF-dz} следует коммутационное
соотношение
\begin{equation}\label{z_dz_sl_2}
z\,dz\,=\,q^{-2}\,dz\,z.
\end{equation}
Применяя дифференциал $d$, получаем:
\begin{equation}
\label{dz_dz_sl_2}dz\wedge dz=0.
\end{equation}

Мы видим, что ''естественное'' ковариантное дифференциальное исчисление над
алгеброй $\mathbb{C}[z]_q$ единственно: градуированная алгебра
$\Lambda(\mathbb{C})_q$ порождена элементами $z$, $dz$ степеней
$$\deg(z)=0,\qquad\deg(dz)\,=\,1$$
и определяется соотношениями \eqref{z_dz_sl_2} -- \eqref{dz_dz_sl_2}, а
структура $U_q\mathfrak{sl}_2$-модульной алгебры описывается равенствами
\eqref{KF-dz} -- \eqref{E-dz}.

Разумеется, нужно доказать существование определяемого равенствами
\eqref{z_dz_sl_2} -- \eqref{dz_dz_sl_2}, \eqref{KF-dz} -- \eqref{E-dz}
ковариантного дифференциального исчисления. Доказательство будет приведено
в \itemiiiах \ref{FODC}, \ref{diff_univ} в существенно большей общности.
Отметим, что дифференциальное исчисление на квантовой комплексной плоскости
было найдено в \cite[стр. 1158]{SchSch}.

\bigskip

Перейдем к $*$-алгебрам. Дифференциальное исчисление над $*$-алгеброй $F$
называют $*$-исчислением \IND{$*$-дифференциальное исчисление} \cite[стр.
462]{KlSch}, если задано антилинейное инволютивное отображение
$*:\Omega\to\Omega$, совпадающее с инволюцией $*$ на $\Omega_0=F$ и такое,
что при всех $\omega_1\in\Omega_m$, $\omega_2\in\Omega_n$,
$\omega\in\Omega$
\begin{equation}\label{def-calc}
(\omega_1\wedge\omega_2)^*\,=\,\,(-1)^{mn}\;\omega_2^*\wedge\omega_1^*,
\qquad(d\omega)^*\,=d\omega^*.
\end{equation}


\bigskip

Рассмотрим частный случай $A=U_q\mathfrak{su}_{1,1}$, $F={\rm
Pol}(\mathbb{C})_q$. Из \eqref{z_dz_sl_2} -- \eqref{dz_dz_sl_2} вытекают
соотношения
\begin{equation}\label{z*_dz*_sl_2}
z^*\,dz^*\,=\,q^2\,dz^*\,z^*,\qquad dz^*\,\wedge\,dz^*=0.
\end{equation}
Предположим, что ${\rm Pol}(\mathbb{C})_q$-подбимодули, порожденные
элементами $dz$, $dz^*$, имеют нулевое пересечение. Дифференцируя
\eqref{Pol_l}, получаем:
\begin{equation}\label{diff-forms_sl_2}
z^*\,dz\,=\,q^2 dz\,z^*,\qquad dz^*\,z\,=\,q^2\, z\,dz^*,
\end{equation}
\begin{equation}\label{dz-dz^*_sl_2}
dz^*\wedge dz\,=\,-q^2\,dz\wedge dz^*.
\end{equation}
Из \eqref{U_q_zstar_sl_2} следуют равенства
\begin{equation}\label{U_q_dzstar_sl_2}
K^{\pm 1}dz^*=q^{\mp 2}dz^*,\quad Ez^*=0,\quad
Fz^*=-q^{5/2}(dz^*\,z^*+z^*\,dz^*).
\end{equation}



\bigskip

Мы получили полный список соотношений и можем приступить к построению
дифференциальных исчислений над
 ${\rm Pol}(\mathbb{C})_q$. Необходимые проверки корректности определений отложим
до \itemiiiа \ref{pol_diff}, где они будут проведены в существенно большей
общности.

Рассмотрим биградуированную алгебру
$$
\Omega(\mathbb{C})_q=\newoplus\limits_{i,j\in\{0,1\}}\,
\Omega^{i,j}(\mathbb{C})_q
$$
с образующими $z$, $z^*$, $dz$, $dz^*$ бистепеней
$$\deg(z)=\deg(z^*)=(0,0),\qquad\deg(dz)=(1,0),\qquad\deg(dz^*)=(0,1)$$
и определяющими соотношениями \eqref{Pol_l}, \eqref{z_dz_sl_2},
\eqref{diff-forms_sl_2}, \eqref{z*_dz*_sl_2}, \eqref{dz_dz_sl_2},
\eqref{dz-dz^*_sl_2}:
\begin{equation*}
z^*\,z\,=\,q^2z\,z^*\,+\,1-q^2,
\end{equation*}
\begin{equation*}
z\,dz\,=\,q^{-2}\,dz\,z,\quad z^*\,dz\,=\,q^2 dz\,z^*,\quad
dz^*\,z\,=\,q^2\,z\,dz^*,\quad dz^*\,z^*\,=\,q^{-2}\,z^*dz^*,
\end{equation*}
\begin{equation*}
dz\wedge dz=0,\quad dz^*\,\wedge\,dz^*=0,\quad dz^*\wedge
dz\,=\,-q^2\,dz\wedge dz^*.
\end{equation*}
Наделим ее градуировкой
$\Omega^k(\mathbb{C})_q=\oplus_{i+j=k}\,\Omega^{i,j}(\mathbb{C})_q$.

Каждый из ${\rm Pol}(\mathbb{C})_q$-бимодулей $\Omega^{i,j}(\mathbb{C})_q$
является как свободным левым, так и свободным правым ${\rm
Pol}(\mathbb{C})_q$-модулем ранга 1. В частности, для любого элемента
$\omega\in \Omega(\mathbb{C})_q$ существует и единственно разложение
\begin{equation}\label{expansion_forms}
\omega\, =\,f_{0,0}\,+\, dz\cdot f_{1,0}\,+\,f_{0,1}\cdot dz^*\,+\, dz\cdot
f_{1,1}\cdot dz^*,
\end{equation}
где $f_{i,j}\in {\rm Pol}(\mathbb{C})_q$.

Структура $U_q\mathfrak{sl}_2$-модульной алгебры в $\Omega(\mathbb{C})_q$
вводится равенствами \eqref{U_q_zstar_sl_2}, \eqref{q-C}, \eqref{KF-dz},
\eqref{E-dz}, \eqref{U_q_dzstar_sl_2}, а дифференциал $d$ -- его значениями
на образующих ($d:\,z\mapsto dz$,\ $d:\,z^*\mapsto dz^*$) и требованием
$d^2=0$. Инволюция $*$ и градуировка определяются очевидным образом.

 Пара $(\Omega(\mathbb{C})_q,d)$ является ковариантным $*$-исчислением над
$U_q\mathfrak{su}_{1,1}$-модульной алгеброй ${\rm Pol}(\mathbb{C})_q$.

Дифференциал $d$ единственным образом разлагается в сумму
 $d=\partial + \bar{\partial}$ опреаторов
 $$
 \partial: \Omega^{i,j}(\mathbb{C})_q \to \Omega^{i+1,j}(\mathbb{C})_q,
 \qquad
 \bar{\partial}: \Omega^{i,j}(\mathbb{C})_q \to
 \Omega^{i,j+1}(\mathbb{C})_q,
 $$
причем
\begin{equation}\label{partial_sl_2}
\partial^2=0,\qquad \bar{\partial}^2=0,\qquad
\partial\bar{\partial}\, +\, \bar{\partial}\partial\,=0,
\end{equation}
как следует из равенства $d^2=0$.

Полагая
$$\Omega^{0,*}(\mathbb{C})_q=\Omega^{0,0}(\mathbb{C})_q\oplus
\Omega^{0,1}(\mathbb{C})_q,\qquad\Omega^{*,\,0}(\mathbb{C})_q=
\Omega^{0,0}(\mathbb{C})_q\oplus\Omega^{1,0}(\mathbb{C})_q,
$$
получаем еще два ковариантных дифференциальных исчисления
\begin{equation}\label{two_calculi}
(\Omega^{0,*}(\mathbb{C})_q,\bar{\partial}),\qquad
(\Omega^{*,\,0}(\mathbb{C})_q,\partial).
\end{equation}
Разумеется, они не являются $*$-исчислениями, и они неизоморфны, поскольку
${\rm Ker}\,\bar{\partial}\ne{\rm Ker}\,\partial$. Отметим, что оператор
$\bar{\partial}$ играет важную роль в комплексном анализе \cite{Hermander}.


\subsubsection{Дифференциальные формы с финитными
коэффициентами.}\label{finite-forms_sl_2}

Начнем с $q$-аналога векторного пространства дифференциальных форм с
обобщенными коэффициентами.

Рассмотрим векторное пространство $\Omega(\mathbb{D})'_q$ сумм, аналогичных
\eqref{expansion_forms}
\begin{equation}\label{gen_expansion_forms}
\omega=f_{0,0}+\,dz\cdot f_{1,0}+\,f_{0,1}\cdot dz^*+\,dz\cdot f_{1,1}\cdot
dz^*,\qquad f_{i,j}\in\mathscr{D}(\mathbb{D})'_q,
\end{equation}
с коэффициентами из пространства $\mathscr{D}(\mathbb{D})'_q$ обобщенных
функций в квантовом круге, см. \itemiii\ \ref{finite_sl_2}. Наделим
$\Omega(\mathbb{D})'_q$ топологией с помощью топологии в
$\mathscr{D}(\mathbb{D})'_q$ и изоморфизма векторных пространств, при
котором $\omega\mapsto(f_{0,0},f_{1,0},f_{0,1},f_{1,1})$. Введенное в
\itemiiiе \ref{pol-forms_sl_2} векторное пространство
$\Omega(\mathbb{C})_q$ является плотным линейным подмногообразием в
$\Omega(\mathbb{D})'_q$. Другими словами, $\Omega(\mathbb{D})'_q$ --
пополнение $\Omega(\mathbb{C})_q$ в соответствующей топологии. Инволюция
$*$, действие $U_q\mathfrak{sl}_2$ и структура
$\Omega(\mathbb{C})_q$-бимодуля переносятся по непрерывности с
$\Omega(\mathbb{C})_q$ на $\Omega(\mathbb{D})'_q$, как и в \itemiiiе
\ref{finite_sl_2}.

Операторы $\partial$, $\bar{\partial}$, $d$ также продолжаются по
непрерывности с $\Omega(\mathbb{C})_q$ на $\Omega(\mathbb{D})'_q$:
\begin{equation}\label{ext_d_sl_2}
\partial\,\psi(y)=-dz\,\frac{\psi(y)-\psi(q^2\,y)}{y-q^2\,y}\,z^*,\quad
\bar{\partial}\,\psi(y)=-z\,\frac{\psi(y)-\psi(q^2\,y)}{y-q^2\,y}\,dz^*,\end{equation}
где $y=1-zz^*$ и $\psi$ -- функция на множестве $q^{2\mathbb{Z}_+}$.
(Равенство \eqref{ext_d_sl_2} достаточно получить в частном случае
$\psi(y)=y^n$, что легко сделать с помощью \eqref{comm_disc}.)

Например, для введенного равенствами \eqref{rel3_sl_2}, \eqref{rel1_sl_2}
элемента $f_0$
\begin{equation}\label{df0_sl_2}
d\,f_0\,=\,-\frac1{1-q^2}\,dz\,f_0\,z^*\;-\;\frac1{1-q^2}\,z\,f_0\,dz^*.
\end{equation}

 Разумеется, $\Omega(\mathbb{D})'_q$ является $U_q\mathfrak{sl}_2$-модульным
$\Omega(\mathbb{C})_q$-бимодулем, дифференциал $d$ -- эндоморфизм
$U_q\mathfrak{sl}_2$-модуля $\Omega(\mathbb{D})'_q$, и выполняются
равенства \eqref{Leibnitz_form}, \eqref{partial_sl_2}.

 Частные производные обобщенной функции $f \in \Omega(\mathbb{D})'_q$
 введем равенствами
 \begin{equation}\label{partial_sl_2_gen}
 \partial\,f\,=\, dz\, \frac{\partial\, f}{\partial z}, \qquad
\bar{\partial}\,f\,=\,dz^* \,\frac{\partial\, f}{\partial z^*}.
 \end{equation}

\bigskip

Плотное линейное подмногообразие
$$
\Omega(\mathbb{D})_q=\{f_{0,0}+\,dz\cdot f_{1,0}+\,f_{0,1}\cdot
dz^*+\,dz\cdot f_{1,1}\cdot
dz^*\in\Omega(\mathbb{D})'_q\,|\,f_{i,j}\in\mathscr{D}(\mathbb{D})_q\}
$$
будем называть пространством дифференциальных форм с финитными
коэффициентами в квантовом круге. \IND{пространство ! дифференциальных форм
с финитными коэффициентами в квантовом круге} Как следует из предложения
\ref{def-finite}, $\omega\in\Omega(\mathbb{D})_q$, если и только если
$z^{*N}\cdot\omega=\omega\cdot z^N=0$ при некотором $N\in\mathbb{N}$.

 Операция умножения дифференциальных форм переносится по непрерывности
 с $\Omega(\mathbb{C})_q$ на $\Omega(\mathbb{D})_q$.
 Из $f_0z=z^*f_0=0$ следуют равенства
 $$
 df_0\;z=-f_0\;dz=-dz\;f_0,\qquad z^*\;df_0=-dz^*\;f_0=-f_0\;dz^*.
$$
С их помощью нетрудно доказать,
 что алгебра $\Omega(\mathbb{D})_q$ порождается своим подпространством
$\mathscr{D}(\mathbb{D})_q \oplus d\mathscr{D}(\mathbb{D})_q$. Пара
$(\Omega(\mathbb{C})_q,d)$ является ковариантным \hbox{$*$-исчислением} над
 $U_q\mathfrak{su}_{1,1}$-модульной алгеброй $\mathscr{D}(\mathbb{D})_q$
 финитных функций в квантовом круге.

  Используя то, что $\Omega(\mathbb{D})'_q$ является
$U_q\mathfrak{sl}_2$-модульным $\Omega(\mathbb{C})_q$-бимодулем, и
соображения непрерывности, наделим
 $\Omega(\mathbb{D})'_q$ структурой
 $U_q\mathfrak{sl}_2$-модульного
 $\Omega(\mathbb{C})_q \oplus \Omega(\mathbb{D})_q$-бимодуля.

  Введем в рассмотрение дифференциальную форму,
  являющуюся $q$-аналогом фундаментальной $2$-формы, ассоциированной
  с инвариантной кэлеровой метрикой в круге \cite[стр. 140]{KobayashiNomidzu2},
  \cite[стр. 120]{GriffHarr1}.

\begin{proposition}\label{Kahler-form}
Дифференциальная форма
\begin{equation}\label{(1,1)-form}
\omega_0\,=\,\frac1{2i}\,y^{-2}\;dz\wedge dz^*\in\Omega(\mathbb{D})'_q,
\end{equation}
является $U_q\mathfrak{sl}_2$-инвариантной,
$\partial\omega_0=\bar{\partial}\omega_0=0$, и
\begin{equation}\label{central_form}
\omega\wedge\omega_0\,=\,\omega_0\wedge\omega
\end{equation}
для всех $\omega\in\Omega(\mathbb{C})_q\oplus\Omega(\mathbb{D})_q$.
\end{proposition}

{\bf Доказательство.} Неочевидна лишь $U_q\mathfrak{sl}_2$-инвариантность
$(1,1)$-формы $\omega_0$. Разумеется, $K^{\pm 1}\omega_0=\omega_0$.
Достаточно доказать одно из равенств $E\omega_0=0$, $F\omega_0=0$. Но, как
следует из \eqref{EF-psi}, $E (y^{-2})= q^{\frac{1}{2}}(1+q^2)y^{-2}z$.
Кроме того,
$$ E(dz\wedge dz^*)=d(Ez)\wedge dz^*=
d(-q^{\frac{1}{2}}z^2)\wedge dz^*=-q^{\frac{1}{2}}(1+q^2)zdz\wedge dz^*.
$$
Значит, $E(2i\omega_0)\,=\,q^{\frac{1}{2}}(1+q^2)y^{-2}zdz\wedge dz^* -
q^{\frac{1}{2}}(1+q^2)y^{-2}zdz\wedge dz^*=0$.\hfill $\square$

\begin{remark}\label{Kahler-log}
 Как следует из \eqref{ext_d_sl_2}, дифференциальная форма $\omega_0$  с
точностью до числового множителя равна $\partial\bar{\partial}\,{\rm
log}\,y$, ср. \cite[стр. 122]{GriffHarr1}, то есть функция $\varphi={\rm
log}\,y$ является $q$-аналогом Кэлерова потенциала.
\end{remark}

\bigskip
Линейный функционал
\begin{equation}\label{Lesbegue_measure}
\int\limits_{\mathbb{D}_q}\,f\,d\mu\,=\, \int\limits_{\mathbb{D}_q}\,f\cdot
y^2d\nu,\qquad f \in \mathscr{D}(\mathbb{D})_q
\end{equation}
является $q$-аналогом интеграла по нормированной мере Лебега в $\mathbb{D}$.
 Воспользуемся разложением
$$ \Omega(\mathbb{D})_q= \Omega^{0,0}(\mathbb{D})_q \oplus
\Omega^{1,0}(\mathbb{D})_q \oplus \Omega^{0,1}(\mathbb{D})_q\oplus
\Omega^{1,1}(\mathbb{D})_q.$$

\begin{proposition}\label{inv_int_for_forms}
Линейный функционал
\begin{equation} \label{inv_int_1_1}
\int\limits_{\mathbb{D}_q} \,f\, dz\wedge dz^*\; \stackrel{\rm
def}{=}\;-2\pi\, i\int\limits_{\mathbb{D}_q}\,f\,d\mu
\end{equation}
 на $\Omega^{1,1}(\mathbb{D})_q$ является $U_q\mathfrak{sl}_2$-инвариантным.
\end{proposition}

{\bf Доказательство.} Из предложения \ref{Kahler-form} следует, что
отображение
$$
\mathscr{D}(\mathbb{D})_q\to\Omega^{1,1}(\mathbb{D})_q,\qquad f\mapsto
f\,\omega_0
$$
является изоморфизмом $U_q\mathfrak{sl}_2$-модулей. Остается воспользоваться
предложением \ref{explicit_int_sl_2} об $U_q\mathfrak{sl}_2$-инвариантности
интеграла $\nu$ на $\mathscr{D}(\mathbb{D})_q$. \hfill $\square$

\bigskip

Отметим, что линейный функционал \eqref{Lesbegue_measure} допускает
продолжение по непрерывности на ${\rm Pol}(\mathbb{C})_q$.

 \begin{corollary} Линейный функционал\ \eqref{inv_int_1_1} на
 $\Omega^{1,1}(\mathbb{C})_q$
 является $U_q\mathfrak{sl}_2$-инвариантным.
\end{corollary}

\bigskip

Рассмотрим двусторонний идеал $J\subset{\rm Pol}({\mathbb C})_q$,
порожденный элементом $1-zz^*$, и коммутативную фактор-алгебру ${\rm
Pol}({\mathbb C})_q/J$. Будем отождествлять ее с алгеброй полиномов Лорана
и называть образ $f|_{\mathbb{T}}$ элемента $f\in{\rm Pol}({\mathbb C})_q$
при каноническом гомоморфизме его сужением на единичную окружность
$\mathbb{T}$.

Определим интеграл на $\Omega^{1,0}({\mathbb C})_q$ равенством
$$
\int\limits_{\mathbb{T}}dz\,f\;=\;2\pi i
\int\limits_{\mathbb{T}}(z\cdot f)|_{\mathbb{T}}\,d\nu,\qquad f\in
\mathrm{Pol}({\mathbb C})_q,
$$
где $\int\limits_{\mathbb{T}}\psi\,d\nu\stackrel{\rm def}{=}\int
\limits_0^{2 \pi}\psi(e^{i \theta})\frac{d\theta}{2\pi}.$ Следующее
утверждение является $q$-аналогом формулы Грина.
\begin{proposition}\label{like_Green}
Для всех $\psi \in \Omega^{1,0}({\mathbb C})_q$
\begin{equation}\label{(A.1)}
\int \limits_{\mathbb{D}_q}\overline{\partial}\psi=\int
\limits_{\mathbb{T}}\psi.
\end{equation}
\end{proposition}

 {\bf Доказательство}. Пусть
 $\psi \in \Omega^{1,0}({\mathbb C})_q$,
\begin{equation}\label{(A.3)}
\psi=dz \left(\sum_{m>0}z^m
\psi_m(y)+\psi_0(y)+\sum_{m>0}\psi_{-m}(y)z^{*m}\right).
\end{equation}
Можно ограничиться случаем $\psi=dz\, \psi_{-1}(y)z^*$, так как слагаемые
$dz\,\psi_j(y)\,dz^*$ дают нулевой вклад в \eqref{(A.1)} при $j \neq -1$.

 Из \eqref{ext_d_sl_2} следует равенство
 $$\overline{\partial}(dz\,\psi_{-1}(y)\,z^*)\;=\;dz \left(
\psi_{-1}(y)-q^{-2} \frac{\psi_{-1}(q^{-2}y)-\psi_{-1}(y)}{q^{-2}y-y}(1-y)
\right)dz^*.$$ Положим $g(y)=\psi_{-1}(y)-q^{-2}
\frac{\psi_{-1}(q^{-2}y)-\psi_{-1}(y)}{q^{-2}y-y}(1-y)$.
 Если $\psi_{-1}(0)=0$, то, как легко
 показать,  $\sum
\limits_{n \in{\mathbb Z}_+}g(q^{2n})q^{2n}=0$. Значит, в этом случае
 $\int
\limits_{\mathbb{D}_q}\overline{\partial}\psi=\int
\limits_{\mathbb{T}}\psi=0$, и предложение \ref{like_Green} доказано для
(1,0)-форм из линейного подпространства коразмерности 1. Остается убедиться
в его справедливости при $\psi=dz \, z^*$:
$$
\int\limits_{\mathbb{D}_q}\bar{\partial}(dz \, z^*) =
 -\int_{\mathbb{D}_q}dz\wedge dz^* = 2\pi i,\qquad
 \int\limits_{\mathbb{T}} dz\,z^* = 2\pi i.
\eqno \square
$$

\medskip
Линейные функционалы в обеих части доказанного равенства продолжаются по
непрерывности на $\Omega^{1,0}(\mathbb{C})_q\oplus
\Omega^{1,0}(\mathbb{D})_q$.

\begin{corollary}\label{to_conjugate}\ \ $\int
\limits_{\mathbb{D}_q}\overline{\partial}\psi=0$ для всех $\psi \in
\Omega^{1,0}(\mathbb{D})_q$.
\end{corollary}

\subsubsection{Инвариантные эрмитовы метрики и инвариантный
лапласиан.}\label{Hodges}

Дифференциальный оператор
\begin{equation}\label{c_Laplace}
\Box=-(1-|z|^2)^2\frac{\partial}{\partial
z}\frac{\partial}{\partial\overline{z}}
\end{equation}
в единичном круге называется инвариантным лапласианом. \IND{инвариантный !
лапласиан} В выборе знака в \eqref{c_Laplace} мы следуем традициям
комплексного анализа \cite[стр. 159]{Beals}.

Напомним доказательство $SU_{1,1}$-инвариантности этого оператора $\Box$.
Векторные пространства гладких финитных функций и гладких $(0,1)$-форм с
финитными коэффициентами можно наделить $SU_{1,1}$-инвариантными
положительными эрмитовыми формами
\begin{equation*}
(f_1,f_2)=\int\limits_\mathbb{D}\overline{f}_2f_1\,
\frac{d\operatorname{Re}z\wedge d\operatorname{Im}z}{(1-|z|^2)^2},\quad
(f_1d\overline{z},f_2d\overline{z})=\int\limits_\mathbb{D}\overline{f}_2f_1\,
d\operatorname{Re}z\wedge d\operatorname{Im}z.
\end{equation*}

Из $SU_{1,1}$-инвариантности дифференциального оператора
$\overline{\partial}$ следует $SU_{1,1}$-инвариантность формально
сопряженного дифференциального оператора $ \overline{\partial}^{\, *}:
fd\overline{z}\mapsto -(1-|z|^2)^2\;\frac{\partial f}{\partial z}$.
  Значит, оператор $\Box=\overline{\partial}^{\, *}\overline{\partial}$
  также $SU_{1,1}$-инвариантен.
Воспользуемся аналогичными рассуждениями в случае квантового круга.

\medskip В \itemiiiе \ref{IntStar} рассматривались инвариантные
полуторалинейные формы. Введем более общие понятия, заменив поле
комплексных чисел $A$-модульной алгеброй $F$. Мы следуем неопубликованной
работе Шклярова о $q$-аналогах кэлеровых многообразий.

Рассмотрим $*$-алгебру $F$. Если $V$ является левым $F$-модулем, то
$\overline{V}$ наделяется структурой правого $F$-модуля
$$vf\stackrel{\rm def}{=}f^*v,\qquad f\in F,\;v\in V,$$
а если $V$ является правым $F$-модулем, то $\overline{V}$ наделяется
структурой левого $F$-модуля:
$$fv\stackrel{\rm def}{=}vf^*,\qquad f\in F,\;v\in V.$$
Рассмотрим $F$-бимодуль $V$. Полуторалинейное отображение $\mathbf{h}:
V\times V\to F$ назовем эрмитовой метрикой, \IND{эрмитова метрика} если оно
вещественно
$$\mathbf{h}(v_2,v_1)\,=\,\mathbf{h}(v_1,v_2)^*,\qquad v_1,v_2\in V,$$
и отвечающее ему линейное отображение
\begin{equation}\label{hermit_form}
{\mathbf h}:\;\overline{V}\otimes_{F}V\;\to\;F
\end{equation}
является морфизмом $F$-бимодулей.

Эрмитову метрику назовем инвариантной, \IND{эрмитова метрика !
инвариантная} если $V$ является $A$-мо\-дульным $F$-бимодулем и отображение
\eqref{hermit_form} -- морфизмом $A$-модулей.

Эрмитову метрику назовем неотрицательной, \IND{эрмитова метрика !
неотрицательная} если $\mathbf{h}(v,v)\ge 0$ для всех $v\in V$, и
положительной, \IND{эрмитова метрика ! положительная} если, кроме
того, $\mathbf{h}(v,v)\ne 0$ при $v\ne 0$.

Приведем примеры эрмитовых метрик. Очевидно,
$\Omega^{0,0}(\mathbb{D})_q$, $\Omega^{1,0}(\mathbb{D})_q$,
$\Omega^{0,1}(\mathbb{D})_q$, $\Omega^{1,1}(\mathbb{D})_q$
являются как свободными правыми, так и свободными левыми
$\mathscr{D}(\mathbb{D})_q$-модулями с образующими $1$, $dz$,
$dz^*$, $dz\wedge dz^*$ соответственно.

\begin{proposition}\label{hermitian_metrics}  Рассмотрим $*$-алгебру Хопфа
$U_q\mathfrak{su}_{1,1}$, $U_q\mathfrak{su}_{1,1}$-модульную $*$-алгебру
$\mathscr{D}(\mathbb{D})_q$ и $U_q\mathfrak{su}_{1,1}$-модульный
$\mathscr{D}(\mathbb{D})_q$-бимодуль $\Omega^{0,1}(\mathbb{D})_q$. Следующие
полуторалинейные отображения являются положительными инвариантными
эрмитовыми метриками:
\begin{equation}\label{metric_00}
\mathbf{h}_{0,0}(f_1\,,\,f_2)\;=\;f_2^*f_1,
\end{equation}
\begin{equation}\label{metric_01}
\mathbf{h}_{0,1}(dz^*\,f_1\,,\,dz^*\,f_2)\;=\;f_2^*\,(1-zz^*)^2\,f_1,
\end{equation}
где $f_1,f_2 \in \mathscr{D}(\mathbb{D})_q$.
\end{proposition}

{\bf Докзательство.} Инвариантность эрмитовых метрик \eqref{metric_00},
\eqref{metric_01} следует из предложения \ref{Kahler-form} и ковариантности
$*$-исчисления $(\Omega(\mathbb{D})_q,d)$. Положительность вытекает из
положительности линейного оператора
$$T_F(1-zz^*):e_n\mapsto q^{2n}\,e_n,\qquad n\in\mathbb{Z}_+,$$
поскольку фоковское представление осуществляет изоморфизм $*$-алгебры
$\mathscr{D}(\mathbb{D})_q$ и $*$-алгебры финитных линейных операторов в
$\mathcal{H}$. \hfill $\square$

\medskip

\begin{corollary}\label{inv_forms_sl_2}
Эрмитовы формы
\begin{equation}\label{new-inner-product}
(f_1,f_2)\,=\, \int\limits_{\mathbb{D}_q} \,f_2^*f_1\,d\nu,\qquad
(dz^*\,f_1,dz^*\,f_2)\,=\, \int\limits_{\mathbb{D}_q}
\,f_2^*\,(1-zz^*)^2\,f_1\,d\nu
\end{equation}
в $\mathscr{D}(\mathbb{D})_q$, $\Omega^{0,1}(\mathbb{D})_q$ являются
инвариантными и положительными.
\end{corollary}

\begin{remark}\label{old-inner-product-remark}\ Как следует из определений,
\begin{equation}\label{old-inner-product}
(f_1dz^*\,,\,f_2dz^*)\,=\, \int\limits_{\mathbb{D}_q}
\,f_2^*\,f_1\,\,(1-zz^*)^2\,d\nu,\qquad f_1,f_2 \in
\mathscr{D}(\mathbb{D})_q.
\end{equation}
\end{remark}

\bigskip
Рассмотрим линейный оператор $\bar{\partial}$ из предгильбертова
пространства $\mathscr{D}(\mathbb{D})_q$ в предгильбертово пространство
$\Omega^{0,1}(\mathbb{D})_q$.

\begin{proposition}\label{conjugate_d}
Линейный оператор
$$
\bar{\partial}^*:\;\Omega^{0,1}(\mathbb{D})'_q\to\mathscr{D}(\mathbb{D})'_q,
\qquad\bar{\partial}^*:\;dz^*\,f\mapsto-q^2\,(1-zz^*)^2\,
\frac{\partial\,f}{\partial z}
$$
является сопряженным к $\bar{\partial}$, точнее, для любых
$f_1,f_2\in\mathscr{D}(\mathbb{D})_q$
\begin{equation}\label{conj_d_bar}
(\bar{\partial}f_1\,,\,dz^*f_2)\,=\left(f_1\,,\,-q^2\,(1-zz^*)^2\,
\frac{\partial f_2}{\partial z}\right).
\end{equation}
\end{proposition}

{\bf Доказательство.} Используя равенство
$\bar{\partial}(f_2^*)=(\partial\,f_2)^*$ и следствие \ref{to_conjugate},
получаем
\begin{multline*}
0=\int\limits_{\mathbb{D}_q}\,\bar{\partial}\,(f_2^*\,dz\,f_1)\,=\,
\int\limits_{\mathbb{D}_q}\,\bar{\partial}(f_2^*)\wedge dz\,f_1\;-\;
\int\limits_{\mathbb{D}_q}\,f_2^*\,dz\wedge\bar{\partial}f_1=
\\ =\int\limits_{\mathbb{D}_q}\,(\partial f_2)^*\wedge dz\,f_1-
\int\limits_{\mathbb{D}_q}\,f_2^*\,dz\,\wedge\bar{\partial}f_1=
\int\limits_{\mathbb{D}_q}\,\left(\frac{\partial f_2}{\partial
z}\right)^*dz^*\wedge
dzf_1-\int\limits_{\mathbb{D}_q}\,f_2^*dz\wedge\bar{\partial}f_1
\end{multline*}
для любых $f_1, f_2 \in \mathscr{D}(\mathbb{D})_q$.
 Значит,
$$\int\limits_{\mathbb{D}_q}\,(dz^*f_2)^*\wedge\bar{\partial}f_1\,=\, \int\limits_{\mathbb{D}_q}\,
\left(\frac{\partial f_2}{\partial z}\right)^*\,dz^*\wedge dz\,f_1.
$$
Следовательно,
$$\int\limits_{\mathbb{D}_q}\,f_2^*dz \wedge dz^*
\frac{\partial f_1}{\partial z^*} \,= -q^2\,\int\limits_{\mathbb{D}_q}\,
\left(\frac{\partial f_2}{\partial z}\right)^*\,dz\wedge dz^*\,f_1.
$$
Остается воспользоваться равенствами  \eqref{(1,1)-form},
\eqref{central_form} и \eqref{inv_int_1_1}. \hfill $\square$

\begin{corollary}\label{inv_L}
 Линейный оператор
\begin{equation}\label{invariant_laplacian}
\square_q=\bar{\partial}^*\bar{\partial}= -q^2\,(1-zz^*)^2
\,\frac{\partial\,}{\partial z}\frac{\partial\,}{\partial z^*}
\end{equation}
в $\mathscr{D}(\mathbb{D})_q$ является эндоморфизмом
$U_q\mathfrak{sl}_2$-модуля $\mathscr{D}(\mathbb{D})_q$.
\end{corollary}

{\bf Доказательство.} Инвариантность эрмитовой формы $(\square_q\,f_1,f_2)$
в $\mathscr{D}(\mathbb{D})_q$ вытекает из предложения \ref{conjugate_d},
следствия \ref{inv_forms_sl_2} и из того, что оператор $\bar{\partial}$
является морфизмом $U_q\mathfrak{sl}_2$-модулей. \hfill $\square$

\medskip

Линейный оператор $\square_q$ допускает продолжение по непрерывности до
эндоморфизма $U_q\mathfrak{sl}_2$-модуля $\mathscr{D}(\mathbb{D})'_q$, и мы
будем называть его инвариантным лапласианом в квантовом круге.
\IND{инвариантный ! лапласиан ! в квантовом круге} Введем в рассмотрение
пополнение $L^2(d\nu)_q$ векторного пространства
$\mathscr{D}(\mathbb{D})_q$ по норме
$\|f\|=\sqrt{\int_{\mathbb{D}}f^*fd\nu}$. Как вытекает из предложения
\ref{conjugate_d}, инвариантный лапласиан является симметрическим линейным
оператором в гильбертовом пространстве $L^2(d\nu)_q$ и, следовательно,
допускает замыкание \cite[стр. 355, 356]{Dunford2}. Мы сохраним для
замыкания обозначение $\square_q$. Как будет показано в \itemiiiе
\ref{boundness_Laplace}, оператор $\square_q$ самосопряжен и ограничен.

\subsubsection{Радиальная часть инвариантного
лапласиана.}\label{radial_part}

Рассмотрим однопараметрическую группу линейных операторов
$$
\sum_{j\in \mathbb{Z}_+} z^j \psi(y)\;+\; \sum_{k\in \mathbb{N}}
\psi(y)z^{*k}\mapsto \sum_{j\in \mathbb{Z}_+} e^{ij\varphi}z^j \psi(y)\;+\;
\sum_{k\in \mathbb{N}} e^{-ik\varphi}\,\psi(y) z^{*k}, \quad \varphi\in
\mathbb{R}
$$
в $\mathscr{D}(\mathbb{D})_q'$, см. \eqref{expansion_ll}.
 Cужение  линейного оператора
$\square_q$ на подпространство ее инвариантов обозначим $\square_q^{(0)}$ и
назовем радиальной частью оператора $\square_q$.
 Пусть
\begin{equation}\label{D-op}
\mathcal{D}_q\psi(x)=\frac{\psi(q^{-1}x)\,-\,\psi(qx)}{q^{-1}x\,-\,qx}.
\end{equation}
 Наделим векторное пространство функций на множестве
$q^{-2\mathbb{Z}_+}$ топологией поточечной сходимости и рассмотрим
непрерывный линейный оператор $ \mathcal{D}_qx(1-q^{-1}x)\mathcal{D}_q$ в
этом пространстве.

\begin{proposition}\label{explicit_square}
Для любой функции $\psi(x)$ на множестве $q^{-2\mathbb{Z}_+}$
\begin{equation}\label{expl_sq}
\Box_q^{(0)}\psi(x)=\mathcal{D}_qx(1-q^{-1}x)\mathcal{D}_q\psi(x),
\end{equation}
где в левой части равенства $x=y^{-1}$.
\end{proposition}

 {\bf Доказательство.} Достаточно рассмотреть частный случай функции
$\psi(x)$ с конечным носителем. Пусть $$\int \limits_1^\infty
\psi(t)d_{q^{-2}}t \stackrel{\rm def}{=}(q^{-2}-1)\sum_{m=0}^\infty
\psi(q^{-2m})q^{-2m}.$$ Очевидно,
  $\|\psi(x)\|^2=q^2 \int \limits_1^\infty|\psi(t)|^2d_{q^{-2}}t$.
Достаточно доказать, что
\begin{equation}\label{2.5.5}
  \|\overline{\partial}\psi(x)\|^2=q^2 \int
  \limits_1^\infty(\mathcal{D}_qt(1-q^{-1}t)\mathcal{D}_q
\psi(t))\overline{\psi(t)}d_{q^{-2}}t.
\end{equation}
Из \eqref{ext_d_sl_2} следует равенство
$
 \bar{\partial}\,\psi(x)\,=\,
 -z\,\frac{\psi(x)-\psi(q^{-2}\,x)}{x^{-1}-q^2x^{-1}}\,dz^*.
$ для любой функции $\psi(x)$ с конечным носителем в $q^{-2{\mathbb Z}_+}$.

 С помощью  \eqref{old-inner-product} получаем
$$
\|\overline{\partial}\psi(x)\|^2=(1-q^2)^{-1}\sum_{m=0}^\infty
\left|\psi(q^{-2m})-\psi(q^{-(2m+2)})\right|^2(1-q^{2m+2})q^{-2m}=$$
$$=q^{-2}(q^{-2}-1)^{-2}\int \limits_1^\infty
\left|\psi(t)-\psi(q^{-2}t)\right|^2(1-q^2t^{-1})d_{q^{-2}}t=$$
$$=-\int \limits_1^\infty
\left|\frac{\psi(t)-\psi(q^{-2}t)}{t-q^{-2}t}\right|^2 \cdot
t(1-q^{-2}t)d_{q^{-2}}t.
$$
 Доказательство завершается с помощью $q$-аналога
формулы ''интегрирования по частям'' (\ref{by_part}). \hfill $\square$

\medskip

В дальнейшем используется следующий базис в пространстве функций с конечными
носителями на множестве $q^{-2\mathbb{Z}_+}$:
\begin{equation}\label{f_j}
f_j(x)=\left\{\begin{array}{cc}1, &\quad\text{при}\quad x=q^{-2j}
\\ 0, &\quad \text{при}\quad x\ne q^{-2j}\end{array}\right.,\qquad
j\in\mathbb{Z}_+.
\end{equation}

Рассмотрим элемент $C_q$ центра алгебры $U_q\mathfrak{sl}_2$, введенный в
\itemiiiе \ref{SL_2_sl_2}.
\begin{corollary}\label{C_q-Laplace}
Для всех $f \in \mathscr{D}(\mathbb{D})_q$ имеем $C_qf\,=\,-q \square_q f$.
\end{corollary}

{\bf Доказательство.} Сравниваемые линейные операторы перестановочны с
действием $U_q\mathfrak{sl}_2$ в $\mathscr{D}(\mathbb{D})_q$. Пусть
$x=y^{-1}$. Из предложения \ref{t2.3.9} следует, что достаточно доказать
равенство $C_qf_0\,=\,-q\square_qf_0$. Используя (\ref{EF-psi}), получаем
равенство
$$C_q f_0=-\frac1{q^{-1}-q}f_0+\frac{q^2}{q^{-1}-q}f_1,$$
а используя предложение \ref{explicit_square} -- равенство
$$
\square_qf_0=\mathcal{D}_qx(1-q^{-1}x)\mathcal{D}_qf_0=\frac{f_0}{1-q^2}
-\frac{q^2f_1}{1-q^2}.\eqno\square
$$

\medskip

\begin{corollary} Для любой функции $\psi(t)$ на $q^{-2{\mathbb Z}_+}$
\begin{equation}\label{2.5.7}
C_q \psi(x)\;=\; -q\mathcal{D}_qx(1-q^{-1}x)\mathcal{D}_q \psi(x),
\end{equation}
где в левой части равенства $x=y^{-1}$.
\end{corollary}

\medskip

Пусть $L^2(1,\infty)_q$ -- гильбертово пространство функций $\psi(x)$ на
$q^{-2\mathbb{Z}_+}$, для которых $
\|\psi\|=\sqrt{\int\limits_1^\infty|\psi(x)|^2d_{q^{-2}}\,x} <\infty. $
Докажем лемму, из которой мгновенно следует, что оператор
\begin{equation}\label{def_L}
L=\mathcal{D}_qx(1-q^{-1}x)\mathcal{D}_q
\end{equation}
в $L^2(1,\infty)_q$ корректно определен, ограничен и самосопряжен.

\begin{lemma}\label{Jacobi_matr}
Векторы
$$e_j=(q^{-2}-1)^{-\frac12}\,q^j\,f_j,\qquad j\in\mathbb{Z}_+,$$
образуют ортонормированный базис в $L^2(1,\infty)_q$, и
\begin{equation}\label{Jacoby}
Le_j= \begin{cases} b_{j-1}e_{j-1}+a_je_j+b_je_{j+1},& j\in \mathbb{N},\\
a_je_j+b_je_{j+1},& j=0
      \end{cases},
\end{equation}
где
$$
a_j=\frac{1+q^2-2q^{2(j+1)}}{(1-q^2)^2},\qquad
b_j=-q\,\frac{1-q^{2(j+1)}}{(1-q^2)^2}.
$$
\end{lemma}

{\bf Доказательство.} Первое утверждение очевидно, а второе вытекает из
равенства
\begin{equation}\label{L_f_j} Lf_j=-\frac{1-q^{2j}}{(1-q^2)^2}f_{j-1}+\frac{1+q^2-2q^{2(j+1)}}{(1-q^2)^2}f_j
-q^2\frac{1-q^{2(j+1)}}{(1-q^2)^2}f_{j+1},
\end{equation}
которое нетрудно получить с помощью \eqref{def_L}. (В \eqref{L_f_j}
подразумевается, что $f_{-1}=0$.) \hfill $\square$

\medskip
Очевидно, линейная оболочка множества $\{L^kf_0\,|\, k\in\mathbb{Z}_+\}$
плотна в $L^2(1,\infty)_q$. Значит \cite[стр. 279]{AG}, ограниченный
самосопряженный линейный оператор $L$ унитарно эквивалентен оператору
умножения $f(\lambda)\mapsto\lambda f(\lambda)$ в гильбертовом пространстве
$L^2(d\mu)$, где $d\mu(\lambda)$ -- борелевская мера с компактным носителем
на вещественной оси.

Зафиксируем эту меру и оператор $U:L^2(1,\infty)_q\to L^2(d\mu)$,
осуществляющий унитарную эквивалентность, требованием $U f_0=q^{-2}-1$.

Найдем меру $d\mu$ и явный вид унитарного оператора $U$.

Будем считать известными основные обозначения теории $q$-спе\-ци\-аль\-ных
функций, см.
\itemiii\ \ref{q-series}. Пусть $d\sigma(\rho)$ -- абсолютно непрерывная
мера на $\left[0,\frac{\pi}{2\log(q^{-1})}\right]$, определяемая равенством
\begin{equation}\label{Pl}
\frac{d\sigma(\rho)}{d\rho}=\frac1{\pi}\;\frac{\log\frac1{q}}{1-q^2}\;
\frac{d\rho}{c(-\frac12+i\rho)c(-\frac12-i\rho)},
\end{equation}
где
\begin{equation}\label{c_l}
c(l)=\dfrac{(q^{2(l+1)};q^2)_\infty^2}
{(q^{2(l+1)};q^2)_\infty(q^2;q^2)_\infty}\,=\,
\frac{\Gamma_{q^2}(2l+1)}{(\Gamma_{q^2}(l+1))^2}.
\end{equation}

Отметим, что $c(l)$ является $q$-аналогом $c$-функции Хариш-Чандры, см.
\itemiii\ \ref{c-function_sl_2}.

\medskip

\begin{proposition}\label{Fourier_disc}
\begin{itemize}
\item[1.]  Рассмотрим функцию
\begin{equation}\label{Phi_l}
\Phi_l(x)=\sideset{_3}{_2}{\mathop{\phi}}
\left[\begin{array}{c}x,q^{-2l},q^{2(l+1)}\\
q^2,0\end{array};q^2,q^2\right],\qquad l\in\mathbb{C}
\end{equation}
на множестве $q^{-2\mathbb{Z}_+}$. Она равна $1$ при $x=1$, и
\begin{equation}\label{eigen_L}
\mathcal{D}_qx(1-q^{-1}x)\mathcal{D}_q\Phi_l(x)=\lambda(l)\Phi_l(x),
\end{equation}
где
\begin{equation}\label{Laplace_eigenvalue}
\lambda(l)=\dfrac{(1-q^{-2l})(1-q^{2(l+1)})}{(1-q^2)^2}.
\end{equation}
\item[2.] Функция $\Phi_{-\frac12+i\rho}(x)$ на $q^{-2\mathbb{Z}_+}$ {\bf
вещественна при вещественных} $\rho$, и интеграл $\int\limits_1^\infty
f(x)\Phi_{-\frac12+i\rho}(x)d_{q^{-2}}x$ принадлежит $L^2(d\sigma)$, если
$f(x)$ -- функция с конечным носителем на множестве $q^{-2\mathbb{Z}_+}$.
Линейный оператор
\begin{equation}\label{lin_op}
f(x)\mapsto\int\limits_1^\infty f(x)\Phi_{-\frac12+i\rho}(x)d_{q^{-2}}x
\end{equation}
допускает продолжение по непрерывности до унитарного оператора из
$L^2(1,\infty)_q$ в $L^2(d\sigma)$.
\end{itemize}
\end{proposition}

{\bf Доказательство.} Векторы $P_j=Ue_j$, $j\in\mathbb{Z}_+$, гильбертова
пространства $L^2(d\mu)$ являются ортогональными полиномами переменной
$\lambda$
$$
\int P_i(\lambda)P_j(\lambda)d\mu(\lambda)=\delta_{ij},\qquad
i,j\in\mathbb{Z}_+.
$$
Как следует из определения оператора $U$ и из леммы \ref{Jacobi_matr},
\begin{multline}\label{P_rec}
\lambda P_j(\lambda)=-q\frac{1-q^{2(j+1)}}{(1-q^2)^2}P_{j+1}(\lambda)+
\left(\frac{1+q^2}{(1-q^2)^2}-\frac{2q}{(1-q^2)^2}q^{2j+1}\right)
P_j(\lambda)-
\\ -q\frac{1-q^{2j}}{(1-q^2)^2}P_{j-1}(\lambda),\qquad j\in\mathbb{Z}_+,
\end{multline}
\begin{equation}\label{init_P0}
P_0(\lambda)=Ue_0=(q^{-2}-1)^{-\frac12}Uf_0=(q^{-2}-1)^{\frac12}.
\end{equation}
В \eqref{P_rec} подразумевается, что $P_{-1}=0$:
\begin{equation}\label{init_P1}
\lambda P_0(\lambda)=-\frac{q}{1-q^2}P_1(\lambda)+\frac1{1-q^2}P_0(\lambda).
\end{equation}
Рекуррентные соотношения \eqref{P_rec} и начальные данные \eqref{init_P0},
\eqref{init_P1} однозначно определяют ортогональные полиномы $P_j(\lambda)$,
$j\in\mathbb{Z}_+$. Сравним их с ортогональными полиномами Аль-Салама-Чихары
$Q_j(x;q,q|q^2)$ (см. \eqref{bible_3.8.1}). Из \eqref{bible_3.8.2} вытекают
следующие рекуррентные соотношения для ''нормированных'' полиномов
$\widetilde{Q}_j\overset{\mathrm{def}}{=}\dfrac1{(q^2;q^2)_j}Q_j$:
\begin{equation}\label{Q_rec}
x\widetilde{Q}_j=\dfrac12\left(1-q^{2(j+1)}\right)\widetilde{Q}_{j+1}+
q^{2j+1}\widetilde{Q}_j+\dfrac12\left(1-q^{2j}\right)\widetilde{Q}_{j-1},
\end{equation}
где подразумевается, что $\widetilde{Q}_{-1}=0$. Из явной формулы
\eqref{bible_3.8.1}
$$
Q_j(x;q,q|q^2)=\frac{(q^2;q^2)_j}{q^j}\sideset{_3}{_2}{\mathop{\phi}}
\left[\begin{array}{c}q^{-2j},qt,qt^{-1}\\ q^2,0\end{array};q^2,q^2\right],
$$
где $x=\dfrac12(t+t^{-1})$, получаем ''начальные данные''
\begin{equation}\label{init_Q0}
\widetilde{Q}_j(x)=1.
\end{equation}
Сравнивая \eqref{P_rec}, \eqref{init_P0} с \eqref{Q_rec}, \eqref{init_Q0},
видим, что полиномы $P_j$ и $\widetilde{Q}_j$ связаны заменой переменной
\begin{equation}\label{change_x}
\lambda=\frac{1+q^2}{(1-q^2)^2}-\frac{2q}{(1-q^2)^2}x.
\end{equation}
Так как $x=\dfrac12(t+t^{-1})$, то
$\lambda=\frac{(1-qt)(1-qt^{-1})}{(1-q^2)^2},$
$$
P_j(\lambda)=(q^{-2}-1)^{\frac12}q^{-j}\sideset{_3}{_2}{\mathop{\phi}}
\left[\begin{array}{c}q^{-2j},qt,qt^{-1}\\ q^2,0\end{array};q^2,q^2\right],
$$
$$
Uf_j=(q^{-2}-1)\cdot q^{-2j}\sideset{_3}{_2}{\mathop{\phi}}
\left[\begin{array}{c}q^{-2j},qt,qt^{-1}\\ q^2,0\end{array};q^2,q^2\right].
$$
Подстановка $t=q^{-2i\rho}$ дает $Uf_j=\int\limits_1^\infty
f_j(x)\Phi_{-\frac12+i\rho}(x)d_{q^{-2}}x$.
 Значит, для любой функции
$f(x)$ на $q^{-2\mathbb{Z}_+}$, имеющей конечный носитель,
\begin{equation}\label{U_f}
Uf=\int\limits_1^\infty f(x)\Phi_{-\frac12+i\rho}(x)d_{q^{-2}}x.
\end{equation}
Завершим доказательство предложения \ref{Fourier_disc}. Очевидно,
$$
\Phi_l(1)=1,\qquad\operatorname{Im}\Phi_{-\frac12+i\rho}(q^{-2j})=0,\qquad
\rho\in\mathbb{R},\;j\in\mathbb{Z}_+.
$$
Равенство \eqref{eigen_L} следует из \eqref{U_f} и из того, что
$$U\mathcal{D}_qx(1-q^{-1}x)\mathcal{D}_qf=\lambda Uf.$$
Остальные утверждения нетрудно получить из \eqref{U_f} с помощью замены
переменной \eqref{change_x} и соотношений ортогональности
\eqref{bible_3.8.2} для полиномов Аль-Салама-Чихары. \hfill $\square$

\begin{remark} \label{classical_limit} Так как
$\lim\limits_{q\to 1}\,\Gamma_q(x)\,=\,\Gamma(x)$ при $-x\notin
\mathbb{Z}_+$, то в пределе $q\to 1$
 $$
 \frac{d\sigma(\rho)}{d\rho}\,
   =\,\frac{1}{2\pi}\,
 \frac{\Gamma(\frac{1}{2}+i\rho)^2}{\Gamma(2i\rho)^2}\;
 \frac{\Gamma(\frac{1}{2}-i\rho)^2}{\Gamma(-2i\rho)^2}.
$$
Используя хорошо известные свойства гамма-функции
$$\Gamma(x)\Gamma(1-x)=\frac{\pi}{\sin \pi x},\qquad
\Gamma(1/2+x)\Gamma(1/2-x)=\frac{\pi}{\cos \pi x},$$ получаем равенство
 $$
 \frac{d\sigma(\rho)}{d\rho}\,
   =\,\frac{1}{2\pi}\,
 \left(\frac{\pi}{{\rm ch}(\pi \rho)}\right)^2\;
 \frac{2\rho\; {\rm \sh (2\pi \rho)}}{\pi}=2\rho\,{\rm th}(\pi \rho),
$$
что согласуется с известной формулой Мелера-Фока \cite[стр. 232]{AG_2}.
\end{remark}

\begin{remark} Приведенное выше простое  доказательство предложения
\ref{Fourier_disc} нам сообщил Эрик Кёлинк. Первоначальное доказательство,
см. \cite{SSV_disc4}, опирается на известные результаты спектральной теории
самосопряженных операторов в гильбертовом пространстве \cite[глава 10,
параграф 6]{Dunford2}. Используемый в этом доказательстве метод обладает
очень широкой областью применимости и является совершенно стандартным.
\end{remark}

\subsubsection{Ограниченность и обратимость инвариантного
лапласиана.}\label{boundness_Laplace}

Как следует из предложения \ref{Fourier_disc}, радиальная часть
$\square_q^{(0)}$ инвариантного лапласиана унитарно эквивалентна оператору
умножения на
\begin{equation}\label{l-rho}
\lambda\left(-\frac12\,+\,i\rho\right)\,=\,
\frac{1\,-\,2q\cos(2\log(q^{-1})\rho)\,+\,q^2}{(1\,-\,q^2)^2}
\end{equation}
в гильбертовом пространстве $L^2(d\sigma(\rho))$ функций на отрезке
$\left[0,\frac{\pi}{2\log(q^{-1})}\right]$. На этом отрезке функция
\eqref{l-rho} монотонно возрастает от $(1+q)^{-2}$ до $(1-q)^{-2}$. Значит,
оператор $\square_q^{(0)}$ ограничен и обратим:
$$(1+q)^{-2}\le\;\square_q^{(0)}\;\le(1-q)^{-2}.$$

Введем в рассмотрение серию $*$-представлений алгебры
$U_q\mathfrak{su}_{1,1}$ и с ее помощью получим такие же оценки для всего
оператора $\square_q$. Из этих оценок будут следовать ограниченность и
обратимость инвариантного лапласиана.

\bigskip

Рассмотрим векторное пространство $\mathcal{E}$ с базисом
$\{e_k\}_{k\in\mathbb{Z}}$ и комплексное число $l$, для которого
\begin{equation}\label{Im_l}
-\frac{\pi}{2\log(q^{-1})}\le\,{\rm Im}\,l\,\le\,\frac{\pi}{2\log(q^{-1})}.
\end{equation}
Нетрудно показать, что равенства
\begin{equation*}\label{principal_spherical_series}
\pi_l(E)\,e_k=\,q^k\,\frac{q^{-(k-l)}-q^{k-l}}{q^{-1}-q}\,e_{k+1},\quad
\pi_l(F)\,e_k=\,-q^{-(k-1)}\,\frac{q^{-(k+l)}-q^{k+l}}{q^{-1}-q}\,e_{k-1},
\end{equation*}
$$\pi_l(K^{\pm 1})e_k\,=\,q^{\pm 2k}\,e_k,\qquad k\in\mathbb{Z},$$
определяют представление $\pi_l$ алгебры $U_q\mathfrak{sl}_2$ и
\begin{equation}\label{2.5.8}
\pi_l(C_{q})=\frac{(q^{-l}-q^l)(q^{-(l+1)}-q^{l+1})}{(q^{-1}-q)^2}\,=\,
-q\dfrac{(1-q^{-2l})(1-q^{2(l+1)})}{(1-q^2)^2}.
\end{equation}

Это множество представлений алгебры $U_q\mathfrak{sl}_2$, рассматриваемых с
точностью до эквивалентности, называют сферической основной серией.
\IND{основная серия ! сферическая}

 Представления $\pi_l$
естественно реализуются в пространствах однородных обобщенных функций
 на квантовом конусе, см. \itemiiiе \ref{q-cone}.

Введем $U_q\mathfrak{sl}_2$-модуль $V^{(l)}$, отвечающий представлению
$\pi_l$. Из определений следует

\begin{lemma}\label{univ_pi_l} Если $l$ не принадлежит $\mathbb{Z}$,
то $U_q\mathfrak{sl}_2$-модуль $V^{(l)}$ прост и может быть описан с
помощью образующей $e_0$ и определяющих соотношений
$$K^{\pm 1}e_0\,=\,e_0,\qquad C_q
e_0\,=\,\frac{(q^{-l}-q^l)(q^{-(l+1)}-q^{l+1})}{(q^{-1}-q)^2}\, e_0.$$
\end{lemma}

\bigskip Последнее утверждение леммы означает биективность морфизма
$U_q\mathfrak{sl}_2$-модулей
$$U_q\mathfrak{sl}_2/\mathcal{J}_l\to V^{(l)},\qquad 1 \mapsto e_0,$$
где $\mathcal{J}_l$ -- наименьший левый идеал, содержащий элементы
$$
K^{\pm 1}-1,\qquad
C_{q}-\frac{(q^{-l}-q^l)(q^{-(l+1)}-q^{l+1})}{(q^{-1}-q)^2}.
$$

\begin{remark}\label{equiv_sl_2}
Отметим, что края полосы $-\frac{\pi}{2\log(q^{-1})}\le\,{\rm
Im}\,l\,\le\,\frac{\pi}{2\log(q^{-1})}$ можно считать склеенными, поскольку
\begin{equation}\label{trivial_iso}
\pi_{(s-i\frac{\pi}{2\log(q^{-1})})}=\pi_{(s+i\frac{\pi}{2\log(q^{-1})})},
\qquad s\in\mathbb{R}.
\end{equation}
Из леммы \ref{univ_pi_l} следует, что представления $\pi_l$ и $\pi_{-1-l}$
эквивалентны, если $l\notin\mathbb{Z}$. Нетрудно показать, что это описание
пар эквивалентных представлений сферической основной серии является
исчерпывающим.
\end{remark}

\bigskip Пусть $\mathrm{Re}\,l\,=\,-\frac12$. Наделим векторное пространство
$\mathcal{E}=V^{(l)}$
структурой предгильбертова пространства так, чтобы базис
$\{e_k\}_{k\in\mathbb{Z}}$ стал ортонормированным. Из определений следует,
что $\pi_l$ -- $*$-представление $*$-алгебры Хопфа
$U_q\mathfrak{su}_{1,1}$. Это сферическая основная унитарная серия.
\IND{основная серия! сферическая ! унитарная}

\bigskip

\begin{proposition}\label{p1.3.1}
Для всех ненулевых $f \in \mathscr{D}(\mathbb{D})_q$
\begin{equation}\label{1.3.1}
(1+q)^{-2}\,\leq\;
\frac{(\overline{\partial}\,f\,,\,\overline{\partial}\,f)}{(f,f)}
 \leq\;  (1-q)^{-2}.
\end{equation}
\end{proposition}

\medskip

{\bf Доказательство.} Рассмотрим меру $d\sigma(\rho)$ на отрезке
$\left[0,\frac{\pi}{2\log(q^{-1})}\right]$, введенную в предыдущем
\itemiiiе. Тензорное произведение предгильбертовых пространств $V=\mathcal{E}\otimes
L^2(d\sigma(\rho))$ является $U_q\mathfrak{sl}_2$-модулем:
 $$\xi: v\otimes f \mapsto \pi_{-\frac{1}{2}+i\rho}(\xi)\,v\otimes\,f,\qquad
  \xi \in U_q\mathfrak{sl}_2,\, v \in \mathcal{E},\, f\in L^2(d\sigma).$$
В силу следствия \ref{C_q-Laplace}, достаточно предъявить такой
изометрический линейный оператор $\mathcal{I}:\mathscr{D}(\mathbb{D})_q
\hookrightarrow V$, что $\mathcal{I}\, C_q\,f\,=\,C_q\, \mathcal{I}\,f$ при
всех $f \in \mathscr{D}(\mathbb{D})_q$.

Как было показано, см. предложение \ref{def-rel-D(D)},
$U_q\mathfrak{sl}_2$-модуль $\mathscr{D}(\mathbb{D})_q$ порождается
элементом $f_0$ и соотношения $K^{\pm 1}f_0=f_0$ являются определяющими.
Значит, существует и единствен морфизм $U_q\mathfrak{sl}_2$-модулей
$$
\mathcal{I}:\mathscr{D}(\mathbb{D})_q \hookrightarrow V,\qquad
\mathcal{I}:\,f_0 \mapsto \, e_0\otimes (q^{-2}-1).
$$
Остается доказать его изометричность, то есть равенство нулю инвариантной
эрмитовой формы
$
 \left<f_1,f_2\right>\,\stackrel{\rm def}{=}\,
 (\mathcal{I}f_1,\mathcal{I}f_2)-\int \limits_{\mathbb{D}_q}f_2^*f_1d \nu
$ на $\mathscr{D}(\mathbb{D})_q$.

Напомним соглашения об обозначениях: $l=-\frac12+i\rho$, $x=(1-zz^*)^{-1}$.
Из равенства $ \mathcal{I}\,C_q\,f\,=\,C_q\,\mathcal{I}\,f$, где
$f\in\mathscr{D}(\mathbb{D})_q$, и из результатов предыдущего
\itemiiiа следует, что для любого элемента вида $\psi(x)$ алгебры
$\mathscr{D}(\mathbb{D})_q$
$$
\mathcal{I}\psi\,=\,e_0\otimes\left(\int\limits_1^\infty
\psi(x)\Phi_{-\frac12+i\rho}(x)d_{q^{-2}}x\right),
$$
и для любых двух элементов вида $\psi_1(x)$, $\psi_2(x)$ алгебры
$\mathscr{D}(\mathbb{D})_q$
\begin{equation*}\langle\psi_1,\psi_2\rangle\,
=\,0.
\end{equation*}
Из последнего равенства и из ортогональности весовых подпространств
относительно инвариантной эрмитовой формы $\langle\cdot,\cdot\rangle$
следует, что каждый элемент вида $\psi(x)$ принадлежит ее ядру. В
частности, это ядро содержит образующую $f_0$\ рассматриваемого
$U_q\mathfrak{sl}_2$-модуля $\mathscr{D}(\mathbb{D})_q$. \hfill $\square$

\medskip

\begin{corollary}\label{boundness_sl_2}\ \
 $(1+q)^{-2}\leq\; \square_q\; \leq (1-q)^{-2}$.
\end{corollary}

\bigskip

Получим разложение $L^2(d \nu)_q$ в прямой интеграл гильбертовых
пространств \cite[стр. 450-451]{BrRob}). Пусть $\overline{V^{(l)}}$ --
пополнение предгильбертова пространства $V^{l}$, а
$\overline{V}$ -- пополнение предгильбертова пространства $V$, то есть\\
$\overline{V}=\newoplus\int\limits_{0}^{\frac{\pi}{2\log(q^{-1})}}
\overline{V^{(-1/2+i\rho)}}\,d\sigma(\rho)$.

\begin{proposition}\label{direct_int_sl_2}\ \  Изометрия $\mathcal{I}$
продолжается по непрерывности до унитарного оператора
$\overline{\mathcal{I}}$, отображающего $L^2(d \nu)_q$ на $\overline{V}$.
\end{proposition}

{\bf Доказательство.} Достаточно показать, что линейное многообразие
$\mathcal{I}\cdot\mathscr{D}(\mathbb{D})_q$ плотно в гильбертовом
пространстве $\overline{V}$. Заметим, что замыкание $L$ этого линейного
многообразия инвариантно относительно оператора умножения на непрерывную
функцию \eqref{l-rho}, вещественную на отрезке
$\left[0,\frac{\pi}{2\log(q^{-1})}\right]$ и разделяющую его точки. Значит,
$L$ -- общее инвариантное подпространство операторов умножения на функции
из $L^\infty(d\sigma)$. Следовательно, достаточно доказать сюръективность
морфизмов $U_q\mathfrak{sl}_2$-модулей
$$
\mathcal{I}_l: \mathscr{D}(\mathbb{D})_q\to V^{(l)},\qquad
\mathcal{I}_l:f_0 \mapsto e_0
$$
при $\mathrm{Re}\, l =-\frac{1}{2}$, см. \cite[стр. 452-453]{BrRob}.
Остается воспользоваться простотой $U_q\mathfrak{sl}_2$-модулей $V^{(l)}$
для нецелых $l$, см. лемму \ref{univ_pi_l}. \hfill $\square$

\subsubsection{От функций к сечениям линейных
расслоений.}\label{bundles_sl_2}

В классическом случае $q=1$ при построении комплексов Дольбо важную роль
играют пространства дифференцируемых сечений голоморфных векторных
расслоений \cite[стр. 83]{Wells}. Приведем примеры их $q$-аналогов.

Пусть $\lambda\in\mathbb{R}$ и $\mathbb{C}[z]_{q,\lambda}$ -- бимодуль над
алгеброй $\mathbb{C}[z]_q$ с образующей $v_\lambda$ и определяющим
соотношением
\begin{equation}\label{v_lambda}
z\,v_\lambda\,=\,q^{-\lambda}\,v_\lambda\,z.
\end{equation}

\begin{proposition}\label{hol_part_sl_2}
Существует и единственна такая структура $U_q\mathfrak{sl}_2$-модуля в
$\mathbb{C}[z]_{q,\lambda}$, что
\begin{equation}\label{U_q_sl_2-bundles}
 K^{\pm 1}v_\lambda\,=\,q^{\pm \lambda}v_\lambda,\qquad
 E\,v_\lambda\,=\,
 -q^{\frac{1}{2}}\,\frac{1-q^{2\lambda}}{1-q^2}\,zv_\lambda,
 \qquad F\,v_\lambda\,=\,0
\end{equation}
и $\mathbb{C}[z]_{q,\lambda}$ является $U_q\mathfrak{sl}_2$-модульным
$\mathbb{C}[z]_q$-бимодулем.
\end{proposition}

{\bf Доказательство.} Единственность очевидна. При доказательстве
существования можно ограничиться частным случаем $\lambda \in \mathbb{Z}$,
поскольку коэффициенты в \eqref{v_lambda}, \eqref{U_q_sl_2-bundles} являются
полиномами Лорана переменной $u=q^\lambda$ (подобные соображения
использовались в \itemiiiе \ref{sL_2_mod_algebras}).

Рассмотрим $U_q\mathfrak{sl}_2$-модульную алгебру $\mathbb{C}[t_1^{\pm
1},t_2^{\pm 1}]_q$. Вложение
$$
\mathbb{C}[z]_q\hookrightarrow\mathbb{C}[t_1^{\pm 1},t_2^{\pm 1}]_q,\qquad
z\mapsto t_2^{-1}t_1,
$$
наделяет $\mathbb{C}[t_1^{\pm 1},t_2^{\pm 1}]_q$ структурой
$U_q\mathfrak{sl}_2$-модульного $\mathbb{C}[z]_q$-бимодуля.

Равенство
$(t_2^{-1}t_1)\,t_2^{-\lambda}=q^{-\lambda}\,t_2^{-\lambda}\,(t_2^{-1}t_1)$
позволяет ввести в рассмотрение гомоморфизм $\mathbb{C}[z]_q$-бимодулей
$$
\mathbb{C}[z]_{q,\lambda}\to\mathbb{C}[t_1^{\pm 1},t_2^{\pm 1}]_q,\qquad
v_\lambda\mapsto t_2^{-\lambda}.
$$
Он инъективен, так как множество $\{t_1^jt_2^k\}_{j,k\in\mathbb{Z}}$
является базисом векторного пространства $\mathbb{C}[t_1^{\pm 1},t_2^{\pm
1}]_q$. Остается заметить, что
$$
K^{\pm 1}t_2^{-\lambda}\,=\,q^{\pm\lambda}t_2^{-\lambda},\quad
E\,t_2^{-\lambda}\,=\,-q^{\frac12}\,\frac{1-q^{2\lambda}}{1-q^2}\,
(t_2^{-1}t_1)t_2^{-\lambda},\quad F\,t_2^{-\lambda}\,=\,0. \eqno \square
$$

\medskip

Рассмотрим $U_q\mathfrak{sl}_2$-модуль ${\rm
Pol}(\mathbb{C})_{q,\lambda}\stackrel{\rm def}{=}
\mathbb{C}[z]_{q,\lambda}\otimes \mathbb{C}[\bar{z}]_q$. Наделим его
структурой $U_q\mathfrak{sl}_2$-модульного ${\rm
Pol}(\mathbb{C})_q$-бимодуля, применяя универсальную $R$-матрицу так же,
как в \itemiiiе \ref{Pol_C_q}. Полученный ${\rm
Pol}(\mathbb{C})_q$-бимодуль можно задать с помощью образующей $v_\lambda$
и определяющих соотношений
\begin{equation}\label{Pol_lambda}
 z\,v_\lambda\,=\, q^{-\lambda}\,v_\lambda\,z,\qquad
 z^*\,v_\lambda\,=\,q^\lambda\, v_\lambda\,z^*,
 \end{equation}
а действие $U_q\mathfrak{sl}_2$ описывается, по-прежнему, равенствами
\eqref{U_q_sl_2-bundles}.

Продолжение по непрерывности приводит к $U_q\mathfrak{sl}_2$-модульному
${\rm Pol}(\mathbb{C})_q\oplus\mathscr{D}(\mathbb{D})_q$-бимодулю, который
можно описать с помощью образующей $v_\lambda$, соотношений
\eqref{Pol_lambda} и дополнительного соотношения
$\psi(y)v_\lambda=v_\lambda\psi(y)$, где $\psi$ -- функция с конечным
носителем на множестве $q^{2\mathbb{Z}_+}$.

  Введем обозначение
$
 \mathscr{D}(\mathbb{D})_{q,\lambda}\,\stackrel{\rm def}{=}\,
 \mathscr{D}(\mathbb{D})_q\,v_\lambda\,=\,
 v_\lambda\, \mathscr{D}(\mathbb{D})_q
$
для полученного $U_q\mathfrak{sl}_2$-модульного
$\mathscr{D}(\mathbb{D})_q$-бимодуля.

Напомним, что $\mathscr{D}(\mathbb{D})_q'$ является
$\mathscr{D}(\mathbb{D})_q$-бимодулем и $y^\lambda \in
\mathscr{D}(\mathbb{D})_q'$. Воспользуемся понятием эрмитовой метрики,
введенным в \itemiiiе \ref{Hodges}.

\begin{proposition}\label{metric_for_bundle}
Полуторалинейное отображение
\begin{equation*}\label{metric_lambda}
{\mathbf h}_\lambda:\,\mathscr{D}(\mathbb{D})_{q,\lambda}\times
\mathscr{D}(\mathbb{D})_{q,\lambda}\to \mathscr{D}(\mathbb{D})_q,\;\;
{\mathbf h}_\lambda:\,
 v_\lambda\,f_1\,\times\,v_\lambda\,f_2 \mapsto f_2^*\,(1-zz^*)^\lambda\,f_1
\end{equation*}
является инвариантной положительной эрмитовой метрикой.
\end{proposition}

 Доказательство будет приведено в \itemiiiе \ref{finite_w_0SU_11}.

\medskip

\begin{corollary}\label{new-form-bundle}
Эрмитова форма
$$
(v_\lambda\,f_1\,,\,v_\lambda\,f_2)\,=\,\int\limits_{\mathbb{D}_q}\,
f_2^*\;(1-zz^*)^\lambda\;f_1\,d\nu,\qquad
f_1,f_2\in\mathscr{D}(\mathbb{D})_q,
$$
в $\mathscr{D}(\mathbb{D})_{q,\lambda}$ является положительной и
инвариантной.
\end{corollary}

\begin{remark}\label{old_form_bundle} Нетрудно доказать равенство
\begin{equation}\label{old-form-bundle}
(f_1\,v_\lambda\,,\,f_2\,v_\lambda)\,=\,
 \int\limits_{\mathbb{D}_q}\, f_2^*\; f_1\,(1-zz^*)^\lambda\;d\nu,\qquad
f_1,f_2 \in \mathscr{D}(\mathbb{D})_q.
\end{equation}
\end{remark}

\medskip  Пусть  $\lambda > 1$. Рассмотрим положительный линейный
функционал $\int\limits_{\mathbb{D}_q}\cdot d\nu_\lambda$ на $*$-алгебре
${\rm Pol}(\mathbb{C})_q$, определяемый равенством
$$ \int\limits_{\mathbb{D}_q}\,f\,d\nu_\lambda=
 \frac{1-q^{2(\lambda-1)}}{1-q^2}\,
 \int\limits_{\mathbb{D}_q}\,f\cdot(1-zz^*)^\lambda d\nu
 =(1-q^{2(\lambda-1)})\,{\rm tr}\,T_F(f\cdot(1-zz^*)^{\lambda-1}).$$
Числовой множитель перед интегралом выбран так, что
$\int\limits_{\mathbb{D}_q}\,1\,d\nu_\lambda\,=\,1$.

Пусть $L^2(d\nu_\lambda)_q$ -- пополнение ${\rm Pol}(\mathbb{C})_q$ по
норме $\|f\|_\lambda=\left(\int\limits_{\mathbb{D}_q}\,f^*f\,d\nu_\lambda
\right)^{\frac12}$ и $\square_{q,\lambda}$ -- инвариантный лапласиан в
$L^2(d\nu_\lambda)_q$, определяемый так же, как в рассмотренном выше
частном случае $\lambda=0$. Отметим, что $L^2(d\nu_\lambda)_q$ естественно
вложено в пространство обобщенных функций $\mathscr{D}(\mathbb{D})'_q$.

 Задача о разложении
 по собственным функциям радиальной части  оператора $\square_{q,\lambda}$
  в несколько  ином контексте решена в диссертации Гроенвельта \cite{Groenevelt}.


\subsubsection{Дополнение о $q$-специальных функциях.}
\label{q-series}\label{Askey-scheme}

Одно из основных свойств гипергеометрического ряда
$\sum\limits_{k=0}^\infty c_kz^k$ состоит в том, что отношение
$\dfrac{c_{k+1}}{c_k}$ является рациональной функцией переменной $k$
\cite[стр. 183]{Bateman1}. В случае базисных гипергеометрических рядов эти
отношения являются рациональными функциями переменной $q^k$, где
$q\in\mathbb{C}$, $q\ne 0$. Более точно, пусть
$a_1,a_2,\ldots,a_r;b_1,b_2,\ldots,b_s\in\mathbb{C}$. Базисный
гипергеометрический ряд \IND{базисный гипергеометрический ряд} определяется
равенством \cite[стр. 23]{GR}
\begin{multline}\label{phi_r_s}
\sideset{_r}{_s}{\mathop{\phi}}\left[\begin{array}{c}a_1,a_2,\ldots,a_r\\
b_1,b_2,\ldots,b_s\end{array};q,z\right]=
\\ =\sum_{k=0}^\infty
\frac{(a_1;q)_k(a_2;q)_k\cdots(a_r;q)_k}{(q;q)_k(b_1;q)_k\cdots(b_s;q)_k}
\left[(-1)^kq^{\frac{k(k-1)}2}\right]^{1+s-r}z^k,
\end{multline}
где
$$
(a;q)_k=
\begin{cases}
(1-a)(1-aq)\cdots(1-aq^{k-1}), & k\in\mathbb{N},\\ 1, & k=0.
\end{cases}
$$
Другое обозначение этого же ряда $\sideset{_r}{_s}{\mathop{\phi}}
(a_1,a_2,\ldots,a_r;b_1,b_2,\ldots,b_s;q,z)$.

В дальнейшем предполагается, что $q\in(0,1)$, $b_1,b_2,\ldots,b_s\notin
q^{-2\mathbb{Z}_+}$.

Ряд $\sideset{_r}{_s}{\mathop{\phi}}$ является полиномом переменной $z$
(обрывается), если хотя бы одно из чисел $a_1,a_2,\ldots,a_r$ принадлежит
множеству $q^{-2\mathbb{Z}_+}$. В противном случае возникает вопрос о
сходимости этого степенного ряда. В наиболее интересном случае $r=s+1$
радиус сходимости равен единице. \footnote{Хорошо известный \cite[стр.
183]{Bateman1} обобщенный
 гипергеометрический ряд $\sideset{_{s+1}}{_s}{\mathop{F}}$
получается из базисного гипергеометрического ряда
$\sideset{_{s+1}}{_s}{\mathop{\phi}}$ подстановкой $a_j=q^{\alpha_j}$,
$b_j=q^{\beta_j}$ и формальным предельным переходом $q\to 1$.}

\bigskip

Приведем некоторые сведения о $q$-специальных функциях, следуя \cite{GR}.
Введем обозначение
$$(a;q)_\infty=\prod\limits_{n=0}^\infty(1-aq^n).$$
Одним из простейших базисных гипергеометрических рядов является
$q$-биномиальный ряд \IND{ряд ! $q$-биномиальный}
$$
\sideset{_1}{_0}{\mathop{\phi}}(a;-;q,z)=
\sum_{n=0}^\infty\frac{(a;q)_n}{(q;q)_n}z^n.
$$
Нетрудно доказать \cite[стр. 26]{GR}, что при $|z|<1$
\begin{equation}\label{qbinom}
\sum_{n=0}^\infty\frac{(a;q)_n}{(q;q)_n}z^n=
\frac{(az;q)_\infty}{(z;q)_\infty}.
\end{equation}
В частности, при $a=0$ и $|z|<1$
\begin{equation}\label{qexp}
e_q(z)\overset{\mathrm{def}}{=}\sum_{n=0}^\infty\frac1{(q;q)_n}z^n=
\frac1{(z;q)_\infty}.
\end{equation}
Ряд $e_q(z)$ называют $q$-экспоненциальным, \IND{ряд ! $q$-экспоненциальный}
поскольку в топологии покоэффициентной сходимости
$$
\lim_{q\to 1}e_q((1-q)z)=\lim_{q\to
1}\sum_{n=0}^\infty\frac{(1-q)^n}{(q;q)_n}z^n=
\sum_{n=0}^\infty\frac{z^n}{n!}=e^z.
$$
Известен еще один $q$-аналог экспоненциального ряда \cite[стр. 285]{GR}
\begin{equation}\label{qExp}
E_q(z)\overset{\mathrm{def}}{=}
\sum_{n=0}^\infty\frac{q^{\frac{n(n-1)}{2}}}{(q;q)_n}z^n=(-z;q)_\infty.
\end{equation}

\bigskip

Приведем полученные Гейне \cite[стр. 29, 30]{GR} $q$-аналог формулы Эйлера
\IND{$q$-аналог ! формулы Эйлера}
\begin{equation}\label{Heine}
\sideset{_2}{_1}{\mathop{\phi}}(a,b;c;q,z)=
\frac{(b;q)_\infty(az;q)_\infty}{(c;q)_\infty(z;q)_\infty}
\sideset{_2}{_1}{\mathop{\phi}}\left(c/b,z;az;q,b\right)
\end{equation}
и $q$-аналог формулы Гаусса \IND{$q$-аналог ! формулы Гаусса}
\begin{equation}\label{Gauss}
\sideset{_2}{_1}{\mathop{\phi}}(a,b;c;q,c/(ab))=
\frac{(c/a;q)_\infty(c/b;q)_\infty}{(c;q)_\infty(c/(ab);q)_\infty}.
\end{equation}
В первой из них $|b|<1$, $|z|<1$, а во второй -- $|b|<1$, $|c/(ab)|<1$.

Интеграл Джексона \IND{интеграл Джексона} вводится равенством
$$\int\limits_0^af(x)d_qx=a(1-q)\sum_{k=0}^\infty f(aq^k)q^k.$$
По аналогии с $B$-функцией Эйлера \IND{$B$-функция Эйлера}
$$
B(x,y)=\int\limits_0^1t^{x-1}(1-t)^{y-1}dt,\qquad
\operatorname{Re}x>0,\;\operatorname{Re}y>0,
$$
определяется ее $q$-аналог \IND{$q$-аналог ! $B$-функции Эйлера}
$$
B_q(x,y)=\int\limits_0^1t^{x-1}\frac{(tq;q)_\infty}{(tq^y;q)_\infty}d_qt,
\qquad\operatorname{Re}x>0,\;y\notin-\mathbb{Z}_+.
$$
$q$-Аналог $\Gamma$-функции Эйлера \IND{$q$-аналог ! $\Gamma$-функции
Эйлера} $\Gamma(x)$ вводится равенством
$$
\Gamma_q(x)=\frac{(q;q)_\infty}{(q^x;q)_\infty}(1-q)^{1-x}.
$$
Можно доказать \cite[стр. 39]{GR}, что $\lim_{q\to
1}\Gamma_q(x)=\Gamma(x)$,
$$B_q(x,y)=\frac{\Gamma_q(x)\Gamma_q(y)}{\Gamma_q(x+y)}.$$
Формула \eqref{Heine} может быть записана следующим образом \cite[стр.
40]{GR}:
$$
\sideset{_2}{_1}{\mathop{\phi}}(q^a,q^b;q^c;q,z)=
\frac{\Gamma_q(c)}{\Gamma_q(b)\Gamma_q(c-b)}\int\limits_0^1t^{b-1}
\frac{(tzq^a;q)_\infty(tq;q)_\infty}{(tz;q)_\infty(tq^{c-b};q)_\infty}d_qt,
$$
что является $q$-аналогом формулы Эйлера --- интегрального представления
для гипергеометрической функции, см. \cite[стр. 72]{Bateman1}.

Из определения функции $\Gamma_q(x)$ следует, что
$$
\Gamma_q(n+1)=1(1+q)(1+q+q^2)\cdots(1+q+\cdots+q^{n-1}).
$$
Пусть
$$\binom{a}{b}_q=\frac{\Gamma_q(a+1)}{\Gamma_q(b+1)\,\Gamma_q(a-b+1)}.$$
При $n,k\in\mathbb{Z}_+$, $n\ge k$ получаем $q$-биномиальные коэффициенты
\begin{equation}\label{bin_coeff}
\binom{n}{k}_q=\frac{(q;q)_n}{(q;q)_k(q;q)_{n-k}},
\end{equation}
называемые также полиномами Гаусса. Они играют важную роль в теории
разбиений и в комбинаторике \cite{Andrews}.\footnote{Например, если
$\mathbb{F}$ -- конечное поле с числом элементов $q$, то векторное
пространство $\mathbb{F}^n$ состоит из $q^n$ элементов, а множество его
$k$-мерных подпространств -- из $\binom{n}{k}_q$ элементов \cite[стр.
218]{Andrews}.}

Нетрудно доказать, что в алгебре над полем $\mathbb{Q}(q)$, определяемой
двумя образующими $a,b$ и коммутационным соотношением $ab=qba$, имеет место
равенство
\begin{equation}\label{qNewton}
(a+b)^n=\sum_{k=0}^n\binom{n}{k}_qb^ka^{n-k}.
\end{equation}
Оно означает, что $(a+b)^n$ является производящей ''функцией''
$q$-биномиаль\-ных коэффициентов.

Используя \eqref{Heine}, нетрудно получить формулу суммирования
Пфаффа-Заальшютца \cite[стр. 286]{GR}: \IND{формула суммирования
Пфаффа-Заальшютца}
\begin{equation}\label{Pfaff}
\sideset{_3}{_2}{\mathop{\phi}}\left[\begin{array}{c}a,b,q^{-n}\\
c,abc^{-1}q^{1-n}\end{array};q,q\right]=
\frac{\left(\frac{c}{a};q\right)_n\left(\frac{c}{b};q\right)_n}
{(c;q)_n\left(\frac{c}{ab};q\right)_n},\qquad n\in\mathbb{Z}_+.
\end{equation}
В \cite{GR} собраны формулы суммирования, аналогичные \eqref{Pfaff}, и
формулы преобразования базисных гипергеометрических рядов, аналогичные
\eqref{Heine}. Например, \cite[стр. 45 и 290 -- 291]{GR}:
\begin{equation}\label{page45}
\sideset{_2}{_1}{\mathop{\phi}}\left[\begin{array}{c}q^{-n},b\\
c\end{array};q,z\right]=\frac{\left(\frac{c}{b};q\right)_n}{(c;q)_n}
\sideset{_3}{_2}{\mathop{\phi}}
\left[\begin{array}{c}q^{-n},b,\frac{bzq^{-n}}{c}\\
\frac{bq^{1-n}}{c},0\end{array};q,q\right],
\end{equation}
\begin{multline}\label{first_3_2}
\sideset{_3}{_2}{\mathop{\phi}}\left[\begin{array}{c}a,b,c\\
d,e\end{array};q,\frac{de}{abc}\right]=
\frac{\left(\frac{e}{a};q\right)_\infty\left(\frac{de}{bc};q\right)_\infty}
{(e;q)_\infty\left(\frac{de}{abc};q\right)_\infty}
\sideset{_3}{_2}{\mathop{\phi}}
\left[\begin{array}{c}a,\frac{d}{b},\frac{d}{c}\\
d,\frac{de}{bc}\end{array};q,\frac{e}{a}\right]=
\\ =\frac{(b;q)_\infty\left(\frac{de}{ab};q\right)_\infty
\left(\frac{de}{bc};q\right)_\infty}{(d;q)_\infty
(e;q)_\infty\left(\frac{de}{abc};q\right)_\infty}
\sideset{_3}{_2}{\mathop{\phi}}
\left[\begin{array}{c}\frac{d}{b},\frac{e}{b},\frac{de}{abc}\\
\frac{de}{ab},\frac{de}{bc}\end{array};q,b\right],
\end{multline}
\begin{multline}\label{second_3_2}
\sideset{_3}{_2}{\mathop{\phi}}\left[\begin{array}{c}q^{-n},b,c\\
d,e\end{array};q,q\right]=\frac{\left(\frac{de}{bc};q\right)_n}{(e;q)_n}
\left(\frac{bc}{d}\right)^n\sideset{_3}{_2}{\mathop{\phi}}
\left[\begin{array}{c}q^{-n},\frac{d}{b},\frac{d}{c}\\
d,\frac{de}{bc}\end{array};q,q\right]=
\\ =\frac{\left(\frac{e}{c};q\right)_n}{(e;q)_n}c^n
\sideset{_3}{_2}{\mathop{\phi}}\left[\begin{array}{c}q^{-n},c,\frac{d}{b}\\
d,\frac{cq^{1-n}}{e}\end{array};q,\frac{bq}{e}\right],
\end{multline}
\begin{equation}\label{third_3_2}
\sideset{_3}{_2}{\mathop{\phi}}\left[\begin{array}{c}q^{-n},b,c\\
d,e\end{array};q,\frac{deq^n}{bc}\right]=
\frac{\left(\frac{e}{c};q\right)_n}{(e;q)_n}\sideset{_3}{_2}{\mathop{\phi}}
\left[\begin{array}{c}q^{-n},c,\frac{d}{b}\\
d,\frac{cq^{1-n}}{e}\end{array};q,q\right].
\end{equation}

\bigskip

Введем $q$-аналоги операторов дифференцирования \IND{$q$-аналог ! операторов
дифференцирования}
\begin{equation}\label{D_pm}
D_q^-f(z)=\frac{f(z)-f(qz)}{z-qz},\qquad
D_q^+f(z)=\frac{f(q^{-1}z)-f(z)}{z-qz}.
\end{equation}

Имеет место следующий $q$-аналог формулы интегрирования по частям:
\IND{$q$-аналог ! формулы интегрирования по частям}
\begin{multline}\label{by_part}
\int\limits_0^a(D_q^-u)(x)v(x)d_qx=
\\ =u(a)v(q^{-1}a)-u(0)v(0)-\int\limits_0^au(x)(D_q^+v)(x)d_qx,
\end{multline}
где $a\in\mathbb{C}$ и функции $u,v$ непрерывны на компакте
$\mathfrak{M}_a=\{0\}\cup(aq^{-1})q^{\mathbb{Z}_+}$.

 Из определений следует, что формальный ряд $e_q((1-q)z)$ является
решением $q$-разностного уравнения $D_q^-y=y$, равным $1$ в точке $z=0$.

Формальный ряд $\sideset{_2}{_1}{\mathop{\phi}}(a,b;q,z)$ удовлетворяет
$q$-разностному уравнению \cite[стр. 44]{GR}
\begin{multline}\label{q_diff}
z(c-abqz)(D_q^-)^2y(z)+\left(\frac{1-c}{1-q}+
\frac{(1-a)(1-b)-(1-abq)}{1-q}z\right)\,D_q^-y(z)-
\\ -\frac{(1-a)(1-b)}{(1-q)^2}y(z)=0.
\end{multline}
Замена $a,b,c$ на $q^a, q^b, q^c$ и формальный предельный переход $q\to 1$
приводят к дифференциальному уравнению
$$z(1-z)\,\frac{d^2}{dz^2}y(z)+(c-(a+b+1)z)\,\frac{dy(z)}{dz}-ab\,y(z)=0.$$

Решению $q$-разностных уравнений с помощью базисных гипергеометрических
рядов и $q$-аналогам задачи Штурма-Лиувилля посвящены главы 4 -- 5 книги
\cite{Exton}.

\bigskip

Перейдем к $q$-гипергеометрическим полиномам одной переменной.

В середине 80х годов Эски и Вильсон в работе \cite{AskeyMAMS} ввели
четырехпараметрическое семейство полиномов, обладающих многими
замечательными свойствами, -- полиномы Эски-Вильсона. \IND{полиномы !
Эски-Вильсона} Они определяются с помощью алгебры ''инвариантных''
полиномов Лорана
$$
\mathbb{C}[t,t^{-1}]^\mathrm{inv}=
\{f(t)\in\mathbb{C}[t,t^{-1}]|\:f(t^{-1})=f(t)\}
$$
и канонического изоморфизма
\begin{equation}\label{to_tore}
\mathbb{C}[x]\overset{\simeq}{\to}\mathbb{C}[t,t^{-1}]^\mathrm{inv},\qquad
f(x)\mapsto f\left(\frac{t+t^{-1}}2\right).
\end{equation}
Именно \cite[стр. 212]{GR},
\begin{multline}\label{AW}
p_n(x;a,b,c,d|q)=
\\ =\frac{(ab;q)_n(ac;q)_n(ad;q)_n}{a^n}
\sideset{_4}{_3}{\mathop{\phi}}
\left[\begin{array}{c}q^{-n},abcdq^{n-1},at,at^{-1}\\
ab,ac,ad\end{array};q,q\right].
\end{multline}
Полиномы Эски-Вильсона являются решениями $q$-разностного уравнения второго
порядка \cite[стр. 225]{GR}
\begin{multline*}
D_q(\sqrt{1-x^2}w(x;aq^\frac12,bq^\frac12,cq^\frac12,dq^\frac12|q)D_qp_n)+
\\ +\lambda_n\sqrt{1-x^2}w(x;a,b,c,d|q)p_n=0,
\end{multline*}
где $\lambda_n=4q^{-n+1}(1-q^n)(1-abcdq^{n-1})(1-q)^{-2}$, для
$w(x;a,b,c,d|q)$ имеется явная формула \cite[стр. 180]{GR}, а $D_q$
определяется с помощью канонического изоморфизма \eqref{to_tore}:
$$
D_q:\mathbb{C}[t,t^{-1}]^\mathrm{inv}\to\mathbb{C}[t,t^{-1}]^\mathrm{inv},
\qquad D_q:f(t)\mapsto\frac{f(q^\frac12 t)-f(q^{-\frac12}t)}{q^\frac12
t-q^{-\frac12}t}.
$$

Отметим, что переход от полиномов переменной $x$ к инвариантным полиномам
Лорана переменной $t$ является шагом по направлению к теории симметрических
ортогональных многочленов нескольких переменных, построенной в работах
И.~Макдональда \cite[глава VI]{Mac95} и Т.~Коорнвиндера \cite{KoornBC}.
Введенные в \cite{KoornBC} полиномы являются обобщениями полиномов
Эски-Вильсона.

Следуя обзору \cite{Koekoek}, приведем известные результаты о двух
семействах ортогональных полиномов, являющихся частными либо предельными
случаями полиномов Эски-Вильсона.

 Простейшими представителями первого семейства
являются полиномы Аль-Салама-Чихары, а второго -- малые $q$-полиномы Якоби.
Первое семейство связано с теорией представлений $*$-алгебры Хопфа
$U_q\mathfrak{su}_{1,1}$, а второе -- с теорией представлений
\hbox{$*$-алгебры} Хопфа $U_q\mathfrak{su}_2$.

В дальнейшем используются обозначения $$t=e^{i\theta}, \qquad
x=\cos\theta,$$
$$(a_1,a_2,\ldots,a_m;q)_k=(a_1;q)_k(a_2;q)_k\cdots(a_m;q)_k,
\qquad k \in \mathbb{N} \cup \{\infty\}.$$

Полиномами Аль-Салама-Чихары \IND{полиномы ! Аль-Салама-Чихары} называют
\cite[стр. 25]{AskeyMAMS} полиномы степени $n$ переменной $x=\cos\theta$,
зависящие от параметров $a$, $b$, $q$ и определяемые равенством
\begin{equation}\label{bible_3.8.1}
Q_n(x;a,b|q)=\frac{(ab;q)_n}{a^n}\sideset{_3}{_2}{\mathop{\phi}}
\left[\begin{array}{c}q^{-n},ae^{i\theta},ae^{-i\theta}\\
ab,0\end{array};q,q\right].
\end{equation}
Если $a$, $b$ либо вещественные числа, либо комплексно сопряженные числа, и
$|a|<1$, $|b|<1$, то полиномы $Q_n(x;a,b|q)$ удовлетворяют следующим
соотношениям ортогональности:
\begin{equation}\label{bible_3.8.2}
\frac1{2\pi}\int\limits_{-1}^1Q_m(x;a,b|q)
Q_n(x;a,b|q)\frac{w(x;a,b|q)dx}{\sqrt{1-x^2}}=
\frac{\delta_{mn}}{(q^{n+1},abq^n;q)_\infty},
\end{equation}
где
\begin{multline}\label{Al_w}
w(x;a,b|q)=\left|\frac{(e^{2i\theta};q)_\infty}
{(ae^{i\theta},be^{i\theta};q)_\infty}\right|^2=
\\ =\frac{h(x,1)h(x,-1)h(x,q^\frac12)h(x,-q^{\frac12})}{h(x,a)h(x,b)},
\end{multline}
$$
h(x,\alpha)=(\alpha e^{i\theta},\alpha
e^{-i\theta};q)_\infty=\prod_{k=0}^\infty\left(1-2\alpha
xq^k+\alpha^2q^{2k}\right).
$$
Имеются рекуррентные соотношения для этих полиномов:
\begin{equation}\label{bible_3.8.4}
2xQ_n(x)=Q_{n+1}(x)+(a+b)q^nQ_n(x)+(1-q^n)(1-abq^{n-1})Q_{n-1}(x).
\end{equation}
Они являются решениями $q$-разностного уравнения
\begin{multline}\label{bible_3.8.6}
(1-q)^2D_q\left[\widetilde{w}\left(\left.x;aq^{\frac12},bq^{\frac12}\right|
q\right)D_qy(x)\right]+
\\ +4q^{-n+1}(1-q^n)\widetilde{w}(x;a,b|q)y(x)=0,
\end{multline}
где $\widetilde{w}(x;a,b|q)=w(x;a,b|q)(1-x^2)^{-\frac{1}{2}}$. Производящая
функция для полиномов Аль-Салама-Чихары имеет следующий вид:
\begin{equation}\label{bible_3.8.13}
\frac{(at,bt;q)_\infty}{(e^{i\theta}t,e^{-i\theta}t;q)_\infty}=
\sum_{n=0}^\infty\frac{Q_n(x;a,b|q)}{(q;q)_n}t^n,\qquad x=\cos\theta.
\end{equation}
Полиномы Аль-Салама-Чихары изучались многими авторами. Отметим работы
\cite{AskeyIsmailMAMS}, \cite{KoelinkGraf}. В последней из них получен
$q$-аналог формулы сложения Графа для функций Бесселя \cite[стр.
53]{Bateman2}.

Полиномы Аль-Салама-Чихары являются частными случаями непрерывных дуальных
$q$-полиномов Хана \cite{Gupta}
\begin{equation}\label{bible_3.3.1}
p_n(x;a,b,c|q)=\frac{(ab,ac;q)_n}{a^n}\sideset{_3}{_2}{\mathop{\phi}}
\left[\begin{array}{c}q^{-n},ae^{i\theta},ae^{-i\theta}\\
ab,ac\end{array};q,q\right],
\end{equation}
где $x=\cos\theta$ (для перехода к полиномам Аль-Салама-Чихары достаточно
положить $c$ равным $0$). Если все три числа $a$, $b$, $c$ вещественны,
либо одно из них вещественно, а два других комплексно сопряжены и $|a|<1$,
$|b|<1$, $|c|<1$, то имеют место соотношения ортогональности, обобщающие
\eqref{bible_3.8.2}:
\begin{equation}\label{bible_3.3.2}
\frac1{2\pi}\int\limits_{-1}^1p_m(x;a,b,c|q)
p_n(x;a,b,c|q)\frac{w(x;a,b,c|q)dx}{\sqrt{1-x^2}}=h_n\delta_{mn},
\end{equation}
где
\begin{multline}
w(x;a,b,c|q)=\left|\frac{(e^{2i\theta};q)_\infty}
{(ae^{i\theta},be^{i\theta},ce^{i\theta};q)_\infty}\right|^2=
\\ =\frac{h(x,1)h(x,-1)h(x,q^\frac12)h(x,-q^{\frac12})}{h(x,a)h(x,b)h(x,c)},
\end{multline}
$$h_n=\frac1{(q^{n+1},abq^n,acq^n,bcq^n;q)_\infty}.$$

Имеются рекуррентные соотношения
\begin{equation}\label{bible_3.3.4}
2x\widetilde{p}_n(x)=
A_n\widetilde{p}_{n+1}(x)+[a+a^{-1}-(A_n+C_n)]\widetilde{p}_n(x)+
C_n\widetilde{p}_{n-1}(x),
\end{equation}
где
$$
\widetilde{p}_n(x)=\frac{a^np_n(x;a,b,c|q)}{(ab,ac;q)_n}=
\sideset{_3}{_2}{\mathop{\phi}}
\left[\begin{array}{c}q^{-n},ae^{i\theta},ae^{-i\theta}\\
ab,ac\end{array};q,q\right],
$$
$$A_n=a^{-1}(1-abq^n)(1-acq^n),\qquad C_n=a(1-q^n)(1-bcq^{n-1}).$$
Полиномы $p_n$ являются решениями следующего $q$-разностного уравнения:
\begin{multline}\label{bible_3.3.6}
(1-q)^2D_q\left[\widetilde{w}\left(\left.x;aq^{\frac12},bq^{\frac12},
cq^{\frac12}\right|q\right)D_qy(x)\right]+
\\ +4q^{-n+1}(1-q^n)\widetilde{w}(x;a,b,c|q)y(x)=0,
\end{multline}
где $\widetilde{w}(x;a,b,c|q) = w(x;a,b,c|q)(1-x^2)^{-\frac{1}{2}}$.
Имеются производящие функции
\begin{flalign}
\frac{(ct;q)_\infty}{(e^{i\theta}t;q)_\infty}\sideset{_2}{_1}{\mathop{\phi}}
\left[\begin{array}{c}ae^{i\theta},be^{i\theta}\\
ab\end{array};q,e^{-i\theta}t\right] &=
\sum_{n=0}^\infty\frac{p_n(x;a,b,c|q)}{(ab,q;q)_n}t^n,
&x=\cos\theta,\label{bible_3.3.13}
\\ \frac{(bt;q)_\infty}{(e^{i\theta}t;q)_\infty}
\sideset{_2}{_1}{\mathop{\phi}}
\left[\begin{array}{c}ae^{i\theta},ce^{i\theta}\\
ac\end{array};q,e^{-i\theta}t\right] &=
\sum_{n=0}^\infty\frac{p_n(x;a,b,c|q)}{(ac,q;q)_n}t^n,
&x=\cos\theta,\label{bible_3.3.14}
\\ \frac{(at;q)_\infty}{(e^{i\theta}t;q)_\infty}
\sideset{_2}{_1}{\mathop{\phi}}
\left[\begin{array}{c}be^{i\theta},ce^{i\theta}\\
bc\end{array};q,e^{-i\theta}t\right] &=
\sum_{n=0}^\infty\frac{p_n(x;a,b,c|q)}{(bc,q;q)_n}t^n, &
x=\cos\theta.\label{bible_3.3.15}
\end{flalign}

\medskip Известны обобщения приведенных выше результатов на случай
полиномов Эски-Вильсона (см. \eqref{AW}). Соотношения ортогональности
справедливы при более общих предположениях, чем обсуждалось выше. В общем
случае мера включает как абсолютно непрерывную, так и дискретную
составляющую \cite[стр. 216]{GR}.

\bigskip Перейдем ко второму из интересующих нас семейств
$q$-гипергеомет\-ри\-ческих ортогональных полиномов.

Малые $q$-полиномы Якоби \IND{малые $q$-полиномы Якоби} были введены В.
Ханом \cite[стр. 29]{Hahn}.
\begin{equation}\label{bible_3.12.12}
p_n(x;a,b|q)=\sideset{_2}{_1}{\mathop{\phi}}
\left[\begin{array}{c}q^{-n},abq^{n+1}\\ aq\end{array};q,qx\right].
\end{equation}
Соотношения ортогональности для них имеют вид \cite[стр. 110]{AnAs}.
\begin{multline}\label{bible_3.12.2}
\sum_{k=0}^\infty\frac{(bq;q)_k}{(q;q)_k}(aq)^kp_m(q^k;a,b|q)p_n(q^k;a,b|q)=
\\ =\frac{(abq^2;q)_\infty}{(aq;q)_\infty}
\frac{(1-abq)(aq)^n}{(1-abq^{2n+1})}\frac{(q,bq;q)_n}{(aq,abq;q)_n}
\delta_{mn},
\end{multline}
где \ \ $0<a<q^{-1}$, $b<q^{-1}$, а рекуррентные соотношения --
\begin{multline}\label{bible_3.12.3}
-xp_n(x;a,b|q)=
\\ =A_np_{n+1}(x;a,b|q)-(A_n+C_n)p_n(x;a,b|q)+C_np_{n-1}(x;a,b|q),
\end{multline}
где
$$
A_n=q^n\frac{(1-aq^{n+1})(1-abq^{n+1})}{(1-abq^{2n+1})(1-abq^{2n+2})},\qquad
C_n=aq^n\frac{(1-q^n)(1-bq^n)}{(1-abq^{2n})(1-abq^{2n+1})}.
$$
Малые $q$-полиномы Якоби являются решениями $q$-разностного уравнения,
вытекающего из \eqref{q_diff}, и равносильного уравнению
\begin{multline}\label{bible_3.12.5}
q^{-n}(1-q^n)(1-abq^{n+1})x\,y(x)= \\
=B(x)y(qx)-[B(x)+D(x)]y(x)+D(x)y(q^{-1}x),
\end{multline}
где
$$B(x)=a(bqx-1),\qquad D(x)=x-1.$$
Производящая функция имеет вид
\begin{multline}\label{bible_3.12.11}
\sideset{_0}{_1}{\mathop{\phi}}\left[\begin{array}{c}-\\
aq\end{array};q,aqxt\right]\sideset{_2}{_1}{\mathop{\phi}}
\left[\begin{array}{c}x^{-1},0\\ bq\end{array};q,xt\right]=
\\ =\sum_{n=0}^\infty\frac{(-1)^nq^{\frac{n(n-1)}{2}}}{(bq,q;q)_n}
p_n(x;a,b|q)t^n.
\end{multline}

\bigskip

Если носитель меры $d\mu$ является конечным множеством и состоит из $N+1$
точки, то размерность векторного пространства $L^2(d\mu)$ равна $N+1$. В
этом случае ортогональный базис получают ортогонализацией степеней $1,x,
\ldots,x^N$. Важным примером таких ортогональных полиномов служат введенные
в работе \cite{Hahn} $q$-полиномы Хана. \IND{$q$-полиномы Хана} Они
определяются равенством \cite[стр. 204]{GR}
\begin{multline}\label{bible_3.6.1}
Q_n(q^{-x};\alpha,\beta,N|q)=
\\ =\sideset{_3}{_2}{\mathop{\phi}}
\left[\begin{array}{c}q^{-n},\alpha\beta q^{n+1},q^{-x}\\
\alpha q,q^{-N}\end{array};q,q\right],\qquad n=0,1,\ldots,N
\end{multline}
и являются полиномами переменной $q^{-x}$. Малые $q$-полиномы Якоби могут
быть получены из $q$-полиномов Хана с помощью предельного перехода:
\begin{equation}\label{bible_4.6.1}
\lim_{N\to\infty}Q_n(q^{x-N};\alpha,\beta,N|q)=p_n(q^x;\alpha,\beta|q).
\end{equation}
Соотношения ортогональности
\begin{multline}\label{bible_3.6.2}
\sum_{x=0}^N\frac{(aq,q^{-N};q)_x}{(q,\beta^{-1}q^{-N};q)_x}(\alpha\beta
q)^{-x}Q_m(q^{-x};\alpha,\beta,N|q)Q_n(q^{-x};\alpha,\beta,N|q)=
\\ =\frac{(\alpha\beta q^2;q)_N}{(\beta q;q)_N(\alpha q)^N}
\frac{(q,\alpha\beta q^{N+2},\beta q;q)_n}{(\alpha q,\alpha\beta
q,q^{-N};q)_n}\frac{(1-\alpha\beta q)(-\alpha q)^n}{(1-\alpha\beta
q^{2n+1})}q^{\frac{n(n-1)}{2}-Nn}\delta_{mn}
\end{multline}
имеют место при $0<\alpha<q^{-1}$ и $0<\beta<q^{-1}$, или $\alpha>q^{-N}$ и
$\beta>q^{-N}$. Рекуррентные соотношения имеют вид
\begin{multline}\label{bible_3.6.3}
-(1-q^{-x})Q_n(q^{-x})=
\\ =A_nQ_{n+1}(q^{-x})-(A_n+C_n)Q_n(q^{-x})+C_nQ_{n-1}(q^{-x}),
\end{multline}
где
$$Q_n(q^{-x})=Q_n(q^{-x};\alpha,\beta,N|q),$$
$$
A_n=\frac{(1-q^{n-N})(1-\alpha q^{n+1})(1-\alpha\beta
q^{n+1})}{(1-\alpha\beta q^{2n+1})(1-\alpha\beta q^{2n+2})},
$$
$$
C_n=-\frac{\alpha q^{n-N}(1-q^n)(1-\alpha\beta q^{n+N+1})(1-\beta
q^n)}{(1-\alpha\beta q^{2n})(1-\alpha\beta q^{2n+1})}.
$$
$q$-Полиномы Хана являются решениями $q$-разностного уравнения
\begin{multline}\label{bible_3.6.5}
q^{-n}(1-q^n)(1-\alpha\beta q^{n+1})y(x)=
\\ =B(x)y(x+1)-[B(x)+D(x)]y(x)+D(x)y(x-1),
\end{multline}
где $y(x)=Q_n(q^{-x};\alpha,\beta,N|q)$,
$$
B(x)=(1-q^{x-N})(1-\alpha q^{x+1}),\qquad D(x)=\alpha
q(1-q^x)(\beta-q^{x-N-1}).
$$

При всех $x=0,1,2,\ldots,N$
\begin{multline}\label{bible_3.6.11}
\sideset{_1}{_1}{\mathop{\phi}}\left[\begin{array}{c}q^{-x}\\
\alpha q\end{array};q,\alpha qt\right]\sideset{_2}{_1}{\mathop{\phi}}
\left[\begin{array}{c}q^{x-N},0\\ \beta q\end{array};q,q^{-x}t\right]=
\\ =\sum_{n=0}^N\frac{(q^{-N};q)_n}{(\beta
q,q;q)_n}Q_n(q^{-x};\alpha,\beta,N|q)t^n.
\end{multline}

\medskip

Приведенные выше результаты о $q$-полиномах Хана обобщаются на случай
$q$-полиномов Рака \cite[стр. 202]{GR} -- полиномов, образы которых при
гомоморфизме алгебр
$$
\mathbb{C}[x]\to\mathbb{C}[u,u^{-1}],\qquad f(x)\mapsto
f(u^{-1}+\gamma\delta qu)
$$
равны
$$
\sideset{_4}{_3}{\mathop{\phi}}
\left[\begin{array}{c}q^{-n},\alpha\beta q^{n+1},u^{-1},\gamma\delta qu\\
\alpha q,\beta\delta q,\gamma q\end{array};q,q\right],\qquad
n=0,1,2,\ldots,N.
$$
$q$-Полиномы Рака связаны с полиномами Эски-Вильсона \eqref{AW} подстановкой
\begin{equation}\label{bible_after_3.2.14}
\alpha=abq^{-1},\quad\beta=cdq^{-1},\quad\gamma=adq^{-1},\quad\delta=ad^{-1},
\quad u=a^{-1}e^{-i\theta}.
\end{equation}

\subsection{Интегральные  представления}

\subsubsection{Алгебры ядер интегральных
операторов.}\label{integral_operators}

В этом \itemiiе рассматриваются $q$-аналоги непрерывных линейных
операторов, играющих важную роль в гармоническом анализе, либо в
комплексном анализе. Как следует из теоремы Лорана Шварца о ядре \cite[стр.
156]{Hermander1}, такие операторы являются интегральными, а их ядра
являются обобщенными функциями. Наша типичная задача-- зная $q$-аналог
непрерывного линейного оператора, найти $q$-аналог ядра этого интегрального
оператора.

В настоящем \itemiiiе вводятся вспомогательные понятия. Прежде всего
покажем, что $A$-ин\-вариантный интеграл может быть использован при
построении морфизмов $A$-модулей.

Рассмотрим $A$-модульные алгебры $F_1$, $F_2$. Каждому линейному функционалу
$\nu:F_2\to\mathbb{C}$ и каждому элементу $\mathscr{K}\in F_1\otimes F_2$
сопоставим линейный интегральный оператор
$$
K:F_2\to F_1,\qquad
K:f\mapsto(\operatorname{id}\otimes\nu)(\mathscr{K}(1\otimes f)).
$$
Элемент $\mathscr{K}$ называют ядром этого интегрального оператора
\IND{ядро интегрального оператора} и используют следующую сокращенную
форму записи:
$$Kf=\int\mathscr{K}(1\otimes f)\,d\nu=\int\mathscr{K}fd\nu.$$

\begin{proposition}\label{intertw_kernels}
Рассмотрим $A$-модульные алгебры $F_1, F_2$ и $A$-ин\-вариантный интеграл
$\nu: F_2\to \mathbb{C}$. Предположим, что билинейная форма
\begin{equation}\label{bilinear_new}
f'\times f''\mapsto\nu(f'\cdot f''),\qquad f',f''\in F_2,
\end{equation}
на $F_2$ невырождена. Интегральный оператор с ядром $\mathscr{K}$ является
морфизмом $A$-модулей, если и только если это ядро $A$-инвариантно.
\end{proposition}

{\bf Доказательство.} 1. Если ядро $\mathscr{K}$ инвариантно, то линейное
отображение
$$
F_2 \rightarrow F_1 \otimes F_2, \qquad f \mapsto \mathscr{K} (1 \otimes f),
$$
а, следовательно, и интегральный оператор с ядром $\mathscr{K}$ являются
морфизмами $A$-модулей.

  2. Предположим, что интегральный оператор с ядром   $\mathscr{K}$ является
морфизмом $A$-модулей. Тогда, как следует из \eqref{int_by_parts} и из
невырожденности билинейной формы \eqref{bilinear_new}, это ядро
удовлетворяет системе уравнений
\begin{equation} \label{kernel_equation}
(a \otimes 1)\, \mathscr{K}=(1 \otimes S^{-1}(a))\, \mathscr{K}, \qquad a
\in A,
\end{equation}
из которой вытекает его инвариантность, поскольку
 $(\mathrm{id}\otimes S^{-1})\Delta^{\mathrm{op}}=\mathrm{id}$.\hfill $\square$

\bigskip
Изучение инвариантных ядер и отвечающих им интегральных операторов является
одной из основных задач некоммутативного гармонического анализа \cite{GKT}.
При построении $A$-инвариантных ядер будет существенно использоваться
следующее

\begin{proposition}\label{IntKer}
Пусть $F_1^\mathrm{op}$ -- алгебра, получаемая из $F_1$ заменой умножения на
противоположное. $A$-инвариантные ядра образуют подалгебру алгебры
$F_1^\mathrm{op}\otimes F_2$.
\end{proposition}

{\bf Доказательство.} Достаточно показать, что произведение
$A$-ин\-вариантных элементов алгебры ядер является $A$-инвариантным ядром.
Сопоставим ядрам
$$
\mathscr{K}'=\sum\limits_ia_i'\otimes
b_i',\qquad\mathscr{K}''=\sum\limits_ja_j''\otimes b_j''
$$
 следующие линейные операторы в $F_1\otimes F_2$:
\begin{equation}\label{op}
\mathcal{K}': f_1\otimes f_2\mapsto\sum\limits_if_1a_i'\otimes
b_i'f_2,\quad \mathcal{K}'': f_1\otimes f_2\mapsto\sum\limits_j
f_1a_j''\otimes b_j''f_2.
\end{equation}
 Из $A$-инвариантности ядер $\mathscr{K}'$,
$\mathscr{K}''$ следует, что операторы \eqref{op}, а, значит, и их
произведение $\mathcal{K}'\mathcal{K}''$, отображающее
 $ f_1\otimes
f_2$ в $\sum\limits_i\sum\limits_jf_1a_j''a_i'\otimes b_i'b_j''f_2, $
являются эндоморфизмами $A$-модуля $F_1\otimes F_2$. Применяя этот
эндоморфизм к $A$-инварианту $1\otimes 1$, получаем $A$-инвариант, равный
произведению $\mathscr{K}''\mathscr{K}'$ в алгебре $F_1^\mathrm{op}\otimes
F_2$. \hfill $\square$

\bigskip В оставшейся части \itemiiiа рассматриваются
   $*$-алгебра Хопфа $A$, $A$-модульные $*$-алгебры $F_1, F_2$,
  вещественный инвариантный интеграл  $\nu: F_2 \rightarrow \mathbb{C}$,
  и предполагается  невырожденность  билинейной формы
  $$
  F_2 \times F_2 \rightarrow \mathbb{C}, \qquad f' \times f'' \mapsto
  \nu ( f'\; f'').
  $$

Легко доказать единственность такой инволюции $*$ в алгебре ядер
$F_1^{op}\otimes F_2$, что при всех $a\in A$, $f\in F_2$, $\mathscr{K}\in
F_1^{op}\otimes F_2$
\begin{equation}\label{real_invol}
\int\mathscr{K}^*(1\otimes f)d\nu=\left(\int\mathscr{K}(1\otimes
f^*)d\nu\right)^*.
\end{equation}
В классическом случае $q=1$ равенство \eqref{real_invol} означает
вещественность интегрального оператора с вещественным ядром.

Инволюция со свойством \eqref{real_invol} допускает сужение на подалгебру
инвариантных ядер, что нетрудно показать, используя предложение
\ref{intertw_kernels}. Действительно, если $\mathscr{K}\in F_1^{op}\otimes
F_2$ -- инвариантное ядро, то при всех $a\in A$, $f\in F_2$ имеют место
равенства
$$
\int\mathscr{K}^*(1\otimes
af)d\nu=\left(\int\mathscr{K}(1\otimes(af)^*)d\nu\right)^*=
\left(\int\mathscr{K}(1\otimes S(a)^*f^*)d\nu\right)^*=
$$
$$=\left(S(a)^*\int\mathscr{K}(1\otimes f^*)d\nu\right)^*=
a\left(\int\mathscr{K}(1\otimes
f^*)d\nu\right)^*=a\int\mathscr{K}^*(1\otimes f)d\nu.
$$
Значит, если интегральный оператор с ядром $\mathscr{K}$ является морфизмом
$A$-модулей, то этим свойством обладает и интегральный оператор с ядром
$\mathscr{K}^*$.

   Инволюцию, удовлетворяющую требованию \eqref{real_invol},
  часто удается найти, используя явный вид инвариантного интеграла. Например,
  если алгебра $F_2$ обладает
  конечномерным точным $*$-представлением $T$ и для некоторого  обратимого
  самосопряженного элемента $Q \in F_2$  имеет место равенство
  \begin{equation}\label{explicit_star_tr}\nu(f)= \operatorname{tr}(T(f)Q),
  \qquad f \in F_2,
  \end{equation}
  то, как следует из определений,
  $\mathscr{K}^*=(1 \otimes Q)^{-1} \;
  \mathscr{K}^{*\otimes *}\; (1 \otimes Q)$-- требуемая инволюция.

  Отметим, что полученная формула верна при
  весьма общих предположениях, поскольку  конечномерность представления
  $T$ использована в ее доказательстве  лишь для
  перестановки сомножителей под знаком следа.

\subsubsection{Функция Грина для оператора
$\mathbf{\square_q}$.}\label{Green_function}

Начнем с классического случая $q=1$. Для широкого класса функций $f$
решение уравнения Пуассона $\square u=f$ в единичном круге $\mathbb{D}$
имеет вид
$$u(z)=\int\limits_{\mathbb{D}}\, G(z,\zeta)\,f(\zeta)\,d\nu(\zeta),$$
где функция Грина $G$ равна
\begin{equation}\label{1.2.8}
G(z,\zeta)=-{\rm ln}(|z-\zeta|^2)\,+\,{\rm ln}(|1-z \overline{\zeta}|^2).
\end{equation}

Равенство \eqref{1.2.8} доказывается методом отражений: первое слагаемое --
вклад основного источника, а второе -- вклад источника, полученного
отражением относительно единичной окружности, то есть инверсией \cite[стр.
429]{Vladimirov}.

\bigskip

Для перехода к квантовому кругу нужны квантовые аналоги интегральных
операторов и их ядер.

Рассмотрим алгебру ${\rm Pol}({\mathbb C})_q^{\rm op}$, получающуюся из
${\rm Pol}({\mathbb C})_q$ заменой умножения на противоположное. В этой
алгебре $zz^*=q^2z^*z+1-q^2$. Стартуя с ${\rm Pol}({\mathbb C})_q^{\rm op}$,
повторим все те рассуждения, которые использовались в \itemiiе
\ref{finite_sl_2} при построении алгебры $\mathscr{D}(\mathbb{D})_q$
финитных функций в квантовом круге и $\mathscr{D}(\mathbb{D})_q$-бимодуля
$\mathscr{D}(\mathbb{D})^\prime_q$ обобщенных функций в квантовом круге. В
итоге получим алгебру $\mathscr{D}(\mathbb{D})_q^{\rm op}$ и
$\mathscr{D}(\mathbb{D})_q^{\rm op}$-бимодуль
$\mathscr{D}(\mathbb{D})_q^{\prime\, {\rm op}}$.

 Введем обозначение ${\rm Pol}({\mathbb
C}\times{\mathbb C})_q$ для алгебры ${\rm Pol}({\mathbb C})_q^{\rm
op}\otimes{\rm Pol}({\mathbb C})_q$ и обозначение $z,z^*,\zeta,\zeta^*$ для
ее образующих $z\otimes 1$, $z^*\otimes 1$, $1\otimes z$, $1\otimes z^*$.
Будем использовать обозначения $y,\,\eta$ для элементов
$$(1-z^*z)\otimes 1,\qquad \,1 \otimes (1-zz^*).$$
 Снова повторяя рассуждения, использованные при построении
$\mathscr{D}(\mathbb{D})_q$ и $\mathscr{D}(\mathbb{D})_q^\prime$, но
применительно к алгебре ${\rm Pol}({\mathbb C}\times{\mathbb C})_q$,
получим алгебру $\mathscr{D}(\mathbb{D}\times\mathbb{D})_q$ финитных
функций на декартовом произведении квантовых кругов, а также
$\mathscr{D}(\mathbb{D} \times \mathbb{D})_q$-бимодуль
$\mathscr{D}(\mathbb{D}\times\mathbb{D})_q^\prime$ обобщенных функций.

Назовем элементы $\mathscr{D}(\mathbb{D}\times\mathbb{D})_q^\prime$
обобщенными ядрами, \IND{обобщенные ! ядра} поскольку с их помощью так же,
как в предыдущем \itemiiiе, вводятся интегральные операторы. Именно,
обобщенному ядру $\mathscr{K}$ отвечает линейный оператор из
$\mathscr{D}(\mathbb{D})_q$ в $\mathscr{D}(\mathbb{D})_q^{\prime\,{\rm
op}}$, определяемый равенством
\begin{equation}\label{intoper}
K\,f\,=\,\int_{{\mathbb{D}}_q}\,\mathscr{K}\,\cdot f\,d\nu=({\rm id}\otimes
\nu)(K(1\otimes f)).
\end{equation}

\begin{remark}\label{why_op_?}
Используя противоположное умножение в первом тензорном сомножителе в ${\rm
Pol}({\mathbb C})_q^{\rm op}\otimes{\rm Pol}({\mathbb C})_q$, мы следуем
соглашениям предыдущего \itemiiiа. Это позволяет ввести в рассмотрение
алгебру $U_q\mathfrak{sl}_2$-инвариантных ядер, \IND{алгебра !
$U_q\mathfrak{sl}_2$-инвариантных ядер} см. предложение \ref{IntKer}.
\end{remark}

Из определений вытекает

\begin{lemma}
Для любого множества $\{\psi_{ij}(y,\eta)\}_{i,j\in{\mathbb Z}_+}$ функций
на \hbox{$q^{2{\mathbb Z}_+}\times q^{2{\mathbb Z}_+}$} ряд\ \
$\sum_{i,j\in\mathbb{Z}_+}z^{*i}\zeta^{i}\,\psi_{ij}(y,\eta)\,
z^{j}\zeta^{*j}$\ \ сходится в
$\mathscr{D}(\mathbb{D}\times\mathbb{D})_q^\prime$.
\end{lemma}

\begin{example} Для всех $m \ge 0$ существует обобщенное
ядро
\begin{equation}\label{1.3.3_new}
G_m\,=\,(1-\zeta \zeta^*)^m(q^{-2(m-1)}z^*\zeta;q^2)_m^{-1}\cdot(q^{-2m}z
\zeta^*;q^2)_m^{-1}\cdot(1-z^*z)^m,
\end{equation}
где $(t;q^2)_m=\prod\limits_{j=0}^{m-1}\, (1-tq^{2j})$ и неявно
используется равенство
$$
\frac{1}{(t;q^2)_m}\,=\,
 \sum_{n=0}^\infty \frac{(q^{2m};q^2)_n}{(q^2;q^2)_n}\,t^n,
$$
вытекающее из \eqref{qbinom}.
\end{example}

\medskip
 Следующее утверждение будет доказано в \itemiiiе \ref{proof_inv_G_m}.
\begin{lemma}\label{inv_G_m} $G_m$  является $U_q \mathfrak{sl}_2$-инвариантным
элементом $U_q \mathfrak{sl}_2$-модуля $\mathscr{D}(\mathbb{D} \times
\mathbb{D})_q'$.
\end{lemma}

\medskip

Для того, чтобы сформулировать основной результат этого \itemiiiа,
обратимся к классическому случаю $q=1$ и воспользуемся разложением функции
Грина (\ref{1.2.8}):
\begin{multline}\label{Green_expan}
-\ln\frac{|z-\zeta|^2}{|1-z\overline{\zeta}|^2}=
-\ln\left(1-\frac{(1-|z|^2)(1-|\zeta|^2)}{|1-z\overline{\zeta}|^2}\right)=
\\ =\sum_{m=1}^\infty\frac1{m}
\left(\frac{(1-|z|^2)(1-|\zeta|^2)}{|1-z\overline{\zeta}|^2}\right)^m.
\end{multline}

Переходя к формальному пределу в \eqref{1.3.3_new}, получаем
$$
\lim_{q \to 1}\,G_m=\left(\frac{(1-|z|^2)(1-|\zeta|^2)}{|1-z
\overline{\zeta}|^2}\right)^m.
$$
Следующий результат является $q$-аналогом разложения \eqref{Green_expan}.
Как утверждает следствие \ref{boundness_sl_2}, инвариантный лапласиан
$\square_q$ в $L^2(d \nu)_q$ ограничен и имеет ограниченный обратный.

\medskip

\begin{proposition}\label{t1.3.4}
Сужение ограниченного линейного оператора $\square_q^{-1}$ в $L^2(d \nu)_q$
на плотное линейное подмногообразие $\mathscr{D}(\mathbb{D})_q$ является
интегральным оператором с ядром
$$G_q=\sum\limits_{m=1}^\infty\frac{q^{-2}-1}{q^{-2m}-1}G_m,$$
то есть
$\square_q^{-1}f=\int\limits_{\mathbb{D}_q}G_q(z,\zeta)f(\zeta)d\nu$ при
всех $f\in\mathscr{D}(\mathbb{D})_q$.
\end{proposition}

{\bf Доказательство.} Сравниваемые линейные операторы являются морфизмами
$U_q\mathfrak{g}$-модулей. Кроме того, $U_q \mathfrak{sl}_2\cdot
f_0=\mathscr{D}(\mathbb{D})_q$, как утверждает предложение \ref{t2.3.9}.
Значит, достаточно доказать равенства
\begin{equation}\label{what_we_need_1}
 \int \limits_{{\mathbb{D}}_q}{G}_q\,f_0\,d\nu
\,=\,
 (1-q^2)\sum_{m=1}^{\infty}\frac{q^{-2}-1}{q^{-2m}-1}y^m,
\end{equation}
\begin{equation}\label{what_we_need_2}
\square_q^{-1}f_0\,=\,
(1-q^2)\sum_{m=1}^{\infty}\frac{q^{-2}-1}{q^{-2m}-1}y^m,
\end{equation}
 из которых следует равенство
\begin{equation*}\label{what_we_need}
\square_q^{-1}f_0=\int\limits_{\mathbb{D}_q}G_q\,f_0\,d \nu.
\end{equation*}

 Справедливость  \eqref{what_we_need_1} вытекает из определения
ядра $G_q$ и равенств
$$ z^*\cdot f_0=0,\qquad
\psi(y)\cdot f_0=\psi(1), \qquad \int
\limits_{{\mathbb{D}}_q}\,z^k\,f_0\,d\nu=(1-q^2)\delta_{k,0}.$$

Остается доказать \eqref{what_we_need_2}. Пусть $x=y^{-1}$. Как следует из
обратимости оператора $\square_q$ существует и единствен такой элемент
$\psi\in L^2(d\nu)_q$, что $\square_q\psi\,=\,f_0$, то есть существует и
единственна такая функция $\psi(x)$ на множестве $q^{-2\mathbb{Z}_+}$, что
\begin{equation}\label{expl_sq_new}
 x(1-q^{-1}x)\mathcal{D}_q \psi(x)=q^{-1}\,-\,q,\qquad
 \sum\limits_{j=0}^\infty |\psi(q^{-2j})|^2\cdot q^{-2j} < \infty,
\end{equation}
см. предложение \ref{explicit_square}. Как следует из \eqref{expl_sq_new},
$$
\psi(x)\,=\, \psi(q^{-2}x)\,+\, (q^{-2}-1)^2\,
 \frac{q^{4}x^{-1}}{1-q^{2}x^{-1}},
$$
$$\psi(x)\,=\, (q^{-2}-1)^2\,q^2\,\sum\limits_{j=1}^\infty
 \frac{q^{2j}x^{-1}}{1-q^{2j}x^{-1}}
\,=\, (q^{-2}-1)^2\,q^2\,\sum\limits_{m=1}^\infty
 \frac{q^{2m}}{1-q^{2m}}x^{-m}
$$
 при всех $x \in q^{-2\mathbb{Z}_+}$.
Значит, $\psi(x)=
(1-q^2)\sum\limits_{m=1}^{\infty}\frac{q^{-2}-1}{q^{-2m}-1}x^{-m}$. \hfill
$\square$

\subsubsection{Квантовые конусы.}\label{q-cone}

 Напомним, что алгебра $\mathbb{C}[SL_2]_q$ регулярных функций на
квантовой группе $SL_2$ порождена матричными элементами $t_{ij}$ векторного
представления алгебры $U_q\mathfrak{sl}_2$, см. \eqref{q2-dim_rep}.
  В \itemiiiе \ref{SL_2_sl_2} алгебра
${\mathbb C}[SL_2]_q$ была наделена структурой
$U_q\mathfrak{sl}_2$-модульной алгебры. Как следует из определений,
\begin{equation}\label{r_reg1}
\begin{pmatrix}    K^{\pm 1}t_{11}, &  K^{\pm 1}t_{12} \\
     K^{\pm 1}t_{21}, &  K^{\pm 1}t_{22}\end{pmatrix}
  =
\begin{pmatrix}
t_{11}, & t_{12} \\
t_{21}, & t_{22}\end{pmatrix}\cdot
\begin{pmatrix} q^{\pm 1}, & 0\\ 0,
 & q^{\mp 1}\end{pmatrix}
=
\begin{pmatrix}
    q^{\pm 1}\, t_{11}, & q^{\mp 1}\, t_{12} \\
 q^{\pm 1}\,   t_{21}, & q^{\mp 1}\, t_{22}\end{pmatrix},
\end{equation}

\begin{equation}\label{r_reg2}
\begin{pmatrix}    Et_{11}, &  Et_{12} \\
     Et_{21}, &  Et_{22}\end{pmatrix}
  =
\begin{pmatrix}
t_{11}, & t_{12} \\
t_{21}, & t_{22}\end{pmatrix}\cdot
\begin{pmatrix} 0, & q^{-\frac{1}{2}}\\ 0, & 0\end{pmatrix}
=\begin{pmatrix}
    0, & q^{-\frac{1}{2}}\, t_{11} \\
    0, & q^{-\frac{1}{2}}\, t_{21}\end{pmatrix},
\end{equation}

\begin{equation}\label{r_reg3}
 \begin{pmatrix}    Ft_{11}, &  Ft_{12} \\
     Ft_{21}, &  Ft_{22}\end{pmatrix}
  =
\begin{pmatrix}
t_{11}, & t_{12} \\
t_{21}, & t_{22}\end{pmatrix}\cdot
\begin{pmatrix} 0, & 0\\ q^{\frac{1}{2}}, & 0\end{pmatrix}
=\begin{pmatrix}
   q^{\frac{1}{2}}\, t_{12}, & 0 \\
   q^{\frac{1}{2}}\, t_{22}, & 0
  \end{pmatrix}.
\end{equation}

\bigskip
 Введем квантовые аналоги пространств  функций на
конусе
$$\widetilde{\Xi}=\{(t_{11},t_{12})
\in{\mathbb C}^2\quad|\quad |t_{11}|=|t_{12}|\}.$$

Алгебра $\mathbb{C}[SL_2]_q$ может быть определена с помощью образующих
$t_{ij}$ и соотношений, все из которых однородны, кроме \eqref{q-det_sl_2},
см. предложение \ref{su_2_relations}. Сделаем в определяющих соотношениях
подстановку
\begin{equation}\label{subst}
t_{ij}=q^{-N}\cdot t_{ij}^{(N)},\qquad N\in{\mathbb N}.
\end{equation}
Формальный предельный переход $N\to+\infty$ не влияет на однородные
соотношения, а \eqref{q-det_sl_2} переходит в
$t_{11}t_{22}-qt_{12}t_{21}=0$.

Введем обозначение $\mathbb{C}[\widetilde{\Xi}]_q$ для алгебры,
определяемой новым множеством соотношений между образующими $t_{ij}$.

\begin{remark} Равенство  $\deg(t_{ij})=1$ наделяет
алгебру $\mathbb{C}[SL_2]_q$ фильтрацией, а алгебру
$\mathbb{C}[\widetilde{\Xi}]_q$-- градуировкой. Описанный выше переход от
$\mathbb{C}[SL_2]_q$ к $\mathbb{C}[\widetilde{\Xi}]_q$ -- это
``приземленное`` описание стандартного перехода от фильтрованной алгебры к
присоединенной градуированной алгебре.
\end{remark}

Введем инволюцию $*$ в $\mathbb{C}[\widetilde{\Xi}]_q$, полагая
\begin{equation}\label{inv_cone}
t_{11}^*=-t_{22},\qquad t_{12}^*=-qt_{21}.
\end{equation}
Из определений следует, что $*$-алгебра $\mathbb{C}[\widetilde{\Xi}]_q$ с
действием $U_q\mathfrak{su}_{1,1}$, определяемым равенствами \eqref{r_reg1}
-- \eqref{r_reg3}, является $U_q\mathfrak{su}_{1,1}$-модульной.

Положим $x=-qt_{12}t_{21}$. Докажем, что каждый элемент алгебры
$\mathbb{C}[\widetilde{\Xi}]_q$ единственным образом разлагается в сумму
\begin{equation}\label{canon_expan_cone}
f=\sum\limits_{ad=bc=0}t_{11}^at_{12}^b\psi_{abcd}(x)t_{21}^ct_{22}^d,
\end{equation}
где $a,b,c,d\in\mathbb{Z}_+$ и $\psi_{abcd}$ -- полиномы одной переменной.
Существование такого разложения очевидно.

Рассмотрим предгильбертово пространство с ортонормированным базисом
$\{e_j\}_{j \in \mathbb{Z}_+}$ и следующее $*$-представление $\pi$ алгебры
$\mathbb{C}[\widetilde{\Xi}]_q$:
$$
\pi(t_{11})e_j\,=\,q^{-(j+1)}\,e_{j+1}, \quad \pi(t_{12})e_j\,=\,q^{-j}e_j,
$$
$$
 \pi(t_{21})e_j\,=\,-q^{-(j+1)}e_j,\quad
 \pi(t_{22})e_j=
-q^{-j}\,e_{j-1}.
$$
 С его помощью
линейная независимость элементов множества
$$
 \{t_{11}^a t_{12}^b t_{21}^c,
\,t_{12}^b t_{21}^c t_{22}^d\,|\, a, b, c, d \in \mathbb{Z}_+\}
$$
доказывается так же, как в \itemiiiе \ref{SL_2_sl_2}. Отсюда следует
единственность разложения \eqref{canon_expan_cone}.

Нетрудно показать, что $x^{\mathbb{Z}_+}$ является множеством Оре и что
отвечающая ему локализация $\mathbb{C}[\widetilde{\Xi}]_{q,x}$ алгебры
$\mathbb{C}[\widetilde{\Xi}]_q$ наследует структуру
$U_q\mathfrak{su}_{1,1}$-модульной алгебры. Элементы алгебры
$\mathbb{C}[\widetilde{\Xi}]_{q,x}$ имеют вид \eqref{canon_expan_cone}, где
$\psi_{abcd}$ -- полиномы Лорана одной переменной.

 Рассмотрим $U_q\mathfrak{su}_{1,1}$-модульные подалгебры $\mathbb{C}[\Xi]_q \subset
\mathbb{C}[\widetilde{\Xi}]_q$, \
 $\mathbb{C}[\Xi]_{q,x} \subset
\mathbb{C}[\widetilde{\Xi}]_{q,x}$ инвариантов однопараметрической группы
автоморфизмов
$$
t_{1j} \mapsto e^{i\varphi}t_{1j},\quad t_{2j} \mapsto e^{-i\varphi}t_{2j},
\qquad \varphi \in \mathbb{R}/(2\pi \mathbb{Z}), \; j=1,2.
$$
Это квантовые аналоги алгебр функций на двумерном конусе. В разложении
\eqref{canon_expan_cone} их элементов равны нулю слагаемые с\ \ $a+b\neq
c+d$.

\bigskip

Наделим алгебру $\mathbb{C}[\Xi]_{q,x}$ градуировкой, полагая
$$\deg(x^{\pm 1})=\pm 1,\qquad\deg(t_{1i}t_{2j})=1,\quad i,j=1,2.$$
Из равенств \eqref{K_psi} -- \eqref{F_psi} следует, что однородные
компоненты являются $U_q\mathfrak{sl}_2$-подмодулями
$U_q\mathfrak{sl}_2$-модуля $\mathbb{C}[\Xi]_{q,x}$. Докажем, что они
изоморфны $U_q\mathfrak{sl}_2$-модулям основной сферической серии. Прежде
всего, отметим, что элементы $t_{12}$, $t_{21}$, $t_{11}$, $t_{22}$ алгебры
$\mathbb{C}[\widetilde{\Xi}]_{q,x}$ обратимы. Следовательно, определен и
обратим элемент $z=t_{12}^{-1}t_{11}\in\mathbb{C}[\Xi]_{q,x}$.

\begin{proposition}\label{in_cone} Пусть $l\in \mathbb{Z}$.
Линейная оболочка $ M^{(l)}$ множества $\{z^k\,x^l\;|\; k\in \mathbb{Z}\}$
является $U_q\mathfrak{sl}_2$-модулем, изоморфным
$U_q\mathfrak{sl}_2$-модулю $V^{(l)}$ основной сферической серии.
\end{proposition}

{\bf Доказательство.} Как нетрудно показать, для всех $l$ и всех полиномов
Лорана $f(z)$
$$K^{\pm 1}\,\left(f(z)x^l\right)\,=\,f(q^{\pm 2}z)x^l,$$
$$
E\left(f(z)x^l\right)\,=\,
-\,q^{\frac12}z^2\frac{f(z)-f(q^2z)}{z-q^2z}\,x^l\,+\,
q^{-\frac{3}{2}}\,\frac{q^{-2l}-1}{q^{-2}-1}zf(q^2z)x^l,
$$
$$
F\left(f(z)x^l\right)\,=\,
q^{\frac12}\,\frac{f(q^{-2}z)-f(z)}{q^{-2}z-z}x^l\,+\,
q^{\frac52}\,\frac{1-q^{2l}}{1-q^2}z^{-1}f(z)x^l.
$$
Эти равенства позволяют найти явный вид морфизма
$U_q\mathfrak{sl}_2$-модулей
\begin{equation}\label{embed_to_cone}
j^{(l)}:\,V^{(l)}\to M^{(l)},\qquad j^{(l)}:\,e_0\mapsto x^l.
\end{equation}
Именно,
\begin{equation}\label{embed_explicit}
j^{(l)}\,e_k\,=\,(-q^{-(l+\frac12)})^k\,z^kx^l,\qquad k\in\mathbb{Z}.
\end{equation}
Остается воспользоваться биективностью $j^{(l)}$.\hfill$\square$

\bigskip

От требования $l\in \mathbb{Z}$ легко освободиться, введя в рассмотрение
обобщенные функции на квантовых конусах. Наметим требуемые построения. Они
основаны на том, что множеством собственных значений оператора $\pi(x)$
является $q^{-2\mathbb{Z}}$.

Пусть $\mathscr{D}(\widetilde{\Xi})'_q$ -- векторное пространство
формальных рядов \eqref{canon_expan_cone} с коэффициентами из пространства
функций на множестве $q^{-2\mathbb{Z}}$. Наделим
$\mathscr{D}(\widetilde{\Xi})'_q$ слабейшей из топологий, в которых
непрерывны линейные функционалы $\mathscr{L}_{a,b,c,d,n}$:
\begin{equation}\label{topology_SU_1_1}
\mathscr{L}_{a,b,c,d,n}(f)\;=\;\psi_{abcd}(q^{-2n}),\qquad ad=bc=0,\;n\in
\mathbb{Z}.
\end{equation}
По непрерывности $\mathscr{D}(\widetilde{\Xi})'_q$ наделяется структурой
$U_q\mathfrak{sl}_2$-модульного $\mathbb{C}[\widetilde{\Xi}]_q$-бимодуля.

Для обобщенных функций на квантовом конусе можно ввести понятие комплексной
степени однородности. \IND{комплексная степень однородности} Пусть
$l\in\mathbb{C}$. Обобщенную функцию $f\in\mathscr{D}(\widetilde{\Xi})'_q$
назовем однородной степени $l$, если коэффициенты $\psi_{abcd}(x)$ ее
разложения \eqref{canon_expan} удовлетворяют следующим требованиям:
$$
\psi_{abcd}(q^2x)=q^{2(l-\frac{a+b+c+d}{2})}\psi(x),\qquad x\in
q^{-2\mathbb{Z}}.
$$
Эквивалентно, пусть $\alpha$ - это автоморфизм
$\mathbb{C}[\widetilde{\Xi}]_q$, заданный равенствами
$\alpha(t_{ij})=qt_{ij}$, где $i,j=1,2$, а $\overline{\alpha}$ - его
продолжение по непрерывности на $\mathscr{D}(\widetilde{\Xi})_q^\prime$.
Элемент $\psi \in \mathscr{D}(\widetilde{\Xi})_q^\prime$ - однородный
степени $l$, если $\overline{\alpha}(\psi)=q^{2l} \psi$.

Из равенств \eqref{K_psi} -- \eqref{F_psi} следует, что обобщенные функции
фиксированной степени однородности образуют $U_q\mathfrak{sl}_2$-подмодуль.
Важно отметить, что в определении подпространства $M^{(l)}$ и в
доказательстве предложения \ref{in_cone} практически не используется
целочисленность $l$.

Используя предложение \ref{in_cone} и равенство \eqref{2.5.8}, нетрудно
доказать, что
\begin{equation}\label{2.5.8_new}
C_{q}v\,=\, -q\dfrac{(1-q^{-2l})(1-q^{2(l+1)})}{(1-q^2)^2}v, \qquad v\in
M^{(l)}
\end{equation}
при всех $l\in \mathbb{C}$.

\bigskip
\begin{proposition}\label{homogeneous_int} Линейный функционал
на подпространстве однородных обобщенных функции степени $-1$, определяемый
равенством
\begin{equation}\label{int_cone}
\int\limits_{\widetilde{\Xi}_q}\,f\, d\eta\;=\;\psi_{0000}(1),
\end{equation}
является $U_q\mathfrak{sl}_2$-инвариантным.
\end{proposition}

{\bf Доказательство.} Нетрудно убедиться в том, что
$$\int\limits_{\widetilde{\Xi}_q}(K^{\pm 1}-1)x^{-1}\,d\eta=0,$$
$$ \int\limits_{\widetilde{\Xi}_q} E(\,t_{12}x^{-2}t_{22}\,)\,d\eta\,=\,
\int\limits_{\widetilde{\Xi}_q} F(\,t_{11}x^{-2}t_{21}\,)\,d\eta\,=\,0.$$
Отсюда легко следует $U_q\mathfrak{sl}_2$-инвариантность рассматриваемого
линейного функционала. \hfill $\square$

\bigskip
Элемент $f \in \mathscr{D}(\widetilde{\Xi})'_q$ назовем финитной функцией
на квантовом конусе, если ${\rm card}\{(a,b,c,d,n)\;|\;
\mathscr{L}_{a,b,c,d,n}(f)\neq 0\} < \infty$. Введем обозначение
$\mathscr{D}(\widetilde{\Xi})_q$ для векторного пространства финитных
функций на квантовом конусе. Оно является $U_q\mathfrak{sl}_2$-модульным
$\mathbb{C}[\widetilde{\Xi}]_q$-бимодулем и по непрерывности наделяется
структурой $U_q\mathfrak{su}_{1,1}$-модульной $*$-алгебры.

\subsubsection{Собственные функции оператора
$\mathbf{\square_q}$.}\label{eigen_functions}

В классическом случае $q=1$ функция
$$
u(z)=\int\limits_{\partial\mathbb{D}}
\left(\frac{1-|z|^2}{(1-z\overline\zeta)(1-\overline
z\zeta)}\right)^{l+1}f(\zeta)d\nu,\qquad d\nu=\frac{d\theta}{2\pi},
$$
является собственной функцией инвариантного лапласиана $\square$ при весьма
широких предположениях о функции $f$. Именно,
$$\square u=-\left( (l+\frac12)^2-\frac14\right)u,$$
см. \cite[стр. 51]{Helg1}.

Желая получить квантовый аналог этого результата, введем квантовые аналоги
степеней ядра Пуассона и пространств функций на единичной окружности
$\partial\mathbb{D}$. Вложение алгебр
$$
\mathbb{C}[z,z^{-1}]\hookrightarrow \mathbb{C}[\widetilde{\Xi}]_{q,x},
\qquad z\mapsto t_{12}^{-1}t_{11},
$$
позволяет наделить алгебру полиномов Лорана переменной $z$ структурой
$U_q\mathfrak{su}_{1,1}$-мо\-дуль\-ной алгебры c инволюцией $z^*\,=\,z^{-
1}$. Введем обозначение
$\mathbb{C}[\partial\mathbb{D}]_q\,=\,(\mathbb{C}[z,z^{-1}],*)$ для
полученной $U_q\mathfrak{su}_{1,1}$-модульной алгебры. Это -- квантовый
аналог алгебры регулярных функций на единичной окружности
$\partial\mathbb{D}$.

Наделим $\mathbb{C}[\partial\mathbb{D}]_q$ топологией покоэффициентной
сходимости, то есть слабейшей из топологий, в которых непрерывны линейные
функционалы
\begin{equation}\label{l_j}
l_j:\;\sum\limits_{n\in\mathbb{Z}}\,c_nz^n\mapsto c_j,\qquad
j\in\mathbb{Z}.
\end{equation}
Действие $U_q\mathfrak{sl}_2$ и инволюция $*$ допускают продолжение по
непрерывности с $\mathbb{C}[\partial\mathbb{D}]_q$ на пополнение --
топологическое векторное пространство формальных рядов
$\sum\limits_{n\in\mathbb{Z}}\,c_nz^n$ с комплексными коэффициентами и
топологией покоэффициентной сходимости. Введем обозначение
$\mathbb{C}[[\partial\mathbb{D}]]_q$ для полученного
$U_q\mathfrak{sl}_2$-модульного
$\mathbb{C}[\partial\mathbb{D}]_q$-бимодуля.

 Наконец, наделим
${\rm Pol}(\mathbb{C})_q^{\rm op}\otimes \mathbb{C}[\partial \mathbb{D}]_q$
слабейшей из топологий, в которых непрерывны линейные функционалы
$l''_{m,n}\otimes l_j$, см. \eqref{l_two_prime}, \eqref{l_j}.
 Пополнение
$\mathscr{D}(\mathbb{D}\times
\partial \mathbb{D})'_q$ можно
отождествить с векторным пространством формальных рядов $\sum\limits_{n\in
\mathbb{Z}}\,c_n\,z^n$ с
 коэффициентами $c_n$ из $\mathscr{D}(\mathbb{D})'_q$.

Введем $q$-аналоги степеней $P^\gamma$ ядра Пуассона
$P=\dfrac{1-|z|^2}{|1-z\overline{\zeta}|^2}$, используя обозначения
\itemiiiа \ref{Green_function}. В квантовом случае $q\in(0,1)$ мы заменяем
    $P^\gamma$ элементом
\begin{equation}\label{1.5.2}
P_\gamma=(q^2z^*\zeta;q^2)_{-\gamma}\;(z \zeta^*;q^2)_{-\gamma}
(1-z^*z)^\gamma\;,
\end{equation}
где ядра $(z \zeta^*;q^2)_{-\gamma}$, $(q^2z^*\zeta;q^2)_{-\gamma}\in
\mathscr{D}(\mathbb{D}\times\partial\mathbb{D})'_q$, определяемые
равенствами
\begin{equation}\label{qbinom_P1}
(z \zeta^*;q^2)_{-\gamma}\,=\,\sum_{n=0}^\infty
\frac{(q^{2\gamma};q^2)_n}{(q^2;q^2)_n}(q^{-2\gamma}z\zeta^*)^n\,=\,
\frac{(z\zeta^*;q^2)_{\infty}}{(q^{-2\gamma}z \zeta^*;q^2)_{\infty}},
\end{equation}
\begin{equation}\label{qbinom_P2}
(q^2z^*\zeta;q^2)_{-\gamma}\,=\,\sum_{n=0}^\infty
\frac{(q^{2\gamma};q^2)_n}{(q^2;q^2)_n}(q^{2-2\gamma}z^*\zeta)^n\,=\,
\frac{(q^2z^*\zeta;q^2)_{\infty}}{(q^{-2\gamma+2}z^*\zeta;q^2)_{\infty}},
\end{equation}
являются $q$-аналогами биномов $(1-z\overline{\zeta})^{-\gamma}$ и
$(1-\overline{z}\zeta)^{-\gamma}$.

\medskip  Основной результат этого \itemiiiа состоит в следующем.

\begin{proposition}\label{t1.5.1}
При всех $l\in\mathbb{C}$, $f\in{\mathbb C}[\partial\mathbb{D}]_q$ элемент
\begin{equation}\label{sobst}
u=\int\limits_{0}^{2\pi}P_{l+1}(z,e^{i\theta})f(e^{i\theta})
\frac{d\theta}{2\pi}
\end{equation}
принадлежит $\mathscr{D}(\mathbb{D})_q^\prime$, и
\begin{equation}\label{eigenvectors_sl_2}
\square_qu=\alpha(l)u,
\end{equation}
где $\alpha(l)=(1-q^{-2l})(1-q^{2l+2})/(1-q^2)^2$.
\end{proposition}

{\bf Доказательство.} Как вытекает из следствия \ref{C_q-Laplace},
достаточно доказать равенство $C_q u=-q\alpha(l)u$.

Используя \eqref{2.5.8_new}, нетрудно доказать что в
$\mathscr{D}(\widetilde{\Xi})'_q$ имеет место равенство
$C_q\,\left(f(z)x^l\right)\,=\,-q\alpha(l)\left(f(z)x^l\right)$ для любого
полинома Лорана $f(z)$ и всех $l\in \mathbb{C}$.

В \itemiiiе \ref{q-cone} были введены $U_q\mathfrak{sl}_2$-модули $M^{(l)}$
элементов степени однородности $l$ на квантовом аналоге конуса $\Xi$. Для
завершения доказательства предложения достаточно предъявить морфизм
$U_q\mathfrak{sl}_2$-модулей $M^{(l)}\to \mathscr{D}(\mathbb{D})_q^\prime$,
отображающий $f(z)x^l$ в $u$.

\medskip
Наделим векторное пространство ${\rm Pol}(\mathbb{C})_q^{\rm op}\otimes
\mathbb{C}[\Xi]_q$ слабейшей из топологий, в которых непрерывны линейные
функционалы $l''_{m,n}\otimes \mathscr{L}_{a,b,c,d,n}$ (см.
\eqref{l_two_prime}, \eqref{topology_SU_1_1}). Следуя обозначениям
\itemiiiа \ref{Green_function}, полагаем $\xi=1\otimes x$.

В \itemiiiе \ref{proof_inv_G_m} будет доказана
$U_q\mathfrak{sl}_2$-инвариантность элемента
$$\mathscr{K}_l=\xi^l \sum_{j=0}^\infty
\frac{(q^{-2l};q^2)_j}{(q^2;q^2)_j}(q^{2(l+1)}z^*\zeta)^j \cdot
\sum_{m=0}^\infty \frac{(q^{-2l};q^2)_m}{(q^2;q^2)_m}(q^{2l}z
\zeta^*)^m(1-z^*z)^{-l}$$ пополнения пространства
 ${\rm Pol}(\mathbb{C})_q^{\rm op}\otimes \mathbb{C}[\Xi]_q$,
 см. замечание \ref{cone-cone},
а в \itemiiiе \ref{q-cone} доказана $U_q\mathfrak{sl}_2$-инвариантность
линейного функционала
$$
 \eta:\, M^{(-1)}\to \mathbb{C},\qquad \eta:
\left(\sum_{m=-\infty}^\infty a_m z^m \right)x^{-1}\,\mapsto\,a_0.
$$
 Произведение
$\mathscr{K}_{-l-1}(1\otimes f(z)x^l)$ по второму тензорному сомножителю
принадлежит $M^{(-1)}$. Значит, корректно определен интегральный оператор
$$
M^{(l)}\to\mathscr{D}(\mathbb{D})_q^\prime,\qquad f \mapsto ({\rm id}
\otimes \eta)\left(\mathscr{K}_{-l-1}(1\otimes f(z)x^l)\right).
$$
Из приведенных выше результатов и равенств \eqref{kernel_equation},
\eqref{l_new} следует, что он является морфизмом
$U_q\mathfrak{sl}_2$-модулей:
$$
({\rm id}\otimes\eta)((a\otimes 1)\mathscr{K}_{-l-1}(1\otimes f(z)x^l))=
 ({\rm id}\otimes\eta)
((1\otimes S^{-1}(a))\mathscr{K}_{-l-1}(1\otimes f(z)x^l)),$$ $$((1\otimes
S^{-1}(a))\mathscr{K}_{-l-1}(1\otimes f(z)x^l))=
 ({\rm id}\otimes\eta)(\mathscr{K}_{-l-1}(1 \otimes a(f(z)x^l))).
$$
 Остается заметить, что $({\rm id}
\otimes \eta)\left(\mathscr{K}_{-l-1}(1\otimes fx^l)\right)=u$.\hfill
$\square$

\begin{remark}
Полагая $\gamma=1$ в \eqref{1.5.2}, получаем $q$-аналог ядра Пуассона
$P=\dfrac{1-|z|^2}{|1-z\overline{\zeta}|^2}$. При $\gamma=1$ интегральный
оператор \eqref{sobst} отображает
$f=\sum\limits_{n=-\infty}^\infty\,c_ne^{in\theta}$ в
$u=\sum\limits_{n=1}^\infty\,c_n\,z^n\;+\;c_0\;+\;
\sum\limits_{n=0}^{\infty}\,c_{-n}\,z^{*n}$ и допускает продолжение по
непрерывности до изометрического линейного оператора из пространства
$C(\partial\mathbb{D})$ непрерывных функций на окружности в пространство
$C(\mathbb{D})_q$ непрерывных функций в квантовом круге. (Последнее
утверждение нетрудно доказать используя унитарную дилатацию сжатия
$T_F(z)$, см. \itemiii\ \ref{Nagy-Foias}.) Разумеется, $\square_q u=0$,
$u_{|\;\partial \mathbb{D}}=f$.
\end{remark}

\subsubsection{Гармонический анализ.}\label{Fourier_sl_2}

Начнем с хорошо известных результатов, относящихся к классическому случаю
$q=1$ \cite[стр. 48, 49]{Helg1}. Неевклидово преобразование Фурье
\IND{неевклидово преобразование Фурье функции в круге} функции $u$ в круге
$\mathbb{D}$ определяется равенством
\begin{equation}\label{Noneu_Fourier}
\widetilde{u}(e^{i\theta},\rho)=\int\limits_{\mathbb{D}}
P^{\frac{1}{2}-i\rho}(z,e^{i\theta})u(z)d\nu,
\end{equation}
где $P=\frac{1-|z|^2}{(1-z\overline\zeta)(1-\overline z\zeta)}$ -- ядро
Пуассона. Обратное преобразование имеет вид
\begin{equation}\label{inverse_Fourier}
u(z)=\int\limits_0^\infty\left(\int\limits_0^{2\pi}P^{\frac12+i\rho}
(z,e^{i\theta})\widetilde{u}(e^{i\theta},\rho)\frac{d\theta}{2\pi}\right)
d\sigma_{\rm classic}(\rho),
\end{equation}
где $d\sigma_{\rm classic}(\rho)=2\rho\,{\rm th}(\pi\rho)\,d\rho$
(несущественное отличие этой меры Планшереля от приведенной в литературе
объясняется выбором используемой нами инвариантной меры $d\nu=\frac{d{\rm
Re}z\;d{\rm Im}z}{\pi(1-|z|)^2}$ в $\mathbb{D}$).

Неевклидово преобразование Фурье является унитарным линейным оператором из
$L^2(d\nu)$ в $L^2(d\sigma)\otimes L^2(\frac{d\theta}{2\pi})$.

 Отметим, см. \itemiii\ \ref{eigen_functions}, что равенство
\eqref{inverse_Fourier} доставляет разложение $u(z)$ по собственным функциям
$u_\rho(z)=\int\limits_0^{2\pi} P^{\frac{1}{2}+i\rho}(z,e^{i\theta})
 \widetilde{u}(e^{i\theta},\rho)\frac{d\theta}{2\pi}$ оператора $\square$.
\ \ Получим $q$-аналог этого разложения. Пусть
$$
P_{l+1}^t=(q^2z^*\zeta;q^2)_{-l-1}(z \zeta^*;q^2)_{-l-1}(1-\zeta
\zeta^*)^{1+l},$$
 см. равенства \eqref{qbinom_P1},\eqref{qbinom_P2}.

\begin{proposition}\label{t1.6.1}
\hspace{-.5em}. Рассмотрим борелевскую меру $d\sigma$ на
$\left[0,\frac{\pi}{2\log(q^{-1})}\right]$, введенную равенством
\eqref{Pl}. Интегральные операторы
$$
u(z)\;\mapsto\;\int\limits_{\mathbb{D}_q}P_{\frac12-i\rho}^t
(e^{i\theta},\zeta)u(\zeta)d\nu,
$$
\begin{equation}\label{qF_sl_2}
\widetilde{u}(e^{i\theta},\rho)\mapsto\int\limits_0^{\pi/h}\,
\left(\int\limits_0^{2\pi}P_{\frac12+i\rho}(z,e^{i\theta})
\widetilde{u}(e^{i\theta},\rho)\frac{d\theta}{2\pi}\right)\,d\sigma(\rho)
\end{equation}
продолжаются по непрерывности с плотных линейных подмногообразий
$$
\mathscr{D}(\mathbb{D})_q\subset
L^2(d\nu)_q,\quad\mathbb{C}[\partial\mathbb{D}]_q\otimes
C^\infty\left[0,\frac{\pi}{2\log(q^{-1})}\right]\subset
L^2\left(\frac{d\theta}{2\pi}\right)\otimes L^2(d\sigma)
$$
до взаимно обратных унитарных операторов.
\end{proposition}

{\bf Доказательство.} Пусть $-\frac{\pi}{2\log(q^{-1})}\le\,{\rm Im}\,
l\,\le\,\frac{\pi}{2\log(q^{-1})}$. Рассмотрим вложение векторных
пространств
$$
\mathbb{C}[\partial\mathbb{D}]\to\mathscr{D}(\Xi)'_q,\qquad f(z)\mapsto
f(z)x^l.
$$
Наделим $\mathbb{C}[\partial\mathbb{D}]_q$ структурой
$U_q\mathfrak{sl}_2$-модуля так, чтобы это вложение стало морфизмом
$U_q\mathfrak{sl}_2$-модулей, и введем для полученного
$U_q\mathfrak{sl}_2$-модуля обозначение
$\mathbb{C}[\partial\mathbb{D}]^{(l)}_q$.

Как следует из предложения \ref{in_cone}, $ V^{(-1/2+i\rho)}\cong
\mathbb{C}[\partial \mathbb{D}]^{(-1/2+i\rho)}_q$. Эта реализация основной
унитарной серии представлений $*$-алгебры $U_q\mathfrak{su}_{1,1}$ в
пространстве функций на окружности хорошо известна в классическом случае
$q=1$. Она позволяет заменить введенный
 в \itemiiiе \ref{boundness_Laplace}  унитарный оператор
\begin{equation}
\overline{\mathcal{I}}:\,L^2(d \nu)_q \to \newoplus
\int\limits_{0}^{\frac{\pi}{2\log(q^{-1})}}
\overline{V^{(-1/2+i\rho)}}\,d\sigma(\rho)
\end{equation}
унитарным оператором
\begin{equation}
\overline{\mathcal{J}}:\,L^2(d \nu)_q \to \newoplus
\int\limits_{0}^{\frac{\pi}{2\log(q^{-1})}}
L_2\left(\frac{d\theta}{2\pi}\right)\,d\sigma(\rho).
\end{equation}

Докажем, что $\overline{\mathcal{J}}\,u=\int
\limits_{\mathbb{D}_q}P_{\frac12-i\rho}^t(z,\zeta)u(\zeta)d\nu$ для всех
$u\in\mathscr{D}(\mathbb{D})_q$. Это утверждение является следствием двух
лемм. Введем линейный оператор
$$
j_\rho:\,\mathscr{D}(\mathbb{D})_q\to
\mathbb{C}[\partial\mathbb{D}]^{(-1/2+i\rho)}_q,\qquad
j_\rho:f\mapsto\int\limits_{\mathbb{D}_q}P_{\frac12-i\rho}^t(z,\zeta)
f(\zeta)d\nu.
$$

\begin{lemma} $j_\rho f_0=1-q^2$.
\end{lemma}

\smallskip

{\bf Доказательство леммы.} Из определений вытекают равенства
$$
P_l^t(z,\zeta)=\sum_{j>0}\zeta^j \cdot \psi_j
+\psi_0+\sum_{j>0}\psi_{-j}\cdot\zeta^*
$$
 и $\psi_0=1$. Остается отметить, что
ненулевой вклад в интеграл $j_\rho f_0$ дает только одно слагаемое
$\psi_0(z,\xi)$ и что $\displaystyle \int \limits_{\mathbb{D}_q}\, f_0\,d
\nu=1-q^2$. \hfill $\square$

\medskip

\begin{lemma}\label{irho}
Линейный оператор $j_\rho:\mathscr{D}(\mathbb{D})_q\to{\mathbb C}[\partial
\mathbb{D}]^{(-\frac1+i \rho)}_q$ является морфизмом
$U_q\mathfrak{sl}_2$-модулей.
\end{lemma}

\smallskip

{\bf Доказательство леммы.} Рассмотрим интегральный оператор
$$
i_\rho:{\mathbb C}[\partial\mathbb{D}]^{(-\frac12+i\rho)}_q\to
\mathscr{D}(\mathbb{D})_q',\qquad i_\rho:f\mapsto
\int\limits_0^{2\pi}P_{\frac12+i\rho}(z,e^{i\theta})f(e^{i\theta})\frac{d
\theta}{2\pi}.
$$
Как видно из доказательства предложения \ref{t1.5.1}, $i_\rho$ является
морфизмом $U_q \mathfrak{sl}_2$-модулей.

Наделим $U_q \mathfrak{su}(1,1)$-модули $\mathscr{D}(\mathbb{D})_q$,
${\mathbb C}[\partial \mathbb{D}]^{(-\frac12+i\rho)}_q$ инвариантными
скалярными произведениями
$$
\mathscr{D}(\mathbb{D})_q \times \mathscr{D}(\mathbb{D})_q \to{\mathbb
C},\qquad f_1 \times f_2 \mapsto \int \limits_{\mathbb{D}_q}f_2^*f_1d \nu,
$$
$$
{\mathbb C}[\partial\mathbb{D}]^{(l)}_q\times{\mathbb C}[\partial
\mathbb{D}]^{(l)}_q\to{\mathbb C},\qquad f_1\times f_2\mapsto\int
\limits_{\partial\mathbb{D}}f_2^*f_1\frac{d\theta}{2\pi}.
$$

Интегральный оператор с ядром \hbox{$K=\sum \limits_ik_i''\otimes k_i'$}
сопряжен к интегральному оператору с ядром $K^t=\sum \limits_ik_i^{'*}
\otimes k_i^{''*}$. Следовательно, $j_ \rho=i_ \rho^*$, и $j_ \rho$--
морфизм $U_q \mathfrak{sl}_2$-модулей, поскольку этим свойством обладает
отображение $i_\rho$. \hfill $\square$

\medskip Завершим доказательство предложения \ref{t1.6.1}.
 Как видно из доказательства леммы
\ref{irho}, линейный оператор, сопряженный к $\mathcal{J}$, является
интегральным оператором \eqref{qF_sl_2}. Для доказательства предложения
\ref{t1.6.1} остается воспользоваться унитарностью оператора $\mathcal{J}$,
см. \itemiii\ \ref{boundness_Laplace}. \hfill $\square$

\bigskip

Как следует из результатов \itemiiiа \ref{eigen_functions}, предложение
\ref{t1.6.1} доставляет разложение элемента $u$ по собственным функциям
$$
u_\rho\,=\,\int\limits_0^{2\pi}P_{\frac12+i\rho}(z,e^{i\theta})
\widetilde{u}(e^{i\theta},\rho)\frac{d\theta}{2\pi}
$$
оператора $\square_q$, где
$\widetilde{u}=\int\limits_{\mathbb{D}_q}P_{\frac12-i\rho}^t
(e^{i\theta},\zeta)u(\zeta)d\nu$ является $q$-аналогом неевклидова
преобразования Фурье. \IND{$q$-аналог ! неевклидова преобразования Фурье}

\subsubsection{Операторы Теплица в пространствах
Бергмана.}\label{weighted_Bergman_sl_2}

Пусть $\lambda>1$. Как объясняется в \itemiiiе \ref{bundles_sl_2} (см.
также \eqref{U_q_sl_2-bundles}, \eqref{EF-q-diff}), равенства
\begin{equation}\label{E_lambda_sl_2}
Ef(z)\,=\,-\,q^{1/2}z^2\frac{f(z)-f(q^2\,z)}{z-q^2\,z}\;-q^{1/2}\,
\frac{1-q^{2\lambda}}{1-q^2}zf(q^2z),
\end{equation}
\begin{equation}\label{F_lambda_sl_2}
Ff(z)\,=\,q^{1/2-\lambda}\,\frac{f(q^{-2}\,z)-f(z)}{q^{-2}\,z-z},\qquad
K^{\pm 1}f(z)\,=q^{\pm\lambda}\,f(q^{\pm 1}z)
\end{equation}
наделяют алгебру полиномов $\mathbb{C}[z]$ структурой унитаризуемого
$U_q\mathfrak{su}_{1,1}$-модуля. Инвариантное скалярное произведение
определяется с помощью вложения в гильбертово пространство
$L^2(d\nu_\lambda)_q$, введенное в \itemiiiе \ref{bundles_sl_2}:
$$
(f_1,f_2)_\lambda\,=\,\frac{1-q^{2(\lambda-1)}}{1-q^2}\,
\int\limits_{\mathbb{D}_q}\,f_2^*f_1(1-zz^*)^\lambda\,d\nu\,=\,
\int\limits_{\mathbb{D}_q}f_2^*f_1d\nu_\lambda.
$$
Это -- $q$-аналог сферической голоморфной дискретной серии. \IND{$q$-аналог
! сферической голоморфной дискретной серии}

\bigskip

Замыкание $\mathbb{C}[z]$ в $L^2(d\nu_\lambda)_q$, является $q$-аналогом
взвешенного пространства Бергмана \IND{$q$-аналог ! взвешенного
пространства Бергмана} \cite[стр. 2]{HedKorZhu} и обозначается
$L^2_a(d\nu_\lambda)_q$.

\begin{lemma}\label{l1.7.1}
\hspace{-.5em}. Мономы $\{z^n\}_{n \in{\mathbb Z}_+}$ попарно ортогональны
в $L^2_a(d\nu_\lambda)_q$, и
\begin{equation}\label{Klimek-Lesn}
 \| z^n \|^2_\lambda\,=\,\frac{(q^2;q^2)_n}{(q^{2 \lambda };q^2)_n}.
\end{equation}
\end{lemma}

\smallskip
{\bf Доказательство}. Попарная ортогональность мономов $z^n$ очевидна.
Используя свойства $q$-бета функции и $q$-гамма функции, приведенные в
\itemiiiе \ref{q-series}, получаем:
$$\frac{ 1-q^2}{1-q^{2(\lambda-1)}}\,
\| z^n \|_\lambda^2=\int \limits_0^1(q^2y;q^2)_n \cdot
 y^{\lambda-2}d_{q^2}y=B_{q^2}(\lambda-1,n+1)=
$$
$$
=\frac{\Gamma_{q^2}(\lambda-1)\cdot\Gamma_{q^2}(n+1)}
{\Gamma_{q^2}(n+\lambda)}=\frac{1-q^2}{1-q^{2(\lambda-1)}}
\frac{(q^2;q^2)_n}{(q^{2\lambda};q^2)_n}.\eqno \square
$$

\bigskip
Равенство \eqref{Klimek-Lesn} получено Климеком и Лесневским в
\cite{KlimLesn}, и его классический $q\to 1$ предел хорошо известен
\cite[стр. 4]{HedKorZhu}.

\medskip

Пусть $P_\lambda$ -- ортогональный проектор в гильбертовом пространстве
$L^2(d\nu_\lambda)_q$ на подпространство $L^2_a(d\nu_\lambda)_q$. Каждому
формальному ряду $\mathscr{K}$ переменной $z\otimes\zeta^*$ сопоставим
интегральный оператор из $\mathscr{D}(\mathbb{D})_q$ в
$\mathscr{D}(\mathbb{D})'_q$: $Kf(z)=\int\limits_{\mathbb{D}_q}
\mathscr{K}(z,\zeta^*)f(\zeta)d\nu_\lambda(\zeta)$.

 Как следует из леммы \ref{l1.7.1} и $q$-биномиальной формулы \eqref{qbinom},
 $$
P_\lambda\,f\,=\, \int\limits_{\mathbb{D}_q}
\mathscr{K}_\lambda(z,\zeta^*)f(\zeta)d\nu_\lambda(\zeta),\qquad f\in
\mathscr{D}(\mathbb{D})_q,
 $$
где
\begin{equation}\label{q-Bergman_sl_2}
\mathscr{K}_\lambda\,=\,
 \sum\limits_{m=0}^\infty \frac{(q^{2\lambda};q^2)_m}
{(q^2;q^2)_m}\,(z\otimes\zeta^*)^m\,=\,
 (q^{2\lambda}\,z\otimes
\zeta^*;q^2)_\infty\;(z\otimes \zeta^*; q^2)_\infty^{-1}
\end{equation}
($\mathscr{K}_\lambda$ является $q$-аналогом взвешенного ядра Бергмана
$(1-z\overline{\zeta})^{-\lambda}$, см. \cite[стр. 5]{HedKorZhu})).

\bigskip

Элемент $f\in\mathscr{D}(\mathbb{D})_q'$ назовем ограниченной функцией в
квантовом круге, \IND{функция ! ограниченная в квантовом круге} если
обобщенный линейный оператор $T_F(f)$ ограничен, и положим
$\|f\|_\infty=\|T_F(f)\|$. Ограниченные функции в квантовом круге образуют
алгебру фон Неймана $L^\infty(\mathbb{D})_q$, естественно изоморфную
алгебре всех ограниченных линейных операторов в пространстве представления
$T_F$, см. \itemiii\ \ref{finite_sl_2}.

Перейдем к операторам Теплица. Как следует из определения класса $S_2$
операторов Гильберта-Шмидта, см. \itemiii\ \ref{S_p},
\begin{equation}\label{like_sigma_2}
L^2(d\nu_\lambda)_q\,=\,\left\{\psi\in\mathscr{D}(\mathbb{D})_q'\;\left|\;
T_F(\psi)T_F\left(y^{\frac{\lambda-1}2}\right)\in S_2\right.\right\},
\end{equation}
\begin{equation}\label{norm_like_sigma_2}
\|\psi\|_\lambda^2\,=\,(1-q^{2(\lambda-1)})\,
\left\|T_F(\psi)T_F\left(y^{\frac{\lambda-1}2}\right)\right\|_2.
\end{equation}

Рассмотрим ограниченную функцию $f\in L^\infty(\mathbb{D})_q$ в квантовом
круге. Так как $S_2$ является двусторонним идеалом алгебры ограниченных
линейных операторов, то
$$
\|f\psi\|_\lambda\le\|f\|_\infty\,\|\psi\|_\lambda,\qquad
\psi\in\mathscr{D}(\mathbb{D})_q.
$$
Значит, линейный оператор $\psi\mapsto P_\lambda(f\psi)$ допускает
продолжение по непрерывности до ограниченного оператора в
$L^2(d\nu_\lambda)_q$ с нормой, не превосходящей $\|f\|_\infty$. Сужение
этого линейного оператора на подпространство $L^2_a(d\nu_\lambda)_q$
называется оператором Теплица \IND{оператор ! Теплица} с (контравариантным)
символом $f$ и обозначается $\mathscr{T}^{(\lambda)}_f$, либо $\hat{f}$.
Очевидно, $\|\hat{f}\|\le\|f\|_\infty$ и $\hat{f}^*=\widehat{f^*}$.

\bigskip

\begin{example} Покажем, что при всех $k,n\in \mathbb{Z}_+$
$$\widehat{y^k}\; z^n\;=\; q^{2kn}
\frac{(q^{2(\lambda-1)};q^2)_k}{(q^{2(\lambda+n)};q^2)_k} \,z^n.$$
Действительно, для любой функции $f(y)\in L^\infty(\mathbb{D})_q$
\begin{equation}\label{Toeplitz_y}
\widehat{f}\,z^n\,=\,\frac{(f(y)z^n,z^n)_\lambda}{(z^n,z^n)_\lambda}\,z^n
\,=\,\frac{\int \limits_0^1\,f(q^{2n}y)\,(q^2y;q^2)_n \;
 y^{\lambda-2}d_{q^2}y}{\int \limits_0^1\,(q^2y;q^2)_n \;
 y^{\lambda-2}d_{q^2}y}\,z^n.
\end{equation}
Рассуждая так же, как при доказательстве леммы \ref{l1.7.1}, получаем
$$
\widehat{f}\,z^n\,=\,
 q^{2kn}\frac{B_{q^2}(\lambda+k-1,n+1)}{B_{q^2}(\lambda-1,n+1)}\,z^n=
 q^{2kn}\frac{\Gamma_{q^2}(\lambda+k-1)\Gamma_{q^2}(\lambda+n)}
 {\Gamma_{q^2}(\lambda+k+n)\Gamma_{q^2}(\lambda-1)}\,z^n.
$$
Остается воспользоваться определением $q$-гамма-функции
$$
\frac{\Gamma_{q^2}(\lambda+k-1)\Gamma_{q^2}(\lambda+n)}
 {\Gamma_{q^2}(\lambda+k+n)\Gamma_{q^2}(\lambda-1)}\,=\,
 \frac{(q^{2(\lambda+k+n)};q^2)_\infty\,(q^{2(\lambda-1)};q^2)_\infty}
 {(q^{2(\lambda+k-1)};q^2)_\infty\,(q^{2(\lambda+n)};q^2)_\infty}
 \,=\, \frac{(q^{2(\lambda-1)};q^2)_k}{(q^{2(\lambda+n)};q^2)_k}.
$$
\end{example}

\smallskip
\begin{example}\label{c1.7.2} Покажем, что
\begin{equation}\label{1.7.2}
\widehat{z}^*\widehat{z}=q^2 \widehat{z}\widehat{z}^*+
 1-q^2+q^{2(\lambda-1)}\, \frac{1-q^2}{1-q^{2(\lambda-1)}}\,
 (1-\widehat{z}\widehat{z}^*) (1-\widehat{z}^*\widehat{z}).
\end{equation}
Действительно,
\begin{equation}\label{for_covariant_symbols}\widehat{z}\,z^n\,=\,z^{n+1},\qquad
\widehat{z}^*\,z^n\,=\,
\begin{cases}
\frac{1-q^{2n}}{1-q^{2(\lambda-1)+2n}}z^{n-1},\; &\; n \neq 0,\\
 0,\; &\; n=0.                          \end{cases}
\end{equation}
Следовательно,
$$
(1-\widehat{z}\widehat{z}^*)^{-1}\,z^m\,=\, \frac{q^{-2m}-
 q^{2(\lambda-1)}}{1-q^{2(\lambda-1)}}z^m,\quad
(1-\widehat{z}^*\widehat{z})^{-1}\,z^m\,=\,
 \frac{q^{-2m-2}-q^{2(\lambda-1)}}{1-q^{2(\lambda-1)}}z^m.
$$
Значит,
$$
(1-\widehat{z}\widehat{z}^*)^{-1}=q^2(1-\widehat{z}^*\widehat{z})^{-1}-
q^{2(\lambda-1)}\dfrac{1-q^2}{1-q^{2(\lambda-1)}},
$$
чторавносильно \eqref{1.7.2}. Коммутационное соотношение \eqref{1.7.2} было
получено в \cite{KlimLesn}.
\end{example}

\medskip \begin{example}\label{Toepltz_symbol}
Пусть $f \in {\rm Pol}(\mathbb{C})_q$. Если
 $f=\sum\limits_{j,k \in \mathbb{Z}_+}\, a_{j,k}z^{*\,j}z^k$, то\\
 $\widehat{f}=\sum\limits_{j,k \in \mathbb{Z}_+}\,
 a_{j,k}\widehat{z}^{*\,j}\widehat{z}^k$.
\end{example}

\begin{example}\label{Toeplitz-f_0}
\begin{equation}\label{f_0_Toeplitz}
\widehat{f_0}: z^j \mapsto \begin{cases} 1-q^{2(\lambda-1)},\;&\; j=0,\\
                       0,\;&\; j \ne 0 \end{cases}.
\end{equation}
\end{example}

\bigskip
Очевидным образом вводятся операторы Ханкеля с ограниченными символами:
\IND{оператор ! Ханкеля}
$$
\mathcal{H}^{(\lambda)}_f:L^2_a(d\nu_\lambda)_q\to
L^2(d\nu_\lambda)_q\ominus
L^2_a(d\nu_\lambda)_q,\qquad\mathcal{H}^{(\lambda)}_f:\,\psi\to
(I-\,P_\lambda)(f\psi).
$$
Отметим, что в классическом случае $q=1$ изучению операторов Теплица и
операторов Ханкеля во взвешенных пространствах Бергмана посвящена обширная
литература.

\subsubsection{Ковариантные символы линейных
операторов.}\label{covariant_symbols} Начнем с нестрогих наводящих
соображений. Пусть $A$-- ограниченный линейный оператор в
$L^2_a(d\nu_\lambda)_q$ и
$$ A=\sum\limits_{i,j=0}^\infty a_{ij}\widehat{z}^i \widehat{z}^{*j}=
 \sum\limits_{i,j=0}^\infty b_{ij} \widehat{z}^{*j} \widehat{z}^i.
$$
Функции в квантовом круге $\sum\limits_{i,j=0}^\infty a_{ij} z^i z^{*j}$,
$\sum\limits_{i,j=0}^\infty b_{ij} z^{*j} z^i$ хотелось бы назвать
ковариантным и контравариантным символами оператора $A$ соответственно.
Отображение, сопоставляющее контравариантному символу линейного опрератора
его ковариантный символ-- это преобразование Березина.

К сожалению, даже в классическом случае $q=1$ эти построения трудно,
уточнив, сделать строгими, что вынуждает использовать менее наглядные
определения ковариантного и контравариантного символов. Их связь с
приведенными выше эвристическими соображениями описывают предложение
\ref{p5.6.6} и пример \ref{Toepltz_symbol}.

\bigskip
 Рассмотрим градуированный $U_q\mathfrak{sl}_2$-модуль
$\mathbb{C}[z]_{q,\lambda}$, введенный в \itemiiiе \ref{bundles_sl_2},
двойственный градуированный $U_q\mathfrak{sl}_2$-модуль
$\mathbb{C}[z]_{q,\lambda}^*$ и $U_q\mathfrak{sl}_2$-модульную алгебру
${\rm End}(\mathbb{C}[z]_{q,\lambda})$ всех линейных операторов в
$\mathbb{C}[z]_{q,\lambda}$. Очевидно,
$\mathbb{C}[z]_{q,\lambda}\otimes\mathbb{C}[z]_{q,\lambda}^*\hookrightarrow
{\rm End}(\mathbb{C}[z]_{q,\lambda})$. Подалгебру
$$
\underline{{\rm End}}(\mathbb{C}[z]_{q,\lambda})\stackrel{\rm def}{=}
\mathbb{C}[z]_{q,\lambda}\otimes\mathbb{C}[z]_{q,\lambda}^*.
$$
будем называть алгеброй финитных линейных операторов в
$\mathbb{C}[z]_{q,\lambda}$. \IND{алгебра ! финитных линейных операторов в
$\mathbb{C}[z]_{q,\lambda}$}

\begin{lemma}\label{p5.1.1} Если $f \in \mathscr{D}(\mathbb{D})_q$, то
интеграл
 $I_{m,j}\,=\,\int\limits_{\mathbb{D}_q}\, z^{*\,m}\,f z^j\, d \nu_ \lambda$
равен нулю для почти всех $m,j \in{\mathbb Z}_+$ (т.е. для всех, за
исключением конечного числа).
\end{lemma}

 {\bf Доказательство.} Для любого
$f\in \mathscr{D}(\mathbb{D})_q$ равенство $z^{*N}f=fz^N=0$ имеет место при
всех достаточно больших $N \in{\mathbb N}$. Значит, $I_{mj}=0$ при
$\max(m,j)\ge N$. \hfill $\square$

\medskip

Из леммы следует, что оператор Теплица с финитным символом является
финитным линейным оператором в $\mathbb{C}[z]_{q,\lambda}$.

\begin{proposition}\label{p5.1.2}
Линейный оператор
$$
\mathscr{T}^{(\lambda)}:\,\mathscr{D}(\mathbb{D})_q\to\underline{{\rm
End}}(\mathbb{C}[z]_{q,\lambda}),\qquad\mathscr{T}^{(\lambda)}:\,f\mapsto
\widehat{f},
$$
является морфизмом $U_q\mathfrak{sl}_2$-модулей.
\end{proposition}

\smallskip

{\bf Доказательство.} Как показано в \itemiiiе \ref{bundles_sl_2},
${\mathbb C}[z]_{q,\lambda}$ является $U_q\mathfrak{sl}_2$-модульным
$\mathbb{C}[z]_q$-модулем и эрмитова форма $(\cdot,\cdot)_\lambda$
инвариантна. Следовательно, линейное отображение
$$
\mathscr{D}(\mathbb{D})_q\otimes\mathbb{C}[z]_{q,\lambda}\to
\mathbb{C}[z]_{q,\lambda},\qquad f\otimes\psi\mapsto P_{q,\lambda}(f\psi),
$$
является морфизмом $U_q\mathfrak{sl}_2$-модулей. Значит, отвечающий ему
элемент тензорного произведения
$$
\mathbb{C}[z]_{q,\lambda}\otimes(\mathscr{D}(\mathbb{D})_q\otimes
\mathbb{C}[z]_{q,\lambda})^*\cong\underline{{\rm
End}}(\mathbb{C}[z]_{q,\lambda})\otimes\mathscr{D}(\mathbb{D})_q^*
$$
$U_q\mathfrak{sl}_2$-инвариантен. \hfill $\square$

Изучение сопряженного к $\mathscr{T}^{(\lambda)}$ линейного оператора начнем
с описания сопряженных векторных пространств.

Пусть $\mathbb{C}[[z]]_{q,\lambda}$ -- векторное пространство формальных
степенных рядов переменной $z$ с топологией покоэффициентной сходимости и
действием $U_q\mathfrak{sl}_2$, полученным продолжением по непрерывности с
$\mathbb{C}[z]_{q,\lambda}$. Линейные операторы из
$\mathbb{C}[z]_{q,\lambda}$ в $\mathbb{C}[[z]]_{q,\lambda}$ будем называть
обобщенными линейными операторами в $\mathbb{C}[z]_{q,\lambda}$ и
использовать обозначение
$$\overline{{\rm End}}(\mathbb{C}[z]_{q,\lambda})\,=\,
{\rm Hom}(\mathbb{C}[z]_{q,\lambda},\mathbb{C}[[z]]_{q,\lambda}).
$$
Очевидно, $\overline{{\rm End}}(\mathbb{C}[z]_{q,\lambda})$ является
$U_q\mathfrak{sl}_2$-модульным $\underline{{\rm
End}}(\mathbb{C}[z]_{q,\lambda})$-бимодулем.

Следуя традиции, введем в рассмотрение инвариантный интеграл
$\mathrm{tr}_q$ на $U_q\mathfrak{sl}_2$-модульной алгебре $\underline{{\rm
End}}(\mathbb{C}[z]_{q,\lambda})$ с помощью равенства
$$ S^2(\xi)\; =\; K^{-1} \xi K,\qquad \xi \in U_q\mathfrak{sl}_2,
$$
 см. \itemiiiы \ref{IntStar},
\ref{inv_int_sl_2}.

Докажем лемму, которая позволит отождествить $\overline{{\rm
End}}(\mathbb{C}[z]_{q,\lambda})$ с $U_q\mathfrak{sl}_2$-модулем,
сопряженным к $\underline{{\rm End}}(\mathbb{C}[z]_{q,\lambda})$, а
$\mathscr{D}(\mathbb{D})_q'$ -- с $U_q\mathfrak{sl}_2$-модулем, сопряженным
к $\mathscr{D}(\mathbb{D})_q$.

\begin{lemma}\label{to_conjugate_new}
1. Линейное отображение \\
$i_\lambda:\overline{\mathrm{End}}(\mathbb{C}[z]_{q,\lambda})\to
\underline{\mathrm{End}}(\mathbb{C}[z]_{q,\lambda})^*$, определяемое
равенством
$$
i_\lambda(A)(X)\,=\,
\frac{1-q^2}{1-q^{2(\lambda-1)}}\cdot {\rm tr}_q(A\,X),
\qquad X\in\underline{{\rm End}}(\mathbb{C}[z]_{q,\lambda}),
$$
является изоморфизмом $U_q\mathfrak{sl}_2$-модулей.

2. Линейное отображение
 $i: \mathscr{D}(\mathbb{D})_q' \to (\mathscr{D}(\mathbb{D})_q)^*$,
 определяемое равенством
$$
 i(f)(\psi)\,=\,\int\limits_{\mathbb{D}_q}\,f\psi\,d\nu,\qquad
  \psi \in \mathscr{D}(\mathbb{D})_q,
$$
является изоморфизмом $U_q\mathfrak{sl}_2$-модулей.
\end{lemma}

{\bf Доказательство.} Линейный функционал
\begin{equation}\label{mu_lambda_sl_2}
\int A\,d\mu_\lambda=\frac{1-q^2}{1-q^{2(\lambda-1)}}\cdot{\rm tr}_q\,A
\end{equation}
является инвариантным интегралом на $\underline{{\rm
End}}(\mathbb{C}[z]_{q,\lambda})$, а линейный функционал
$\int\limits_{\mathbb{D}_q}\,f\,d\nu$ -- инвариантным интегралом на
$\mathscr{D}(\mathbb{D})_q$. Остается вспомнить определение двойственного
$U_q\mathfrak{sl}_2$-модуля и воспользоваться формулой интегрирования по
частям \eqref{by_parts}. \hfill $\square$

\bigskip С помощью  доказанной леммы введем в рассмотрение
 сопряженный линейный оператор
$$
(\mathscr{T}^{(\lambda)})^*:\,
\overline{{\rm End}}(\mathbb{C}[z]_{q,\lambda}) \to
\mathscr{D}(\mathbb{D})_q',
$$
являющийся морфизмом
$U_q\mathfrak{sl}_2$-модулей.

Обобщенную функцию $(\mathscr{T}^{(\lambda)})^*\,A$ назовем ковариантным
символом обобщенного линейного оператора \IND{ковариантный символ
обобщенного линейного оператора} $A\in\overline{{\rm
End}}(\mathbb{C}[z]_{q,\lambda})$. Это понятие использовалось Березиным в
другом контексте, и нам предстоит установить соответствие определений.
Прежде всего, получим интегральное представление матричных элементов
оператора Теплица.

\begin{lemma}\label{p5.1.5}
Пусть $\widehat{f}:\mathbb{C}[z]_{q,\lambda}\to\mathbb{C}[z]_{q,\lambda}$
-- оператор Теплица с символом $f\in\mathscr{D}(\mathbb{D})_q$ и
$\{\widehat{f}_{m,j}\}$ -- его матричные элементы в базисе $\{z^j\}_{j\in
\mathbb{Z}_+}$. Тогда
\begin{equation}\label{5.1.2}
\widehat{f}_{mj}=\frac{1-q^{2(\lambda-1)}}{1-q^2}
\int\limits_{\mathbb{D}_q}P_{z,m,j}\,f(z)\,d\nu(z),
\end{equation}
где
\begin{equation}\label{5.1.3}
P_{z,m,j}=
\frac{(q^{2\lambda};q^2)_m}{(q^2;q^2)_m}q^{2j}z^j(1-zz^*)^{\lambda}z^{*m}.
\end{equation}
\end{lemma}

\smallskip

 {\bf Доказательство.} Так как
$$(z^m,z^l)_\lambda=\frac{(q^2;q^2)_m}{(q^{2\lambda};q^2)_m}\delta_{ml},
\qquad m,l \in{\mathbb Z}_+,$$ то
$$\widehat{f}_{mj}=
\frac{(f\,z^j,z^m)_\lambda}{(z^m,z^m)_\lambda}=
\frac{(q^{2\lambda};q^2)_m}{(q^2;q^2)_m}\int
\limits_{\mathbb{D}_q}z^{*m}\,f\,z^j\,d \nu_\lambda=$$
$$=\frac{1-q^{2(\lambda-1)}}{1-q^2}\cdot
 \frac{(q^{2\lambda};q^2)_m}
{(q^2;q^2)_m}\int
 \limits_{\mathbb{D}_q}z^{*m} \,f\,z^j(1-zz^*)^{\lambda}d
\nu.$$ Как следует из явного вида инвариантного интеграла,
$$\int
\limits_{\mathbb{D}_q}f_1f_2(1-zz^*)d \nu(z)=\int
\limits_{\mathbb{D}_q}f_2f_1(1-zz^*)d \nu(z)$$ для любых $f_1\in
\mathscr{D}(\mathbb{D})_q^\prime$, $f_2\in \mathscr{D}(\mathbb{D})_q$.
Значит,
\begin{equation}\label{5.1.4}
  \widehat{f}_{mj}=\frac{1-q^{2(\lambda-1)}}{1-q^2}\cdot
   \frac{(q^{2\lambda};q^2)_m}{(q^2;q^2)_m}
   \int\limits_{\mathbb{D}_q}(1-zz^*)z^j(1-zz^*)^{\lambda-1}
   z^{*m}\,f\,d \nu.
\end{equation}
 Остается воспользоваться коммутационным соотношением
$$(1-zz^*)z=q^2z(1-zz^*). \eqno \square$$

 Матрица $P_z=(P_{z,m,
 j})_{m,j \in{\mathbb Z}_+}$ является
$q$-аналогом матрицы ортогонального проектора на одномерное
подпространство, порожденное собственным вектором оператора, сопряженного к
умножению на $z$, то есть вектором переполненной системы, используемой
Березиным \cite{Ber_CMP}.


\begin{proposition}\label{Berezin_2.1}
Пусть $A \in \overline{\operatorname{End}}(\mathbb{C}[z]_{q,\lambda})$ и
$\{a_{mj}\}$ -- матричные элементы $A$ в базисе
$\{z^j\}_{j\in\mathbb{Z}_+}:\; A z^j=\sum\limits_{m=0}^\infty a_{mj}z^m$.
Тогда
\begin{equation}\label{Berezin_2.2}
(\mathscr{T}^{(\lambda)})^*(A) =\operatorname{tr}_q(A\, P_z),
\end{equation}
где $\operatorname{tr}_q(A\,P_z)\;=\; \sum\limits_{j,m=0}^\infty
a_{jm}P_{z,j,m}q^{-2j}$ и ряд сходится в топологии пространства
$\mathscr{D}(\mathbb{D})_q'$.
\end{proposition}

{\bf Доказательство.} Как следует из \eqref{5.1.2},
\begin{multline*}
\int\limits\,A\,
\left(\mathscr{T}^{(\lambda)}\psi\right)d\mu_\lambda=
\frac{1-q^2}{1-q^{2(\lambda-1)}}
\operatorname{tr}_q\left(A\,\mathscr{T}^{(\lambda)}\psi\right)=
\\ =\sum_{j,m=0}^\infty a_{jm}\left(\int\limits_{\mathbb{D}_q}P_{z,m,j}\psi
d\nu\right)q^{-2j}=\int\limits_{\mathbb{D}_q}\left(\sum_{j,m=0}^\infty
a_{jm}P_{z,m,j}q^{-2j}\right)\psi d\nu.\tag*{$\square$}
\end{multline*}

\begin{example}\label{cov_symbol_f_0_hat} Рассмотрим оператор $A$,
 для которого
 $Az^j=\begin{cases}1,\;&\;j=0,\\ 0, \;&\; j \ne 0\end{cases}$.
 Его ковариантный символ равен $(1-zz^*)^\lambda$.
\end{example}

\medskip Важно отметить, что в работе Березина \cite{Ber_CMP}
ковариантный символ линейного оператора вводится равенством, аналогичным
\eqref{Berezin_2.2}.

\bigskip В оставшейся части \itemiiiа описаны основные свойства
ковариантных символов обобщенных линейных операторов.

 Напомним, что операторы
$\widehat{z}$, $\widehat {z}^*$ действуют в градуированном векторном
пространстве ${\mathbb C}[z]_{q,\lambda}$, причем ${\rm
deg}(\widehat{z})=+1$, ${\rm deg}(\widehat{z}^*)=-1$. Следовательно, для
любой числовой матрицы $(a_{ij})_{i,j \in{\mathbb Z}_+}$ ряд
\begin{equation}\label{5.6.1}
A=\sum\limits_{i,j=0}^\infty a_{ij}\widehat{z}^i \widehat{z}^{*j}
\end{equation}
 сходится в сильной операторой топологии пространства
$$\overline{\operatorname{End}}(\mathbb{C}[z]_{q,\lambda})
\,=\,{\rm Hom}_{\mathbb C}({\mathbb C}[z]_{q,\lambda},{\mathbb
C}[[z]]_{q,\lambda}).$$

Аналогично доказывается, что для любой матрицы $(c_{ij})_{i,j \in{\mathbb
Z}_+}$ ряд
\begin{equation}\label{5.6.1_infty}
f=\sum\limits_{i,j=0}^\infty c_{ij} z^i z^{*j}
\end{equation}
 сходится в топологическом векторном пространстве
$\mathscr{D}(\mathbb{D})_q '$, канонически изоморфном топологическому
векторному пространству обобщенных линейных операторов в $\mathcal{H}$.

\medskip Из  \eqref{for_covariant_symbols} вытекает следующее
утверждение, см. \cite{KlimLesn}.

\begin{lemma}\label{p5.6.1} Если
$A=\sum\limits_{i,j=0}^\infty a_{ij}\widehat{z}^i \widehat{z}^{*j}$, то
$$A\,z^n=\displaystyle \sum \limits_{m=0}^\infty b_{mn}z^m,\;\;\text{где}\;\; b_{mn}=\displaystyle\sum
\limits_{j=0}^{\min(m,n)}
\frac{\textstyle(q^{2n};q^{-2})_{n-j}}{\textstyle
 (q^{2(\lambda+n-1)};q^{-2})_{n-j}}a_{m-j,n-j}.$$
\end{lemma}

\medskip

Из этой леммы следует

\begin{proposition}\label{c5.6.2}
Для любого $A\in \overline{\operatorname{End}}(\mathbb{C}[z]_{q,\lambda})$
существует и единственно разложение (\ref{5.6.1}).
\end{proposition}

Аналогично доказывается
\begin{proposition}\label{c5.6.2_infty} Для
любой обобщенной функции $f\in \mathscr{D}(\mathbb{D})_q'$ существует и
единственно разложение (\ref{5.6.1_infty}). Линейное отображение $
 (c_{ij})_{i,j \in{\mathbb Z}_+} \mapsto
\sum\limits_{i,j=0}^\infty c_{ij} z^i z^{*j} $ является изоморфизмом
векторного пространства матриц $(c_{ij})_{i,j \in{\mathbb Z}_+}$ с
топологией поэлементной сходимости и топологического векторного пространства
$\mathscr{D}(\mathbb{D})_q'$.
\end{proposition}

\begin{example}\label{f_0_hat_expansin}
Рассмотрим оператор Теплица
$$
\widehat{f_0}\,z^j=\;
\left\{\begin{array}{cll}1-q^{2(\lambda-1)} &,& j=0 \\
0 &,& j\ne 0\end{array}\right.
$$
с символом $f_0$, см. \eqref{f_0_Toeplitz}. Докажем равенство
\begin{equation}\label{5.6.3}
\widehat{f}_0=\left(1-q^{2(\lambda-1)}\right)
\sum_{k=0}^\infty\frac{(q^{-2\lambda};q^2)_k}{(q^2;q^2)_k}\,q^{2\lambda
k}\widehat{z}^k\widehat{z}^{*k}.
\end{equation}

 Перейдем от равенства операторов к равенству их матричных элементов
в базисе $\{z^n \}_{n \in{\mathbb Z}_+}$. Внедиагональные элементы равны
нулю. Сравнение диагональных элементов дает
\begin{equation}\label{5.6.4}
  \sum_{k=0}^j \frac{(q^{-2\lambda};q^2)_k}{(q^2;q^2)_k}\cdot
\frac{(q^{2j};q^{-2})_k}{(q^{2(\lambda+j-1)};q^{-2})_k}\cdot
 q^{2\lambda k}=\delta_{j,0}.
\end{equation}
 Достаточно рассмотреть случай $j>0$. Умножая обе части равенства
(\ref{5.6.4}) на
$\frac{\textstyle(q^{2(\lambda+j-1)};
q^{-2})_j}{\textstyle(q^{2j};q^{-2})_j},$
получаем
$$\sum_{k=0}^j \frac{(q^{-2\lambda};q^2)_k}{(q^2;q^2)_k}\cdot
 \frac{(q^{2\lambda};q^2)_{j-k}}{(q^2;q^2)_{j-k}}\cdot q^{2\lambda k}=0.$$
 То есть,
$$\sum_{k+m=j} \frac{(q^{-2\lambda};q^2)_k}{(q^2;q^2)_k}\cdot
 q^{2\lambda k}\cdot \frac{(q^{2\lambda};q^2)_m}{(q^2;q^2)_m}=0.$$
 Остается ввести в рассмотрение $q$-биномиальные ряды (см. \eqref{qbinom})
$$a(t)=\sum_{k=0}^\infty
 \frac{(q^{-2\lambda};q^2)_k}{(q^2;q^2)_k}\cdot q^{2\lambda k}\cdot
 t^k=\frac{(t;q^2)_\infty}{(q^{2\lambda}t;q^2)_\infty},$$
$$
b(t)=\sum_{m=0}^\infty\frac{(q^{2\lambda};q^2)_m}{(q^2;q^2)_m}
 \cdot
t^m=\frac{(q^{2\lambda}t;q^2)_\infty}{(t;q^2)_\infty},
$$
и заметить, что $a(t)\cdot b(t)=1$.
\end{example}

\begin{example}\label{f_0_expansin} Аналогично доказывается,
что в $\mathscr{D}(\mathbb{D})_q'$
\begin{equation}\label{f0_exp}
  f_0 =  \sum\limits_{j=0}^\infty \frac{(-1)^j
q^{j(j-1)}}{(q^2; q^2)_j} z^j z^{*j}.
\end{equation}
\end{example}

\begin{example} Напомним, что $y=1-zz^*$ и что $y^\lambda \in
\mathscr{D}(\mathbb{D})_q'$ при всех $\lambda$. Получим разложение
\begin{equation}\label{5.6.5}
 (1-zz^*)^\lambda\,=\,\sum_{k=0}^\infty
 \frac{(q^{-2\lambda};q^2)_k}{(q^2;q^2)_k}q^{2\lambda k}z^kz^{*k}
\end{equation}
в $\mathscr{D}(\mathbb{D})_q'$. Переходя к операторам фоковского
представления и приравнивая их матричные элементы в базисе $\{T_F(z^m)e_0
\}_{m\in \mathbb{Z}_+}$ (достаточно рассматривать диагональные матричные
элементы), получаем
$$
\sum\limits_{k=0}^j \frac{(q^{-2\lambda};q^2)_k}{(q^2;q^2)_k}\cdot
(q^{2j};q^{-2})_k \cdot q^{2\lambda k}=q^{2\lambda j},\qquad
j\in\mathbb{Z}_+
$$
то есть
$$
\sum\limits_{k+m=j}\frac{(q^{-2\lambda};q^2)_k}{(q^2;q^2)_k}
\cdot q^{2\lambda k}\cdot\frac1{(q^2;q^2)_m}=
\frac{q^{2\lambda j}}{(q^2;q^2)_j},\qquad j \in \mathbb{Z}_+.
$$
Остается ввести в рассмотрение производящие функции
$$
a(t)=\frac{(t;q^2)_\infty}{(q^{2\lambda}t;q^2)_\infty},\qquad b(t)=
\frac1{(t;q^2)_\infty},\qquad c(t)=\frac1{(q^{2\lambda}t;q^2)_\infty}
$$
и применить $q$-биномиальное разложение \eqref{qbinom} к обеим частям
очевидного равенства $a(t)b(t)=c(t)$.
\end{example}

\bigskip Покажем, что элемент $\widehat{f_0}$ порождает
$\overline{\operatorname{End}}(\mathbb{C}[z]_{q,\lambda})$ как
топологический $U_q \mathfrak{sl}_2$-модуль.

\begin{lemma}\label{p5.6.4}
$U_q \mathfrak{sl}_2\,
\widehat{f_0}=\underline{\operatorname{End}}(\mathbb{C}[z]_{q,\lambda})$.
\end{lemma}

\smallskip

{\bf Доказательство.} Используя \eqref{for_covariant_symbols},
\eqref{f_0_Toeplitz}, нетрудно получить равенство
$$
\widehat{z}^i\widehat{f}_0\widehat{z}^{*n}:z^j\mapsto
\left(1-q^{2(\lambda-1)}\right)\;
\frac{\textstyle(q^2;q^2)_n}{\textstyle(q^{2\lambda};q^2)_n}\;
\delta_{jn}z^i,\qquad i,j,n\in\mathbb{Z}_+.
$$
Значит, линейные операторы
$\{\widehat{z}^i\widehat{f}_0\widehat{z}^{*n}\}_{i,n\in\mathbb{Z}_+}$
порождают векторное пространство
$\underline{\operatorname{End}}(\mathbb{C}[z]_{q,\lambda})$. Остается
показать, что все они принадлежат $U_q\mathfrak{sl}_2\widehat{f}_0$. Для
этого достаточно повторить доказательство предложения \ref{t2.3.9}, заменив
образующие ${z}$, ${z}^*$, ${f}_0$ на $\widehat{z}$, $\widehat{z}^*$,
$\widehat{f}_0$. Действительно, формулы, описывающие действие $K^{\pm 1}$,
$E$, $F$ на $\widehat{f}_0$, отличаются от формул \eqref{lemma_EFf0} лишь
числовыми множителями:
$$
K^{\pm 1}\widehat{f_0}=\widehat{f_0},\quad
E\widehat{f_0}=c'\widehat{z}\cdot\widehat{f_0},\quad
F\widehat{f_0}=c''\widehat{f_0}\cdot\widehat{z}^*,\qquad c',c''\ne 0.
$$
Эти соотношения вытекают из предложения \ref{p5.1.2}, равенств
\eqref{lemma_EFf0} и равенств $ {\mathbb
C}\widehat{zf_0}=\mathbb{C}\widehat{z}\widehat{f_0}$,
$\mathbb{C}\widehat{f_0z^*}=\mathbb{C}\widehat{f_0}\widehat{z^*},$ которые
легко следуют из
$$
{\rm Im}\,\widehat{zf_0}=\mathbb{C}z,\quad{\rm
Im}\,\widehat{f_0z^*}={\mathbb C}\cdot1,\quad{\rm Ker}\,\widehat{f_0}={\rm
Ker}\,\widehat{zf_0}=(\mathbb{C}\cdot1)^\perp. \eqno \square
$$

\medskip

\begin{proposition}\label{p5.6.6} Для любой числовой матрицы
 $(a_{jk})_{j,k \in \mathbb{Z}_+}$ ковариантный символ оператора
$\sum\limits_{j,k=0}^\infty a_{jk}\widehat{z}^j \widehat{z}^{*k}$ равен
$\sum\limits_{j,k =0}^\infty a_{jk}z^jz^{*k}$.
\end{proposition}

\smallskip

 {\bf Доказательство.} Нам известен ковариантный символ оператора Теплица
 $\widehat{f_0}$, см. пример \ref{cov_symbol_f_0_hat}.
Используя \eqref{5.6.3}, \eqref{5.6.5}, получаем доказываемое утверждение в
этом частном случае.

Остается воспользоваться леммой \ref{p5.6.4} и тем, что линейный оператор
$\left(\mathscr{T}^{(\lambda)}\right)^*$ является морфизмом
$U_q\mathfrak{sl}_2$-модулей. \hfill $\square$

\medskip

\begin{corollary}\label{c5.6.7}
 Отображение
$
\left(\mathscr{T}^{(\lambda)}\right)^*:\;
\overline{\operatorname{End}}(\mathbb{C}[z]_{q,\lambda})\to
\mathscr{D}(\mathbb{D})_q',
$
сопоставляющее оператору его ковариантный символ, биективно.
\end{corollary}

\subsubsection{Преобразование Березина.}\label{Berezin_transform}

Рассмотрим линейный оператор $\mathscr{B}_\lambda$ из
$\mathscr{D}(\mathbb{D})_q$ в $\mathscr{D}(\mathbb{D})'_q$, определяемый
равенством $\mathscr{B}_\lambda\,=\,\left(\mathscr{T}^{(\lambda)}\right)^*
\mathscr{T}^{(\lambda)}$. Он сопоставляет финитной функции $f$ ковариантный
символ оператора Теплица $\hat{f}$ с контравариантным символом $f$.
Оператор $\mathscr{B}_\lambda$ называется преобразованием Березина,
\IND{преобразование ! Березина} \cite[стр. 566]{UU}. Из полученных в
предыдущем
\itemiiiе результатов следует, что $\mathscr{B}_\lambda$ -- морфизм
$U_q\mathfrak{sl}_2$-модулей.

\begin{example}\label{ex_Berezin}
 $\mathscr{B}_\lambda\,f_0=(1-q^{2(\lambda-1)})(1-zz^*)^\lambda$,
см. примеры \ref{Toeplitz-f_0}, \ref{cov_symbol_f_0_hat}.
\end{example}

\medskip

В книге \cite[стр. 29]{HedKorZhu}, посвященной теории пространств Бергмана
в обычном круге, преобразование Березина определяется равенством
\begin{equation}\label{classical_Bergman}
\mathscr{B}_\lambda f(z)\,=\,(\lambda-1)\,\int\limits_{\mathbb{D}}\,
\left(\frac{(1-|z|^2)(1-|\zeta|^2)}{|1-z\bar{\zeta}|^2}\right)^\lambda\,
f(\zeta)\,d\nu(\zeta)
\end{equation}
(авторы используют параметр $\alpha=\lambda-2$ вместо параметра Уоллака
$\lambda$). Покажем, что в случае квантового круга преобразование Березина
также является интегральным оператором, и найдем его ядро. Тем самым будет
получен $q$-аналог равенства \eqref{classical_Bergman}.

\begin{proposition}\label{p5.2.5}
$ (\mathscr{B}_\lambda f)(z)=\int
\limits_{\mathbb{D}_q}b_\lambda(z,\zeta)\,f\,(\zeta)d \nu(\zeta),$ где
$b_\lambda\in \mathscr{D}(\mathbb{D}\times \mathbb{D})_q'$ и
\begin{equation}\label{q-kernel}
b_\lambda(z,\zeta)=\frac{1-q^{2(\lambda-1)}}{1-q^2} (1-\zeta
\zeta^*)^\lambda(q^2z^*\zeta;q^2)_{-\lambda} \cdot(z
\zeta^*;q^2)_{-\lambda}(1-z^*z)^\lambda.
\end{equation}
\end{proposition}

\smallskip

{\bf Доказательство.} Рассмотрим линейный интегральный оператор
$$
\widetilde{\mathscr{B}}_\lambda:\mathscr{D}(\mathbb{D})_q\to
\mathscr{D}(\mathbb{D})_q',\qquad\widetilde{\mathscr{B}}_\lambda:\,f\,
\mapsto\int\limits_{\mathbb{D}_q}\,b_\lambda(z,\zeta)f(\zeta)d\nu(\zeta),
$$
с ядром $b_\lambda$. Оно с точностью до постоянного множителя равно
$U_q\mathfrak{sl}_2$-инвариантному ядру $ k_{22}^{-\lambda}\cdot
k_{11}^{-\lambda}$, см. \itemiii\ \ref{inv_G_m}. Значит,
$\widetilde{\mathscr{B}}_\lambda$ -- морфизм $U_q\mathfrak{sl}_2$-модулей.
Заметим, что $\mathscr{B}_\lambda$ обладает этим же свойством, и что $f_0$
порождает $U_q \mathfrak{sl}_2$-модуль $\mathscr{D}(\mathbb{D})_q$.
Следовательно, равенство
$\widetilde{\mathscr{B}}_\lambda\,f=\mathscr{B}_\lambda\,f$ достаточно
доказать в частном случае $f=f_0$. Но $\widetilde{\mathscr{B}}_\lambda
f_0=(1-q^{2(\lambda-1)})(1-zz^*)^\lambda=\mathscr{B}_\lambda\,f_0$, см.
пример \ref{ex_Berezin}. \hfill $\square$

Получим разложение оператор-функции $\mathscr{B}_\lambda$ по степеням
$t=q^{2(\lambda-1)}$. Очевидно, $t \to 0$ при $\lambda \to +\infty$. В
классическом случае $q=1$ асимптотическое разложение преобразования
Березина при $\lambda \to +\infty$ хорошо известно \cite{Moreno, Cahen,UU,
CGR, Englis, ArazyOrsted}.

Пусть $f_n$ -- элемент $\mathscr{D}(\mathbb{D})_q$, определяемый равенством
$$
T_F(f_n)e_j=
\begin{cases}
e_j,\;&\;j=n,
\\ 0,\;&j\ne n.
\end{cases}
$$

Следующая лемма доказывается с помощью спектрального разложения
$1-zz^*\;=\;\sum_{n=0}^\infty q^{2n}f_n $ элемента $y=1-zz^*$ и равенства
(\ref{5.1.3}).

\medskip

\begin{lemma}\label{l5.3.1}
Для всех $m,j \in{\mathbb Z}_+$ в пространстве
$\mathscr{D}(\mathbb{D})_q^\prime$ имеет место разложение
\begin{equation}\label{5.3.1}
  P_{z,m,j}=\frac{(q^{2\lambda};q^2)_m}{(q^2;q^2)_m}\sum_{n=0}^\infty
q^{2(\lambda-1) n}\cdot P_{z,m,j}^{(n)},
\end{equation}
где $P_{z,m,j}^{(n)}=q^{2(j+n)}z^j \cdot f_n \cdot z^{*m} \in
\mathscr{D}(\mathbb{D})_q$.
\end{lemma}

\bigskip Пусть $f,\psi \in \mathscr{D}(\mathbb{D})_q$. Рассмотрим интеграл
 $\int\limits_{\mathbb{D}_q}\psi^* \, (\mathscr{B}_\lambda f)\, d \nu$ как
функцию переменной $t=q^{2(\lambda-1)}$. Из предложения \ref{Berezin_2.1} и
леммы \ref{l5.3.1} следует аналитичность этой функции в интервале $(0,1)$.
Получаем

\medskip

\begin{proposition}\label{p5.3.2}
Существует и единственна такая последовательность морфизмов $U_q
\mathfrak{sl}_2$-модулей $\mathscr{B}^{(n)}:\mathscr{D}(\mathbb{D})_q \to
\mathscr{D}(\mathbb{D})_q^\prime$, $n \in{\mathbb Z}_+$, что
\begin{equation}\label{5.3.2}
  \mathscr{B}_\lambda\, f=\sum_{n=0}^\infty q^{2(\lambda-1)\, n}\;
\mathscr{B}^{(n)}\, f
\end{equation}
 для всех
$f\in \mathscr{D}(\mathbb{D})_q$.
\end{proposition}

\medskip

 Перейдем к доказательству того, что линейные операторы
 $\mathscr{B}^{(n)}$ отображают $\mathscr{D}(\mathbb{D})_q$ в себя и
 принадлежат подалгебре, порожденной $\square_q$.

Пусть
\begin{multline}\label{5.3.3}
p_j(x)=
\\ =\sum_{k=0}^j\frac{(q^{-2j};q^2)_k}{(q^2;q^2)^2_k}q^{2k}\cdot
\prod_{i=0}^{k-1}\left(1\,+\,q^{2i}\left((1-q^2)^2x-1-q^2\right)+q^{4i+2}
\right).
\end{multline}

\medskip

\begin{lemma}\label{l5.3.3}
$p_j(\square_q)f_0=q^{2j}\; f_j$ для всех $j \in{\mathbb Z}_+$.
\end{lemma}

\smallskip

{\bf Доказательство.} Воспользуемся результатами
 \itemiiiа \ref{radial_part}. Для всех $l \in{\mathbb C}$
базисный гипергеометрический ряд
$$
\Phi_l=_3\!\!\Phi_2 \left[\begin{array}{c}(1-zz^*)^{-1},q^{-2l},q^{2(l+1)}\\
q^2,0\end{array};q^2,q^2\right]
$$
сходится в $\mathscr{D}(\mathbb{D})_q^\prime$, и
$$\square_q \Phi_l=\frac{(1-q^{-2l})(1-q^{2(l+1)})}{(1-q^2)^2}\;\Phi_l.$$
Достаточно показать, что для всех $l \in{\mathbb C}$
\begin{equation}\label{5.3.4}
q^{-2j}\cdot\int \limits_{\mathbb{D}_q}\Phi_l^*\cdot p_j(\square_q)\, f_0d
\nu=\int \limits_{\mathbb{D}_q}\Phi_l^*\,f_jd \nu.
\end{equation}
Заменяя $l$ на $\overline{l}$, приходим к выводу, что (\ref{5.3.4})
эквивалентно равенству
$$p_j \left(\frac{(1-q^{-2l})(1-q^{2(l+1)})}{(1-q^2)^2}\right)=
_3\!\!\Phi_2\left[\frac{q^{-2j},q^{-2l},q^{2(l+1)};q^2;q^2}{q^2,0}\right].
$$
Полагая для краткости $u=q^{2l+2}+q^{-2l}$, запишем его в
следующем виде:

$$
p_j\left(\frac{1+q^2-u}{(1-q^2)^2}\right)=_3\!\!\Phi_2
\left[\frac{q^{-2j},q^{-2l},q^{2(l+1)};q^2;q^2}{q^2,0}\right].
$$
Из определения $\,_3\Phi_2$ получаем:
$$
_3\!\!\;\Phi_2\left[\begin{array}{c} q^{-2j},q^{-2l},q^{2(l+1)}\\
q^2,0\end{array};q^2,q^2\right]=
$$
$$
=\sum_{k=0}^j\frac{(q^{-2j};q^2)_k}{(q^2;q^2)^2_k}\cdot
\prod_{i=0}^{k-1}\left(1-q^{2i}\cdot q^{2(l+1)})(1-q^{2i}\cdot
q^{-2l})\right)\cdot q^{2k}=
$$
$$
=\sum_{k=0}^j \frac{(q^{-2j};q^2)_k}{(q^2;q^2)^2_k}\cdot
\prod_{i=0}^{k-1}\left(1-q^{2i}\cdot u+q^{4i+2}\right)\cdot q^{2k}.
$$
Остается воспользоваться определением полиномов $p_j$. \hfill $\square$

\bigskip

Следующее утверждение существенно уточняет утверждение предложения
\ref{p5.3.2}.

\begin{proposition}\label{p5.3.4}
Для всех $f\in \mathscr{D}(\mathbb{D})_q$ и $\lambda > 1$ имеет место
разложение:
\begin{equation}\label{5.3.5}
 \mathscr{B}_\lambda f=(1-q^{2(\lambda-1)})\sum_{j \in{\mathbb
Z}_+}q^{2(\lambda - 1)\, j}\; p_j(\square_q) f
\end{equation}
в пространстве $\mathscr{D}(\mathbb{D})_q^\prime$.
\end{proposition}

\smallskip

 {\bf Доказательство.} В частном случае $f=f_0$ требуемое утверждение
следует из леммы \ref{l5.3.3}, поскольку
$$
 \mathscr{B}_\lambda f_0=(1-q^{2(\lambda-1)})\,
 \sum \limits_{k \in{\mathbb Z}_+} q^{2(\lambda-1)\,k}\,f_k,
$$
см. пример \ref{ex_Berezin}. Остается учесть, что $f_0$ порождает $U_q
\mathfrak{sl}_2$-модуль $\mathscr{D}(\mathbb{D})_q$, а операторы
$\mathscr{B}_\lambda$, $\square_q$ являются морфизмами $U_q
\mathfrak{sl}_2$-модулей. \hfill $\square$

\medskip

\begin{corollary}\label{c5.3.5}
\begin{equation}\label{5.3.6}
  \mathscr{B}^{(n)}=\left \{\begin{array}{lcl}I &,& n=0 \\
                                      p_n(\square_q)-p_{n-1}(\square_q) &,&
                                      n \in{\mathbb N}\end{array}\right.
\end{equation}
\end{corollary}

\bigskip
\begin{remark}\label{koelink}
Ерик Коелинк обратил наше внимание на то, что полиномы $p_j$ лишь аффинными
заменами переменных отличаются от ортогональных полиномов Аль-Салама-Чихары
\eqref{bible_3.8.1}.
\end{remark}

\bigskip
Как следует из предложения \ref{p5.2.5}, существуют такие элементы
$\varphi_{j,k,\lambda} \in {\rm Pol}(\mathbb{C})_q$, что
$$
\mathscr{B}_\lambda f = \sum\limits_{j,k \in \mathbb{Z}_+}
\,\left(\int\limits_{\mathbb{D}_q} \varphi_{j,k,\lambda}\,f\,d\nu_\lambda
\right)\, z^jz^{*\,k},\qquad f\in \mathscr{D}(\mathbb{D})_q.
$$
 Значит,
преобразование Березина продолжается по непрерывности до линейного
отображения из $L^1(d\nu_\lambda)_q$ в $\mathscr{D}(\mathbb{D})_q'$.
Получим разложение обобщенной функции $\mathscr{B}_\lambda\,f$ для
$f\in{\rm Pol}(\mathbb{C})_q$, аналогичное \eqref{5.3.5}.

\begin{corollary}\label{p5.3.4_new}
Для всех $f\in{\rm Pol}(\mathbb{C})_q$ и $\lambda>1$ имеет место
разложение:
\begin{equation}\label{5.3.5_new}
\mathscr{B}_\lambda
f=\left(1-q^{2(\lambda-1)}\right)\sum_{j\in\mathbb{Z}_+}q^{2(\lambda-1)\,
j}\;p_j(\square_q)f
\end{equation}
в пространстве $\mathscr{D}(\mathbb{D})_q'$.
\end{corollary}

{\bf Доказательство.} Напомним, что алгебра $L^\infty(\mathbb{D})_q$
естественно изоморфна алгебре всех ограниченных линейных операторов в
гильбертовом пространстве. Наделим $L^\infty(\mathbb{D})_q$ ультраслабой
топологией \cite[стр. 76]{BrRob} и будем использовать именно эту топологию.
Достаточно показать, что для некоторой последовательности непрерывных
линейных операторов
$$
B^{(n)}:L^\infty(\mathbb{D})_q\to\mathscr{D}(\mathbb{D})_q',\qquad n\in
\mathbb{Z}_+,
$$
имеет место разложение $\mathscr{B}_\lambda\,f=\sum_{n=0}^\infty
q^{2(\lambda-1)\,n}\;B^{(n)}\,f$ при всех $f\in L^\infty(\mathbb{D})_q$.

Рассуждая так же, как при доказательстве леммы \ref{p5.1.5} и используя
ядерность оператора $(T_F(1-zz^*))^{\lambda-1}$ при $\lambda>1$, получаем
интегральное представление матричных элементов оператора Теплица
$\widehat{f}$ с ограниченным символом $f$:
$$
\widehat{f}_{mj}=\left(1-q^{2(\lambda-1)}\right)\,
\frac{(q^{2\lambda};q^2)_m}{(q^2;q^2)_m}\;{\rm
tr}\left(T_F(z^j(1-zz^*)^{\lambda-1}z^{*m})T_F(f)\right).
$$
Отметим, что оператор-функция
$$
\left(1-q^{2(\lambda-1)}\right)\,\frac{(q^{2\lambda};q^2)_m}{(q^2;q^2)_m}\;
T_F\left(z^j(1-zz^*)^{\lambda-1}z^{*m}\right)
$$
разлагается в сходящийся при $\lambda>1$ ряд по степеням
$q^{2(\lambda-1)}$, и коэффициенты этого степенного ряда являются ядерными
операторами. Остается воспользоваться равенством
$$
\int\limits_{\mathbb{D}_q}\,\psi^*\,(\mathscr{B}_\lambda\,f)\,d\nu\,
=\,\frac{1-q^2}{1-q^{2(\lambda-1)}}\cdot{\rm
tr}_q(\widehat{f}\,\widehat{\psi}),\qquad f\in
L^\infty(\mathbb{D})_q,\;\psi\in\mathscr{D}(\mathbb{D})_q
$$
(ср. с леммой \ref{to_conjugate_new}) и учесть, что оператор
$\widehat{\psi}$ имеет лишь конечное число ненулевых матричных элементов в
базисе $\{z^m\}_{m\in\mathbb{Z}_+}$. \hfill $\square$

\medskip Получим $q$-аналог известного разложения преобразования Березина в
бесконечное произведение \cite{Ber_CMP, UU}. Напомним обозначения
$y=1-zz^*$, $[\lambda]_q=\frac{q^\lambda -q^{-\lambda}}{q-q^{-1}}$.

Используя \eqref{explicit_square}, нетрудно показать, что
$$
\square_q\,y^\lambda\,=\,-q^{-1}\,[\lambda]_q[\lambda-1]_q\,y^\lambda
+[\lambda]_q^2\,y^{\lambda + 1}.
$$
Применяя $\square_q$ к обеим частям равенства
$\mathscr{B}_\lambda\,f_0=(1-q^{2(\lambda-1)})(1-zz^*)^\lambda$, получаем:
$\left(\square_q\mathscr{B}_\lambda+q^{-1}[\lambda]_q[\lambda-1]_q
(\mathscr{B}_\lambda-\mathscr{B}_{\lambda+1})\right)\,f_0=0$. Значит,
$$
\left(1+\frac{q}{[\lambda]_q[\lambda-1]_q}\,\square_q\right)\,
\mathscr{B}_\lambda\,f_0=\,\mathscr{B}_{\lambda+1}f_0,
$$
\begin{equation}\label{prod_f_0}
\mathscr{B}_\lambda\,f_0\,=\,\prod\limits_{j=0}^\infty\,
\left(1+\frac{q}{[\lambda+j]_q[\lambda+j-1]_q}\,\square_q\right)^{-1}\,f_0
\end{equation}
(операторы $1 +\frac{q}{[\lambda+j]_q[\lambda+j-1]_q}\,\square_q$ обратимы
в $L^2(d\nu)_q$ при всех $\lambda>1$, поскольку
$\square_q=\bar{\partial}^*\bar{\partial}\ge 0$).

Так как $\square_q$ и $\mathscr{B}_q$ являются морфизмами
$U_q\mathfrak{sl}_2$-модулей и $f_0$ порождает $U_q\mathfrak{sl}_2$-модуль
$\mathscr{D}(\mathbb{D})_q$, то из \eqref{prod_f_0} вытекает равенство
\begin{equation}\label{prod_Berezin}
\mathscr{B}_\lambda\,=\,\prod\limits_{j=0}^\infty\,
\left(1+\frac{q}{[\lambda+j]_q[\lambda+j-1]_q}\,\square_q\right)^{-1}
\end{equation}
($\mathscr{B}_\lambda:\mathscr{D}(\mathbb{D})_q\to\mathscr{D}(\mathbb{D})_q'$,
и произведение операторов в правой части \eqref{prod_Berezin} сходится на
каждом векторе из $\mathscr{D}(\mathbb{D})_q$).

\subsubsection{Квантование по Березину}\label{Berezin-quantization_sl_2}

Рассмотрим унитальную алгебру $A$ и $\mathbb{C}[[t]]$-модуль $A[[t]]$
формальных рядов с коэффициентами из $A$. Деформацией алгебры $A$
\IND{деформация алгебры} называют $\mathbb{C}[[t]]$-билинейное умножение
$$*_t:\,A[[t]]\times A[[t]]\,\to\,A[[t]]$$
вида $f_1*_tf_2=f_1f_2\,+\,\sum\limits_{k=1}^\infty t^k\;m_k(f_1,f_2)$, где
$m_k$ -- билинейные отображения из $A\times A$ в $A$. В определение
включается требование ассоциативности умножения $*_t$ и наличие единицы,
точнее, требование, чтобы естественное вложение $A\hookrightarrow A[[t]]$
было гомоморфизмом комплексных унитальных алгебр: $m_k(1,f)=m_k(f,1)=0$.
Это позволяет отождествить алгебру $A$ с ее образом при вложении в
$A[[t]]$.

Если алгебра $A$ является $U_q\mathfrak{su}_{1,1}$-модульной $*$-алгеброй,
то $A[[t]]$ наследует эту структуру (подразумевается, что $t^*=t$).
Деформацию называют $U_q\mathfrak{su}_{1,1}$-модульной, \IND{деформация
алгебры ! $U_q\mathfrak{su}_{1,1}$-модульная} если
$$(f_1*_tf_2)^*=f_2^*\,*_t\,f_1^*,\qquad f_1,f_2\in A[[t]],$$
и все отображения
$$A\otimes A\to A,\qquad f_1\otimes f_2\mapsto m_k(f_1,f_2),$$
являются морфизмами $U_q\mathfrak{sl}_2$-модулей.

Обратимся к частному случаю $A={\rm Pol}(\mathbb{C})_q$. Говорят, что
деформация алгебры ${\rm Pol}(\mathbb{C})_q$ является деформацией с
разделяющимися переменными, \IND{деформация алгебры ! с разделяющимися
переменными} если $f_1*_tf\,=\,f_1f$ и $f*_tf_2=ff_2$ для всех
$f_1\in\mathbb{C}[z]_q$, $f_2\in\mathbb{C}[\bar{z}]_q$, $f\in{\rm
Pol}(\mathbb{C})_q$, см. \cite{Karabegov96}.

\bigskip
 Коммутационное соотношение \eqref{1.7.2} и подразумеваемая всюду  далее в этом
\itemiiiе подстановка $t=q^{2(\lambda-1)}$
 позволяют получить вещественную
$U_q\mathfrak{sl}_2$-инвариантную деформацию с разделяющимися переменными
алгебры ${\rm Pol}(\mathbb{C})_q$, как показывает следующее утверждение.

\begin{proposition}\label{def_def}
1. Существует и единственна такая деформация c разделяющимися переменными
алгебры ${\rm Pol}(\mathbb{C})_q$, что
\begin{equation}\label{KliLesn}
z^* \, *_t z \;=\; q^2 zz^*\, +\, 1-q^2\; +\;
 \frac{(1-q^2)t}{1-q^2t}\,(1-z^*z)(1-zz^*).
\end{equation}
2. Эта деформация является $U_q\mathfrak{su}_{1,1}$-инвариантной.
\end{proposition}

Нетрудно доказать это предложение, исходя из определений. В дальнейшем мы
получим более точный результат методом Березина.

Как показано заметке Шклярова в сборнике \cite{Vaks_lectures}, всякая
$U_q\mathfrak{su}_{1,1}$-инвариантная деформация с разделяющимися
переменными алгебры ${\rm Pol}(\mathbb{C})_q$ может быть получена из
описанной выше деформации формальной заменой $ t:= \sum\limits_{j=1}^\infty
a_jt^j$, где $a_j \in \mathbb{R}$. Этот результат согласуется с описанием
инвариантных деформаций с разделяющимися переменными коммутативных алгебр
\cite{MB-Neumaier}.

\bigskip Как будет показано,  для деформации, описанной в предложении \ref{def_def},
 билинейные операторы $m_k$ являются $q$-би\-диффе\-рен\-циаль\-ными.
Наряду с введенными в \itemiiiе \ref{finite-forms_sl_2} ''левыми'' частными
производными
 $$
 \frac{\partial^{(l)}\, f}{\partial z}\,=\,
 \frac{\partial\, f}{\partial z}, \qquad
 \frac{\partial^{(l)}\, f}{\partial z^*}\,=\,
 \frac{\partial\, f}{\partial z^*}
 $$
 нам понадобятся ''правые'' частные производные,
 определяемые равенствами
\begin{equation*}\label{partial_right_sl_2}
 \partial\,f\,=\, \frac{\partial^{(r)}\, f}{\partial z}\, dz, \qquad
\bar{\partial}\,f\,=\,\frac{\partial^{(r)}\, f}{\partial z^*}\, dz^*.
 \end{equation*}

\bigskip

Рассмотрим алгебру линейных операторов в векторном пространстве
 ${\rm Pol}(\mathbb{C})_q=\mathbb{C}[z]_q  \otimes \mathbb{C}[\bar{z}]_q$
и ее гомоморфизм $A\mapsto \widetilde{A}$ в алгебру линейных операторов в
${\rm Pol}(\mathbb{C})_q\otimes{\rm Pol}(\mathbb{C})_q$, определяемый с
помощью изоморфизма векторных пространств
$$
i:\,\mathbb{C}[\bar{z}]_q\otimes\mathbb{C}[z]_q\stackrel{\approx}{\to}{\rm
Pol}(\mathbb{C})_q,\qquad i:\,g_1\otimes g_2\mapsto g_1g_2$$ равенством
$\widetilde{A}={\rm id}_{\mathbb{C}[z]_q}\otimes(i^{-1}\,A\;i)\otimes{\rm
id}_{\mathbb{C}[\bar{z}]_q}$.

\begin{example}
Нетрудно показать, что
$$
\square_q(f(z^*)\,g(z))\,=\,-q^2\,\frac{\partial^{(r)}\,f}{\partial z^*}\;
\left(1-\left(1+q^{-2}\right)z^*z+q^{-2}z^{*2}z^2\right)\;
\frac{\partial^{(l)}\,g}{\partial z}
$$
для всех $f\in\mathbb{C}[\bar{z}]_q$, $g\in\mathbb{C}[z]_q$. Значит,
$$
\widetilde{\square}_q\,=\,-q^{-2}\left(1-\left(1+q^{-2}\right)z^*\otimes
z+q^{-2}z^{*2}\otimes z^2\right)\;\frac{\partial^{(r)}}{\partial
z^*}\otimes\frac{\partial^{(l)}}{\partial z},
$$
где умножение производится в алгебре ядер ${\rm Pol}(\mathbb{C})_q^{\rm
op}\otimes{\rm Pol}(\mathbb{C})_q$.
\end{example}

\begin{proposition}\label{t1.7.3}
\hspace{-.5em}. Для всех $f_1,f_2 \in{\rm Pol}({\mathbb C})_q$
$$
f_1*_t f_2=(1-t)\sum_{j=0}^\infty\; t^j \cdot
m(p_j(\stackrel{\sim}{\square}_q)(f_1 \otimes f_2)),
$$
где $p_j$ -- полиномы (\ref{5.3.3}).
\end{proposition}

\medskip {\bf Доказательство предложений
 \ref{def_def}, \ref{t1.7.3}.} Будем писать $\mathbb{C}[z]_q$ вместо $\mathbb{C}[z]_{q,\lambda}$
 в тех случаях, когда структура $U_q\mathfrak{sl}_2$-модуля не играет
 существенной роли.

Пусть $Q:{\rm Pol}({\mathbb C})_q[[t]]\to{\rm End}({\mathbb C}[z]_q)[[t]]$
-- такой непрерывный в $t$-адической топологии морфизм
$\mathbb{C}[[t]]$-модулей, что
$$
Q:z^jz^{*k}\mapsto\widehat{z}^j\widehat{z}^{*k},\qquad j,k\in\mathbb{Z}_+.
$$
Это отображение сопоставляет полиному $f\in{\rm
Pol}(\mathbb{C})_q$ линейный оператор, ковариантный символ
которого равен $f$, см. предложение \ref{p5.6.6}. Умножение в
${\rm End}({\mathbb C}[z]_q)[[t]]$ ассоциативно, и, следуя
Березину \cite{Ber_CMP}, мы перенесем его в ${\rm Pol}({\mathbb
C})_q[[t]]$ с помощью отображения $Q$.

\begin{lemma}\label{lemma_new_sl_2}
Существует единственное такое ${\mathbb C}[[t]]$-билинейное отображение
$*_t:{\rm Pol}({\mathbb C})_q[[t]]\times{\rm Pol}({\mathbb C})_q[[t]]
 \to{\rm Pol}({\mathbb C})_q[[t]]$, что
\begin{equation}\label{for_*}
 Q(f_1*_t f_2)=(Qf_1)\cdot(Qf_2)
 \end{equation}
 для всех
$f_1,f_2 \in{\rm Pol}({\mathbb C})_q[[t]]$.
\end{lemma}

{\bf Доказательство леммы.} Единственность следует из инъективности
отображения $Q$, см. следствие \ref{c5.6.7}. Предъявим умножение $*_t$,
удовлетворяющее требованию \eqref{for_*}, на образе вложения
$$
\mathbb{C}[\bar{z}]_q\times\mathbb{C}[z]_q\hookrightarrow {\rm
Pol}(\mathbb{C})_q[[t]]\times{\rm Pol}(\mathbb{C})_q[[t]],\qquad f\times g
\mapsto (1\otimes f)\times (g\otimes 1).
$$
 Сопоставим каждой паре
$(f(z^*),g(z))\in\mathbb{C}[\bar{z}]_q\times\mathbb{C}[z]_q$ операторы
$\check{f}$, $\check{g}$ с ковариантными символами $f(z^*)$, $g(z)$. Они
сами и их произведение являются аналитическими функциями переменной
$t=q^{2(\lambda-1)}$ со значениями в банаховой алгебре ограниченных
линейных операторов. Ковариантный символ оператора
$\check{f}\cdot\check{g}$ равен
$$
\mathscr{B}_\lambda(f(z^*)g(z))=
(1-q^{2(\lambda-1)})\sum_{j\in\mathbb{Z}_+}q^{2(\lambda-1)\,j}\;
p_j(\square_q)(f(z^*)g(z)),
$$
см. следствие \ref{5.3.5_new} и пример \ref{Toepltz_symbol}. Остается
заменить $q^{2(\lambda-1)}$ на $t$:
$$
f(z^*)*_t g(z)\stackrel{\rm def}{=}(1-t)\sum_{j\in\mathbb{Z}_+}t^j\;
p_j(\square_q)(f(z^*)g(z)),\quad
f(z^*)\in\mathbb{C}[\bar{z}]_q,\;g(z)\in\mathbb{C}[z]_q.
$$

 Используя предложение \ref{p5.6.6}, продолжим $*_t$ на
 $$ {\rm Pol}(\mathbb{C})_q\times  {\rm Pol}(\mathbb{C})_q \subset
 {\rm Pol}({\mathbb C})_q[[t]]\times{\rm Pol}({\mathbb C})_q[[t]],
 $$
 с сохранением свойства \eqref{for_*}:
 $$
 f_1(z)f_2(z^*) *_t g_1(z)g_2(z^*)=(1-t)\sum_{j \in{\mathbb Z}_+}
 t^j f_1(z)\cdot p_j(\square_q)(f_2(z^*)g_1(z))\cdot g_2(z^*),
 $$
 где $f_1(z), g_1(z) \in \mathbb{C}[z]_q$,\,
 $f_2(z^*), g_2(z^*) \in \mathbb{C}[\bar{z}]_q$. Остается продолжить
 $*_t$ на
 ${\rm Pol}({\mathbb C})_q[[t]]\times{\rm Pol}({\mathbb C})_q[[t]]$
 по $\mathbb{C}[[t]]$-линейности и непрерывности в $t$-адической
 топологии. Лемма доказана. \hfill $\square$

\medskip

Рассмотрим отображение $*_t$, описанное в лемме. Из его построения и
равенства \eqref{1.7.2} вытекает \eqref{KliLesn}. Ассоциативность $*_t$ и
$U_q\mathfrak{su}_{1,1}$-инвариантность этой деформации, а также то, что
она является деформацией с разделяющимися переменными, проверяется, исходя
из определений. Отсюда следуют утверждения предложений \ref{def_def},
\ref{t1.7.3}. \hfill $\square$

\bigskip
\begin{remark}
При построении деформации методом Березина функции $f$ сопоставляется
оператор Теплица с ковариантным символом $f$. При другом подходе, так
называемом квантовании Теплица-Березина, функции $f$ сопоставляется
оператор Теплица с контравариантным символом $f$, см., например,
\cite{KlimLesn_disc,BLU,Englis05}.
\end{remark}

\subsubsection{Дополнение о классах Неймана-Шэттена.}\label{S_p}

Пусть $A$ -- ограниченный неотрицательный линейный оператор в сепарабельном
гильбертовом пространстве $H$. Сумма ${\rm tr}(A)=\sum\limits_j(Ae_j,e_j)$
не зависит от выбора ортонормированного базиса $\{e_j\}$ в $H$, причем для
конечности следа необходимо, чтобы оператор $A$ был компактен \cite[стр.
126]{GKr}.

Пусть $p\ge 1$. Класс Неймана-Шэттена \IND{класс Неймана-Шэттена} $S_p$
образуют компактные линейные операторы в $H$, для которых
$\|A\|_p\,=\,({\rm tr}(|A|^p))^{\frac1{p}}<\infty$. Здесь
$|A|=\sqrt{A^*A}$. След -- это некоммутативный аналог интеграла, а классы
$S_p$ -- некоммутативные аналоги пространств $L^p$. Некоммутативным
аналогом алгебры $L^\infty$ является $\mathcal{B}(H)$.

Классы $S_p$ можно ввести по-другому. Cингулярные числа оператора $A$
\IND{сингулярные числа оператора} определим равенством
$s_j(A)\,=\,\inf\{\|A-B\|\;|\;{\rm rk}(B)\le j-1\}$, \cite[стр. 48]{GKr}.
Здесь $j\in\mathbb{N}$. Можно показать, что
$\|A\|_p\,=\,\left(\sum\limits_{j=1}^\infty s_j^p(A)\right)^{\frac1{p}}$.
Получаем определение классов $S_p$ в терминах сингулярных чисел \cite[стр.
121]{GKr}. Эти числа $s_j(A)$ можно определить иначе. Они являются
собственными значениями оператора $|A|$, взятыми с учетом кратности и
упорядоченными по убыванию \cite[стр. 46]{GKr}.

Нетрудно показать \cite[стр. 19]{Zhu}, что
$\|A\|_p^p=\sup\{\sum\limits_j|(Ae_j,e_j)|^p\;| \; \{e_j\}\}$, где $\{e_j\}$
пробегает множество всех ортонормированных базисов гильбертова пространства
$H$. Значит, функционал $\|\cdot\|_p$ является нормой. Другое доказательство
см. в \cite[стр. 107]{GKr}.
 Нормированные пространства $S_p$ являются полными \cite[стр.
108]{GKr}.

Как следует из определений, $S_p$ -- двусторонний идеал алгебры
ограниченных линейных операторов в $H$ и $\|A\|_p=\|U_1AU_2\|_p$ для любых
унитарных операторов $U_1$, $U_2$.

Операторы из $S_1$ называются ядерными, и след продолжается по
непрерывности с двустороннего идеала конечномерных линейных операторов на
$S_1$. Операторы из $S_2$ называются операторами Гильберта-Шмидта,
\IND{оператор ! Гильберта-Шмидта} и
$\|A\|_2=\sum\limits_{j,k}|(Ae_j,e_k)|^2$ для любого ортонормированного
базиса $\{e_n\}$ в $H$.

Пусть $p>1$, $q>1$ и $\frac1{p}+\frac1{q}=1$. Каждому ограниченному
линейному оператору $A\in S_p$ отвечает непрерывный линейный функционал
$l_A(X)={\rm tr}(AX)$ на $S_q$, и его норма равна $\|A\|_p$. Это общий вид
непрерывного линейного функционала на $S_q$ \cite[стр. 170]{GKr}.

Аналогично, каждому ограниченному линейному оператору $A$ в $H$ отвечает
непрерывный линейный функционал $l_A(X)={\rm tr}(AX)$ на $S_1$ и его норма
равна $\|A\|$. Это общий вид непрерывного линейного функционала на $S_1$
\cite[стр. 166]{GKr}. Важно отметить, что $S_1$ -- пространство всех
ультраслабо непрерывных линейных функционалов на $\mathcal{B}(H)$
\cite[стр. 77]{BrRob}.

\bigskip

Хорошо известны интерполяционные теоремы в шкале банаховых пространств
$S_p$ \cite[стр. 140-145]{GKr2}. В частности, $S_p$ можно получить из
$\mathcal{B}(H)$ и $S_1$ методом комплексной интерполяции:
\begin{equation}\label{compl_interpol}
S_p\,=\,[\mathcal{B}(H),S_1,]_{p^{-1}}.
\end{equation}

Поясним сказанное. Рассмотрим векторное пространство ограниченных и
непрерывных в вертикальной полосе $0\le{\rm Re}\,z\le 1$ оператор-функций
со значениями в $\mathcal{B}(H)$. Пусть $\mathcal{A}$ -- подпространство
голоморфных при $0<{\rm Re}\,z<1$ оператор-функций, сужения которых на
прямую ${\rm Im}\,z=1$ непрерывны и ограничены в $S_1$. Нетрудно показать
\cite[стр. 288]{ KreinPetuninSemenov}, что $\mathcal{A}$ с нормой
$$
\|f\|_\mathcal{A}=\max\{\sup\limits_{y\in\mathbb{R}}\|f(iy)\|\,,\,
\sup\limits_{y\in\mathbb{R}}\|f(1+iy)\|_1\}
$$
является банаховым пространством.

Оказывается \cite[стр. 31]{Zhu}, операторы из $S_p$ и только они являются
значениями оператор-функций из $\mathcal{A}$ в точке $p^{-1}$, причем
$$
\|A\|_p=\inf\{\|f\|_\mathcal{A}\;|\; f \in \mathcal{A}\;\&\;f(p^{-1})=A\}
$$
для всех $A \in S_p$.

   Отсюда вытекает следующий результат
   \cite[стр. 309]{KreinPetuninSemenov}. Пусть линейный оператор
  $T$ в банаховом пространстве $\mathcal{B}(H)$ отображает $S_1$ в себя
  и \hbox{$\|T\|_{\mathcal{B}(H)} <\infty$}, $\|T\|_{S_1} <\infty$. Тогда
 $T$ отображает каждое из пространств $S_p$ в себя и
норма $\|T\|_{S_p}$ логарифмически выпукло зависит от $\theta=p^{-1}$.

\subsection{О ядрах сплетающих интегральных операторов}\label{intertw_sl_2}

\subsubsection{ Вложение $\operatorname{Pol}(\mathbb{C})_q\hookrightarrow
\mathbb{C}[w_0SU_{1,1}]_{q,x}$.}\label{canon_embed_sl_2}
\label{canon_example_sl_2}

Единичный круг $\mathbb{D}$ является $SU_{1,1}$-пространством, и в категории
$SU_{1,1}$-пространств
\begin{equation}\label{naive_embedding}
\mathbb{D}\cong U_1\setminus SU_{1,1},
\end{equation}
где $U_1=\left\{\mathrm{diag}(e^{i\varphi},e^{-i\varphi})|\:
\varphi\in\mathbb{R}/(2\pi\mathbb{Z})\right\}$. Изоморфизм
\eqref{naive_embedding} доставляет вложение $*$-алгебры полиномов
$\mathrm{Pol}(\mathbb{C})$ в $*$-алгебру рациональных функций на
вещественной аффинной алгебраической группе $SU_{1,1}$, при котором
\begin{equation}\label{embed_SU11}
z\mapsto t_{11}^{-1}t_{12}.
\end{equation}
Обладает ли это вложение $q$-аналогом? Для того, чтобы ответ оказался
утвердительным, постановку вопроса следует немного изменить. Пусть
$w_0=\begin{pmatrix}0, & -1\\ 1, & 0\end{pmatrix}$. Группу $SU_{1,1}$
заменим ее главным однородным пространством
\begin{equation}\label{su_1_1-orbit}
w_0\,SU_{1,1}=\left\{\left.\left(\begin{array}{cl}t_{11}, & \quad t_{12}\\
t_{21}, & \quad t_{22}\end{array}\right)\in
SL_2\quad\right|\quad\overline{t}_{11}=-t_{22},\:\overline{t}_{12}=
-t_{21}\right\}.
\end{equation}
а вложение \eqref{embed_SU11} -- вложением
\begin{equation}\label{embed_hat}
z\mapsto t_{12}^{-1}t_{11}.
\end{equation}

\bigskip
Введем квантовый аналог алгебры регулярных функций на $SU_{1,1}$-орбите
\eqref{su_1_1-orbit}. Наделим ${\mathbb C}[SL_2]_q$ инволюцией
\begin{equation}\label{inv}
t_{11}^*=-t_{22},\qquad t_{12}^*=-qt_{21}
\end{equation}
и введем обозначение $\mathbb{C}[w_0SU_{1,1}]_q=(\mathbb{C}[SL_2]_q,*)$ для
полученной $*$-алгебры. Одним из доводов в пользу сделанного выбора
инволюции является то, что в формальном пределе $q \to 1$ получается система
уравнений $\bar{t}_{11}=-t_{22}$, $\bar{t}_{12}=-t_{21}$, выделяющая
вещественное
 $w_0SU_{1,1}$ из $SL_2$. Более
содержательным доводом служит приведенное ниже предложение \ref{*-hat}.

 Напомним, что алгебра $\mathbb{C}[SL_2]_q$ регулярных функций на
квантовой группе $SL_2$ порождена матричными элементами $t_{ij}$ векторного
представления алгебры $U_q\mathfrak{sl}_2$, см. \eqref{q2-dim_rep}, и что
  в \itemiiiе \ref{SL_2_sl_2} алгебра
${\mathbb C}[SL_2]_q$ была наделена структурой
$U_q\mathfrak{sl}_2$-модульной алгебры.

\bigskip Следующее утверждение вытекает из определений.
\begin{proposition}\label{*-hat}
 $*$-Алгебра
$\mathbb{C}[w_0SU_{1,1}]_q$ является $U_q\mathfrak{su}_{1,1}$-модульной, то
есть
\begin{equation}\label{*_*}
(\xi f)^*=(S(\xi))^*f^*,\qquad \xi\in U_q\mathfrak{su}_{1,1},\;
f\in\mathbb{C}[w_0SU_{1,1}]_q.
\end{equation}
\end{proposition}

\medskip Как и в классическом случае $q=1$, см. \eqref{embed_hat}, алгебра регулярных
 функций слишком мала для наших целей, и ее нужно расширить. Воспользуемся тем, что
элемент
\begin{equation}\label{def_x_SL2}
x=-qt_{12}t_{21}
\end{equation}
самосопряжен и квазикоммутирует со всеми образующими $t_{ij}$:
$$
t_{11}\,x=q^2\,x\,t_{11},\qquad t_{22}\,x=q^{-2}\,x\,t_{22},\qquad
t_{12}\,x=x\,t_{12}, \qquad t_{21}\,x=x\,t_{21}.
$$
Значит, мультипликативное множество $x^{\mathbb{Z}_+}$ удовлетворяет условию
Оре. Как отмечалось в \itemiiiе \ref{SL_2_sl_2}, алгебра
$\mathbb{C}[SL_2]_q$ не имеет делителей нуля. Следовательно, она вложена в
свою локализацию $\mathbb{C}[w_0SU_{1,1}]_{q,x}$ по мультипликативному
множеству $x^{\mathbb{Z}_+}$. Самосопряженность $x$ позволяет продолжить
инволюцию $*$ на $\mathbb{C}[w_0SU_{1,1}]_{q,x}$ и доказать, что
естественное вложение
$\mathbb{C}[w_0SU_{1,1}]_{q,x}\hookrightarrow\mathbb{C}[w_0SU_{1,1}]_{q,x}$
является гомоморфизмом $*$-алгебр.

Очевидно, элементы $t_{12}$, $t_{21}$ обратимы в
$\mathbb{C}[w_0SU_{1,1}]_{q,x}$:
$$t_{12}^{-1}=t_{21}\,x,\qquad t_{21}^{-1}=t_{12}\,x.$$

\begin{proposition}\label{*_local}
\begin{itemize}
\item[1.] Существует и единственно продолжение \\ структуры
    $U_q\mathfrak{su}_{1,1}$-модульной алгебры с
    $\mathbb{C}[w_0SU_{1,1}]_q$ на $\mathbb{C}[w_0SU_{1,1}]_{q,x}$.

\item[2.] Для любых $\xi\in U_q\mathfrak{su}_{1,1}$,
$f\in\mathbb{C}[w_0SU_{1,1}]_{q,x}$ существует и единствен такой полином
Лорана $p_{f,\xi}$ одной переменной с коэффициентами из
$\mathbb{C}[w_0SU_{1,1}]_{q,x}$, что
\begin{equation}\label{p_f}
\xi(f\cdot x^n)=p_{f,\xi}(q^n)\cdot x^n
\end{equation}
при всех $n\in\mathbb{Z}$.
\end{itemize}
\end{proposition}

{\bf Доказательство.} Утверждения о единственности очевидны. Исходя из
\eqref{r_reg1} -- \eqref{r_reg3}, легко показать, что для любого полинома
$\psi$ одной переменной
\begin{equation}\label{K_psi}
K^{\pm 1}\psi(x)=\psi(x),
\end{equation}
\begin{equation}\label{E_psi}
E\psi(x)=-q^{1/2}t_{11}\frac{\psi(q^{-2}x)-\psi(x)}{q^{-2}x-x}t_{21}=
-q^{1/2}t_{11}t_{21}\frac{\psi(q^{-2}x)-\psi(x)}{q^{-2}x-x},
\end{equation}
\begin{equation}\label{F_psi}
F\psi(x)=-q^{3/2}t_{12}\frac{\psi(q^{-2}x)-\psi(x)}{q^{-2}x-x}t_{22}=
-q^{3/2}t_{12}t_{22}\frac{\psi(x)-\psi(q^2x)}{x-q^2x}.
\end{equation}
Как следует из этих равенств, для любых $\xi\in U_q\mathfrak{su}_{1,1}$,
$f\in\mathbb{C}[w_0SU_{1,1}]_q$ существует такой полином Лорана $p_{f,\xi}$
с коэффициентами из $\mathbb{C}[w_0SU_{1,1}]_q$, что равенство \eqref{p_f}
имеет место при всех $n\in\mathbb{Z}_+$. Полином Лорана $p_{f,\xi}$
однозначно определяется своими значениями в точках геометрической
прогрессии $q^\mathbb{N}$. Из его единственности следует что равенство
\eqref{p_f} корректно определяет действие элементов $\xi\in
U_q\mathfrak{su}_{1,1}$ в пространстве $\mathbb{C}[w_0SU_{1,1}]_{q,x}$,
поскольку
$$p_{f,\xi}(q^n)x=p_{f\cdot x,\xi}(q^{n-1}),\qquad n\in\mathbb{Z}.$$
Нетрудно доказать, что полученное отображение
$$
\pi:U_q\mathfrak{su}_{1,1}\to
\operatorname{End}(\mathbb{C}[w_0SU_{1,1}]_{q,x})
$$
линейно и является гомоморфизмом алгебр, а также доказать согласованность
инволюций \eqref{*_*}. Действительно, требуемые равенства нетрудно свести к
равенствам полиномов Лорана, заменив входящие в них элементы
$f\in\mathbb{C}[w_0SU_{1,1}]_{q,x}$ последовательностями $f_n=f\cdot x^n$,
$n\in\mathbb{Z}_+$. Остается заметить, что каждое из требуемых равенств
полиномов Лорана справедливо в точках геометрической прогрессии
$\{q^n\}_{n\in\mathbb{Z}_+}$ при достаточно больших $n\in\mathbb{N}$.

   Наконец, получим равенство
\begin{equation}\label{Umod_extension}
 \xi(f_1f_2)=\sum_i \xi'_i(f_1)\, \xi''_i(f_2),\qquad
 f_1, f_2 \in \mathbb{C}[w_0SU_{1,1}]_{q,x},\; \xi \in
  U_q\mathfrak{sl}_2,
\end{equation}
  где подразумевается, что $\Delta \xi=\sum\limits_i \xi'_i \otimes \xi''_i$.
Так же, как прежде, доказывается существование и единственность полинома
Лорана $p_{f,\xi}$ двух переменных, для которого при всех $m, n$
$$ \xi(x^m f x^n)\,=\,x^m\, p_{f,\xi}(q^m,q^n)\, x^n.$$
 Если числа $m,n \in \mathbb{Z}_+$ достаточно велики, то
\begin{equation}\label{two_sided_x}
 x^{-m}\,\xi(x^mf_1f_2x^n)\, x^{-n}= \sum_i (x^{-m}\xi'_i(x^mf_1))\,
(\xi''_i(f_2x^n)\,x^{-n}),
\end{equation}
поскольку алгебра $\mathbb{C}[SL_2]_q$ является
$U_q\mathfrak{sl}_2$-модульной. Приходим к равенству двух полиномов Лорана
во всех точках $\{(q^m,q^n)\}$ с достаточно большими $m,n$. Значит, эти
полиномы Лорана равны, и \eqref{two_sided_x} имеет место при всех $m,n \in
\mathbb{Z}$. Полагая $m=n=0$, получаем \eqref{Umod_extension}. \hfill
$\square$

\begin{proposition}\label{*_embed}
Отображение $\mathcal{I}:z\mapsto t_{12}^{-1}t_{11}$ единственным образом
продолжается до вложения $U_q\mathfrak{su}_{1,1}$-модульных алгебр
\begin{equation}\label{disc_embed}
\mathcal{I}:\operatorname{Pol}(\mathbb{C})_q\hookrightarrow
\mathbb{C}[w_0SU_{1,1}]_{q,x}.
\end{equation}
\end{proposition}

{\bf Доказательство.} Единственность рассматриваемого гомоморфизма
$*$-алгебр очевидна, а его существование вытекает из следующих равенств в
алгебре $\mathbb{C}[w_0SU_{1,1}]_{q,x}$:
$$(t_{12}^{-1}t_{11})^*=q^{-1}t_{22}t_{21}^{-1},\quad
(t_{12}^{-1}t_{11})(t_{22}t_{21}^{-1})-
q^{-2}(t_{22}t_{21}^{-1})(t_{12}^{-1}t_{11})=q-q^{-1}.$$

Докажем инъективность гомоморфизма $\mathcal{I}$. Рассмотрим
предгильбертово пространство с ортонормированным базисом
$\{e_j\}_{j\in\mathbb{Z}_+}$ и представление $\pi_+$ алгебры
$\mathbb{C}[SL_2]_q$, введенное в \itemiiiе \ref{SL_2_sl_2}. Это
$*$-представление алгебры $\mathbb{C}[w_0SU_{1,1}]_q$. Сохраним обозначение
$\pi_+$ для его продолжения до $*$-представления локализации
$\mathbb{C}[w_0SU_{1,1}]_{q,x}$.

Из определений следует, что представление $\pi_+\cdot\mathcal{I}$ алгебры
$\operatorname{Pol}(\mathbb{C})_q$ эквивалентно ее фоковскому
представлению. Значит, представление $\pi_+\circ\mathcal{I}$ является
точным, см. \itemiii\ \ref{Fock_l}. Следовательно, гомоморфизм
$\mathcal{I}$ инъективен.

Докажем, что $\mathcal{I}$ является морфизмом $U_q\mathfrak{sl}_2$-модулей.
Так как \hbox{$*$-алгебры} $\operatorname{Pol}(\mathbb{C})_q$ и
$\mathbb{C}[w_0SU_{1,1}]_{q,x}$ являются
$U_q\mathfrak{su}_{1,1}$-модульными, то равенство $\mathcal{I}(\xi
f)=\xi\mathcal{I}(f)$ достаточно получить для $\xi\in
U_q\mathfrak{su}_{1,1}$, $f\in\mathbb{C}[z]_q$. Остается воспользоваться
вложением $U_q\mathfrak{sl}_2$-модульных алгебр
$$
\mathbb{C}[t_1,t_2]_q\to\mathbb{C}[SL_2]_q,\qquad t_1\mapsto t_{11},\;t_2
\mapsto t_{12}
$$
и описанным в \itemiiiе \ref{sL_2_mod_algebras} вложением $\mathbb{C}[z]_q$
в локализацию алгебры $\mathbb{C}[t_1,t_2]_q$. \hfill $\square$

\bigskip Описанное вложение $U_q\mathfrak{su}_{1,1}$-модульных алгебр \\
$\operatorname{Pol}(\mathbb{C})_q\hookrightarrow
\mathbb{C}[w_0SU_{1,1}]_{q,x}$ будем называть каноническим вложением.

\subsubsection{Инвариантные интегралы.}\label{finite_w_0SU_11}

В силу следствия \ref{PBW_SL_2}, каждый элемент $f \in
\mathbb{C}[w_0SU_{1,1}]_q$ единственным образом разлагается в сумму
\begin{equation}\label{canon_expan}
f=\sum\limits_{ad=bc=0} t_{11}^at_{12}^b\psi_{abcd}(x)t_{21}^ct_{22}^d,
\end{equation}
где $a,b,c,d \in \mathbb{Z}_+$ и $\psi_{abcd}$-полиномы одной переменной.

Построения \itemiiiа \ref{finite_sl_2} и равенство
\begin{equation}\label{x_y}
\mathcal{I}(y)\;=\;x^{-1}
\end{equation}
наводят на мысль ввести в рассмотрение векторное пространство
$\mathscr{D}(w_0SU_{1,1})'_q$ формальных рядов \eqref{canon_expan} с
коэффициентами из пространства функций на множестве $q^{-2\mathbb{Z}_+}$ и
наделить его слабейшей из топологий, в которых непрерывны линейные
функционалы $\mathcal{L}_{a,b,c,d,n}$:
\begin{equation}\label{topology_SU_1_1_new}
\mathcal{L}_{a,b,c,d,n}(f)\;=\;\psi_{abcd}(q^{-2n}),\qquad ad=bc=0,\;n\in
\mathbb{Z}_+.
\end{equation}

Полиномы плотны в пространстве функций на $q^{-2\mathbb{Z}_+}$ с топологией
поточечной сходимости. Значит, $\mathbb{C}[w_0SU_{1,1}]_q$ -- плотное
линейное подмногообразие в $\mathscr{D}(w_0SU_{1,1})'_q$. Очевидно,
инволюция $*$ продолжается по непрерывности на
$\mathscr{D}(w_0SU_{1,1})'_q$.

Из равенств \eqref{K_psi} -- \eqref{F_psi} и
$$
(t_{11}^at_{12}^b\,\psi(x)\,t_{21}^c)\,t_{22}\;=\;
q^{b+c}t_{11}^{a-1}t_{12}^b\,(1-x)\psi(q^2x)\,t_{21}^c,
$$
$$
t_{11}\,(t_{12}^b\psi(x)\,t_{21}^ct_{22}^d)\;=
\;q^{b+c}t_{12}^b\,(1-x)\psi(q^2x)\,t_{21}^ct_{22}^{d-1}
$$
следует, что структура $U_q\mathfrak{sl}_2$-модульного
$\mathbb{C}[w_0SU_{1,1}]_q$-бимодуля также допускает продолжение по
непрерывности с $\mathbb{C}[w_0SU_{1,1}]_q$ на
$\mathscr{D}(w_0SU_{1,1})'_q$.

Обобщенную функцию $f\in\mathscr{D}(w_0SU_{1,1})'_q$ назовем финитной,
\IND{финитная обобщенная функция из $f\in\mathscr{D}(w_0SU_{1,1})'_q$} если
$${\rm card}\{(a,b,c,d,n)\;|\;\mathcal{L}_{a,b,c,d,n}(f)\ne 0\}<\infty.$$
Векторное пространство $\mathscr{D}(w_0SU_{1,1})_q$ финитных функций
является подмодулем $U_q\mathfrak{sl}_2$-модульного
$\mathbb{C}[w_0SU_{1,1}]_q$-бимодуля $\mathscr{D}(w_0SU_{1,1})'_q$. Более
того, на $\mathscr{D}(w_0SU_{1,1})_q$ переносится по непрерывности
структура $U_q\mathfrak{su}_{1,1}$-модульной алгебры.

Описанное в предыдущем \itemiiiе вложение $U_q\mathfrak{su}_{1,1}$-модульных
алгебр
 $\mathcal{I}: {\rm Pol}(\mathbb{C})_q \hookrightarrow
 \mathbb{C}[w_0SU_{1,1}]_{q,x}$
 допускает продолжение по непрерывности до вложения
 $U_q\mathfrak{su}_{1,1}$-модульных алгебр
 $\mathcal{I}: \mathscr{D}(\mathbb{D})_q \hookrightarrow
 \mathscr{D}(w_0SU_{1,1})_q$, а также до вложения пространств
 обобщенных функций $\mathscr{D}(\mathbb{D})'_q \hookrightarrow
 \mathscr{D}(w_0SU_{1,1})'_q$.

\bigskip В формулировке следующего утверждения участвует $U_q\mathfrak{sl}_2$-модульная
 алгебра $\mathscr{D}(\widetilde{\Xi})_q$ финитных
 функций на квантовом конусе, введенная в \itemiiiе \ref{q-cone}.
\begin{proposition}\label{inv_int_x}  Линейный функционал
на $\mathscr{D}(w_0SU_{1,1})_q$, определяемый равенством $
\int\limits_{(w_0SU_{1,1})_q}\,f\, d\nu\;=\;
(1-q^2)\,\sum\limits_{n=0}^\infty \psi_{0000}(q^{-2n})q^{-2n}, $ является
положительным и $U_q\mathfrak{sl}_2$-инвариантным. Линейный функционал на
$\mathscr{D}(\widetilde{\Xi})_q$, определяемый равенством $
\int\limits_{\widetilde{\Xi}_q}\,f\, d\nu\;=\;
(1-q^2)\,\sum\limits_{n=-\infty}^\infty \psi_{0000}(q^{-2n})q^{-2n}, $
является положительным и $U_q\mathfrak{sl}_2$-ин\-вариантным.
\end{proposition}

{\bf Доказательство.} Ограничимся первым утверждением. Второе доказывается
аналогично.

Воспользуемся представлениями $\pi_+,\rho$, введенными в
\itemiiiе \ref{SL_2_sl_2}. Представление $\pi_+\otimes \rho$ продолжается
до точного $*$-представления алгебры $\mathscr{D}(w_0SU_{1,1})_q$.
Положительность рассматриваемого линейного функционала вытекает из
равенства
$$\int\limits_{(w_0SU_{1,1})_q}\,f\, d\nu\;=\;(1-q^2)\, {\rm
tr} (\pi_+\otimes \rho) (f\cdot x),\qquad f \in \mathscr{D}(w_0SU_{1,1})_q
$$
(след корректно определен потому, что
 матрица оператора $(\pi_+\otimes \rho) (f\cdot x)$
в базисе $\{e_j \otimes z^k\}_{j \in \mathbb{Z}_+,\,k \in \mathbb{Z}}$
имеет конечное число ненулевых элементов на главной диагонали).

Пусть $\psi(x)$ -- функция с конечным носителем на множестве
$q^{-2\mathbb{Z}_+}$. Используя \eqref{r_reg1} -- \eqref{r_reg3}, нетрудно
показать, что
$$\int\limits_{(w_0SU_{1,1})_q}(K^{\pm 1}-1)\psi(x)\,d\nu=0,$$
$$
\int\limits_{(w_0SU_{1,1})_q}E(\,t_{12}\psi(x)t_{22}\,)\,d\nu\,=\,
\int\limits_{(w_0SU_{1,1})_q}F(\,t_{11}\psi(x)t_{21}\,)\,d\nu\,=\,0.
$$
Утверждение об $U_q\mathfrak{sl}_2$-инвариантности рассматриваемого
линейного функционала легко следует из этих равенств. \hfill $\square$

\bigskip Завершим \itemiii\ доказательством предложения
\ref{metric_for_bundle}.

 Инвариантность полуторалинейнго отображения
\begin{equation*}
{\mathbf h}_\lambda:\,\mathscr{D}(\mathbb{D})_{q,\lambda}\times
\mathscr{D}(\mathbb{D})_{q,\lambda}\to \mathscr{D}(\mathbb{D})_q,\;\;
{\mathbf h}_\lambda:\,
 v_\lambda\,f_1\,\times\,v_\lambda\,f_2 \mapsto f_2^*\,(1-zz^*)^\lambda\,f_1
\end{equation*}
достаточно установить в частном случае $\lambda \in \mathbb{Z}$. С помощью
вложения
 $$
 \mathcal{I}_\lambda:\,\mathscr{D}(\mathbb{D})_{q,\lambda} \hookrightarrow
 \mathscr{D}(w_0SU_{1,1})_q, \qquad \mathcal{I}_\lambda:\,
 v_\lambda\,f \mapsto t_{12}^{-\lambda}\; \mathcal{I}(f),
 $$
 доказательство инвариантности
сводится к равенству $
\mathcal{I}(y^\lambda)=(t_{12}^{*})^{-\lambda}t_{12}^\lambda, $
 поскольку $\mathscr{D}(w_0SU_{1,1})'_q$ является
$U_q\mathfrak{sl}_2$-модульным $\mathscr{D}(w_0SU_{1,1})_q$-бимодулем.

 Положительность при $\lambda \in \mathbb{R}$ устанавливается без труда, см.
доказательство предложения \ref{hermitian_metrics}.

\subsubsection{Доказательство леммы \ref{inv_G_m}.} \label{proof_inv_G_m}

Введем алгебры ядер
 $$\mathbb{C}[w_0SU_{1,1}\times w_0SU_{1,1}]_q\stackrel{\rm def}{=}
 \mathbb{C}[w_0SU_{1,1}]_q^{\rm op}\otimes \mathbb{C}[w_0SU_{1,1}]_q,$$
$$\mathbb{C}[w_0SU_{1,1}\times \widetilde{\Xi}]_q\stackrel{\rm def}{=}
 \mathbb{C}[w_0SU_{1,1}]_q^{\rm op}\otimes \mathbb{C}[\widetilde{\Xi}]_q,\quad
\mathbb{C}[\widetilde{\Xi}\times \widetilde{\Xi}]_q\stackrel{\rm def}{=}
\mathbb{C}[\widetilde{\Xi}]_q^{\rm op}\otimes \mathbb{C}[\widetilde{\Xi}]_q
 $$
 и  обозначения $t_{ij},\tau_{ij}$ для их образующих $t_{ij}\otimes 1,1
\otimes t_{ij}$ соответственно.

Из определений следует

\begin{lemma}\label{trivial_lemma}
Ядра
$$
k_{11}=t_{11}\tau_{22}-qt_{12}\tau_{21};\qquad
k_{12}=-q^{-1}t_{11}\tau_{12}+t_{12}\tau_{11};
$$
$$
k_{21}=t_{21}\tau_{22}-qt_{22}\tau_{21};\qquad
k_{22}=-q^{-1}t_{21}\tau_{12}+t_{22}\tau_{11}
$$
$U_q\mathfrak{sl}_2$-инвариантны.
\end{lemma}

\medskip

\begin{proposition} В $\mathbb{C}[w_0SU_{1,1}\times w_0SU_{1,1}]_q$
выполняются следующие соотношения:
\begin{equation}
\left \{\begin{array}{lcl}k_{11}k_{12}=q^{-1}k_{12}k_{11}&,& \qquad
k_{21}k_{22}=q^{-1}k_{22}k_{21},\\ k_{11}k_{21}=q^{-1}k_{21}k_{11}&,& \qquad
k_{12}k_{22}=q^{-1}k_{22}k_{12},\\ k_{12}k_{21}=k_{21}k_{12}&,&
\qquad k_{11}k_{22}-k_{22}k_{11}=(q^{-1}-q)k_{12}k_{21}\\
k_{22}k_{11}-qk_{12}k_{21}=1 \end{array}\right.
\end{equation}
\end{proposition}

\smallskip

{\bf Доказательство.} Рассмотрим отображение
\begin{equation}\label{pi_SL_2}
\mathbb{C}[SL_2]_q\to\mathbb{C}[SL_2]_q^{\rm
op}\otimes\mathbb{C}[SL_2]_q,\qquad f\mapsto({\rm id}\otimes S)\Delta(f),
\end{equation}
где $\Delta$ -- коумножение, а $S$ -- антипод в алгебре Хопфа
$\mathbb{C}[SL_2]_q$. Это отображение является антиавтоморфизмом и
$t_{ij}\mapsto k_{ij}$. \hfill $\square$

\begin{corollary}\label{for_cone}
В $\mathbb{C}[w_0SU_{1,1}\times \widetilde{\Xi}]_q$ и в
$\mathbb{C}[\widetilde{\Xi}\times \widetilde{\Xi}]_q$ выполняются следующие
соотношения:
\begin{equation}
\left \{\begin{array}{lcl}k_{11}k_{12}=q^{-1}k_{12}k_{11}&,& \qquad
k_{21}k_{22}=q^{-1}k_{22}k_{21},\\ k_{11}k_{21}=q^{-1}k_{21}k_{11}&,&
\qquad k_{12}k_{22}=q^{-1}k_{22}k_{12},\\ k_{12}k_{21}=k_{21}k_{12}&,&
\qquad k_{11}k_{22}-k_{22}k_{11}=(q^{-1}-q)k_{12}k_{21}\\
k_{22}k_{11}-qk_{12}k_{21}=0.\end{array}\right.
\end{equation}
\end{corollary}

\medskip

Введем в рассмотрение элементы $x=t_{12}^*t_{12}$,
$\xi=\tau_{12}\tau_{12}^*$ алгебр ядер $\mathbb{C}[w_0SU_{1,1}\times
w_0SU_{1,1}]_q$, $\mathbb{C}[w_0SU_{1,1}\times \widetilde{\Xi}]_q$,
$\mathbb{C}[\widetilde{\Xi}\times \widetilde{\Xi}]_q$
 и локализации $\mathbb{C}[w_0SU_{1,1}\times w_0SU_{1,1}]_{q,x\xi}$,
 $\mathbb{C}[w_0SU_{1,1}\times \widetilde{\Xi}]_{q,x\xi}$,
 $\mathbb{C}[\widetilde{\Xi}\times \widetilde{\Xi}]_{q,x\xi}$ этих алгебр
по мультипликативным множествам $(x\xi)^{\mathbb{Z}_+}$. Пусть
$$
z=qt_{12}^{-1}t_{11},\quad z^*=t_{22}t_{21}^{-1},\quad \zeta=q
\tau_{11}\tau_{12}^{-1},\quad \zeta^*=\tau_{21}^{-1}\tau_{22}.
$$

\medskip
\begin{lemma}\label{to_powers}   В
$\mathbb{C}[w_0SU_{1,1}\times w_0SU_{1,1}]_{q,x\xi}$,
 $\mathbb{C}[w_0SU_{1,1}\times \widetilde{\Xi}]_{q,x\xi}$,
 $\mathbb{C}[\widetilde{\Xi}\times \widetilde{\Xi}]_{q,x\xi}$
при всех $l \in \mathbb{Z}$
\begin{equation}\label{kl2}
k_{22}^lk_{11}^l=q^{-2l}\xi^l \sum_{j=0}^\infty
\frac{(q^{-2l};q^2)_j}{(q^2;q^2)_j}(q^{2(l+1)}z^*\zeta)^j \sum_{m=0}^\infty
\frac{(q^{-2l};q^2)_m}{(q^2;q^2)_m}(q^{2l}z \zeta^*)^m x^l.
\end{equation}
\end{lemma}

  {\bf Доказательство.} С помощью равенств
$$
t_{12}\tau_{21}(z \zeta^*)=q^2(z \zeta^*)t_{12}\tau_{21};\qquad (\zeta
z^*)t_{12}\tau_{21}=q^2t_{21}\tau_{12}(\zeta z^*),
$$
получаем:
\begin{equation}\label{kl1}
\begin{gathered}
k_{11}^l=(-qt_{12}\tau_{21}(1-q^{-2}z
\zeta^*))^l=(-qt_{12}\tau_{21})^l(q^{-2}z \zeta^*;q^{-2})_l,
\\ k_{22}^l=((1-z^*\zeta)(-q^{-1}t_{21}\tau_{12}))^l=
(z^*\zeta;q^{-2})_l(-q^{-1}t_{21}\tau_{12})^l.
\end{gathered}
\end{equation}
Значит,
\begin{multline*}
k_{22}^lk_{11}^l=
(z^*\zeta;q^{-2})_l(t_{12}t_{21})^l(\tau_{12}\tau_{21})^l(q^{-2}z
\zeta^*;q^{-2})_l=
\\ =q^{-2l}\xi^l(q^2z^*\zeta;q^2)_l(z\zeta^*;q^2)_lx^l.
\end{multline*}
Остается воспользоваться равенством
\begin{equation}\label{tqn}
(t;q)_n=\sum_{j=0}^n \frac{(q^{-n};q)_j}{(q;q)_j}q^{j(n+1)}t^j
\end{equation}
(см. \eqref{qbinom}) и его следствиями
$$
(q^2z^*\zeta;q^2)_l= \sum_{n=0}^\infty
\frac{(q^{-2l};q^2)_n}{(q^2;q^2)_n}q^{2(l+1)n}(z^*\zeta)^n,
$$
$$
\qquad (z \zeta^*;q^2)_l= \sum_{n=0}^\infty
\frac{(q^{-2l};q^2)_n}{(q^2;q^2)_n}q^{2ln}(z \zeta^*)^n.\eqno \square
$$

\bigskip Получим аналитическое
продолжение ядра $k_{22}^lk_{11}^l$ по параметру $l$.

  Для этого нам понадобится топологическое векторное пространство обобщенных
ядер. Напомним, что в \itemiiiе \ref{finite_w_0SU_11} векторное
пространство $\mathbb{C}[w_0SU_{1,1}]_q$ было наделено слабейшей из
топологий, в которых непрерывны линейные функционалы
$\mathcal{L}_{a,b,c,d,n}$, введенные равенством
\eqref{topology_SU_1_1_new}. Рассуждая аналогично, наделим векторное
пространство $\mathbb{C}[w_0SU_{1,1}\times w_0SU_{1,1}]_q$ слабейшей из
топологий, в которых непрерывны линейные функционалы
$\mathcal{L}_{a',b',c',d',n'}\otimes \mathcal{L}_{a'',b'',c'',d'',n''}$, и
введем в рассмотрение пополнение, называемое пространством обобщенных ядер
и обозначаемое $\mathscr{D}(w_0SU_{1,1}\times w_0SU_{1,1})'_q$.

 Так же, как в \itemiiiе \ref{finite_w_0SU_11}, получаем описание этого
 пополнения в терминах формальных рядов и продолжение на него по
 непрерывности действия алгебры $U_q\mathfrak{sl}_2$.
Алгебра финитных ядер
 $$ \mathscr{D}(w_0SU_{1,1}\times w_0SU_{1,1})_q
 \,\stackrel{\rm def}{=}\, \mathscr{D}(w_0SU_{1,1})_q^{\rm op} \otimes
 \mathscr{D}(w_0SU_{1,1})_q$$
 является плотным линейным подмногообразием в
  $\mathscr{D}(w_0SU_{1,1}\times  w_0SU_{1,1})'_q$.

\medskip
 Покажем, что $k_{22}^lk_{11}^l$
 является обобщенным  ядром. Действительно, как следует из (\ref{kl2}),
$$k_{22}^lk_{11}^l=q^{-2l}\sum_{j=0}^\infty \sum_{m=0}^\infty
\frac{(q^{-2l};q^2)_j}{(q^2;q^2)_j}
\frac{(q^{-2l};q^2)_m}{(q^2;q^2)_m}q^{2j(l+1)+2ml}z^{*j}z^mx^l \xi^l \zeta^j
\zeta^{*m}=$$
$$=q^{-2l}\sum_{j<m}\frac{(q^{-2l};q^2)_j}{(q^2;q^2)_j}
\frac{(q^{-2l};q^2)_m}{(q^2;q^2)_m}q^{2j(l+1)+2ml}(z^{*j}z^j)z^{m-j}x^l
\xi^l(\zeta^j \zeta^{*j})\zeta^{*(m-j)}+$$
$$+q^{-2l}\sum_{j=0}^\infty
\frac{(q^{-2l};q^2)_j^2}{(q^2;q^2)_j^2}q^{4jl+2j}(z^{*j}z^j)x^l
\xi^l(\zeta^j \zeta^{*j})+$$
$$+q^{-2l}\sum_{m<j}\frac{(q^{-2l};q^2)_j}{(q^2;q^2)_j}
\frac{(q^{-2l};q^2)_m}{(q^2;q^2)_m}q^{2j(l+1)+2ml}z^{*(j-m)}(z^{*m}z^m)x^l
\xi^l \zeta^{j-m}(\zeta^m \zeta^{*m}).$$

\medskip \stepcounter{theorem}

Займемся аналитическим продолжением. Рассмотрим вектор-функ\-цию
$\mathscr{K}(\lambda)$ комплексной переменной $\lambda$ со значениями в
$\mathscr{D}(w_0SU_{1,1}\times w_0SU_{1,1})'_q$. Назовём такую
вектор-функцию полиномиальной, если для любой финитной функции
$\psi\in\mathscr{D}(w_0SU_{1,1}\times w_0SU_{1,1})'_q$ интеграл
$\int\limits_{(w_0SU_{1,1})_q}\int\limits_{(w_0SU_{1,1})_q}\,
\mathscr{K}(\lambda)\,\psi\, d\nu \,d\nu$ является полиномом Лорана
переменной $\lambda$.

\medskip
\begin{proposition}\label{p3.6.7_1}
 Существует единственная
такая полиномиальная вектор-функция $\mathscr{K}(\lambda)$ со значениями в
$\mathscr{D}(w_0SU_{1,1}\times w_0SU_{1,1})'_q$, что
$\mathscr{K}(q^{2l})=k_{22}^lk_{11}^l$ для всех $l \in{\mathbb Z}_+$.
\end{proposition}

\smallskip

{\bf Доказательство.} Единственность рациональной функции $f(\lambda)$ с
заданными значениями в точках множества $q^{2\mathbb{Z}_+}$, очевидна.

Рассмотрим финитную функцию $\psi$ и заменим $k_{22}^lk_{11}^l$ в интеграле
$\int\limits_{(w_0SU_{1,1})_q}\int\limits_{(w_0SU_{1,1})_q}\,
\mathscr{K}(\lambda)\,\psi\, d\nu \,d\nu$ правой частью равенства
(\ref{kl2}). Только конечное число слагаемых из (\ref{kl2}) дают ненулевой
вклад в этот интеграл. Остаётся заметить, что каждое из этих слагаемых
является полиномом Лорана переменной $q^{2l}$. \hfill $\square$

\medskip

Сохраним обозначение $k_{22}^lk_{11}^l$ для обобщенного ядра, описанного в
предыдущем предложении. Как следует из леммы \ref{trivial_lemma} и
предложения \ref{IntKer}, это ядро $U_q\mathfrak{sl}_2$-инвариантно при
$l\in\mathbb{Z}_+$, а, следовательно, и в общем случае.

Утверждение леммы \ref{inv_G_m} вытекает из инвариантности ядра
$k_{22}^lk_{11}^l$ при $l=-m$ и леммы \ref{to_powers}, поскольку
\\$q^{-2m}\,k_{22}^{-m}k_{11}^{-m}=(1-\zeta \zeta^*)^m \cdot$
$$\cdot \sum_{j=0}^\infty
\frac{(q^{2m};q^2)_j}{(q^2;q^2)_j}(q^{2(-m+1)}z^*\zeta)^j \cdot
\sum_{n=0}^\infty \frac{(q^{2m};q^2)_n}{(q^2;q^2)_n}(q^{-2m}z \zeta^*)^n
(1-z^*z)^m=G_m.$$

\medskip Отметим, что в этом доказательстве неявно использованы каноническое вложение
$U_q\mathfrak{sl}_2$-модульных алгебр и вложение
$U_q\mathfrak{sl}_2$-модулей $\mathscr{D}(\mathbb{D}\times \mathbb{D})'_q
\to \mathscr{D}(w_0SU_{1,1}\times w_0SU_{1,1})'_q$, получаемое продолжением
по непрерывности.

\begin{remark}\label{cone-cone} Заменяя в предыдущих рассуждениях $\mathbb{C}[w_0SU_{1,1}\times
w_0SU_{1,1}]_q$ на $\mathbb{C}[w_0SU_{1,1}\times \widetilde{\Xi}]_q$, либо
на $\mathbb{C}[\widetilde{\Xi}\times \widetilde{\Xi}]_q$, получаем
 утверждение об
$U_q\mathfrak{sl}_2$-инвариантности при всех $l \in \mathbb{C}$ обобщенных
функций
$$
\mathscr{K}_l=\xi^l \sum_{j=0}^\infty
\frac{(q^{-2l};q^2)_j}{(q^2;q^2)_j}(q^{2(l+1)}z^*\zeta)^j \sum_{m=0}^\infty
\frac{(q^{-2l};q^2)_m}{(q^2;q^2)_m}(q^{2l}z \zeta^*)^m x^l
$$
на квантовых аналогах декартовых произведений $w_0SU_{1,1}\times
\widetilde{\Xi}$ и $\widetilde{\Xi}\times \widetilde{\Xi}$.
\end{remark}

\subsubsection{Стандартные сплетающие операторы и ${\mathbf c}$-функция
Хариш-Чандры.}\label{c-function_sl_2}

В классическом случае $q=1$ к $c$-функции Хариш-Чандры можно прийти, изучая
операторы, сплетающие представления основной серии. Опишем $q$-аналог этого
подхода.

Напомним обозначение $ (a;q^2)_l=(a;q^2)_\infty/ (aq^{2l};q^2)_\infty$.

Рассмотрим $U_q\mathfrak{sl}_2$-модули
$$
\mathbb{C}[\partial\mathbb{D}]^{(l)}_q,\qquad
-\frac{\pi}{2\log(q^{-1})}\le\,{\rm Im}\,l\,\le\,\frac{\pi}{2\log(q^{-1})},
$$
основной сферической серии, введенные в \itemiiiе \ref{Fourier_sl_2}.

\begin{proposition}\label{intertw_op_sl_2} При всех
$l\notin \{-1,-2,-3,\ldots\}$ линейный оператор
\begin{equation}\label{def_I_l}
 I_l:\, f(e^{i\varphi})\mapsto \int\limits_0^{2\pi}
 (e^{i(\varphi-\theta)};q^2)_{l}\;
(q^2e^{i(\theta-\varphi)};q^2)_{l}\;f(e^{i\theta})\;\frac{d\theta}{2\pi}
\end{equation}
является морфизмом $U_q\mathfrak{sl}_2$-модулей
 $$I_l: \mathbb{C}[\partial \mathbb{D}]^{(-l-1)}_q \to
 \mathbb{C}[\partial \mathbb{D}]^{(l)}_q.$$
\end{proposition}

{\bf Доказательство.} Можно ограничиться случаем ${\rm Re} l>0$. Требуемое
утверждение легко получить, повторяя доказательство предложения
\ref{t1.5.1}. Действительно, во-первых, имеют место равенства
\eqref{qbinom_P1}, \eqref{qbinom_P2}, и обобщенная функция $\mathscr{K}_l$
 на квантовом конусе $\Xi\times \Xi$ является $U_q\mathfrak{sl}_2$
 инвариантной, см. замечание \ref{cone-cone}. Во-вторых,
 действие алгебры $U_q\mathfrak{sl}_2$
перенесено в $\mathbb{C}[\partial \mathbb{D}]^{(-l-1)}_q$ с квантового
конуса при помощи отображения $f\mapsto f\,x^{-l-1}$. В-третьих,
произведение $\mathscr{K}_l\,(1\otimes f\,x^{-l-1})$ имеет степень
однородности $-1$ по второму тензорному сомножителю. В-четвертых, линейный
функционал \hbox{$\eta: \left(\sum_{m=-\infty}^\infty a_m z^m
\right)x^{-1}\,\mapsto\,a_0$} на подпространстве обобщенных функций степени
однородности $-1$ на квантовом конусе является
$U_q\mathfrak{sl}_2$-инвариантным.
 \hfill $\square$

\medskip Следуя традиции, см., например, \cite{JoWall}, назовем линейный
оператор $I_l$ в векторном пространстве $\mathbb{C}[\partial \mathbb{D}]$
стандартным сплетающим оператором.

Напомним, что $c$-функция Хариш-Чандры \IND{$q$-аналог ! $c$-функции
Хариш-Чандры} была введена равенством \eqref{c_l}:
\begin{equation*}
c(l)=\dfrac{(q^{2(l+1)};q^2)_\infty^2}
{(q^{2(l+1)};q^2)_\infty(q^2;q^2)_\infty}\,=\,
\frac{\Gamma_{q^2}(2l+1)}{(\Gamma_{q^2}(l+1))^2}.
\end{equation*}

\medskip

 Функция $c(l)$ неявно зависит от параметра $q$
 и в формальном
 пределе  $q\to 1$ несущественно отличается от обычной $c$-функции
 Хариш-Чандры \cite[стр. 55]{Helg1}, как следует из  формулы удвоения
 Лежандра \cite[стр. 19]{Bateman1}
 $$\Gamma(2l+1)=2^{2l}\,\pi^{-1/2}\,\Gamma(l+\frac{1}{2})\,
  \Gamma(l+1)$$
(ее $q$-аналог приведен в \cite[стр. 38]{GR}).

\medskip
\begin{proposition}\label{coeff_I_l}
При всех $l\notin \{-1,-2,-3,\ldots\}$, $n \in \mathbb{Z}$
\begin{equation}\label{c_a_sl_2}
I_l\, z^n \,=\, c(l)a(l,n)\, z^n,
\end{equation}
 где
\begin{equation}\label{a-l-n}
a(l,n)\,=\,q^{-(2l+1)n}\,\begin{cases}
 \prod\limits_{j=0}^{n-1}
       \frac{q^{-(n-l-j-1)}-q^{n-l-j-1}}{q^{-(n+l-j)}-q^{n+l-j}},
       &\qquad n>0,\\
  1,   &\qquad n=0,\\
  \prod\limits_{j=0}^{n-1}
       \frac{q^{-(n+l+j+1)}-q^{n+l+j+1}}{q^{-(n-l+j)}-q^{n-l+j}},
  &\qquad n<0.
       \end{cases}
\end{equation}
\end{proposition}

{\bf Доказательство.} Достаточно рассмотреть случай $l\notin \mathbb{Z}$.
Как следует из леммы \ref{univ_pi_l}, при $l\notin \mathbb{Z}$ существует и
единствен такой изоморфизм $U_q\mathfrak{sl}_2$-модулей $A_l:
\mathbb{C}[\partial \mathbb{D}]^{(-l-1)}_q \to
 \mathbb{C}[\partial \mathbb{D}]^{(l)}_q$, что $A_l\; 1=1$.
Используя явный вид изоморфизмов этих $U_q\mathfrak{sl}_2$-модулей и
модулей основной сферической серии, нетрудно показать, что
 $A_l\;z^n=a(l,n) z^n$ при всех $n\in \mathbb{Z}$.

Изоморфизм $U_q\mathfrak{sl}_2$-модулей
 $\mathbb{C}[\partial \mathbb{D}]^{(-l-1)}_q\to
 \mathbb{C}[\partial \mathbb{D}]^{(l)}_q$
   единствен
с точностью до числового множителя. Значит, из предложения
 \ref{intertw_op_sl_2} следует равенство
 $I_l\,z^n\,=\,{\rm const}(l)\,a(l,n)\,z^n$.
 Остается показать, что
 ${\rm const}(l)=c(l)$, то есть доказать равенство $I_l\;1=c(l)$.
 Но, как следует из \eqref{def_I_l},
$$
I_l\; 1\,=\, \sum\limits_{k=o} \frac{(q^{-2l};q^2)_k(q^{-2l};q^2)_k}
{(q^2;q^2)_k(q^2;q^2)_k} q^{2(2l+1)j}\,=\,
\sideset{_2}{_1}{\mathop{\phi}}\left[\begin{array}{c}{q^{-2l},q^{-2l}}\\
q^2\end{array};q^2,q^{2(2l+1)}\right].
$$
Используя $q$-аналог \eqref{Gauss} формулы Гаусса и определение \eqref{c_l}
функции $c(l)$, получаем:
$$
I_l\cdot 1\,=\,\frac{(q^{2(l+1)};q^2)_\infty(q^{2(l+1)};q^2)_\infty}
{(q^2;q^2)_\infty(q^{2(2l+1)};q^2)_\infty}\,=\,c(l).\eqno \square
$$

\medskip
\begin{corollary}\label{complementary_series_sl_2}
 При $-1 <l< 0$ эрмитова форма
\begin{equation}\label{h-10} (f_1,f_2)_l=\int\limits_0^{2\pi}\,
\overline{f_2(e^{i\theta})}(I_{-l-1}f_1)(e^{i\theta})\,\frac{d\theta}{2\pi}
\end{equation}
в $\mathbb{C}[\partial \mathbb{D}]^{(l)}_q$ является положительной и
$U_q\mathfrak{su}_{1,1}$-инвариантной.
\end{corollary}

{\bf Доказательство.} Из определения $U_q\mathfrak{sl}_2$-модулей $
\mathbb{C}[\partial \mathbb{D}]^{(-l-1)}_q$, $\mathbb{C}[\partial
\mathbb{D}]^{(l)}_q$ и из предложения \ref{homogeneous_int} вытекает
$U_q\mathfrak{su}_{1,1}$-инвариантность полуторалинейной формы $
 \mathbb{C}[\partial \mathbb{D}]^{(-l-1)}_q\times
   \mathbb{C}[\partial \mathbb{D}]^{(l)}_q \to \mathbb{C},
$
 определяемой  равенством
 $$ \langle f_1, f_2 \rangle_l\,=\, \int\limits_0^{2\pi}\,
 \overline{f_2(e^{i\theta})}\,f_1(e^{i\theta})\,\frac{d\theta}{2\pi}.$$
 Другими словами,
$
 \langle \xi\, f_1, f_2 \rangle_l\,=\,\langle f_1, \xi^*\,f_2 \rangle_l$ при всех
 $\xi \in U_q\mathfrak{su}_{1,1}$.
Значит, эрмитова форма \eqref{h-10} является
$U_q\mathfrak{su}_{1,1}$-инвариантной. Остается воспользоваться равенствами
\eqref{coeff_I_l},\eqref{a-l-n} для доказательства ее положительности.
\hfill $\square$

\bigskip

Дополнительной серией \IND{дополнительная серия} называют множество простых
унитаризуемых $U_q\mathfrak{su}_{1,1}$-модулей
$\mathbb{C}[\partial\mathbb{D}]^{(l)}_q$ с $-\frac12<l<0$.

Аналогично доказывается унитаризуемость $U_q\mathfrak{su}_{1,1}$-модулей
$\mathbb{C}[\partial\mathbb{D}]^{(l)}_q$ при ${\rm
Im}l=\pm\frac{\pi}{2\log(q^{-1})}$. Сферической странной серией
\IND{сферичесая странная серия} называют множество простых унитаризуемых
$U_q\mathfrak{su}_{1,1}$-модулей $\mathbb{C}[\partial\mathbb{D}]^{(l)}_q$ с
${\rm Im}l\,=\,\frac{\pi}{2\log(q^{-1})}$ и ${\rm Re}\,l>0$.

\begin{remark}\label{complete_list_sl_2}
Эти серии вместе с основной унитарной сферической серией
$$
\mathbb{C}[\partial \mathbb{D}]^{(l)}_q,\quad\qquad{\rm
Re}\,l=-\frac12\quad\&\quad 0\le\,{\rm Im}l\le\,\frac{\pi}{2\log(q^{-1})},
$$
и тривиальным $U_q\mathfrak{su}_{1,1}$-модулем $\mathbb{C}$ образуют полный
список простых сферических унитаризуемых $U_q\mathfrak{su}_{1,1}$-модулей,
см. работу Бурбана и Климыка \cite{BurKli}, где приведен полный список
простых весовых унитаризуемых $U_q\mathfrak{su}_{1,1}$-модулей.
\end{remark}

\bigskip

Покажем, что, как и в классическом случае $q=1$, $c$-функция Хариш-Чандры
возникает при изучении асимптотики
 $|z| \to 1$ собственных функций инвариантного Лапласиана $\square_q$.
Это объясняет появление $c$-функции в формуле для ядра интегрального
оператора $(\lambda-\square_q)^{-1}$ и, следовательно \cite[стр.
77]{Dunford2}, в формуле \eqref{Pl} для меры Планшереля.

Во-первых, как показано в \itemiiiе \ref{radial_part}, функция
\begin{equation*}
\Phi_l(x)=\sideset{_3}{_2}{\mathop{\phi}}
\left[\begin{array}{c}x,q^{-2l},q^{2(l+1)}\\
q^2,0\end{array};q^2,q^2\right],\qquad x=(1-zz^*)^{-1},
\end{equation*}
является собственной функцией оператора $\square_q^{(0)}$, отвечающей
собственному значению \eqref{Laplace_eigenvalue}. Кроме того,
$\Phi_l(1)=1$. Во-вторых, как показано в \itemiiiе \ref{eigen_functions},
обобщенная функция $\psi(x)=\int\limits_0^{2\pi}
P_{l+1}(z,e^{i\theta})f(e^{i\theta})\frac{d\theta}{2\pi}$ является
собственной функцией оператора $\square_q^{(0)}$ с тем же собственным
значением и равна 1 при $x=1$. Значит,
$$
\Phi_l(x)\;=\;\int\limits_0^{2\pi}
P_{l+1}(z,e^{i\theta})f(e^{i\theta})\frac{d\theta}{2\pi}.
$$
Это -- $q$-аналог хорошо известного интегрального представления сферической
функции \cite[стр.54]{Helg1}.

\begin{proposition}\label{l1.5.4} 1. При
${\rm Re}\,l>- \frac{1}{2}$
$$ \lim \limits_{x \in q^{-2{\mathbb Z}_+},\, x \to \infty}x^{-l}
 \Phi_l(x)\;=\;\frac{\Gamma_{q^2}(2l+1)}{\Gamma_{q^2}^2(l+1)}.
$$
2. При ${\rm Re}\,l<-\frac{1}{2}$
$$
\lim \limits_{x \in q^{-2{\mathbb Z}_+},\,x \to
\infty}x^{l+1}\Phi_l(x)\;=\;\frac{\Gamma_{q^2}(-2l-1)}{\Gamma_{q^2}^2(-l)}.
$$
\end{proposition}

{\bf Доказательство.} Соотношение $\Phi_l(x)=\Phi_{-1-l}(x)$ позволяет
ограничиться случаем ${\rm Re}\,l>-\frac{1}{2}$. Заменяя в \eqref{page45}
$q$, $b$, $c$, $z$ на $q^2$, $q^{-2l}$, $q^{-2l-2n}$, $q^{2l+2}$
соответственно, получаем:
$$\Phi_l(q^{-2n})=\frac{(q^{-2l-2n};q^2)_n}{(q^{-2n};q^2)_n}\cdot
\sideset{_2}{_1}{\mathop{\phi}}\left[\begin{array}{c}q^{-2n},q^{-2l}\\
q^{-2(l+n)}\end{array};q^2,q^{2(l+1)}\right] \sim
$$
$$
\sim q^{-2nl}\frac{(q^{2(l+1)};q^2)_\infty}{(q^2;q^2)_\infty}\cdot \,
 \sideset{_1}{_0}{\mathop{\phi}}\left[\begin{array}{c}q^{-2l}\\
-\end{array};q^2,q^{2(2l+1)}\right].
$$
Согласно $q$-биномиальной формуле \eqref{qbinom},
$$ \Phi_l(q^{-2n}) \sim
q^{-2nl}\frac{(q^{2(l+1)};q^2)_\infty}{(q^2;q^2)_\infty}\cdot
\frac{(q^{2(l+1)};q^2)_\infty}{(q^{2(2l+1)};q^2)_\infty},
$$
и остается воспользоваться определением $q$-гамма-функции. \hfill $\square$

\medskip Это доказательство нам сообщил Л. И. Корогодский.

\newpage

\section{ОСНОВЫ КВАНТОВОЙ ТЕОРИИ ОГРАНИЧЕННЫХ СИММЕТРИЧЕСКИХ ОБЛАСТЕЙ}
\label{basic_theory}

\subsection{Сводка   известных результатов о квантовых универсальных обертывающих
алгебрах}\label{UE}

\subsubsection{Алгебры Дринфельда-Джимбо.}\label{DrJ}

Одним из самых известных результатов математики 20-го столетия является
полученная Э.~Картаном классификация простых комплексных алгебр Ли
\cite{Helg, Hum_Lie}. С точностью до изоморфизма они определяются своими
матрицами Картана $\mathbf{a}=(a_{ij})_{i,j=1,2,\ldots,l}$, полный список
которых приведен в \cite{Bou4-6}. Каждая такая матрица $\mathbf{a}$
неразложима, является целочисленной и $a_{ii}=2$, $a_{ij}<0$ при $i\ne j$.
Кроме того, для некоторых чисел $d_i>0$, $i=1,2,\ldots,l$, матрица
$\mathbf{b}=(d_i a_{ij})_{i,j=1,2,\ldots,l}$ положительно определена.

 В дальнейшем
 числа $d_i$ будут выбираться натуральными и взаимно простыми.

   Простую комплексную алгебру Ли
$\mathfrak{g}$ с матрицей Картана $\mathbf{a}$ можно описать с помощью
 образующих и соотношений. Приведем такое описание применительно к
универсальной обертывающей алгебре
 $U\mathfrak{g}$ \cite[стр. 51]{Jant}. Ее можно задать  образующими
$\{H_i,E_i,F_i\}_{i=1,2,\ldots,l}$ и следующими определяющими
соотношениями:
$$
H_iH_j-H_jH_i=0,\quad E_iF_j-F_jE_i=\delta_{ij}H_i,
$$
$$
H_iE_j-E_jH_i=a_{ij}E_j,\quad H_iF_j-F_jH_i=-a_{ij}F_j
$$
при $i,j=1,2,\ldots,l$;
$$
\sum\limits_{m=0}^{1-a_{ij}}(-1)^m\binom{1-a_{ij}}{m}
E_i^{1-a_{ij}-m}\; E_j\; E_i^m=0,
$$
$$
\sum\limits_{m=0}^{1-a_{ij}}(-1)^m\binom{1-a_{ij}}{m}
F_i^{1-a_{ij}-m}\; F_j\; F_i^m=0
$$
при $i\neq j$.

Универсальная обертывающая алгебра $U\mathfrak{g}$ является алгеброй Хопфа
с коумножением $\Delta$:
$$
\Delta(H_i)=H_i\otimes 1+1\otimes H_i,\quad\Delta(E_i)=E_i\otimes1+1\otimes
E_i,\quad\Delta(F_i)=F_i\otimes1+1\otimes F_i,
$$
коединицей $\varepsilon$ и антиподом $S$:
$$\varepsilon(H_i)=\varepsilon(E_i)=\varepsilon(F_i)=0,$$
$$S(H_i)=-H_i,\quad S(E_i)=-E_i,\quad S(F_i)=-F_i,\quad i=1,2,\ldots,l.$$
Алгебру Ли $\mathfrak{g}$ можно отождествить с подалгеброй Ли, порожденной
множеством $\{H_i,E_i,F_i\}_{i=1,2,\ldots,l}$.

\bigskip Перейдем к понятию квантовой универсальной обертывающей алгебры,
\IND{квантовая ! универсальная обертывающая алгебра} следуя обозначениям
работ \cite{Lust1, DaDeCon, Rosso, Jant, BrGood}. Эта алгебра Хопфа введена
В.~Дринфельдом и М.~Джимбо в середине 80-х годов в существенно большей
общности, чем описано ниже.

Пусть ${\mathbf a}=(a_{ij})_{i,j=1,2,\ldots,l}$ -- матрица Картана простой
комплексной алгебры Ли и $q\in\mathbb{C}$ -- такое ненулевое комплексное
число, что $q^{2d_i}\ne 1$ при всех $1=1,2,\ldots,l$. Алгебра
$U_q\mathfrak{g}$ определяется своими образующими
$\{K_i,K_i^{-1},\,E_i,F_i\}_{i=1,2,\ldots,l}$ и соотношениями:
$$K_iK_j=K_jK_i,\quad K_iK_i^{-1}=K_i^{-1}K_i=1,$$
\begin{equation}\label{first_line}
K_iE_j=q_i^{a_{ij}}E_jK_i,\quad K_iF_j=q_i^{-a_{ij}}F_jK_i,
\end{equation}
\begin{equation}\label{first_line_new}
E_iF_j-F_jE_i=\delta_{ij}\,\frac{K_i-K_i^{-1}}{q_i-q_i^{-1}},
\end{equation}
\begin{equation}\label{Serre1}
\sum\limits_{m=0}^{1-a_{ij}}(-1)^m\begin{bmatrix}1-a_{ij}\\ m\end{bmatrix}_{q_i}
E_i^{1-a_{ij}-m}E_jE_i^m=0,
\end{equation}
\begin{equation}\label{Serre2}
\sum\limits_{m=0}^{1-a_{ij}}(-1)^m
\begin{bmatrix}1-a_{ij}\\ m\end{bmatrix}_{q_i}
F_i^{1-a_{ij}-m}F_jF_i^m=0,
\end{equation}
где $q_i=q^{d_i}$, $1\le i\le l$,
\begin{equation}\label{qbinom_VG}
\begin{bmatrix}m\\ n\end{bmatrix}_q=\frac{[m]_q!}{[n]_q![m-n]_q!},\quad
[n]_q!=[n]_q\ldots[2]_q[1]_q,\quad [n]_q=\frac{q^n-q^{-n}}{q-q^{-1}}.
\end{equation}

Доказывается \cite[глава 4]{Jant} существование и единственность такой
структуры алгебры Хопфа в $U_q\mathfrak{g}$, что при всех $i=1,2, \cdots, l$
\begin{equation}\label{comult_D}
  \Delta(E_i)=E_i\otimes 1+K_i\otimes E_i,\, \,
  \Delta(F_i)=F_i\otimes K_i^{-1}+1\otimes F_i, \, \,
  \Delta(K_i)=K_i\otimes K_i,
\end{equation}
$$
\varepsilon(E_i)=\varepsilon(F_i)=0,\quad \varepsilon(K_i)=1,
$$
\begin{equation}\label{antipode_DJ}
  S(E_i)=-K_i^{-1}E_i,\quad
  S(F_i)=-F_iK_i,\quad
  S(K_i)=K_i^{-1}.
\end{equation}

В частном случае $l=1$ получаем алгебру Хопфа $U_q\mathfrak{sl}_2$, в
терминах которой была описана симметрия квантового круга в главе
\ref{quantum_disc_part}. Как и в этой главе, мы в дальнейшем ограничимся
частным случаем $q\in(0,1)$. Только такие $q$ возникают в интересующих нас
задачах некоммутативного комплексного анализа.

Как и в случае $l=1$, определяющие соотношения для $U\mathfrak{g}$ могут
быть получены из соотношений Дринфельда--Джимбо подстановкой
\begin{equation}\label{KH}
  K_i=q_i^{H_i},\quad
  K_i^{-1}=q_i^{-H_i},\quad
  i=1,2,\ldots,l
\end{equation}
и формальным предельным переходом $q \rightarrow 1$.

\subsubsection{Весовые $U_q\mathfrak{g}$-модули.}\label{Weight}

Напомним некоторые понятия теории алгебр Ли \cite{Kac}. Рассмотрим простую
комплексную алгебру Ли $\mathfrak{g}$ ранга $l$, заданную образующими
$\{H_i, E_i, F_i\}_{i=1,2,\ldots,l}$ и хорошо известными определяющими
соотношениями, см. предыдущий \itemiii\ .

Линейную оболочку $\mathfrak{h}$ попарно коммутирующих элементов \\
$\{H_j\}_{j=1,2,\ldots,l}$ называют картановской подалгеброй,
\IND{картановская подалгебра} элементы биортогонального базиса
$\{\overline{\omega}_j\}_{j=1,2,\ldots,l}$ в $\mathfrak{h}^*$--
фундаментальными весами, \IND{фундаментальный вес} а абелеву группу
$P=\sum\limits_{i=1}^l\mathbb{Z}\,\overline{\omega}_i$ -- решеткой весов.
\IND{решетка весов}

Вещественная линейная оболочка фундаментальных весов обозначается
$\mathfrak{h}_\mathbb{R}^*$.
 Каждому  $\lambda \in \mathfrak{h}_\mathbb{R}^*$ отвечает
строка $(\lambda(H_1),\lambda(H_2),\ldots,\lambda(H_l))$ числовых отметок
веса $\lambda$, чем устанавливается изоморфизм
$\mathfrak{h}_\mathbb{R}^*\cong \mathbb{R}^l$.

Линейный функционал $\lambda$ на $\mathfrak{h}$ называют весом
представления $\pi$ алгебры Ли $\mathfrak{g}$, \IND{вес ! представления
алгебры Ли} если для некоторого ненулевого вектора $v$
$$\pi(\xi)v=\lambda(\xi)v,\qquad\xi\in\mathfrak{h}.$$
Участвующий в этом равенстве общий собственный вектор $v$ операторов
$\pi(\xi)$ называют весовым вектором веса $\lambda$. \IND{весовой ! вектор}

Ненулевые веса присоединенного представления
$\operatorname{ad}_{\xi}=[\xi,\cdot]$ алгебры Ли $\mathfrak{g}$ называются
корнями, \IND{корень алгебры Ли} а отвечающие им весовые подпространства --
корневыми подпространствами. \IND{корневое подпространство} Линейные
функционалы
\begin{equation}\label{alpha_definition}
\alpha_j\,=\,\sum\limits_{i=1}^l\;a_{ij}\,\overline{\omega}_i\qquad
j=1,2,\ldots,l
\end{equation}
являются корнями, поскольку $[H_i,E_j]\,=\,a_{ij}E_j $. Их называют
простыми корнями, \IND{корень алгебры Ли ! простой} и, как следует из
обратимости матрицы Картана, они образуют базис векторного пространства
$\mathfrak{h}^*_{\mathbb{R}}$.

Все корни принадлежат группе $Q=\sum\limits_{i=1}^l\mathbb{Z}\alpha_i$,
называемой решеткой корней. \IND{решетка корней} Используются обозначения
$$
P_+\,=\,\sum\limits_{j=1}^l\,\mathbb{Z}_+\,\overline{\omega}_j,\qquad
Q_+\,=\,\sum\limits_{j=1}^l\,\mathbb{Z}_+\,\alpha_j.
$$

 Векторное пространство
   $\mathfrak{h}^*_{\mathbb{R}}$ наделяют
 отношением частичного порядка, полагая
$\lambda'\ge \lambda''$ если и только если неотрицательны все координаты
вектора $\lambda'-\lambda''$ в базисе $\alpha_1,\alpha_2,\ldots,\alpha_l$.

Пусть $\mathfrak{n}^+$, $\mathfrak{n}^-$ -- подалгебры Ли, порожденные
множествами $\{E_1,E_2,\ldots,E_l\}$ и $\{F_1,F_2,\ldots,F_l\}$,
соответственно. Подпространства $\mathfrak{n}^\pm$ являются
$\operatorname{ad}_\xi$-инвариантными для всех $\xi\in\mathfrak{h}$.
Разложению
$\mathfrak{g}=\mathfrak{n}^-\oplus\mathfrak{h}\oplus\mathfrak{n}^-$
отвечает разложение множества корней $\Phi$ в объединение подмножества
$\Phi^+$ положительных корней и подмножества $\Phi^-$ отрицательных корней.
Очевидно, $\Pi=\{\alpha_1,\alpha_2,\cdots,\alpha_l\}\subset\Phi^+$.

Можно показать, что все корневые подпространства одномерны и что среди
корней есть наибольший.

\bigskip Перейдем от алгебр Ли к группам. Рассмотрим простую аффинную
алгебраическую группу $G$ с алгеброй Ли $\mathfrak{g}$. Решетка весов ее
рациональных представлений канонически изоморфна подгруппе $\mathscr{L}$,
заключенной между $Q$ и $P$: $Q\subset\mathscr{L}\subset P$. Пара
$(\mathfrak{g},\mathscr{L})$ определяет группу $G$ с точностью до
изоморфизма \cite[стр. 312]{Hum}. Другими словами, $G$ определяется с
точностью до изоморфизма парой $(\mathbf{a}, \mathscr{L})$, где
$\mathbf{a}$ -- матрица Картана алгебры Ли $\mathfrak{g}$.

\medskip Перейдем к квантовым группам.
Модуль $V$ над алгеброй $U_q\mathfrak{g}$ назовем $\mathscr{L}$-весовым,
\IND{весовой ! модуль} если
\begin{equation}\label{weight}
V=\newoplus\limits_{\mathbf{\lambda}\in\mathscr{L}}V_{\mathbf{\lambda}},
\qquad V_\lambda=\left\{v\in V\;\left|\;K_iv=q_i^{\lambda_i}v,\quad
i=1,2,\ldots,l\right.\right\},
\end{equation}
где $\mathbf{\lambda}=\sum\limits_{j=1}^l\lambda_j\overline{\omega}_j$.
Ненулевое слагаемое $V_{\mathbf{\lambda}}$ в этом разложении называют
весовым подпространством веса $\lambda$. \IND{весовое подпространство}
Полная подкатегория конечномерных $\mathscr{L}$-весовых
$U_q\mathfrak{g}$-модулей (наделенная некоторыми дополнительными
структурами) играет роль квантовой алгебраической группы $G$ \IND{квантовая
! алгебраическая группа} \cite{BrGood}.

\begin{example}\label{q-ad}
Присоединенное представление алгебры Хопфа \IND{присоединенное
представление алгебры Хопфа} $U_q\mathfrak{g}$ определяется равенством
${\rm ad}(\xi)\,\eta\,=\,\sum\limits_j\xi'_j\eta S(\xi'')$, где
$\Delta(\xi)=\sum\limits_j\xi'_j\otimes\xi_j''$, и наделяет
$U_q\mathfrak{g}$ структурой весового $U_q\mathfrak{g}$-модуля:
$U_q\mathfrak{g}=\newoplus\limits_{\lambda\in Q}(U_q\mathfrak{g})_\lambda$.
\end{example}

\medskip

Требование $\mathscr{L}\subset P$ может оказаться слишком ограничительным.
Ослабим его. Как и в \itemiiiе \ref{RM_sl_2}, весовым
$U_q\mathfrak{g}$-модулем будем называть $U_q\mathfrak{g}$-модуль,
удовлетворяющий условию (\ref{weight}) с
$\mathscr{L}=\mathfrak{h}^*_{\mathbb{R}}\cong\mathbb{R}^l$.
\footnote{Конечномерные весовые $U_q\mathfrak{g}$-модули обычно называют
$U_q\mathfrak{g}$-модулями типа $1$ \cite{Jant}.} Равенством
$H_i|_{V_{\mathbf{\lambda}}}=\lambda_i$ вводятся линейные операторы $H_i$ в
$V$, для которых
$$K_i^{\pm1}v=q_i^{\pm H_i}v,\qquad v\in V.$$

Полная подкатегория весовых $U_q\mathfrak{g}$-модулей замкнута относительно
перехода к подмодулям,\footnote{См., например, \cite{Kac}, предложение 1.5}
фактор-модулям и тензорным произведениям. Каждый $U_q\mathfrak{g}$-модуль
обладает наибольшим весовым подмодулем, что позволяет ввести понятие
двойственного весового $U_q\mathfrak{g}$-модуля. \IND{двойственный !
весовой модуль} Именно, если
$V=\bigoplus\limits_{\mathbf{\lambda}}V_{\mathbf{\lambda}}$, то
$V^*\stackrel{\operatorname{def}}{=}
\bigoplus\limits_{\mathbf{\lambda}}(V_{\mathbf{\lambda}})^*$.

\bigskip

Каждому элементу $\lambda=\sum\limits_{i=1}^l n_i\alpha_i$ решетки корней
$Q$ сопоставим элемент $$K_\lambda =K_1^{n_1}K_2^{n_2}\cdots K_l^{n_l}$$
подалгебры $U_q\mathfrak{h}$. Из определения алгебры $U_q\mathfrak{g}$
следуют равенства
 $$K_\lambda K_j K_\lambda^{-1}=K_j,\quad
 K_\lambda E_jK_\lambda^{-1}=q_j^{\lambda(H_j)}E_j,\quad
 K_\lambda F_jK_\lambda^{-1}=q_j^{-\lambda(H_j)}F_j.$$

Покажем, что автоморфизм $S^2$ алгебры $U_q\mathfrak{g}$ является
внутренним. Как следует из (\ref{antipode_DJ}),
$$
S^2(K_j^{\pm1})=K_j^{\pm1},\quad S^2(E_j)=q_j^{-2}E_j,\quad
S^2(F_j)=q_j^2F_j,\quad j=1,2,\ldots,l.
$$
Значит,
\begin{equation}\label{inner_S2}
S^2(\xi)=K_{-2\rho} \, \xi K_{-2\rho}^{-1},\qquad \xi\in U_q\mathfrak{g},
\end{equation}
где
 $$\rho\stackrel{\rm def}{=}\frac{1}{2}\sum\limits_{\beta\in\Phi^+}\beta=
 \sum\limits_{i=1}^l\overline{\omega}_i$$
 (последнее равенство
хорошо известно, см. \cite[стр. 210]{Bou4-6}).

Рассмотрим весовой $U_q\mathfrak{g}$-модуль $V$ и отвечающее ему
представление $\pi$. Предложение \ref{tr_g} и равенство \eqref{inner_S2}
позволяют заключить, что линейный функционал
 \begin{equation}\label{q-tr-def}
 \operatorname{tr}_q A
\stackrel{\operatorname{def}}{=}
 \operatorname{tr}(A\pi(K_{-2\rho}))
 \end{equation}
на $U_q\mathfrak{g}$-модульной алгебре $\operatorname{End}_f(V)$ линейных
операторов конечного ранга в $V$ является инвариантным интегралом.

 Линейный функционал $\operatorname{tr}_q$
называют $q$-следом. Требование \eqref{group_like}, очевидно, выполняется:
$\Delta(K_{-2\rho})=K_{-2\rho} \otimes K_{-2\rho}$,\ \
$\varepsilon(K_{-2\rho})=1$.


\subsubsection{Мономиальные базисы.}\label{PBW}

Рассмотрим алгебру Хопфа $U_q\mathfrak{g}$ и подалгебры Хопфа
$U_q\mathfrak{h}$, $U_q\mathfrak{n}^+$, $U_q\mathfrak{n}^-$, порожденные
множествами $\{K_i^{\pm1}\}_{i=1,2,\ldots,l}$, $\{E_i\}_{i=1,2,\ldots,l}$,
$\{F_i\}_{i=1,2,\ldots,l}$, соответственно.

Построим мономиальные базисы векторных пространств $U_q\mathfrak{n}^\pm$.
Это приведет к базисам $U_q\mathfrak{g}$, поскольку мономы $\{K_1^{j_1}
K_2^{j_2} \cdots K_l^{j_l}\}_{j_1,j_2,\ldots,j_l\in\mathbb{Z}}$ образуют
базис $U_q\mathfrak{h}$ и линейные отображения
$$
U_q\mathfrak{n}^-\otimes U_q\mathfrak{h}\otimes U_q\mathfrak{n}^+\to
U_q\mathfrak{g},\quad u^-\otimes u^0\otimes u^+\mapsto u^-u^0u^+,
$$
$$
U_q\mathfrak{n}^+\otimes U_q\mathfrak{h}\otimes U_q\mathfrak{n}^-\to
U_q\mathfrak{g},\quad u^+\otimes u^0\otimes u^-\mapsto u^+u^0u^-
$$
являются изоморфизмами векторных пространств \cite[стр. 66]{Jant}.

Как следует из теоремы Пуанкаре-Биркгофа-Витта \cite[стр. 28]{Bou1-3}, в
классическом случае $q=1$ мономиальные базисы векторных пространств
$U\mathfrak{n}^\pm$ можно получить с помощью базисов в $\mathfrak{n}^\pm$,
например, с помощью базисов корневых векторов. Получим $q$-аналоги этих
базисов.

\medskip Прежде всего, введем в рассмотрение группу Вейля.
Наделим $\mathfrak{h}^*_{\mathbb{R}}$ скалярным произведением
\begin{equation}\label{inv_inner_prod}
(\alpha_i,\alpha_j)=d_ia_{ij},\qquad i,j=1,2,\ldots,l.
\end{equation}
Как следует из \eqref{alpha_definition},
\begin{equation}\label{inner_omega} (\overline{\omega}_i, \alpha_j)\,=\,
d_i\delta_{ij},\qquad i,j=1,2,\ldots,l.
\end{equation}
Введенное скалярное произведение доставляет изоморфизм
\begin{equation}\label{h_hstar}
\mathfrak{h}_\mathbb{R} \stackrel{\approx}{\to}
\mathfrak{h}_\mathbb{R}^*,\qquad
 H_i \mapsto \frac{2\alpha_i}{(\alpha_i,\alpha_i)},
\end{equation}
позволяющий отождествить $\mathfrak{h}_\mathbb{R}$ c
$\mathfrak{h}_\mathbb{R}^*$.

 Линейные операторы
$s_i:\lambda\mapsto\lambda-\lambda(H_i)\alpha_i$ в полученном евклидовом
пространстве являются отражениями относительно гиперплоскостей,
ортогональных простым корням $\alpha_i$. Действительно, из равенства
\hbox{$\alpha_i(H_i)=2$} следует, что $s_i^2=\operatorname{id}$, а из
равенств
$$
 (s_i(\alpha_j),\alpha_k)\,=\,d_ja_{jk}-a_{ij}d_ia_{ik}\,=\,
 d_ja_{jk}-a_{ik}d_ja_{ji}\,=\,(\alpha_j,s_i(\alpha_k))
$$
вытекает ортогональность операторов $s_i$.

Группой Вейля \IND{группа Вейля} называют группу $W$ обратимых линейных
операторов в $\mathfrak{h}^*_{\mathbb{R}}$, порожденную простыми
отражениями $s_1,s_2,\ldots,s_l$.

\begin{remark}\label{abstract_root_system}
Нетрудно показать, что конечное подмножество $\Phi\subset E$ является
(абстрактной) системой корней, \IND{система корней (абстрактная)} то есть
порождает $\mathfrak{h}^*_{\mathbb{R}}$ и
\begin{itemize}
\item[1.]\quad{$s_\alpha(\beta)\in\Phi$ при всех $\alpha,\beta\in\Phi$;}

\item[2.]\quad{$s_\alpha(\beta)-\beta\in\mathbb{Z}\alpha$ при всех
$\alpha,\beta\in\Phi$,} где
$s_\alpha(\lambda)=\lambda-2\frac{(\lambda,\alpha)\alpha}{(\alpha,\alpha)}$.
\end{itemize}
\end{remark}

Определяющий список соотношений между образующими $s_i$ группы $W$ имеет
вид \cite[стр. 116, 184]{Bou4-6}
$$s_i^2=1,\quad (s_is_j)^{m_{ij}}=1,\quad i\ne j,$$
где
$$
m_{ij}=
\begin{cases}
2, & \text{если}\quad a_{ij}a_{ji}=0
\\ 3, & \text{если}\quad a_{ij}a_{ji}=1
\\ 4, & \text{если}\quad a_{ij}a_{ji}=2
\\ 6, & \text{если}\quad a_{ij}a_{ji}=3.
\end{cases}
$$
Рассмотрим элемент $w_0$ максимальной длины группы $W$. Выберем его
приведенное разложение
\begin{equation}\label{reduced}
w_0=s_{i_1}s_{i_2}s_{i_3}\ldots s_{i_M},
\end{equation}
то есть разложение, для которого $M$ равно длине $w_0$. Ему отвечает
отношение линейного порядка на множестве положительных корней $\Phi^+$
\begin{equation}\label{positive_roots}
\beta_1=\alpha_{i_1},\quad\beta_2=s_{i_1}(\alpha_{i_2}),\quad\ldots,\quad
\beta_M=s_{i_1}s_{i_2}\ldots s_{i_{M-1}}(\alpha_{i_M}),
\end{equation}
см. \cite[стр. 196 -- 197]{Bou4-6}. Чтобы получить $q$-аналоги корневых
векторов $E_{\beta_1},E_{\beta_2},\ldots,E_{\beta_M}$, отвечающих этим
корням, воспользуемся введенными Люстигом в \cite{Lust1} автоморфизмами
$T_1,T_2,\ldots,T_l$ алгебры $U_q\mathfrak{g}$. \IND{автоморфизмы Люстига
алгебры $U_q\mathfrak{g}$} Их действие на образующие описывается следующими
равенствами:
\begin{equation}\label{reflection}
T_i(K_j)=K_jK_i^{-a_{ij}},
\end{equation}
$$T_i(E_i)=-F_iK_i,\qquad T_i(F_i)=-K_i^{-1}E_i,$$
$$
T_i(E_j)=\sum\limits_{r+s=-a_{ij}}\frac{(-1)^sq_i^{-r}}{[s]_{q_i}![r]_{q_i}!}
E_i^sE_jE_i^r,\quad i\ne j,
$$
$$
T_i(F_j)=\sum\limits_{r+s=-a_{ij}}\frac{(-1)^s q_i^r}{[s]_{q_i}![r]_{q_i}!}
F_i^rF_jF_i^s,\quad i\ne j.
$$
Рассмотрим присоединенное представление алгебры Хопфа $U_q\mathfrak{g}$. Из
\eqref{reflection} следует, что автоморфизмы $T_i$ переставляют весовые
подпространства: $T_i(U_q\mathfrak{g})_{\boldsymbol{\lambda}}=
(U_q\mathfrak{g})_{s_i(\boldsymbol{\lambda})},$ где
$$
(U_q\mathfrak{g})_{\boldsymbol{\lambda}}=\left\{\xi\in
U_q\mathfrak{g}\left|\:K_i\xi
K_i^{-1}=q_i^{\lambda_i}\xi\right.\right\},\qquad
\boldsymbol{\lambda}=(\lambda_1,\lambda_2,\ldots,\lambda_l)\in\mathbb{Z}^l.
$$
Как показал Люстиг, эти автоморфизмы удовлетворяют соотношениям кос
$$
\underbrace{T_iT_jT_i\cdots}_{m_{ij}}=\underbrace{T_jT_iT_j\cdots}_{m_{ij}},
\qquad i\ne j
$$
(число сомножителей в каждой из частей равенства равно $m_{ij}$).

Определим ''$q$-аналоги корневых векторов'' равенствами
$$
E_{\beta_k}=
\begin{cases}
T_{i_1}T_{i_2}T_{i_3}\ldots T_{i_{k-1}}(E_{i_k}), & k\ge 2
\\ E_{i_1}, & k=1
\end{cases},
$$
$$
F_{\beta_k}=
\begin{cases}
T_{i_1}T_{i_2}T_{i_3}\ldots T_{i_{k-1}}(F_{i_k}), & k\ge 2
\\ F_{i_1}, & k=1
\end{cases}.
$$
Следующий результат также принадлежит Люстигу.

\begin{proposition}\label{PBW-basis}
Мономы
\begin{equation}\label{basisE}
E_{\beta_1}^{j_1}E_{\beta_2}^{j_2}\ldots E_{\beta_M}^{j_M},\qquad
j_1,j_2,\ldots,j_M\in\mathbb{Z}_+,
\end{equation}
образуют базис в $U_q\,\mathfrak{n}^+$, а мономы
\begin{equation}\label{basisF}
F_{\beta_M}^{j_M}F_{\beta_{M-1}}^{j_{M-1}}\ldots F_{\beta_1}^{j_1},\qquad
j_1,j_2,\ldots,j_M\in\mathbb{Z}_+,
\end{equation}
образуют базис в $U_q\,\mathfrak{n}^-$.
\end{proposition}

В заключение приведем коммутационные соотношения, полученные Левендорским и
Сойбельманом \cite[стр. 67-68]{DeConProc}.

\begin{proposition}\label{Levend}
\begin{enumerate}
                 \item При всех $i<j$
\begin{equation}\label{Levend_1}
E_{\beta_i}E_{\beta_j}-q^{(\beta_i,\beta_j)}E_{\beta_j}E_{\beta_i}=
\sum\limits_{\mathbf{m}\in\mathbb{Z}_+^M}C'_\mathbf{m}(q)\cdot E^\mathbf{m},
\end{equation}
\begin{equation}\label{Levend_2}
F_{\beta_i}F_{\beta_j}-q^{-(\beta_i,\beta_j)}F_{\beta_j}F_{\beta_i}=
\sum\limits_{\mathbf{m}\in\mathbb{Z}_+^M}C''_\mathbf{m}(q)\cdot
F^\mathbf{m},
\end{equation}
где $\mathbf{m}=(m_1,\ldots,m_M)$,
$$
E^\mathbf{m}=E_{\beta_1}^{m_1}E_{\beta_2}^{m_2}\ldots E_{\beta_M}^{m_M},
\qquad F^\mathbf{m}=F_{\beta_M}^{m_M}F_{\beta_{M-1}}^{m_{M-1}}\ldots
F_{\beta_1}^{m_1},
$$
и коэффициенты $C'_\mathbf{m}(q)$, $C''_\mathbf{m}(q)$ могут быть отличны от
нуля лишь при $m_1=m_2=\ldots=m_i=0$ и $m_j=m_{j+1}=\ldots=m_M=0$.
             \item $C'_\mathbf{m}(q), C''_\mathbf{m}(q)\in
\mathbb{Q}[q,q^{-1}]$, то есть $C'_\mathbf{m}(q)$ и $C''_\mathbf{m}(q)$
являются полиномами Лорана с рациональными коэффициентами.
\end{enumerate}
\end{proposition}


Используя предложение \ref{Levend}, нетрудно наделить алгебру
$U_q\mathfrak{g}$ такой фильтрацией, что в присоединенной градуированной
алгебре
$$K_iK_j=K_jK_i,\quad K_iK_{i-1}=K_{i-1}K_i=1,\quad i=1,2,\ldots,l,$$
$$E_{\beta_i}E_{\beta_j}=E_{\beta_j}E_{\beta_i},\quad i=1,2,\ldots,M,$$
$$
E_{\beta_i}E_{\beta_j}=q^{(\beta_i,\beta_j)}E_{\beta_j}E_{\beta_i},\quad
F_{\beta_i}F_{\beta_j}=q^{-(\beta_i,\beta_j)}F_{\beta_j}F_{\beta_i},\quad
i<j.
$$

\begin{corollary}\label{nodivisors}\cite{DeConKac}
Алгебра $U_q\mathfrak{g}$ является областью целостности.
\end{corollary}

Следующее утверждение вытекает из предложения \ref{PBW-basis} и из
квазиперестановочности образующих присоединенной градуированной алгебры.

\begin{corollary}\label{PBWbasis-new}
Мономы
$$
E_{\beta_M}^{j_M}E_{\beta_{M-1}}^{j_{M-1}}\ldots E_{\beta_1}^{j_1},\qquad
j_1,j_2,\ldots,j_M\in\mathbb{Z}_+,
$$
образуют базис в $U_q\mathfrak{n}^+$, а мономы
$$
F_{\beta_1}^{j_1}F_{\beta_2}^{j_2}\ldots F_{\beta_M}^{j_M},\qquad
j_1,j_2,\ldots,j_M\in\mathbb{Z}_+,
$$
образуют базис в $U_q\mathfrak{n}^-$.
\end{corollary}

\subsubsection{Универсальная $R$-матрица.}\label{RM} Приведем обобщения
результатов \itemiiiа \ref{RM_sl_2}, используя обозначения этого
\itemiiiа.

Рассмотрим подалгебры Хопфа $U_q\mathfrak{b}^+$, $U_q\mathfrak{b}^-$,
порожденные множествами $\{K_i^{\pm1},E_i\}_{i=1,2,\ldots,l}$,
$\{K_i^{\pm1},F_i\}_{i=1,2,\ldots,l}$, соответственно. Полную подкатегорию
весовых локально \hbox{$U_q\mathfrak{b}^+$-конечномерных}
\hbox{$U_q\mathfrak{g}$-модулей} обозначим $\mathcal{C}^+$, \IND{категории
$\mathcal{C}^+$ и $\mathcal{C}^-$} а полную подкатегорию весовых локально
\hbox{$U_q\mathfrak{b}^-$-конечномерных} \hbox{$U_q\mathfrak{g}$-модулей}--
$\mathcal{C}^-$.

\begin{proposition}\label{Pairing}(\cite[стр. 946]{Tanisaki},
\cite[стр. 114]{Jant}). Существует и единственна такая билинейная форма
$(\cdot, \cdot)$: $U_q\mathfrak{b}^+ \times U_q\mathfrak{b}^-\rightarrow
\mathbb{C}$, что,
\begin{itemize}
 \item[1.\ \ ]
 $(\xi'\xi'',\eta)=(\xi''\otimes \xi', \Delta(\eta))$\ \  и\ \
$(\xi,\eta'\eta'')=(\Delta(\xi),\eta'\otimes\eta'')$ при всех
\hbox{$\xi,\xi',\xi''\in U_q\mathfrak{b}^+$}, $\eta,\eta',\eta''\in
U_q\mathfrak{b}^-$,
 \item[2.\ \ ] $ (K_\mu,K_\nu)=q^{-(\mu,\nu)}$ при
всех $\mu, \nu \in Q$,
 \item[3.\ \ ] $ {(E_i,K_j^{\pm1})=(K_i^{\pm1},F_j)=0}$,
$ {(E_i,F_j)=\delta_{ij}\frac{1}{q_i^{-1}-q_i}} $ при всех
$i,j=1,2,\ldots,l$.
\end{itemize}
\end{proposition}

\bigskip

Первое условие означает, что билинейная форма $(\cdot,\cdot)$ осуществляет
спаривание алгебр Хопфа $(U_q\mathfrak{b}^-)^{cop}$ и $U_q\mathfrak{b}^+$.
Из предложения \ref{R-nondegenerate} следует его невырожденность.

Напомним стандартное обозначение для весовых подпространств
$$
(U_q\mathfrak{n}^{\pm})_{\mathbf{\lambda}}=\left\{\xi\in
U_q\mathfrak{n}^{\pm}
\;\left|\;\operatorname{ad}(K_i)\,\xi=q_i^{\lambda_i}\,\xi,\quad
i=1,2,\ldots,l\right.\right\},
$$
где, как обычно, $\mathbf{\lambda}=(\lambda_1,\lambda_2,\ldots,\lambda_l)$.
Очевидно,
$$
U_q\mathfrak{n}^+=\newoplus\limits_{\lambda\in
Q_+}(U_q\mathfrak{n}^+)_{\lambda},\quad
U_q\mathfrak{n}^-=\newoplus\limits_{\lambda\in
-Q_+}(U_q\mathfrak{n}^-)_{\lambda}.
$$
Из определения билинейной формы $(\cdot,\cdot)$ следует, что ее сужение на
$(U_q\mathfrak{n}^-)_{\mu}\times(U_q\mathfrak{n}^+)_{\lambda}$ равно нулю
при $\lambda+\mu\ne 0$. Следующее утверждение доказано в \cite{Jant}
сначала в предположении о трансцендентности $q$, а потом для всех чисел
$q$, не являющихся корнями из единицы, см. \cite[стр. 105]{Jant}.

\begin{proposition}\label{R-nondegenerate}
Сужение билинейной формы $(\cdot,\cdot)$ на
$(U_q\mathfrak{b}^-)_{-\lambda}\times(U_q\mathfrak{b}^+)_{\lambda}$
невырождено при всех $\lambda\in Q_+$.
\end{proposition}

Пусть $\{\xi_i\}_{i=1}^{\infty}$ -- базис весовых векторов в
$U_q\mathfrak{n}^+$, а $\{\eta_i\}_{i=1}^{\infty}$ -- биортогональный
базис весовых векторов в $U_q\mathfrak{n}^-$: $(\eta_i,\xi_i)=\delta_{ij}$.

Если $V'$, $V''$ -- весовые $U_q\mathfrak{g}$-модули и либо $V'\in
\mathcal{C}^+$, либо $V''\in \mathcal{C}^-$, то корректно определен
линейный оператор
$$
\Theta_{V',V''}:V'\otimes V''\rightarrow V'\otimes V'',\quad
\Theta_{V',V''}:v'\otimes v''\mapsto\sum\limits_{i=1}^{\infty}\xi_i
v'\otimes\eta_i v'' $$ (при всех $v'\otimes v''\in V'\otimes V''$ ряд
$\sum\limits_{i=1}^{\infty}\xi_iv'\otimes\eta_iv''$ обрывается, и его сумма
не зависит от выбора базиса весовых векторов $\{\xi_i\}_{i=1}^{\infty}$).

В этих предположениях рассмотрим линейный оператор в $V'\otimes V''$,
определенный равенством
$$
R_{V',V''}=\Theta_{V',V''}\cdot
q^{-\sum\limits_{i,j=1}^l c_{ij}d_id_jH_i\otimes H_j},
$$
где $\mathbf{c}=(c_{ij})_{i,j=1,2,\ldots,l}$ -- матрица, обратная к
${\mathbf{b}=(d_ia_{ij})_{i,j=1,2,\ldots,l}}$.

Для линейных операторов $\check{R}_{V'V''}=\sigma_{V'V''}\; R_{V',V''}$
доказываются утверждения, дословно повторяющие предложения
\ref{braiding1_sl_2}, \ref{braiding2_sl_2} и означающие, что
$\check{R}_{V'V''}$ является $q$-аналогом оператора перестановки тензорных
сомножителей. Приведем эти утверждения.

\begin{proposition}\label{braiding1}
  1. Рассмотрим такие весовые $U_q\mathfrak{g}$-модули $V'$, $V''$, что либо
$V'\in \mathcal{C}^+$, либо $V''\in \mathcal{C}^-$. Линейный оператор
$$\check{R}_{V'V''}: V'\otimes V''\rightarrow V''\otimes V'$$ обратим и
является морфизмом
$U_q\mathfrak{g}$-модулей. \\
  2. Если $f':V'\rightarrow W'$, $f'':V''\rightarrow W''$ -- морфизмы весовых
$U_q\mathfrak{g}$-модулей, и либо $V',W'\in \mathcal{C}^+$, либо $V'',W''\in
\mathcal{C}^-$, то
$$
(f''\otimes f') \cdot \check{R}_{V',V''}=\check{R}_{W',W''}\cdot(f'\otimes
f'').
$$
\end{proposition}

\begin{proposition}\label{braiding2}

  1. Если $V$, $V'$, $V''$ -- весовые $U_q\mathfrak{g}$-модули и
   либо $V',V''\in \mathcal{C}^+$, либо $V\in \mathcal{C}^-$, то
$$
\check{R}_{V'\otimes V'',V}=(\check{R}_{V',V}\otimes
\operatorname{id}_{V''})(\operatorname{id}_{V'}\otimes \check{R}_{V'',V}).
$$
  2. Если $V$, $V'$, $V''$ -- весовые $U_q\mathfrak{g}$-модули и либо
  $V\in \mathcal{C}^+$, либо $V',V''\in \mathcal{C}^-$, то
$$
\check{R}_{V,V'\otimes V''}=(\operatorname{id}_{V'}\otimes
\check{R}_{V,V''})(\check{R}_{V,V'}\otimes \operatorname{id}_{V''}).
$$
  3. Для любого весового $U_q\mathfrak{g}$-модуля $V$
$$
\check{R}_{V,\mathbb{C}}=\check{R}_{\mathbb{C},V}=\operatorname{id}_V.
$$
\end{proposition}

\bigskip

Перейдем к обобщению формулы \eqref{Rmatrix_sl_2} для универсальной
$R$-матрицы алгебры Хопфа $U_q\mathfrak{sl}_2$. Выберем приведенное
разложение элемента $w_0\in W$ максимальной длины группы Вейля.
Воспользуемся участвующими в формулировке предложения \ref{PBW-basis}\ \
$q$-аналогами корневых векторов $E_{\beta_i}$, $F_{\beta_i}$, отвечающими
этому приведенному разложению. Пусть
$$q_{\beta}=q^{\frac{(\beta,\beta)}2},\qquad\beta\in\Phi^+.$$
Рассмотрим произведение
\begin{equation}\label{Rmatrix}
R=\prod\limits_{\beta\in\Phi^+}^{\curvearrowleft}
\exp_{q_\beta^2}\left((q_\beta^{-1}-q_\beta)E_\beta\otimes F_\beta\right)
q^{-\sum\limits_{i,j=1}^l c_{ij}d_id_jH_i\otimes H_j},
\end{equation}
сомножители которого записаны в порядке убывания номеров корней $\beta$:
\medskip \begin{multline*}
 \exp_{q_{\beta_M}^2}\left((q_{\beta_M}^{-1}-q_{\beta_M})E_{\beta_M}\otimes
 F_{\beta_M}\right)\exp_{q_{\beta_{M-1}}^2}
 \left((q_{\beta_{M-1}}^{-1}-q_{\beta_{M-1}})E_{\beta_{M-1}}\otimes
F_{\beta_{M-1}}\right)\\ \cdots
 \exp_{q_{\beta_1}^2}
\left((q_{\beta_1}^{-1}-q_{\beta_1})E_{\beta_1}\otimes
F_{\beta_1}\right)q^{-\sum\limits_{i,j=1}^lc_{ij}d_id_jH_i\otimes H_j}.
\end{multline*}

\begin{proposition}
\label{MultFormula}\cite{LevSoib},\cite{KR}
 Пусть $V'$,
$V''$ -- весовые $U_q\mathfrak{g}$-модули и либо $V'\in \mathcal{C}^+$,
 либо $V''\in \mathcal{C}^-$. Линейный оператор в $V'\otimes V''$, определяемый
правой частью равенства \eqref{Rmatrix}, равен $R_{V',V''}$.
\end{proposition}

\medskip
 Важные дополнительные
свойства универсальной \hbox{$R$-матрицы} приведены в работах С.~Хорошкина и
В.~Толстого \cite{KhorTol1, KhorTol2}.

\begin{note}\label{t_0}
Из равенств
$$
\alpha_i(H_j)=a_{ij},\qquad (\alpha_i,\alpha_j)=d_i a_{ij}=d_j a_{ji}
$$
и из определения матрицы $\mathbf{c}=(c_{ij})_{i,j=1,2,\ldots,l}$ следует,
что элемент
$$
t_0=\sum\limits_{i,j=1}^l c_{ij}d_id_jH_i\otimes H_j
$$
является каноническим элементом в $\mathfrak{h}\otimes \mathfrak{h}$,
отвечающим скалярному произведению $(\cdot,\cdot)$ в $\mathfrak{h}^*$:
\begin{equation}\label{t_canon}
 (\lambda,\mu)=\lambda\otimes\mu(t_0),\qquad \lambda,\mu\in\mathfrak{h}^*.
\end{equation}
Следовательно, для любых весовых $U_q\mathfrak{g}$-модулей $V'$, $V''$ и
любых их весовых векторов $v'\in V_{\lambda}'$, $v''\in V_\mu''$ имеет место
равенство
\begin{equation}\label{fraction}
  q^{-t_0}(v'\otimes v'')=q^{-(\lambda,\mu)} v'\otimes v''.
\end{equation}
\end{note}
Если $\{I_j\}_{j=1,2,\ldots,l}$ -- ортогональный базис в $\mathfrak{h}
\cong \mathfrak{h}^*$, то
\begin{equation}\label{t0_new}
 t_0=\sum\limits_k \frac{I_k\otimes I_k}{(I_k,I_k)}.
\end{equation}

\bigskip

Рассуждая так же, как в \itemiiiе \ref{RM_sl_2}, можно перейти от
$\mathbb{C}$ и числового параметра $q$ к топологической алгебре
$\mathbb{C}[[h]]$ с $h$-адической топологией.

Топологическая алгебра Хопфа $U_h\mathfrak{g}$ определяется своими
образующими $\{H_i,X_i^\pm\}_{i=1,2,\ldots,l}$ и соотношениями:
$$
H_iH_j-H_jH_i=0,\quad H_iX_j^\pm-X_j^\pm H_i=\pm a_{ij}X_j^\pm,
$$
$$
X_i^+X_j^- - X_j^-X_i^+
=\delta_{ij}\frac{e^{d_ihH_i/2}-e^{-d_ihH_i/2}}{e^{d_ih/2}-e^{-d_ih/2}},
$$
$$
\sum\limits_{s=0}^{1-a_{ij}}(-1)^s\begin{bmatrix}1-a_{ij}\\
s\end{bmatrix}_{e^{d_ih/2}}X_i^\pm X_j^\pm X_i^\pm =0,\quad i\ne j,
$$
$$
\Delta(H_i)=H_i\otimes1+1\otimes H_i,\quad \Delta(X_i^\pm)=X_i^\pm\otimes
e^{d_ihH/4}+e^{-d_ihH/4}\otimes X_i^\pm,
$$
$$
\varepsilon(H_i)=\varepsilon(X_i^\pm)=0,\quad S(H_i)=-H_i,
$$
$$
S(X_i^\pm)=-e^{\pm d_ih/2}X_i^\pm,\quad i=1,2,\ldots,l.
$$

\begin{note}\label{Change}
Эта алгебра Хопфа над $\mathbb{C}[[h]]$ связана с алгеброй
$U_q\mathfrak{g}$ над полем $\mathbb{C}$ следующими ``заменами образующих и
параметра деформации'':
\begin{equation}\label{change1}
  q=e^{-h/2},\quad K_i^{\pm1}=e^{\mp d_ihH_i/2},
\end{equation}
\begin{equation}\label{change2}
  E_i=X_i^+ e^{-d_ihH_i/4},\quad F_i=e^{d_ihH_i/4}X_i^-.
\end{equation}
\end{note}

\bigskip
Так же, как в случае $U_q\mathfrak{g}$, определяются спаривание
$(\cdot,\cdot)$ и формальный ряд $R$ по переменной $h$ с коэффициентами из
$U\mathfrak{g}\otimes U\mathfrak{g}$, называемый универсальной
$R$-матрицей. Применяя к коэффициентам этого ряда стандартные вложения
$
i_{1,2},i_{1,3},i_{2,3}: U\mathfrak{g}\otimes U\mathfrak{g}\rightarrow
U\mathfrak{g}\otimes U\mathfrak{g}\otimes U\mathfrak{g},
$
получаем формальные ряды $R^{1,2}$, $R^{1,3}$, $R^{2,3}$. Аналоги
предложений \ref{braiding1}, \ref{braiding2} вытекают из следующих основных
свойств универсальной $R$-матрицы:
\begin{equation}\label{R_property1}
\Delta^{\mathrm{op}}(\xi)=R\cdot\Delta(\xi)\cdot R^{-1},\quad \xi\in
U_q\mathfrak{g},
\end{equation}
\begin{equation}\label{R_property2}
 (\Delta\otimes \operatorname{id})(R)=R^{1,3}\cdot R^{2,3},
\end{equation}
\begin{equation}\label{R_property3}
 (\operatorname{id}\otimes\Delta)(R)=R^{1,3}\cdot R^{1,2},
\end{equation}
где $\Delta^{\mathrm{op}}$ -- противоположное коумножение, ср. с
\eqref{R_property1_sl_2}-\eqref{R_property3_sl_2}. Произведение
(\ref{Rmatrix}) называют мультипликативной формулой для универсальной
$R$-матрицы \cite{LevSoib, KR}.

Так же, как как в
\itemiiiе \ref{RM_sl_2}, вводится антиавтоморфизм $\tau$ алгебры
$U_h\mathfrak{g}$, определяемый равенствами
 $$\tau(H_i)\,=\,H_i,\qquad \tau(E_i)\,=\, e^{-\frac{hd_i H_i}{2}}F_i,\qquad
 \tau(F_i)\,=\,  E_ie^{\frac{hd_i H_i}{2}},$$
 и доказываются равенства, обобщающие
\eqref{S_theta_sl_2}:
\begin{equation}\label{S_theta}
(S\otimes S)R=R,\qquad(\tau\otimes\tau)R=R^{21},
\end{equation}
где $R^{21}$ отличается от $R$ перестановкой тензорных сомножителей. См.
статью Дринфельда \cite[стр. 34, 37]{Drinf2}.

Из \cite[предложение 4.2]{Drinf2} следует, что при замене коумножения на
противоположное универсальная $R$-матрица заменяется на обратную.

\subsubsection{Модули Верма.}\label{HChV}

Пусть $\mathbf{\lambda}=\sum\limits_{j=1}^l\lambda_j\,\overline{\omega}_j
\in \mathfrak{h}_\mathbb{R}^*$ и $\chi_\lambda$ -- одномерное представление
алгебры $U_q\mathfrak{b}^+$ в векторном пространстве $\mathbb{C}$:
\begin{equation}\label{chi_lambda_new}
\chi_\lambda:K_i^{\pm 1}\mapsto
q_i^{\pm\lambda_i},\qquad\chi_\lambda:E_i\mapsto 0,\qquad i=1,2,\ldots,l.
\end{equation}

Рассмотрим одномерный $U_q\mathfrak{b}^+$-модуль $\mathbb{C}_\lambda$,
отвечающий этому представлению, и определим модуль Верма \IND{модуль !
Верма} $M(\mathbf{\lambda})$ над алгеброй $U_q\mathfrak{g}$ так же, как в
классическом случае $q=1$:\ \
$M(\lambda)=U_q\mathfrak{g}\otimes_{U_q\mathfrak{b}^+}\mathbb{C}_\lambda.$

\begin{remark}\label{complex_lambda} Требование вещественности
числовых отметок $\lambda_j$ веса $\lambda$ можно было бы отбросить, но в
интересующих нас задачах оно всегда выполняется.
\end{remark}

Пусть $i_\lambda$ -- естественное вложение $U_q\mathfrak{b}^+$-модулей:
$$
i_\lambda:\mathbb{C}_\lambda\to M(\lambda),\qquad i_\lambda:z\mapsto
1\otimes z.
$$
Для любого $U_q\mathfrak{g}$-модуля $M$ и любого морфизма
$U_q\mathfrak{b}^+$-модулей $f:\mathbb{C}_\lambda\to M$ существует и
единствен морфизм $U_q\mathfrak{g}$-модулей $\widetilde{f}:M(\lambda)\to
M$, для которого $f=\widetilde{f}\cdot i_\lambda$. Это свойство модуля
Верма $M(\lambda)$ называют его универсальностью.

Пусть $v(\lambda)=1\otimes 1$. Модуль $M(\lambda)$ можно задать с помощью
образующей $v(\lambda)$ и определяющих соотношений:
$$
E_iv(\lambda)=0,\qquad K_i^{\pm
1}v(\lambda)=q_i^{\pm\lambda_i}v(\lambda),\qquad i=1,2,\ldots,l.
$$

Весовые векторы
\begin{equation}\label{PBV-Verma}
v_J(\lambda)=F_{\beta_M}^{j_M}F_{\beta_{M-1}}^{j_{M-1}}\ldots
F_{\beta_1}^{j_1}v(\lambda),\qquad J=(j_1,j_2,\ldots,j_M)\in\mathbb{Z}_+^M,
\end{equation}
образуют базис векторного пространства $M(\lambda)$ в силу предложения
\ref{PBW-basis}. Значит, $M(\lambda)$ -- весовой $U_q\mathfrak{g}$-модуль с
множеством весов $\lambda-Q_+$ и такими же размерностями весовых
подпространств, как в классическом случае $q=1$. Именно,
$$\dim M(\lambda)_{\lambda-\mu}=
\operatorname{card}\left\{(j_1,j_2,\ldots,j_M)\in\mathbb{Z}_+^M\left|
\:\sum\limits_{i=1}^Mj_i\beta_i=\mu\right.\right\}.
$$
Напомним, что множество весов является частично упорядоченным. Как следует
из определений, $M(\lambda)$ является $U_q\mathfrak{g}$-модулем со старшим
весом $\lambda$ и одномерным старшим весовым подпространством
$\mathbb{C}v(\lambda)$.

Наделим $U_q\mathfrak{g}$ инволюцией
\begin{equation}\label{aster1}
(K_j^{\pm 1})^\star=K_j^{\pm 1},\quad E_j^\star=K_jF_j,\quad
F_j^\star=E_jK_j^{-1},\quad j=1,2,\ldots,l,
\end{equation}
обобщающей \eqref{U_q_su_2}. Полученная $*$-алгебра Хопфа является
$q$-аналогом $*$-алгебры Хопфа
$U_q\mathfrak{g}_0\otimes_{\mathbb{R}}\mathbb{C}$, где $\mathfrak{g}_0$ --
компактная вещественная форма простой комплексной алгебры Ли $\mathfrak{g}$
\cite[стр. 24]{Bou9}.

Напомним что мы не включили требование неотрицательности в определение
эрмитовой формы и что условие $U_q\mathfrak{g}$-инвариантности эрмитовой
формы в $M(\lambda)$ имеют вид
$$
(\xi v',v'')=(v',\xi^\star v''),\qquad v',v''\in M(\lambda),\quad\xi\in
U_q\mathfrak{g},
$$
см. \itemiiiы \ref{IntStar} и \ref{Hodges}.

В $M(\lambda)$ существует и единственна инвариантная эрмитова форма
$(\cdot,\cdot)$, для которой $(v(\lambda),v(\lambda))=1$, как легко
показать, используя изоморфизм векторных пространств
$U_q\mathfrak{n}^-\otimes U_q\mathfrak{h}\otimes U_q\mathfrak{n}^+\approx
U_q\mathfrak{g}$, описанный в \itemiiiе\ \ref{PBW}. Этот изоморфизм
позволяет ввести в рассмотрение линейный оператор
$$
{\bf h}:\xi\eta\zeta\mapsto\varepsilon(\xi)\eta\varepsilon(\zeta),\qquad
\xi\in U_q\mathfrak{n}^-,\;\eta\in U_q\mathfrak{h},\;\zeta\in
U_q\mathfrak{n}^+
$$
из $U_q\mathfrak{g}$ в $U_q\mathfrak{h}$ и получить равенство $(\xi
v(\lambda),\eta v(\lambda))=\chi_\lambda({\bf h}(\eta^\star\xi))$, см.
\eqref{chi_lambda_new}.

 Весовые подпространства $M(\lambda)_\mu$ модуля Верма $M(\lambda)$ попарно
ортогональны. Каждый его собственный подмодуль является весовым, ортогонален
$\mathbb{C}v(\lambda)$ и, следовательно, содержится в ядре инвариантной
эрмитовой формы. Значит, это ядро $K(\lambda)$ является наибольшим
нетривиальным подмодулем модуля Верма $M(\lambda)$, фактормодуль
$L(\lambda)=M(\lambda)/K(\lambda)$ прост и наделен невырожденной эрмитовой
формой, для которой
\begin{equation*}\label{Hermit-L}
(\xi v',v'')=(v',\xi^\star v''),\qquad v',v''\in L(\lambda),\quad\xi\in
U_q\mathfrak{g}.
\end{equation*}

\begin{note}\label{uniqueness}
Нетрудно доказать, что $K(\lambda)$ является единственным подмодулем
конечной коразмерности $U_q\mathfrak{g}$-модуля $M(\lambda)$.
Действительно, пусть $\mathrm{dim}(M(\lambda)/\mathscr{K})<\infty$ для
некоторого подмодуля $\mathscr{K}\subset K(\lambda)$. Центральные характеры
$U_q\mathfrak{g}$-модулей $M(\lambda)/\mathscr{K}$ и $L(\lambda)$ равны.
Следовательно, $U_q\mathfrak{g}$-модуль $M(\lambda)/\mathscr{K}$ кратен
$L(\lambda)$, см. \itemiiiы \ref{U_ext}, \ref{loc_finite}. Кратность равна
единице, поскольку $\dim (M(\lambda)/\mathscr{K})_\lambda \leq \dim
M(\lambda)_\lambda=1$. Значит, $\mathscr{K}=K(\lambda)$.
\end{note}

Очевидно, $L(\lambda)=M(\lambda)$, если и только если матрица
$\mathbf{g}(\lambda,\mu)$ попарных скалярных произведений (матрица Грама)
векторов
$$
\left\{v_{J}(\lambda)\left|\: J=(j_1,j_2,\ldots,j_M)\in\mathbb{Z}_+^M \quad
\& \quad \sum\limits_{i=1}^M j_i\beta_i=\mu\right.\right\}
$$
невырождена при всех $\mu\in Q_+$. Имеется $q$-аналог известного результата
Шаповалова \cite{Shap}.

\begin{proposition}(\cite[стр. 129]{Jak3},\cite[стр. 484]{DeConKac})\label{DetShap}
Детерминант матрицы Грама $\mathbf{g}(\lambda,\mu)$ с точностью до
ненулевого числового множителя равен
$$
\chi_\lambda\left(\prod\limits_{\beta\in\Phi^+} \prod\limits_{m=1}^\infty
\left([m]_{q_\beta}\,\frac{K_\beta\,
q^{(\rho-\frac{m}2\beta,\beta)}-K_\beta^{-1}
q^{-(\rho-\frac{m}2\beta,\beta)}} {q^{\frac{(\beta,\beta)}2}
-q^{-\frac{(\beta,\beta)}2}}\right)^{P(\mu-m\beta)} \right),
$$
где, как и прежде, $q_{\beta}=q^{\frac{(\beta,\beta)}{2}}$ и
$[m]_{q_{\beta}}=\frac{q_{\beta}^m-q_{\beta}^{-m}}{q_{\beta}-q_{\beta}^{-1}}$.
\end{proposition}

\begin{remark}
Если $\mathbf{\lambda}=\sum\limits_j \lambda_j \overline{\omega}_j$,
 $\beta=\sum\limits_in_i\alpha_i$, то, как следует из (\ref{inner_omega}),
 $$ \chi_\lambda( K_\beta^{\pm 1})=q^{\pm \sum\limits_i d_i n_i \lambda_i}=
   q^{\pm \sum\limits_{i,j} (\alpha_i, \overline{\omega}_j) n_i \lambda_j}=
   q^{\pm (\lambda,\beta)}.
 $$
Это позволяет записать детерминант матрицы Грама следующим образом
$$
\operatorname{const} \cdot \prod\limits_{\beta\in\Phi^+}
\prod\limits_{m=1}^\infty \left([m]_{q_\beta}\, \frac{
  q^{(\lambda +\rho-\frac{m}2\beta,\beta)}-
q^{-(\lambda + \rho-\frac{m}2\beta,\beta)}} {q^{\frac{(\beta,\beta)}2}
-q^{-\frac{(\beta,\beta)}2}}\right)^{P(\mu-m\beta)},
$$
где $\operatorname{const} \neq 0$.
\end{remark}

\begin{corollary}\label{IrrVerma}
Модуль Верма $M(\lambda)$ прост, если и только если
$\frac{2(\lambda+\rho,\beta)}{(\beta,\beta)}\notin\mathbb{N}$ при всех
$\beta\in\Phi^+$.
\end{corollary}

\subsubsection{Центральные элементы алгебры \boldmath
$U_q^\mathrm{ext}\mathfrak{g}$.}\label{U_ext}

В \itemiiiе \ref{Weight} каждому $\lambda \in Q$ был сопоставлен элемент
\begin{equation}\label{K_ll}
K_\lambda=K_1^{a_1}K_2^{a_2}\cdots K_l^{a_l}, \qquad
\mathbf{\lambda}=\sum_{i=1}^l a_i\alpha_i.
\end{equation}
Очевидно,
\begin{equation}\label{relation_1}
K_\lambda K_\mu = K_{\lambda+\mu},\qquad K_0 = 1.
\end{equation}
Так как $(\alpha_i,\alpha_j)=d_ia_{ij}$, то
\begin{equation} \label{relation_2}
\\ K_\lambda E_j = q^{(\lambda,\alpha_j)}E_jK_\lambda, \qquad
K_\lambda F_j = q^{-(\lambda,\alpha_j)}F_jK_\lambda.
\end{equation}

Решения некоторых задач теории квантовых групп упрощаются, если использовать
элементы $K_\lambda$, отвечающие всем $\lambda\in P$. В этом случае числа
$a_i$ в \eqref{K_ll} являются рациональными, но не обязательно целыми, что
приводит к необходимости расширить алгебру Хопфа $U_q\mathfrak{g}$.

Следуя обзорам \cite[стр. 281]{KlSch}, \cite[стр. 73]{Letz}, введем в
рассмотрение алгебру $U_q^\mathrm{ext}\mathfrak{g}$, \IND{алгебра !
$U_q^\mathrm{ext}\mathfrak{g}$} определяемую своими образующими
$$E_i,F_i,\quad i=1,2,\ldots,l;\qquad K_\lambda,\quad\lambda\in P,$$
и соотношениями, получаемыми из соотношений для $U_q\mathfrak{g}$ заменой
\eqref{first_line} на \eqref{relation_1}, \eqref{relation_2}. Очевидно,
$U_q^\mathrm{ext}\mathfrak{g}$ является алгеброй Хопфа:
$$
\triangle(K_\lambda)=K_\lambda\otimes K_\lambda,\qquad
S(K_\lambda)=K_{-\lambda},\qquad\varepsilon(K_\lambda)=1,
$$
а $U_q\mathfrak{g}$-- ее подалгеброй Хопфа.

Так же, как в \itemiiiе \ref{Weight}, вводится понятие весового
$U_q^\mathrm{ext}\mathfrak{g}$-модуля $V$:
$$
V=\newoplus_\mu V_\mu,\qquad V_\mu=\left\{v\in V\left|\:K_\lambda
v=q^{(\lambda,\mu)}v, \; \; \lambda \in P\right.\right\}.
$$
Очевидно, категория весовых $U_q^\mathrm{ext}\mathfrak{g}$-модулей
канонически изоморфна категории весовых $U_q\mathfrak{g}$-модулей. Допуская
вольность речи, можно сказать, что алгебрам $U_q^\mathrm{ext}\mathfrak{g}$ и
$U_q\mathfrak{g}$ отвечает одна и та же квантовая группа.

Опишем центр $Z(U_q^\mathrm{ext}\mathfrak{g})$ алгебры
$U_q^\mathrm{ext}\mathfrak{g}$. Линейная оболочка
$U_q^\mathrm{ext}\mathfrak{h}$ множеств $\{K_\lambda\}_{\lambda\in P}$, и
линейная оболочка $U_q^\mathrm{even}\mathfrak{h}$ множеств
$\{K_\lambda\}_{\lambda\in 2P}$, являются коммутативными подалгебрами.
Линейное отображение
$$
U_q\mathfrak{n}^-\otimes U_q^\mathrm{ext}\mathfrak{h}\otimes
U_q\mathfrak{n}^+\to U_q^\mathrm{ext}\mathfrak{g},\qquad
\xi\otimes\eta\otimes\zeta\mapsto\xi\eta\zeta
$$
биективно (ср. с аналогичным результатом в \itemiiiе \ref{PBW}). Это
позволяет ввести в рассмотрение линейный оператор
$$
\gamma':Z(U_q^\mathrm{ext}\mathfrak{g})\to
U_q^\mathrm{ext}\mathfrak{h},\qquad
\gamma':\xi\eta\zeta\mapsto\varepsilon(\xi)\eta\varepsilon(\zeta),
$$
где $\xi\in U_q\mathfrak{n}^-$, $\eta\in U_q^\mathrm{ext}\mathfrak{h}$,
$\zeta\in U_q\mathfrak{n}^+$. Так же, как в классическом случае $q=1$
\cite[стр. 268]{Dix}, доказывается, что отображение $\gamma'$ является
гомоморфизмом алгебр.

Напомним теоретико-представленческий смысл этого гомоморфизма. Пусть $z \in
Z(U_q^\mathrm{ext}\mathfrak{g})$. Алгебра $U_q^\mathrm{ext}\mathfrak{h}$
естественно вкладывается в алгебру функций на $\mathfrak{h}_\mathbb{R}^*$:
$$K_\mu: \lambda \mapsto q^{(\mu,\lambda)},\qquad \mu \in P.$$ При этом вложении элементу
$\gamma'(z)$ отвечает функция, описывающая действие $z$ в модуле Верма
$M(\lambda)$. Введем изоморфизм "сдвига на полусумму положительных корней"
$$
\tau:U_q^\mathrm{ext}\mathfrak{h}\overset{\approx}{\to}
U_q^\mathrm{ext}\mathfrak{h},\qquad\tau:K_i^{\pm 1}\mapsto q_i^{\mp
1}K_i^{\pm 1}.
$$

Следующее утверждение является $q$-аналогом известной теоремы Хариш-Чандры
об описании центра алгебры $U\mathfrak{g}$ \cite[стр. 269]{Dix} и использует
естественное действие группы $W$ в $U_q^\mathrm{ext}\mathfrak{h}$:
$$wK_\lambda\,\stackrel{\rm def}{=}\,K_{w^{-1}\lambda},\qquad
 \lambda\in P,\;\;w\in W.$$

\begin{proposition}\cite[стр. 951]{Tanisaki} \label{center}
 Ограничение гомоморфизма\\
$\gamma\overset{\mathrm{def}}{=}\tau\gamma'$ на
$Z(U_q^\mathrm{ext}\mathfrak{g})$ инъективно, и его образ совпадает с
подалгеброй всех $W$-инвариантных элементов алгебры
$U_q^\mathrm{even}\mathfrak{h}$.
\end{proposition}

Опишем принадлежащий В. Дринфельду метод построения элементов центра алгебры
$U_q^\mathrm{ext}\mathfrak{g}$.

Из (\ref{t_canon}) следует, что изоморфизм, обратный к (\ref{h_hstar}),
 имеет следующий вид:
\begin{equation}\label{H_mu}
\mathfrak{h}^* \stackrel{\approx}{_\rightarrow} \mathfrak{h},\qquad
\mu\mapsto\operatorname{id}\otimes\mu(t_0)
\end{equation}
Пусть $H_\mu$-- его значение на элементе $\mu$. Используя обозначения
\itemiiiа \ref{Weight} и $\eqref{h_hstar}$, получаем
$$ H_{\alpha_i}=d_i H_i,\qquad
K_{\alpha_i}^{\pm 1}=K_i^{\pm 1}=q_i^{\pm H_i}=q^{\pm H_{\alpha_i}}.
$$
Значит,
\begin{equation}\label{H_omega_i}
K_{\overline{\omega}_i}^{\pm 1}=q^{\pm H_{\overline{\omega}_i}},\qquad
i=1,2,\ldots,l.
\end{equation}

\begin{remark}\label{id_pi}
Пусть $V$ -- конечномерный весовой $U_q^\mathrm{ext}\mathfrak{g}$-модуль и
$\pi_V$ -- отвечающее ему представление алгебры
$U_q^\mathrm{ext}\mathfrak{g}$. Как следует из \eqref{H_mu},
\eqref{H_omega_i}, элемент $\operatorname{id}\otimes\pi_V(q^{-t_0})$
алгебры $U_q^\mathrm{ext}\mathfrak{g}\otimes\operatorname{End}V$ корректно
определен. Отвечающий ему линейный оператор из $V$ в
$U_q^\mathrm{ext}\mathfrak{g}\otimes V$ отображает весовой вектор $v_\mu$
веса $\mu$ в $K_{-\mu}\otimes v_\mu$.
\end{remark}

Пусть $R$ -- универсальная R-матрица \eqref{Rmatrix}, и $R_{21}$ отличается
от $R$ перестановкой тензорных сомножителей:
$$
R_{21}=\prod_{\beta\in\Phi^+}^\curvearrowleft
\exp_{q_\beta^2}((q_\beta^{-1}-q_\beta)F_\beta\otimes E_\beta)q^{-t_0}.
$$

С каждым конечномерным весовым $U_q^\mathrm{ext}\mathfrak{g}$-модулем $V$
свяжем линейный функционал
$$
\operatorname{tr}_q:\operatorname{End}V\to\mathbb{C},\qquad
\operatorname{tr}_q:A\mapsto\operatorname{tr}(A\pi_V(K_{-2\rho})),
$$
где $\pi_V$ -- представление алгебры $U_q^\mathrm{ext}\mathfrak{g}$,
отвечающее $V$. Так же, как в \itemiiiе \ref{Weight}, доказывается, что этот
линейный функционал $U_q^\mathrm{ext}\mathfrak{g}$-инвариантен.

\begin{proposition}\label{C_V}
Элемент $C_V=(\operatorname{id}\otimes\operatorname{tr}_q\pi_V)(R_{21}R)$
 принадлежит центру алгебры $U_q^\mathrm{ext}\mathfrak{g}$.
\end{proposition}

{\bf Доказательство.} Нужный результат известен в контексте формальных рядов
переменной $h$. Именно, в \cite[стр. 32]{Drinf2} показано, что элемент
$\operatorname{id}\otimes l(R_{21}R)$ принадлежит центру, если линейный
функционал $l$ на алгебре Хопфа выбран так, что $l(\xi\eta)=l(\eta
S^2(\xi))$ для всех элементов $\xi$, $\eta$ этой алгебры. Остается положить
$l(\xi)=\operatorname{tr}_q\pi_V(\xi)$. Для перехода от числового параметра
$q\in(0,1)$ к формальному параметру $h$ в равенстве
$$C_V\cdot\xi=\xi\cdot C_V,\qquad\xi\in U_q^\mathrm{ext}\mathfrak{g},$$
достаточно разложить обе части этого равенства по базису
Пуанкаре-Биркгофа-Витта в $U_q^\mathrm{ext}\mathfrak{g}$ (см. \itemiii\
\ref{PBW}). После подстановки $q=e^{-h/2}$ каждый из коэффициентов
разложения становится аналитической функцией на полуоси $h>0$ и в
окрестности точки $h=0$ (см. предложение \ref{Levend}). Остается доказать,
что равны нулю все производные этих функций в точке $h=0$. Но это следует из
справедливости равенства $C_V\xi=\xi C_V$ в контексте формальных рядов.
$\square$

\medskip

Следующий результат доставляет описание образа элемента $C_V$ при
гомоморфизме $\gamma'$.

\begin{proposition}\label{center_char}
Пусть $M(\lambda)$ -- модуль Верма со старшим весом $\lambda$ и $V$ --
конечномерный весовой $U_q^\mathrm{ext}\mathfrak{g}$-модуль. Тогда
\begin{equation}\label{C_v_image}
C_V|_{M(\lambda)}=\sum_{\mu\in P}(\dim V_\mu)\, q^{-2(\mu,\lambda+\rho)},
\end{equation}
где $\rho$ -- полусумма положительных корней.
\end{proposition}

{\bf Доказательство.} Равенство (\ref{C_v_image}) в контексте формальных
рядов известно \cite[стр. 39]{Drinf2}. Переход от числового параметра
$q\in(0,1)$ к формальному параметру $h$ осуществляется так же, как при
доказательстве предыдущего предложения. \hfill $\square$

\begin{corollary}\label{image_} Для любого
 конечномерного весового $U_q^\mathrm{ext}\mathfrak{g}$-модуля $V$
\begin{equation}\label{image_new}
\gamma(C_V)=\sum\limits_{\mu\in P}(\dim V_\mu)K_{-2\mu}.
\end{equation}
\end{corollary}

{\bf Доказательство.} Правая часть \eqref{image_new} равна
$$
\sum_{\mu\in P}(\dim V_\mu)q^{-2(\mu\otimes(\lambda+\rho))t_0}=\sum_{\mu\in
P}(\dim V_\mu)q^{-2(\lambda+\rho)(H_\mu)}.\eqno\square
$$

\medskip

Рассмотрим конечномерный простой весовой $U_q\mathfrak{g}$-модуль
$L(\overline{\omega}_i)$ со старшим весом $\overline{\omega}_i$,
$i=1,2,\ldots,l$. Сохраним для отвечающего ему
$U_q^\mathrm{ext}\mathfrak{g}$-модуля то же обозначение.

\begin{corollary}\label{generators}
Множество $\{C_{L(\overline{\omega}_i)}\}_{i=1,2,\ldots,l}$ порождает центр
\\ $Z(U_q^\mathrm{ext}\mathfrak{g})$ алгебры $U_q^\mathrm{ext}\mathfrak{g}$.
\end{corollary}

{\bf Доказательство} проводится с помощью стандартных рассуждений \cite[стр.
125]{Jant} с использованием следствия \ref{image_} и предложения
\ref{center}. Именно, достаточно доказать, что следующие подалгебры алгебры
$U_q^\mathrm{ext}\mathfrak{h}$ совпадают:
\begin{description}
\item[]$A_1$.\ \ подалгебра, порожденная элементами $\sum\limits_{\mu}\dim
L(\overline{\omega}_i)K_{-2\mu}$, $i=1,2,\ldots,l$;

\item[]$A_2$.\ \ подалгебра, порожденная элементами $\sum\limits_{\mu}\dim
L(\lambda)_\mu K_{-2\mu}$, $\lambda\in P_+$;

\item[]$A_3$.\ \ подалгебра, порожденная элементами $\sum\limits_{\mu\in
W\lambda}K_{-2\mu}$, где $\lambda\in P_+$;

\item[]$A_4$.\ \ алгебра $W$-инвариантов линейной оболочки элементов
$K_{-2\lambda}$, $\lambda\in P_+$;

\item[]$A_5$.\ \ $\gamma(Z(U_q^\mathrm{ext}\mathfrak{g}))$.
\end{description}

Равенство $A_1\,=\,A_2$ имеет место в силу того, что кольцо конечномерных
весовых представлений $U_q^\mathrm{ext}\mathfrak{g}$ канонически изоморфно
кольцу конечномерных представлений алгебры $U\mathfrak{g}$ и,
следовательно, порождается фундаментальными представлениями \cite[стр.
175]{Bou7-8}.

Для доказательства равенства $A_2\,=\,A_3$ достаточно воспроизвести
рассуждения, приведенные в \cite[стр. 125]{Jant}, причем индукцию удобно
вести по $(\lambda,\lambda)$, где $\lambda \in P_+$, используя то, что для
любого $\lambda \in P_+$ и любого веса $\mu\in P_+$ представления
$L(\lambda)$, отличного от $\lambda$, имеет место неравенство
$(\mu,\mu)<(\lambda,\lambda)$ \cite[стр. 234]{KV}.

Равенство $A_3\,=\,A_4$ следует из \cite[стр. 232]{Bou4-6}, а равенство
$A_5\,=\,A_6$ вытекает из предложения \ref{center}. \hfill $\square$

\medskip

Нетрудно доказать, что элементы $C_{L(\overline{\omega}_i)}$,
$i=1,2,\cdots,l$, алгебраически независимы и, следовательно, центр
$Z(U_q^{\operatorname{ext}}\mathfrak{g})$ алгебры
$U_q^{\operatorname{ext}}\mathfrak{g}$ изоморфен алгебре полиномов $l$
переменных, как и в классическом случае $q=1$.

\begin{remark} \label{center_Uq}
 Центр $Z(U_q\mathfrak{g})$
 алгебры $U_q\mathfrak{g}$ устроен несколько сложнее и
  является линейной оболочкой элементов
 $C_{L(\lambda)}$, $\lambda \in P_+ \cap \frac{1}{2}Q$,
 см. \cite[стр. 109]{Jant}.
  Нетрудно доказать,  что для любых
  $\mu_1, \mu_2 \in P_+$ найдется элемент
  $z \in Z(U_q\mathfrak{g})$, разделяющий $U_q\mathfrak{g}$-модули
  $L(\mu_1),L(\mu_2)$, см. \cite[стр. 125-126]{Jant}.
\end{remark}

\subsubsection{Редукция к случаю $q=1$.}\label{reduction}

Опишем в самых общих чертах метод сведения задач об
$U_q\mathfrak{g}$-модулях к соответствующим задачам об
$U\mathfrak{g}$-модулях. Данный метод особенно эффективен при
дополнительном предположении о трансцендентности $q$ \cite[главы
5-7]{Jant}, а переход от этого частного случая к общему случаю требует
дополнительных усилий \cite[глава 8]{Jant}.


Равенство \eqref{fraction} показывает, что в формулах теории
$U_q\mathfrak{g}$-модулей могут возникать дробные степени $q$. Пусть $s$--
порядок группы $P/Q$. Как следует из \eqref{inner_omega},
$(\lambda,\mu)\in\frac1{s}\mathbb{Z}$ при всех $\lambda,\mu\in P$.
Рассмотрим поле $\mathbb{C}(t^{1/s})$ рациональных функций переменной
$t^{1/s}$ и алгебру Хопфа $U_t\mathfrak{g}$ \IND{алгебра ! Хопфа !
$U_t\mathfrak{g}$} над этим полем, определенную теми же образующими и
соотношениями, что $U_q\mathfrak{g}$, но с заменой числового параметра $q$
переменной $t$. Эта же замена позволяет получить элементы
$E_{\beta_j},\;F_{\beta_j},\;j=1,2,\ldots,M$, алгебры $U_t\mathfrak{g}$,
модули Верма $M(\lambda)$ над этой алгеброй, их простые фактормодули
$L(\lambda)$ и линейные операторы $\check{R}_{L(\lambda),L(\mu)}$ при
$\lambda,\mu\in P$. Последнее условие будет предполагаться выполненным во
всех утверждениях настоящего \itemiiiа.

Пусть $\mathscr{A}=\mathbb{Q}[t^{1/s},t^{-1/s}]$ -- алгебра полиномов
Лорана переменной $t^{1/s}$ с рациональными коэффициентами.

Основное наблюдение состоит в том, что многие свойства
$U_t\mathfrak{g}$-модулей могут быть сформулированы и доказаны с помощью
содержащихся в них выделенных $\mathscr{A}$-модулей. Действительно,
рассмотрим $\mathscr{A}$-подалгебру $U_{\mathscr{A}}$ в $U_q\mathfrak{g}$,
порожденную элементами $K_i^{\pm
1},E_i,F_i,{h_i=\frac{K_i-K_i^{-1}}{t_i-t_i^{-1}}}$, ${i=1,2,\ldots,l}$,
где $t_i=t^{d_i}$. Известно \cite{DeConKac}, что определяющий список
соотношений между этими образующими можно получить, заменив в стандартном
списке соотношений для $U_t\mathfrak{g}$ равенство
${E_iF_j-F_jE_i=\delta_{ij}\frac{K_i-K_i^{-1}}{t_i-t_i^{-1}}}$ на
${E_iF_j-F_jE_i=\delta_{ij}h_i}$ и добавив соотношение
$$(t_i-t_i^{-1})h_i=K_i-K_i^{-1},\qquad i=1,2,\ldots,l.$$
Разумеется, $\mathscr{A}$-подалгебра $U_{\mathscr{A}}$ наследует структуру
алгебры Хопфа:
$$
\Delta(h_i)=h_i\otimes K_i+K_i^{-1}\otimes h_i,\quad S(h_i)=-h_i,\quad
\varepsilon(h_i)=0,\quad i=1,2,\ldots,l.
$$
Векторы $F_{\beta_\mu}^{j_\mu}F_{\beta_{\mu-1}}^{j_{\mu-1}}\ldots
F_{1}^{j_1}v_\lambda$ порождают $\mathscr{A}$-модуль
$M_{\mathscr{A}}(\lambda)=U_{\mathscr{A}}\cdot v_{\lambda}$. Образ
$M_{\mathscr{A}}(\lambda)$ при каноническом эпимоморфизме $M(\lambda)\to
L(\lambda)$ обозначим $L_{\mathscr{A}}(\lambda)$. Очевидно,
$M_{\mathscr{A}}(\lambda)$ и $L_{\mathscr{A}}(\lambda)$ являются
$U_{\mathscr{A}}$-модулями.

Опишем два гомоморфизма алгебр, первый из которых позволяет переносить
свойства $U\mathfrak{g}$-модулей на $U_{\mathscr{A}}$-модули, а второй ---
свойства $U_{\mathscr{A}}$-модулей на $U_q\mathfrak{g}$-модули. Рассмотрим
гомоморфизм $\mathbb{Q}$-алгебр
$$\mathbb{Q}[t^{1/s},t^{-1/s}]\to\mathbb{C},\qquad t^{\pm 1/s}\mapsto 1.$$
Существует и единственно такое его продолжение до гомоморфизма
\hbox{$\mathbb{Q}$-алгебр} $U_{\mathscr{A}}\to U\mathfrak{g}$, что
$$
K_i^{\pm 1}\mapsto 1,\quad E_i\mapsto E_i,\quad F_i\mapsto F_i,\quad
h_i\mapsto H_i,\qquad i=1,2,\ldots,l.
$$
Этот гомоморфизм позволяет использовать классическую теорию
$U\mathfrak{g}$-модулей при изучении $U_{\mathscr{A}}$-модулей. С другой
стороны, рассмотрим гомоморфизм $\mathbb{Q}$-алгебр
$$
\mathbb{Q}[t^{1/s},t^{-1/s}]\to\mathbb{C},\qquad t^{\pm 1/s}\mapsto q^{\pm
1/s}.
$$
Существует и единственно такое его продолжение до гомоморфизма
\hbox{$\mathbb{Q}$-алгебр} $U_{\mathscr{A}}\to U_q\mathfrak{g}$, что
$$
K_i^{\pm 1}\mapsto K_i^{\pm 1},\quad E_i\mapsto E_i,\quad F_i\mapsto F_i,
\quad h_i\mapsto\frac{K_i-K_i^{-1}}{q_i-q_i^{-1}},\qquad i=1,2,\ldots,l.
$$
Этот гомоморфизм инъективен при трансцендентных $q$ и позволяет переносить
свойства $U_{\mathscr{A}}$-модулей на $U_q\mathfrak{g}$-модули.

Например, известно \cite[стр. 34]{Drinf2}, что
\begin{equation}\label{antianti}
  S\otimes S(R)=R.
\end{equation}
Введя обозначения $q=e^{-h/2}$, $\Theta=Rq^{\mathrm{t}_0}$, этому равенству
можно придать вид
\begin{equation}\label{antitheta}
  S\otimes S(\Theta) = q^{-\mathrm{t}_0} \Theta q^{\mathrm{t}_0}.
  \end{equation}
 Используя переходы от $U_q\mathfrak{g}$ к
$U_{\mathscr{A}}\mathfrak{g}$ и от $U_{\mathscr{A}}\mathfrak{g}$ к
$U_h\mathfrak{g}$, нетрудно получить с помощью (\ref{antitheta}) равенство
$$
R_{V_1,V_2}^{*}\; (l_1 \otimes l_2)\;=\;R_{V_1^*,V_2^*}\; (l_1 \otimes
l_2), \qquad l_1 \in V_1^*,\; l_2 \in V_2^*
$$
для широкого класса весовых $U_q\mathfrak{g}$-модулей $V_1, V_2$.

\subsubsection{Локально конечномерные $U_q\mathfrak{g}$-модули.}
\label{loc_finite}

Следующие ут\-верж\-дения хорошо известны \cite{Jant}.

\begin{proposition}\label{irreducible}
1) Простые весовые $U_q\mathfrak{g}$-модули $L(\lambda)$, $\lambda\in
P_{+}$, конечномерны и попарно неизоморфны.
\\ 2) Каждый простой весовой конечномерный $U_q\mathfrak{g}$-модуль
изоморфен одному из $U_q\mathfrak{g}$-модулей $L(\lambda)$, $\lambda\in
P_{+}$.
\end{proposition}

Сохраним обозначение $v(\lambda)$ для образа одноименной образующей модуля
Верма $M(\lambda)$ при каноническом эпиморфизме $M(\lambda)\to L(\lambda)$.

\begin{proposition}\label{listl}
Пусть
$\mathbf{\lambda}=\sum\limits_{j=1}^l\,\lambda_j\,\overline{\omega}_j\in
P_+$. Элемент $v(\lambda)$ порождает $U_q\mathfrak{g}$-модуль $L(\lambda)$,
и следующие соотношения являются определяющими:
\begin{equation}\label{complete_list}
E_i\;v(\lambda)=0,\quad K_i^{\pm 1}\;v(\lambda)=q_i^{\pm\lambda_i}\;\;
v(\lambda),\quad F_i^{\lambda_i+1}\;v(\lambda)=0,
\end{equation}
где $i=1,2,\ldots,l$.
\end{proposition}

\medskip

\begin{remark} \label{second-def} Соотношения \eqref{complete_list}
доставляют второе, эквивалентное определение $U_q \mathfrak{g}$-модуля
$L(\lambda)$ при $\lambda \in P_+$.
\end{remark}

\begin{proposition}(ср. \cite[ стр. 76-77]{Jant})\label{separate}
 Для любого ненулевого вектора $\xi \in U_q\mathfrak{g}$
найдется такой элемент $\lambda \in P_{+}\bigcap Q$, что $\xi\,L(\lambda)
\neq 0$.
\end{proposition}

\begin{proposition}\label{qWeyl}
(\cite{Jant}, стр. 81) Множество весов $L(\lambda)$ и их кратности не
изменяются при переходе от $U\mathfrak{g}$ к $U_q\mathfrak{g}$, то есть от
классического случая к квантовому.
\end{proposition}

\begin{proposition}\label{semisimple}
(\cite{Jant}, стр. 82) Каждый весовой конечномерный $U_q\mathfrak{g}$-модуль
полупрост, то есть отвечающее ему представление алгебры $U_q\mathfrak{g}$
вполне приводимо.
\end{proposition}

Из предложений \ref{semisimple}, \ref{irreducible}, \ref{qWeyl} получаем

\begin{corollary}\label{Clebsh_Gordan}
При всех $\lambda, \mu \in P_{+}$
$$
L(\lambda) \otimes L(\mu) \approx\sum \limits_{\nu \in P_{+}}
c^{\nu}_{\lambda\mu} L(\nu),
$$
 и кратности $c^{\nu}_{\lambda\mu}$ вхождения $L(\lambda)$ в $L(\lambda) \otimes L(\mu)$
 равны их значениям в классическом случае $q=1$.
\end{corollary}

Пусть $P^{\nu}_{\lambda\mu}$ -- проектор в $L(\lambda)\otimes L(\mu)$ на
изотипическую компоненту, кратную $L(\nu)$, параллельно сумме остальных
изотипических компонент.

\begin{proposition}\label{spectr_R}
(\cite{Drinf2}, стр. 39) При всех $\lambda,\mu \in P_{+}$
\begin{equation}\label{Resh_Dr}
  \check{R}_{L(\mu),L(\lambda)} \check{R}_{L(\lambda),L(\mu)}=
  \newoplus\limits_{\nu \in P_{+}} q^{(\lambda,\lambda+2\rho)+(\mu,\mu+2\rho)-
  (\nu,\nu+2\rho)} P^{\nu}_{\lambda\mu}.
\end{equation}
\end{proposition}

\medskip

Пусть $\mathcal{C}$ -- полная подкатегория всех весовых локально
конечномерных $U_q\mathfrak{g}$-модулей. \IND{категория ! $\mathcal{C}$}

\begin{proposition}\label{semisimple1}
Каждый $U_q\mathfrak{g}$-модуль $V\in \mathcal{C}$ единственным образом
разлагается в прямую сумму $V=\newoplus\limits_{\lambda \in P_+}
V(\lambda)$ своих подмодулей $V(\lambda)$, кратных $L(\lambda)$.
\end{proposition}

{\bf Доказательство.} Пусть $\lambda\in P_{+}$. Рассмотрим векторное
пространство
$$
X(\lambda)=\left\{v\in V\left|\:E_jv=0,\quad K_j^{\pm
1}v=q_j^{\pm\lambda_j}v,\;\;j=1,2,\ldots,l\right.\right\}
$$
и порожденный им $U_q\mathfrak{g}$-модуль $V(\lambda)=U_q\mathfrak{g}\cdot
X(\lambda)$. Достаточно доказать, что это наибольший подмодуль
$U_q\mathfrak{g}$-модуля $V$, кратный $L(\lambda)$, и что
$V=\bigoplus\limits_{\lambda\in P_{+}}V(\lambda)$. Для любого ненулевого
вектора $v_\lambda\in X(\lambda)$ существует и единственно вложение
$U_q\mathfrak{g}$-модулей
$$L(\lambda)\hookrightarrow V,\qquad v(\lambda)\mapsto v_\lambda.
$$
(единственность очевидна, а существование легко доказать, используя
универсальность модуля Верма $M(\lambda)$ и замечание \ref{uniqueness}).

Значит, каждому базису векторного пространства $X(\lambda)$ отвечает
разложение $V(\lambda)$ в прямую сумму модулей, кратных $L(\lambda)$.
Разложение $V=\newoplus\limits_{\lambda\in P_{+}} V(\lambda)$ вытекает из
приведенных выше свойств конечномерных $U_q\mathfrak{g}$-модулей. \hfill
$\square$

\begin{corollary}\label{semisimple2}
Каждый $U_q\mathfrak{g}$-модуль категории $\mathcal{C}$ полупрост.
\end{corollary}

\bigskip

Категория $\mathcal{C}$ замкнута относительно перехода к подмодулям,
фактормодулям и тензорным произведениям.

Рассмотрим функтор $F$ из категории всех $U_q\mathfrak{g}$-модулей в
категорию $\mathcal{C}$, сопоставляющий каждому $U_q\mathfrak{g}$-модулю
$V$ его наибольший подмодуль категории $\mathcal{C}$. Действие на морфизмы
определяется очевидным образом. Двойственный объект $V^*$ в категории
$\mathcal{C}$ -- это результат применения функтора $F$ к
$U_q\mathfrak{g}$-модулю всех линейных функционалов на $V$.

Алгебра $U_q\mathfrak{g}$ с присоединенным действием является
$U_q\mathfrak{g}$-модульной алгеброй, и $F(U_q\mathfrak{g})$-- ее
собственная $U_q\mathfrak{g}$-модульная подалгебра.

\begin{proposition}\label{Noether}
(\cite{Jo}, глава 7) Алгебры $U_q\mathfrak{g}$ и $F(U_q\mathfrak{g})$
являются как левыми, так и правыми нетеровыми алгебрами.
\end{proposition}

 \bigskip Пусть $\lambda \in P_+$.
 Результаты \itemiiiов \ref{PBW}, \ref{HChV}
  и настоящего
 \itemiiiа позволяют  естественным образом ввести
непрерывное над $(0,1]$ и аналитическое в $(0,1)$ векторное
 расслоение $\mathcal{E}_\lambda$ со слоями, изоморфными
$L(\lambda)$.

 Опишем вкратце построение расслоения $\mathcal{E}_\lambda$. В
\itemiiiе \ref{HChV} была введена эрмитова форма
 $(\cdot,\cdot)$ в $L(\lambda)$ и базисы  векторного пространства
$M(\lambda)$,
 отвечающие приведенным разложениям элемента $w_0\in W$.

Для любой точки $q_0\in (0,1]$ выберем подмножество
 $\{v_j\}_{j=1,2,\ldots,\dim L(\lambda)}$ элементов одного из выделенных
 базисов так, чтобы матрица Грама
 $((v_i,v_j))_{i,j=1,2,\ldots,\dim L(\lambda)}$ была невырожденной.
 По непрерывности она невырождена при всех $q$, достаточно близких
 к $q_0$, и мы получаем тривиальное расслоение с требуемыми свойствами
над окрестностью точки $q_0$, поскольку $\dim L(\lambda)$ не зависит от
$q\in (0,1]$.

 Элементы
матрицы Грама являются непрерывными в $(0,1]$ и аналитическими в $(0,1)$
функциями. Следовательно, во-первых, матрицы операторов $E_i, F_i, H_i$
непрерывны в $(0,1]$ и аналитичны в $(0,1)$, во-вторых, такими же
свойствами обладают матрицы перехода, заданные на пересечениях
окрестностей. Наконец, полученное векторное расслоение на полуинтервале
$(0,1]$ в существенном не зависит от сделанных при его построении выборов.

Используя $\mathcal{E}_\lambda$, легко доказать положительную
определенность инвариантной эрмитовой формы в $L(\lambda)$, введенной в
\itemiiiе \ref{HChV}. Действительно, эта эрмитова форма положительно
определена в пределе $q \rightarrow 1$, аналитична в интервале $q \in
(0,1)$ и невырождена при всех $q \in (0,1)$, как следует из простоты
$U_q\mathfrak{g}$-модуля $L(\lambda)$. Получаем

\begin{proposition}\label{star_rep}
Эрмитова форма $(\cdot,\cdot)$ в $L(\lambda)$ положительно определена.
\end{proposition}

Следующее утверждение вытекает из полупростоты $U_q\mathfrak{g}$-модулей
категории $\mathcal{C}$ и из предложения \ref{star_rep}.

\begin{corollary}\label{star_C}
Каждый $U_q\mathfrak{g}$-модуль категории $\mathcal{C}$ унитаризуем.
\end{corollary}

\begin{remark}\label{A-module}
Пусть $\lambda\in P_+$ и $\mu$-- один из весов $U\mathfrak{g}$-модуля
$L(\lambda)$. Как следует из рассмотрений предыдущего \itemiiiа, имеется
$\mathbb{Q}[t^\frac{1}{s},t^{-\frac{1}{s}}]$-модуль, доставляющий весовые
подпространства модулей $U_q\mathfrak{g}$ веса $\mu$ при специализации
$t^\frac{1}{s}=q^\frac{1}{s}$, $q\in (0,1]$.
\end{remark}

\subsection{Сводка известных результатов об алгебрах функций на компактных
квантовых группах}\label{func_alg}
\subsubsection{Алгебраические квантовые группы.}\label{alg_gr} Как
упоминалось в \itemiiiе \ref{Weight}, простая комплексная аффинная
алгебраическая группа $G$ определяется с точностью до изоморфизма парой
$(\mathbf{a}, \mathscr{L})$, где $\mathbf{a}$ -- матрица Картана ее алгебры
Ли, а $L$-- решетка, заключенная между решетками весов и корней:
$Q\subset\mathscr{L}\subset P$. Подгруппа $\mathscr{L}$ является решеткой
весов рациональных представлений группы $G$.

 Наметим построение алгебры Хопфа \IND{алгебра ! Хопфа !
$\mathbb{C}[G]$} регулярных функций на $G$. Положим
$\mathscr{L}_+=\mathscr{L}\cap P_+$. Подпространства
$(\operatorname{End}\,L(\lambda))^*$ вложены в алгебру Хопфа
$(U\mathfrak{g})^\star$ и линейно независимы. Можно показать, что
подалгебра Хопфа
\begin{equation}\label{PW}
\mathbb{C}[G]=\newoplus\limits_{\lambda\in\mathscr{L}_+}
(\operatorname{End}\,L(\lambda))^*
\end{equation}
является конечнопорожденной областью целостности и, следовательно, является
алгеброй регулярных функций на связной аффинной алгебраической группе.
\IND{алгебра ! Хопфа ! регулярных функций на связной аффинной
алгебраической группе}

Это построение неявно использует редуктивность алгебры Ли $\mathfrak{g}$ и,
как следствие, наличие у нее компактной вещественной формы. Равенство
\eqref{PW} является разложением Петера-Вейля для соответствующей компактной
вещественной алгебраической группы \cite[параграф 5.2]{VinbOn}.

 Об описании алгебры
$\mathbb{C}[G]$ с помощью образующих и соотношений см. работу Попова
\cite{Popov_ESI}.

\bigskip

Введем в рассмотрение квантовый аналог алгебры Хопфа $\mathbb{C}[G]$.
\IND{квантовый аналог ! алгебры Хопфа $\mathbb{C}[G]$} Каждому
представлению $\pi:U_q\mathfrak{g}\to\operatorname{End}\,V$ и каждому
линейному функционалу ${\varphi:\operatorname{End}\,V\to\mathbb{C}}$
сопоставим линейный функционал
\hbox{$\varphi\cdot\pi:U_q\mathfrak{g}\to\mathbb{C}$}. Такие линейные
функционалы называют матричными элементами представления $\pi$. Если
представление $\pi$ конечномерно и неприводимо, то сопряженный линейный
оператор $(\operatorname{End}V)^*\to(U_q\mathfrak{g})^*$ инъективен, как
следует из теоремы Бернсайда \cite[стр. 177]{CurtisReiner}, \cite[стр.
100]{Zhelob}.

Рассмотрим полную подкатегорию всех $\mathscr{L}$-весовых конечномерных
$U_q\mathfrak{g}$-модулей. Если $V$, $V'$, $V''$ принадлежат этой категории,
то и $V'\oplus V''$, $V''\otimes V'$, $V^*$ принадлежат ей. Тривиальный
$U_q\mathfrak{g}$-модуль также принадлежит этой категории. Значит,
подмножество векторного пространства линейных функционалов на
$U_q\mathfrak{g}$, образованное матричными элементами $\mathscr{L}$-весовых
конечномерных представлений $U_q\mathfrak{g}$, является алгеброй Хопфа. Ее
обозначают $\mathbb{C}[G]_q$ и называют алгеброй регулярных функций на
квантовой группе $G$, \IND{алгебра ! регулярных функций ! на квантовой
группе} \IND{алгебра ! $\mathbb{C}[G]_q$} см. \cite[стр. 96]{KorSoib}.
\footnote{Ее также называют координатным кольцом квантовой группы $G$ и
обозначают $\mathscr{O}_q(G)$.}

 Следующий результат
является аналогом теоремы Петера--Вейля.

\begin{proposition}\label{PeterWeyl}
  $\mathbb{C}[G]_q \cong
  \newoplus\limits_{\lambda\in \mathscr{L}_+}
   (\operatorname{End}\,L(\lambda))^*$.
\end{proposition}

{\bf Доказательство.} Каждый конечномерный весовой $U_q\mathfrak{g}$-модуль
изоморфен сумме подмодулей, кратных $L(\lambda)$ с $\lambda\in P_+$, см.
\itemiii\ \ref{loc_finite}. Следовательно, $\mathbb{C}[G]_q$ является
    суммой подпространств $(\operatorname{End}\,L(\lambda))^*$, естественно
вложенных в $\mathbb{C}[G]_q$. Остается доказать, что эта сумма является
прямой.

Пусть $U_q\mathfrak{g}^{\rm op}$ -- алгебра Хопфа, отличающаяся от
$U_q\mathfrak{g}$ заменой умножения на противоположное. Наделим
$\mathbb{C}[G]_q$ структурой $U_q\mathfrak{g}$-бимодуля, задав действие
$U_q\mathfrak{g}$ и действие $U_q\mathfrak{g}^{\rm op}$ следующим образом:
\begin{equation}\label{Rreg_Lreg}
R_{\operatorname{reg}}(\xi):\;f(\eta)\;\mapsto\;f(\eta\xi),\qquad
L_{\operatorname{reg}}(\xi):\;f(\eta)\;\mapsto\;f(\xi\eta),
\end{equation}
где $f\in\mathbb{C}[G]_q$, $\xi,\eta\in U_q\mathfrak{g}$.

Подпространства $(\operatorname{End}\,L(\lambda))^* \hookrightarrow
\mathbb{C}[G]_q$ являются простыми конечномерными попарно неизоморфными
$U_q\mathfrak{g}^{\rm op} \otimes U_q\mathfrak{g}$-модулями. Следовательно,
естественное отображение их внешней прямой суммы в $\mathbb{C}[G]_q$
инъективно (см. замечание \ref{center_Uq}, позволяющее при доказательстве
инъективности заменить некоммутативную алгебру $U_q\mathfrak{g}$ ее
центром). $\hfill\square$

\medskip
\begin{remark}
\label{bimodule} $\mathbb{C}[G]_q$ является
  $U_q\mathfrak{g}^{\rm op} \otimes U_q\mathfrak{g}$-модульной алгеброй.
\end{remark}

\begin{remark} В дальнейшем изложении представление $R_{\operatorname{reg}}$
играет существенно большую роль, чем $L_{\operatorname{reg}}$. Будем писать
$\xi f$ вместо $R_{\operatorname{reg}}(\xi)f$, в тех случаях, когда это не
приводит к недоразумениям.
\end{remark}

Введем присоединенное представление $\operatorname{Ad}_\xi$ алгебры
$U_q\mathfrak{g}$ в пространстве $\mathbb{C}[G]_q$, полагая
$$
\operatorname{Ad}_\xi:f(\eta)\mapsto\sum\limits_if(S^{-1}(\xi'_i)\,\eta\,
\xi''_i),\quad\xi,\eta\in U_q\mathfrak{g},
$$
где подразумевается, что $\Delta(\xi)=\sum\limits_i\xi'_i\otimes\xi''_i$.
\IND{представление ! присоединенное ! алгебры $U_q\mathfrak{g}$ в
пространстве $\mathbb{C}[G]_q$}

Пусть $\pi_\lambda$ --- представление алгебры $U_q\mathfrak{g}$, отвечающее
простому $\mathscr{L}$-весовому конечномерному $U_q\mathfrak{g}$-модулю
$L(\lambda)$. Элемент
\begin{equation}\label{chi_lambda}
\chi_\lambda(\xi)=\operatorname{tr}\,\pi_\lambda(\xi\cdot
K_{-2\rho}),\quad\xi\in U_q\mathfrak{g},
\end{equation}
алгебры $\mathbb{C}[G]_q$ называют характером представления $\pi_\lambda$.
\IND{характер представления}

Из результатов \itemiiiа \ref{loc_finite} следует, что $\dim {\mathrm
Hom}_{U_q\mathfrak{g}} ({\mathrm End} (L(\lambda)^*, \mathbb{C})\,=\,1$,
как и в классическом случае $q=1$. Значит, из предложения \ref{tr_g} и из
разложения Петера-Вейля вытекает

\begin{corollary}\label{characters}
Характеры $\chi_\lambda$ с $\lambda\in\mathscr{L}_+$ принадлежат
подпространствам
$(\operatorname{End}\,L(\lambda))^*\hookrightarrow\mathbb{C}[G]_q$,
являются инвариантами присоединенного действия
$$
\operatorname{Ad}_\xi
\chi_\lambda\;=\;\varepsilon(\chi_\lambda)\,f,\qquad\xi\in U_q\mathfrak{g}
$$
и образуют базис векторного пространства инвариантов.
\end{corollary}

\begin{remark} 1. Как следует из предложения \ref{separate}, естественное
спаривание алгебр Хопфа
\begin{equation}\label{pairing}
  \mathbb{C}[G]_q\times U_q\mathfrak{g}\rightarrow \mathbb{C},\qquad
  f\times\xi\mapsto f(\xi)
\end{equation}
невырождено. Значит, $U_q\mathfrak{g} \hookrightarrow \mathbb{C}[G]_q^*$.
Это аналог хорошо известного в теории вещественных групп Ли вложения
универсальной обертывающей алгебры в пространство обобщенных функций с
носителями в единице группы, см. \cite[стр. 167]{Kir}.

2. Алгебра $\mathbb{C}[G]_q$ является конечно порожденной. Действительно,
нетрудно доказать \cite[стр. 61]{BrGood}, что полугруппа $\mathscr{L}_+$
порождена конечным числом элементов и что наименьшая подалгебра, содержащая
 отвечающие им подпространства $(\operatorname{End}\,L(\lambda))^*$,
 совпадает с
$\mathbb{C}[G]_q$.
\end{remark}

\bigskip

Фундаментальная группа $\pi_1(G)$ естественно изоморфна $P/\mathscr{L}$
\cite[стр. 204]{VinbOn}. В дальнейшем, если не оговорено противное,
предполагается, что $\mathscr{L}=P$, то есть что группа $G$ односвязна.

\bigskip
\begin{remark}\label{to-L_S}
В классическом случае $q=1$ каждому подмножеству
$\mathbb{S}\subset\{1,2,\ldots,l\}$ отвечают подалгебры Ли
$\mathfrak{q}_\mathbb{S}^\pm\subset\mathfrak{g}$, порожденные множеством
$\{H_j,E_j\}_{j=1,2,\ldots,l}\;\cup\;\{F_i\}_{i\in \mathbb{S}}$, и
множеством $\{H_j,E_j\}_{j=1,2,\ldots,l}\;\cup\;\{F_i\}_{i\in\mathbb{S}}$
соответственно \cite[стр. 291]{Hum}.
 Построения настоящего \itemiiiа в полной мере применимы к
связной редуктивной подгруппе $L_\mathbb{S}\subset G$ с алгеброй Ли
 $\mathfrak{l}_\mathbb{S}=
 \mathfrak{q}_\mathbb{S}^+ \cap \mathfrak{q}_\mathbb{S}^-$.
Получаем алгебру Хопфа $\mathbb{C}[L_\mathbb{S}]_q$, являющуюся квантовым
аналогом алгебры Хопфа $\mathbb{C}[L_\mathbb{S}]$ регулярных функций на
$L_\mathbb{S}$.
\end{remark}

\bigskip
 В ряде задач теории квантовых групп трудно или
невозможно указать квантовый аналог замкнутой подгруппы, но удается найти
квантовый аналог алгебраического $G$-пространства со стабилизатором точки,
изоморфным этой подгруппе. Примером служит квантовый конус, введенный в
\itemiiiе \ref{q-cone}. Приведем несколько хорошо известных результатов
о подгруппах аффинных алгебраических групп и об алгебраических
$G$-пространствах. Разумеется, используемой топологией будет топология
Зарисского и все $G$-модули будут предполагаться конечномерными.

   Гомоморфный образ аффинной алгебраической группы является замкнутой
подгруппой.

Фактором $G/H$ называется однородное пространство $X$ с отмеченной точкой
$x$, стабилизатором которой является $H$, удовлетворяющее следующему
требованию универсальности. Для любого однородного пространства $Y$ с
отмеченной точкой $y$ и стабилизатором точки, равным $H$, существует и
единствен морфизм однородных пространств $X\to Y$, отображающий $x$ в $y$.
Фактор существует и единствен с точностью до изоморфизма.

Каждая замкнутая подгруппа $H\subset G$ является стабилизатором одномерного
подпространства некоторого $G$-модуля, что наделяет орбиту
  структурой квазипроективного алгебраического многообразия, изоморфного
$G/H$ \cite[стр. 135-146]{Hum}.

В дальнейшем все подгруппы предполагаются замкнутыми. Однородное
пространство $G/H$ аффинно если и только если подгруппа $H$ редуктивна
\cite[стр. 216]{VinbPopov}. Подгруппа $H\subset G$ называется
параболической, \IND{подгруппа ! параболическая} если $G/H$-- проективное
алгебраическое многообразие.

 Рассмотрим  проективизацию $G$-орбиты старшего вектора
конечномерного простого $G$-модуля. Как показал Костант, она замкнута и ее
однородное координатное кольцо является квадратичной алгеброй. То есть,
идеал полиномов, равных нулю на $G$-орбите старшего вектора, порожден
 элементами второй степени. Подробное доказательство  см., например,
 в \cite{GorodRudakov}.

Рассмотрим $G$-пространство $X$ и точку $x\in X$. Орбита $G x$ является
неособым подмногообразием, а ее граница-- объединением орбит строго меньшей
размерности \cite[стр. 108]{Hum}. В частности, существует замкнутая
$G$-орбита в замыкании орбиты $G x$.

Дальнейшие результаты на эту тему см. в обзорах \cite{Springer_AG},
\cite{VinbPopov} и в лекциях \cite{Steinberg}.

\subsubsection{Коммутационные соотношения.}\label{CommRelations}

Рассмотрим конечномерный $\mathscr{L}$-весовой $U_q\mathfrak{g}$-модуль
$V$. Каждой паре $v\in V$, $l\in V^*$ сопоставим линейный функционал на
$U_q\mathfrak{g}$:
\begin{equation}\label{matrix_elements}
c_{l,v}^V(\xi)\stackrel{\operatorname{def}}{=}l(\xi v),\quad \xi\in
U_q\mathfrak{g}
\end{equation}
Алгебра $\mathbb{C}[G]_q$ обладает конечным множеством образующих вида
$c_{l,v}^{L(\Lambda)}$, причем линейный функционал $l$ и вектор $v$ можно
считать весовыми.

 В дальнейшем используется сокращенное обозначение
$c_{\lambda,\mu}^\Lambda$ для элемента $c_{l,v}^{L(\Lambda)}$, где $v\in
L(\Lambda)_\mu$ , $l\in(L(\Lambda)_{\lambda})^*$, если это сокращение не
приводит к недоразумениям. Подразумевается, что
$(L(\Lambda)_{\lambda})^*\hookrightarrow L(\Lambda)^*$.

Коммутационные соотношения между элементами вида $c_{\lambda,\mu}^\Lambda$
можно получить, используя морфизмы $U_q\mathfrak{g}$-модулей
$\check{R}_{V,V'}: V \otimes V' \rightarrow V' \otimes V$, описанные в
\itemiiiе \ref{RM}. В самом деле, пусть $V$, $V'$ -- конечномерные
\hbox{$\mathscr{L}$-весовые} $U_q\mathfrak{g}$-модули и $\pi$, $\pi'$ --
отвечающие им представления алгебры $U_q\mathfrak{g}$. Тогда
\begin{equation}\label{QISM}
  (\pi'\otimes\pi)(\xi)\; \check{R}_{V,V'} =
  \check{R}_{V,V'}\; (\pi\otimes\pi')(\xi),\qquad
  \xi\in U_q\mathfrak{g}.
\end{equation}
Остается перейти от равенства операторов (\ref{QISM}) к равенству их
матричных элементов:
\begin{equation}\label{comm_BG}
  c_{l'\otimes l,\, \check{R}(v\otimes v')}^{V'\otimes V} =
  c_{\check{R}^*(l'\otimes l),\, v\otimes v'}^{V\otimes V'},
\end{equation}
где $\check{R} = \check{R}_{V,V'}$ и $v\in V$, $l\in V^*$, $v'\in V'$,
$l'\in(V')^*$. Равенства (\ref{comm_BG}) являются коммутационными
соотношениями между элементами вида $c_{l,v}^V$, $c_{l',v'}^{V'}$, поскольку
матричные элементы тензорных произведений представлений являются
произведениями матричных элементов сомножителей:
$$
c_{l\otimes l',v\otimes v'}^{V\otimes V'}=c_{l,v}^V\cdot
c_{l',v'}^{V'},\quad c_{l'\otimes l,v'\otimes v}^{V'\otimes V}=c_{l',v'}^{V'}\cdot c_{l,v}^V.
$$

Описанный выше подход к выводу коммутационных соотношений с помощью
универсальной $R$-матрицы принадлежит {В.~Дринфельду} \cite{Drinf1} и
восходит к квантовому методу обратной задачи рассеяния, развитому
 {Л.~Д.~Фадеевым} и его сотрудниками \cite{STF}.
Следующее утверждение доказывается с помощью коммутационных соотношений,
вытекающих из (\ref{comm_BG}).

\begin{proposition}\label{dom_noether}
(\cite[стр. 78]{BrGood}, \cite[стр. 266]{Jo})
Алгебра $\mathbb{C}[G]_q$ является (как правой, так и левой)
нетеровой областью целостности.
\end{proposition}

Напомним стандартное обозначение $w_0$ элемента максимальной длины группы
Вейля $W$.
\begin{proposition}\label{Rquasicom}
  (\cite[стр. 98]{KorSoib}). Пусть
$$\Lambda,\Lambda' \in P_+ ,\qquad
\lambda,\lambda',\mu,\mu'\in P$$
 и выполнены два
условия:
 $$1).\ \ \mu = \Lambda \ \ \text{или} \ \ \mu'= w_0\Lambda',\qquad
 2).\ \ \lambda = w_0\Lambda \ \text{ или} \ \ \lambda'=\Lambda'.$$
 Тогда в алгебре $\mathbb{C}[G]_q$
\begin{equation}\label{quasicom}
  c_{\lambda',\mu'}^{\Lambda'}\; c_{\lambda,\mu}^\Lambda =
  q^{(\mu,\mu')-(\lambda,\lambda')}\;
  c_{\lambda,\mu}^\Lambda\; c_{\lambda',\mu'}^{\Lambda'}.
\end{equation}
\end{proposition}

{\bf Доказательство.} Как следует из (\ref{QISM}), (\ref{comm_BG}),
(\ref{fraction}), достаточно доказать равенства
\begin{equation}\label{Rq1}
R_{L(\Lambda),L(\Lambda')}(v\otimes v')=q^{-t_0}v\otimes v',\qquad
v\otimes v'\in L(\Lambda)_\mu\otimes L(\Lambda')_{\mu'},
\end{equation}
$$
\left(R_{L(\Lambda),L(\Lambda')}\right)^*(l\otimes l')=(q^{-t_0})^*l\otimes l',\qquad l\otimes l'\in(L(\Lambda)_\lambda)^*\otimes(L(\Lambda')_{\lambda'})^*.
$$
Последнее из них означает, что
\begin{equation}\label{Rq2}
R_{L(\Lambda)^*,L(\Lambda')^*}(l\otimes l')=q^{-t_0}l\otimes l',\qquad
l\otimes l'\in(L(\Lambda)_\lambda)^*\otimes(L(\Lambda')_{\lambda'})^*,
\end{equation}
поскольку $(S\otimes S)R=R$, см. \eqref{S_theta}. Отметим, что старшие веса
$U_q\mathfrak{g}$-модулей $L(\Lambda)$, $L(\Lambda')$, $L(\Lambda)^*$,
$L(\Lambda')^*$ равны $\Lambda$, $\Lambda'$, $-w_0\Lambda$, $-w_0\Lambda'$
соответственно, а их младшие веса равны $w_0\Lambda$, $w_0\Lambda'$,
$-\Lambda$, $-\Lambda'$.

Воспользуемся разложением $U_q\mathfrak{n}^+\otimes U_q\mathfrak{n}^-$ в
сумму весовых подпространств:
$$ U_q\mathfrak{n}^+\otimes U_q\mathfrak{n}^-\;=\;
 \newoplus\limits_{\mu \in Q_+}\newoplus\limits_{\nu \in Q_+}
 (U_q\mathfrak{n}^+)_\mu \otimes (U_q\mathfrak{n}^-)_{-\nu}.$$
Из явного вида универсальной $R$-матрицы (\ref{Rmatrix}) следует, что в
пространстве формальных рядов
 \begin{equation}\label{triangle}
  Rq^{\mathbf{t}_0}-1\in \newoplus\limits_{\stackrel{\nu\in Q_+}{\nu\neq 0}}
  (U_q\mathfrak{n}^+)_\nu \otimes (U_q\mathfrak{n}^-)_{-\nu}.
  \end{equation}

  Из \eqref{fraction} следует, что равенство
(\ref{Rq1}) выполнено, если $\mu=\Lambda$ или $\mu'=w_0\Lambda'$, а
равенство (\ref{Rq2}) выполнено, если $\lambda=w_0\Lambda$ или
$\lambda'=\Lambda'$. $\hfill\square$

\begin{corollary}\label{subalgebra}
Подалгебра, порожденная элементами
\begin{equation}\label{gen_sub}
c_{w_0\Lambda,\Lambda}^\Lambda,\;\;c_{\Lambda,w_0\Lambda}^\Lambda,\qquad
\Lambda\in\mathscr{L}_+,
\end{equation}
коммутативна, и ее образующие (\ref{gen_sub}) квазикоммутируют со всеми
элементами $c_{\lambda',\mu'}^{\Lambda'}\in\mathbb{C}[G]_q$:
$$
c_{\lambda',\mu'}^{\Lambda'}\;c_{w_0\Lambda,\Lambda}^\Lambda=
q^{(\Lambda,\mu')-(w_0\Lambda,\lambda')}\;c_{w_0\Lambda,\Lambda}^\Lambda\;
c_{\lambda',\mu'}^{\Lambda'},
$$
$$
c_{\Lambda,w_0\Lambda}^\Lambda\;c_{\lambda',\mu'}^{\Lambda'}=
q^{(\mu',w_0\Lambda)-(\lambda',\Lambda)}\;c_{\lambda',\mu'}^{\Lambda'}\;
c_{\Lambda,w_0\Lambda}^\Lambda.
$$
\end{corollary}

\medskip

\begin{note}\label{infinite_dim_R} Отбросим требование $\lambda \in P$ и
 заменим конечномерные весовые $U_q\mathfrak{g}$-модули $L(\lambda)$
 модулями Верма
 $M(\lambda)$. Матричные элементы $c_{\lambda,\mu}^\Lambda$,
 $c_{\lambda',\mu'}^{\Lambda'}$ отвечающих им представлений
 являются элементами алгебры
 $(U_q\mathfrak{g})^*$. Как видно из доказательства предложения
 \ref{Rquasicom}, коммутационные соотношения (\ref{quasicom}) справедливы и
 в этом более общем случае, если $\mu=\Lambda$ и $\lambda'=\Lambda'$.
\end{note}

\subsubsection{Компактные квантовые
группы.}\label{regular_K}\label{compact}

Получим квантовый аналог компактной формы $K$ комплексной редуктивной
группы $G$ \cite{VinbOn}. Вещественная алгебра Ли называется компактной,
\IND{компактная вещественная алгебра Ли} если она обладает инвариантной
относительно присоединенного действия положительно определенной
невырожденной симметрической билинейной формой. Такие подалгебры сопряжены
\cite[стр. 26]{Bou9}. Примером может служить вещественная подалгебра Ли
простой комплексной алгебры Ли $\mathfrak{g}$, порожденная множеством
$\{iH_j,\,E_j-F_j,\,i(E_j+F_j)\}_{j=1,2,\ldots,l}$, см. \cite[стр.
139]{VinbGorbOn}.

 Рассмотрим $*$-алгебру Хопфа $(U_q\mathfrak{g}, \star)$, введенную в
\itemiiiе \ref{HChV} равенствами (\ref{aster1}) и
двойственную
 $*$-алгебру Хопфа $(U_q\mathfrak{g})^\star$:
$$
f^\star(\xi) \stackrel{\rm def}{=} \overline{f((S(\xi))^\star)},\qquad
\xi\in U_q\mathfrak{g},\; f\in(U_q\mathfrak{g})^\star,
$$
см. (\ref{invol1}).

Напомним, что ее элементами являются матричные элементы конечномерных
представлений алгебры $U_q\mathfrak{g}$, не обязательно весовых. Подалгебра
Хопфа $\mathbb{C}[G]_q \hookrightarrow (U_q\mathfrak{g})^\star$ наследует
эту инволюцию, поскольку $f^*$ является матричным элементом представления
$\pi^*$, если $f$ --- матричный элемент конечномерного весового
представления $\pi$. (В дальнейшем мы уточним последнее утверждение, см.
\eqref{aster_c}.)

Следуя Дринфельду, будем называть $*$-алгебру Хопфа
\begin{equation}\label{compact_form}
\mathbb{C}[K]_q\stackrel{\operatorname{def}}{=}(\mathbb{C}[G]_q,\star)
\end{equation}
алгеброй регулярных функций на вещественной аффинной алгебраической группе
$K$. \IND{алгебра ! регулярных функций ! на вещественной аффинной
алгебраической группе}

\begin{remark}\label{compact_group}
В классическом случае $q=1$ инволюция $\star$ в $\mathbb{C}[G]$ вводится
точно так же, и
\begin{equation}\label{compact_form_classic}
K\,=\,\left\{\,g\in G\;\left|\;f^\star(g)=\overline{f(g)}\;\;\text{при
всех}\;\;f\in\mathbb{C}[G]\,\right.\right\}.
\end{equation}
\end{remark}

\bigskip В \itemiiiе \ref{loc_finite} конечномерные весовые
$U_q\mathfrak{g}$-модули $L(\Lambda)$ были наделены инвариантными
положительно определенными эрмитовыми формами $(\cdot,\cdot)$, для которых
$(v(\lambda),v(\lambda))=1$. Выберем ортонормированный базис весовых
векторов $\left\{v_{\mu,j}^\Lambda\right\}$ в каждом из конечномерных
весовых подпространств $L(\Lambda)_\mu$. Если это подпространство
$L(\Lambda)_\mu$ одномерно, то вместо $v_{\mu,1}^\Lambda$ будем для
краткости писать $v_\mu^\Lambda$.

Введем обозначение
$$
c_{\lambda,i;\mu,j}^\Lambda(\xi)\stackrel{\operatorname{def}}{=}\left(\xi
v_{\mu,j}^\Lambda, v_{\lambda,i}^\Lambda\right)
$$
для матричных элементов представлений $\pi_\Lambda$, отвечающих
$U_q\mathfrak{g}$-модулям $L(\Lambda)$. Желая обеспечить соответствие с
обозначениями раздела \ref{CommRelations}, вместо
$c_{\lambda,i;\mu,j}^\Lambda$ будем писать $c_{\lambda,\mu}^\Lambda$, если
это не приводит к недоразумениям.

Множество $\left\{c_{\lambda,i;\mu,j}^\Lambda\right\}$ является базисом
весовых векторов $U_q\mathfrak{g}^{\rm op}\otimes U_q\mathfrak{g}$-модуля
$\mathbb{C}[K]_q$:
\begin{equation}\label{rweight}
R_{\operatorname{reg}}(K_m^{\pm1})\,c_{\lambda,i;\mu,j}^\Lambda\;=\;
q_m^{\pm\mu_m}\,c_{\lambda,i;\mu,j}^\Lambda,
\end{equation}
\begin{equation}\label{lweight}
L_{\operatorname{reg}}(K_m^{\pm1})\,c_{\lambda,i;\mu,j}^\Lambda\;=\;
q_m^{\pm\lambda_m}\,c_{\lambda,i;\mu,j}^\Lambda,
\end{equation}
где $\lambda=\sum\limits_k\lambda_k\,\overline{\omega}_k$,
$\mu=\sum\limits_k\mu_k\,\overline{\omega}_k$.

Из определений следует, что коумножение $\Delta$ и коединица $\varepsilon$
действуют на матричные элементы $c_{\lambda,i;\mu,j}^\Lambda$ следующим
образом:
\begin{equation}\label{comult_c}
\Delta\left(c_{\lambda,i;\mu,j}^\Lambda\right)=
\sum\limits_{\nu,m}c_{\lambda,i;\nu,m}^\Lambda\otimes
c_{\nu,m;\mu,j}^\Lambda,
\end{equation}
\begin{equation}\label{counit_c}
\varepsilon\left(c_{\lambda,i;\mu,j}^\Lambda\right)=
\delta_{\lambda,\mu}\delta_{i,j}.
\end{equation}

Для того, чтобы найти явный вид элементов
$S\left(c_{\lambda,i;\mu,j}^\Lambda\right)$,
$\left(c_{\lambda,i;\mu,j}^\Lambda\right)^\star$,
наделим $U_q\mathfrak{g}$ {\bf антилинейной} инволюцией $\tau$
\begin{equation}\label{theta1}
\tau(K_j^{\pm1})=S(K_j^{\pm1}),\quad\tau(E_j)=q_jS(E_j),\quad
\tau(F_j)=q_j^{-1}S(F_j),
\end{equation}
где $j=1,2,\ldots,l$. Каждому $\Lambda\in P_+$ отвечает
$U_q\mathfrak{g}$-модуль $L(\Lambda)$ и $*$-представление $\pi_\Lambda$ в
конечномерном гильбертовом пространстве $L(\Lambda)$. Пусть
$\pi_\Lambda^\tau$ -- представление алгебры $U_q\mathfrak{g}$ в
$L(\Lambda)$, определяемое равенством
$$
(\pi_\Lambda^\tau(\xi)v_1,v_2)=(v_1,\pi_\Lambda(\tau(\xi))v_2),\qquad
\xi\in U_q\mathfrak{g},\;v_1,v_2\in L(\Lambda).
$$
Так как $\star\tau=\tau\star$, то
$$
(\pi_\Lambda^\tau(\xi)v_1,v_2)=(v_1,\pi_\Lambda(\tau(\xi))v_2)=
(\pi_\Lambda((\tau(\xi))^\star)v_1,v_2)=
$$
$$
=(\pi_\Lambda(\tau(\xi^\star))v_1,v_2)=(v_1,\pi_\Lambda^\tau(\xi^\star)v_2).
$$
Значит, $\pi_\Lambda^\tau$ является $*$-представлением со старшим весом
$-w_0\Lambda$.

   Следующее утверждение несущественно отличается от
аналогичного результата Сойбельмана \cite[стр. 100]{KorSoib}.

\begin{proposition}\label{invol_c}
Пусть $U:L(\Lambda)\to L(-w_0\Lambda)$ -- единственный с точностью до
числового множителя унитарный оператор, для которого
$\pi_\Lambda^\tau(\xi)=U^{-1}\,\pi_{-w_0\Lambda}(\xi)\,U$ при всех $\xi\in
U_q\mathfrak{g}$. В $*$-алгебре $\mathbb{C}[K]_q$
\begin{equation}\label{antipode_c}
S\left(c_{\lambda,i;\mu,j}^\Lambda\right)=q^{(\mu-\lambda,\rho)}\cdot
\widetilde{c}_{-\mu,j;-\lambda,i}^{\;-w_0\Lambda},
\end{equation}
\begin{equation}\label{aster_c}
\left(c_{\lambda,i;\mu,j}^\Lambda\right)^\star=q^{(\lambda-\mu,\rho)}\cdot
\widetilde{c}_{-\lambda,i;-\mu,j}^{\;-w_0\Lambda},
\end{equation}
где $\left\{c_{\lambda,i;\mu,j}^\Lambda\right\}$,
$\left\{\widetilde{c}_{\lambda,i;\mu,j}^{\;-w_0\Lambda}\right\}$ --
матричные элементы представлений $\pi_\Lambda$, $\pi_{-w_0\Lambda}$ в
базисах $\left\{v_{\mu,j}^\Lambda\right\}$,
$\left\{U\,v_{-\mu,j}^\Lambda\right\}$.
\end{proposition}

{\bf Доказательство.} Рассмотрим автоморфизм $\xi\mapsto K_\rho^{-1}\xi
*_\rho$ алгебры $U_q\mathfrak{g}$. Он корректно определен, см.
\itemiii\ \ref{Weight}. Рассуждая так же, как при доказательстве
(\ref{inner_S2}), и используя (\ref{theta1}), получаем равенство $
S(\xi)=K_\rho^{-1}\tau(\xi)K_\rho$ для любого элемента $\xi$ наименьшей
вещественной подалгебры, содержащей все образующие $K_i^{\pm1}$, $E_i$,
$F_i$, $i=1,2,\ldots, l$. Значит, для таких элементов $\xi$
$$
c_{\lambda,i;\mu,j}^\Lambda(S(\xi))=
q^{(\mu-\lambda,\rho)}c_{\lambda,i;\mu,j}^\Lambda(\tau(\xi)).
$$
Из полученного равенства вытекает (\ref{antipode_c}). При всех $\xi\in
U_q\mathfrak{g}$
$$
c_{\lambda,i;\mu,j}^\Lambda(S(\xi))=
\overline{c_{\mu,j;\lambda,i}^\Lambda((S(\xi))^\star)} =
(c_{\mu,j;\lambda,i}^\Lambda)^\star(\xi).
$$
Значит,
\begin{equation}\label{S_aster}
S\left(c_{\lambda,i;\mu,j}^\Lambda\right)=
\left(c_{\mu,j;\lambda,i}^\Lambda\right)^\star,
\end{equation}
и равенство (\ref{aster_c}) следует из (\ref{antipode_c}),
(\ref{S_aster}). $\hfill\square$


Следующие равенства вытекают из определения антипода алгебры Хопфа и из
(\ref{comult_c}), (\ref{counit_c}), (\ref{S_aster}):
\begin{equation}\label{unitary_1}
\sum\limits_{\mu,j}c_{\lambda,i;\mu,j}^\Lambda\;
\left(c_{\lambda,i;\mu,j}^\Lambda\right)^\star=1,
\end{equation}
\begin{equation}\label{unitary_2}
\sum\limits_{\lambda,i}\left(c_{\lambda,i;\mu,j}^\Lambda\right)^\star\;
c_{\lambda,i;\mu,j}^\Lambda=1.
\end{equation}
Из них следует, что для любого элемента $f\in\mathbb{C}[K]_q$
\begin{equation}\label{infty_norm}
\|f\|_\infty\stackrel{\rm def}{=}\sup\limits_{T}\|T(f)\|<\infty,
\end{equation}
где $T$ пробегает множество классов унитарной эквивалентности
$*$~-~представлений алгебры $\mathbb{C}[K]_q$ в гильбертовых пространствах.
Построим точное $*$~-~представление алгебры $\mathbb{C}[K]_q$. Тем самым
будет доказано, что $\|f\|_\infty$ является нормой, а $\mathbb{C}[K]_q$ с
этой нормой -- нормированная алгебра.

\begin{note}\label{irrep_only}
Хорошо известно \cite[стр.~53]{Dix2}, что та же самая норма получается,
если в \eqref{infty_norm} использовать не все $*$-представления алгебры
$\mathbb{C}[K]_q$, а только неприводимые $*$-представления.
\end{note}

Введем в рассмотрение линейный функционал
$$
\nu:\mathbb{C}[K]_q\to\mathbb{C},\qquad\nu:f\mapsto\int\limits_{K_q}f\,d\nu,
$$
равный нулю на гиперплоскости $\newoplus\limits_{\lambda\ne
0}(\mathrm{End}\,L(\lambda))^*\subset\mathbb{C}[K]_q$ и единице на
единичном элементе.
Очевидно,
\begin{equation}\label{invar}
(\nu\otimes\operatorname{id})\Delta(f)=
(\operatorname{id}\otimes\nu)\Delta(f)=\nu(f)\cdot 1,\qquad
f\in\mathbb{C}[K]_q,
\end{equation}
и $\nu$ является единственным с точностью до числового множителя
$(U_q\mathfrak{g})^{\rm op}\otimes U_q\mathfrak{g}$-инвариантным линейным
функционалом. Это квантовый аналог интеграла по мере Хаара на $K$.
\IND{квантовый аналог ! интеграла по мере Хаара}

Из равенства \eqref{inner_S2} следует, что представления
$\pi_\Lambda(\xi)$
и $\pi_\Lambda(S^2(\xi))$ алгебры $U_q\mathfrak{g}$ в $L(\Lambda)$
эквивалентны. Именно, линейный оператор
$$
F_\Lambda: L(\Lambda)\rightarrow L(\Lambda),\quad F_\Lambda:
v_{\mu,j}^\Lambda\mapsto q^{-2(\mu,\rho)}\ v_{\mu,j}^\Lambda,
$$
является сплетающим:
$$\pi_\Lambda(S^2(\xi))F_\Lambda = F_\Lambda \
\pi_\Lambda(\xi), \qquad \xi\in U_q\mathfrak{g}.$$ Как объясняется в
\cite[стр. 26]{SoibVak92}, зная явный вид сплетающего оператора, нетрудно
получить следующие соотношения ортогональности.
\begin{proposition}\label{ort_matrix}
\begin{equation}\label{ort_eq}
\int\limits_{K_q} (c_{\lambda,i;\mu,j}^\Lambda)^\star\;
c_{\lambda',i';\mu',j'}^{\Lambda'} d\nu = \frac{\delta_{\Lambda,\Lambda'}
\delta_{\lambda,\lambda'} \delta_{i,i'}\delta_{\mu,\mu'} \delta_{j,j'}
q^{2(\lambda,\rho)}}{\operatorname{dim}_q L(\Lambda)},
\end{equation}
где
$$
\operatorname{dim}_q L(\Lambda) \stackrel{\operatorname{def}}{=}
\sum\limits_{\lambda \in P} \dim L(\Lambda)_\lambda \;
q^{2(\lambda,\rho)}.
$$
\end{proposition}

\begin{corollary}
Полуторалинейная форма
$$
(\psi_1,\psi_2)=\int\limits_{K_q}\psi_2^\star\psi_1 d\nu,\quad
\psi_1,\psi_2\in\mathbb{C}[K]_q,
$$
положительно определена.
\end{corollary}

  Следующее утверждение вытекает из \eqref{ort_eq}, \eqref{S_aster} и из определения
  характера \eqref{chi_lambda}, согласно которому
  $ \chi_\Lambda=\sum\limits_{\lambda ,i}\;
  c_{\lambda,i;\lambda,i}^\Lambda \;q^{-2(\lambda,\rho)}
  $.

\begin{corollary}\label{projectors}
  Проектор $P_\Lambda$ в $\mathbb{C}[K]_q$ на
  $(\operatorname{End}\,L(\Lambda))^*$ параллельно
  $\bigoplus\limits_{\lambda\neq\Lambda}
  (\operatorname{End}\,L(\lambda))^*$ является
интегральным оператором с ядром
$$
\mathscr{P}_\Lambda = \operatorname{dim}_q L(\Lambda) \cdot
(\operatorname{id}\otimes S) \Delta(\chi_\Lambda).
$$
\end{corollary}

\medskip

Пополняя $\mathbb{C}[K]_q$ по норме
$\|\psi\|_2 = \left(\int\limits_{K_q} \psi^\star\psi\;
d\nu\right)^{1/2},$
получаем гильбертово пространство $L^2(d\nu)_q$, которое является
квантовым аналогом пространства $L^2(d\nu)$. \IND{квантовый аналог !
пространства $L^2(d\nu)$}

\begin{remark}\label{norm_character}
Так как $w_0\rho=-\rho$, то
$$
\rm{dim}_qL(-w_0\Lambda)\;=\;\sum\limits_{\lambda\in P}
\dim(L(-w_0\Lambda)_{-w_0\lambda})q^{2(\lambda,\rho)}\;=\;\rm{dim}_q
L(\Lambda),
$$
где $w_0$-элемент максимальной длины группы Вейля $W$.
\end{remark}

Из \eqref{unitary_2} и из предложения \ref{dom_noether} вытекает

\begin{proposition}
Представление $*$-алгебры $\mathbb{C}[K]_q$ в предгильбертовом
пространстве
$\mathbb{C}[K]_q$, которое каждому элементу $f$ сопоставляет линейный
оператор умножения слева на $f$,
является ее точным $*$-представлением, и все операторы этого
представления
ограничены.
\end{proposition}

\begin{corollary}
 Продолжая по непрерывности операторы рассматриваемого представления
 на все пространство $L^2(d\nu)_q$, получаем точное $*$-представление
 алгебры $\mathbb{C}[K]_q$
 в этом гильбертовом пространстве.
\end{corollary}

\bigskip

Пополнение $C(K)_q$ $*$-алгебры $\mathbb{C}[K]_q$ по норме $\|f\|_\infty$
является $C^*$~-~алгеброй и называется алгеброй непрерывных функций на
квантовой группе $K$. \IND{алгебра ! непрерывных функций на квантовой
группе}

Рассмотрим пространственное тензорное произведение $C^*$-алгебр
$C(K)_q\otimes C(K)_q$ \cite[стр.~71]{Dix2}. Из определений следует, что
структура биалгебры по непрерывности переносится с плотной подалгебры
$\mathbb{C}[K]_q$ на $C(K)_q$:
$$
\Delta: C(K)_q\rightarrow C(K)_q\otimes C(K)_q;\qquad \varepsilon:
C(K)_q\rightarrow \mathbb{C}.
$$

Инвариантный интеграл $\nu$ допускает продолжение по непрерывности на
$C(K)_q$. Действительно, $\|f\|_2\leq \|f\|_\infty$, поскольку из
определения $\|\cdot\|_\infty$ вытекает неравенство $\|f\varphi\|_2\leq
\|f\|_\infty\,\|\varphi\|_2$, и $\|1\|_2=1$. Следовательно,
$|\int\limits_{K_q}\,f\,d\nu|\le\|f\|_2\,\|1\|_2\,=\,\|f\|_2\le
\|f\|_\infty$.

 Очевидно, равенства (\ref{invar}) имеют место при всех $f\in C(K)_q$.

\medskip

\begin{remark}
Подмножество
$\left\{c_{\lambda,i;\mu,j}^\rho\right\}\subset\mathbb{C}[K]_q$ и
$C^*$-алгебра $C(K)_q$ удовлетворяют всем требованиям, включенным С.
Вороновичем в определение компактной матричной псевдогруппы \cite[стр.
616]{Woron1}.
\end{remark}

\subsubsection{Пример: алгебра функций на квантовой группе $SU_2$.}\label{SU_2}
\label{int_SU2}

Рассмотрим алгебру Хопфа $U_q\mathfrak{sl}_2$, $*$-алгебру Хопфа
$U_q\mathfrak{su}_2=(U_q\mathfrak{sl}_2,\star)$ и алгебру Хопфа
$\mathbb{C}[SL_2]_q$, введенные в \itemiiiах \ref{q-sl_2},
\ref{star-Hopf_algebras} и \ref{SL_2_sl_2}.

Наделяя $\mathbb{C}[SL_2]_q$ инволюцией по двойственности, см. предыдущий
\itemiii, получаем $*$-алгебру Хопфа $\mathbb{C}[SU_2]_q$, являющуюся квантовым
аналогом алгебры регулярных функций на квантовой группе $SU_2$.
\IND{$*$-алгебра ! Хопфа ! $\mathbb{C}[SU_2]_q$} Из определений вытекает
равенство
\begin{equation}\label{ast2}
\left(
\begin{array}{cc}
t_{11}^\star, & t_{12}^\star
\\ t_{21}^\star, & t_{22}^\star
\end{array}
\right)=
\left(
\begin{array}{cc}
t_{22}, & -qt_{21}
\\ -q^{-1}t_{12}, & t_{11}
\end{array}
\right).
\end{equation}

Используя описание алгебры $\mathbb{C}[SL_2]_q$ в терминах образующих и
соотношений, приведенное в \itemiiiе \ref{SL_2_sl_2}, найдем неприводимые
$*$-представления алгебры $\mathbb{C}[SU_2]_q$.

\begin{example}\label{one-dim}
Алгебра $\mathbb{C}[SU_2]_q$ обладает одномерными
$*$-пред\-став\-ле\-ни\-я\-ми
$$
\chi_\varphi(t_{11})=e^{i\varphi},\quad
\chi_\varphi(t_{12})=\chi_\varphi(t_{21})=0,\quad
\chi_\varphi(t_{22})=e^{-i\varphi},\quad\text{где}\quad
\varphi\in\mathbb{R}/(2\pi\mathbb{Z}).
$$
\end{example}

\begin{example}\label{infinite_dim_rep}
  Пусть $\{e_j\}_{j=0}^\infty$ -- стандартный базис гильбертова пространства
$l^2(\mathbb{Z}_+)$. Следующие равенства определяют неприводимое
бесконечномерное $*$-представление $\Pi$ алгебры $\mathbb{C}[SU_2]_q$ в
$l^2(\mathbb{Z}_+)$:
$$ \Pi(t_{11})e_j=\begin{cases}
\sqrt{1-q^{2j}}\,e_{j-1},\;&\; j>0,\\ 0,\;&\;j=0
\end{cases},
 \qquad \Pi(t_{12})e_j=q^{j+1}e_j,$$
 $$ \Pi(t_{22})e_j=\sqrt{1-q^{2(j+1)}}\,e_{j+1}, \qquad \qquad \qquad
  \Pi(t_{21})e_j=q^j e_j.
$$
\end{example}

\medskip

\begin{remark}\label{inequal} При всех $f\in C(SU_2)_q$
\begin{equation}\label{norms}
\|\Pi(f)\|\ge\max\limits_{\varphi\in\mathbb{R}/(2\pi\mathbb{Z})}|\chi_\varphi(f)|.
\end{equation}
Действительно, пусть $\mathcal{B}$ -- алгебра всех ограниченных операторов
в $l^2(\mathbb{Z}_+)$, $\mathcal{K}$ -- ее двусторонний идеал компактных
операторов и $j:\mathscr{B}\to\mathcal{B}/\mathcal{K}$ -- каноническая
факторизация. Тогда в алгебре Калкина $\mathcal{B}/\mathcal{K}$
$$
\|j\cdot\Pi(f)\|=\max\limits_{\varphi\in\mathbb{R}/(2\pi\mathbb{Z})}|\chi_\varphi(f)|,\qquad
f\in\mathbb{C}[SU_2]_q.
$$
Остается воспользоваться тем, что $\|j\|\le 1$.
\end{remark}

\bigskip

Первое утверждение следующего предложения использует коумножение в
$\mathbb{C}[SU_2]_q$ и почти очевидно. Второе утверждение доказывается
методом, описанным в \itemiiiе \ref{Fock_l}. Роль коммутационных соотношений
\eqref{comm_disc} при этом играют коммутационные соотношения
$$
t_{11}\,x\,=\,q^2\,x\,t_{11},\quad t_{12}\,x\,=\,x\,t_{12},\quad
t_{21}\,x\,=\,x\,t_{21},\quad t_{22}\,x\,=\,q^{-2}\,x\,t_{22},
$$
где $x=t_{12}t_{12}^\star$.

\begin{proposition}\label{list_SU2}
(\cite{VakSoib88}, стр. 5)\ \ \ $*$-Представления алгебры
$\mathbb{C}[SU_2]_q$
\begin{equation}\label{list_rep}
\chi_\varphi,\qquad \Pi\otimes \chi_\varphi,
\end{equation}
где $\varphi\in\mathbb{R}/(2\pi\mathbb{Z})$, неприводимы и попарно
неэквивалентны.
 Каждое неприводимое $*$~-~представление алгебры $\mathbb{C}[SU_2]_q$
унитарно эквивалентно одному из представлений (\ref{list_rep}).
\end{proposition}

\begin{corollary}
В $C^*$-алгебре $C(SU_2)_q$
\begin{equation}\label{norm_SU2}
\|f\|=\max\limits_{\varphi\in\mathbb{R}/(2\pi\mathbb{Z})}\|(\Pi\otimes\chi_\varphi)(f)\|,\qquad
f\in\mathbb{C}[SU_2]_q.
\end{equation}
\end{corollary}

\medskip

В классическом случае $q=1$ спектром $\operatorname{spec}(x)$ элемента
$x=|t_{12}|^2$ банаховой алгебры $C(SU_2)$ является отрезок $[0,1]$. В
квантовом случае спектр $x$ оказывается счетным множеством.

\begin{corollary}
  В $C^*$-алгебре $C(SU_2)_q$
\begin{equation}\label{spectr_x}
  \operatorname{spec}(x)=q^{2\mathbb{N}} \ \cup \ \{0\}.
\end{equation}
\end{corollary}

\bigskip
Известен явный вид инвариантного интеграла, равного единице на единичном
элементе.
\begin{proposition}
(\cite{VakSoib88}, стр. 11). Для всех $f\in \mathbb{C}[SU_2]_q$
\begin{equation}\label{intSU2}
  \int\limits_{(SU_2)_q} f d\nu = \frac{q^{-2} - 1}{2\pi}
\int\limits_0^{2\pi}\operatorname{tr}((\Pi\otimes \chi_\varphi)(f\cdot x))
d\varphi.
\end{equation}
\end{proposition}

\medskip Перейдем к существенно более общим результатам.

\subsubsection{Неприводимые $*$-представления алгебры $\mathbb{C}[K]_q$ и
инвариантный интеграл.}\label{IrrepK}

Построение неприводимых $*$-представлений алгебры $\mathbb{C}[K]_q$
осуществляется сведением к частному случаю $K=SU_2$. Выберем
$i\in\{1,2,\ldots,l\}$ и рассмотрим вложение $*$-алгебр Хопфа
$$
U_{q_i}\mathfrak{sl}_2 \hookrightarrow U_q\mathfrak{g},\qquad K^{\pm1}
\mapsto K_i^{\pm1},\quad E \mapsto E_i,\quad F \mapsto F_i.
$$
Его образ будем обозначать $(U_{q_i} \mathfrak{sl}_2)_i$\ .

Сопряженный линейный оператор является гомоморфизмом $*$-алгебр Хопфа
$\mathbb{C}[K]_q\to\mathbb{C}[SU_2]_{q_i}$, поскольку сужение
конечномерного весового представления алгебры $U_q\mathfrak{g}$ на ее
подалгебру $(U_{q_i}\mathfrak{sl}_2)_i$ является конечномерным весовым
представлением этой подалгебры. Возникают гомоморфизмы $*$-алгебр Хопфа
$$\psi_i:\mathbb{C}[K]_q\to\mathbb{C}[SU_2]_{q_i},\qquad i=1,2,\ldots,l,$$
и неприводимые $*$-представления алгебры $\mathbb{C}[K]_q$:
$$
\Pi_i=\Pi\circ\psi_i,\qquad\qquad\chi_{\varphi,i}=
\chi_\varphi\circ\psi_i,\quad\varphi\in\mathbb{R}/(2\pi\mathbb{Z}).
$$

Получен $l$-мерный тор $\mathbb{T}^l=(\mathbb{R}/(2\pi\mathbb{Z}))^l$
одномерных $*$-представлений
$$
\chi_{\varphi}=\chi_{\varphi_1,1}\otimes\chi_{\varphi_2,2}\otimes\ldots
\otimes\chi_{\varphi_l,l},\quad
\varphi=(\varphi_1,\varphi_2,\ldots,\varphi_l)\in\mathbb{T}^l.
$$
Естественно ожидать, что неприводимы и некоторые из тензорных произведений
$\Pi_{i_1}\otimes\Pi_{i_2}\otimes\ldots\otimes\Pi_{i_N}$. Подчеркнем, что
операторы представления ограничены (см. (\ref{unitary_1}),
(\ref{unitary_2})), и нас интересуют тензорные произведения в категории
гильбертовых пространств. Неприводимость в этом контексте означает
отсутствие общих нетривиальных инвариантных подпространств, а не отсутствие
общих нетривиальных инвариантных линейных подмногообразий операторов
представления.

\begin{proposition}\label{Soib_list}
  (\cite{KorSoib}, стр. 121).
Пусть $w=s_{i_1}s_{i_2}\ldots s_{i_N}$ -- приведенное разложение элемента
$w\in W$ и $\varphi \in \mathbb{T}^l$. Рассмотрим $*$-представление алгебры
$\mathbb{C}[K]_q$
\begin{equation}\label{tensor_rep}
  \Pi_{i_1} \otimes \Pi_{i_2} \otimes \ldots \otimes \Pi_{i_N} \otimes
\chi_\varphi.
\end{equation}
1. $*$-Представление (\ref{tensor_rep}) неприводимо и каждое неприводимое
\hbox{$*$-представление} алгебры $\mathbb{C}[K]_q$ унитарно эквивалентно
одному из представлений (\ref{tensor_rep}).
\\ 2. Два $*$-представления (\ref{tensor_rep}) унитарно эквивалентны, если и
только если совпадают отвечающие им элементы $w\in W$ и элементы
$(\mathbb{R}/(2\pi\mathbb{Z}))^l$.
\end{proposition}

\medskip
Используя (\ref{norm_SU2}) и предложение \ref{Soib_list}, получаем

\begin{corollary}
Пусть $w_0=s_{i_1}s_{i_2}\ldots s_{i_M}$ -- приведенное разложение элемента
максимальной длины и
$$
\Pi_{w_0} = \Pi_{i_1} \otimes \Pi_{i_2} \otimes \ldots \otimes \Pi_{i_M}.
$$
Тогда для любого элемента $f\in \mathbb{C}[K]_q$
\begin{equation}\label{tore}
  \|f\|_\infty = \max\limits_{\varphi \in \mathbb{T}^l}
  \|(\Pi_{w_0}\otimes\chi_\varphi)(f)\|.
\end{equation}
\end{corollary}

\medskip

Из предложения \ref{Soib_list} также следует, что $*$-представление
$\Pi_w$, $w\in W$, алгебры $\mathbb{C}[K]_q$ определяется элементом $w$ с
точностью до унитарной эквивалентности, то есть не зависит от выбора
приведенного разложения.

\bigskip

Как показал Стокман в \cite{Stokman_list}, используя результаты работы
\cite{DijStok1}, аналогичное описание $*$-представлений можно получить для
 важного класса однородных пространств квантовых групп.

Обратимся к классическому случаю $q=1$. Рассмотрим подмножество
 $\mathbb{S} \subset \{1,2,\ldots l\}$ и
алгебраический тор $\mathbb{T}\subset G$
  с алгеброй Ли
 $$\{\xi \in  \mathfrak{h}
\;|\; \alpha_j(\xi)=0\;\;\text{при}\;\; j \in \mathbb{S}\},$$
 см. \cite[стр. 133]{VinbOn}.
  Централизатор $L_\mathbb{S} \subset G$ этого тора является
связной редуктивной группой
 с алгеброй Ли $\mathfrak{l}_\mathbb{S}\subset
\mathfrak{g}$, порожденной множеством
 $\{ F_i, E_j\;|\; i \in \mathbb{S}\} \cup \{H_j\;|\; j=1,2,\ldots l\}$
 \cite[стр. 202, 254]{Hum}.
 Значит, $L_\mathbb{S}$-инвариантные  функции на $G$ выделяются требованием
 их $U\mathfrak{l}$-инвариантности, и оправданы  следующие построения.

Пусть $\mathbb{S}\subset \{1,2,\ldots,l\}$ и $U_q\mathfrak{l}_\mathbb{S}
\subset U_q\mathfrak{g}$ -- подалгебра Хопфа, порожденная множеством
$$
\{K_i^{\pm 1}\}_{i=1,2,\ldots,l}\;\; \cup\;\; \{E_j, F_j\}_{j\in
\mathbb{S}}.
$$
Рассмотрим подалгебру
$$ \mathbb{C}[L_\mathbb{S}\backslash G]_q
 \stackrel{\operatorname{def}}{=} \left\{
f\in \mathbb{C}[G]_q\; |\; L_{\operatorname{reg}}(\xi)f =
\varepsilon(\xi)f,\quad \xi\in U_q\mathfrak{l}_\mathbb{S} \right\}.
$$
Она наследует структуру локально конечномерной $U_q\mathfrak{g}$-модульной
алгебры, поскольку
$$
L_{\operatorname{reg}}(\xi)\cdot R_{\operatorname{reg}}(\eta) =
R_{\operatorname{reg}}(\eta)\cdot L_{\operatorname{reg}}(\xi),\qquad \xi,
\eta \in U_q\mathfrak{g}.
$$

 Подалгебра $\mathbb{C}[L_\mathbb{S}\backslash G]_q$
наследует инволюцию $\star$. Известно \cite[стр. 474]{Stokman_list},
\cite{HeckenbergerKolb04}, что она порождается множеством
 $\{c_{\overline{\omega}_{i};\lambda,j}^{\overline{\omega}_{i}}
 (c_{\overline{\omega}_{i};\mu,k} ^{\overline{\omega}_{i}})^\star\;
 |\; i\in \mathbb{S}\}$.

 Займемся представлениями алгебры  с инволюцией
 $\mathbb{C}[L_\mathbb{S}\backslash G]_q$.
Пусть $W_\mathbb{S}\subset W$ -- подгруппа, порожденная множеством
 $\{s_i\}_{i\in \mathbb{S}}$, и
$$^\mathbb{S}W=\left\{w\in W\,|\,l(w s_i)>l(w), \quad i\in
\mathbb{S}\right\}.$$


Следующий результат получен Костантом в \cite{Kostant_1961}.

\begin{lemma}\label{w_s}
(\cite{Hum2}, стр. 19) 1. Каждый элемент $w\in W$ единственным образом
разлагается в произведение
$$ w\,=\,^\mathbb{S}w\cdot\,w_{\mathbb{S}},\qquad \quad
w_{\mathbb{S}}\in W_{\mathbb{S}},\;\;\; ^{\mathbb{S}}w\in\,
^{\mathbb{S}}W;$$
 2.\ \ \  $^{\mathbb{S}}w$ -- единственный элемент наименьшей длины в
$\{w\,W_\mathbb{S}\}$;
\\ 3.\ \ \ $l(w)=l(w_{\mathbb{S}})+l(^\mathbb{S}w)$.
\end{lemma}

\begin{proposition}\label{Stokman_list}
  (\cite{Stokman1}, стр. 480)
Если $w\in\,^\mathbb{S}W$, то сужение $*$-представлений $\Pi_w$
  на подалгебру $\mathbb{C}[L_\mathbb{S}\backslash G]_q$ неприводимо. Каждое
неприводимое $*$-представление алгебры $\mathbb{C}[L_\mathbb{S}\backslash
G]_q$ унитарно эквивалентно одному и только одному из этих сужений.
\end{proposition}

\bigskip Равенство (\ref{intSU2}) доставляет явный вид инвариантного
интеграла на $\mathbb{C}[SU_2]_q$. Приведем обобщение этого результата на
случай $*$-алгебры Хопфа $\mathbb{C}[K]_q$.


При всех $\Lambda\in P_+$, $w\in W$ весовое подпространство
$L(\Lambda)_{w\Lambda}$ одномерно. Значит, нормированный базисный вектор
$v_{w\Lambda}^\Lambda \in L(\Lambda)_{w\Lambda} $ определен с точностью до
числового множителя, по модулю равного единице,\footnote{Напомним, что
инвариантная эрмитова форма в $L(\Lambda)$, введенная в \itemiiiе
\ref{HChV}, положительно определена при всех $\lambda \in P_+$ (см.
\itemiii\ \ref{loc_finite}).} и однозначно определен элемент
$$ x_{\Lambda, w} = c_{\Lambda,\; w\Lambda}^\Lambda \cdot
(c_{\Lambda,\;w\Lambda}^\Lambda)^\star$$ алгебры $\mathbb{C}[K]_q$.

Для элемента $w\in W$ и его приведенного разложения $w=s_{i_1}s_{i_2}\ldots
s_{i_N}$ введем элементы
$$
w_j = s_{i_{j+1}}s_{i_{j+2}}\ldots s_{i_N},\quad j=0,1,\ldots,N-1;\quad
w_N=1.
$$

\begin{lemma}\label{x_w}
(\cite{ReshYak}, стр. 10). Ограниченные линейные операторы
$$(\Pi_w\otimes\chi_\varphi)(x_{\Lambda,w}),\qquad\Lambda\in P_+,$$
в $l^2(\mathbb{Z}_+)^{\otimes N}$ попарно коммутируют, и
$$
(\Pi_w\otimes\chi_t)e_{k_1}\otimes e_{k_2}\otimes\ldots\otimes e_{k_N}=
\prod\limits_{j=1}^Nq^{2(k_j+1)(w_i\Lambda,\alpha_{i_j})}e_{k_1}\otimes
e_{k_2}\otimes\ldots\otimes e_{k_N},
$$
где $\{e_{k_1}\otimes e_{k_2}\otimes\ldots\otimes e_{k_N}\}$ -- стандартный
базис гильбертова пространства $l^2(\mathbb{Z}_+)^{\otimes N}$.
\end{lemma}

\begin{corollary}\label{x_w0}
Неотрицательный линейный оператор \\
$(\Pi_{w_0}\otimes\chi_\varphi)(x_{\rho,w_0})$ является ядерным и имеет
нулевое ядро.
\end{corollary}

\begin{proposition}\label{intK}
(\cite{ReshYak}, стр. 11). Для любого элемента $f\in \mathbb{C}[K]_q$
\begin{equation}\label{RYak}
  \int\limits_{K_q} f d\nu =
  \operatorname{const}(q) \int\limits_{\mathbb{T}^l}
   \operatorname{tr}((\Pi_{w_0} \otimes \chi_\varphi)
(f\cdot x_{\rho,w_0})) d\varphi,
\end{equation}
где $d\varphi$ -- мера Хаара на группе $\mathbb{T}^l$.
\end{proposition}

\begin{note}
Явный вид положительного числового множителя $\operatorname{const}(q)$
приведен в работе \cite{ReshYak}.
\end{note}

Инвариантный интеграл допускает продолжение по непрерывности с
$\mathbb{C}[K]_q$ на $C(K)_q$.

\begin{corollary}\label{faith_nu}
Состояние $\nu$ является точным, то есть для всех ненулевых элементов $f\in
C(K)_q$
$$
\int\limits_{K_q} f^\star f d\nu > 0.
$$
\end{corollary}

{\bf Доказательство.} Все $*$-представления алгебры $\mathbb{C}[K]_q$ будем
считать продолженными по непрерывности на $C(K)_q$. Пусть $f \in C(K)_q$ и
$\int\limits_{K_q} f^\star f d\nu = 0$. Из непрерывности оператор-функции
$(\Pi_{w_0} \otimes \chi_\varphi)(f)$ на торе $\mathbb{T}^l$ и отсутствия
ядра у оператора $(\Pi_{w_0} \otimes \chi_\varphi)(x_0)\geq 0$ следует, что
$\Pi_{w_0} \otimes \chi_\varphi(f)= 0$ при всех $\varphi \in \mathbb{T}^l$.
Продолжая обе части равенства (\ref{tore}) по непрерывности на $C(K)_q$,
получаем:
$$
\|f\|_\infty = \max\limits_{\varphi\in \mathbb{T}^l} \|(\Pi_{w_0} \otimes
\chi_\varphi)(f)\| = 0,\qquad f\in C(K)_q.
$$
Значит, $f=0$. $\hfill \square$

\subsection{Квантовые векторные пространства и модули
Хариш-Чандры}\label{Vect_and_HCh}
\subsubsection{Обобщенные модули Верма.}\label{GVerma}

Так же, как в \itemiiiе \ref{DrJ}, рассмотрим матрицу Картана $(a_{ij})$,
$i,j=1,2,\ldots,l$, простой комплексной алгебры Ли $\mathfrak{g}$ и
отвечающие ей алгебры Хопфа $U\mathfrak{g}$, $U_q\mathfrak{g}$ с
образующими
$$
\{H_i,E_i,F_i\},\quad i=1,2,\ldots,l;\qquad\{K_i,K_i^{-1},E_i,F_i\},\quad
i=1,2,\ldots,l,
$$
соответственно.

Выберем подмножество $ \mathbb{S} \subset \{1,2,\ldots,l\}$ и наделим
алгебру Ли $\mathfrak{g}$ градуировкой с помощью элемента $H_\mathbb{S} \in
\mathfrak{h}$, определяемого равенством
$$
\alpha_j(H_\mathbb{S}) = \left\{
\begin{array}{ll}
  2, &\quad j  \notin \mathbb{S} \\
  0, &\quad j \in \mathbb{S}.
\end{array}\right.
$$
Именно, $$\mathfrak{g}=\newoplus\limits_j\mathfrak{g}_j,\qquad
\mathfrak{g}_j=\{\xi\in\mathfrak{g}|
\:\mathrm{ad}_{H_\mathbb{S}}\xi=2j\xi\}.$$ Очевидно,
$$\mathfrak{g}=\mathfrak{u}^-\oplus\mathfrak{l}\oplus\mathfrak{u}^+,$$ где
$\mathfrak{l}=\mathfrak{g}_0$,
$\mathfrak{u}^-=\bigoplus\limits_{j<0}\mathfrak{g}_j$,
$\mathfrak{u}^+=\bigoplus\limits_{j>0}\mathfrak{g}_j$ и
$\mathfrak{q}^\pm=\mathfrak{l}\oplus\mathfrak{u}^\pm$. Разумеется,
$\mathfrak{q}^\pm \supset \mathfrak{b}^\pm$, то есть подалгебры Ли
$\mathfrak{q}^\pm$ являются параболическими.

В следующем разделе мы наделим алгебру $\mathbb{C}[\mathfrak{u}^-]$
полиномов на $\mathfrak{u}^-$ структурой $U\mathfrak{g}$-модульной алгебры и
получим ее квантовый аналог $\mathbb{C}[\mathfrak{u}^-]_q$. Введем
необходимые определения.

Прежде всего, наделим каждый весовой $U_q\mathfrak{g}$-модуль $V$
 градуировкой
$$
V = \newoplus\limits_r V_r,\quad V_r = \{v \in V\ |\ H_\mathbb{S} v = 2rv\}.
$$

 Далее, рассмотрим подалгебры Хопфа
 $U_q\mathfrak{q}^+$, $U_q\mathfrak{q}^-$,
 $U_q\mathfrak{l}$, $U_q\mathfrak{l}_\mathrm{ss}$,
порожденные множествами
$$
\{K_i^{\pm 1}, E_i\}_{i=1,2,\ldots,l}\; \cup\; \{F_j\}_{j \in \mathbb{S}},
\qquad \{K_i^{\pm 1}, F_i\}_{i=1,2,\ldots,l}\;\cup\; \{E_j\}_{j \in
\mathbb{S}},
$$
$$
\{K_i^{\pm 1}\}_{i=1,2,\ldots,l}\;\cup\;\{E_j, \; F_j\}_{j \in \mathbb{S}},
\qquad \{K_j^{\pm 1},\; E_j,\; F_j\}_{j \in \mathbb{S}}.
$$
Отметим, что подалгебра Хопфа $U_q\mathfrak{l}_\mathrm{ss}\subset
U_q\mathfrak{l}$ является квантовым аналогом универсальной обертывающей
алгебры $U\mathfrak{l}_\mathrm{ss}$ ''полупростой части''
 $\mathfrak{l}_\mathrm{ss}=[\mathfrak{l},\mathfrak{l}]$
алгебры Ли $\mathfrak{l}$.

Очевидно, $U_q\mathfrak{l}=U_q\mathfrak{q}^+\cap U_q \mathfrak{q}^-$, и
корректно определены сюръективные гомоморфизмы алгебр Хопфа
$$
\pi^+:U_q\mathfrak{q}^+\to U_q\mathfrak{l},\qquad\pi^-:U_q\mathfrak{q}^-\to
U_q\mathfrak{l},
$$
для которых
\begin{equation}\label{pi_1}
\pi^\pm|_{U_q\mathfrak{l}}=\mathrm{id},\qquad\pi^+(E_j)=\pi^-(F_j)=0,\qquad
j\notin\mathbb{S}.
\end{equation}
Эти гомоморфизмы $\pi^\pm$ зачастую способны заменить квантовые аналоги
подалгебр Хопфа $U\mathfrak{u}^\pm\subset U\mathfrak{g}$. Например,
$q$-аналоги пространств $\mathfrak{u}^\pm$-инвариантов
$U_q\mathfrak{g}$-модуля $V$ можно ввести равенствами \IND{$q$-аналоги !
пространств $\mathfrak{u}^\pm$-инвариантов}
$$
V^{\mathfrak{u}^+}\stackrel{\operatorname{def}}{=}\{v\in
V|\:\operatorname{Ker}\pi^+\cdot v=0\},\qquad
V^{\mathfrak{u}^-}\stackrel{\operatorname{def}}{=}\{v\in
V|\:\operatorname{Ker}\pi^-\cdot v=0\}.
$$

\bigskip
Перейдем к обобщенным модулям Верма, см. \cite{Lep}. Пусть
$$
P_+^\mathbb{S}=\left\{\left.\mathbf{\lambda}=\sum\limits_{j=1}^l\lambda_j\,
\overline{\omega}_j\in\mathfrak{h}_\mathbb{R}^*\;\right|\;
\lambda_j\in\mathbb{Z}_+\text{\ при\ }j\in\mathbb{S}\right\}.
$$
Каждому $\lambda\in P_+^\mathbb{S}$ сопоставим простой весовой
$U_q\mathfrak{l}$-модуль $L(\mathfrak{l},\lambda)$ с образующей
$v(\mathfrak{l},\lambda)$ и определяющими соотношениями
$$
K_j^{\pm 1}v(\mathfrak{l},\lambda)=q_j^{\pm\lambda_j}
v(\mathfrak{l},\lambda),\qquad j=1,2,\ldots,l;
$$
$$
E_jv(\mathfrak{l},\lambda)=0,\qquad F_j^{\lambda_j+1}
v(\mathfrak{l},\lambda)=0,\qquad j\in\mathbb{S},
$$
ср. с \eqref{complete_list}. Используя гомоморфизмы
$\pi^\pm:U_q\mathfrak{q}^\pm\to U_q\mathfrak{l}$, получаем
$U_q\mathfrak{q}^\pm$-модули $L(\mathfrak{q}^\pm,\lambda)$ в том же самом
векторном пространстве.

Обобщенным модулем Верма со старшим весом $\lambda\in P_+^\mathbb{S}$
называют $U_q\mathfrak{g}$-модуль \IND{обобщенный модуль Верма}
$$
N(\mathfrak{q}^+,\lambda)=U_q\mathfrak{g}\otimes_{U_q\mathfrak{q}^+}
L(\mathfrak{q}^+,\lambda).
$$

Построим базис векторного пространства $N(\mathfrak{q}^+, \lambda)$ с
помощью базиса Пуанкаре-Биркгофа-Витта алгебры $U_q\mathfrak{n}^-$,
отвечающего специальным образом выбранному приведенному разложению элемента
$w_0 \in W$ максимальной длины.

Рассмотрим подгруппу $W_\mathbb{S}\subset W$, порожденную множеством
простых отражений $\{s_i\}_{i\in \mathbb{S}}$, и элемент наименьшей длины
$w_0^\mathbb{S}$ в классе смежности $W_\mathbb{S}\,w_0$. Как следует из
леммы \ref{w_s}, такой элемент $w_0^\mathbb{S}$ единствен (именно,
$w_0^\mathbb{S}\,=\, (^\mathbb{S}w_0)^{-1})$ и
$$w_0=w_{0,\mathbb{S}}\cdot\; w_0^\mathbb{S},$$
где $w_{0,\mathbb{S}}$-- элемент максимальной длины группы $W_\mathbb{S}$.

 Выберем приведенные разложения
\begin{equation}\label{decomp_w}
w_{0,\mathbb{S}}=s_{i_1}s_{i_2}\ldots s_{i_{M'}},\qquad
w_0^\mathbb{S}=s_{i_{M'+1}}s_{i_{M'+2}}\ldots s_{i_M}.
\end{equation}
Их конкатенация $w_0 = s_{i_1}s_{i_2}\ldots s_{i_M}$ является приведенным
разложением элемента $w_0$. Из предложения \ref{PBW-basis} следует, что
$U_q\mathfrak{g}$ является свободным правым модулем над алгеброй
$U_q\mathfrak{q}^+$ с базисом
$$
F_{\beta_M}^{j_M}F_{\beta_{M-1}}^{j_{M-1}}\ldots
F_{\beta_{M'+1}}^{j_{M'+1}},\qquad
(j_M,j_{M-1},\ldots,j_{M'+1})\in\mathbb{Z}_+^{M-M'}.
$$
Значит, имеет место

\begin{proposition}\label{N_basis}
Пусть $\mathbf{\lambda} \in P_+^\mathbb{S}$ и $\{v_1, v_2, \ldots, v_d\}$ --
базис векторного пространства $L(\mathfrak{q}^+, \lambda)$. Тогда векторы
\begin{equation}\label{basis_N}
F_{\beta_M}^{j_M} F_{\beta_{M-1}}^{j_{M-1}} \ldots
F_{\beta_{M'+1}}^{j_{M'+1}} v_i,\qquad j_k\in\mathbb{Z}_+,\; \;
i\in\{1,2,\ldots,d\},
\end{equation}
образуют базис однородных элементов градуированного векторного пространства
$N(\mathfrak{q}^+, \lambda)$.
\end{proposition}

\begin{note}\label{homogen}
Однородность элементов базиса (\ref{basis_N}) вытекает из того, что элементы
$F_\beta$ являются весовыми векторами $U_q\mathfrak{g}$-модуля
$U_q\mathfrak{g}$ с весами, равными их классическим значениям, а все векторы
$v_i$ имеют одну и ту же степень однородности.
\end{note}

Из предложений \ref{N_basis}, \ref{qWeyl} вытекает

\begin{proposition}\label{dim_N}
Однородные компоненты градуированного векторного пространства
$N(\mathfrak{q}^+, \lambda)$ являются конечномерными весовыми
$U_q\mathfrak{l}$-модулями. Кратности весов и размерности этих
$U_q\mathfrak{l}$-модулей равны их классическим значениям, то есть значениям
при $q=1$.
\end{proposition}

Отметим, что $N(\mathfrak{q}^+, \lambda)$ является локально $U_q
\mathfrak{b}^+$-конечномерным модулем, где $\mathfrak{b}^+ \subset
\mathfrak{g}$ -- стандартная борелевская подалгебра. Так же, как предложение
\ref{N_basis}, доказывается

\begin{lemma}\label{finite_gen}Однородные компоненты
$U_q\mathfrak{l}$-бимодулей $U_q\mathfrak{q}^\pm$
являются как свободными левыми, так и свободными правыми
$U_q\mathfrak{l}$-модулями конечного ранга.
\end{lemma}

Следующее свойство универсальности обобщенного модуля Верма \IND{обобщенный
модуль Верма ! свойство универсальности} $N(\mathfrak{q}^+,\lambda)$
вытекает из его определения и использует естественное вложение
$U_q\mathfrak{q}^+$-модулей
\begin{equation}\label{embed_L}
L(\mathfrak{q}^+,\lambda)\hookrightarrow U_q\mathfrak{g}
\otimes_{U_q\mathfrak{q}^+}L(\mathfrak{q}^+,\lambda),\qquad v\mapsto 1
\otimes v.
\end{equation}

\begin{lemma}\label{induction}
 Рассмотрим $U_q\mathfrak{g}$-модуль $V$. Каждый морфизм $U_q
\mathfrak{q}^+$-модулей \hbox{$f: L(\mathfrak{q}^+, \lambda) \rightarrow V$}
единственным образом продолжается до морфизма \ $U_q\mathfrak{g}$-модулей
$$N(\mathfrak{q}^+, \lambda) \rightarrow V,\qquad \xi \otimes v \mapsto \xi
f(v).$$
\end{lemma}

\bigskip Пусть $v(\mathfrak{q}^+,\lambda)=1\otimes v(\mathfrak{l},\lambda)$
-- старший вектор $N(\mathfrak{q}^+,\lambda)$. Так же, как в случае $q=1$
\cite{Lep}, доказывается

\begin{proposition}\label{Lep}
Рассмотрим модуль Верма $M(\lambda)$ со старшим весом $\lambda\in
P_+^\mathbb{S}$ и стандартный эпиморфизм $U_q\mathfrak{g}$-модулей
\begin{equation}\label{M-N}
p_\lambda:M(\lambda)\to N(\mathfrak{q}^+,\lambda),\qquad
p_\lambda:v(\lambda)\mapsto v(\mathfrak{q}^+,\lambda).
\end{equation}
Ядро $\mathrm{Ker}\,p_\lambda$ является наименьшим подмодулем
$U_q\mathfrak{g}$-модуля $M(\lambda)$, содержащим подмножество
$\left\{F_i^{\lambda_i+1}v(\lambda)\right\}_{i\in\mathbb{S}}$.
\end{proposition}

\begin{note}\label{proof_Lep} Пусть $L\subset M(\lambda)$ -- подмодуль,
порожденный множеством
$\left\{F_i^{\lambda_i+1}v(\lambda)\right\}_{i\in\mathbb{S}}$. В \cite{Lep}
равенство $L=\mathrm{Ker}\,p_\lambda$ выводится из расщепимости точной
последовательности $U\mathfrak{g}$-модулей
\begin{equation}\label{exact_Lep}
0\to\operatorname{Ker}p_\lambda/L\to M(\lambda)/L\to
M(\lambda)/\mathrm{Ker}\,p_\lambda\to 0
\end{equation}
и из того, что $U\mathfrak{g}$-модуль $M(\lambda)/L$ порождается своим
старшим вектором. В свою очередь, расщепимость последовательности
(\ref{exact_Lep}) следует из леммы \ref{induction} об универсальности
${N(\mathfrak{q}^+,\lambda)\cong M(\lambda)/\mathrm{Ker}\,p_\lambda}$ и из
коммутативности диаграммы
$$
\xymatrix{U\mathfrak{q}^+\cdot v(\lambda)\ar[r]^\approx\ar@{_(->}[d] &
U\mathfrak{q}^+\cdot v(\mathfrak{q}^+,\lambda)\ar@{_(->}[d]
\\ M(\lambda)\ar[r]^{p_\lambda} & N(\mathfrak{q}^+,\lambda)}
$$
в категории $U\mathfrak{q}^+$-модулей.
\end{note}

\medskip

\begin{corollary}\label{second_def}
Обобщенный модуль Верма $N(\mathfrak{q}^+,\lambda)$ канонически изоморфен
$U_q\mathfrak{g}$-модулю с образующей $v(\mathfrak{q}^+,\lambda)$ и
определяющими соотношениями
\begin{equation}\label{def_rel}
\begin{gathered}
E_j\,v(\mathfrak{q}^+,\lambda)=0,\quad K_j^{\pm 1}\,v(\mathfrak{q}^+,
\lambda)=q_j^{\pm\lambda_j}\,v(\mathfrak{q}^+,\lambda),\qquad
j=1,2,\ldots,l;
\\ F_j^{\lambda_j+1}\,v(\mathfrak{q}^+,\lambda)=0,\qquad j\in\mathbb{S}.
\end{gathered}
\end{equation}
\end{corollary}

\medskip

\begin{note}\label{uniqueness2}
Из определений следует существование наибольшего собственного подмодуля
$K(\mathfrak{q}^+, \lambda) \varsubsetneqq N(\mathfrak{q}^+, \lambda)$. Если
$\lambda \in P_+$, то $\dim(N(\mathfrak{q}^+, \lambda) / K(\mathfrak{q}^+,
\lambda)) < \infty$ и $\dim(N(\mathfrak{q}^+, \lambda) / K) = \infty$ для
любого собственного подмодуля $K \subsetneqq K(\mathfrak{q}^+, \lambda)$.
Последнее утверждение доказывается так же, как аналогичное утверждение в
\itemiiiе \ref{HChV}.
\end{note}

В заключение отметим, что без существенных изменений в доказательствах можно
получить аналогичные результаты для $U_q\mathfrak{g}$-модулей
$$
N(\mathfrak{q}^-,\lambda)=U_q\mathfrak{g}\otimes_{U_q \mathfrak{q}^-}
L(\mathfrak{q}^-,w_{0,\mathbb{S}}\lambda),\qquad\lambda\in -P_+^\mathbb{S},
$$
-- обобщенных модулей Верма с младшим весом $\lambda$. Здесь, как и прежде,
$w_{0,\mathbb{S}}$ -- элемент максимальной длины подгруппы $W_\mathbb{S}$
группы Вейля $W$.

\subsubsection{Квантовое  векторное пространство
$\mathfrak{u}^-$.}\label{vect}

Перейдем к построению $U_q\mathfrak{g}$-модульной алгебры
$\mathbb{C}[\mathfrak{u}^-]_q$. \IND{$U_q\mathfrak{g}$-модульная алгебра !
$\mathbb{C}[\mathfrak{u}^-]_q$} Используемый метод восходит к Дринфельду
\cite{Drinf1} и состоит в предварительном построении двойственной коалгебры.
В основе этого метода лежит следующий факт, вытекающий из определений.

\begin{proposition}\label{dual}
Рассмотрим алгебру Хопфа $A^{\mathrm{cop}}$, отличающуюся от алгебры Хопфа
$A$ заменой коумножения на противоположное. Если $C$ является
$A^{\mathrm{cop}}$-модульной коалгеброй, то сопряженное векторное
пространство $C^*$ является $A$-модульной алгеброй.
\end{proposition}

\begin{remark}\label{why_cop_?}
Следует пояснить причину замены коумножения на противоположное. Очевидно,
векторное пространство, сопряженное к коалгебре, является алгеброй, причем
тензорные сомножители при переходе к сопряженным пространствам не
переставляются. Напротив, в категории $U_q\mathfrak{g}$-модулей $
V_1^*\otimes V_2^* \hookrightarrow (V_2\otimes V_1)^*$, то есть тензорные
сомножители переставляются. Вводя противоположное коумножение, мы обходим
это несоответствие.
\end{remark}
Приведем пример $U_q\mathfrak{g}^{\mathrm{cop}}$-модульной коалгебры.

\begin{proposition}\label{coalg_N}
Отображения
$$
\Delta^+: v(\mathfrak{q}^+, 0) \mapsto v(\mathfrak{q}^+, 0) \otimes
v(\mathfrak{q}^+, 0),\quad \varepsilon^+: v(\mathfrak{q}^+, 0) \mapsto 1
$$
единственным образом продолжаются до морфизмов
$U_q\mathfrak{g}^{\mathrm{cop}}$-модулей
$$
\Delta^+: N(\mathfrak{q}^+, 0) \rightarrow N(\mathfrak{q}^+, 0) \otimes
N(\mathfrak{q}^+, 0),\quad \varepsilon^+:N(\mathfrak{q}^+, 0) \rightarrow
\mathbb{C}.
$$
Эти морфизмы наделяют $N(\mathfrak{q}^+, 0)$ структурой
$U_q\mathfrak{g}^{\mathrm{cop}}$-модульной коалгебры.
\end{proposition}

{\bf Доказательство.} Единственность морфизма $\Delta^+$ следует из того,
что $v(\mathfrak{q}^+,0)$ порождает $U_q\mathfrak{g}$-модуль
$N(\mathfrak{q}^+,0)$. Его существование вытекает из равенств
$$
E_j(v(\mathfrak{q}^+,0)\otimes v(\mathfrak{q}^+,0))=(K_j^{\pm
1}-1)(v(\mathfrak{q}^+,0)\otimes v(\mathfrak{q}^+,0))=0,\quad
j=1,2,\ldots,l,
$$
$$
F_j(v(\mathfrak{q}^+,0)\otimes v(\mathfrak{q}^+,0))=0,\quad j\in\mathbb{S},
$$
которые следуют из (\ref{comult_D}) и (\ref{def_rel}). Существование и
единственность морфизма $\varepsilon$ очевидны. Все нужные свойства морфизма
$\Delta^+$, кроме коассоциативности, вытекают из определений, а
коассоциативность достаточно проверить на образующей $v(\mathfrak{q}^+, 0)$:
$$
(v(\mathfrak{q}^+,0)\otimes v(\mathfrak{q}^+,0))\otimes
v(\mathfrak{q}^+,0)=v(\mathfrak{q}^+,0)\otimes(v(\mathfrak{q}^+,0)\otimes
v(\mathfrak{q}^+,0)). \eqno\square
$$

\medskip

Как следует из предложения \ref{dual}, алгебра
$\mathbb{C}[[\mathfrak{u}^-]]_q$ всех линейных функционалов на
$N(\mathfrak{q}^+,0)$ является $U_q\mathfrak{g}$-модульной. Это $q$-аналог
$U\mathfrak{g}$-модульной алгебры $\mathbb{C}[[\mathfrak{u}^-]]$ формальных
рядов на комплексном векторном пространстве
$\mathfrak{u}^-\cong\mathfrak{g}/\mathfrak{q}^+$. Для перехода от формальных
рядов к полиномам достаточно заменить в предыдущих построениях тензорную
категорию векторных пространств тензорной категорией градуированных
векторных пространств: прямая сумма и тензорное произведение градуированных
векторных пространств наделяются стандартными градуировками:
$$
(V'\oplus V'')_r=V'_r\oplus V''_r,\qquad(V'\otimes V'')_r= \newoplus
\limits_{r'+r''=r}(V'_{r'}\otimes V''_{r''}).
$$
Отметим, что каждому градуированному векторному пространству
\begin{equation}\label{grade}
V=\newoplus\limits_r V_r,\qquad V_r=\{v\in V|\:\deg v=r\}
\end{equation}
отвечает градуированное векторное пространство
$V^*\stackrel{\mathrm{def}}{=}\newoplus\limits_r V_r^*$. Оно совпадает с
пространством всех линейных функционалов на $V$, если и только если число
ненулевых слагаемых в (\ref{grade}) конечно.

Биалгебра $A$ называется градуированной, \IND{биалгебра ! градуированная}
если она является градуированным векторным пространством, причем умножение
$m:A\otimes A\to A$, и коумножение $\Delta:A\to A\otimes A$ уважают
градуировки, то есть являются линейными операторами нулевой степени.
Аналогично вводятся понятия градуированной алгебры и градуированного модуля,
градуированной коалгебры и градуированного комодуля, а также понятие
градуированной алгебры Хопфа. \IND{алгебра ! градуированная}\IND{модуль !
градуированный}\IND{коалгебра ! градуированная}\IND{комодуль !
градуированный}\IND{алгебра ! Хопфа ! градуированная}

В предыдущем \itemiiiе каждый весовой $U_q\mathfrak{g}$-модуль был наделен
градуировкой с помощью элемента $H_\mathbb{S}$.

Будем обобщенные модули Верма $N(\mathfrak{q}^+,\lambda)$ считать объектами
тензорной категории $U_q\mathfrak{g}^\mathrm{cop}$-модулей, а двойственные
им градуированные $U_q\mathfrak{g}$-модули -- объектами тензорной категории
$U_q\mathfrak{g}$-модулей.

Отметим, что $U_q\mathfrak{g}$, $U_q\mathfrak{g}^{\rm cop}$ --
градуированные алгебры Хопфа
$$
\deg(K_j^{\pm 1})=0,\quad\deg(E_j)=
\begin{cases}
1, & j\in\mathbb{S},
\\ 0,& j\notin\mathbb{S}
\end{cases},\quad
\deg(F_j)=
\begin{cases}
-1, & j\in\mathbb{S},
\\ 0,& j\notin\mathbb{S},
\end{cases}
$$
а $N(\mathfrak{q}^+,0)$ -- градуированная
$U_q\mathfrak{g}^{\operatorname{cop}}$-модульная коалгебра. Двойственную
градуированную $U_q\mathfrak{g}$-модульную алгебру
$$
\mathbb{C}[\mathfrak{u}^-]_q=\newoplus\limits_{j\ge 0}
\mathbb{C}[\mathfrak{u}^-]_{q,j},\qquad\mathbb{C}[\mathfrak{u}^-]_{q,j}
\stackrel{\mathrm{def}}{=}N(\mathfrak{q}^+,0)_{-j}^*
$$
назовем алгеброй полиномов на квантовом комплексном векторном пространстве
$\mathfrak{u}^-$. \IND{алгебра ! полиномов ! на квантовом комплексном
векторном пространстве $\mathfrak{u}^-$} Структура
$U_q\mathfrak{g}$-модульной алгебры вводится по двойственности:
\begin{equation}\label{dual_cop}
\xi f(v)\;=\;f(S^{-1}(\xi)v),\qquad\xi\in U_q\mathfrak{g},\;
f\in\mathbb{C}[\mathfrak{u}^-]_q,\;v\in N(\mathfrak{q}^+,0),
\end{equation}
где $S$ -- антипод алгебры Хопфа $U_q\mathfrak{g}$. Такие
$U_q\mathfrak{g}$-модульные алгебры рассматривал Кебэ \cite{Kebe_CR, Kebe}.

\begin{proposition}\label{integr}
Алгебра $\mathbb{C}[\mathfrak{u}^-]_q$ является областью целостности.
\end{proposition}

{\bf Доказательство.}
Рассмотрим естественный сюръективный гомоморфизм коалгебр
\begin{equation}\label{sur_N}
  (U_q \mathfrak{b}^-)^{\rm cop} \rightarrow N(\mathfrak{q}^+,0) \rightarrow 0,
  \qquad \xi \mapsto \xi\,v(\mathfrak{q}^+,0).
\end{equation}
Из предложений \ref{Pairing}, \ref{R-nondegenerate} следует, что сопряженный
(в категории градуированных векторных пространств) линейный оператор
является вложением алгебр $\mathbb{C}[\mathfrak{u}^-]_q \hookrightarrow U_q
\mathfrak{b}^+$. Остается воспользоваться следствием \ref{nodivisors},
согласно которому алгебра $U_q\mathfrak{g}$ является областью
целостности.\hfill $\square$

\medskip Из предложения \ref{dim_N} следует, что однородные компоненты
$\mathbb{C}[\mathfrak{u}^-]_{q,j}$ градуированного векторного пространства
$\mathbb{C}[\mathfrak{u}^-]_q$ являются $U_q\mathfrak{l}$-модулями той же
размерности, что в классическом случае $q=1$:
\begin{equation}\label{dim_j}
  \dim\,\mathbb{C}[\mathfrak{u}^-]_{q,j} =
  \dbinom{j + \dim\,\mathfrak{u}^- - 1}{\dim\,\mathfrak{u}^- - 1} =
  \dbinom{j + \dim\,\mathfrak{u}^- - 1}{j}.
\end{equation}

\bigskip

Рассмотрим простейший пример $\mathfrak{g}=\mathfrak{sl}_2$. В этом случае
$N(\mathfrak{q}^+,0)$ --- модуль Верма $M(0)$ с нулевым старшим весом, и
множество $\{F^j v(0)\}_{j\in\mathbb{Z}_+}$ является базисом векторного
пространства $M(0)$.

Очевидно, $M(0)=\newoplus\limits_{j\ge 0} M(0)_{-j}$,
$M(0)_{-j}=\mathbb{C}F^jv(0)$. Единицей $U_q\mathfrak{sl}_2$-модульной
алгебры $\mathbb{C}[\mathfrak{u}^-]_q$ является линейный функционал, равный
единице на $v(0)$ и нулю на $\newoplus\limits_{j\ne 0}\,M(0)_j$. Рассмотрим
элемент $z\in M(0)^*$, определяемый равенством
$$
z(S(F^j)v(0))=
\left\{
\begin{array}{rl}
q^{1/2}, & \;\quad j=1
\\ 0, &\;\quad j\ne 1,
\end{array}
\right.
$$
где $S$ -- антипод алгебры Хопфа $U_q \mathfrak{sl}_2$. Индукцией по $j$
нетрудно доказать, что $z^j(F^jv(0))\ne 0$. Значит, $z^j$ -- ненулевой
элемент одномерного пространства $M(0)_{-j}^*\subset M(0)^*$, и алгебра
$\mathbb{C}[\mathfrak{u}^-]_q$ изоморфна алгебре $\mathbb{C}[z]$ полиномов
одной переменной. Структура $U_q\mathfrak{sl}_2$-модульной алгебры в
$\mathbb{C}[\mathfrak{u}^-]_q$ полностью определяется действием образующих
$K^{\pm 1}$, $E$, $F$ алгебры $U_q\mathfrak{sl}_2$ на образующую $z$.
Используя \eqref{dual_cop} и явный вид антипода $S$, нетрудно получить
равенства
\begin{equation}\label{U_q_z}
K^{\pm 1}z=q^{\pm 2}z,\qquad Fz=q^{1/2},\qquad Ez=-q^{1/2}z^2
\end{equation}
(множитель $-q^{1/2}$ в последнем из них найден с помощью соотношения
$(EF-FE)z=\frac{K-K^{-1}}{q-q^{-1}}z$).

\bigskip

Рассмотренный выше простейший пример допускает обобщение в нескольких
направлениях. Случай минимальной параболической подалгебры $\mathbb{S} =
\varnothing$, $\mathfrak{q}^+ = \mathfrak{b}^+$ изучался А.~Джозефом и его
сотрудниками \cite{Jos}. В противоположном случае максимальной
параболической подалгебры множество $\mathbb{S}$ содержит все элементы
$1,2,\ldots,l$, кроме одного, обозначаемого в дальнейшем $l_0$. Пусть
$\alpha_1, \alpha_2, \ldots, \alpha_l$ -- простые корни алгебры Ли
$\mathfrak{g}$ (см. раздел \ref{DrJ}) и $\delta = \sum\limits_{j=1} n_j
\alpha_j$ -- ее старший корень. Наиболее интересны случаи $n_{l_0}=1$ и
$n_{l_0}=2$. Первый из них приводит к $q$-аналогам предоднородных векторных
пространств коммутативного параболического типа (см. \itemiii\
\ref{prehomogeneous}), а второй случай связан с дискретной серией
Бореля-Зибенталя \cite{EPWW,Fraj1} представлений вещественных редуктивных
групп.

\bigskip

 \begin{note}\label{plus}
Обобщенный модуль Верма $N(\mathfrak{q}^-,0) = \bigoplus\limits_{j\geq 0}
N(\mathfrak{q}^-,0)_j$ наделяется структурой градуированной
$U_q\mathfrak{g}$-модульной коалгебры точно так же, как
$N(\mathfrak{q}^-,0)$: $\Delta v(\mathfrak{q}^-,0)=v(\mathfrak{q}^-,0)
\otimes v(\mathfrak{q}^-,0)$. Двойственное градуированное векторное
пространство
$$
\mathbb{C}[\mathfrak{u}^+]_q = \newoplus\limits_{j\geq 0}
\mathbb{C}[\mathfrak{u}^+]_{q,-j},\qquad \mathbb{C}[\mathfrak{u}^+]_{q,-j} =
N(\mathfrak{q}^-,0)_j^*
$$
является $U_q\mathfrak{g}$-модульной алгеброй ''полиномов на квантовом
комплексном векторном пространстве'' $\mathfrak{u}^+$. Как будет показано в
\itemiiiе \ref{Pol}, алгебры $\mathbb{C}[\mathfrak{u}^+]_q$ и
$\mathbb{C}[\mathfrak{u}^-]_q$ антиизоморфны. Значит,
$\mathbb{C}[\mathfrak{u}^+]_q$ -- область целостности.
\end{note}

\subsubsection{$(U_q\mathfrak{g},*)$-модульная $*$-алгебра
$\mathrm{Pol}(\mathfrak{u}^-)_q$.}\label{braid_algebra}\label{Pol}

В \itemiiiе \ref{Pol_C_q} использование универсальной $R$-матрицы позволило
получить $U_q\mathfrak{su}_{1,1}$-модульную алгебру ${\rm
Pol}(\mathbb{C})_q$. Покажем, что те же построения позволяют получить
существенно более общий результат.

Рассмотрим комплексное векторное пространство
$\mathfrak{u}^-\oplus\mathfrak{u}^+$ и алгебру
$\mathbb{C}[\mathfrak{u}^-\oplus\mathfrak{u}^+]$ полиномов на
$\mathfrak{u}^-\oplus\mathfrak{u}^+$. Получим ее $q$-аналог, наделив
структурой $U_q\mathfrak{g}$-модульной алгебры $U_q\mathfrak{g}$-модуль
\IND{$q$-аналог ! алгебры полиномов на $\mathfrak{u}^-\oplus\mathfrak{u}^+$}
$$
\mathbb{C}[\mathfrak{u}^-\oplus\mathfrak{u}^+]_q
\stackrel{\operatorname{def}}{=}\mathbb{C}[\mathfrak{u}^-]_q\otimes
\mathbb{C}[\mathfrak{u}^+]_q.
$$

Воспользуемся тем, что $U_q\mathfrak{g}$-модули
$\mathbb{C}[\mathfrak{u}^\pm]_q$ являются весовыми и
$U_q\mathfrak{g}$-модуль $\mathbb{C}[\mathfrak{u}^-]_q$ локально
$U_q\mathfrak{b}^-$-конечномерен, а $U_q\mathfrak{g}$-модуль
$\mathbb{C}[\mathfrak{u}^+]_q$ --- локально
$U_q\mathfrak{b}^+$-конечномерен. Как объясняется в \itemiiiе \ref{RM}, это
позволяет корректно определить морфизм $U_q\mathfrak{g}$-модулей
$$
\check{R}=\check{R}_{\mathbb{C}[\mathfrak{u}^+]_q
\mathbb{C}[\mathfrak{u}^-]_q}:\mathbb{C}[\mathfrak{u}^+]_q\otimes
\mathbb{C}[\mathfrak{u}^-]_q\to\mathbb{C}[\mathfrak{u}^-]_q\otimes
\mathbb{C}[\mathfrak{u}^+]_q.
$$
Наделим $\mathbb{C}[\mathfrak{u}^-\oplus\mathfrak{u}^+]_q$ структурой
алгебры, задав умножение
$$
{m:\mathbb{C}[\mathfrak{u}^-\oplus\mathfrak{u}^+]_q\otimes
\mathbb{C}[\mathfrak{u}^-\oplus\mathfrak{u}^+]_q\to
\mathbb{C}[\mathfrak{u}^-\oplus\mathfrak{u}^+]_q}
$$
с помощью умножений
$$
m^\pm:\mathbb{C}[\mathfrak{u}^\pm]_q\otimes\mathbb{C}[\mathfrak{u}^\pm]_q\to
\mathbb{C}[\mathfrak{u}^\pm]_q,\qquad m^\pm:f_1\otimes f_2\mapsto f_1 f_2.
$$
в алгебрах $\mathbb{C}[\mathfrak{u}^\pm]_q$. Именно, положим
\begin{equation}\label{mult_pm}
m=(m^-\otimes m^+)(\mathrm{id}_{\mathbb{C}[\mathfrak{u}^-]_q}\otimes
\check{R}\otimes\mathrm{id}_{\mathbb{C}[\mathfrak{u}^+]_q}),
\end{equation}
см. аналогичное определение в \itemiiiе \ref{Pol_C_q}.

Из предложения \ref{braiding1} вытекают следующие равенства морфизмов
\hbox{$U_q\mathfrak{g}$-модулей}:
$$
(m^- \otimes \mathrm{id}) \check{R}_{\mathbb{C}[\mathfrak{u}^+]_q,\,
\mathbb{C}[\mathfrak{u}^-]_q \otimes \mathbb{C}[\mathfrak{u}^-]_q} =
\check{R}_{\mathbb{C}[\mathfrak{u}^+]_q \mathbb{C}[\mathfrak{u}^-]_q}
(\mathrm{id} \otimes m^-),
$$
$$
(\mathrm{id} \otimes m^+) \check{R}_{\mathbb{C}[\mathfrak{u}^+]_q \otimes
\mathbb{C}[\mathfrak{u}^+]_q,\, \mathbb{C}[\mathfrak{u}^-]_q} =
\check{R}_{\mathbb{C}[\mathfrak{u}^+]_q \mathbb{C}[\mathfrak{u}^-]_q} (m^+
\otimes \mathrm{id}),
$$
а из предложения \ref{braiding2} -- равенства
$$
\check{R}_{\mathbb{C}[\mathfrak{u}^+]_q,\, \mathbb{C}[\mathfrak{u}^-]_q
\otimes \mathbb{C}[\mathfrak{u}^-]_q} =
(\mathrm{id}_{\mathbb{C}[\mathfrak{u}^-]_q} \otimes
\check{R}_{\mathbb{C}[\mathfrak{u}^+]_q \mathbb{C}[\mathfrak{u}^-]_q}) \cdot
(\check{R}_{\mathbb{C}[\mathfrak{u}^+]_q \mathbb{C}[\mathfrak{u}^-]_q}
\otimes \mathrm{id}_{\mathbb{C}[\mathfrak{u}^-]_q}),
$$
$$
\check{R}_{\mathbb{C}[\mathfrak{u}^+]_q \otimes
\mathbb{C}[\mathfrak{u}^+]_q,\, \mathbb{C}[\mathfrak{u}^-]_q} =
(\check{R}_{\mathbb{C}[\mathfrak{u}^+]_q \mathbb{C}[\mathfrak{u}^-]_q}
\otimes \mathrm{id}_{\mathbb{C}[\mathfrak{u}^+]_q}) \cdot
(\mathrm{id}_{\mathbb{C}[\mathfrak{u}^+]_q} \otimes
\check{R}_{\mathbb{C}[\mathfrak{u}^+]_q \mathbb{C}[\mathfrak{u}^-]_q}).
$$
C их помощью, используя ассоциативность умножений $m^\pm$ в алгебрах
$\mathbb{C}[\mathfrak{u}^\pm]_q$, нетрудно доказать ассоциативность
умножения $m$ в $\mathbb{C}[\mathfrak{u}^- \oplus \mathfrak{u}^+]_q$.
Остальные утверждения следующего предложения вытекают из предложения
\ref{braiding1}.

\begin{proposition}\label{assoc}
Умножение $m$ наделяет $U_q\mathfrak{g}$-модуль
$\mathbb{C}[\mathfrak{u}^-\oplus \mathfrak{u}^+]_q$
 структурой $U_q\mathfrak{g}$-модульной алгебры с единицей $1 \otimes 1$.
\end{proposition}

\begin{note}\label{subalgebras}
Как следует из предложения \ref{braiding2}, вложения векторных пространств
$$
\mathbb{C}[\mathfrak{u}^-]_q \hookrightarrow \mathbb{C}[\mathfrak{u}^-
\oplus \mathfrak{u}^+]_q,\quad f\mapsto f \otimes 1,
$$
$$
\mathbb{C}[\mathfrak{u}^+]_q \hookrightarrow \mathbb{C}[\mathfrak{u}^-
\oplus \mathfrak{u}^+]_q,\quad f \mapsto 1 \otimes f,
$$
являются гомоморфизмами $U_q\mathfrak{g}$-модульных алгебр.
\end{note}


Рассмотрим градуированнную алгебру Ли $\mathfrak{g} =
\newoplus\limits_i \mathfrak{g}_i$,  введенную в \itemiiiе \ref{GVerma},
и ее инволютивный автоморфизм $\theta$:
$$
\theta|_{\mathfrak{g}_i} = \left\{
\begin{array}{rl}
   1, &\;\quad i \in 2\mathbb{Z} \\
  -1, &\;\quad i \in 2\mathbb{Z}+1.
\end{array}
\right.
$$
Рассмотрим продолжение этого автоморфизма до автоморфизма универсальной
обертывающей алгебры $U\mathfrak{g}$. Его $q$-аналогом служит автоморфизм
$\theta_q$ алгебры Хопфа $U_q\mathfrak{g}$, определяемый равенствами
$$
\theta_q(K_i^{\pm1})=K_i^{\pm1},\qquad\theta_q(E_i)=\left\{
\begin{array}{rl}
E_i, & \;i\in\mathbb{S}
\\ -E_i, &\; i\notin\mathbb{S}
\end{array}
\right.,\qquad\theta_q(F_i)=\left\{
\begin{array}{rl}
F_i, &\; i\in\mathbb{S},
\\ -F_i, &\; i\notin\mathbb{S}.
\end{array}
\right.
$$

В \itemiiе \ref{func_alg} при изучении компактных квантовых групп
существенно использовалась $*$-алгебра Хопфа $(U_q\mathfrak{g}, \star)$,
где $\star$ -- инволютивный антилинейный антиавтоморфизм (\ref{aster1}).
Такую же роль в теории некомпактных квантовых групп играет $*$-алгебра
Хопфа $(U_q\mathfrak{g},
*)$, где
\begin{equation}\label{star_1}
  * = \theta_q \cdot \star = \star \cdot \theta_q.
  \end{equation}
То есть,
$$
(K_j^{\pm1})^* = K_j^{\pm1},\qquad j=1,2,\ldots,l,
$$
$$
E_j^* = \left\{
\begin{array}{rl}
   K_j F_j, & \; j \in \mathbb{S} \\
  -K_j F_j, & \; j \notin \mathbb{S}
\end{array}
\right.,\qquad F_j^* = \left\{
\begin{array}{rl}
   E_j K_j^{-1}, &\; j \in \mathbb{S}, \\
  -E_j K_j^{-1}, &\; j \notin \mathbb{S}.
\end{array}
\right.
$$
Отметим, что инволюции такого типа впервые рассматривались в
\cite[стр.~188]{RTF}.

Наделим $\mathbb{C}[\mathfrak{u}^- \oplus \mathfrak{u}^+]_q$ структурой
$(U_q\mathfrak{g}, *)$-модульной алгебры, используя естественные вложения
$$\mathbb{C}[\mathfrak{u}^\pm]_q \hookrightarrow
 \mathbb{C}[\mathfrak{u}^-]_q \otimes \mathbb{C}[\mathfrak{u}^+]_q,
  \qquad f_- \mapsto f_- \otimes 1,\quad f_+ \mapsto 1 \otimes f_+.$$
\begin{proposition}\label{star_Pol}
Существует единственный такой инволютивный антилинейный антиавтоморфизм $*$
алгебры $\mathbb{C}[\mathfrak{u}^- \oplus \mathfrak{u}^+]_q$, что,
во-первых, $*: \mathbb{C}[\mathfrak{u}^\pm]_q \rightarrow
\mathbb{C}[\mathfrak{u}^\mp]_q$ и, во-вторых,
\begin{equation}\label{*_*_new}
  (\xi f)^* = (S(\xi))^* f^*, \qquad \quad \qquad
  \xi \in U_q\mathfrak{g},\; f \in \mathbb{C}[\mathfrak{u}^- \oplus
\mathfrak{u}^+]_q.
\end{equation}
\end{proposition}

{\bf Доказательство.} Начнем с единственности. Рассмотрим инволюцию с
описанными в формулировке предложения свойствами и ее сужения на подалгебры
$\mathbb{C}[\mathfrak{u}^\pm]_q$. Из (\ref{*_*_new}) и из определения
градуировок в $\mathbb{C}[\mathfrak{u}^\pm]_{q,j}$ следует, что
$$
*:\mathbb{C}[\mathfrak{u}^\pm]_{q,j}\to\mathbb{C}[\mathfrak{u}^\mp]_{q,-j}.
$$
Это позволяет ввести в рассмотрение сопряженные антилинейные операторы
${*:N(\mathfrak{q}^\pm,0)\to N(\mathfrak{q}^\mp,0)}$, определяемые
равенством
\begin{equation}\label{conjugate}
f^*(v)=\overline{f(v^*)},\qquad f\in\mathbb{C}[\mathfrak{u}^\pm]_q,\;\,
v\in N(\mathfrak{q}^\pm,0).
\end{equation}
Отметим, что $*:N(\mathfrak{q}^\pm,0)_j\to N(\mathfrak{q}^\mp,0)_{-j}$. Так
как $1^*=1$ в $\mathbb{C}[\mathfrak{u}^-\oplus\mathfrak{u}^+]_q$, то
\begin{equation}\label{v_v}
*:v(\mathfrak{q}^\pm,0)\mapsto v(\mathfrak{q}^\mp, 0).
\end{equation}
Используя (\ref{*_*_new})--(\ref{v_v}) и $S*S*=\operatorname{id}$, получаем
цепочку равенств
$$ (\xi
f)^*(v(\mathfrak{q}^\pm,0))=((S(\xi))^*f^*)(v(\mathfrak{q}^\pm,0)),
$$
$$
\overline{(\xi
f)(v(\mathfrak{q}^\mp,0))}=f^*(S((S(\xi))^*)v(\mathfrak{q}^\pm,0)),
$$
$$
\overline{(\xi f)(v(\mathfrak{q}^\mp,0))}=f^*(\xi^* v(\mathfrak{q}^\pm,0)),
$$
$$
\overline{f(S(\xi)v(\mathfrak{q}^\mp,0))}=\overline{f((\xi^*v(\mathfrak{q}^\pm,0))^*)},
$$
$$
f((S^{-1}(\xi))^*v(\mathfrak{q}^\mp,0))=f((\xi v(\mathfrak{q}^\pm))^*).
$$
Следовательно,
\begin{equation}\label{def_star1}
*:\xi v(\mathfrak{q}^\pm,0)\mapsto(S^{-1}(\xi))^*v(\mathfrak{q}^\mp,0),
\qquad\xi\in U_q\mathfrak{g}.
\end{equation}
Инволюции
$*:\mathbb{C}[\mathfrak{u}^\pm]_q\to\mathbb{C}[\mathfrak{u}^\mp]_q$,\ \
${*:N(\mathfrak{q}^\pm,0)\to N(\mathfrak{q}^\mp,0)}$, однозначно
определяются равенствами (\ref{conjugate}), (\ref{def_star1}). Значит,
инволюция с предписанными свойствами в
$\mathbb{C}[\mathfrak{u}^-\oplus\mathfrak{u}^+]_q$ единственна.

Ее построение будет осуществлено с помощью тех же рассуждений, проведенных
в обратном порядке. Равенство (\ref{def_star1}) позволяет определить
антилинейные операторы $*: N(\mathfrak{q}^\pm, 0) \rightarrow
N(\mathfrak{q}^\mp, 0)$ (корректность следует из определения
$U_q\mathfrak{g}$-модулей $N(\mathfrak{q}^\pm, 0)$). Равенство
(\ref{conjugate}) определяет антилинейные операторы
 $*:\mathbb{C}[\mathfrak{u}^\pm]_q \rightarrow
\mathbb{C}[\mathfrak{u}^\mp]_q$.
 Остается ввести антилинейный оператор в
${\mathbb{C}[\mathfrak{u}^- \oplus \mathfrak{u}^+]_q =
\mathbb{C}[\mathfrak{u}^-]_q \otimes \mathbb{C}[\mathfrak{u}^+]_q}$,
полагая
$$
(f_- \otimes f_+)^* = f_+^* \otimes f_-^*,\qquad f_{\pm} \in
\mathbb{C}[\mathfrak{u}^\pm]_q,
$$
и доказать, что он наделяет $\mathbb{C}[\mathfrak{u}^- \oplus
\mathfrak{u}^+]_q$ структурой $(U_q\mathfrak{g}, *)$-модульной алгебры.

Во-первых, антилинейные операторы $*:N(\mathfrak{q}^\pm,0)\to
N(\mathfrak{q}^\mp, 0)$ являются антигомоморфизмами коалгебр потому, что
отображение $\xi\mapsto(S^{-1}(\xi))^*$ -- антиавтоморфизм коалгебры
$U_q\mathfrak{g}$. Во-вторых, эти антилинейные операторы являются взаимно
обратными потому, что $*S^{-1}*S^{-1}=\mathrm{id}$. Значит, антилинейные
операторы
$*:\mathbb{C}[\mathfrak{u}^\pm]_q\to\mathbb{C}[\mathfrak{u}^\mp]_q$ являются
взаимно обратными антилинейными антигомоморфизмами алгебр. В-третьих, как
следует из \eqref{S_theta},
\begin{equation}\label{new_R}
S\otimes S(R)=R,\qquad R^{*\otimes *}=R^{21},
\end{equation}
где $R^{21}$ отличается от $R$ лишь перестановкой тензорных сомножителей.
Используя \eqref{new_R}, нетрудно доказать, что
$$
(f_1f_2)^*=f_2^*f_1^*,\qquad f_1,f_2\in\mathbb{C}[\mathfrak{u}^-\oplus
\mathfrak{u}^+]_q,
$$
см. раздел 16 статьи \cite{SV}. Значит, $*$ является инволюцией
$U_q\mathfrak{g}$-модульной алгебры
$\mathbb{C}[\mathfrak{u}^-\oplus\mathfrak{u}^+]_q$. Остается показать, что
инволюции в $U_q\mathfrak{g}$ и в
$\mathbb{C}[\mathfrak{u}^-\oplus\mathfrak{u}^+]_q$ согласованы
$$
(\xi f)^*=(S(\xi))^*f^*,\qquad\xi\in U_q\mathfrak{g},\;\,f\in
\mathbb{C}[\mathfrak{u}^-\oplus\mathfrak{u}^+]_q.
$$
Начнем с частного случая $f\in\mathbb{C}[\mathfrak{u}^\pm]_q$. Используя
равенства (\ref{conjugate}), (\ref{def_star1}) и \hbox{$S*S=*$},
 пройдем
предыдущую цепочку равенств в обратном направлении:
$$
f((\xi v)^*)=f((S^{-1}(\xi))^*v^*),\qquad\overline{f(S(\xi)v^*)}=
\overline{f((\xi^*v)^*)},
$$
$$
\overline{(\xi f)(v^*)}=f^*(\xi^*v),\qquad\overline{(\xi f)(v^*)}=
f^*(S((S(\xi))^*)v),
$$
$$
(\xi f)^*(v) = ((S(\xi))^* f^*)(v),
$$
где $v\in N(\mathfrak{q}^\pm,0)$, $f\in N(\mathfrak{q}^\pm,0)^*$. Общий
случай сводится к рассмотренному частному случаю. Действительно, пусть
$f=f_-f_+$, где $f_{\pm}\in\mathbb{C}[\mathfrak{u}^\pm]_q$. Тогда
$(\xi(f_-f_+))^*=\sum_i(\xi''_i f_+)^*(\xi'_if_-)^*$, где
$\Delta(\xi)=\sum_i\xi'_i\otimes\xi''_i$, и
$$
(S(\xi))^*(f_-f_+)^*=(S(\xi))^*(f_+^*f_-^*)=\sum_i((S(\xi''_i))^*f_+^*)
(S(\xi'_i)^* f_-^*). \eqno\square
$$

\medskip

Наличие антилинейного антиизоморфизма $*: \mathbb{C}[\mathfrak{u}^-]_q
\rightarrow \mathbb{C}[\mathfrak{u}^+]_q$ позволяет считать
$U_q\mathfrak{g}$-модульную алгебру $\mathbb{C}[\mathfrak{u}^+]_q$
квантовым аналогом ($q$-аналогом) алгебры антиголоморфных полиномов на
$\mathfrak{u}^-$.

Введем обозначение $\mathrm{Pol}(\mathfrak{u}^-)_q =
(\mathbb{C}[\mathfrak{u}^- \oplus \mathfrak{u}^+]_q, *)$ для полученной
$(U_q\mathfrak{g}, *)$-модульной алгебры. Пример квантового круга
$\mathfrak{g}=\mathfrak{sl}_2$ показывает, что
$\mathrm{Pol}(\mathfrak{u}^-)_q$ является $q$-аналогом алгебры полиномов на
овеществленном векторном пространстве $\mathfrak{u}^-$.

\bigskip
 Векторное пространство
$\mathrm{Pol}(\mathfrak{u}^-)_q$ допускает разложение
\begin{equation}\label{expan1}
  \mathrm{Pol}(\mathfrak{u}^-)_q = \newoplus\limits_{i,j=0}^\infty
\mathbb{C}[\mathfrak{u}^-]_{q,i} \cdot \mathbb{C}[\mathfrak{u}^+]_{q,-j},
\end{equation}
причем $\mathbb{C}[\mathfrak{u}^-]_{q,0} \cdot
\mathbb{C}[\mathfrak{u}^+]_{q,0} = \mathbb{C} \cdot 1$. Другими словами,
для любого элемента $f \in \mathrm{Pol}(\mathfrak{u}^-)_q$ существует и
единственно разложение
\begin{equation}\label{expan2}
  f = \langle f \rangle \cdot 1 +
   \sum\limits_{\{ (i,j) \in \mathbb{Z}_+^2|(i,j)\neq (0,0)\}} f_{i,j},\qquad
  f_{i,j} \in \mathbb{C}[\mathfrak{u}^-]_{q,i} \cdot
\mathbb{C}[\mathfrak{u}^+]_{q,-j},
\end{equation}
где $\langle \cdot \rangle$ -- определяемый последним равенством линейный
функционал на $\mathrm{Pol}(\mathfrak{u}^-)_q$.

\begin{proposition}\label{inner1}
Полуторалинейная форма $(f_1,f_2) \stackrel{\operatorname{def}}{=} \langle
f_2^* f_1 \rangle$
 в $\mathbb{C}[\mathfrak{u}^-]_q$ невырождена и
\begin{equation}\label{Re}
  (f_1, f_2) = \overline{(f_2, f_1)},\qquad
  f_1,f_2 \in \mathbb{C}[\mathfrak{u}^-]_q.
\end{equation}
\end{proposition}

{\bf Доказательство.} Прежде всего, докажем равенство
\begin{equation}\label{inner_new}
\langle f_2^*f_1\rangle=f_2^*\otimes f_1(R(v(\mathfrak{q}^-,0)\otimes
v(\mathfrak{q}^+,0))),\qquad f_1,f_2\in\mathbb{C}[\mathfrak{u}^-]_q,
\end{equation}
где $R$ -- универсальная $R$-матрица (\ref{Rmatrix}). Напомним, что
$$
N(\mathfrak{q}^-,0)\otimes N(\mathfrak{q}^+,0)=
\newoplus\limits_{i,j=0}^\infty N(\mathfrak{q}^-,0)_i\otimes
N(\mathfrak{q}^+,0)_{-j}.
$$
Наделим конечномерные векторные пространства $N(\mathfrak{q}^-,0)_i\otimes
N(\mathfrak{q}^+,0)_{-j}$ обычной топологией, а $N(\mathfrak{q}^-,0)\otimes
N(\mathfrak{q}^+,0)$ -- слабейшей из топологий, в которых непрерывны все
''проекторы'' на $N(\mathfrak{q}^-,0)_i\otimes N(\mathfrak{q}^+,0)_{-j}$
параллельно $\newoplus\limits_{i'\ne i,j'\ne
j}\;N(\mathfrak{q}^-,0)_{i'}\otimes N(\mathfrak{q}^+,0)_{-j'} $.
Пространство формальных рядов вида
$$
\sum\limits_{i,j} f_{i,j},\qquad f_{i,j}\in{N(\mathfrak{q}^-,0)_i\otimes
N(\mathfrak{q}^+,0)_{-j}}
$$
%
является пополнением введенного выше топологического векторного пространства
$N(\mathfrak{q}^-,0)\otimes N(\mathfrak{q}^+,0)$ и обозначается
$N(\mathfrak{q}^-,0)\widehat{\otimes} N(\mathfrak{q}^+,0)$. Спаривание
$(\mathbb{C}[\mathfrak{u}^+]_q\otimes\mathbb{C}[\mathfrak{u}^-]_q)\times
(N(\mathfrak{q}^-,0)\otimes N(\mathfrak{q}^+,0))\to\mathbb{C}$ продолжается
по непрерывности до спаривания
$(\mathbb{C}[\mathfrak{u}^+]_q\otimes\mathbb{C}[\mathfrak{u}^-]_q)\times
(N(\mathfrak{q}^-,0)\widehat{\otimes}N(\mathfrak{q}^+,0))\to\mathbb{C}$.
Если
$$
\varphi=\sum\limits_i\varphi_i,\quad\varphi_i\in
\mathbb{C}[\mathfrak{u}^-]_{q,i},\qquad\quad\psi=\sum\limits_i\psi_j,\quad
\psi_j\in\mathbb{C}[\mathfrak{u}^+]_{q,-j},
$$
то $\varphi_0=\varphi(v(\mathfrak{q}^+,0))$,
$\psi_0=\psi(v(\mathfrak{q}^-,0))$. Значит, из определения умножения в
алгебре $\mathbb{C}[\mathfrak{u}^-\oplus\mathfrak{u}^+]_q$ и из равенства
$S\otimes S(R)=R$ вытекает \eqref{inner_new}.

Завершим доказательство предложения. Действие $U_q\mathfrak{g} \otimes
U_q\mathfrak{g}$ в $N(\mathfrak{q}^-,0) \otimes N(\mathfrak{q}^+,0)$ и
антилинейный оператор $v_1 \otimes v_2 \mapsto v_2^* \otimes v_1^*$ в этом
пространстве допускают продолжение по непрерывности на $N(\mathfrak{q}^-,0)
\widehat{\otimes} N(\mathfrak{q}^+,0)$.

Как следует из \eqref{inner_new}, равенство (\ref{Re}) вытекает из
$$
 R^{21}\, (v(\mathfrak{q}^+,0)\otimes v(\mathfrak{q}^-,0))=
 \left(R\, (v(\mathfrak{q}^-,0)\otimes v(\mathfrak{q}^+,0))\right)^{*\otimes
*}.
$$
В свою очередь, последнее равенство следует из (\ref{conjugate}),
(\ref{def_star1}) и (\ref{new_R}).

Невырожденность полуторалинейной формы $(\cdot,\cdot)$ в
$\mathbb{C}[\mathfrak{u}^-]_q$ вытекает из равенства (\ref{inner_new}),
мультипликативной формулы для инвариантной $R$-матрицы (\ref{Rmatrix}) и из
явного вида базисов в $N(\mathfrak{q}^\pm, 0)$, полученных в \itemiiiе
\ref{GVerma}. \hfill $\square$

\subsubsection{Категория $C(\mathfrak{g},\mathfrak{l})_q$ и
эрмитово-симметрические пары.}\label{Hermit}

Рассмотрим полную подкатегорию $C(\mathfrak{g},\mathfrak{l})_q$ категории
$U_q\mathfrak{g}$-модулей, образованную $P$-весовыми локально
$U_q\mathfrak{l}$-конечномерными модулями. \IND{категория !
$C(\mathfrak{g},\mathfrak{l})_q$} Она замкнута относительно операций
$\oplus$, $\otimes$, перехода к подмодулям и фактор-модулям. Очевидно,
$\operatorname{Pol}(\mathfrak{u}^-)_q\in C(\mathfrak{g},\mathfrak{l})_q$.

В \itemiiiе \ref{GVerma} каждый весовой $U_q\mathfrak{g}$-модуль был наделен
градуировкой с помощью элемента $H_\mathbb{S}\in\mathfrak{h}$. Из
определения этой градуировки следует, что алгебры Хопфа
$U_q\mathfrak{q}^\pm$ являются градуированными:
$$
U_q\mathfrak{q}^+=\newoplus\limits_{j=0}^\infty(U_q\mathfrak{q}^+)_j,\qquad
U_q\mathfrak{q}^-=\newoplus\limits_{j=0}^\infty(U_q\mathfrak{q}^-)_{-j},
$$
причем $(U_q\mathfrak{q}^+)_0=(U_q\mathfrak{q}^-)_0=U_q\mathfrak{l}$.

\begin{proposition}\label{gamma}
Каждый $U_q\mathfrak{g}$-модуль $V$ обладает наибольшим подмодулем
$V_\mathfrak{l}$ категории $C(\mathfrak{g},\mathfrak{l})_q$.
\end{proposition}

{\bf Доказательство.} Очевидно, каждый $U_q\mathfrak{g}$-модуль $V$
обладает наибольшим $P$-весовым подмодулем $V_\mathfrak{h}$. Остается
доказать, что подпространство $V'=\{v\in
V_\mathfrak{h}|\:\dim(U_q\mathfrak{l}\,v)<\infty\}$ является подмодулем
$U_q\mathfrak{g}$-модуля $V_\mathfrak{h}$. Как следует из леммы
\ref{finite_gen}, для любого $j$ существует такое конечное подмножество
$\{\eta_1^\pm,\eta_2^\pm,\cdots\}$, что
$(U_q\mathfrak{q}^\pm)_j=\sum\limits_{k=1}^{N(j)}\eta_k^\pm
U_q\mathfrak{l}$. Значит,
$$
\dim(U_q\mathfrak{l}\;\xi v)\le\dim\left((U_q\mathfrak{q}^\pm)_j\,v\right)=
\dim\left(\sum\limits_{k=1}^{N(j)}\eta_k^\pm U_q\mathfrak{l}\,v\right)<
\infty
$$
при всех $\xi\in(U_q\mathfrak{q}^\pm)_j$, $v\in V'$; то есть $\xi v\in V'$.
\hfill $\square$

\bigskip

Получаем функтор $\Gamma$ из категории всех $U_q\mathfrak{g}$-модулей в
категорию $C(\mathfrak{g},\mathfrak{l})_q$, сопоставляющий каждому
$U_q\mathfrak{g}$-модулю $V$ его подмодуль $V_{\mathfrak{l}}$, а каждому
морфизму $V\rightarrow U$ --- его сужение на $V_\mathfrak{l}$. В теории
представлений вещественных редуктивных групп функторы такого типа называют
функторами Цакермана \cite{KV}. \IND{функтор Цакермана}

\medskip

Пусть $V\in C(\mathfrak{g},\mathfrak{l})_q$. Векторное пространство всех
линейных функционалов на $V$ является $U_q\mathfrak{g}$-модулем. Применяя к
нему функтор $\Gamma$, получаем $U_q\mathfrak{g}$-модуль $V^\mathrm{c}\in
C(\mathfrak{g},\mathfrak{l})_q$, который будем называть двойственным
объектом категории $C(\mathfrak{g},\mathfrak{l})_q$. \IND{двойственный !
объект категории $C(\mathfrak{g},\mathfrak{l})_q$} Каждому морфизму
$f:V_1\to V_2$ в категории $C(\mathfrak{g},\mathfrak{l})_q$ очевидным
образом сопоставляется двойственный морфизм $f^*:V_2^\mathrm{c}\to
V_1^\mathrm{c}$. \IND{двойственный ! морфизм категории
$C(\mathfrak{g},\mathfrak{l})_q$} Отметим, что в частном случае
$\mathbb{S}=\varnothing$ категория $C(\mathfrak{g},\mathfrak{l})_q$
совпадает с категорией $P$-весовых $U_q\mathfrak{g}$-модулей, и введенное
понятие двойственности согласуется с описанным ранее в \itemiiiе
\ref{Weight}.

Получим разложение $U_q\mathfrak{g}$-модуля категории
$C(\mathfrak{g},\mathfrak{l})_q$ в сумму его
\hbox{$U_q\mathfrak{l}$-изотипических} компонент. Следующее утверждение
доказывается так же, как предложение \ref{semisimple1}.

\begin{proposition}\label{decomp}
Пусть $V\in C(\mathfrak{g},\mathfrak{l})_q$. Тогда в категории
$U_q\mathfrak{l}$-модулей $V=\bigoplus\limits_{\lambda\in\ P_+^\mathbb{S}}
V(\lambda)$, где $U_q\mathfrak{l}$-модули $V(\lambda)$ кратны
$L(\mathfrak{l},\lambda)$.
\end{proposition}

Следуя \cite{HechtSchmid}, \cite{KashSchmid}, конечно порожденный модуль
$V\in C(\mathfrak{g},\mathfrak{l})_q$ будем называть модулем Хариш-Чандры,
\IND{модуль ! Хариш-Чандры} если кратности вхождения в $V$ простых
$U_q\mathfrak{l}$-модулей $L(\mathfrak{l},\lambda)$ конечны:
$$
\dim_{U_q\mathfrak{l}}(L(\mathfrak{l,\lambda}),V)<\infty,\qquad\lambda\in
P_+^\mathbb{S},
$$
(другими словами, если $\dim\,V(\lambda)<\infty$ при всех $\lambda\in
P_+^\mathbb{S}$).

\begin{remark} Изучение модулей Хариш-Чандры является важным направлением
исследований в теории представлений. При $q=1$ имеются отдельные работы,
посвященные частному случаю $\mathbb{S}=\varnothing$. Например,
классификация простых модулей Хариш-Чандры в этом случае получена Матье
\cite{Mathieu}. Обширная литература посвящена модулям Хариш-Чандры,
связанным с представлениями вещественных редуктивных групп \cite{Schmid}.
Классификация простых модулей Хариш-Чандры в этом
случае является одним из наиболее известных достижений теории \cite[стр.
741, 747]{KV}. Задача о классификации простых модулей Хариш-Чандры в общем
случае обсуждается в статьях Пенкова, Цакермана и Сергановой \cite{PZ1, PZ2,
PS}.
\end{remark}

Важно отметить, что далеко не все построения теории вещественных редуктивных
групп переносятся на квантовый случай. Одна из причин в том, что подалгебра
Ли $\mathfrak{g}_0$ совпадает с подалгеброй неподвижных элементов инволюции
$\theta$, введенной в \itemiiiе \ref{Pol}, если и только если
\begin{equation}\label{short_grad}
\mathfrak{g}=\mathfrak{g}_{-1}\oplus\mathfrak{g}_0\oplus\mathfrak{g}_1.
\end{equation}

 В дальнейшем, желая обеспечить такой обрыв градуировки, будем
предполагать, что выбран простой корень $\alpha_{l_0}$, входящий в
разложение старшего корня $\delta$ алгебры Ли $\mathfrak{g}$ с
коэффициентом 1, и что $$\mathbb{S}=\{1,2,\ldots,l\}\setminus\{l_0\}.$$
Используя таблицы, приведенные в \cite{Bou4-6}, легко найти все такие
простые корни, см.
\itemiii\ \ref{prehomogeneous}.

По традиции в этом важном частном случае используются обозначения
$$
\mathfrak{k}=\mathfrak{g}_0,\qquad\mathfrak{p}^{\pm}=\mathfrak{g}_{\pm 1},
\qquad\mathfrak{q}^\pm=\mathfrak{k}\oplus\mathfrak{p}^\pm,
$$
и пара $(\mathfrak{g},\mathfrak{k})$ называется эрмитовой симметрической
\cite{Helg}. \IND{эрмитово симметрическая пара}

\begin{remark}\label{long_root}
1. Мы заменили вещественные подалгебры Ли, рассматриваемые в теории
эрмитовых симметрических пространств \cite{Helg}, их комплексификациями.

2. Хотя алгебра $\mathfrak{k}$ не является простой, для $U_q\mathfrak{k}$
имеют место результаты теории квантовых универсальных обертывающих алгебр,
приведенные в \itemiiе \ref{UE}.

3. Скалярное произведение в $\mathfrak{h}_\mathbb{R}$ выбрано в \itemiiiе
\ref{PBW} так, что квадраты длин коротких корней равны 2, а, как нетрудно
показать, $\alpha_{l_0}$ является длинным корнем. Значит, скалярное
произведение в $\mathfrak{h}_\mathbb{R}$, выбранное для $\mathfrak{g}$,
подходит и для $\mathfrak{k}$.

\end{remark}

\bigskip В \itemiiiе \ref{Pol} была введена эрмитова форма $(f_1,f_2)=\langle
f_2^*f_1\rangle$ в $\mathbb{C}[\mathfrak{p}^-]_q$.
\begin{proposition}\label{positive_herm}
Если $(\mathfrak{g},\mathfrak{k})$ -- эрмитово-симметрическая пара, то
$(f,f)>0$ для всех ненулевых элементов $f\in\mathbb{C}[\mathfrak{p}^-]_q$.
\end{proposition}

{\bf Доказательство.} Из $U_q\mathfrak{h}$-инвариантности формы
$(\cdot,\cdot)$ следует, что однородные компоненты
$\mathbb{C}[\mathfrak{p}^-]_{q,j}$ попарно ортогональны. Остается доказать
положительность сужений рассматриваемой эрмитовой формы на подпространства
$\mathbb{C}[\mathfrak{p}^-]_{q,j}$.

Рассмотрим базис векторного пространства $N(\mathfrak{q}^+,0)_{-j}$
$$
\left\{F_{\beta_M}^{j_M}F_{\beta_{M-1}}^{j_{M-1}}\ldots
F_{\beta_{M'+1}}^{j_{M'+1}}v(\mathfrak{q}^+,0)\left|\:
\sum\limits_{k=M'+1}^Mj_k=j\right.\right\},
$$
введенный в \itemiiiе \ref{GVerma}, и биортогональный базис в
$\mathbb{C}[\mathfrak{p}^-]_{q,j}$. Из предложения \ref{inner1} и
результатов \itemiiiа \ref{PBW} следует, что матрица эрмитовой формы
$(q^{-1}-q)^{-j}(\cdot,\cdot)$ в этом базисе невырождена и непрерывна по
переменной $q$ в интервале $(0,1)$. Используя равенства \eqref{inner_new},
\eqref{Rmatrix} и предложение \ref{inner1}, нетрудно доказать, что при
$f\in\mathbb{C}[\mathfrak{p}^-]_{q,j}$ и $q\to 1$:
$$(q^{-1}-q)^{-j}(f,f)\sim$$
$$
\sim f^*\otimes f\left(\prod\limits_{k>M'}^{\curvearrowleft}
\exp_{q_{\beta_k}^2}\left(\frac{q_{\beta_k}^{-1}-q_{\beta_k}}{q^{-1}-q}
E_{\beta_k}\otimes F_{\beta_k}\right)(v(\mathfrak{q}^-,0)\otimes
v(\mathfrak{q}^+,0))\right).
$$
Значит, матрица эрмитовой формы имеет предел при $q\to 1$. Остается
доказать положительность предельной матрицы. Но для эрмитово-симметрической
пары $(\mathfrak{g},\mathfrak{k})$ в пределе $q=1$ имеют место соотношения
\eqref{conjugate}, \eqref{def_star1}, а при $k>M'$ --- равенство
$E_{\beta_k}^*=-F_{\beta_k}$. Значит, предельная эрмитова форма
несущественными положительными множителями отличается от суммы квадратов
модулей коэффициентов ряда Тейлора функции $f$ в нуле. \hfill $\square$

\bigskip

 Особого внимания заслуживает изучение модулей Хариш-Чандры,
отвечающих эрмитовым симметрическим парам $(\mathfrak{g},\mathfrak{k})$, и
построение геометрических реализаций таких модулей Хариш-Чандры. На этом
пути будут получены содержательные результаты некоммутативного комплексного
анализа в квантовых ограниченных симметрических областях.

\subsubsection{Квантовые аналоги подалгебр $U\mathfrak{u}^{\pm}\subset
U\mathfrak{g}$.}\label{Add_Uq}

Результаты этого и следующего \itemiiiов получены Шкляровым, ср.
\cite{Shkl}. В \itemiiiе \ref{GVerma} рассматривалась градуированная
алгебра Ли $\mathfrak{g}$:
$\mathfrak{g}=\mathfrak{u}^-\newoplus\mathfrak{l}\newoplus\mathfrak{u}^+$.
В настоящем \itemiiiе будут введены $q$-аналоги универсальных обертывающих
алгебр $U\mathfrak{u}^{\pm}$. Точнее, мы рассмотрим $U_q\mathfrak{g}$ с
присоединенным действием подалгебры $U_q\mathfrak{l}$ и построим вложения
\hbox{$U_q\mathfrak{l}$-модульных} алгебр
$i^\pm:\mathbb{C}[\mathfrak{u}^\mp]_q\hookrightarrow U_q\mathfrak{g}$.
После этого положим $U_q\mathfrak{u}^\pm\overset{\mathrm{def}}{=}
i^\pm\mathbb{C}[\mathfrak{u}^\mp]_q$.

 Приведем доводы в пользу того, что построенные
таким образом подалгебры $U_q\mathfrak{u}^\pm$ можно считать $q$-аналогами
универсальных обертывающих алгебр $U\mathfrak{u}^\pm$. Во-первых, эти
подалгебры являются $U_q\mathfrak{l}$-модульными. Во-вторых,
$U_q\mathfrak{u}^\pm\subset U_q\mathfrak{b}^\pm$. В-третьих, образующие
$U_q\mathfrak{u}^\pm$ в формальном пределе $q\to1$ переходят в образующие
подалгебр $U\mathfrak{u}^\pm$, см. предложение \ref{obraz}. Наконец, как
будет показано в следующем \itemiiiе, с помощью подалгебр
$U_q\mathfrak{u}^\pm$ можно получить $q$-аналог разложения
$U\mathfrak{g}=U\mathfrak{u}^-\cdot U\mathfrak{l}\cdot U\mathfrak{u}^+$.

    Воспользуемся  сокращенными обозначениями Свидлера для коумножения:
$\Delta(\xi)=\xi_1\otimes\xi_2$, ср. \cite[стр. 55]{Kassel_QG}. Двойные
индексы означают повторное применение коумножения, например:
$(1\otimes\Delta)\Delta(\xi)=\xi_1\otimes\xi_{21}\otimes\xi_{22}$. В этих
обозначениях коассоциативность коумножения записывается следующим образом:
$\xi_{1}\otimes
\xi_{21}\otimes\xi_{22}=\xi_{11}\otimes\xi_{12}\otimes\xi_{2}.$

Начнем с построения вложения $i^+$. Рассмотрим обобщенный модуль Верма
$N(\mathfrak{q}^+,0)=\bigoplus\limits_{j=0}^\infty N(\mathfrak{q}^+,0)_{-j}$
со стандартной градуировкой и векторное пространство
$U_q\mathfrak{g}\widehat{\otimes}N(\mathfrak{q}^+,0)$ формальных рядов
$$
\sum_{j=0}^\infty c_j,\qquad c_j\in U_q\mathfrak{g}\otimes
N(\mathfrak{q}^+,0)_{-j}.
$$
Очевидно, $U_q\mathfrak{g}\otimes N(\mathfrak{q}^+,0)\hookrightarrow
U_q\mathfrak{g}\widehat{\otimes}N(\mathfrak{q}^+,0)$. Пусть
$$
\mathscr{K}=\sum_{j=0}^\infty a_j\otimes b_jv(\mathfrak{q}^+,0)\in
U_q\mathfrak{g}\widehat{\otimes}N(\mathfrak{q}^+,0),
$$
где $\sum\limits_{j=0}^\infty a_j\otimes b_j$ -- универсальная R-матрица.
Предполагается, что последнее разложение получено из мультипликативной
формулы \eqref{Rmatrix} умножением рядов для $q$-экспонент, а приведенное
разложение элемента $w_0$, выбрано так, как в \itemiiiе \ref{GVerma}.
Сомножитель $q^{-t_0}$ можно опустить, поскольку $v(\mathfrak{q}^+,0)$ --
весовой вектор нулевого веса.

\begin{proposition}\label{4.3_4.4}
Линейный оператор
$$
i^+:\mathbb{C}[\mathfrak{u}^-]_q\to U_q\mathfrak{g}, \qquad
i^+:f\mapsto(\operatorname{id}\otimes f)(\mathscr{K})
$$
является вложением $U_q\mathfrak{l}$-модульных алгебр.
\end{proposition}

{\bf Доказательство.} Докажем инъективность $i^+$. Воспользуемся тем, что
образующая $v(\mathfrak{q}^+,0)$ является $U_q\mathfrak{l}$-инвариантом:
\begin{equation}\label{v0inv}
\xi v(\mathfrak{q}^+,0)=\varepsilon(\xi)\cdot
v(\mathfrak{q}^+,0),\qquad\xi\in U_q\mathfrak{l}.
\end{equation}
Так как для выбранного в \itemiiiе \ref{GVerma} приведенного разложения
элемента $w_0$ элементы $F_{\beta_{1}},F_{\beta_{2}},\ldots F_{\beta_{M'}}$
принадлежат $U_q\mathfrak{l}$, то
\begin{multline}\label{k}
\mathscr{K}=\sum {c}_{j_M,j_{M-1},\ldots,j_{M'+1}}\cdot
\\ \cdot E_{\beta_M}^{j_M}E_{\beta_{M-1}}^{j_{M-1}}\ldots
E_{\beta_{M'+1}}^{j_{M'+1}}\otimes
F_{\beta_M}^{j_M}F_{\beta_{M-1}}^{j_{M-1}}\ldots
F_{\beta_{M'+1}}^{j_{M'+1}}v(\mathfrak{q}^+,0),
\end{multline}
причем все коэффициенты
$$
{c}_{j_M,j_{M-1},\ldots,j_{M'+1}},\qquad
(j_M,j_{M-1},\ldots,j_{M'+1})\in\mathbb{Z}_+^{M-M'},
$$
отличны о нуля, см. \eqref{Rmatrix}. Остается воспользоваться тем, что
элементы $E_{\beta_M}^{j_M}E_{\beta_{M-1}}^{j_{M-1}}\ldots
E_{\beta_{M'+1}}^{j_{M'+1}}$ линейно независимы, а элементы \\
$F_{\beta_M}^{j_M}F_{\beta_{M-1}}^{j_{M-1}}\ldots
F_{\beta_{M'+1}}^{j_{M'+1}}v(\mathfrak{q}^+,0)$ образуют базис в
$N(\mathfrak{q}^+,0)$, см. следствие \ref{PBWbasis-new} и предложение
\ref{N_basis}.

Покажем, что $i^+$ -- гомоморфизм алгебр. Пусть
$f_1,f_2\in\mathbb{C}[\mathfrak{u}^-]_q$. Тогда
\begin{multline*}
i^+(f_1f_2)=(\operatorname{id}\otimes f_1f_2)\left(\sum_{j=0}^\infty
a_j\otimes b_jv(\mathfrak{q}^+,0)\right)=\\
=(\operatorname{id}\otimes f_1\otimes f_2)\left(\sum_{j=0}^\infty
a_j\otimes\Delta^+(b_jv(\mathfrak{q}^+,0))\right),
\end{multline*}
где $\Delta^+$--- коумножение, введенное в \itemiiiе \ref{vect}. Из его
определения вытекает равенство
\begin{multline*}
(\operatorname{id}\otimes f_1\otimes f_2)\left(\sum_{j=0}^\infty a_j\otimes
\Delta^+(b_jv(\mathfrak{q}^+,0))\right)=
\\ =(\operatorname{id}\otimes f_1\otimes f_2)\left(\sum_{j=0}^\infty
a_j\otimes\Delta^{\rm cop}(b_j)(v(\mathfrak{q}^+,0)\otimes
v(\mathfrak{q}^+,0))\right).
\end{multline*}
Как следует из свойства \eqref{R_property3} универсальной R-матрицы,
$$
\sum_{j=0}^\infty a_j\otimes\Delta^{\rm op}(b_j)=\sum_{i,j=0}^\infty
a_ia_j\otimes b_i\otimes b_j.
$$
Остается заметить, что
$$
(\operatorname{id}\otimes f_1\otimes f_2)\left(\sum_{i,j=0}^\infty
a_ia_j\otimes b_iv(\mathfrak{q}^+,0)\otimes b_jv(\mathfrak{q}^+,0)\right)=
i^+(f_1)\,i^+(f_2).
$$

Для того, чтобы завершить доказательство предложения, покажем, что $i^+$
является морфизмом $U_q\mathfrak{l}$-модулей. Отождествим пространство
$U_q\mathfrak{g}\widehat{\otimes}N(\mathfrak{q}^+,0)$ с
$\mathrm{Hom}(\mathbb{C}[\mathfrak{u}^-]_q,U_q\mathfrak{g})$, сопоставляя
элементу $\sum\limits_j\eta_j\otimes v_j\in
U_q\mathfrak{g}\widehat{\otimes}N(\mathfrak{q}^+,0)$ линейный оператор $
f\mapsto(\operatorname{id}\otimes f)\left(\sum\limits_j\eta_j\otimes
v_j\right)$. Снабдим тензорное произведение
$U_q\mathfrak{g}\widehat{\otimes}N(\mathfrak{q}^+,0)$ структурой
$U_q\mathfrak{l}$-модуля, используя противоположное коумножение
$\Delta^{op}:U_q\mathfrak{l}\to U_q\mathfrak{l}\otimes U_q\mathfrak{l}$:
$$
\xi\left(\sum\limits_j\eta_j\otimes
v_j\right)=\sum\limits_j\operatorname{ad}(\xi_2)\eta_j\otimes\xi_1v_j.
$$

\begin{lemma} Если элемент $U_q\mathfrak{l}$-модуля
   $U_q\mathfrak{g}\widehat{\otimes}N(\mathfrak{q}^+,0)$ является
 $U_q\mathfrak{l}$-инвариантным, то отвечающий ему линейный оператор из
$\mathbb{C}[\mathfrak{u}^-]_q$ в $U_q\mathfrak{g}$-- морфизм
$U_q\mathfrak{l}$-модулей.
\end{lemma}

{\bf Доказательство леммы.} Пусть $A$ -- алгебра Хопфа, а $V_1,V_2$ --
некоторые $A$-модули. Каноническое отображение
$$
V_1^*\otimes V_1\otimes V_2\to V_2,\qquad f_1\otimes v_1\otimes v_2\mapsto
f_1(v_1)v_2
$$
является морфизмом $A$-модулей. Получаем отображение из пространства
 $A$-инвариантов  $(V_1\otimes V_2)^A$ в $\mathrm{Hom}_{A}(V_1^*,V_2)$.
Остается применить эти соображения, полагая $A=U_q\mathfrak{l}$,
$V_1=N(\mathfrak{q}^+,0)$, $V_2=U_q\mathfrak{g}$. \hfill $\square$

\medskip

Таким образом, для завершения доказательства предложения \ref{4.3_4.4}
достаточно показать, что элемент $\mathscr{K}$ является
$U_q\mathfrak{l}$-инвариантом. Пусть \hbox{$\xi\in U_q\mathfrak{l}$}. Тогда
\begin{multline*}
\xi(\mathscr{K})=\xi\left(\sum\limits_ja_j\otimes
b_jv(\mathfrak{q}^+,0)\right)=\sum\limits_j\operatorname{ad}(\xi_2)a_j
\otimes\xi_1b_jv(\mathfrak{q}^+,0)=
\\ =\sum\limits_j\xi_{21}a_jS(\xi_{22})\otimes\xi_1b_jv(\mathfrak{q}^+,0).
\end{multline*}
В силу коассоциативности
$$
\sum\limits_j\xi_{21}a_jS(\xi_{22})\otimes\xi_1b_jv(\mathfrak{q}^+,0)=
\sum\limits_j\xi_{12}a_jS(\xi_{2})\otimes\xi_{11}b_jv(\mathfrak{q}^+,0).
$$
Согласно свойству \eqref{R_property1} универсальной R-матрицы,
$$
\sum\limits_j\xi_{12}a_j\otimes\xi_{11}b_j=
\sum\limits_ja_j\,\xi_{11}\otimes b_j\xi_{12}.
$$
Значит,
$\sum\limits_j\xi_{12}a_jS(\xi_{2})\otimes\xi_{11}b_jv(\mathfrak{q}^+,0)=
\sum\limits_ja_j\,\xi_{11}S(\xi_{2})\otimes
b_j\xi_{12}v(\mathfrak{q}^+,0)$. Из $U_q\mathfrak{l}$-инвариантности
образующей $v(\mathfrak{q}^+,0)$ (см. \eqref{v0inv}) вытекает равенство
$$
\sum\limits_ja_j\xi_{11}S(\xi_2)\otimes
b_j\xi_{12}v(\mathfrak{q}^+,0)=\sum\limits_ja_j\xi_{11}S(\xi_2)\otimes
b_j\varepsilon(\xi_{12})v(\mathfrak{q}^+,0).
$$
По определению коединицы и антипода в алгебре Хопфа
$$
\sum\limits_ja_j\,\xi_{11}S(\xi_2)\otimes
b_j\varepsilon(\xi_{12})v(\mathfrak{q}^+,0)=
\sum\limits_ja_j\xi_{11}\varepsilon(\xi_{12})S(\xi_2)\otimes
b_jv(\mathfrak{q}^+,0)=
$$
$$
=\sum\limits_ja_j\xi_{1}S(\xi_2)\otimes
b_jv(\mathfrak{q}^+,0)=\varepsilon(\xi)\sum\limits_ja_j\otimes
b_jv(\mathfrak{q}^+,0).
$$
Следовательно, элемент $\mathscr{K}$ является
$U_q\mathfrak{l}$-инвариантным, и предложение \ref{4.3_4.4} полностью
доказано. $\hfill$ $\square$

\bigskip

Введем обозначение $U_q\mathfrak{u}^+$ для образа линейного оператора $i^+$.
\begin{corollary}\label{U_qu}
$U_q\mathfrak{u}^+$ является $U_q\mathfrak{l}$-модульной подалгеброй в
$U_q\mathfrak{g}$.
\end{corollary}

Из формулы \eqref{k} и следующих за ней рассуждений вытекает

\begin{proposition}\label{bas1} Элементы
\begin{equation}\label{bu+}
E_{\beta_M}^{j_M}E_{\beta_{M-1}}^{j_{M-1}}\ldots
E_{\beta_{M'+1}}^{j_{M'+1}},\qquad
j_{M'+1},\ldots,j_{M-1},j_M\in\mathbb{Z}_+,
\end{equation}
образуют базис векторного пространства $U_q\mathfrak{u}^+$.
\end{proposition}

\medskip Отсюда и из предложения \ref{Levend} получаем
\begin{proposition}\label{obraz}
$U_q\mathfrak{u}^+$ -- подалгебра с единицей в $U_q\mathfrak{g}$,
порожденная элементами
$E_{\beta_{M'+1}},\ldots,E_{\beta_{M-1}},E_{\beta_{M}}$.
\end{proposition}

\medskip Последнее утверждение и результаты \itemiiiа \ref{PBW} позволяют
переставить сомножители в \eqref{bu+}.
\begin{corollary}\label{bas2} Элементы
\begin{equation}\label{u+b}
E_{\beta_{M'+1}}^{j_{M'+1}}\ldots
E_{\beta_{M-1}}^{j_{M-1}}E_{\beta_M}^{j_M},\qquad
j_{M'+1},\ldots,j_{M-1},j_M\in\mathbb{Z}_+,
\end{equation}
образуют базис векторного пространства $U_q\mathfrak{u}^+$.
\end{corollary}

Обратимся к построению $q$-аналога универсальной обертывающей алгебры
$U\mathfrak{u}^-$.

Напомним, что $U_q\mathfrak{g}$ является $U_q\mathfrak{g}$-модульной
алгеброй по отношению к присоединенному действию, и инволюция \eqref{star_1}
наделяет $U_q\mathfrak{g}$ структурой $*$-алгебры Хопфа. Следующее
утверждение хорошо известно.

\begin{proposition}\label{*_add}
$*$-Алгебра $(U_q\mathfrak{g},*)$ является $(U_q\mathfrak{g},*)$-модульной
алгеброй
$$(\operatorname{ad}(\xi)\eta)^*=\operatorname{ad}((S(\xi))^*)\eta^*, \qquad
\xi, \eta \in U_q\mathfrak{g}.$$
\end{proposition}

{\bf Доказательство.} Очевидно
\begin{equation}\label{*_ad_1}
(\operatorname{ad}(\xi)\eta)^*=\left(\xi_1\eta
S(\xi_2)\right)^*=(S(\xi_2))^*\eta^*\xi_1^*.
\end{equation}

С другой стороны, $*$ является автоморфизмом, а $S$ --- антиавтоморфизмом
коалгебры $U_q\mathfrak{g}$. Следовательно,
$\triangle((S(\xi))^*)=(S(\xi_2))^*\otimes(S(\xi_1))^*.$ Значит,
\begin{equation}\label{*_ad_2}
\operatorname{ad}((S(\xi))^*)\eta^*=(S(\xi_2))^*\eta^*S((S(\xi_1))^*).
\end{equation}
Так как $S((S(\xi))^*)=\xi^*$ для любого $\xi\in U_q\mathfrak{g}$, то правая
часть \eqref{*_ad_1} совпадает с правой частью \eqref{*_ad_2}. \hfill
$\square$

\bigskip Из последнего предложения немедленно следует, что подалгебра
\begin{equation}\label{def_u_-}
U_q\mathfrak{u}^-\overset{\mathrm{def}}{=} \{\xi^*|\:\xi\in
U_q\mathfrak{u}^+\}
\end{equation}
алгебры $U_q\mathfrak{g}$ является $U_q\mathfrak{l}$-модульной.

\begin{remark}\label{u_g_minus}
Рассмотрим линейный оператор
$$
i^-:\mathbb{C}[\mathfrak{u}^+]_q\to U_q\mathfrak{g},\qquad
i^-:~f\mapsto(i^+(f^*))^*. \footnote{Здесь инволюция внутри скобок -- это
определенная в предложении \ref{star_Pol} инволюция в алгебре
$\mathbb{C}[\mathfrak{u}^-\oplus\mathfrak{u}^+]_q$, а инволюция вне скобок
-- это инволюция \eqref{star_1} в $U_q\mathfrak{g}$.}
$$
 Он является вложением $U_q\mathfrak{l}$-модульных алгебр, как следует из
предложений \ref{*_add}, \ref{4.3_4.4}, и
 $i^-\mathbb{C}[\mathfrak{u}^+]_q=U_q\mathfrak{u}^-$.
\end{remark}

\medskip Операторы $i^\pm$ позволяют построить вложение
$U_q\mathfrak{l}$-модулей
$$
i:\mathbb{C}[\mathfrak{u}^-]_q\otimes\mathbb{C}[\mathfrak{u}^+]_q\to
U_q\mathfrak{g},\quad i:f_1\otimes f_2\mapsto i^+(f_1)\cdot i^-(f_2).
$$
Это вложение $i$ уважает\, инволюции в
$\operatorname{Pol}(\mathfrak{u}^-)_q=
\mathbb{C}[\mathfrak{u}^-]_q\otimes\mathbb{C}[\mathfrak{u}^+]_q$ и в
$U_q\mathfrak{g}$:
\begin{equation}\label{respect*}
(i(f))^*=i(f^*),\qquad f\in\operatorname{Pol}(\mathfrak{u}^-)_q,
\end{equation}
но не является гомоморфизмом алгебр.

\medskip
В оставшейся части \itemiiiа мы будем предполагать, что
$(\mathfrak{g},\mathfrak{k})$ -- эрмитово-симметрическая пара и использовать
обозначения $\mathfrak{k}$, $\mathfrak{p}^{\pm}$ вместо $\mathfrak{l}$,
$\mathfrak{u}^{\pm}$, следуя \itemiiiу \ref{Hermit}. Напомним, что в этом
случае $\mathbb{S}=\{1,2,\ldots,l\}\backslash\{l_0\}$,
$$
\mathfrak{k}=\mathfrak{g}_0,\quad \mathfrak{p}^{\pm}=\mathfrak{g}_{\pm 1},
\quad \mathfrak{q}^\pm=\mathfrak{g}_0 \oplus \mathfrak{g}_{\pm 1}.
$$

\begin{proposition}\label{hom_gen}
Если $(\mathfrak{g},\mathfrak{k})$ -- эрмитово-симметрическая пара, то
подпространство $\mathbb{C}[\mathfrak{p}^-]_{q,1}$ порождает
подалгебру $\mathbb{C}[\mathfrak{p}^-]_q$, а подпространство
$\mathbb{C}[\mathfrak{p}^+]_{q,-1}$ -- подалгебру
$\mathbb{C}[\mathfrak{p}^+]_q$.
\end{proposition}

{\bf Доказательство.} Достаточно показать, что
$\mathbb{C}[\mathfrak{p}^-]_{q,1}$ порождает $\mathbb{C}[\mathfrak{p}^-]_q$.
Это подпространство и эту алгебру можно заменить их образами при вложении
$i^+$ в $U_q\mathfrak{g}$. Остается воспользоваться тем, что элементы
\eqref{bu+} образуют базис в $i^+\mathbb{C}[\mathfrak{p}^-]_q$, а
подпространство $i^+\mathbb{C}[\mathfrak{p}^-]_{q,1}$ содержит все
''корневые векторы''
\begin{equation}\label{root_hermit}
E_{\beta_{M'+1}},E_{\beta_{M'+2}},\ldots,E_{\beta_{M-1}},E_{\beta_M},
\end{equation}
поскольку пара $(\mathfrak{g},\mathfrak{k})$ является
эрмитово-симметрической, а в этом случае векторы \eqref{root_hermit} имеют
степень однородности 1. \hfill $\square$

\medskip
  Подчеркнем, что в  эрмитово-симметрическом случае используется
  обозначение $U_q\mathfrak{p}^\pm$ вместо $U_q\mathfrak{u}^\pm$.

\begin{proposition}\label{Jakobsen_def}
Алгебра $U_q\mathfrak{p}^+$ является наименьшей
$U_q\mathfrak{k}$-мо\-дуль\-ной унитальной подалгеброй в $U_q\mathfrak{g}$,
содержащей элемент $E_{l_0}$.
\end{proposition}

{\bf Доказательство.} Пусть $\mathcal{U}^+$ -- наименьшая
$U_q\mathfrak{k}$-модульная подалгебра с единицей алгебры
$U_q\mathfrak{g}$, содержащая $E_{l_0}$. Требуется доказать, что
\hbox{$\mathcal{U}^+=U_q\mathfrak{p}^+$}.

Докажем включение $\mathcal{U}^+\subset U_q\mathfrak{p}^+$. Для этого
достаточно показать, что $E_{l_0}\in U_q\mathfrak{p}^+$. Но среди корней
$\beta_{M'+1},\ldots,\beta_{M-1},\beta_{M}$ есть выделенный простой корень
$\alpha_{l_0}$. Пусть, для определенности, $\beta_{N}=\alpha_{l_0}$,
$M'+1\le N\le M$. Элементы $E_{\beta_{N}}$ и $E_{l_0}$ отличаются лишь
(ненулевым) множителем, поскольку оба они принадлежат пересечению весового
подпространства $(U_q\mathfrak{g})_{\alpha_{l_0}}$ с подалгеброй
$U_q\mathfrak{n}^+$, а это пересечение одномерно (см. предложение
\ref{PBW-basis}). Остается воспользоваться тем, что $E_{\beta_{N}}\in
U_q\mathfrak{p}^+$, как следует из предложения \ref{bas1}.

Перейдем к доказательству включения $\mathcal{U}^+\supset
U_q\mathfrak{p}^+$. Напомним, что в \itemiiiе \ref{GVerma} были наделены
градуировкой все весовые $U_q\mathfrak{g}$-модули и, в частности,
$U_q\mathfrak{g}$ с присоединенным действием. Подпространства
$\mathcal{U}^+$, $U_q\mathfrak{p}^+$ наследуют градуировку, и их однородные
компоненты являются $U_q\mathfrak{k}$-подмодулями.
 Так как  $\mathcal{U}^+ \subset U_q\mathfrak{p}^+$, то
 $(\mathcal{U}^+)_1 \subset (U_q\mathfrak{p}^+)_1$.
  Остается доказать, неравенство
  $$\dim(\mathcal{U}^+)_1 \geq \dim(U_q\mathfrak{p}^+)_1.$$
Действительно, из него будет следовать включение $(\mathcal{U}^+)_1 \supset
(U_q\mathfrak{p}^+)_1$ и требуемое утверждение $\mathcal{U}^+\supset
U_q\mathfrak{p}^+$, поскольку $U_q\mathfrak{k}$-подмодуль
$(U_q\mathfrak{p}^+)_1$ порождает алгебру $U_q\mathfrak{p}^+$.

Прежде всего, отметим, что размерность $(U_q\mathfrak{p}^+)_1$ равна ее
значению в предельном случае $q=1$ и равна $M-M'$, поскольку элементы
$E_{\beta_{M'+1}},\ldots,E_{\beta_{M-1}},E_{\beta_{M}}$ имеют степень
однородности 1 и образуют базис в $(U_q\mathfrak{p}^+)_1$ (см. предложение
\ref{obraz}).

С другой стороны, $E_{l_0}$ является однородным элементом степени 1, и,
следовательно,
$(\mathcal{U}^+)_1\supset\operatorname{ad}(U_q\mathfrak{k})E_{l_0}$. Кроме
того, $E_{l_0}$ -- особый вектор $U_q\mathfrak{k}$-модуля
$\operatorname{ad}(U_q\mathfrak{k})E_{l_0}$, то есть
$$
\mathrm{ad}(F_j)\,E_{l_0}\,=\,F_jE_{l_0}K_j\,-\,E_{l_0}F_jK_j\,=\,0,\qquad
j\ne l_0,
$$
и весовой вектор веса $\alpha_{l_0}$. Из предложения \ref{qWeyl} следует,
что размерности простых весовых модулей не меняются при переходе от
классических универсальных обертывающих алгебр к квантовым.

Значит, как и при $q=1$, $$\dim(\operatorname{ad}(U_q\mathfrak{k})E_{l_0})
\geq \dim((U_q\mathfrak{p}^+)_1).$$ Следовательно,
$$
\dim(\mathcal{U}^+)_1\geq \dim(\operatorname{ad}(U_q\mathfrak{k})E_{l_0})
 \geq \dim(U_q\mathfrak{p}^+)_1. \eqno\square
 $$

\bigskip

Описание алгебры $U_q\mathfrak{p}^+$, приведенное в предложении
\ref{Jakobsen_def}, заимствовано из работы Якобсена \cite{Jak-Hermit}, где
оно использовалось в качестве определения квантового аналога универсальной
обертывающей алгебры $U\mathfrak{p}^+$ в эрмитово симметрическом случае.
Одной из основных целей работы \cite{Jak-Hermit} было построение аналога
алгебры $\mathbb{C}[\mathfrak{p}^-]$ полиномов на векторном пространстве
$\mathfrak{p}^-$, и этот аналог {\it определяется} в \cite{Jak-Hermit}
равенством $\mathbb{C}[\mathfrak{p}^-]_q=U_q\mathfrak{p}^+$. В
эрмитово-симметрическом случае такое определение может быть мотивировано
следующим образом: при $q=1$ форма Киллинга на $\mathfrak{g}$ устанавливает
изоморфизм $U\mathfrak{k}$-модулей $\mathfrak{p}^-\cong(\mathfrak{p}^+)^*$,
и, следовательно, изоморфизм $U\mathfrak{k}$-модульных алгебр
$\mathbb{C}[\mathfrak{p}^-]\cong S(\mathfrak{p}^+)$. Близкий подход к
построению квантового аналога алгебры $\mathbb{C}[\mathfrak{p}^-]$
используется в \cite{KMT, Baldoni_Frajria}.

Предложение \ref{Jakobsen_def} и изоморфизм
$i^+:\mathbb{C}[\mathfrak{p}^-]_q\rightarrow U_q\mathfrak{p}^+$, введенный в
этом \itemiiiе, по существу, устанавливают эквивалентность нашего подхода и
подхода работы \cite{Jak-Hermit} к построению аналога алгебры
$\mathbb{C}[\mathfrak{p}^-]$.

\subsubsection{Квантовый аналог разложения \boldmath
$U\mathfrak{g}=U\mathfrak{u}^-\cdot U\mathfrak{l}\cdot
U\mathfrak{u}^+$.}\label{decomp_Dima}

Из равенства
$\mathfrak{g}=\mathfrak{u}^-\oplus\mathfrak{l}\oplus\mathfrak{u}^+$ и из
теоремы Пуанкаре-Биркгофа-Витта следует, что имеет место разложение
\begin{equation}\label{111}
U\mathfrak{g}=U\mathfrak{u}^-\cdot U\mathfrak{l}\cdot U\mathfrak{u}^+.
\end{equation}
Точнее говоря, линейное отображение
$$
m:U\mathfrak{u}^-\otimes U\mathfrak{l}\otimes U\mathfrak{u}^+ \rightarrow
U\mathfrak{g},\qquad m: \xi \otimes \eta \otimes \zeta \mapsto \xi\eta\zeta,
$$
определяемое с помощью умножения в алгебре $U\mathfrak{g}$, является
биективным. Именно так термин ''разложение'' будет пониматься в настоящем
\itemiiiе.

Получим квантовый аналог \eqref{111}. Точнее, покажем, что имеют место шесть
разложений
\begin{multline}\label{fact}
U_q\mathfrak{g}=U_q\mathfrak{u}^-\cdot U_q\mathfrak{l}\cdot
U_q\mathfrak{u}^+=U_q\mathfrak{u}^-\cdot U_q\mathfrak{u}^+\cdot
U_q\mathfrak{l}=U_q\mathfrak{l}\cdot U_q\mathfrak{u}^-\cdot
U_q\mathfrak{u}^+\\=U_q\mathfrak{l}\cdot U_q\mathfrak{u}^+\cdot
U_q\mathfrak{u}^-=U_q\mathfrak{u}^+\cdot U_q\mathfrak{u}^-\cdot
U_q\mathfrak{l} =U_q\mathfrak{u}^+\cdot U_q\mathfrak{l}\cdot
U_q\mathfrak{u}^-.
\end{multline}

\begin{proposition}\label{factor+}
Имеют место разложения
\begin{equation}\label{factor++}
U_q\mathfrak{q}^+=U_q\mathfrak{l}\cdot U_q\mathfrak{u}^+, \qquad
U_q\mathfrak{q}^+= U_q\mathfrak{u}^+\cdot U_q\mathfrak{l}.
\end{equation}
\end{proposition}

{\bf Доказательство.} Получим первое из них. Используя разложение
$U_q\mathfrak{g}=U_q\mathfrak{n}^-\cdot U_q\mathfrak{h}\cdot
U_q\mathfrak{n}^+ $ (см. \itemiii\ \ref{PBW}), получаем:
$$
U_q\mathfrak{q}^+=(U_q\mathfrak{q}^+\cap U_q\mathfrak{n}^-)\cdot
U_q\mathfrak{h}\cdot U_q\mathfrak{n}^+,\quad
U_q\mathfrak{l}=(U_q\mathfrak{l}\cap U_q\mathfrak{n}^-)\cdot
U_q\mathfrak{h}\cdot(U_q\mathfrak{l}\cap U_q\mathfrak{n}^+).
$$
Так как $U_q\mathfrak{u}^+\subset U_q\mathfrak{n}^+$ и
$U_q\mathfrak{q}^+\cap U_q\mathfrak{n}^-=U_q\mathfrak{l}\cap
U_q\mathfrak{n}^-$, то достаточно получить разложение
$
U_q\mathfrak{n}^+=(U_q\mathfrak{l}\cap U_q\mathfrak{n}^+)\cdot
U_q\mathfrak{u}^+.
$
Остается воспользоваться тем, что элементы
\begin{equation}
E_{\beta_1}^{j_1}E_{\beta_{2}}^{j_{2}}\ldots E_{\beta_{M'}}^{j_{M'}},\quad
j_{1},j_{2},\ldots,j_{M'}\in\mathbb{Z}_+,
\end{equation}
образуют базис векторного пространства $U_q\mathfrak{l}\cap
U_q\mathfrak{n}^+$, а также следствием \ref{bas2} и предложением
\ref{PBW-basis}.

Второе разложение в \eqref{factor++} доказывается аналогично, нужно лишь
начать с разложения $ U_q\mathfrak{g}=U_q\mathfrak{n}^+\cdot
U_q\mathfrak{h}\cdot U_q\mathfrak{n}^- $, а в конце воспользоваться
предложением \ref{bas1}, следствием \ref{PBWbasis-new}, и тем, что элементы
\begin{equation}
E_{\beta_{M'}}^{j_{M'}}\ldots E_{\beta_2}^{j_2}E_{\beta_{1}}^{j_{1}},\qquad
j_{1},j_{2},\ldots,j_{M'}\in\mathbb{Z}_+,
\end{equation}
образуют базис векторного пространства $U_q\mathfrak{l}\cap
U_q\mathfrak{n}^+$. \hfill $\square$

\medskip

Применяя инволюцию $*$ к \eqref{factor++}, получаем

\begin{corollary}\label{factor-}
Имеют место разложения
\begin{equation}\label{factor--}
U_q\mathfrak{q}^-=U_q\mathfrak{l}\cdot U_q\mathfrak{u}^-, \qquad
U_q\mathfrak{q}^-= U_q\mathfrak{u}^-\cdot U_q\mathfrak{l}.
\end{equation}
\end{corollary}

\begin{proposition}\label{factor+++}
Имеют место разложения
\begin{equation}\label{factor++++}
U_q\mathfrak{g}=U_q\mathfrak{q}^-\cdot U_q\mathfrak{u}^+,\qquad
U_q\mathfrak{g}= U_q\mathfrak{u}^+\cdot U_q\mathfrak{q}^-.
\end{equation}
\end{proposition}

{\bf Доказательство} этого утверждения легко получить, анализируя
доказательство предложения \ref{factor+}. Действительно,
$$
U_q\mathfrak{g}=U_q\mathfrak{n}^-\cdot U_q\mathfrak{h}\cdot
U_q\mathfrak{n}^+=U_q\mathfrak{n}^-\cdot U_q\mathfrak{h}\cdot
(U_q\mathfrak{l}\cap U_q\mathfrak{n}^+)\cdot
U_q\mathfrak{u}^+=U_q\mathfrak{q}^-\cdot U_q\mathfrak{u}^+.
$$
 Аналогично,
$$
U_q\mathfrak{g}=U_q\mathfrak{n}^+\cdot U_q\mathfrak{h}\cdot
U_q\mathfrak{n}^-=U_q\mathfrak{u}^+\cdot(U_q\mathfrak{l}\cap
U_q\mathfrak{n}^+)\cdot U_q\mathfrak{h}\cdot
U_q\mathfrak{n}^-=U_q\mathfrak{u}^+\cdot U_q\mathfrak{q}^-.\eqno\square
$$

Применяя инволюцию $*$ к (\ref{factor++++}), получаем
\begin{corollary}\label{factor---}
Имеют место разложения
\begin{equation}\label{factor----}
U_q\mathfrak{g}=U_q\mathfrak{q}^+\cdot U_q\mathfrak{u}^-,\qquad
U_q\mathfrak{g}= U_q\mathfrak{u}^-\cdot U_q\mathfrak{q}^+.
\end{equation}
\end{corollary}

\medskip
Объединяя разложения в \eqref{factor++++} и \eqref{factor----} с
разложениями в \eqref{factor++} и \eqref{factor--}, мы приходим к
результату, анонсированному в начале \itemiiiа.
\begin{proposition}\label{6razl}
Имеют место разложения \eqref{fact}.
\end{proposition}

\subsubsection{Квадратичная алгебра
$\mathbb{C}[\mathfrak{p}^-]_q$.}\label{quadratic_first}

В \itemiiiе \ref{Add_Uq} было показано, что в эрмитово-симметрическом
случае подпространство $\mathbb{C}[\mathfrak{p}^-]_{q,1}$ порождает алгебру
$\mathbb{C}[\mathfrak{p}^-]_q$. Найдем определяющие соотношения.

Рассмотрим линейное отображение
$$ m:\,\mathbb{C}[\mathfrak{p}^-]_q
\otimes \mathbb{C}[\mathfrak{p}^-]_q \rightarrow
\mathbb{C}[\mathfrak{p}^-]_q, \qquad m:\,f_1 \otimes f_2 \mapsto f_1f_2.$$ и
универсальную $R$-матрицу алгебры Хопфа $U_q\mathfrak{g}$. Из определения
умножения $m$ следует равенство
\begin{equation}\label{comm_general}
 \varphi \psi=m\check{R}_{\mathbb{C}[\mathfrak{p}^-]_q,\mathbb{C}[\mathfrak{p}^-]_q}
 (\varphi \otimes \psi),\qquad \varphi, \psi \in \mathbb{C}[\mathfrak{p}^-]_q,
\end{equation}
означающее, что $\mathbb{C}[\mathfrak{p}^-]_q$ является коммутативной
алгеброй в сплетенной тензорной категории $\mathcal{C}^-$, см. \itemiii\
\ref{braided_categories_sl_2}.

Пусть $\mathcal{L}$ -- ядро линейного отображения
$$
\mathbb{C}[\mathfrak{p}^-]_{q,1}\otimes\mathbb{C}[\mathfrak{p}^-]_{q,1}\to
\mathbb{C}[\mathfrak{p}^-]_{q,2},\qquad\varphi\otimes\psi\mapsto\varphi\psi-
m\check{R}_{\mathbb{C}[\mathfrak{p}^-]_q,\mathbb{C}[\mathfrak{p}^-]_q}
(\varphi\otimes\psi).
$$
В классическом предельном случае $q=1$ это подпространство всех
антисимметрических тензоров.

Рассмотрим подалгебру $U_q\mathfrak{k}_\mathrm{ss}$, порожденную множеством
\\ $\{K^{\pm 1}_j,\,E_j,\,F_j\}_{j\ne l_0}$.

Как следует из \eqref{comm_general}, $m\mathcal{L}=0$. Другими словами,
элементы подпространства $\mathcal{L}$ являются квадратичными
соотношениями.
 Опишем это подпространство в терминах
универсальной $R$-матрицы алгебры Хопфа $U_q\mathfrak{k}_{\rm ss}$ и
отвечающего ей
 морфизма \hbox{$U_q\mathfrak{k}$-модулей}
 $$\widetilde{R}_{\mathbb{C}[\mathfrak{p}^-]_{q,1},\mathbb{C}[\mathfrak{p}^-]_{q,1}}:
 \mathbb{C}[\mathfrak{p}^-]_{q,1}\otimes
\mathbb{C}[\mathfrak{p}^-]_{q,1}\rightarrow
 \mathbb{C}[\mathfrak{p}^-]_{q,1}\otimes
 \mathbb{C}[\mathfrak{p}^-]_{q,1},$$
  определяемого действием универсальной $R$-матрицы
 с последующей ''наивной'' перестановкой тензорных
 сомножителей.

\begin{lemma}\label{eigenvectors_in_L}  $\mathcal{L}$
содержит все собственные векторы с отрицательными собственными
 значениями оператора
$\widetilde{R}_{\mathbb{C}[\mathfrak{p}^-]_{q,1},
\mathbb{C}[\mathfrak{p}^-]_{q,1}}$.
\end{lemma}

{\bf Доказательство.} Пусть $\mathcal{L}'$ -- спектральное подпространство
линейного оператора $\widetilde{R}_{\mathbb{C}[\mathfrak{p}^-]_{q,1},
\mathbb{C}[\mathfrak{p}^-]_{q,1}}$, отвечающее отрицательной полуоси
$(-\infty,0)$. Очевидно, $\mathcal{L}$ и $\mathcal{L}'$ являются
$U_q\mathfrak{k}$-подмодулями, и нам достаточно доказать, что
$\mathcal{L}\supset\mathcal{L}'$.

 В классическом случае $q=1$ кратности вхождения
простых весовых $U\mathfrak{k}$-модулей в $(\mathfrak{p}^-)^* \otimes
(\mathfrak{p}^-)^*$ не превосходят 1, поскольку весовые подпространства
$U\mathfrak{k}$-модуля $\mathfrak{p}^-$ одномерны \cite[стр.
591--592]{Zhelob}. Значит, согласно предложению \ref{qWeyl}, эти кратности
равны 1 при всех $q \in (0,1)$. Следовательно, подпространства
$\mathcal{L}$ и $\mathcal{L}'$, определяются своими
$U_q\mathfrak{k}$-спектрами, то есть множествами старших весов своих
простых $U_q\mathfrak{k}$-подмодулей.

При $q=1$ эти множества старших весов совпадают. Остается проследить за их
зависимостью от параметра $q\in (0,1]$. Как следует из \eqref{Resh_Dr},
спектр оператора $\widetilde{R}_{\mathbb{C}[\mathfrak{p}^-]_{q,1},
\mathbb{C}[\mathfrak{p}^-]_{q,1}}$ лежит на вещественной оси и не содержит
0.
 Значит, аналитическая зависимость спектрального проектора, отвечающего
 отрицательной полуоси, вытекает из аналитической зависимости оператора
$\widetilde{R}_{\mathbb{C}[\mathfrak{p}^-]_{q,1},
\mathbb{C}[\mathfrak{p}^-]_{q,1}}$.

Проследим за зависимостью от параметра $q$ операторов
$$E_j,F_j,\quad j\ne l_0,\qquad H_i,\quad i=1,2,\ldots,l,$$
подпространства $\mathcal{L}$ и оператора
$\widetilde{R}_{\mathbb{C}[\mathfrak{p}^-]_{q,1},
\mathbb{C}[\mathfrak{p}^-]_{q,1}}$. Для этого выберем базис весовых
векторов в $\mathbb{C}[\mathfrak{p}^-]_q$, биортогональный описанному в
предложении \ref{N_basis}. Получаем базисы в
$\mathbb{C}[\mathfrak{p}^-]_{q,1}$,
$\mathbb{C}[\mathfrak{p}^-]_{q,1}\otimes\mathbb{C}[\mathfrak{p}^-]_{q,1}$,
$\mathbb{C}[\mathfrak{p}^-]_{q,2}$. Требуемое утверждение
$\mathcal{L}\supset\mathcal{L}'$ вытекает из того, что матричные элементы
операторов $m$ и $\widetilde{R}_{\mathbb{C}[\mathfrak{p}^-]_{q,1},
\mathbb{C}[\mathfrak{p}^-]_{q,1}}$ в рассматриваемых базисах аналитически
зависят от $q\in(0,1]$. \hfill $\square$

   \begin{proposition}\label{new_relations} $\mathcal{L}$ является
   собственным подпространством  линейного оператора
 $\widetilde{R}_{\mathbb{C}[\mathfrak{p}^-]_{q,1},
 \mathbb{C}[\mathfrak{p}^-]_{q,1}}$,
 отвечающим его единственному отрицательному собственному значению
$-q^{\frac{4}{(H_\mathbb{S},H_\mathbb{S})}}$. Кратность этого
 собственного значения равна
 $\frac{\dim \mathfrak{p}^-(\dim \mathfrak{p}^--1)}{2}$.
    \end{proposition}

{\bf Доказательство.} Размерность пространства $\mathcal{L}$ не превосходит
$\frac{\dim\mathfrak{p}^-(\dim \mathfrak{p}^--1)}{2}$, поскольку
 $\dim \mathbb{C}[\mathfrak{p}^-]_{q,2}\;=\;
  \frac{\dim\mathfrak{p}^-(\dim \mathfrak{p}^- +1)}{2}$.

Остается воспользоваться леммой \ref{eigenvectors_in_L} и тем, что у
оператора $\widetilde{R}_{\mathbb{C}[\mathfrak{p}^-]_{q,1},
\mathbb{C}[\mathfrak{p}^-]_{q,1}}$ размерность собственного
подпространства, отвечающего собственному значению
$-q^{\frac{4}{(H_\mathbb{S},H_\mathbb{S})}}$, больше или равна
$\frac{\dim\mathfrak{p}^-(\dim\mathfrak{p}^--1)}2$. Докажем последнее
неравенство.

Рассмотрим сужение линейного оператора
$\check{R}_{N(\mathfrak{q}^+,-\alpha_{l_0})^*,
N(\mathfrak{q}^+,-\alpha_{l_0})^*}$ на
$N(\mathfrak{q}^+,-\alpha_{l_0})_{-1}^*\otimes
N(\mathfrak{q}^+,-\alpha_{l_0})_{-1}^*$. Нетрудно показать, что это сужение
отличается от интересующего нас линейного оператора
$\widetilde{R}_{\mathbb{C}[\mathfrak{p}^-]_{q,1},
\mathbb{C}[\mathfrak{p}^-]_{q,1}}$ лишь числовым множителем
$q^{-\frac4{(H_\mathbb{S},H_\mathbb{S})}}$.
(Для сравнения универсальных $R$-матриц алгебр Хопфа $U_q\mathfrak{g}$ и
$U_q\mathfrak{k}_\mathrm{ss}$ используется приведенное разложение элемента
$w_0$, описанное в \itemiiiе \ref{GVerma}, равенства \eqref{Rmatrix},
\eqref{t0_new} и согласованность скалярных произведений в картановских
подалгебрах, см. замечание \ref{long_root}.)

Доказательство предложения \ref{new_relations} завершается ссылкой на
равенство
$$\left(\check{R}_{N(\mathfrak{q}^+,\lambda)^*
N(\mathfrak{q}^+,\mu)^*}\right)^*\;=\;
\left(\check{R}_{N(\mathfrak{q}^+,\lambda)
N(\mathfrak{q}^+,\mu)}\right)^{-1}, $$
см. \itemiii\ \ref{RM},
 и на следующие вспомогательные результаты, которые
будут получены в
 \itemiiiе \ref{diff_univ}. Во-первых, образы морфизмов $U_q\mathfrak{g}$-модулей
$$
N(\mathfrak{q}^+,w\rho-\rho) \rightarrow N(\mathfrak{q}^+,-\alpha_{l_0})
\otimes N(\mathfrak{q}^+,-\alpha_{l_0}), \qquad w \in\, W^{\mathbb{S}}\,
\&\, l(w)=2
$$
пересекаются с $N(\mathfrak{q}^+,-\alpha_{l_0})_{-1}\otimes
N(\mathfrak{q}^+,-\alpha_{l_0})_{-1}$ по подпространствам, сумма которых
имеет размерность $\frac{\dim \mathfrak{p}^-(\dim \mathfrak{p}^--1)}{2}$.
Во-вторых, на образе любого такого морфизма линейный оператор
$\check{R}_{N(\mathfrak{q}^+,-\alpha_{l_0}),N(\mathfrak{q}^+,-\alpha_{l_0})}
$ равен -1. \hfill $\square$

\bigskip

\begin{remark}\label{monomial_bases}\label{concrete_bases}
Рассмотрим квадратичную алгебру, определяемую пространствами образующих
$\mathbb{C}[\mathfrak{p}^-]_{q,1}$ и соотношений
$\mathcal{L}\subset\mathbb{C}[\mathfrak{p}^-]_{q,1}\otimes
\mathbb{C}[\mathfrak{p}^-]_{q,1}$. Это фактор-алгебра
$F=T(\mathbb{C}[\mathfrak{p}^-]_{q,1})/J$ тензорной алгебры
$T(\mathbb{C}[\mathfrak{p}^-]_{q,1})$ по ее двустороннему идеалу $J$,
порожденному $\mathcal{L}\subset T(\mathbb{C}[\mathfrak{p}^-]_{q,1})$.
Получим мономиальный базис векторного пространства $F$. Как и в
классическом случае $q=1$, веса $U_q\mathfrak{k}$-модуля
$\mathbb{C}[\mathfrak{p}^-]_{q,1}$ имеют вид $-\alpha_{l_0}-\sum_{i\ne
l_0}n_i\alpha_i$ и все весовые подпространства одномерны. Линейно
упорядочим множество весов этого $U_q\mathfrak{k}$-модуля в соответствии с
лексикографической упорядоченностью множества строк
$(-n_1,\;-n_2,\;\cdots,\;-n_{l_0-1},\;-1,\;-n_{l_0+1},\;\cdots,\;-n_l)$ из
коэффициентов разложения по простым корням. Выберем базис
$\{z_1,z_2,\cdots,z_{\dim\mathfrak{p}^-}\}$ весовых векторов пространства
$\mathbb{C}[\mathfrak{p}^-]_{q,1}$ так, как при доказательстве леммы
\ref{eigenvectors_in_L}, и упорядочим его по возрастанию весов. Очевидно,
свободная алгебра $\mathbb{C}\langle
z_1,z_2,\cdots,z_{\dim\mathfrak{p}^-}\rangle$ изоморфна тензорной алгебре
$T(\mathbb{C}[\mathfrak{p}^-]_{q,1})$ и из \eqref{Rmatrix} следует, что
тензоры
\begin{equation*}
z_i\otimes z_j+q^{-\frac4{(H_\mathbb{S},H_\mathbb{S})}}\sum\limits_{k<m}
\widetilde{R}_{ij}^{km}z_k\otimes z_m,\qquad i>j
\end{equation*}
образуют базис Гребнера идеала $J$ (см. \itemiii\ \ref{Grobner}). Значит,
множество
\begin{equation}\label{b_p}
\{z_1^{j_1}z_2^{j_2}z_3^{j_3}\ldots
z_{\mathrm{dim}\,\mathfrak{p}^-}^{j_{\mathrm{dim}\,\mathfrak{p}^-}}
\;|\;j_1<j_2<\ldots<j_{\mathrm{dim}\,\mathfrak{p}^-}\}
\end{equation}
является базисом векторного пространства $F$. Действие образующих $E_i$,
$F_i$, $K_i^{\pm 1}$ в этом базисе описывается матрицами, элементы которых
принадлежат полю рациональных функций $\mathbb{Q}(q)$ и не имеют полюсов
при $q\in(0,1]$. Это следует из определений и из предложения \ref{Levend}.
\end{remark}

\bigskip

Как вытекает из следующего утверждения, та же формула \eqref{b_p} задает
базис векторного пространства $\mathbb{C}[\mathfrak{p}^-]_q$.

   \begin{proposition}\label{quadr_alg} Алгебра $\mathbb{C}[\mathfrak{p}^-]_q$
 является квадратичной с пространством образующих $\mathbb{C}[\mathfrak{p}^-]_{q,1}$
 и пространством соотношений $\mathcal{L}$.
\end{proposition}

{\bf Доказательство.} Как и в замечании \ref{monomial_bases}, рассмотрим
квадратичную алгебру $F$ с пространством образующих
$\mathbb{C}[\mathfrak{p}^-]_{q,1}$ и пространством соотношений
$\mathcal{L}$. Требуется доказать, что естественный гомоморфизм
градуированных алгебр $F\rightarrow \mathbb{C}[\mathfrak{p}^-]_q$
биективен. Сюръективность вытекает из того, что подпространство
$\mathbb{C}[\mathfrak{p}^-]_{q,1}$ порождает алгебру
$\mathbb{C}[\mathfrak{p}^-]_q$, см. предложение \ref{hom_gen}. Остается
воспользоваться равенствами
$$ \dim \mathbb{C}[\mathfrak{p}^-]_{q,j}=\binom{j+\dim \mathfrak{p}^--1}{\dim
\mathfrak{p}^--1}=\dim F_j,$$
 последнее из которых легко получить с помощью базиса \eqref{b_p}. \hfill $\square$

\begin{remark}\label{HeckenbergerKolb} Предложение \ref{hom_gen} и
уточняющее его предложение \ref{quadr_alg} принадлежат
 Геккенберже и Колбу \cite{HeckenbergerKolb04}.
 Их доказательства существенно  отличаются от наших.
\end{remark}


\subsubsection{Дополнение о предоднородных векторных пространствах.}
\label{prehomogeneous}

Следующие определения и результаты подробно изложены в монографии
Рубентхалера \cite{Rub}.

Пусть $V$ -- пространство конечномерного рационального представления $\pi$
связной комплексной аффинной алгебраической группы. Его называют
предоднородным векторным пространством, если существует открытая орбита
$\Omega$. \IND{предоднородное векторное пространство} Она, очевидно,
единственна. Если представление $\pi$ неприводимо, то и предоднородное
векторное пространство $V$ называют неприводимым. \IND{предоднородное
векторное пространство ! неприводимое} Полный список неприводимых
предоднородных векторных пространств, рассматриваемых с точностью до
естественного изоморфизма, получен Сато и Кимурой \cite[стр. 228 --
229]{Kimura}. В дальнейшем будем предполагать рассматриваемую группу
редуктивной.

Предоднородное векторное пространство $V$ называют регулярным, если
замкнутое по Зарисскому множество $V\setminus\Omega$ является
гиперповерхностью. \IND{предоднородное векторное пространство ! регулярное}

Введем часто используемые классы неприводимых предоднородных векторных
пространств. Пусть $\mathfrak{g}$ -- простая комплексная алгебра Ли,
$\{\alpha_1,\alpha_2,\ldots,\alpha_l\}$ -- система ее простых корней, и
$\mathbb{S}\subset\{1,2,\ldots,l\}$.

Рассмотрим односвязную комплексную аффинную алгебраическую группу $G$ с
алгеброй Ли $\mathfrak{g}$. Ее координатное кольцо является алгеброй Хопфа,
порожденной матричными элементами фундаментальных представлений. \cite[стр.
94]{VinbGorbOn}, \cite[стр. 204]{VinbOn}.

В \itemiiiе \ref{GVerma} алгебра Ли была наделена градуировкой
$\mathfrak{g}=\newoplus_i\mathfrak{g}_i$. Пусть $L$ -- подгруппа всех
элементов, сохраняющих введенную градуировку, $L=\{a\in
G|\:\operatorname{Ad}_a\,H_\mathbb{S}=H_\mathbb{S}\}$, где использовано
присоединенное действие $L$ в $\mathfrak{g}$ \cite[стр. 114--116]{Hum}. Эта
подгруппа связна \cite[стр. 227--228]{VinbOn}, \cite[стр. 33]{Hum_conj}.

Как следует из теоремы Винберга \cite[стр. 119]{Rub}, \cite{Vinberg},
векторные пространства $\mathfrak{g}_{\pm 1}$ являются предоднородными
векторными пространствами с группой симметрии $L$. Эти
$\mathfrak{g}_0$-модули просты если и только если подмножество
$\{1,2,\cdots,l\}\setminus\mathbb{S}$ состоит из одного элемента \cite[стр.
126]{Rub}, \cite[стр. 113]{VinbGorbOn}.

Следуя Рубентхалеру, выделим важный класс таких предоднородных векторных
пространств. Предположим, что
\hbox{$\mathbb{S}=\{1,2,\ldots,l\}\setminus\{l_0\}$} и что простой корень
$\alpha_{l_0}$ входит в разложение максимального корня
\begin{equation}\label{max_root}
\delta=\sum_jn_j\alpha_j
\end{equation}
с коэффициентом $n_{l_0}$, равным 1. В этом случае градуировка алгебры Ли
$\mathfrak{g}$ быстро обрывается
$$
\mathfrak{g}=\mathfrak{g}_{-1}\oplus\mathfrak{g}_0\oplus\mathfrak{g}_{+1},
$$
и её однородные компоненты $\mathfrak{g}_{\pm 1}$ являются абелевыми
алгебрами Ли. Их называют предоднородными векторными пространствами
коммутативного параболического типа \cite[стр. 147]{Rub}. Полный список
таких предоднородных векторных пространств легко получить, используя явный
вид разложений \eqref{max_root} для неприводимых систем корней
\cite{Bou4-6}. В обозначениях Бурбаки \cite{Bou4-6} этот список выглядит
следующим образом:
\begin{itemize}
\item[$A_l$] \ \ $l_0\le 1,2,\cdots,l$, комплексные $l_0\times
    (l+1-l_0)$-матрицы;
\item[$C_l$]\ \ $l_0=l$,\ \ комплексные симметрические $l\times
    l$-матрицы;
\item[$D_l$]\ \ $l_0=l$,\ \ комплексные антисимметрические $l\times
    l$-матрицы;
\item[$B_l$]\ \ $l_0=1$;\ \ $D_l$\ \ $l_0=1$;
\item[$E_6$]\ \ $l_0=6$;\ \ $E_7$\ \ $l_0=7$
\end{itemize}
(нетрудно перечислить пары изоморфных предоднородных векторных пространств
этого списка).

Нерегулярными являются следующие предоднородные векторные
  пространства \cite{Muller}:
$$A_l, \; \text{при}\;l_0 \neq \frac{l+1}{2}; \quad \qquad D_l,\;  \text{при нечетных}\; l;
\quad \qquad E_6.$$


\IND{неприводимая ограниченная симметрическая область ! картановский
список}

Элементы этого списка часто изображают в виде диаграммы Дынкина с
выделенной вершиной, отвечающей простому корню $\alpha_{l_0}$.

\newpage
\bigskip

\centerline{ \sl Cписок предоднородных векторных пространств}
\centerline{\sl коммутативного параболического типа}

\bigskip

 \unitlength=1mm \linethickness{.5mm}

\begin{center}
\begin{picture}(111,160)(0,0)

\multiput(3,25)(15,0){6}{\circle*{2}} \put(78,25){\circle{4}}
\put(33,10){\circle*{2}} \put(3,25){\line(1,0){75}}
\put(33,25){\line(0,-1){15}}

\multiput(3,50)(15,0){5}{\circle*{2}} \put(3,50){\circle{4}}
\put(33,35){\circle*{2}} \put(3,50){\line(1,0){60}}
\put(33,50){\line(0,-1){15}}

\multiput(3,70)(15,0){7}{\circle*{2}}
\multiput(108,60)(0,20){2}{\circle*{2}} \put(3,70){\circle{4}}
\dashline[70]{3}(63,70)(78,70) \put(3,70){\line(1,0){60}}
\put(78,70){\line(1,0){15}} \put(93,70){\thicklines \line(3,2){15}}
\put(93,70){\thicklines \line(3,-2){15}}

\multiput(3,90)(15,0){5}{\circle*{2}} \multiput(93,90)(15,0){2}{\circle*{2}}
\put(3,90){\circle{4}} \put(93.5,89.4){\framebox(14,1.3){}}
\dashline[70]{3}(63,90)(93,90) \put(3,90){\line(1,0){60}}
\put(100,89.1){$\big >$}

\multiput(3,110)(15,0){7}{\circle*{2}}
\multiput(108,100)(0,20){2}{\circle*{2}} \put(108,120){\circle{4}}
\dashline[70]{3}(63,110)(78,110) \put(3,110){\line(1,0){60}}
\put(78,110){\line(1,0){15}} \put(93,110){\thicklines \line(3,2){15}}
\put(93,110){\thicklines \line(3,-2){15}}

\multiput(3,135)(15,0){5}{\circle*{2}}
\multiput(93,135)(15,0){2}{\circle*{2}} \put(108,135){\circle{4}}
\put(93.5,134.4){\framebox(14,1.3){}} \dashline[70]{3}(63,135)(93,135)
\put(3,135){\line(1,0){60}} \put(100,134.1){$\big <$}

\multiput(3,150)(15,0){8}{\circle*{2}} \put(63,150){\circle{4}}
\put(3,150){\line(1,0){30}} \put(48,150){\line(1,0){30}}
\put(93,150){\line(1,0){15}} \dashline[70]{3}(33,150)(48,150)
\dashline[70]{3}(78,150)(93,150)


\end{picture}
\end{center}

Каждая неприводимая ограниченная симметрическая область \cite[стр.
342]{Helg} допускает стандартную реализацию (реализацию Хариш-Чандры)
\IND{неприводимая ограниченная симметрическая область ! стандартная
реализация (реализация Хариш-Чандры)} в виде центрально-симметричной
выпуклой области предоднородного векторного пространства коммутативного
параболического типа. Это предоднородное векторное пространство регулярно,
если и только если рассматриваемая ограниченная симметрическая область
имеет трубчатый тип, то есть $w_0\alpha_{_{l_0}}=-\alpha_{_{l_0}}$, см.
\cite[стр. 38]{Arazy}. \IND{предоднородное векторное пространство !
регулярное}

\subsubsection{Дополнение о  некоммутативных базисах
Гребнера.}\label{Grobner}

Изложим основные понятия теории базисов Гребнера для некоммутативных
алгебр, следуя обзору В.~Уфнаровского \cite[стр. 31-34]{Ufn}. Отметим, что
для перехода от некоммутативного случая к существенно более простому
коммутативному случаю \cite{KoLiShi, Arjantsev, Prasolov} достаточно внести
в определение базиса Гребнера небольшие изменения, перечисленные в \cite[
стр. 35,36]{Ufn}.

Выберем ненулевое число $q \in \mathbb{C}$. Основным полем $K$ в этом
\itemiiiе будет наименьшее подполе поля $\mathbb{C}$, содержащее
 $q$. Таким образом, поле $K$ будет либо конечным расширением поля
рациональных чисел, либо его трансцендентным расширением степени 1.

Рассмотрим {\bf конечное} множество $X$, называемое в дальнейшем алфавитом.
\IND{алфавит} Пусть $\langle X \rangle$ -- полугруппа слов, включающая
пустое слово $\varnothing$, c операцией дописывания слов друг за другом.
Формальные линейные комбинации элементов $\langle X \rangle$ образуют
свободную алгебру $K\langle X\rangle$ с единицей $1\cdot\varnothing$.
Каждому {\bf конечному} подмножеству $R\subset K\langle X\rangle$
сопоставим наименьший содержащий его двусторонний идеал $I$ алгебры
$K\langle X\rangle$ и фактор-алгебру
\begin{equation}\label{finite_def}
A=K \langle X \rangle/I
\end{equation}
Алгебры, изоморфные \eqref{finite_def}, называют конечно определенными, а
множества $X$ и $R$ --- множествами образующих и соотношений
соответственно.

\begin{example}\label{S_N}(\cite[стр. 97]{CoxMos}).
Групповая алгебра симметрической группы $S_N$ допускает описание с помощью
образующих $s_i=(i,i+1)$, $i=1,2,\ldots,N-1$, и определяющих соотношений\ \
$s_i^2=1,\quad i=1,2,\ldots,N-1$,
$$(s_i s_j)^2=1,\qquad|i-j|\ne 1,\;\;i,j\in\{1,2,\ldots,N-1\},$$
$$(s_j s_{j+1})^3=1,\qquad j\in\{1,2,\ldots,N-2\}.$$
\end{example}

\bigskip

В дальнейшем алфавит $X$ предполагается линейно упорядоченным, а полугруппа
слов $\langle X\rangle$ -- наделенной отношением порядка
$\operatorname{deglex}$: слова сравниваются по длине, а, если их длины
равны, то лексикографически. \IND{отношение порядка
$\operatorname{deglex}$} Из определения следует, что множество слов
$\langle X\rangle$ является вполне упорядоченным и отношение порядка
согласуется с полугрупповой операцией: если $a_1\ge a_2\ \&\ b_1\ge b_2$,
то $a_1b_1\ge a_2b_2$. Очевидно, $\varnothing$ является наименьшим
элементом полугруппы слов.

Старшим словом элемента \IND{старшее слово}
$$
v=\sum a_w,\qquad w\in K\langle X\rangle,\; a_w\in K,\;w\in\langle X\rangle,
$$
называют наибольшее из слов $w$, входящих в это разложение с ненулевым
коэффициентом $a_w$. Старшее слово элемента $v$ обозначают $L(v)$, а
множество старших слов элементов двустороннего идеала $I$ обозначают
$L(I)$.

Воспользуемся введенным отношением порядка в $\langle X\rangle$ для
построения подпространства, дополняющего $I$ в векторном пространстве
$K\langle X\rangle$. Пусть $N$ -- множество слов, не принадлежащих $L(I)$.
Его называют множеством нормальных слов. \IND{множество нормальных слов}

\begin{proposition}\label{normal_form}
$K \langle X \rangle = KN \oplus I,$
 где $KN$ -- подпространство векторного
пространства $K \langle X \rangle$, порожденное множеством $N$.
\end{proposition}

{\bf Доказательство.} Очевидно, $KN\cap I=0$. Пусть $KN\oplus I\ne K\langle
X\rangle$. Рассмотрим наименьшее слово $w$, которое не принадлежит
$KN\oplus I$. Из $w\notin N$ вытекает существование такого элемента $v\in
I$, что его старшим словом является $w$. Иными словами, $v=aw+v'$, $a\in
K$, и каждое из слов, входящих в $v'$, принадлежит $KN\oplus I$. Значит,
$v'\in KN\oplus I$. Следовательно, $w\in KN\oplus I$, что противоречит
выбору слова $w$. \hfill $\square$

\medskip

Рассмотрим проектор в векторном пространстве $K\langle X\rangle$ на
подпространство $KN$ параллельно $I$:
$$K\langle X\rangle\to KN,\qquad v\to\overline{v}.$$
Его называют оператором приведения к нормальной форме, а элемент
$\overline{v}$ --- нормальной формой элемента $v$. \IND{нормальная форма
элемента $v$} Из предложения \ref{normal_form} следует, что этот оператор
порождает изоморфизм векторных пространств $A=K\langle X\rangle/I$ и $KN$,
полезный, если множество нормальных слов $N$ описано явно.

Предлагается совместить поиск явного описания этого множества с поиском
алгоритма приведения к нормальной форме. Следует отметить, что для
некоторых пар $X,\,R$ такого алгоритма не существует. Примером может
служить групповая алгебра $K[G]$ группы
 $G$ с конечным числом образующих и
определяющих соотношений, для которой знаменитая проблема равенства слов
алгоритмически неразрешима \cite[стр. 103]{Man1}.

Слово $w'$ называют подсловом \IND{подслово} слова $w$, если $w=aw'b$ при
некоторых $a,b\in\langle X\rangle$. Подслово $w'$ является собственным,
\IND{подслово ! собственное} если $w'\ne w$. Множество $M$ называют
двусторонним идеалом полугруппы $S$, \IND{двусторонний идеал полугруппы}
если $SM\subset M$ и \hbox{$MS\subset M$.}

Среди подмножеств, порождающих двусторонний идеал $L(I)$ полугруппы
$\langle X\rangle$, есть наименьшее подмножество $O$. Оно образовано
словами из $L(I)$, не имеющими собственных подслов из $L(I)$. Следуя
Д.~Анику, такие слова называют обструкциями. \IND{обструкция} Так как $O$
порождает двусторонний идеал $L(I)$ полугруппы $\langle X\rangle$, то
множество $G=\{w-\overline{w}\,|\,w\in O\}$ порождает двусторонний идеал
$I$ алгебры $K\langle X\rangle$.

Полученное множество образующих $G$ двустороннего идеала $I$ называют
базисом Гребнера этого идеала \IND{базис Гребнера идеала} (в обзоре
\cite{Ufn} используется термин ``приведенный базис Гребнера''). Разумеется,
базис Гребнера зависит от выбранного отношения линейного порядка в $X$.

\begin{lemma}\label{leader_term}
Обструкция $w \in O$ является старшим словом отвечающего ей элемента базиса
Гребнера: $w = L(w-\overline{w})$.
\end{lemma}

{\bf Доказательство.} Очевидно, $N\cap L(I)=\varnothing$, $N\cup
L(I)=\langle X\rangle$. Следовательно, из $v=v'+v''$, $v'\in KN$, $v''\in
I$ вытекает равенство
$$
L(v)=L(v'+v'')=L(L(v')+L(v''))=\left\{
\begin{array}{rl}
L(v'), & \;\mbox{если }L(v)\in N,
\\ L(v''), & \;\mbox{если }L(v)\in L(I).
\end{array}
\right.
$$
Остается положить $v=w$, $v'=\overline{w}$, $v''=w-\overline{w}$ и
воспользоваться тем, что $w= L(w)\in L(I)$. \hfill $\square$

\medskip

Предположим, что базис Гребнера конечен. Это требование существенно в
некоммутативном случае. Покажем, что существует алгоритм приведения
элемента $v \in K \langle X \rangle$ к нормальной форме и опишем этот
алгоритм.

   Сопоставим каждому элементу $u \in I$ линейный оператор
\begin{equation}\label{u_hat}
\hat{u}: K \langle X \rangle \rightarrow K \langle X \rangle.
\end{equation}
Пусть $u \in I$ и $a_{L(u)} \in K$ -- коэффициент при старшем слове $L(u)$ в
разложении $u=\sum_{w \in \langle X \rangle}a_w w$. Введем линейный оператор
\eqref{u_hat}, задав его действие на элементы базиса $\langle X \rangle$
векторного пространства $K \langle X \rangle$. По определению,
$\hat{u}(w)=w$, если $L(u)$ не является одним из подслов слова $w \in
\langle X \rangle$. В противном случае, можно найти первое вхождение
подслова $L(u)$ в слово $w$, и $\hat{u}(w)$ получается заменой этого
подслова выражением $L(u)-\frac{1}{a_{L(u)}}u$. Другими словами, оператор
$\hat{u}$ осуществляет подстановку
\begin{equation}\label{substitution}
  L(u) \mapsto L(u)-\frac{1}{a_{L(u)}}u.
\end{equation}

\medskip

Пусть $v\in K\langle X\rangle$. Элемент $v'\in K\langle X\rangle$
называется $G$-редукцией элемента $v$, \IND{$G$-редукция элемента $v$}
если, во-первых, $v'$ можно получить из $v$ в результате последовательного
применения операторов множества $\{\hat{u}|u\in G\}$ и, во-вторых, $v'$
является общим неподвижным вектором этих линейных операторов.

\begin{proposition}\label{reduction_Grobner}
$G$-редукция элемента $v \in K \langle X \rangle$ существует, единственна и
равна нормальной форме $\overline{v}$ этого элемента.
\end{proposition}

\begin{corollary}\label{normal}
Слово $w \in \langle X \rangle$ является нормальным, если и только если
среди его подслов нет старших слов элементов базиса Гребнера.
\end{corollary}

\begin{note}\label{comp_alg}
Переход от конечного списка правил подстановки (\ref{substitution}) к
алгоритму приведения к нормальной форме мы не обсуждаем, поскольку в
современных программных системах компьютерной алгебры имеются средства,
позволяющие автоматизировать этот переход.
\end{note}

\begin{example}\label{ex_Gr}\cite[стр. 33]{Ufn}
Рассмотрим алгебру $A$ с образующими $x$, $y$ и определяющим соотношением
$x^2 + y^2 =0$: $A = K \langle x,y \rangle / (x^2 +y^2)$. Упорядочим
образующие: $x>y$. Можно показать, что в этом случае $G=\{x^2+y^2,
xy^2-y^2x\} $ и список правил подстановки имеет вид
$$
x^2 \mapsto -y^2,\quad xy^2 \mapsto y^2x.
$$
Слово $w$ является нормальным, если и только если оно не имеет подслов
$x^2$, $xy^2$:
$
N = \{y^n(xy)^k,\quad y^n(xy)^k x\ |\ n,k \in \mathbb{Z}_+\}.
$
\end{example}

  Итак, имея конечный базис Гребнера, легко получить алгоритм приведения
к нормальной форме. Остается объяснить, как найти конечный базис Гребнера
двустороннего идеала $I=(R)$ в предположении, что он существует. Для этого
предлагается трансформировать множество $R$
$$
R \rightarrow R_1 \rightarrow R_2 \rightarrow \ldots,
$$
используя для перехода от $R_j$ к $R_{j+1}$ операции трех описанных ниже
типов.

\begin{itemize}
\item[ 1]. (Нормировка.) Предположим, что $v \in R_j$ и старшее слово $L(v)$
входит в разложение $v = \sum a_w w$ с коэффициентом $a_{L(v)} \neq 1$.
Нормировка состоит в замене элемента $v$ элементом $a_{L(v)}^{-1}v$,
пропорциональным $v$. \\
\item[2.] (Редукция.) Пусть $v,v' \in R_j$, $v \neq v'$ и $L(v')$ является
подсловом слова $L(v)$: $L(v) = a L(v')b$, $a,b \in \langle X \rangle$.
Редукция состоит в удалении $v$ из множества $R_j$, если $v=av'b$, и в
замене элемента $v$ элементом $v-av'b$, если $v \neq av'b$. \\
\item[3.] Рассмотрим два элемента $u,v \in R_j$ с равными единице
коэффициентами при старших мономах $L(u)$, $L(v)$. Слово $w \in \langle X
\rangle$ называется композицией элементов $u$, $v$, если $L(u)$ является
началом слова $w$, $L(v)$ -- его концом, и эти подслова слова $w$
пересекаются:
$$
w=xyz,\quad L(u)=xy,\quad L(v)=yz,
$$
где $y$ -- непустое слово. Наличие композиции позволяет применить к слову
$w$ как подстановку ${xy \mapsto L(u)-u}$, так и подстановку $yz \mapsto
L(v)-v$. Результаты, вообще говоря, не совпадают. В самом деле, их разность
равна $(L(u)-u)z - x(L(v)-v) = xv-uz$. Операция третьего типа состоит в
присоединении к $R_j$ элемента $xv-uz \in K \langle X \rangle$, если этот
элемент не содержится в $R_j$. Элемент $xv-uz$ называется результатом
композиции.
\end{itemize}

Порядок выполнения операций весьма существен. Операциям первого типа
присваивается высший приоритет: именно их нужно применять, если хотя бы у
одного из элементов множества $R_j$ коэффициент при старшем слове отличен от
$1$. Операциям третьего типа присваивается низший приоритет: такую операцию
выбирают только в том случае, когда операции двух первых типов неприменимы.

\begin{note}\label{cycle}
Может случиться так, что элемент, присоединяемый к $R_j$ операцией третьего
типа, ``редуцируется к нулю'' непосредственно следующими за ней операциями
первых двух типов. Это источник ``зацикливаний'', и для борьбы с ним следует
ввести дополнительный запрет примененять ''бесполезные'' операции третьего
типа.
\end{note}

Если описанный ``алгоритм'' прекратил работу после конечного числа шагов
из-за того, что к полученному множеству нельзя применить ни одну из
описанных выше операций (с учетом замечания \ref{cycle}), то это множество
является базисом Гребнера идеала $I$. Для доказательства достаточно
воспользоваться даймонд-леммой Бергмана \cite[стр. 181]{Bergman} в следующей
ее формулировке (ср. \cite[стр.~33]{Ufn}).

\begin{lemma}\label{diamond_compos}
Пусть $G\subset K\langle X\rangle$ -- множество, к которому неприменимы ни
операции редукции, ни операции нормировки. Если результат любой из
композиций приводится к нулю за конечное число последовательно выполняемых
шагов редукции и нормировки, то $G$ -- базис Гребнера.
\end{lemma}

\subsection{ Пространства функций в квантовых ограниченных симметрических
областях}\label{function_spaces}
\subsubsection{Фоковское представление и финитные функции.}\label{finite}

Рассмотрим комплексное векторное пространство $V\ne\{0\}$. Ограниченная
область $\mathbb{D}$ этого векторного пространства называется
симметрической, \IND{ограниченная симметрическая область} если каждая точка
$z\in\mathbb{D}$ является изолированной неподвижной точкой некоторого
инволютивного автоморфизма $s_z^2=\operatorname{id}$ области $\mathbb{D}$,
рассматриваемой как комплексное многообразие \cite[стр.~342]{Helg},
\cite[стр. 193]{Kor}. Такую область называют приводимой, \IND{ограниченная
симметрическая область ! приводимая} если в категории комплексных
многообразий она изоморфна декартовому произведению ограниченных
симметрических областей. Простейший пример неприводимой ограниченной
симметрической области -- единичный круг
$\mathbb{D}=\{z\in\mathbb{C}\,|\,|z|<1\}$.

Рассмотрим эрмитову симметрическую пару $(\mathfrak{g},\mathfrak{k})$ (см.
\itemiii\ \ref{Hermit}). Такие пары $(\mathfrak{g},\mathfrak{k})$ получаются
комплексификацией ортогональных инволютивных алгебр Ли \cite[стр.
217]{Helg}, \cite[стр. 196]{Kor}. Каждой из них отвечает разложение простой
комплексной алгебры Ли $\mathfrak{g}$:
$$\mathfrak{g} = \mathfrak{p}^- \oplus \mathfrak{k} \oplus \mathfrak{p}^+.
$$

Как показал Хариш-Чандра (\cite[стр. 342--352]{Helg}, \cite[стр.
6]{Satake}), каждая неприводимая ограниченная симметрическая область
изоморфна области $\mathbb{D}$, стандартным образом вложенной в
$\mathfrak{p}^-$. Наша цель -- получить ее $q$-аналог.

\bigskip

Как и в частном случае квантового круга, см. \itemiii\
\ref{inv_int_subsection}, нам нужна $U_q\mathfrak{g}$-модульная алгебра
$\mathscr{D}(\mathbb{D})_q$ финитных функций в квантовой ограниченной
симметрической области $\mathbb{D}$. Прежде всего, введем фоковское
представление \hbox{$*$-алгебры} $\mathrm{Pol}(\mathfrak{p}^-)_q=
(\mathbb{C}[\mathfrak{p}^-\oplus\mathfrak{p}^+]_q, *)$.

Рассмотрим $\mathbb{C}[\mathfrak{p}^- \oplus \mathfrak{p}^+]$-модуль
$\mathcal{H}$ с одной образующей $v_0$ и определяющими соотношениями
$$
f v_0 = \langle f \rangle v_0,\qquad f \in \mathbb{C}[\mathfrak{p}^+]_q,
$$
где $\langle f \rangle$ -- линейный функционал, введенный в разделе
\ref{Pol}. Другими словами,
\begin{equation}\label{Fock1}
f v_0 = 0,\qquad f \in \newoplus\limits_{j=1}^\infty
\mathbb{C}[\mathfrak{p}^+]_{q,-j}.
\end{equation}
Из предположения \ref{hom_gen} следует, что множество определяющих соотношений
(\ref{Fock1}) можно заменить его подмножеством
\begin{equation}\label{Fock2}
  f v_0 = 0,\qquad f \in \mathbb{C}[\mathfrak{p}^+]_{q,-1}.
\end{equation}

\begin{lemma}\label{H_as_vect}
Линейное отображение
\begin{equation}\label{iso_vect}
  \mathbb{C}[\mathfrak{p}^-]_q \rightarrow \mathcal{H},\quad f\mapsto f v_0
\end{equation}
является изоморфизмом векторных пространств.
\end{lemma}

{\bf Доказательство.} Ядро линейного отображения $\mathbb{C}[\mathfrak{p}^-
\oplus \mathfrak{p}^+]_q \rightarrow \mathcal{H}$, $f \mapsto f v_0$ -- это
левый идеал алгебры $\mathbb{C}[\mathfrak{p}^- \oplus \mathfrak{p}^+]_q$,
порожденный $\bigoplus\limits_{j=1}^\infty
\mathbb{C}[\mathfrak{p}^+]_{q,-j}$. Этот идеал совпадает с
$\bigoplus\limits_{i=0}^\infty \bigoplus\limits_{j=1}^\infty
\mathbb{C}[\mathfrak{p}^-]_{q,i} \otimes \mathbb{C}[\mathfrak{p}^+]_{q.-j}$.
Значит, $\mathcal{H} \cong \bigoplus\limits_{i=0}^\infty
\mathbb{C}[\mathfrak{p}^-]_{q,i} = \mathbb{C}[\mathfrak{p}^-]_q$. \hfill
$\square$

\medskip

Из этой леммы и предложения \ref{positive_herm} вытекает
\begin{proposition}\label{preHilbert}
Эрмитова форма
$$
(f_1v_0,f_2v_0)=\langle f_2^*f_1\rangle,\qquad f_1,f_2\in
\mathrm{Pol}(\mathfrak{p}^-)_q,
$$
в $\mathcal{H}$ корректно определена, положительна и
$$
(fv',v'')=(v',f^*v''),\qquad v',v''\in\mathcal{H},\;f\in
\mathrm{Pol}(\mathfrak{p}^-)_q.
$$
\end{proposition}

\bigskip

Последнее предложение позволяет ввести в рассмотрение представление $T_F$
$*$-алгебры $\mathrm{Pol}(\mathfrak{p}^-)_q$ в предгильбертовом
пространстве $\mathcal{H}$. Его называют фоковским представлением, а вектор
$v_0$ -- вакуумным вектором. \IND{фоковское представление ! $*$-алгебры
$\mathrm{Pol}(\mathfrak{p}^-)_q$} \IND{вакуумный вектор}

\medskip

Обратимся к частному случаю $\mathfrak{g}=\mathfrak{sl}_2$. Из
\eqref{comm_disc} следует, что
$$
(z^nv_0,z^nv_0)=((1-q^2y)z^{n-1}v_0,z^{n-1}v_0)=
(1-q^{2n})(z^{n-1}v_0,z^{n-1}v_0)=(q^2;q^2)_n.
$$
Значит, действие операторов $T_F(z)$, $T_F(z^*)$ на элементы
ортонормированного базиса
$\left\{\sqrt{(q^2;q^2)_n}\;z^n\,v_0\right\}_{n\in\mathbb{Z}_+}$
описывается равенствами \eqref{T_F}.

\medskip

\begin{note}\label{irrep_Fock}
Очевидно, фоковское представление неприводимо. Действительно, пусть
$\mathcal{H}_1\subset\mathcal{H}$ -- ненулевое общее инвариантное
подпространство операторов этого представления. Тогда $\mathcal{H}_1$
содержит ненулевой вектор $v$, для которого $T_F(f)v=0$ при всех
$f\in\mathbb{C}[\mathfrak{p}^+]_{q,-j}$, $j\ne 0$. Этот вектор ортогонален
$\bigoplus\limits_{j=0}^\infty\mathbb{C}[\mathfrak{p}^-]_{q,j}v_0$ и,
следовательно, принадлежит $\mathbb{C}v_0$. Значит, $v_0\in\mathcal{H}_1$ и
$\mathcal{H}_1=\mathcal{H}$.
\end{note}

\medskip

Перейдем к построению алгебры $\mathscr{D}(\mathbb{D})_q$. В частном случае
квантового круга эта алгебра введена в \itemiiiе \ref{finite_sl_2}. Наделим
векторное пространство $\mathcal{H}$ градуировкой
$$
\mathcal{H}=\newoplus\limits_{j=0}^\infty\mathcal{H}_j,\qquad
\mathcal{H}_j=\mathbb{C}[\mathfrak{p}^-]_{q,j}v_0.
$$
Подпространства $\mathcal{H}_j\cong\mathbb{C}[\mathfrak{p}^-]_{q,j}$
попарно ортогональны и конечномерны. Из \eqref{dim_j} следует, что
\begin{equation}\label{dim_H_j}
\dim\mathcal{H}_j=\dbinom{j+\dim\mathfrak{p}^- -1}{\dim\mathfrak{p}^--1}.
\end{equation}
Линейный оператор $A$ в $\mathcal{H}$ назовем финитным, \IND{финитный
линейный оператор} если $A\mathcal{H}_j=0$ для всех $j\in\mathbb{Z}_+$,
кроме конечного их числа. Примером может служить одномерный ортогональный
проектор $P_0$ на вакуумное подпространство $\mathcal{H}_0=\mathbb{C}v_0$.
Финитные линейные операторы образуют алгебру, и определяемая ниже алгебра
$\mathscr{D}(\mathbb{D})_q$ окажется ей канонически изоморфной. Введем в
рассмотрение элемент $f_0$, которому при этом каноническом изоморфизме
будет отвечать проектор $P_0$.

Точнее, рассмотрим алгебру с единицей $\mathrm{Fun}(\mathbb{D})_q$,
порожденную своей подалгеброй $\mathrm{Pol}(\mathfrak{p}^-)_q$ и элементом
$f_0$, удовлетворяющими следующим определяющим соотношениям
\begin{equation}\label{rel1}
f_0^2=f_0,
\end{equation}
\begin{equation}\label{rel2}
\psi f_0=\langle\psi\rangle f_0,\quad\psi\in\mathbb{C}[\mathfrak{p}^+]_q;
\qquad f_0\psi=\langle\psi\rangle f_0,\quad
\psi\in\mathbb{C}[\mathfrak{p}^-]_q.
\end{equation}
Из предложения \ref{hom_gen} следует, что множество соотношений
(\ref{rel2}) в этом определении можно заменить его подмножеством
\begin{equation}\label{rel3}
\psi f_0=0,\quad\psi\in\mathbb{C}[\mathfrak{p}^+]_{q,-1};\qquad f_0\psi=0,
\quad\psi\in\mathbb{C}[\mathfrak{p}^-]_{q,1}.
\end{equation}

\medskip

Например, в случае квантового круга элементы $z$, $z^*$, $f_0$ порождают
алгебру $\mathrm{Fun}(\mathbb{D})_q$, а множество определяющих соотношений
включает (\ref{Pol_l}),
\begin{equation}\label{local_060214}
z^*\, f_0 = f_0\, z = 0,\quad f_0^2 = f_0.
\end{equation}

\medskip
Существует и единственно такое продолжение инволюции $*$ с алгебры
$\mathrm{Pol}(\mathfrak{p}^-)_q$ на алгебру $\mathrm{Fun}(\mathbb{D})_q$,
что
 $f_0^* = f_0$. Очевидно, имеет место
\begin{lemma}\label{extension_Fock} Существует и единственно  продолжение фоковского
представления $T_F$ до представления $*$-алгебры
$\mathrm{Fun}(\mathbb{D})_q$, при котором элементу $f_0$ отвечает
ортогональный проектор $P_0$ на вакуумное подпространство.
\end{lemma}

\medskip Для этого продолжения фоковского представления
сохраним прежнее обозначение $T_F$.

Двусторонний идеал алгебры $\mathrm{Fun}(\mathbb{D})_q$, порожденный
элементом $f_0$, будем обозначать $\mathscr{D}(\mathbb{D})_q$ и называть
алгеброй финитных функций в квантовой ограниченной симметрической области.
\IND{алгебра ! финитных функций ! в квантовой ограниченной симметрической
области}

Отметим, что алгебра финитных линейных операторов в $\mathcal{H}$ является
$*$-алгеброй (сопряженный линейный оператор вводится стандартным образом).

\begin{proposition}\label{bijection}  Представление $T_F$ доставляет
 изоморфизм $*$-ал\-гебры $\mathscr{D}(\mathbb{D})_q$ и $*$-алгебры всех финитных
линейных операторов в $\mathcal{H}$.
\end{proposition}

{\bf Доказательство.} $T_F$ -- $*$-представление. Следовательно, достаточно
доказать, что сужение $T_F$ на $\mathscr{D}(\mathbb{D})_q$ является
биективным отображением $\mathscr{D}(\mathbb{D})_q$ на алгебру всех
финитных линейных операторов в $\mathcal{H}$.

Пусть
$\mathscr{D}(\mathbb{D})_{q,i,j}=\mathbb{C}[\mathfrak{p}^-]_{q,i}\cdot
f_0\cdot\mathbb{C}[\mathfrak{p}^+]_{q,-j}$. Если
$f\in\mathscr{D}(\mathbb{D})_{q,i,j}$, то линейный оператор $T_F(f),\qquad
f\in\mathscr{D}(\mathbb{D})_{q,i,j}$, отображает $\mathcal{H}_j$ в
$\mathcal{H}_i$ и равен нулю на $\bigoplus\limits_{k\ne j}\mathcal{H}_k$.
Получаем линейное отображение из $\mathscr{D}(\mathbb{D})_{q,i,j}$ в
$\mathrm{Hom}(\mathcal{H}_j,\mathcal{H}_i)$. Оно сюръективно, в силу
предложения \ref{inner1}, и
$$
\dim\mathscr{D}(\mathbb{D})_{q,i,j}=
\dim\,\mathrm{Hom}(\mathcal{H}_j,\mathcal{H}_i).
$$
Значит, представление $T_F$ устанавливает изоморфизм
\begin{equation}\label{iso_finite}
\mathscr{D}(\mathbb{D})_{q,i,j}=\mathbb{C}[\mathfrak{p}^-]_{q,i}\;f_0\;
\mathbb{C}[\mathfrak{p}^+]_{q,-j}\cong
\mathrm{Hom}(\mathcal{H}_j,\mathcal{H}_i).
\end{equation}
Но $\mathscr{D}(\mathbb{D})_q=
\bigoplus\limits_{i,j=0}^\infty\mathscr{D}(\mathbb{D})_{q,i,j}$, а
векторное пространство финитных линейных операторов совпадает с
$\bigoplus\limits_{i,j=0}^\infty\mathrm{Hom}(\mathcal{H}_j,\mathcal{H}_i)
\hookrightarrow\mathrm{End}\,\mathcal{H}$. \hfill $\square$

\bigskip

Перейдем к доказательству точности фоковского представления алгебры
$\mathrm{Pol}(\mathbb{D})_q$. Так же, как предложение \ref{bijection},
доказывается
\begin{lemma}\label{ball_I_8.7}
При всех $j,k \in \mathbb{Z}_+$ отображение
$$
\mathbb{C}[\mathfrak{p}^-]_{q,k} \cdot \mathbb{C}[\mathfrak{p}^+]_{q,-j}
\rightarrow \mathrm{Hom}(\mathcal{H}_j, \mathcal{H}_k),\qquad f \mapsto
T_F(f)|_{\mathcal{H}_j}
$$
биективно.
\end{lemma}

\begin{proposition}\label{ball_I_8.8}
Фоковское представление алгебры $\mathrm{Pol}(\mathfrak{p}^-)_q$ является
точным.
\end{proposition}
{\bf Доказательство.} Наделим $\mathbb{Z}_+^2$ стандартным отношением
частичного порядка ($(k_1,j_1)\le(k_2, j_2)$, если и только если $k_1\le
k_2\ \&\ j_1\le j_2$), а векторное пространство
$\mathrm{Pol}(\mathfrak{p}^-)_q$ биградуировкой
$$
\mathrm{Pol}(\mathfrak{p}^-)_q=
\newoplus\limits_{k,j=0}^\infty\mathbb{C}[\mathfrak{p}^-]_{q,k}\;
\mathbb{C}[\mathfrak{p}^+]_{q,-j}.
$$
Предположим, что утверждение леммы неверно и $T_F(f)=0$, $f\ne 0$.
Воспользуемся разложением
$$
f=\sum\limits_{k,j=0}^\infty f_{k,j},\qquad f_{k,j}\in
\mathbb{C}[\mathfrak{p}^-]_{q,k}\mathbb{C}[\mathfrak{p}^+]_{q,-j},
$$
и выберем ненулевую компоненту $f_{k_0 j_0}$ с минимальным ''номером''
$(k_0,j_0)$ в этом разложении. Из $T_F(f)=0$ следует, что
$T_F(f_{k_0j_0})|_{\mathcal{H}_{j_0}}=0$, поскольку
$T_F(f_{k_0j_0})|_{\mathcal{H}_{j_0}}$ -- произведение линейного оператора
$T_F(f)|_{\mathcal{H}_{j_0}}$ и проектора на $\mathcal{H}_{k_0}$
параллельно $\bigoplus\limits_{k\ne k_0}\mathcal{H}_k$. Остается заметить,
что
$$
f_{k_0 j_0}\ne 0,\qquad f_{k_0j_0}\in\mathbb{C}[\mathfrak{p}^-]_{q,k_0}
\mathbb{C}[\mathfrak{p}^+]_{q,-j_0},\qquad
T_F(f_{k_0j_0})|_{\mathcal{H}_{j_0}}=0.
$$
Приходим к противоречию с утверждением леммы \ref{ball_I_8.7}. \hfill
$\square$

\medskip

\begin{remark}\label{Fun-faith}
Как следует из доказанных ниже предложения \ref{image} и следствия
\ref{Pol-integrity}, образ подалгебры $\operatorname{Pol}(\mathfrak{p}^-)_q$
в фоковском представлении не содержит ненулевых элементов образа подалгебры
$\mathscr{D}(\mathbb{D})_q$. Значит, фоковское представление алгебры
$\mathrm{Fun}(\mathbb{D})_q=\operatorname{Pol}(\mathfrak{p}^-)_q\oplus
\mathscr{D}(\mathbb{D})_q$ является точным.
\end{remark}

\subsubsection{Фоковское представление и обобщенные функции.}\label{gen}


Рассмотрим введенные в разделе \ref{finite} конечномерные подпространства
$\mathcal{H}_j\subset\mathcal{H}$, $j\in\mathbb{Z}_+$, и их прямое
произведение $\overline{\mathcal{H}}$ -- векторное пространство формальных
рядов с членами из $\mathcal{H}_j$. Отметим, что векторные пространства
$\mathrm{Hom}(\mathcal{H}_j,\mathcal{H}_k)$ конечномерны. Как и в
\itemiiiе \ref{finite_sl_2}, линейные операторы из $\mathcal{H}$ в
    $\overline{\mathcal{H}}$ будем называть обобщенными линейными
операторами в $\mathcal{H}$. \IND{обобщенные ! линейные операторы} Каждый
из них является формальным рядом
\begin{equation}\label{formal}
A=\sum\limits_{k,j=0}^\infty A_{kj},\quad
A_{kj}\in\mathrm{Hom}(\mathcal{H}_j,\mathcal{H}_k),
\end{equation}
и векторное пространство обобщенных линейных операторов в $\mathcal{H}$
наделяется топологией почленной сходимости формальных рядов.

Определяемое ниже топологическое векторное пространство
$\mathscr{D}(\mathbb{D})_q'$ ``обобщенных функций в квантовой ограниченной
симметрической области'' окажется канонически изоморфно топологическому
векторному пространству обобщенных линейных операторов в $\mathcal{H}$, как
и в частном случае квантового круга, см. \itemiii\ \ref{finite_sl_2}.

Напомним, что
$$
\mathrm{Pol}(\mathfrak{p}^-)_q=
\newoplus\limits_{k,j=0}^\infty\mathbb{C}[\mathfrak{p}^-]_{q,k}\cdot
\mathbb{C}[\mathfrak{p}^+]_{q,-j}\cong
\newoplus\limits_{k,j=0}^\infty\mathbb{C}[\mathfrak{p}^-]_{q,k}\otimes
\mathbb{C}[\mathfrak{p}^+]_{q,-j},
$$
причем $\dim(\mathbb{C}[\mathfrak{p}^-]_{q,k}\cdot
\mathbb{C}[\mathfrak{g}^+]_{q,-j})<\infty$. Другими словами, каждый элемент
$f\in\mathrm{Pol}(\mathfrak{p}^-)_q$ единственным образом разлагается в
конечную сумму
\begin{equation}\label{series_f}
f=\sum\limits_{k,j=0}^\infty f_{k,j},\qquad
f_{k,j}\in\mathbb{C}[\mathfrak{p}^-]_{q,k}\cdot
\mathbb{C}[\mathfrak{p}^+]_{q,-j},
\end{equation}
где $\mathbb{C}[\mathfrak{p}^-]_{q,k}\cdot
\mathbb{C}[\mathfrak{g}^+]_{q,-j}$-- линейная оболочка множества
 $\{f'f''\;|\; f'\in \mathbb{C}[\mathfrak{p}^-]_{q,k},\;
 f''\in \mathbb{C}[\mathfrak{g}^+]_{q,-j}\}$.

Рассмотрим векторное пространство $\mathscr{D}(\mathbb{D})_q'$ формальных
рядов вида (\ref{series_f}) с топологией почленной сходимости. Инволюция $*$
продолжается по непрерывности с плотного линейного подмногообразия
$\mathrm{Pol}(\mathfrak{p}^-)_q$ на $\mathscr{D}(\mathbb{D})_q'$
\begin{equation}\label{invol_'}
*:\sum\limits_{k,j=0}^\infty f_{k,j}\mapsto\sum\limits_{k,j=0}^\infty
f_{k,j}^*.
\end{equation}
Элементы пространства $\mathscr{D}(\mathbb{D})_q'$ будем называть
обобщенными функциями в квантовой ограниченной симметрической области.
\IND{обобщенные ! функции в квантовой ограниченной симметрической области}

\begin{proposition}\label{isom_general}
  1. Фоковское представление алгебры $\mathrm{Pol}(\mathfrak{p}^-)_q$
допускает продолжение по непрерывности до линейного отображения
$\mathscr{D}(\mathbb{D})_q'$ в пространство обобщенных линейных операторов в
$\mathcal{H}$.
\\ 2. Это линейное отображение является изоморфизмом топологических
векторных пространств.
\end{proposition}

{\bf Доказательство.}
Введем в рассмотрение линейные операторы
$$
a_{k,j}: \mathrm{Pol}(\mathfrak{p}^-)_q \rightarrow
\mathbb{C}[\mathfrak{p}^-]_{q,k} \cdot
\mathbb{C}[\mathfrak{p}^+]_{q,-j};\quad b_{kj}:
\mathrm{Pol}(\mathfrak{p}^-)_q \rightarrow \mathrm{Hom}(\mathcal{H}_j,
\mathcal{H}_k).
$$
Именно, операторы $a_{k,j}$ введем с помощью разложения (\ref{series_f}),
полагая $a_{kj} (f) = f_{k,j}$, а операторы $b_{kj}$ -- с помощью разложения
(\ref{formal}):
$$
T_F(f) = \sum\limits_{k,j=0}^\infty b_{kj}(f), \quad f \in
\mathrm{Pol}(\mathfrak{p}^-)_q.
$$
Здесь $k,j \in \mathbb{Z}_+$. Наделим векторное пространство
$\mathrm{Pol}(\mathfrak{p}^-)_q$ двумя топологиями. Пусть $\mathscr{T}_a$ --
слабейшая из топологий, в которых непрерывны все операторы $a_{k,j}$, а
$\mathscr{T}_b$ -- слабейшая из топологий, в которых непрерывны все
операторы $b_{k,j}$. Топологии $\mathscr{T}_a$, $\mathscr{T}_b$ индуцированы
вложениями $\mathrm{Pol}(\mathfrak{p}^-)_q$ в $\mathscr{D}(\mathbb{D})_q'$ и
в пространство обобщенных линейных операторов в $\mathcal{H}$
соответственно. Значит, эти топологические векторные пространства являются
пополнениями $\mathrm{Pol}(\mathfrak{p}^-)_q$ в топологиях $\mathscr{T}_a$,
$\mathscr{T}_b$, и нам достаточно доказать, что $\mathscr{T}_a =
\mathscr{T}_b$. Очевидно, $\mathrm{Ker}\,b_{kj} \supset \bigcap\limits_{i
\leq \min(k,j)} \mathrm{Ker}\,a_{k-i,j-i}$. Другими словами, оператор
$b_{kj}(f)$ является линейной комбинацией операторов $a_{k-i,j-i}(f)$ с
коэффициентами
$$
c_{kji}: \mathbb{C}[\mathfrak{p}^-]_{q,k-i}
\mathbb{C}[\mathfrak{g}^+]_{q,-j+i} \rightarrow \mathrm{Hom}(\mathcal{H}_j,
\mathcal{H}_k),
$$
 не зависящими от $f \in \mathrm{Pol}(\mathfrak{p}^-)_q$:
\begin{equation}\label{system}
  b_{kj}(f) = \sum\limits_{i \leq \min(k,j)} c_{k,j,i} \, a_{k-i,j-i} (f).
\end{equation}
Значит, топология $\mathscr{T}_a$ не слабее топологии $\mathscr{T}_b$.
остается заметить, что из леммы \ref{ball_I_8.7} следует обратимость
линейных операторов $c_{k,j,0}$, $k,j \in \mathbb{Z}_+$. Значит,
``треугольную'' систему уравнений (\ref{system}) можно решить относительно
$a_{k,j}(f)$,
и топология $\mathscr{T}_b$ не слабее топологии $\mathscr{T}_a$. \hfill
$\square$

\medskip Отметим, что векторное пространство обобщенных линейных операторов
наделено инволюцией $*$ и полученное выше линейное отображение
$\mathscr{D}(\mathbb{D})_q'$ в это векторное пространство уважает инволюции.

 Из этого предложения и из предложения \ref{bijection} вытекает

\begin{corollary}\label{embed_fun}
Вложение векторного пространства $\mathrm{Pol}(\mathfrak{p}^-)_q$ в
$\mathscr{D}(\mathbb{D})_q'$ канонически продолжается до линейного
отображения векторного пространства $\mathrm{Fun}(\mathbb{D})_q$ в
$\mathscr{D}(\mathbb{D})_q'$ и доставляет вложение
$\mathscr{D}(\mathbb{D})_q$ в $\mathscr{D}(\mathbb{D})_q'$ в виде плотного
линейного подмногообразия.
\end{corollary}

\medskip
Будем отождествлять $\mathscr{D}(\mathbb{D})_q,
\mathrm{Pol}(\mathfrak{p}^-)_q$ с их образами при вложении в
$\mathscr{D}(\mathbb{D})_q'$. Из предложения \ref{isom_general} вытекают

\begin{corollary}\label{fun_bimod}
 Для любых $f_1,f_2 \in \operatorname{Fun}(\mathbb{D})_q$ линейный
оператор
$$
\mathscr{D}(\mathbb{D})_q \rightarrow \mathscr{D}(\mathbb{D})_q, \qquad f
\mapsto f_1 \cdot f \cdot f_2
$$
допускает продолжение по непрерывности до линейного оператора в
$\mathscr{D}(\mathbb{D})_q'$. Тем самым $\mathscr{D}(\mathbb{D})_q'$
наделяется структурой $\operatorname{Fun}(\mathbb{D})_q$-бимодуля.
\end{corollary}

\begin{remark}
Возникающую за счет вложения $\mathrm{Pol}(\mathfrak{p}^-)_q\hookrightarrow
\operatorname{Fun}(\mathbb{D})_q$ структуру
$\mathrm{Pol}(\mathfrak{p}^-)_q$-бимодуля в $\mathscr{D}(\mathbb{D})_q'$
можно получить по-другому -- продолжением по непрерывности с
$\mathrm{Pol}(\mathfrak{p}^-)_q$:
$$
f_1\times f_2:f\mapsto f_1\cdot f\cdot f_2,\qquad
f_1,f_2,f\in\mathrm{Pol}(\mathfrak{p}^-)_q.
$$
\end{remark}

\medskip

\begin{corollary}\label{embed_finite}
1. Каждый элемент $f \in \mathscr{D}(\mathbb{D})_q'$ единственным образом
разлагается в сходящийся ряд
\begin{equation}\label{series_f0}
  f = \sum\limits_{k,j=0}^\infty f_{k,j},\qquad
  f_{k,j} \in\mathbb{C}[\mathfrak{p}^-]_{q,k} \; f_0 \;
\mathbb{C}[\mathfrak{p}^+]_{q,-j}.
\end{equation}
2. Элемент $f \in \mathscr{D}(\mathbb{D})_q'$ принадлежит
$\mathscr{D}(\mathbb{D})_q$, если и только если все члены ряда
(\ref{series_f0}), кроме конечного их числа, равны нулю.
\end{corollary}

В частном случае квантового круга следующее утверждение вытекает из
предложения \ref{def-finite}.

\begin{proposition}\label{image}
Элемент $f \in \mathscr{D}(\mathbb{D})_q'$ принадлежит подмногообразию
$\mathscr{D}(\mathbb{D})_q$, если и только если для некоторого $N$
\begin{equation}\label{eq_finite}
  \mathbb{C}[\mathfrak{p}^+]_{q,-n}\; f =
f\; \mathbb{C}[\mathfrak{p}^-]_{q,n} = 0
\end{equation}
при всех $n \geq N$.
\end{proposition}

{\bf Доказательство.} Если $f \in \mathscr{D}(\mathbb{D})_q$, то при
достаточно больших $n$ равенства (\ref{eq_finite}) выполняются по
определению алгебр $\mathrm{Fun}(\mathbb{D})_q$ и
$\mathscr{D}(\mathbb{D})_q$. С другой стороны, из (\ref{eq_finite}) следует,
что
$$
f_0 \,\mathbb{C}[\mathfrak{p}^+]_{q,-n} \, f = f\,
\mathbb{C}[\mathfrak{p}^-]_{q,n}\, f_0 = 0
$$
при достаточно больших $n$. Значит, согласно предложению \ref{inner1}, равны
нулю все члены ряда (\ref{series_f0}), кроме конечного их числа. Остается
воспользоваться следствием \ref{embed_finite}. \hfill $\square$

\subsubsection{Действие алгебры $U_q\mathfrak{g}$ в пространствах основных и
обобщенных функций.}\label{Uq-mod}

\begin{lemma}\label{finite_deg}
Для любого $\xi \in U_q\mathfrak{g}$ существует такое число $N \in
\mathbb{N}$, что
$$
\xi \left(\mathbb{C}[\mathfrak{p}^-]_{q,k} \;
\mathbb{C}[\mathfrak{p}^+]_{q,-j} \right) \subset \newoplus_{|k'-k| \leq N
\; \& \; \ |j'-j| \leq N} \mathbb{C}[\mathfrak{p}^-]_{q,k'} \;
\mathbb{C}[\mathfrak{p}^+]_{q,-j'}
$$
при всех $k,j \in \mathbb{Z}_+$.
\end{lemma}

{\bf Доказательство.} Множество тех элементов $\xi$, для которых справедливо
утверждение леммы, является подалгеброй и содержит все элементы $E_i$,
$F_i$, $K_i^\pm$, $i=1,2,\ldots,l$. \hfill $\square$

\medskip

  Напомним  что
$\mathrm{Pol}(\mathfrak{p}^-)_q$ является $U_q\mathfrak{g}$-модульным
$\mathrm{Pol}(\mathfrak{p}^-)_q$-бимодулем. Как было показано в предыдущем
\itemiiiе, структура $\mathrm{Pol}(\mathfrak{p}^-)_q$-бимодуля переносится
по непрерывности на $\mathscr{D}(\mathbb{D})_q'$. Из леммы \ref{finite_deg}
следует, что структура $U_q\mathfrak{g}$-модуля также переносится по
непрерывности на $\mathscr{D}(\mathbb{D})_q'$ и что эти структуры
согласованы. Более точно, справедливо следующее утверждение.

\begin{proposition}\label{Uq-action}
1. Для любого $\xi\in U_q\mathfrak{g}$ оператор
\begin{equation}\label{U-action}
\mathrm{Pol}(\mathfrak{p}^-)_q\to\mathrm{Pol}(\mathfrak{p}^-)_q,\quad
f\mapsto\xi f
\end{equation}
допускает продолжение по непрерывности на топологическое векторное
пространство $\mathscr{D}(\mathbb{D})_q'$. Тем самым векторное пространство
$\mathscr{D}(\mathbb{D})_q'$ наделяется структурой
$U_q\mathfrak{g}$-модульного $\mathrm{Pol}(\mathfrak{p}^-)_q$-бимодуля.
\\ 2. Для всех $\xi\in U_q\mathfrak{g}$, $f\in\mathscr{D}(\mathbb{D})_q'$
\begin{equation}\label{*-action}
(\xi f)^*=(S(\xi))^*f^*.
\end{equation}
\end{proposition}

\bigskip

Используя предложение \ref{image}, получаем
\begin{corollary}\label{D-D'}
Линейное подмногообразие $\mathscr{D}(\mathbb{D})_q \subset
\mathscr{D}(\mathbb{D})_q'$ является $U_q\mathfrak{g}$-модульным
$\mathrm{Pol}(\mathfrak{p}^-)_q$-бимодулем.
\end{corollary}

Напомним обозначение, введенное в \itemiiiе \ref{Hermit},\ \ \
$$\mathfrak{q}^\pm=\mathfrak{k} \oplus \mathfrak{p}^\pm.$$

\begin{corollary}\label{submodules}
Рассмотрим $U_q\mathfrak{g}$-модуль $\mathscr{D}(\mathbb{D})_q$. Линейное
подмногообразие $\mathbb{C}[\mathfrak{p}^-]_qf_0$ является его
$U_q\mathfrak{q}^+$-подмодулем, а линейное подмногообразие
$f_0\mathbb{C}[\mathfrak{p}^+]_q$ --- его $U_q\mathfrak{q}^-$-подмодулем:
$$
U_q\mathfrak{q}^+(\mathbb{C}[\mathfrak{p}^-]_qf_0)\subset
\mathbb{C}[\mathfrak{p}^-]_qf_0,\qquad
U_q\mathfrak{q}^-(f_0\mathbb{C}[\mathfrak{p}^+]_q)\subset f_0
\mathbb{C}[\mathfrak{p}^+]_q.
$$
\end{corollary}

{\bf Доказательство.} Финитный линейный оператор в $\mathcal{H} =
\bigoplus\limits_{j=0}^\infty \mathcal{H}_j$ принадлежит образу подалгебры
$\mathbb{C}[\mathfrak{p}^-]_q f_0$ при гомоморфизме $T_F$, если и только
если его ядро содержит $\bigoplus\limits_{j=1}^\infty \mathcal{H}_j$.
Значит, в силу предложения \ref{inner1},
$$
\mathbb{C}[\mathfrak{p}^-]_q f_0 = \left\{f \in \mathscr{D}(\mathbb{D})_q
\left|\, f\left(\newoplus_{j=1}^\infty
\mathbb{C}[\mathfrak{p}^-]_{q,j}\right) = 0 \right. \right\}.
$$
Аналогично доказывается равенство
$$
f_0 \mathbb{C}[\mathfrak{p}^+]_q = \left\{f \in \mathscr{D}(\mathbb{D})_q
\left|\, \left(\newoplus_{j=1}^\infty
\mathbb{C}[\mathfrak{p}^+]_{q,-j}\right) f = 0\right. \right\}.
$$
Остается воспользоваться тем, что $\mathscr{D}(\mathbb{D})_q$ является
$U_q\mathfrak{q}^\pm$-модульной алгеброй и
$$
U_q\mathfrak{q}^+ \left(\newoplus_{j=1}^\infty
\mathbb{C}[\mathfrak{p}^-]_{q,j}\right) \subset \newoplus_{j=1}^\infty
\mathbb{C}[\mathfrak{p}^-]_{q,j};\quad U_q\mathfrak{q}^-
\left(\newoplus_{j=1}^\infty \mathbb{C}[\mathfrak{p}^+]_{q,-j}\right)
\subset
\newoplus_{j=1}^\infty \mathbb{C}[\mathfrak{p}^+]_{q,-j}.
$$
\hfill $\square$

\medskip

\begin{note}\label{geom_real}
Линейное отображение
$$
\mathcal{H}\stackrel{\sim}{_\to}\mathbb{C}[\mathfrak{p}^-]_qf_0,\qquad
\psi v_0 \mapsto \psi f_0,
$$
где $\psi\in\mathbb{C}[\mathfrak{p}^-]_q$, является изоморфизмом
$\mathrm{Pol}(\mathfrak{p}^-)_q$-модулей и доставляет ``геометрическую
реализацию'' фоковского представления. \IND{фоковское представление !
геометрическая реализация} В силу следствия \ref{submodules}, это линейное
отображение наделяет векторное пространство $\mathcal{H}$ структурой
$U_q\mathfrak{q}^+$-модульного $\mathrm{Pol}(\mathfrak{p}^-)_q$-модуля.
Соответствующее представление алгебры $U_q\mathfrak{q}^+$ в $\mathcal{H}$
обозначим $\Gamma$.
\end{note}

\medskip

Покажем, что ${\rm Pol}(\mathbb{C})_{q_{l_0}}$ естественно вкладывается в
${\rm Pol}(\mathfrak{p}^-)_q$.

\begin{lemma}\label{z1}
1. Рассмотрим градуированную $U_q\mathfrak{g}$-модульную алгебру
$\mathbb{C}[\mathfrak{p}^-]_q$. Существует и единствен такой элемент
$z_{\operatorname{low}} \in \mathbb{C}[\mathfrak{p}^-]_{q,1}$ , что
\begin{equation}\label{z_l0}
  F_{l_0} z_{\operatorname{low}} = q_{l_0}^{1/2},\qquad
  F_j z_{\operatorname{low}} = 0,\;
  j \neq l_0.
\end{equation}
2. Имеют место равенства
\begin{equation}\label{z_l0_weight}
  K_i^{\pm 1} z_{\operatorname{low}} =
   q_i^{\pm a_{i l_0}} z_{\operatorname{low}},\quad
  i=1,2,\ldots,l,
\end{equation}
\begin{equation}\label{E_z_l0}
E_{l_0} z_{\operatorname{low}} = -q_{l_0}^{1/2}
z_{\operatorname{low}}^2,
\end{equation}
\medskip
$$  E_{l_0} z_{\operatorname{low}}^* =
q_{l_0}^{-3/2},\quad F_{l_0} z_{\operatorname{low}}^* = -q^{5/2}
z_{\operatorname{low}}^{*2},\quad K_i^{\pm 1} z_{\operatorname{low}}^* =
q_i^{\mp a_{i l_0}} z_{\operatorname{low}}^*, \; i=1,2,\cdots,l.
$$
\end{lemma}

{\bf Доказательство.}
 1. Как упоминалось в \itemiiiе \ref{prehomogeneous},
$U_q\mathfrak{k}$-модуль $\mathbb{C}[\mathfrak{p}^-]_{q,1}$ прост. Значит,
его младший вектор $z_{\operatorname{low}} \neq 0$ определен с точностью до
ненулевого числового множителя. Докажем, что
\begin{equation}\label{F_l0_zlow}
F_{l_0}z_{\operatorname{low}} \neq 0.
\end{equation}
Действительно, в противном случае
 $U_q\mathfrak{n}^- z_{\operatorname{low}}=\mathbb{C} z_{\operatorname{low}}$.
Следовательно, при всех $\xi \in U_q\mathfrak{n}^-$ линейный функционал
$\xi\, z_{\operatorname{low}} \in \mathbb{C}[\mathfrak{p}^-]_{q,1} $
аннулирует вектор
 $v(\mathfrak{q}^+,0) \in N(\mathfrak{q}^+,0)$.
Значит, $z_{\operatorname{low}}$ аннулирует все векторы из
$U_q\mathfrak{n}^-\, v(\mathfrak{q}^+,0)$. Другими словами,
$z_{\operatorname{low}}=0$, что противоречит выбору элемента
$z_{\operatorname{low}}$.

  Как следует из  \eqref{F_l0_zlow}, элемент $z_{\operatorname{low}}$
можно выбрать так, что $F_{l_0} z_{\operatorname{low}} = q_{l_0}^{1/2}$, и
он является единственным элементом со свойствами (\ref{z_l0}).

 2. Равенство \eqref{z_l0_weight} очевидно, поскольку
младший вес $U_q\mathfrak{k}$-модуля $\mathbb{C}[\mathfrak{p}^-]_{q,1}$
отличается знаком от старшего веса $U_q\mathfrak{k}$-модуля
$N(\mathfrak{q}^+,0)_{-1}$, а старшим вектором последнего является вектор
$F_{l_0} v(\mathfrak{q}^+,0)$ веса $-\alpha_{l_0}$.

Равенство $E_{l_0}z_{\operatorname{low}}=\mathrm{const}\cdot
z_{\operatorname{low}}^2$ вытекает из того, что младшее весовое
подпространство $U_q\mathfrak{k}$-модуля $\mathbb{C}[\mathfrak{p}^-]_{q,2}$
порождено элементом $z_{\operatorname{low}}^2$, а вектор
$E_{l_0}z_{\operatorname{low}}\in\mathbb{C}[\mathfrak{p}^-]_{q,2}$
принадлежит этому весовому подпространству (поскольку его вес равен
$2\alpha_{l_0}$). Равенство $\mathrm{const}=-q_{\operatorname{low}}^{1/2}$
вытекает из того, что
$$
(E_{l_0}F_{l_0}-F_{l_0}E_{l_0})z_{\operatorname{low}}=
\frac{K_{l_0}-K_{l_0}^{-1}}{q_{l_0}-q_{l_0}^{-1}}z_{\operatorname{low}}.
$$
Остальные равенства следуют из уже доказанных и из условия согласованности
инволюций \eqref{*-action}. \hfill $\square$


\begin{lemma}\label{z2}
\begin{equation}\label{like_disc}
  z_{\operatorname{low}}^* z_{\operatorname{low}} =
  q_{l_0}^2 z_{\operatorname{low}} z_{\operatorname{low}}^* + 1 - q_{l_0}^2.
\end{equation}
\end{lemma}

{\bf Доказательство.} Воспользуется мультипликативной формулой \\
\eqref{Rmatrix} для универсальной $R$-матрицы, считая, что приведенное
разложение элемента $w_0\in W$ выбрано так же, как в разделе \ref{GVerma}:
$w_0=w_{0,\mathbb{S}}\cdot\,^\mathbb{S}w_0$, где
$$
w_{0,\mathbb{S}}=s_{i_1}s_{i_2}\ldots s_{i_{M'}},\qquad
\,^\mathbb{S}w_0=s_{i_{M'+1}}s_{i_{M'+2}}\ldots s_{i_M}.
$$
Из определения умножения $m:\mathrm{Pol}(\mathfrak{p}^-)_q^{\otimes
2}\to\mathrm{Pol}(\mathfrak{p}^-)_q$ следует, что с точностью до порядка
тензорных сомножителей элемент $m(z_{\operatorname{low}}^*\otimes
z_{\operatorname{low}})$ равен
$$
\exp_{q_{M'}^2}((q_{M'}^{-1})E_{\beta_{M'}}\otimes F_{\beta_{M'}})
q^{-t_0}z_{\operatorname{low}}^*\otimes z_{\operatorname{low}}.
$$
Нетрудно показать, что
\begin{equation}\label{two}
\beta_{M'}=\alpha_{l_0},\qquad E_{\beta_{M'}}\otimes F_{\beta_{M'}}=
c\cdot E_{l_0}\otimes F_{l_0},\quad c\in\mathbb{C}\backslash\{0\}.
\end{equation}
Из первого равенства и из (\ref{fraction}) следует, что
$$
q^{-t_0}(z_{\operatorname{low}}^*\otimes z_{\operatorname{low}})=
q_{l_0}^{-\frac{H_{l_0}\otimes H_{l_0}}2}(z_{\operatorname{low}}^*\otimes
z_{\operatorname{low}}).
$$
Значит,
\begin{equation}\label{c}
z_{\operatorname{low}}^*z_{\operatorname{low}}=
q_{l_0}^2z_{\operatorname{low}}z_{\operatorname{low}}^*+c(1-q_{l_0}^2),\quad
c\ne 0.
\end{equation}

Остается доказать, что $c=1$. Для этого достаточно подействовать элементом
$F_{l_0} \in U_q\mathfrak{g}$ на обе части равенства
$$
z_{\operatorname{low}}^* z_{\operatorname{low}}^2 = q_{l_0}^4
z_{\operatorname{low}}^2 z_{\operatorname{low}}^* + c(1-q_{l_0}^4)
z_{\operatorname{low}},
$$
что приводит к квадратному уравнению для $c$, оба решения которого ($c=0$,
$c=1$) легко угадываются. \hfill $\square$

\begin{corollary}\label{z4}
Наделим $\mathrm{Pol}(\mathfrak{p}^-)_q$ структурой
$U_{q_{l_0}}\mathfrak{su}_{1,1}$-модульной алгебры с помощью вложения
$$
U_{q_{l_0}}\mathfrak{sl}_2 \hookrightarrow U_q\mathfrak{g},\qquad K^{\pm 1}
\mapsto K_{l_0}^{\pm 1},\quad E \mapsto E_{l_0},\quad F \mapsto F_{l_0}.
$$
Отображение $z \mapsto z_{l_0}$ единственным образом продолжается до
вложения $U_q\mathfrak{su}_{1,1}$-модульных $*$-алгебр
\begin{equation}\label{embed_C}
  \mathrm{Pol}(\mathbb{C})_{q_{l_0}} \hookrightarrow
  \mathrm{Pol}(\mathfrak{p}^-)_q.
\end{equation}
\end{corollary}

\medskip Опишем действие образующих алгебры $U_q\mathfrak{g}$ на элемент
$f_0$.

\begin{proposition}\label{action_f0}
$K_i^{\pm 1} f_0 = f_0$ при $i=1,2,\ldots,l$,
$$
F_j f_0 =
\left\{
\begin{array}{rl}
  -\frac{q_{l_0}^{1/2}}{q_{l_0}^{-2}-1} f_0 z_{\operatorname{low}}^*, &\; j=l_0 \\
  0, & \; j \neq l_0
\end{array}
\right.,\qquad
 E_j f_0 = \left\{
\begin{array}{rl}
  -\frac{q_{l_0}^{1/2}}{1 - q_{l_0}^2} z_{\operatorname{low}} f_0, &\; j=l_0 \\
  0, &\; j \neq l_0
\end{array}
\right.
$$
\end{proposition}

{\bf Доказательство.} Отметим, что $H_\mathbb{S} f_0 = 0$, поскольку в
разложении $$ f_0 = \sum\limits_{j,k=0}^\infty a_{j,k},\qquad a_{j,k} \in
\mathbb{C}[\mathfrak{p}^-]_{q,k} \mathbb{C}[\mathfrak{p}^+]_{q,-j} ,$$ равны
нулю все члены $a_{j,k}$ с $j \neq k$, как следует из определения
представления $T_F$ и из леммы \ref{ball_I_8.7}. Значит, из следствия
\ref{submodules} вытекает $U_q\mathfrak{k}$-инвариантность элемента $f_0$:
\begin{equation}\label{K_weight}
  K_i^{\pm 1} f_0 = f_0,\quad i=1,2,\ldots,l, \qquad
E_j f_0 = F_j f_0 = 0,\quad j \neq l_0.
\end{equation}
Осталось доказать равенства
$
 F_{l_0}\,f_0\,=\,-\frac{q_{l_0}^{1/2}}{q_{l_0}^{-2}-1} f_0
 z_{\operatorname{low}}^*$, $
 E_{l_0}\,f_0\,=\,-\frac{q_{l_0}^{1/2}}{1 - q_{l_0}^2} z_{\operatorname{low}}
 f_0.
$
 Из следствия \ref{submodules} вытекает, что
$F_{l_0} f_0 \in f_0 \mathbb{C}[\mathfrak{p}^+]_q$, $E_{l_0} f_0 \in
\mathbb{C}[\mathfrak{p}^-]_q f_0$. Это утверждение можно уточнить. Из
(\ref{K_weight}) следует, что веса элементов $F_{l_0} f_0$, $E_{l_0} f_0$
равны $-\alpha_{l_0}$ и $\alpha_{l_0}$ соответственно. Но весовые
подпространства таких весов в $f_0 \mathbb{C}[\mathfrak{p}^+]_q$,
$\mathbb{C}[\mathfrak{p}^-]_q f_0$ одномерны и порождены элементами $f_0
z_{\operatorname{low}}^*$, $z_{\operatorname{low}} f_0$. Значит, $ F_{l_0}
f_0 = c^- f_0 z_{\operatorname{low}}^*,\quad E_{l_0} f_0 = c^+
z_{\operatorname{low}} f_0 $.

Остается найти числовые множители $c^-$, $c^+$. Из равенств
$$
f_0 z_{\rm low}= 0,\quad \Delta(F_{l_0}) = F_{l_0} \otimes K_{l_0}^{-1} + 1
\otimes F_{l_0},\quad
 z_{\rm low}^* z_{\rm low}\,=\,
  q_{l_0}^2 z_{\rm low} z_{\rm low}^* + 1 - q_{l_0}^2
$$
следует, что
$$
 0 = F_{l_0}(f_0 z_{\rm low})\, =\,
 (c^- f_0 z_{\rm low}^*)(q_{l_0}^{-2}\,
 z_{\rm low}) + f_0\, q_{l_0}^{1/2} = (c^- q_{l_0}^{-2} (1-q_{l_0}^2)
  + q_{l_0}^{1/2}) f_0.
$$
Значит, $c^- = -q_{l_0}^{1/2} / (q_{l_0}^{-2}-1)$. Используя равенства
$(S(E_{l_0}))^*=q_{l_0}^{-2}F_{l_0}$ и (\ref{*-action}), находим $c^+$.
 \hfill $\square$

\bigskip Отметим, что в частном случае квантового круга
доказанное утверждение было получено ранее в \itemiiiе \ref{finite_sl_2}.

\bigskip
Наделим $\operatorname{Fun}(\mathbb{D})_q$ структурой $U_q\mathfrak{g}$
модуля с помощью разложения
 $\operatorname{Fun}(\mathbb{D})_q=
 \operatorname{Pol}(\mathfrak{p}^-)_q \oplus
 \mathscr{D}(\mathbb{D})_q$.

\begin{proposition}\label{F_alg}
$\mathrm{Fun}(\mathbb{D})_q$ является $U_q\mathfrak{g}$-модульной алгеброй.
\end{proposition}

{\bf Доказательство.} Воспользуемся описанием алгебры
$\mathrm{Fun}(\mathbb{D})_q$ в терминах образующих и соотношений:
образующими служат $f_0$ и ненулевые элементы
$\mathrm{Pol}(\mathfrak{p}^-)_q$, а множество определяющих соотношений
состоит из (\ref{rel1}), (\ref{rel2}) и всех соотношений, имеющихся между
ненулевыми элементами алгебры $\mathrm{Pol}(\mathfrak{p}^-)_q$. Действие
элементов $K_i^{\pm 1}$, $E_i$, $F_i$, $i=1,2,\ldots,l$, алгебры Хопфа
$U_q\mathfrak{g}$ на эти образующие алгебры $\mathrm{Fun}(\mathbb{D})_q$
было описано выше в настоящем \itemiiiе и в \itemiiiах \ref{vect},
\ref{braid_algebra}. Используя это описание, наделим свободную алгебру,
порожденную образующими алгебры $\mathrm{Fun}(\mathbb{D})_q$ структурой
$U_q\mathfrak{g}$-модульной алгебры.

 Остается убедиться в том, что  двусторонний идеал $J$ свободной алгебры,
 отвечающий определяющим соотношениям алгебры $\mathrm{Fun}(\mathbb{D})_q$,
  является  $U_q\mathfrak{g}$-подмодулем. Другими
словами, остается убедиться в том, что действие образующих алгебры
$U_q\mathfrak{g}$ на образующие двустороннего идеала $J$ не выводит из $J$.

Начнем с элемента $f_0^2\,-\,f_0\in J$. Используя предложение
\ref{action_f0}, получаем
$$
K_i^{\pm 1}(f_0 f_0 - f_0) = (K_i^{\pm 1} f_0)(K_i^{\pm 1} f_0) - K_i^{\pm
1} f_0 = f_0 f_0 - f_0,\quad i=1,2,\ldots,l,
$$
$$
F_j (f_0 f_0 - f_0) = E_j (f_0 f_0 - f_0) = 0,\quad j \neq l_0,
$$
$$
F_{l_0} (f_0 f_0 - f_0) =
-\frac{q_{l_0}^{1/2}}{q_{l_0}^{-2}-1}(f_0(z_{\operatorname{low}}^*f_0) +
(f_0 f_0 - f_0) \; z_{\operatorname{low}}^*),
$$
$$
E_{l_0} (f_0 f_0 - f_0) = -\frac{q_{l_0}^{1/2}}{1 -
q_{l_0}^2}(z_{\operatorname{low}}(f_0 f_0 - f_0) + (f_0
z_{\operatorname{low}}) f_0).
$$
Очевидно, правые части полученных равенств принадлежат $J$.

Действие элементов $K_i^{\pm 1}$, $E_i$, $F_i$, $i=1,2,\ldots,l$, на
элементы, отвечающие остальным определяющим соотношениям, не выводит из $J$
потому, что $\mathrm{Pol}(\mathfrak{p}^-)_q$ является
$U_q\mathfrak{g}$-модульной алгеброй, а $\mathscr{D}(\mathbb{D})_q$ --
$U_q\mathfrak{g}$-модульным $\mathrm{Pol}(\mathfrak{p}^-)_q$-бимодулем.
\hfill $\square$

\medskip

 \begin{corollary}\label{D_alg}
$*$-Алгебра $\mathscr{D}(\mathbb{D})_q$ является $(U_q\mathfrak{g},
*)$-модульной.
\end{corollary}

\subsubsection{Существование, единственность и явный вид инвариантного
интеграла на $\mathscr{D}(\mathbb{D})_q$.}\label{InvInt}

В предыдущем \itemiiiе показано, что алгебра $\mathscr{D}(\mathbb{D})_q$
является $U_q\mathfrak{g}$-модульной. Докажем существование инвариантного
интеграла $\nu: \mathscr{D}(\mathbb{D})_q \to \mathbb{C}$ и его
единственность с точностью до числового множителя. Будем следовать
\itemiiiу \ref{inv_int_sl_2}, где эти результаты получены в частном случае
квантового круга.

Рассмотрим введенное в \itemiiiе \ref{Uq-mod} представление $\Gamma$ алгебры
$U_q \mathfrak{q}^+$ в фоковском пространстве $\mathcal{H}$. Каждому
$\lambda =\sum \limits_{i=1}^{l} n_i\alpha_i \in Q$ отвечают элемент
$K_\lambda=K_1^{n_1}K_2^{n_2}\ldots K_l^{n_l}\in U_q \mathfrak{h}$ и
линейный оператор $\Gamma(K_{\lambda})$.

Конечномерность линейных операторов $T_F(f)$ при
$f\in\mathscr{D}(\mathbb{D})_q$ позволяет ввести $q$-след
$$
\mathrm{tr}_qT_F(f)\stackrel{\operatorname{def}}{=}
\mathrm{tr}(T_F(f)\Gamma(K_{-2\rho})),\qquad
f\in\mathscr{D}(\mathbb{D})_q,
$$
где $\rho$ -- полусумма положительных корней алгебры Ли $\mathfrak{g}$.

\begin{proposition}\label{exist_int}
Линейный функционал
\begin{equation}\label{int_D}
\int\limits_{\mathbb{D}_q}fd\nu\stackrel{\operatorname{def}}{=}
\left(1-q^2_{l_0}\right)^{\dim\mathfrak{p}^-}\,\mathrm{tr}_qT_F(f),
\end{equation}
на $\mathscr{D}(\mathbb{D})_q$ является положительным
$U_q\mathfrak{g}$-инвариантным интегралом.
\end{proposition}

{\bf Доказательство} Сужение линейного оператора $\Gamma(K_{-2\rho})$ на
каждое из весовых подпространств $U_q \mathfrak{q}^+$-модуля $\mathcal{H}$
является положительным скалярным оператором. Значит, $\mathrm{tr}_q
T_F(f^*f)=\mathrm{tr}(T_F(f) \Gamma(K_{-2\rho})T_F(f)^*)>0$ для всех $f \in
\mathscr{D}(\mathbb{D})_q$, и положительность линейного функционала
(\ref{int_D}) вытекает из предложения \ref{bijection}. Разумеется,
\begin{equation}\label{real_int}
\int \limits_{\mathbb{D}_q} f^* d\nu = \overline{\int
\limits_{\mathbb{D}_q}f d\nu}.
\end{equation}
Из предположения \ref{tr_g} и равенства (\ref{inner_S2}) следует, что
линейный функционал (\ref{int_D}) является $U_q \mathfrak{q}^+$-инвариантным
интегралом. Его $U_q \mathfrak{g}$-инвариантность вытекает из
$U_q\mathfrak{q}^+$-инвариантности и из равенств (\ref{real_int}),
(\ref{*-action}). $\hfill \square$

\bigskip
Перейдем к доказательству единственности инвариантного интеграла.

\begin{lemma}\label{L7.1} Элемент $f_0 \in \mathscr{D}(\mathbb{D})_q$ порождает
$U_q \mathfrak{b}^+$-модуль $\mathbb{C}[\mathfrak{p}^-]_q f_0$ и $U_q
\mathfrak{b}^-$-модуль $f_0 \mathbb{C}[\mathfrak{p}^+]_q$.
\end{lemma}

{\bf Доказательство} Как следует из (\ref{*-action}), достаточно показать,
что $U_q \mathfrak{b}^+\, f_0 = \mathbb{C}[\mathfrak{p}^-]_q \, f_0$.
Используя предложение \ref{inner1}, легко доказать, что $U_q
\mathfrak{b}^+$-инвариантное спаривание
$$\mathbb{C}[\mathfrak{p}^+]_q \times \mathbb{C}[\mathfrak{p}^-] f_0
 \to \mathbb{C}, \quad
f_+ \times f_- f_0 \mapsto \int \limits_{\mathbb{D}_q} f_+ f_- f_0 d\nu$$
невырождено.

Значит, в категории градуированных $U_q \mathfrak{b}^+$-модулей
$\mathbb{C}[\mathfrak{p}^-]_q f_0 \cong \mathbb{C}[\mathfrak{p}^+]_q^*$. В
этой же категории $N(\mathfrak{q}^-,0) \cong
\mathbb{C}[\mathfrak{p}^+]_q^*$, как следует из определения
$U_q\mathfrak{g}$-модуля $\mathbb{C}[\mathfrak{p}^+]_q$. Значит,
$\mathbb{C}[\mathfrak{p}^-]_q f_0 \cong N(\mathfrak{q}^-,0)$, причем
последний изоморфизм отображает однородную компоненту $\mathbb{C} f_0$ на
однородную компоненту $\mathbb{C} v(\mathfrak{q}^-,0)$. Остается
воспользоваться тем, что вектор $v(\mathfrak{q}^-,0)$ порождает $U_q
\mathfrak {b}^+$-модуль $N(\mathfrak{q}^-,0)$. $\hfill \square$

\medskip Следующее утверждение обобщает предложение \ref{t2.3.9}.
\begin{proposition}\label{generator_f0} $U_q \mathfrak{g}\, f_0
= \mathscr{D}(\mathbb{D})_q$.
\end{proposition}
{\bf Доказательство} Как следует из леммы \ref{L7.1}, достаточно показать,
что
\begin{equation} \label{incl_7.1}
(E_{j_1} E_{j_2}\ldots E_{j_r} f_0) \mathbb{C}[\mathfrak{p}^+]_q \subset
U_q \mathfrak{g}\, f_0
\end{equation}
при всех $r \in \mathbb{Z}_+,\quad j_1,j_2,\ldots,j_r \in
\{1,2,\ldots,l\}$. Воспользуемся методом математической индукции. При $r=0$
включение (\ref{incl_7.1}) следует из леммы \ref{L7.1}. Опишем переход от
$r-1$ к $r$.

Пусть $\varphi=E_{j_2} E_{j_3} \ldots E_{j_r} f_0,\; \psi \in
\mathbb{C}[\mathfrak{p}^+]_q$. По предположению индукции
$$\varphi \psi \in U_q \mathfrak{g}\, f_0, \qquad (K_{j_1} \varphi)(E_{j_1} \psi)
\in U_q \mathfrak{g} f_0.$$ Значит, из предложения \ref{F_alg} следует, что
$$ (E_{j_1} E_{j_2}\ldots E_{j_r} f_0) \psi = (E_{j_1} \varphi ) \psi =
E_{j_1} (\varphi \psi) - (K_{j_1} \varphi )(E_{j_1} \psi),$$ поскольку
$\Delta(E_{j_1})=E_{j_1} \otimes 1 + K_{j_1} \otimes E_{j_1}$.

Остается заметить, что элементы $E_{j_1} (\varphi \psi)$ и $(K_{j_1}
\varphi )(E_{j_1} \psi)$ принадлежат $U_q\mathfrak{g} f_0$ по предположению
индукции.
 $\hfill
\square$

\begin{corollary}\label{unique_int} $U_q \mathfrak{g}$-инвариантный интеграл на
$\mathscr{D}(\mathbb{D})_q$ единствен с точностью до числового множителя.
\end{corollary}

\bigskip Положительность $U_q \mathfrak{g}$-инвариантного интеграла
(\ref{int_D}) позволяет наделить векторное пространство
$\mathscr{D}(\mathbb{D})_q$ нормой $\| f \|_2 = \left(\int \limits
_{\mathbb{D}_q} f^* f d\nu\right)^\frac{1}{2}$. Пополнение
$\mathscr{D}(\mathbb{D})_q$ по этой норме является гильбертовым
пространством. Оно обозначается $L^2(d\nu)_q$ и является $q$-аналогом
пространства $L^2(d\nu)$ квадратично суммируемых по инвариантной мере $d\nu$
функций в ограниченной симметрической области $\mathbb{D}$.

 Как следует из (\ref{int_D}), для любых однородных нормированных векторов
 $v_i \in \mathcal{H}_i$, $v_j \in \mathcal{H}_j$
градуированного предгильбертова пространства
$\mathcal{H}=\newoplus\limits_{k=0}^\infty \mathcal{H}_k$
 имеет место неравенство
 $$ |(T_F(f) v_i, v_j)| \leq \operatorname{const}(i,j) \,\|f\|_2,
   \qquad f \in \mathscr{D}(\mathbb{D})_q.
 $$
Значит, описанное в \itemiiiе \ref{gen} вложение $\mathscr{D}(\mathbb{D})_q
\hookrightarrow \mathscr{D}(\mathbb{D})_q'$ допускает продолжение по
непрерывности до вложения $L^2(d\nu)_q \hookrightarrow
\mathscr{D}(\mathbb{D})_q'$ (инъективность легко следует из определения
гильбертова пространства $L^2(d\nu)_q$ и из явного вида инвариантного
интеграла (\ref{int_D})). В дальнейшем гильбертово пространство
$L^2(d\nu)_q$ будет отождествляться с его образом при вложении в
$\mathscr{D}(\mathbb{D})_q'$.

\bigskip
Для любой обобщенной функции $f \in \mathscr{D}(\mathbb{D})_q'$ в квантовой
ограниченной симметрической области $\mathbb{D}$ корректно определена
полуторалинейная форма $(T_F(f) v',v'')$ в $\mathcal{H}$. Обобщенную функцию
$f\in \mathscr{D}(\mathbb{D})_q'$ будем называть ограниченной, если
$$ \|f \|_\infty = \sup\limits_{v',v'' \in \mathcal{H}\backslash \{0\}}
\frac{(T_F(f)v',v'')}{\|v'\|\|v''\|}<\infty.$$

В этом случае существует и единствен ограниченный линейный оператор
$\overline{T}$ в пополнении $\overline{\mathcal{H}}$ предгильбертова
пространства $\mathcal{H}$, для которого $(\overline{T} v',v'')=(T_F(f)
v',v'')$ при всех $v',v''\in \mathcal{H}$. Ограниченные обобщенные функции с
нормой $\| \cdot\|_\infty$ образуют $C^*$-алгебру $L^\infty(\mathbb{D})_q$,
естественно изоморфную $C^*$-алгебре всех ограниченных линейных операторов в
гильбертовом пространстве $\overline{\mathcal{H}}$. Если $\varphi \in
L^\infty(\mathbb{D})_q, \psi \in L^2(d\nu)_q$, то $\varphi \psi \in
L^2(d\nu)_q$ и
\begin{equation}\label{bound_mult}
\| \varphi \psi\|_2 \leq \| \varphi\|_\infty \|\psi\|_2
\end{equation}

Действительно,
$$
\|\varphi\psi\|_2^2=\left(1-q_{l_0}^2\right)^{\dim\mathfrak{p}^-}
\mathrm{tr}(T_F(\varphi)T_F(\psi)\Gamma(K_{-2\rho})T_F(\psi)^*
T_F(\varphi)^*)\le
$$
$$
\le\|\varphi\|_\infty^2\left(1-q_{l_0}^2\right)^{\dim\mathfrak{p}^-}
\mathrm{tr}(T_F(\psi)\Gamma(K_{-2\rho})T_F(\psi)^*)\le\|\varphi\|_\infty^2
\|\psi\|_2^2.
$$
Сопоставляя каждому элементу $\varphi\in L^\infty(\mathbb{D})_q$ линейный
оператор $\hat{\varphi}:\psi\mapsto\varphi\psi$, получаем естественное
представление $C^*$-алгебры $L^\infty(\mathbb{D})_q$ в гильбертовом
пространстве $L^2(d\nu)_q$.

\subsubsection{Изоморфизм  $U_q\mathfrak{g}$-модулей $\mathscr{D}(\mathbb{D})_q$
и $U_q\mathfrak{g}\otimes_{U_q\mathfrak{k}}\mathbb{C}$
.}\label{rel-for-f_0}

Каждому весовому локально конечномерному $U_q\mathfrak{k}$-модулю $V$
сопоставим $U_q\mathfrak{g}$-модуль
$P(V)=U_q\mathfrak{g}\otimes_{U_q\mathfrak{k}}V$. Как нетрудно показать, он
является проективным объектом категории $C(\mathfrak{g},\mathfrak{k})_q$.
Такие проективные объекты называют стандартными, ср. \cite[стр. 104]{KV}.
\IND{стандартный проективный объект} Разложим
$P(\mathbb{C})=U_q\mathfrak{g}\otimes_{U_q\mathfrak{k}}\mathbb{C}$ в
произведение обобщенных модулей Верма в тензорной категории
$U_q\mathfrak{g}^\mathrm{cop}$-модулей.

\medskip

Вектор $v(\mathfrak{q}^+,0) \otimes v(\mathfrak{q}^-,0)\in
  N(\mathfrak{q}^+,0) \otimes N(\mathfrak{q}^-,0)
$ является $U_q\mathfrak{k}$-инвариантным, что позволяет ввести в
рассмотрение морфизм $U_q\mathfrak{g}$-модулей
\begin{equation}\label{P-to-GVerma}
\mathcal{J}:P(\mathbb{C})\to N(\mathfrak{q}^+,0)\otimes
N(\mathfrak{q}^-,0), \quad \mathcal{J}:1\otimes 1\mapsto
 v(\mathfrak{q}^+,0) \otimes v(\mathfrak{q}^-,0).
\end{equation}

\begin{proposition}\label{P-of-C}
Морфизм $U_q\mathfrak{g}$ модулей \eqref{P-to-GVerma} биективен.
\end{proposition}

Доказательству предпошлем следующую лемму.
\begin{lemma}\label{iso_special_new} Пусть
$P_\mathfrak{h}(\mathbb{C})=U_q\mathfrak{g}\otimes_{U_q\mathfrak{h}}
\mathbb{C}$. В тензорной категории $U_q\mathfrak{g}^{\rm cop}$-модулей
морфизм $U_q\mathfrak{g}$-модулей
$$
\mathcal{J}_\mathfrak{h}:P_\mathfrak{h}(\mathbb{C})\to
 N(\mathfrak{b}^+,0)\otimes N(\mathfrak{b}^-,0), \qquad
\mathcal{J}_\mathfrak{h}:1\otimes 1\mapsto
 v(\mathfrak{b}^+,0)\otimes v(\mathfrak{b}^-,0)
$$
является изоморфизмом.
\end{lemma}

{\bf Доказательство леммы.} Наделим $U_q\mathfrak{h}$-модули
$N(\mathfrak{b}^\pm,0)$ фильтрациями, связанными с весовыми
подпространствами
$$
F_{\mu^+}N(\mathfrak{b}^+,0)=\newoplus_{\substack{\nu^+\in Q_+\\
0\le\nu^+\le\mu^+}}N(\mathfrak{b}^+,0)_{-\nu^+},
\quad
F_{\mu^-}N(\mathfrak{b}^-,0)=\newoplus_{\substack{\nu^-\in Q_+\\
0\le\nu^-\le\mu^-}}N(\mathfrak{b}^-,0)_{\nu^-}.
$$
Тензорное произведение $N(\mathfrak{b}^+,0)\otimes N(\mathfrak{b}^-,0)$
наделим аналогичной фильтрацией
\begin{multline*}
F_{\mu^+,\mu^-}\left(N(\mathfrak{b}^+,0)\otimes N(\mathfrak{b}^-,0)
\right)=
\\ =\newoplus_{\substack{\nu^+\in Q_+\\ 0\le\nu^+\le\mu^+}}\;
\newoplus_{\substack{\nu^-\in Q_+\\ 0\le\nu^-\le\mu^-}}
N(\mathfrak{b}^+,0)_{-\nu^+}\otimes N(\mathfrak{b}^-,0)_{\nu^-}.
\end{multline*}
Наконец, используя присоединенное действие $U_q\mathfrak{h}$ в
$U_q\mathfrak{n}^\pm$, наделим фильтрацией $U_q\mathfrak{h}$-модуль
$P_\mathfrak{h}(\mathbb{C})$:
$$
F_{\mu^+,\mu^-}\left(P_\mathfrak{h}(\mathbb{C})
\right)=\sum_{\substack{\nu^+\in Q_+\\
\nu^+\le\mu^+}}\;\sum_{\substack{\nu^-\in Q_+\\ \nu^-\le\mu^-}}\text{л.о.}
\left((U_q\mathfrak{n}^-)_{-\nu^+}(U_q\mathfrak{n}^+)_{\nu^-}\otimes
1\right).
$$
где $\mu^\pm$ -- элементы решетки весов $P$.

Воспользуемся этими фильтрациями для доказательства леммы. Во-первых,
нетрудно показать, что
\begin{equation}\label{sur_special}
\mathcal{J}_\mathfrak{h}\;F_{\mu^+,\mu^-}
\left(P_\mathfrak{h}(\mathbb{C})\right)\supset
F_{\mu^+,\mu^-}\left(N(\mathfrak{b}^+,0)\otimes N(\mathfrak{b}^-,0)\right).
\end{equation}
Действительно, в противном случае среди пар $\mu^+,\mu^-$, для которых не
выполняется \eqref{sur_special}, есть минимальные (отношение частичного
порядка на множестве таких пар вводится очевидным образом). Приходим к
противоречию, поскольку если
$$
\xi=F_{i_1}^{j_1}F_{i_2}^{j_2}F_{i_3}^{j_3}\ldots,\;
\eta=E_{k_1}^{l_1}E_{k_2}^{l_2}E_{k_3}^{l_3}\ldots,\qquad
i_1,i_2,\ldots;\;k_1,k_2\ldots\in\{1,2,\ldots,l\},
$$
то с точностью до членов младших степеней
$$
\mathcal{J}_\mathfrak{h}(\xi\eta(1\otimes 1))=\mathrm{const}\;\xi
v(\mathfrak{b}^+,0)\otimes\eta v(\mathfrak{b}^-,0),
$$
где $\mathrm{const}$ -- ненулевое число.

Во-вторых, из результатов \itemiiiа \ref{PBW} следует, что подпространства
$F_{\mu^+,\mu^-} \left(P_\mathfrak{h}(\mathbb{C})\right)$,
$F_{\mu^+,\mu^-}\left(
N(\mathfrak{b}^+,0)\otimes
N(\mathfrak{b}^-,0)\right)$ конечномерны, и их
размерности равны. Значит,
$$
\mathcal{J}_\mathfrak{h}:\;F_{\mu^+,\mu^-}
\left(P_\mathfrak{h}(\mathbb{C})\right)
\overset{\approx}{\longrightarrow}F_{\mu^+,\mu^-} \left(
 N(\mathfrak{b}^+,0)\otimes N(\mathfrak{b}^-,0)\right),
$$
и остается воспользоваться равенствами
$\cup_{\mu^+,\mu^-}F_{\mu^+,\mu^-}\left( P_\mathfrak{h}(\mathbb{C})\right)=
P_\mathfrak{h}(\mathbb{C})$,
$$
\bigcup\limits_{\mu^+,\mu^-}F_{\mu^+,\mu^-} \left(
 N(\mathfrak{b}^+,0)\otimes N(\mathfrak{b}^-,0)\right)=
N(\mathfrak{b}^+,0)\otimes N(\mathfrak{b}^-,0).\eqno \square $$

{\bf Доказательство предложения.} Будем рассуждать так же, как при
доказательстве леммы. Напомним, что весовые $U_q\mathfrak{k}$-модули
наделены градуировками, определяемыми с помощью элемента
$H_\mathbb{S}\in\mathfrak{h}$, см. \itemiii \ \ref{Weight}. В частности,
$$
N(\mathfrak{q}^+,0)=\newoplus\limits_{j=0}^\infty
N(\mathfrak{q}^+,0)_{-j},\qquad
N(\mathfrak{q}^-,0)=\newoplus\limits_{j=0}^\infty N(\mathfrak{q}^-,0)_{j}.
$$
Наделим $ N(\mathfrak{q}^+,0)\otimes N(\mathfrak{q}^-,0)$ фильтрацией
$$F_{i^+i^-}\left(
N(\mathfrak{q}^+,0)\otimes N(\mathfrak{q}^-,0)\right)= \newoplus_{0\le
j^+\le i^+}\;\newoplus_{0\le j^-\le i^-} N(\mathfrak{q}^+,0)_{-j^+}\otimes
N(\mathfrak{q}^-,0)_{j^-}. $$
Используя изоморфизм
$$
 U_q\mathfrak{p}^-\otimes U_q\mathfrak{p}^+\otimes
U_q\mathfrak{k}\overset{\approx}{\longrightarrow}U_q\mathfrak{g},\qquad
\xi^-\otimes\xi^+\otimes\eta\mapsto \xi^-\xi^+\eta
$$
и присоединенное действие $U_q\mathfrak{k}$ в $U_q\mathfrak{p}^\pm$,
наделим $U_q\mathfrak{k}$-модуль $P(\mathbb{C})$ фильтрацией
$$
F_{i^+i^-}\left(P(\mathbb{C})\right)= \sum_{0\le j^+\le i^+}\;\sum_{0\le
j^-\le
i^-}\text{л.о.}\left((U_q\mathfrak{p}^-)_{-j^+}(U_q\mathfrak{p}^+)_{j^-}\otimes
1\right).
$$
Из определений следует, что $P(\mathbb{C})=
\cup_{i^+,i^-}F_{i^+i^-}\left(P(\mathbb{C})\right)$,
$$
N(\mathfrak{q}^+,0)\otimes
N(\mathfrak{q}^-,0)=\bigcup\limits_{i^+,i^-}F_{i^+i^-} \left(
N(\mathfrak{q}^+,0)\otimes N(\mathfrak{q}^-,0)\right).
$$
Кроме того, подпространства $ F_{i^+i^-}\left( N(\mathfrak{q}^+,0)\otimes
N(\mathfrak{q}^-,0)\right)$ и $F_{i^+i^-}\left(P(\mathbb{C})\right) $
конечномерны, а их размерности равны, см. \itemiiiы \ref{GVerma} и
\ref{Add_Uq}.

Остается доказать включение
\begin{equation}\label{sur_new}
\mathcal{J} F_{i^+i^-}\left(P(\mathbb{C})\right)\supset F_{i^+i^-}\left(
N(\mathfrak{q}^+,0)\otimes N(\mathfrak{q}^-,0)\right),\qquad
i^+,i^-\in\mathbb{Z}_+.
\end{equation}
 Используя коммутативную диаграмму
$$
\xymatrix{P_{\mathfrak{h}}(\mathbb{C}) \ar[r]^-{\mathcal{J}_\mathfrak{h}}
\ar[d]_\varphi & N(\mathfrak{b}^+,0)\otimes N(\mathfrak{b}^-,0)\ar[d]^\psi
\\ P(\mathbb{C})
\ar[r]^-{\mathcal{J}} & N(\mathfrak{q}^+,0)\otimes N(\mathfrak{q}^-,0)}
$$
в категории $U_q\mathfrak{g}^{\rm cop}$-модулей, нетрудно получить
\eqref{sur_new} из \eqref{sur_special} и из следующих утверждений.

Прежде всего, подпространство $N(\mathfrak{q}^+,0)_{j^+}$ является образом
суммы весовых подпространств $N(\mathfrak{b}^+,0)_{-\nu^+}$, для которых
$\nu^+(H_\mathbb{S})=j^+$, при каноническом эпиморфизме
$N(\mathfrak{b}^+,\lambda^+)\to N(\mathfrak{q}^+,\lambda^+)$. Аналогично,
$N(\mathfrak{q}^-,0)_{j^-}$ -- образ суммы весовых подпространств
$N(\mathfrak{b}^-,0)_{\nu^-}$, для которых $\nu^-(H_\mathbb{S})=j^-$.
Наконец, $ (U_q\mathfrak{n}^-)_{-\nu^+}(U_q\mathfrak{n}^+)_{\nu^-} \otimes
1 \subset (U_q\mathfrak{q}^-)_{-\nu^+}(U_q\mathfrak{q}^+)_{\nu^-} \otimes
1$ и
$$
\mathcal{J}\,\left((U_q\mathfrak{q}^-)_{-\nu^+}(U_q\mathfrak{q}^+)_{\nu^-}
\otimes 1\right) \subset
(U_q\mathfrak{p}^-)_{-j^+}(U_q\mathfrak{p}^+)_{j^-}
(v(\mathfrak{q}^+,0)\otimes v(\mathfrak{q}^-,0)),
$$
если $j^\pm=\nu^\pm(H_\mathbb{S})$.

Первые утверждения дают возможность оценить сверху прообраз подпространства
$F_{i^+i^-}\left(N(\mathfrak{q}^+,0)\otimes N(\mathfrak{q}^-,0)\right)$ при
отображении $\psi$, а последнее -- оценить снизу прообраз
$F_{i^+i^-}(P(\mathbb{C}))$ при отображении $\varphi$. \hfill $\square$

\bigskip

Докажем, что $P(\mathbb{C})\approx \mathscr{D}(\mathbb{D})_q$. Рассмотрим
векторное пространство
$$
\mathrm{Pol}(\mathfrak{p}^-)_q=\mathbb{C}[\mathfrak{p}^-]_q\otimes
\mathbb{C}[\mathfrak{p}^+]_q=N(\mathfrak{q}^+,0)^* \otimes
N(\mathfrak{q}^-,0)^*
$$
и его убывающую фильтрацию подпространствами $\newoplus\limits_{j+k\geq
n}\;
 \mathbb{C}[\mathfrak{p}^-]_{q,j}\,
\mathbb{C}[\mathfrak{p}^+]_{q,-k}$.
 Их аннуляторы в сопряженном
пространстве монотонно возрастают, и из предложения \ref{P-of-C} следует,
что объединение аннуляторов является $U_q\mathfrak{g}$-подмодулем,
естественно изоморфным $P(\mathbb{C})$.

С другой стороны, как легко показать, используя основные свойства алгебры
финитных функций и инвариантного интеграла, это объединение аннуляторов
естественно изоморфно $\mathscr{D}(\mathbb{D})_q$.

Мы получили изоморфизм $U_q\mathfrak{g}$-модулей
$$
P(\mathbb{C})\stackrel{\approx}{\to}\mathscr{D}(\mathbb{D})_q,\qquad
1\otimes 1\mapsto\frac1{\left(1-q^2_{l_0}\right)^{\dim\mathfrak{p}^-}}\,f_0.
$$
Тем самым доказано следующее уточнение предложения \ref{generator_f0},
обобщающее предложение \ref{def-rel-D(D)}.

\begin{proposition}\label{def-rel-D(D)-general}
Элемент $f_0$ -- образующая $U_q\mathfrak{g}$-модуля
$\mathscr{D}(\mathbb{D})_q$, и соотношение
$$\xi f_0=\varepsilon(\xi)f_0,\qquad\xi\in U_q\mathfrak{k},$$
является его определяющим соотношением.
\end{proposition}

\begin{remark}\label{like_HeckenbergerKolb}
Использованная в доказательстве двойственность $U_q\mathfrak{g}$-модульной
алгебры $\mathrm{Pol}(\mathbb{C})_q$ и $U_q\mathfrak{g}$-модульной
коалгебры $P(\mathbb{C})$ позволяет считать $P(\mathbb{C})$ квантовым
аналогом пространства обобщенных функций с носителями в точке
$\mathbf{z}=0$.
\end{remark}

\begin{remark}\label{GVerma-tensor-product}
Пусть $\lambda^+\in P^\mathbb{S}_+$, $\lambda^-\in-P^\mathbb{S}_-$ и
$w_{0,\mathbb{S}}$ -- элемент максимальной длины группы Вейля
$W_\mathbb{S}\subset W$. Рассмотрим обобщенный модуль Верма
$N(\mathfrak{q}^+,\lambda^+)$ со старшим весом $\lambda^+$ и обобщенный
модуль Верма $N(\mathfrak{q}^-,\lambda^-)$ с младшим весом $\lambda^-$.
Доказательство предложения \ref{P-of-C} без существенных изменений
переносится на случай ненулевых весов $\lambda^\pm$ и доставляет изоморфизм
$P(L(\mathfrak{k},w_{0,\mathbb{S}}\lambda^-)\otimes
L(\mathfrak{k},\lambda^+))\stackrel{\approx}{\to}
N(\mathfrak{q}^+,\lambda^+)\otimes N(\mathfrak{q}^-,\lambda^-)$ в тензорной
категории $U_q\mathfrak{g}^\mathrm{cop}$-модулей.
\end{remark}

\subsubsection{Ограниченность операторов фоковского представления
$*$-алгебры $\mathrm{Pol}(\mathfrak{p}^-)_q$.} \label{boundness}

В классическом случае $q=1$ полиномы на векторном пространстве
$\mathfrak{p^-}$ принадлежат $L^\infty(\mathbb{D})$, поскольку область
$\mathbb{D}$ ограничена. Получим аналогичный результат в случае $0<q<1$ в
подтверждение ''ограниченности рассматриваемых нами квантовых областей''.

\bigskip Начнем с доказательства  ограниченности линейного оператора
$T_F(z_{\operatorname{low}})$. В частном случае квантового круга это
очевидно: $\|T_F(z)\|=1$. Рассмотрим вложение $*$-алгебр
$$
i: \mathrm{Pol}(\mathbb{C})_{q_{l_0}} \hookrightarrow
\mathrm{Pol}(\mathfrak{p}^-)_q,\qquad i: z \mapsto z_{\operatorname{low}},
$$
см. следствие \ref{z4}.

\begin{lemma}\label{bound_zl0}1.
Представление $T_F \cdot i$ в предгильбертовом пространстве $\mathcal{H}$
кратно фоковскому представлению $*$-алгебры
$\mathrm{Pol}(\mathbb{C})_{q_{l_0}}$.\\ 2.\ \ \ \
$\|T_F(z_{\operatorname{low}})\| = 1$.
\end{lemma}

{\bf Доказательство.}
 Достаточно доказать первое
 утверждение.
 Пусть $T=T_F \cdot i$. Рассмотрим подпространство
   $\mathcal{L} =\sum_{j=0}^\infty
   (\operatorname{Ker}T(z^*) \cap \mathcal{H}_j)$. Достаточно доказать, что
линейная оболочка $\mathcal{M}$ множества
 $\{T(z^j)v| j \in \mathbb{Z}_+,\; v \in \mathcal{L}\}$ совпадает с $\mathcal{H}$.

 Очевидно,
$ \mathcal{M}=\sum_{j=0}^\infty (\mathcal{M} \cap \mathcal{H}_j)$.
Предположим, что $\mathcal{M} \neq \mathcal{H}$. Тогда $\mathcal{M} \cap
\mathcal{H}_j \neq \mathcal{H}_j$ при некотором $j \in \mathbb{Z}_+$, и,
следовательно, существует ненулевой вектор $v \in \mathcal{H}$,
ортогональный $\mathcal{M}$. Для некоторого $N \in \mathbb{Z}_+$ имеем $
T(z^*)^Nv \neq 0$ и $T(z^*)^{(N+1)}v =0$, поскольку оператор $T(z^*)$
уменьшает степень однородности. Значит, $T(z)^N\, T(z^*)^N v \in
\mathcal{M}$ и $(T(z^*)^Nv,T(z^*)^Nv)=0$, что противоречит выбору $N$.
\hfill $\square$

\begin{corollary}\label{z_zstar}
$\|T_F(z_{\operatorname{low}} z_{\operatorname{low}}^*)\| = 1$.
\end{corollary}

Рассмотрим линейный функционал $\langle \cdot \rangle$ на
$\mathrm{Pol}(\mathfrak{p}^-)_q$, введенный равенством \eqref{expan2}. Пусть
$(\cdot, \cdot)$ -- стандартная эрмитова форма в
$\mathbb{C}[\mathfrak{p}^-]_q$
\begin{equation}\label{herm}
  (\psi_1, \psi_2) = \langle \psi_2^* \psi_1 \rangle, \qquad
  \psi_1, \psi_2 \in \mathbb{C}[\mathfrak{p}^-]_q,
  \end{equation}
  см. \itemiii\  \ref{Pol}.

\begin{lemma}\label{inner}
\begin{itemize}

\item[1.] Стандартная эрмитова форма в $\mathbb{C}[\mathfrak{p}^-]_{q,1}$
положительна, $U_q\mathfrak{k}$-инвариантна и
\begin{equation}\label{inner_zl0w}
  (z_{\operatorname{low}}, z_{\operatorname{low}}) = 1 - q_{l_0}^2.
\end{equation}

\item[2.] Если $\{ z_j \}$ -- ортогональный базис в
$\mathbb{C}[\mathfrak{p}^-]_{q,1}$, то элемент
$$
\Psi = \sum\limits_{j=1}^{\dim\, \mathfrak{p}^-} \frac{z_j z_j^*}{(z_j,
z_j)}
$$
не зависит от выбора этого базиса и $U_q\mathfrak{k}$-инвариантен.

\item[3.] Рассмотрим $U_q\mathfrak{k}$-модуль
    $\mathbb{C}[\mathfrak{p}^-]_{q,1}\cdot
    \mathbb{C}[\mathfrak{p}^+]_{q,-1}$ -- линейную оболочку множества
    $\left\{\psi^-\cdot\psi^+\in
    \operatorname{Pol}(\mathfrak{p}^-)_q\;\left|\;\psi^-\in
    \mathbb{C}[\mathfrak{p}^-]_{q,1}\,,\,\psi^+\in
    \mathbb{C}[\mathfrak{p}^+]_{q, -1}\right.\right\}$. Подпространство
    его $U_q\mathfrak{k}$-инвариантов одномерно и порождается элементом
    $\Psi$.
\end{itemize}
\end{lemma}

{\bf Доказательство.} Положительность стандартной эрмитовой формы следует из
предложения \ref{positive_herm}, а равенство (\ref{inner_zl0w}) -- из
(\ref{like_disc}), (\ref{herm}) и из соотношений $f_0
\mathbb{C}[\mathfrak{p}^-]_{q,1} = \mathbb{C}[\mathfrak{p}^+]_{q, -1} f_0 =
0$.
  Из определения линейного  функционала $\langle \cdot \rangle$ вытекает
его $U_q\mathfrak{k}$-инвариантность, откуда следует
$U_q\mathfrak{k}$-инвариантность эрмитовой формы \eqref{herm}.
 Элемент $\Psi$, очевидно, не зависит от выбора ортогонального базиса $\{
z_j \}$.

  Для доказательства оставшихся утверждений сопоставим каждому элементу
  $\mathscr{K}$ векторного пространства
$\mathbb{C}[\mathfrak{p}^-]_{q,1} \mathbb{C}[\mathfrak{p}^+]_{q, -1}$
линейный оператор умножения на $\mathscr{K}$ в пространстве
$\mathbb{C}[\mathfrak{p}^-]_{q,1} \, f_0$: $$
  \mathbb{C}[\mathfrak{p}^-]_{q,1} \, f_0 \mapsto
\mathbb{C}[\mathfrak{p}^-]_{q,1}\, f_0,\qquad \psi \mapsto \mathscr{K} \psi.
$$ Из предложения \ref{inner1} следует, что полученное линейное отображение
\begin{equation}\label{15.12}
\mathbb{C}[\mathfrak{p}^-]_{q,1} \mathbb{C}[\mathfrak{p}^+]_{q, -1}
\rightarrow \operatorname{End}(\mathbb{C}[\mathfrak{p}^-]_{q,1} f_0)
\end{equation}
биективно, а из $U_q\mathfrak{k}$-инвариантности элемента $f_0$ вытекает,
что это изоморфизм $U_q\mathfrak{k}$-модулей. Значит, оно осуществляет
биекцию пространств \hbox{$U_q\mathfrak{k}$-инвариантов}.
  Элементу $\Psi$ при отображении
(\ref{15.12}) отвечает единичный оператор в
$\mathbb{C}[\mathfrak{p}^-]_{q,1} f_0$. Следовательно, этот элемент является
\hbox{$U_q\mathfrak{k}$-инвариантным}.

 Из простоты $U_q\mathfrak{k}$-модуля
 $\mathbb{C}[\mathfrak{p}^-]_{q,1} \cong \mathbb{C}[\mathfrak{p}^-]_{q,1} f_0$,
 отмечавшейся в разделе
\ref{vect}, следует одномерность подпространства
$U_q\mathfrak{k}$-инвариантов $U_q\mathfrak{k}$-модуля
$\operatorname{End}(\mathbb{C}[\mathfrak{p}^-]_{q,1} f_0)$. Значит,
подпространство $U_q\mathfrak{k}$-инвариантов $U_q\mathfrak{k}$-модуля
$\mathbb{C}[\mathfrak{p}^-]_{q,1} \mathbb{C}[\mathfrak{p}^+]_{q, -1}$
одномерно. Оно содержит ненулевой $U_q\mathfrak{k}$-инвариантный элемент
$\Psi$, и, следовательно, порождается этим элементом. \hfill $\square$

\begin{proposition}\label{bound_psi}
Линейный оператор $T_F(\Psi)$ в предгильбертовом пространстве $\mathcal{H}$
ограничен.
\end{proposition}

{\bf Доказательство.}
 Так же, как в разделах \ref{alg_gr}, \ref{regular_K}, введем в рассмотрение
$*$-алгебру Хопфа $\mathbb{C}[K]_q \hookrightarrow (U_q\mathfrak{k})^*$
матричных элементов \hbox{$P$-весовых} конечномерных
$U_q\mathfrak{k}$-модулей. Пусть $\nu: \mathbb{C}[K]_q \rightarrow
\mathbb{C}$ -- инвариантный интеграл на квантовой группе $K$ с условием
нормировки : $\int\limits_{K_q} 1 d\nu = 1$.

Из локальной конечномерности действия $U_q\mathfrak{k}$ в
$\mathrm{Pol}(\mathfrak{p}^-)_q$ следует, что корректно определено
кодействие
\begin{equation}\label{co_act_Pol}
  \Delta: \mathrm{Pol}(\mathfrak{p}^-)_q \rightarrow
  \mathrm{Pol}(\mathfrak{p}^-)_q
\otimes \mathbb{C}[K]_q,
\end{equation}
см. \itemiii\ \ref{SL_2_sl_2}. Рассмотрим линейный оператор ``усреднения по
действию квантовой группы $K$'':
$$
\operatorname{av}: \mathrm{Pol}(\mathfrak{p}^-)_q
\rightarrow\mathrm{Pol}(\mathfrak{p}^-)_q, \qquad \operatorname{av}: f
\mapsto (\mathrm{id} \otimes \nu) \Delta (f).
$$
Очевидно, $\Delta: (\mathbb{C}[\mathfrak{p}^-]_{q,1}
\mathbb{C}[\mathfrak{p}^+]_{q, -1}) \rightarrow
(\mathbb{C}[\mathfrak{p}^-]_{q,1} \mathbb{C}[\mathfrak{p}^+]_{q, -1})
\otimes \mathbb{C}[K]_q$,
$$
\operatorname{av}: (\mathbb{C}[\mathfrak{p}^-]_{q,1}
\mathbb{C}[\mathfrak{p}^+]_{q, -1}) \rightarrow
(\mathbb{C}[\mathfrak{p}^-]_{q,1} \mathbb{C}[\mathfrak{p}^+]_{q, -1}).
$$
Из $U_q\mathfrak{k}$-инвариантности элемента
$\operatorname{av}(z_{\operatorname{low}} z_{\operatorname{low}}^*)$ и из
леммы \ref{inner} следует, что
\begin{equation}\label{z_psi}
  \operatorname{av}(z_{\operatorname{low}} z_{\operatorname{low}}^*) = \mathrm{const} \cdot \Psi.
\end{equation}

Покажем, что в последнем равенстве $\mathrm{const} \neq 0$. Используя
морфизм $U_q\mathfrak{k}$-модулей \eqref{15.12} и $q$-след $\mathrm{tr}_q:
\mathrm{End}(\mathbb{C}[\mathfrak{p}^-]_{q,1} f_0) \rightarrow \mathbb{C}$,
введенный в \itemiiiе \ref{Weight}, нетрудно найти ненулевой
$U_q\mathfrak{k}$-инвариантный интеграл $\mu$ на
$\mathbb{C}[\mathfrak{p}^-]_{q,1} \, \mathbb{C}[\mathfrak{p}^+]_{q, -1}$ и
доказать, что \hbox{$\mu(z_{\operatorname{low}} z_{\operatorname{low}}^*)
\neq 0$.} Имеем
$$
(\mu \otimes \mathrm{id}) \Delta (\psi) = \mu(\psi) \otimes 1, \qquad \psi
\in \mathbb{C}[\mathfrak{p}^-]_{q,1} \mathbb{C}[\mathfrak{p}^+]_{q, -1}.
$$
Следовательно, $\mu(\operatorname{av}(z_{\operatorname{low}}
z_{\operatorname{low}}^*)) = \mu(z_{\operatorname{low}}
z_{\operatorname{low}}^*) \neq 0$. Значит, $\mathrm{const} \cdot \mu(\Psi)
\neq 0$.

Как следует из \eqref{z_psi}, для доказательства предложения
\ref{bound_psi} достаточно получить неравенство
$\|T_F(\operatorname{av}(z_{\operatorname{low}}z_{\operatorname{low}}^*))\|
\le 1$. Напомним, что однородные компоненты $\mathcal{H}_k$ градуированного
векторного пространства
\hbox{$\mathcal{H}=\bigoplus\limits_{k=0}^\infty\mathcal{H}_k$}
конечномерны и попарно ортогональны. Конечномерные \hbox{$C^*$-алгебры}
$\mathrm{End}(\mathcal{H}_k)$ являются $U_q\mathfrak{k}$-модульными, что
позволяет ввести в рассмотрение линейные операторы
$$
\Delta_k:\mathrm{End}(\mathcal{H}_k)\to\mathrm{End}(\mathcal{H}_k),\quad
\operatorname{av}_k:\mathrm{End}(\mathcal{H}_k)\to
\mathrm{End}(\mathcal{H}_k)
$$
так же, как были введены линейные операторы $\Delta$, $\operatorname{av}$.
Из определений следует, что
$$
\|T_F(\operatorname{av}(z_{\operatorname{low}}z_{\operatorname{low}}^*))\|=
\sup\limits_k\|T_F(\operatorname{av}(z_{\operatorname{low}}
z_{\operatorname{low}}^*))_{|_{\mathcal{H}_k}}\|=\sup\limits_{k}\|
\operatorname{av}_k(T_F(z_{\operatorname{low}}
z_{\operatorname{low}}^*)|_{\mathcal{H}_k})\|,
$$
поскольку линейное отображение
$$
\mathbb{C}[\mathfrak{p}^-]_{q,1}\mathbb{C}[\mathfrak{p}^+]_{q,1}\to
\mathrm{End}(\mathcal{H}_k),\qquad\psi\mapsto
T_F(\psi)|_{\mathcal{H}_k}
$$
является гомоморфизмом $U_q\mathfrak{k}$-модульных алгебр. Из леммы
\ref{bound_zl0} следует, что $\|T_F(z_{\operatorname{low}}
z_{\operatorname{low}}^*)|_{\mathcal{H}_k}\|\le 1$. Остается доказать, что
линейный оператор $\operatorname{av}_k$ ``усреднения по действию квантовой
группы $K$'' не увеличивает норму элемента $C^*$-алгебры
$\mathrm{End}\,\mathcal{H}_k$. Но
$\Delta_k:\mathrm{End}\,\mathcal{H}_k\to\mathrm{End}\,\mathcal{H}_k\otimes
C(K)_q$ является сжатием, как и любой гомоморфизм $C^*$-алгебр. Значит,
доказываемое предложение \ref{bound_psi} вытекает из приведенной ниже леммы
\ref{contr_inv}.

\medskip

Рассмотрим инвариантный интеграл $\nu: C(K)_q \rightarrow \mathbb{C}$ (см.
\itemiii\ \ref{IrrepK}). Пусть $\mathrm{Mat}_n$ -- алгебра комплексных
матриц порядка $n$. Известно, что $\mathrm{Mat}_n$ единственным образом
наделяется структурой $C^*$-алгебры.

\begin{lemma}\label{contr_inv}
При всех $n\in\mathbb{N}$ линейное отображение $C^*$-алгебр
$$\mathrm{id}\otimes\nu:\;\mathrm{Mat}_n\otimes C(K)_q\to\mathrm{Mat}_n$$
является сжимающим.
\end{lemma}

{\bf Доказательство.}
 Начнем со случая $n=1$. Из результатов раздела \ref{compact} следует, что
$C^*$-алгебра $C(K)_q$ обладает представлением в гильбертовом пространстве
$L^2(d\nu)_q$
$$
\widehat{\psi}: f \mapsto \psi f, \qquad \psi \in C(K)_q, \quad f \in
L^2(d\nu)_q,
$$
причем $\nu(\psi) = (\widehat{\psi} 1, 1)$. Значит, $ | \nu(\psi) | = |
(\widehat{\psi} 1, 1) | \leq \| \widehat{\psi} \| \leq \| \psi \|$. В общем
случае рассуждения совершенно аналогичны.

Именно, $C^*$-алгебра $\mathrm{Mat}_n \otimes C(K)_q$ естественным образом
действует в гильбертовом пространстве $\mathbb{C}^n \otimes L^2(d\nu)_q$.
Введем отвечающее этому действию представление $\mathcal{R}$. Пусть
${\mathcal{L} = \mathbb{C}^n \otimes 1 \hookrightarrow \mathbb{C}^n \otimes
L^2(d\nu)_q}$ и $P_\mathcal{L}$ -- ортопроектор на $\mathcal{L}$.
Утверждение леммы вытекает из очевидного неравенства
$$\|\mathrm{id} \otimes \nu (\psi) \| =
\| P_\mathcal{L}
\mathcal{R}(\psi)_{|_{\mathcal{L}}} \| \leq \|\psi\|,
 \qquad \psi \in\mathrm{Mat}_n \otimes C(K)_q. \eqno \square$$

\medskip
\begin{remark} Заключительную часть доказательства леммы
можно заменить ссылкой на известные результаты теории $C^*$-алгебр
\cite[стр. 148, 153]{Arv}, которые позволяют перейти от частного случая
$n=1$ к общему случаю $n \in \mathbb{Z}_+$ при оценке нормы линейного
отображения $\mathrm{id} \otimes \nu$.
\end{remark}

\medskip Предложение \ref{bound_psi} полностью доказано. Из него следует
ограниченность линейных операторов $T_F(f)$ при $f \in
\mathbb{C}[\mathfrak{p}^-]_{q,1}$. Но такие элементы порождают $*$-алгебру
$\operatorname{Pol}(\mathfrak{p}^-)_q$, как вытекает из предложения
\ref{hom_gen}. Значит, имеет место следующее утверждение.

\begin{corollary}\label{bound_D}
Все операторы представления $T_F$ в предгильбертовом пространстве
$\mathcal{H}$ ограничены.
\end{corollary}

Другими словами, $\mathrm{Fun}(\mathbb{D})_q \subset
L^\infty(\mathbb{D})_q$. Замыкание алгебры $\mathrm{Fun}(\mathbb{D})_q$ в
$L^\infty(\mathbb{D})_q$ обозначим $C(\overline{\mathbb{D}})_q$. Эта
$C^*$-алгебра является квантовым аналогом коммутативной $C^*$-алгебры
$C(\overline{\mathbb{D}})$ непрерывных функций в $\overline{\mathbb{D}}$. В
частном случае квантового круга она рассматривалась в \itemiiе
\ref{topology}.

\subsection{Каноническое вложение}\label{canonical_embed}

\subsubsection{Введение.}\label{canon_example}

В дальнейшем все пары $(\mathfrak{g},\mathfrak{k})$ будут предполагаться
эрмитово-симметрическими. Напомним, что в \itemiiiе \ref{Pol} были введены
картановская инволюция $\theta_q$ в $U_q\mathfrak{g}$
$$
\theta_q(K_i^{\pm 1})=K_i^{\pm 1},\qquad\theta_q(E_i)=
\begin{cases}
E_i, & i\ne l_0,
\\ -E_i, & i=l_0,
\end{cases}
\qquad \theta_q(F_i)=
\begin{cases}
F_i, & i\ne l_0,
\\ -F_i, & i=l_0.
\end{cases}
$$
и $*$-алгебра Хопфа $(U_q\mathfrak{g},*)$, где
$*=\star\theta_q=\theta_q\star$. По двойственности (см. \eqref{invol1})
определяется $*$-алгебра Хопфа $((U_q\mathfrak{g})^\star,\sharp)$:
\begin{equation}\label{g_0}
f^\sharp(\xi)=\overline{f((S(\xi))^*)},\qquad
f\in(U_q\mathfrak{g})^*,\;\xi\in U_q\mathfrak{g}.
\end{equation}

Из теоремы Бернсайда следует сюръективность гомоморфизма алгебр
$U_q\mathfrak{g} \rightarrow \operatorname{End}L(\Lambda)$ при всех
 $\Lambda\in P_+$. Значит,
$(\operatorname{End}L(\Lambda))^*\hookrightarrow (U_q\mathfrak{g})^\star$.
Так как
$$
\sharp:(\operatorname{End}L(\Lambda))^*\to
(\operatorname{End}L(-w_0\Lambda))^*,
$$
то инволюция $\sharp$ допускает сужение на подалгебру
$$
\mathbb{C}[G]_q=\sum\limits_{\Lambda \in
P_+}(\operatorname{End}L(\Lambda))^* \subset(U_q\mathfrak{g})^*.
$$

Введем в рассмотрение $*$-алгебру Хопфа
$\mathbb{C}[G_0]_q=(\mathbb{C}[G]_q,\sharp)$. \IND{$*$-алгебра ! Хопфа !
$\mathbb{C}[G_0]_q$} Это $q$-аналог алгебры регулярных функций на
вещественной аффинной алгебраической группе $G_0=\{g\in
G\,|\,f^\sharp(g)=\overline{f(g)},\quad f\in\mathbb{C}[G]\}$, ср. с
\eqref{compact_form}. \IND{$q$-аналог ! алгебры регулярных функций на
вещественной аффинной алгебраической группе} Известно, что $G_0$ действует
автоморфизмами ограниченной симметрической области $\mathbb{D}$,
стандартным образом вложенной в $\mathfrak{p}^-$, см. \cite{Kor}.


 Наличие выделенной точки $0 \in \mathbb{D}$ позволяет каждой функции $f(z)$
 в области $\mathbb{D}$ сопоставить функцию $f(g^{-1}\cdot 0)$ на $G_0$.  Нас
будут интересовать
 квантовые аналоги алгебр
функций в $\mathbb{D}$ и их образы при описанном вложении в алгебру функций
на $G_0$.

Если $\mathfrak{g}=\mathfrak{sl}_2$, то роль группы $G_0$ играет
вещественная аффинная алгебраическая группа $SU_{1,1}$ и, как нетрудно
показать,
\begin{equation}\label{inv_su11}
t_{11}^\sharp=t_{22},\qquad t_{12}^\sharp = q t_{21},
\end{equation}
где $\{t_{ij}\}_{i,j=1,2}$-- стандартные образующие алгебры
$\mathbb{C}[SL_2]_q$. Действительно, рассмотрим автоморфизм $\theta_q^*$
алгебры Хопфа $\mathbb{C}[SL_2]_q$, двойственный автоморфизму
$\theta_q:U_q\mathfrak{sl}_2\to U_q\mathfrak{sl}_2$,
$$
\theta_q:K^{\pm 1}\mapsto K^{\pm 1},\qquad\theta_q:E\mapsto
-E,\qquad\theta_q:F\mapsto -F.
$$
Он действует на образующие $t_{ij}$ следующим образом:
$$
\begin{pmatrix}\theta_q^*(t_{11}), & \; \theta_q^*(t_{12})\\
 \theta_q^*(t_{21}), &\;
\theta_q^*(t_{22})\end{pmatrix}=\begin{pmatrix}t_{11}, & -t_{12}\\ -t_{21},
& t_{22}\end{pmatrix}.
$$
Остается воспользоваться равенством \eqref{star_1}, устанавливающим связь
инволюций $*$ и $\star$, а также равенством \eqref{ast2}, описывающим
действие $\star$ на образующие $t_{ij}$.

Напомним, что каноническое вложение в случае
$\mathfrak{g}=\mathfrak{sl}_2$, то есть $\mathbb{D}=\{z\in \mathbb{C}\,|\,
|z| < 1\}$, подробно описано в
\itemiiiе \ref{canon_example_sl_2}.


\subsubsection{Квантовый аналог элемента $w_0$ группы Вейля.}\label{w_0}

Построение канонического вложения существенно усложняется при переходе от
частного случая $\mathfrak{g}=\mathfrak{sl}_2$ к общему случаю и использует
описанное ниже понятие квантовой группы Вейля.

\begin{proposition}\label{s_bar}
Существует и единствен такой линейный функционал $\overline{s}$ на
$\mathbb{C}[SL_2]_q$, что
\begin{align*}
\overline{s}(t_{22}f)&=\overline{s}(ft_{11})=0,&
\\ \overline{s}(t_{12}f)&=\overline{s}(ft_{12})=q\overline{s}(f),&
\\ \overline{s}(t_{21}f)&=\overline{s}(ft_{21})=-\overline{s}(f)&
\end{align*}
при всех $f\in\mathbb{C}[SL_2]_q$ и $f(1)=1$.
\end{proposition}

{\bf Доказательство.} Пусть $\{e_j\}$, $j\in\mathbb{Z}_+$, -- стандартный
базис гильбертова пространства $l^2(\mathbb{Z}_+)$ и $\Pi$ -- построенное в
\itemiiiе \ref{SU_2} $*$-представление алгебры
$\mathbb{C}[SU_2]_q=(\mathbb{C}[SL_2]_q,\star)$ в этом пространстве.
Линейный функционал $(\Pi(f)e_0,e_0)$ обладает всеми требуемыми свойствами.
Единственность $\overline{s}$ следует из того, что $\mathbb{C}[SL_2]_q$
является линейной оболочкой множества
$$ \{t_{22}^d t_{21}^c t_{12}^b t_{11}^a\;|\; a, b, c, d \in \mathbb{Z}_+\}.
\eqno \square
$$

\medskip
\begin{remark}
Нетрудно доказать, что каждый элемент алгебры $\mathbb{C}[SL_2]_q$
единственным образом разлагается в сумму
$$
f\,=\,\sum_{j=1}^\infty
t_{22}^j\,f_j(t_{12},t_{21})\,+\,f_0(t_{12},t_{21})\,+\,\sum_{j=1}^\infty
f_{-j}(t_{12},t_{21})\,t_{11}^j,
$$
где $f_j$ -- полиномы двух переменных, см. \itemiii\ \ref{int_SU2}. Это дает
второе доказательство существования элемента $\overline{s}$:
$$\overline{s}(f)=f_0(q,-1),\qquad f\in\mathbb{C}[SL_2]_q.$$

\bigskip
Элемент $\overline{s}\in\mathbb{C}[SL_2]_q^*$ является $q$-аналогом элемента
$\left(\begin{smallmatrix}0, & 1\\ -1, & 0\end{smallmatrix}\right)$ группы
$SL_2$, поскольку
\begin{equation}\label{s_bar_eq}
\begin{pmatrix}\overline{s}(t_{11}), & \overline{s}(t_{12})\\
\overline{s}(t_{21}), & \overline{s}(t_{22})\end{pmatrix}=\begin{pmatrix}0,
& q\\ -1, & 0\end{pmatrix}.
\end{equation}
\end{remark}

Наделим векторное пространство $U_q\mathfrak{sl}_2$ слабейшей из топологий,
в которых непрерывны все неприводимые конечномерные весовые представления
$\pi_\lambda: U_q\mathfrak{sl}_2 \to {\rm End}(L(\lambda))$, $\lambda\in
\mathbb{Z}_+\overline{\omega}$.
 Так как эти представления разделяют
элементы алгебры $U_q\mathfrak{sl}_2$, то введенная топология отделима и
пополнение алгебры $U_q\mathfrak{sl}_2$ канонически изоморфно прямому
произведению алгебр $\mathop{\times}\limits_j
\operatorname{End}L(j\overline{\omega})$, то есть алгебре
$\mathbb{C}[SL_2]_q^*$ со $*$-слабой топологией. Следующее утверждение
вытекает из определений.

\begin{lemma}\label{extension_pi} Каждое конечномерное весовое представление
алгебры $U_q\mathfrak{sl}_2$ в векторном пространстве $L$ допускает
продолжение по непрерывности на $\mathbb{C}[SL_2]_q^*$, наделяя $L$
структурой $\mathbb{C}[SL_2]_q^*$-модуля.
\end{lemma}

\bigskip Напомним, что для каждого весового $U_q \mathfrak{sl}_2$-модуля $V$ корректно определены
 линейные операторы $H$, $X^{\pm}$ в $V$, для которых
\begin{equation}\label{changesl_2}
  K^{\pm1}=q^{\pm H},\quad X^+=Eq^{-\frac{1}{2}H},
  \quad X^-=q^{\frac{1}{2}H}F,
\end{equation}
см. \itemiii \ \ref{Weight} и равенства \eqref{change1}, \eqref{change2}.

Пусть $j\in \mathbb{Z}_+$. Следуя традициям квантовой теории углового
момента, введем стандартный базис $\{e_m^j\}_{m=j,j-1,\ldots,-j}\}$
векторного пространства полагаем $L(2j\overline{\omega})$ равенствами
\begin{equation}\label{explicit_action}
\begin{gathered}
He_m^j=2me_m^j,\quad X^+e_m^j=
\begin{cases}
([j-m]_q[j+m+1]_q)^{1/2}e_{m+1}^j, & m\ne j
\\ 0, & m=j
\end{cases}
\\ X^-e_m^j=
\begin{cases}
([j-m+1]_q[j+m]_q)^{1/2}e_{m-1}^j, & m\ne-j
\\ 0, & m=-j
\end{cases},
\end{gathered}
\end{equation}
где, как обычно, $[n]_q=\dfrac{q^n-q^{-n}}{q-q^{-1}}$.

 \begin{lemma}(\cite[стр. 8]{VakSoib88})\label{simple_refl}\label{add_sbar}
 При всех $j,m$,
$$ \bar{s} e_m^j=(-1)^{j-m} q^{j^2+j}q^{-(m^2+m)}e_{-m}^j,$$ и, следовательно,
\begin{equation}\label{refl_weightK}
\overline{s}K\overline{s}^{-1}=K^{-1}.
\end{equation}
\end{lemma}

\bigskip

 Перейдем от частного случая $G=SL_2$ к общему
случаю. Напомним соглашение об односвязности $G$, см. \itemiii\
\ref{alg_gr}. Наделим векторное пространство $\mathbb{C}[G]_q^*\cong
\mathop{\times}\limits_{\lambda \in P_+}\operatorname{End}(L(\lambda))$
топологией прямого произведения, то есть $*$-слабой топологией.
Воспользуемся стандартными вложениями $*$-алгебр Хопфа
$$\varphi_i:\,U_{q_i}\mathfrak{sl}_2 \hookrightarrow U_q\mathfrak{g},$$
$$
\varphi_i:K^{\pm 1}\mapsto K_i^{\pm 1},\quad\varphi_i:E\mapsto
E_i,\quad\varphi_i:F\mapsto F_i,\qquad i=1,2,\ldots,l.
$$

\begin{proposition}\label{embed_Sl2}
Вложения алгебр $\varphi_i:U_{q_i}\mathfrak{sl}_2\hookrightarrow
U_q\mathfrak{g}$ допускают продолжения по непрерывности до вложений алгебр
$\overline{\varphi}_i:\mathbb{C}[SL_2]_{q_i}^*\hookrightarrow
\mathbb{C}[G]_q^*$.
\end{proposition}

{\bf Доказательство.} Существование продолжения вытекает из леммы
\ref{extension_pi}, а его инъективность -- из неравенства
$$
\dim\operatorname{Hom}_{\,U_{q_i}\mathfrak{sl}_2}
(L(k\overline{\omega}),L(k\overline{\omega}_i))>0, \qquad k \in
\mathbb{Z}_+, \eqno \square
$$
где $U_q\mathfrak{g}$-модуль $L(k\overline{\omega}_i)$ рассматривается как
$U_{q_i}\mathfrak{sl}_2$-модуль.

\bigskip

Пусть $\overline{s}_i=\overline{\varphi}_i(\overline{s})$, где
$i=1,2,\ldots,l$. Как известно \cite[стр. 157]{LS}, \cite[стр. 427]{KR},
элементы $\overline{s}_1,\overline{s}_2,\ldots,\overline{s}_l$ алгебры
$\mathbb{C}[G]_q^*$ удовлетворяют соотношениям кос: \IND{соотношения кос}
\begin{equation}\label{braid_rel}
\overline{s}_i\overline{s}_j\overline{s}_i\cdots=
\overline{s}_j\overline{s}_i\overline{s}_j\cdots,
\end{equation}
где число сомножителей в каждой из частей \eqref{braid_rel} равно двум,
если $a_{ij}a_{ji}=0$, трем, если $a_{ij}a_{ji}=1$, четырем, если
$a_{ij}a_{ji}=2$, и шести, если $a_{ij}a_{ji}=3$.\footnote{Последняя
возможность не реализуется в интересующем нас эрмитово-симметрическом
случае.} Пусть $w=s_{i_1}s_{i_2}\cdots s_{i_M}$ -- приведенное разложение
элемента $w\in W$, и
\begin{equation}\label{w_q}
\overline{w}=\overline{s}_{i_1}\overline{s}_{i_2}\cdots\overline{s}_{i_M}.
\end{equation}
Этот элемент алгебры $\mathbb{C}[G]_q^*$ не зависит от выбора приведенного
разложения, как следует из \eqref{braid_rel} и хорошо известного свойства
групп Кокстера \cite[стр. 19]{Bou4-6}.

Мы получили квантовый аналог $\overline{w}$ элемента $w \in W$ и, в
частности, квантовый аналог элемента $w_0$ максимальной длины группы Вейля
$W$.
Используя лемму \ref{simple_refl}, нетрудно доказать

\begin{proposition}[\cite{KorSoib}, стр. 138.]\label{s_general}
Элементы $\overline{s}_i$ алгебры $\mathbb{C}[G]_q^*$ обратимы, и
$$
\overline{s}_iK_j\overline{s}_i^{-1}=K_jK_i^{-a_{ij}},\qquad
i,j=1,2,\ldots,l.
$$
\end{proposition}

\begin{corollary}\label{refl_weight}
Пусть $V=\bigoplus\limits_\lambda V_\lambda$ -- разложение конечномерного
весового $U_q\mathfrak{g}$-модуля $V$ в сумму его весовых подпространств.
Тогда
\begin{equation}\label{w_action}
\overline{w}V_\lambda=V_{w\lambda},\qquad w\in W.
\end{equation}
\end{corollary}

\medskip

Покажем, что равенство
\begin{equation}\label{w0_tilda}
\widetilde{w}_0=\overline{w}_0^{\,-1}\cdot
q^{-\frac12\sum\limits_{i,j=1}^lc_{ij}d_id_jH_iH_j}.
\end{equation}
определяет элемент $\widetilde{w}_0$ алгебры $\mathbb{C}[G]_q^*$.
 Во-первых, из леммы \ref{add_sbar}
вытекает обратимость элементов $\overline{s}_i$, а, значит, и их
произведения $\overline{w}_0$. Во-вторых, $\mathbb{C}[G]_q^*$ -- прямое
произведение конечномерных алгебр $\operatorname{End}L(\lambda)$, и все
$U_q\mathfrak{g}$-модули $L(\lambda)$ являются весовыми. Значит, второй
множитель в \eqref{w0_tilda}, называемый корректирующим картановским
множителем, определяет элемент алгебры $\mathbb{C}[G]_q^*$.

Если
 $\mathfrak{g}=\mathfrak{sl}_2$, то
 $\widetilde{w}_0=\bar{s}^{\,\,-1}\cdot
q^{-\frac{H^2}{4}}$ . Можно показать \cite[стр. 425]{KR}, что
 $$
\widetilde{w}_0^{-1} X^+ \widetilde{w}_0 = -q^{-1}\, X^-,\qquad
\widetilde{w}_0^{-1} X^- \widetilde{w}_0 = -q\, X^+,
$$
\begin{equation}\label{w0_R_sl_2}
\Delta \widetilde{w}_0 =(\widetilde{w}_0\otimes\widetilde{w}_0) R,
\end{equation}
где $\Delta: \mathbb{C}[SL_2]_q^*\to (\mathbb{C}[SL_2]_q^{\otimes 2})^*$ и
элемент $R\in (\mathbb{C}[SL_2]_q^{\otimes 2})^*$ отвечает универсальной
$R$-матрице алгебры Хопфа $U_q\mathfrak{sl}_2$, см. \cite[стр. 425]{KR}.

Аналогично, в общем случае введем элемент $R$ алгебры
$(\mathbb{C}[G]_q^{\otimes 2})^*$, отвечающий универсальной $R$-матрице
алгебры Хопфа $U_q\mathfrak{g}$, ср. \cite[стр. 289]{Jo}.

\begin{proposition}(\cite[стр.
252]{LS})\label{w0_tilda_Prop} В алгебре $(\mathbb{C}[G]_q^{\otimes 2})^*$
имеет место равенство
\begin{equation}\label{w0_R_new}
\Delta \widetilde{w}_0 =(\widetilde{w}_0\otimes\widetilde{w}_0) R
\end{equation}
где $\Delta: \mathbb{C}[G]_q^*\to (\mathbb{C}[G]_q^{\otimes 2})^*$.
\end{proposition}

\begin{remark}\label{to_LS} В работах \cite{LS,KR} используется
$\widetilde{w}_0^{-1}$ вместо $\widetilde{w}_0$.
\end{remark}

Рассмотрим элементы $X_i^\pm$ алгебры $\mathbb{C}[G]_q^*$, определяемые
равенствами
\begin{equation}\label{from_EF_to_X}
X_i^+=E_iq_i^{-\frac{H_i}{2}},\qquad X_i^-=q_i^{\frac{H_i}{2}}F_i,
\end{equation}
см. \eqref{changesl_2} и замечание \ref{Change}.

Следующий результат получен Джозефом в \cite{Jo-R}.
\begin{proposition}\label{from_Jo-R} В алгебре
$\mathbb{C}[G]_q^*$
\begin{equation}\label{w0_Joseph}
\widetilde{w}_0^{-1}\,X_i^{\pm}\,\widetilde{w}_0\,=\, - q_i^{\mp 1}
X_{i'}^{\mp},\qquad i=1,2,\ldots,l,
\end{equation}
где $i'$ определяется равенством $\alpha_{i'}=- w_0 \alpha_i$.
\end{proposition}

\begin{remark}\label{to_Joseph} Выбор корректирующего картановского
множителя в \eqref{w0_tilda} является традиционным, см., например,
\cite[стр. 252]{LS}, и немногим отличается от используемого Джозефом, см.
\cite[стр. 422]{Jo-R}.
\end{remark}


\subsubsection{Фундаментальные представления и специальные
базисы.}\label{true_bases}

При изучении алгебр функций на компактных квантовых группах широко
используются матричные элементы представлений $\pi_\Lambda$ алгебры
$U_q\mathfrak{g}$, отвечающих простым весовым конечномерным
$U_q\mathfrak{g}$-модулям $L(\Lambda)$, $\Lambda\in P_+$, см. \itemiii\
\ref{CommRelations}. При этом неудобства доставляет неопределенность в
выборе ортонормированного базиса весовых векторов в $L(\Lambda)$.

   Даже в частном случае фундаментального веса $\Lambda$
и $\lambda,\mu\in W\Lambda$, $\lambda \neq \mu$, когда весовые
подпространства одномерны, матричный элемент $c^{\Lambda}_{\lambda,\mu}$,
определен лишь с точностью до числового множителя, по модулю равного
единице. Устраним эту неопределенность, уменьшив произвол в выборе
ортонормированного базиса весовых векторов.

\bigskip Пусть $\mathbb{S}=\{1,2,\ldots,l\}\setminus \{l_0\}$. Рассмотрим
подгруппу $W_{\mathbb{S}}\subset W$, порожденную простыми отражениями $s_j$,
$j\neq l_0$. Как известно \cite[стр. 22]{Hum2}, она является стационарной
подгруппой фундаментального веса
 $\overline{\omega}=\overline{\omega}_{l_0}$, то есть
 $ W_{\mathbb{S}}=\{w\in W\,|\, w \overline{\omega}=\overline{\omega}\}$.
 Следовательно, множество $W\overline{\omega}$ естественно изоморфно
 множеству  $^\mathbb{S}W$
элементов наименьшей длины в классах смежности $wW_\mathbb{S}$, $w \in W$,
см. лемму \ref{w_s}.

Напомним, что представление $\pi_{\overline{\omega}}$ $*$-алгебры
$(U_q\mathfrak{g},\star)$ в $L(\overline{\omega})$ унитаризуемо, и
соответствующая эрмитова форма $(\cdot,\cdot)$ однозначно определяется
условием нормировки $(v(\overline{\omega}),v(\overline{\omega}))=1$, см.
\itemiii\ \ref{loc_finite}.

Ортонормированный базис весовых векторов $U_q\mathfrak{g}$-модуля
$L(\overline{\omega})$ будем называть специальным, \IND{специальный базис}
если весовые векторы с весами из $W\overline{\omega}$, равны
$$
v_{w\overline{\omega}}^{\overline{\omega}}=
\frac{\overline{w}^{-1}v(\overline{\omega})}
{\left\|\overline{w}^{-1}v(\overline{\omega})\right\|},\qquad
w\in\,^\mathbb{S}W.
$$
Здесь $\overline{w}\in\mathbb{C}[G]_q^*$ -- квантовый аналог элемента $w$,
введенный в предыдущем \itemiiiе\ (использованы его обратимость и равенство
\eqref{w_action}).

Дадим эквивалентное определение специального базиса, ср. \cite[стр.
153]{LakshResh}. Пусть $w=s_{i_L}s_{i_{L-1}}\cdots s_{i_1}$ -- приведенное
разложение элемента $w\in\, ^\mathbb{S}W$. Положим
$u_j=s_{i_j}s_{i_{j-1}}\cdots s_{i_1}$,
$$
\lambda_j=
\begin{cases}
u_j\overline{\omega}, & j=1,2,\ldots,L,
\\ \overline{\omega}, & j=0,
\end{cases}\qquad \qquad
m_j=\frac{2(\lambda_j,\lambda_{j-1})}{(\lambda_j,\lambda_j)}.
$$
Как известно, $m_j\in\mathbb{N}$ \cite[стр. 93]{Bou4-6}. Следующее
утверждение вытекает из леммы \ref{simple_refl}.
\begin{proposition}\label{standard_def}
Ортонормированный базис весовых векторов является специальным, если и только
если для всех $w\in\, ^\mathbb{S}W$
$$
v_{w\overline{\omega}}^{\overline{\omega}}=
\frac{F_{i_L}^{m_L}F_{i_{L-1}}^{m_{L-1}}\cdots
F_{i_1}^{m_1}v(\overline{\omega})}
{\left\|F_{i_L}^{m_L}F_{i_{L-1}}^{m_{L-1}}\cdots
F_{i_1}^{m_1}v(\overline{\omega})\right\|}.
$$
\end{proposition}

\medskip
\begin{proposition}\label{U_standard}
Для матричных элементов фундаментальных представлений в специальных базисах
имеет место следующее соотношение:
\begin{eqnarray}
S(c_{w'\overline{\omega},w''\overline{\omega}}^{\overline{\omega}}) &=&
(-1)^{l(w')-l(w'')}q^{(\mu-\lambda,\rho)}
c_{-w''\overline{\omega},-w'\overline{\omega}}^{-w_0\overline{\omega}}\;,
\label{antipode_standard}
\\ (c_{w'\overline{\omega},w''\overline{\omega}}^{\overline{\omega}})^\star
&=& (-1)^{l(w')-l(w'')}q^{(\lambda-\mu,\rho)}
c_{-w'\overline{\omega},-w''\overline{\omega}}^{-w_0\overline{\omega}}\;,
\label{aster_standard}
\end{eqnarray}
где $w',w''\in\, ^\mathbb{S}W$, а $\rho$ -- полусумма положительных корней.
\end{proposition}

{\bf Доказательство.} Пусть $\pi_{\overline{\omega}}$,
$\pi_{-w_0\overline{\omega}}$ -- фундаментальные представления алгебры
$U_q\mathfrak{g}$, отвечающие $U_q\mathfrak{g}$-модулям
$L(\overline{\omega})$, $L(-w_0\overline{\omega})$ соответственно.
Рассмотрим инволютивный антилинейный антиавтоморфизм $\tau$ алгебры
$U_q\mathfrak{g}$ и представление $\pi_{\overline{\omega}}^\tau$ этой
алгебры в $L(\overline{\omega})$, описанные в \itemiiiе \ref{regular_K}, где
показано, что существует и единствен унитарный оператор
$$
U:L(\overline{\omega})\to L(-w_0\overline{\omega}),\qquad
U:v_{\overline{\omega}}^{\overline{\omega}}\mapsto
v_{-\overline{\omega}}^{-w_0\overline{\omega}},
$$
осуществляющий эквивалентность представлений $\pi_{\overline{\omega}}^\tau$
и $\pi_{-w_0\overline{\omega}}$. Из предложения \ref{standard_def} и из
определения инволюции $\tau$ вытекает равенство
$$
Uv_{w\overline{\omega}}^{\overline{\omega}}=
(-1)^{l(w)}v_{-w\overline{\omega}}^{-w_0\overline{\omega}},\qquad w\in
^\mathbb{S}W.
$$
Остается воспользоваться предложением \ref{invol_c}. \hfill $\square$

\bigskip

В заключение рассмотрим важный частный случай фундаментального веса
$\overline{\omega}$, которому отвечает простой корень, входящий в
разложение старшего корня с коэффициентом 1. Кроме того, ограничимся
алгебрами Ли с корнями равной длины (сериями $A$, $D$, $E$). В этих
предположениях $\overline{\omega}$ является микровесом, \IND{микровес} то
есть все веса $U_q\mathfrak{g}$-модуля $L(\overline{\omega})$ лежат на его
$W$-орбите \cite[стр. 164]{Bou4-6}. Равенства \eqref{antipode_standard},
\eqref{aster_standard} принимают следующий вид
\begin{eqnarray}
S(c_{\lambda,\mu}^{\overline{\omega}}) &=& (-q)^{(\mu-\lambda,\rho)}\,
c_{-\mu,-\lambda}^{-w_0\overline{\omega}},\qquad\lambda,\mu\in
W\overline{\omega},\label{antipode_ADE}
\\ (c_{\lambda,\mu}^{\overline{\omega}})^\star &=&
(-q)^{(\lambda-\mu,\rho)}\,c_{-\lambda,-\mu}^{-w_0\overline{\omega}},\qquad
\lambda,\mu\in W\overline{\omega},\label{aster_ADE}
\end{eqnarray}
поскольку из $l(s_iw)=l(w)+1$, $s_iw\in\,^\mathbb{S}W$ следует, что
$$
(s_i\,w\overline{\omega}-w\overline{\omega},\rho)=
-\left(\alpha_i,\sum_{j=1}^l\overline{\omega}_j\right)=-1.
$$


\subsubsection{Каноническое вложение $\mathbb{C}[\mathfrak{p}^-]_q
\hookrightarrow\mathbb{C}[X^-_{\mathbb{S}}]_{q,t}$.}\label{w0_G_x}

 По-прежнему предполагаем, что
$\mathbb{S}=\{1,2,\ldots,l\}\setminus\{l_0\},$ и что простой корень
$\alpha_{l_0}$ входит в разложение старшего корня с коэффициентом 1. Введем
в рассмотрение фундаментальный вес
$$\overline{\omega}_{l_1}=-w_0\overline{\omega}_{l_0}$$
 и отвечающее ему фундаментальное представление $\pi_{l_1}$
 алгебры $U_q\mathfrak{g}$
 в векторном пространстве $L(\overline{\omega}_{l_1})$.

Воспользуемся специальным базисом в $L(\overline{\omega}_{l_1})$ и
обозначениями для матричных элементов представлений алгебры
$U_q\mathfrak{g}$, введенными в
 \itemiiiе \ref{CommRelations}.
Рассмотрим следующий матричный элемент представления $\pi_{l_1}$:
\begin{equation}\label{t_def}
t= c_{\overline{\omega}_{l_1},-\overline{\omega}_{l_0}}^
{\overline{\omega}_{l_1}}.
\end{equation}
 Из следствия \ref{subalgebra} вытекает

\begin{proposition}\label{t_tstar}
1. $tt^\star=t^\star t$,

2. Пусть $\Lambda\in P_+$ и $\lambda,\mu$ -- веса $L(\Lambda)$. Матричные
элементы $c_{\lambda,j;\mu,k}^{L(\Lambda)}$ алгебры $\mathbb{C}[G]_q$
квазикоммутируют с $t$, $t^\star$. Именно,
$$
t\,c_{\lambda,j;\mu,k}^{L(\Lambda)}=\mathrm{const}_1\;
c_{\lambda,j;\mu,k}^{L(\Lambda)}\,t,\qquad t^\star\,
c_{\lambda,j;\mu,k}^{L(\Lambda)}=\mathrm{const}_2\;
c_{\lambda,j;\mu,k}^{L(\Lambda)}\,t^\star,
$$
где $\mathrm{const}_1,\mathrm{const}_2\in q^{\frac1s\mathbb{Z}}$ зависят от
$\Lambda$, $\lambda$, $\mu$, $q$ и $s$ -- порядок группы $P/Q$.
\end{proposition}

\medskip

Пусть $x=tt^\star=t^\star t$.

\begin{corollary}\label{quasicom_x}
Если $\Lambda\in P_+$ и $\lambda,\mu$ -- веса $L(\Lambda)$, то
\begin{equation}\label{x_c}
x\,c_{\lambda,j;\mu,k}^{L(\Lambda)}=\mathrm{const}\;
c_{\lambda,j;\mu,k}^{L(\Lambda)}\,x,
\end{equation}
где число $\mathrm{const}$ зависит от $\Lambda$, $\lambda$, $\mu$, $q$ и
принадлежит $q^{\frac1s\mathbb{Z}}$.
\end{corollary}

   Мультипликативное подмножество $x^{\mathbb{Z}_+}$ является множеством
   Орэ,
см. \eqref{x_c}. Пусть $\mathbb{C}[G]_{q,x}$ -- локализация алгебры
$\mathbb{C}[G]_q$ по этому мультипликативному подмножеству. Из предложения
\ref{dom_noether} следует, что
$\mathbb{C}[G]_q\hookrightarrow\mathbb{C}[G]_{q,x}$.

\begin{proposition}\label{*_local_gen}
\begin{itemize}
\item[1.] Существует и единственно продолжение структуры
$U_q\mathfrak{g}$-модульной алгебры с $\mathbb{C}[G]_q$ на
$\mathbb{C}[G]_{q,x}$.

\item[2.] Для любых $\xi\in U_q\mathfrak{g}$, $f\in\mathbb{C}[G]_{q,x}$
существует и единствен такой полином Лорана $p_{f,\xi}$ одной переменной с
коэффициентами из $\mathbb{C}[G]_{q,x}$, что
\begin{equation}\label{p_f_gen}
\xi(f\cdot x^n)=p_{f,\xi}(q^\frac{n}s)\,x^n
\end{equation}
при всех $n\in\mathbb{Z}$.
\end{itemize}
\end{proposition}

{\bf Доказательство.} Начнем со второго утверждения. Используя то, что
$\mathbb{C}[G]_q$ является $U_q\mathfrak{g}$-модульной алгеброй, нетрудно
свести общий случай $\xi\in U_q\mathfrak{g}$ к частному случаю
$\xi\in\{K_j^{\pm 1},E_j,F_j|\:j=1,2,\ldots,l\}$. В этом частном случае
существование полинома Лорана $p_{f,\xi}$ вытекает из следствия
\ref{quasicom_x} и явного вида элементов $\triangle(K_j^{\pm 1})$,
$\triangle(E_j)$, $\triangle(F_j)$, см. \eqref{comult_D}. Остальные
утверждения предложения \ref{*_local_gen} доказываются точно так же, как
соответствующие утверждения предложения \ref{*_local}. \hfill $\square$

\medskip
Получаем хорошо известное утверждение.
\begin{corollary} Мультипликативное подмножество $t^{\mathbb{Z}_+}$ элементов
алгебры $\mathbb{C}[G]_q$ является множеством Орэ, алгебра $\mathbb{C}[G]_q$
вкладывается в локализацию $\mathbb{C}[G]_{q,t}$, и структура
$U_q\mathfrak{g}$-модульной алгебры единственным образом продолжается с
$\mathbb{C}[G]_q$ на $\mathbb{C}[G]_{q,t}$.
\end{corollary}

\begin{remark}\label{Rosenberg}
Другой подход к доказательству того, что структура
$U_q\mathfrak{g}$-модульной алгебры допускает продолжение на локализацию,
описан в работах В.~Лунца и А.~Розенберга \cite{LuntsRosLoc},
\cite{LuntsRosDiff2}.
\end{remark}

Введем в рассмотрение
 матричный элемент
\begin{equation}\label{t_prime}
t'=c_{\overline{\omega}_{l_1},-\overline{\omega}_{l_0}+\alpha_{l_0}}^
{\overline{\omega}_{l_1}}.
\end{equation}
 фундаментального представления $\pi_{l_1}$
алгебры $U_q\mathfrak{g}$ в специальном базисе векторного пространства
$L(\overline{\omega}_{l_1})$. Отметим, что веса $-\overline{\omega}_{l_0}$,
$-\overline{\omega}_{l_0}+\alpha_{l_0}$ принадлежат $W$-орбите старшего
веса $\overline{\omega}_{l_1}$, и, следовательно, элементы $t$, $t'$ не
зависят от выбора специального базиса.

Из определений следует, что при всех $j\in\{1,2,\cdots,l\}$

\begin{gather}
F_j\,t=0,\qquad K^{\pm 1}_j\,t=
\begin{cases}
q_j^{\mp 1}\,t, & j=l_0,
\\ 0, & j\ne l_0.
\end{cases} \qquad
 E_j\,t=
\begin{cases}
t', & j=l_0,
\\ 0, & j\ne l_0.
\end{cases}\label{FHE_t}
\end{gather}

\begin{proposition}\label{hol_embed}
Отображение $i:z_\mathrm{low}\mapsto t^{-1}t'$ единственным образом
продолжается до вложения $U_q\mathfrak{g}$-модульных алгебр
\begin{equation}\label{hom_i}
i:\mathbb{C}[\mathfrak{p}^-]_q\hookrightarrow\mathbb{C}[G]_{q,x}.
\end{equation}
\end{proposition}

{\bf Доказательство.} Рассматриваемая пара $(\mathfrak{g},\mathfrak{k})$
является эрмитово-симметрической. Значит, единственность гомоморфизма
\eqref{hom_i} следует из простоты $U_q\mathfrak{k}$-модуля
$\mathbb{C}[\mathfrak{p}^-]_{q,1}$ (cм. \itemiii \ \ref{vect}) и из
предложения \ref{hom_gen}, согласно которому подпространство
$\mathbb{C}[\mathfrak{p}^-]_{q,1}$ порождает алгебру с единицей
$\mathbb{C}[\mathfrak{p}^-]_q$.

Перейдем к доказательству существования гомоморфизма \eqref{hom_i}.
Рассмотрим $U_q\mathfrak{g}$-модуль $\widetilde{w}_0U_q\mathfrak{g}$,
двойственный $U_q\mathfrak{g}$-модуль $(\widetilde{w}_0U_q\mathfrak{g})^*$
и его наибольший весовой подмодуль
$\mathbb{F}\subset(\widetilde{w}_0U_q\mathfrak{g})^*$, состоящий из
$U_q\mathfrak{b}^-$-конечных элементов. Здесь, как обычно,
$U_q\mathfrak{b}^-$ -- подалгебра Хопфа, порожденная элементами множества
$\{K_j^\pm,F_j\}_{j=1,2,\cdots,l}$. Иначе говоря,
$\mathbb{F}=\bigoplus\limits_{\lambda\in P}\mathbb{F}_\lambda$,
$$
\mathbb{F}_\lambda=\left\{\left.f\in(\widetilde{w}_0U_q\mathfrak{g})^*
\right|\:\dim(U_q\mathfrak{b}^-f)<\infty,\quad H_if=\lambda(H_i)f,\;
i=1,\ldots,l\right\}.
$$
Линейное подмногообразие $\mathbb{C}[G]_q$ плотно в $\mathbb{F}$ в слабой
топологии. Используя \eqref{w0_R_new}, нетрудно доказать, что структура
$U_q\mathfrak{g}$-модульной алгебры допускает продолжение по непрерывности
с $\mathbb{C}[G]_q$ на $\mathbb{F}$. Действительно,
\begin{eqnarray*}
\langle f_1f_2,\widetilde{w}_0\xi\rangle &=& \langle
R_{\mathbb{F}\mathbb{F}}\triangle(\xi)(f_1\otimes
f_2),\widetilde{w}_0\otimes\widetilde{w}_0\rangle,
\\ \langle\xi f,\widetilde{w}_0\eta\rangle &=& \langle
f,\widetilde{w}_0\eta\xi\rangle,\qquad\xi,\eta\in U_q\mathfrak{g},\qquad
f,f_1,f_2\in\mathbb{F}.
\end{eqnarray*}
Здесь $R_{\mathbb{F}\mathbb{F}}$ -- линейный оператор в
$\mathbb{F}\otimes\mathbb{F}$, определяемый действием универсальной
R-матрицы.

Как отмечалось при доказательстве единственности гомоморфизма \eqref{hom_i},
     элемент
$z_\mathrm{low}$ порождает $U_q\mathfrak{k}$-модульную алгебру
$\mathbb{C}[\mathfrak{p}^-]_q$. Значит, существование вложения \eqref{hom_i}
будет доказано, если удастся построить такие вложения
$U_q\mathfrak{g}$-модульных алгебр
$$
i_1:\mathbb{C}[G]_{q,t}\hookrightarrow\mathbb{F},\qquad
i_2:\mathbb{C}[\mathfrak{p}^-]_q\hookrightarrow\mathbb{F},
$$
что $i_1(t^{-1}t')=i_2(z_{\operatorname{low}})$.

Построим вложение $i_1$, продолжив естественное спаривание
$\mathbb{C}[G]_q\times\widetilde{w}_0U_q\mathfrak{g}\to\mathbb{C}$ до
спаривания $\mathbb{C}[G]_{q,t}\times\widetilde{w}_0U_q\mathfrak{g}\to
\mathbb{C}$. В классическом случае $q=1$ такое продолжение возможно
благодаря тому, что $w_0$ принадлежит области $t\ne 0$, в которой регулярны
все функции из $\mathbb{C}[G]_{q,t}$. В квантовом случае будет использовано
то, что
\begin{equation}\label{w0_t}
\langle t,\widetilde{w}_0\rangle\ne 0,
\end{equation}
как следует из \eqref{w_action}.

Пусть $f\in\mathbb{C}[G]_{q,\lambda}$ -- весовой вектор веса $\lambda\in P$
и $c_k(f)=\langle ft^k,\widetilde{w}_0\rangle$. Используя равенства
$$
\left\langle ft^k,\widetilde{w}_0\right\rangle=\left\langle
R_{\mathbb{C}[G]_q\mathbb{C}[G]_q}(ft^{k-1}\otimes
t),\widetilde{w}_0\otimes\widetilde{w}_0\right\rangle,
$$
\eqref{FHE_t}, \eqref{Rmatrix}, \eqref{fraction}, нетрудно доказать, что
числовая последовательность $\{c_k(f)\}$, $k\in\mathbb{Z}_+$, является
решением разностного уравнения первого порядка с коэффициентами, не
зависящими от $f$. Именно,
$$
c_k(f)=
q^{\mathrm{const}_1\cdot(k-1)+\mathrm{const}_2\cdot\lambda}\cdot\langle
t,\widetilde{w}_0\rangle c_{k-1}(t),\qquad k\in\mathbb{N}.
$$
Из \eqref{w0_t} следует, что каждое решение этого уравнения единственным
образом продолжается до его решения $\{c_k\}$, определенного при всех
$k\in\mathbb{Z}$. Получаем ''естественное'' продолжение линейного
функционала $f\mapsto\langle f,\widetilde{w}_0\rangle$ с $\mathbb{C}[G]_q$
на $\mathbb{C}[G]_{q,t}$.

Из предложения \ref{*_local_gen} следует, что подалгебра
$\mathbb{C}[G]_{q,t}\subset\mathbb{C}[G]_{q,x}$ является
$U_q\mathfrak{g}$-модульной. Введем в рассмотрение спаривание
\begin{equation}\label{pairing_t}
\mathbb{C}[G]_{q,t}\times\widetilde{w}_0U_q\mathfrak{g}\to\mathbb{C},\qquad
f\times\widetilde{w}_0\xi\mapsto\langle\xi f,\widetilde{w}_0\rangle
\end{equation}
и отвечающий ему морфизм $U_q\mathfrak{g}$-модулей
$$
i_1:\mathbb{C}[G]_{q,t}\to(\widetilde{w}_0U_q\mathfrak{g})^*,\qquad
i_1:f\mapsto\langle f,\cdot\rangle.
$$
Докажем инъективность $i_1$. Из обратимости элемента
$\widetilde{w}_0\in\mathbb{C}[G]_q^*$ следует, что линейное подмногообразие
$\widetilde{w}_0U_q\mathfrak{g}$ плотно в $\mathbb{C}[G]_q^*$ в слабой
топологии. Значит, $\operatorname{Ker}i_1\cap\mathbb{C}[G]_q=0$. Остается
заметить, что из $i_1(f)=0$ следует равенство $i_1(f\cdot t^j)=0$ при всех
$j\in\mathbb{N}$, поскольку
$$
\langle f\cdot t^j,\widetilde{w}_0\xi\rangle=\left\langle
R_{\mathbb{C}[G]_{q,t}\mathbb{C}[G]_{q,t}}\triangle(\xi)(f\otimes
t^j),\widetilde{w}_0\otimes\widetilde{w}_0\right\rangle.
$$

Отметим, что $i_1\mathbb{C}[G]_{q,t}\subset\mathbb{F}$, поскольку $i_1$
является морфизмом $U_q\mathfrak{g}$-модулей и все элементы
$f\in\mathbb{C}[G]_{q,t}$ являются
$U_q\mathfrak{b}^-$-конечными.\footnote{Последнее утверждение достаточно
доказать для $f\in\mathbb{C}[G]_q$ и $f=t^{-1}$. Но
$\mathbb{C}[G]_q=\bigoplus\limits_{\lambda\in
P_+}\operatorname{End}L(\lambda)$, а $\dim(U_q\mathfrak{b}^-\cdot
t^{-1})=1$, как следует из \eqref{FHE_t}}

Докажем, что полученное вложение
$i_1:\mathbb{C}[G]_{q,t}\hookrightarrow\mathbb{F}$ является гомоморфизмом
алгебр. Достаточно показать, что для каждой пары весовых векторов
$f_1,f_2\in\mathbb{C}[G]_q$ следующее равенство
\begin{equation}\label{eq_f1f2}
\langle(f_1\cdot t^k)(f_2\cdot t^k),\widetilde{w}_0\rangle=\left\langle
R_{\mathbb{C}[G]_{q,t}\mathbb{C}[G]_{q,t}}(f_1\cdot t^k\otimes f_2\cdot
t^k),\widetilde{w}_0\otimes\widetilde{w}_0\right\rangle
\end{equation}
имеет место при всех $k\in\mathbb{Z}$. Так как последовательности чисел в
обеих частях равенства \ref{eq_f1f2} являются линейными комбинациями
последовательностей вида $q^{ak^2+bk}$, $a,b\in\mathbb{Q}$, то достаточно
доказать \eqref{eq_f1f2} при больших $k$. Например, можно ограничиться
частным случаем $f_1,f_2\in\mathbb{C}[G]_q$. В этом случае \eqref{eq_f1f2}
следует из \eqref{w0_R_new}.

Вложение $U_q\mathfrak{g}$-модульных алгебр $i_1$ построено. Перейдем к
построению вложения $i_2$. Пусть
$i_2:\mathbb{C}[\mathfrak{p}^-]_q\to(\widetilde{w}_0U_q\mathfrak{g})^*$ --
линейный оператор, сопряженный к
$$
\mathscr{I}:\widetilde{w}_0U_q\mathfrak{g}\to
N(\mathfrak{q}^+,0),\qquad\mathscr{I}:\widetilde{w}_0\xi\mapsto
S(\xi)v(\mathfrak{q}^+,0).
$$
Из определений следует, что $i_2$ является инъективным морфизмом
$U_q\mathfrak{g}$-модулей. Значит,
$i_2:\mathbb{C}[\mathfrak{p}^-]_q\hookrightarrow\mathbb{F}$, поскольку
$U_q\mathfrak{g}$-модуль $\mathbb{C}[\mathfrak{p}^-]_q$ является локально
$U_q\mathfrak{b}^-$-конечномерным.

Из \eqref{Rmatrix} и из того, что $v(\mathfrak{q}^+,0)$ является старшим
вектором нулевого веса, следуют равенства

\begin{multline*}
\langle\xi(f_1f_2),v(\mathfrak{q}^+,0)\rangle=
\langle\triangle(\xi)(f_1\otimes f_2),v(\mathfrak{q}^+,0)\otimes
v(\mathfrak{q}^+,0)\rangle=
\\ =\langle R_{\mathbb{C}[\mathfrak{p}^-]_q\mathbb{C}[\mathfrak{p}^-]_q}
\triangle(\xi)(f_1\otimes f_2),v(\mathfrak{q}^+,0)\otimes
v(\mathfrak{q}^+,0)\rangle=
\\ =\langle
R_{\mathbb{F}\mathbb{F}}\triangle(\xi)(i_2(f_1)\otimes
i_2(f_2)),\widetilde{w}\otimes\widetilde{w}\rangle=
\langle\xi(i_2(f_1)i_2(f_2),\widetilde{w}\rangle,
\end{multline*}
где $f_1,f_2\in\mathbb{C}[\mathfrak{p}^-]_q$, $\xi\in U_q\mathfrak{g}$.
Значит, $i_2$ является гомоморфизмом алгебр.

Покажем, что элемент $i_1(t^{-1}t')$ принадлежит
$i_2\mathbb{C}[\mathfrak{p}^-]_q$. Достаточно доказать ортогональность
вектора $i_1(t^{-1}t')$ и подпространства
$\operatorname{Ker}\mathscr{I}=\widetilde{w}_0(U_q\mathfrak{q}^+)_0$, где
\hbox{$(U_q\mathfrak{q}^+)_0=
U_q\mathfrak{q}^+\cap\operatorname{Ker}\varepsilon$}, а $\varepsilon$ --
коединица алгебры Хопфа $U_q\mathfrak{q}^+$. Рассмотрим ''регулярное''
представление $L_{\operatorname{reg}}$ алгебры $U_q\mathfrak{g}^{\rm op}$ в
пространстве $\mathbb{C}[G]_q$, см. \itemiii\ \ref{alg_gr}. Так же, как при
доказательстве предложения \ref{*_local_gen}, продолжим операторы
$L_\mathrm{reg}(\xi)$ на $\mathbb{C}[G]_{q,t}$. Тем самым
$\mathbb{C}[G]_{q,t}$ наделяется структурой
$U_q\mathfrak{g}^\mathrm{op}$-модуля. Если $\xi \in (U_q\mathfrak{q}^+)_0$,
то элемент $\eta= \widetilde{w}_0 \cdot \xi \cdot (\widetilde{w}_0)^{-1}$
алгебры $(\mathbb{C}[G]_q)^*$ принадлежит подалгебре $U_q\mathfrak{g}$
 и имеет место равенство
\begin{equation}\label{last_eq_embed}
L_\mathrm{reg}(\eta)(t^{-1}t')=0.
\end{equation}
Действительно, при доказательстве этого равенства элемент $\xi$ можно
считать однородным положительной степени. Наделим
$U_q\mathfrak{g}^{op}$-модуль $V=\{L_\mathrm{reg}(\zeta)(t^{-1}t')| \zeta
\in U_q\mathfrak{g}^{op}\}$ градуировкой, отвечающей подмножеству
  $\mathbb{S}_1=\{1,2, \cdots, l\} \setminus \{l_1\}$. Нетрудно показать,
  что элемент
$L_{\operatorname{reg}}(\eta) (t^{-1}t') \in V$ имеет недопустимую степень
однородности и, следовательно, равен нулю.

 Из \eqref{last_eq_embed} вытекает требуемая ортогональность.

\medskip Доказательство предложения \ref{hol_embed} будет завершено, если
доказать, что элемент $z$, для которого $i_1(t^{-1}t')=i_2(z)$, равен
  $z_\mathrm{low}$.
Напомним, что $\mathbb{S}=\{1,2,\ldots,l\}\setminus\{l_0\}$, и что
$U_q\mathfrak{g}$-модуль $\mathbb{C}[G]_{q,t}$ является весовым.
Следовательно, корректно определен линейный оператор $H_\mathbb{S}$ в
$\mathbb{C}[G]_{q,t}$. Из равенств
\begin{eqnarray*}
H_\mathbb{S}t &=& -\overline{\omega}_{l_0}(H_\mathbb{S})t,\qquad
F_jt=0,\quad j=1,2,\ldots,l,
\\ H_\mathbb{S}t' &=&
(-\overline{\omega}_{l_0}(H_\mathbb{S})+\alpha_{l_0}(H_\mathbb{S}))t',
\qquad F_jt=0,\quad j\ne l_0,
\end{eqnarray*}
следует, что
$$
H_\mathbb{S}(t^{-1}t')=2t^{-1}t',\qquad F_j(t^{-1}t')=0,\quad j\ne l_0.
$$
Значит, элемент $z\in\mathbb{C}[\mathfrak{p}^-]_q$, для которого
$i_1(t^{-1}t')=i_2(z)$, является младшим вектором $U_q\mathfrak{k}$-модуля
$\mathbb{C}[\mathfrak{p}^-]_{q,1}$. Как показывают несложные вычисления,
\begin{equation}\label{t_prime_second}
t'=q_{l_0}^{\frac{1}{2}}E_{l_0} t,
\end{equation}
 $F_{l_0}(t^{-1}t')=q_{l_0}^{\frac{1}{2}}$. Следовательно,
$F_{l_0}(z)=q_{l_0}^{1/2}$, и равенство $z=z_{\operatorname{low}}$ вытекает
из леммы \ref{z1}. \hfill $\square$

\bigskip Полученное вложение
 $i: \mathbb{C}[\mathfrak{p}^-]_q \hookrightarrow \mathbb{C}[G]_{q,x}$
 будем называть каноническим.


\begin{proposition}\label{image_hol}
Образом алгебры $\mathbb{C}[\mathfrak{p}^-]_q$ при каноническом вложении в
$\mathbb{C}[G]_{q,x}$ является подалгебра, порожденная множеством
$\{t^{-1}\,c_{\overline{\omega}_{l_1},\lambda}^{\overline{\omega}_{l_1}}\}$,
где $c_{\overline{\omega}_{l_1},\lambda}^{\overline{\omega}_{l_1}}$ --
элементы ''старшей'' строки матрицы фундаментального представления
$\pi_{l_1}$.
\end{proposition}

{\bf Доказательство.} Пусть $\mathcal{F} \subset \mathbb{C}[G]_{q,x}$ --
подалгебра, порожденная элементами $\{
t^{-1}c_{\overline{\omega}_{l_1},\lambda} ^{\overline{\omega}_{l_1}} \}$. Из
определений следует, что она является $U_q\mathfrak{g}$-модульной
подалгеброй.

Докажем, что $i\mathbb{C}[\mathfrak{p}^-]_q=\mathcal{F}$. Прежде всего,
$i\mathbb{C}[\mathfrak{p}^-]_q\subset\mathcal{F}$, поскольку
$i(z_\mathrm{low})\in\mathcal{F}$ и, следовательно,
$i\mathbb{C}[\mathfrak{p}^-]_{q,1} \subset\mathcal{F}$. Обратное включение
$i\mathbb{C}[\mathfrak{p}^-]_q\supset\mathcal{F}$ нетрудно установить,
используя то, что линейная оболочка $\mathcal{L}$ множества
$\{c_{\overline{\omega}_{l_1},\lambda}^{\overline{\omega}_{l_1}}\}$
совпадает с $U_q\mathfrak{b}^+\cdot t$, поскольку элемент $t$ является
младшим вектором $U_q\mathfrak{g}$-модуля
 $\mathcal{L} \approx L(\overline{\omega}_{l_0})$. \hfill
$\square$

\begin{remark}\label{mult_free} Используя каноническое вложение и известную
теорему Хуа-Шмида (см., например, \itemiii\ \ref{Hua-Schmid}), нетрудно
показать, что
$$
 \dim\operatorname{Hom}_{U_q\mathfrak{k}}
(L(\mathfrak{k},\lambda),L(j\overline{\omega}_{l_0}))\le 1
$$
 при всех  $\lambda\in P_+^\mathbb{S}$, $j\in\mathbb{N}$.
\end{remark}

\bigskip

Пусть $\mathbb{C}[X^-_{\mathbb{S}}]_q\subset\mathbb{C}[G]_q$-- наименьшая
$U_q\mathfrak{g}$-модульная подалгебра, содержащая элемент $t$, а
$\mathbb{C}[X^+_{\mathbb{S}}]_q\subset\mathbb{C}[G]_q$-- наименьшая
$U_q\mathfrak{g}$-модульная подалгебра, содержащая элемент $t^\star$.
Другими словами, подалгебра $\mathbb{C}[X^-_{\mathbb{S}}]_q$ порождена
матричными элементами
$\left\{c_{\overline{\omega}_{l_1};\lambda,j}^{\overline{\omega}_{l_1}}
\right\}$, а подалгебра $\mathbb{C}[X^+_{\mathbb{S}}]_q$-- матричными
элементами
$\left\{c_{w_0\overline{\omega}_{l_0};\lambda,j}^{\overline{\omega}_{l_0}}
\right\}$, см. \eqref{aster_standard}.

Наделим алгебру $\mathbb{C}[X^-_{\mathbb{S}}]_q$ градуировкой, полагая
\hbox{$\deg\left(c_{\overline{\omega}_{l_1};\lambda,j}^{\overline{\omega}_{l_1}}
\right)=1$.} Определение коректно, поскольку эту градуировку можно задать
по-другому, с помощью линейного оператора
$L_\mathrm{reg}(H_{\mathbb{S}_1})$. Именно, $\deg f=j$, если и только если
$L_\mathrm{reg}(H_{\mathbb{S}_1})f=jf$.

\medskip
Из полученных результатов вытекают следующие известные в теории квантовых
групп утверждения.

Во-первых, мультипликативное подмножество $t^{\mathbb{Z}_+}$ алгебры
$\mathbb{C}[X^-_\mathbb{S}]_q$ является множеством Орэ, и эта алгебра
естественно вкладывается в свою локализацию: $\mathbb{C}[X^-_\mathbb{S}]_q
\hookrightarrow \mathbb{C}[X^-_\mathbb{S}]_{q,t}$.
 Во-вторых, подалгебра
 $\mathbb{C}[X^-_{\mathbb{S}}]_q$
является $U_q\mathfrak{g}$-модульной.
 В-третьих,  структура $U_q\mathfrak{g}$-модульной алгебры допускает
 продолжение с
$\mathbb{C}[X^-_{\mathbb{S}}]_q$ на $\mathbb{C}[X^-_{\mathbb{S}}]_{q,t}$, и
такое продолжение единственно.

\medskip
  С помощью предложения \ref{hol_embed} получаем

\begin{corollary}\label{embed_to_flags}
 Отображение $i:z_\mathrm{low}\mapsto t^{-1}t'$ единственным образом
продолжается до вложения $U_q\mathfrak{g}$-модульных алгебр
$i:\mathbb{C}[\mathfrak{p}^-]_q\hookrightarrow
\mathbb{C}[X^-_{\mathbb{S}}]_{q,t}$.
\end{corollary}

\begin{remark}\label{embed-star} Аналогично можно наделить
локализацию $\mathbb{C}[X_\mathbb{S}^+]_{q,t^\star}$ алгебры
$\mathbb{C}[X_\mathbb{S}^+]_q$ по мультипликативному множеству
$(t^\star)^{\mathbb{Z}_+}$ структурой $U_q\mathfrak{g}$-модульной алгебры и
получить вложение $\mathbb{C}[\mathfrak{p}^+]_q\hookrightarrow
\mathbb{C}[X_\mathbb{S}^+]_{q,t^\star}$.
\end{remark}

Обратимся к классическому случаю $q=1$. Будет использован изоморфизм
$U\mathfrak{g}$-модулей $L(\overline{\omega}_{l_1})^*\approx
L(\overline{\omega}_{l_0})$. Рассмотрим конус $X^-_\mathbb{S}\subset
L(\overline{\omega}_{l_0})$, порожденный $G$-орбитой старшего весового
вектора $v(\overline{\omega}_{l_0})\in L(\overline{\omega}_{l_0})$.
 Его проективизация  обсуждалась
в \itemiiiе \ref{alg_gr}. Она является однородным пространством группы $G$,
а стабилизатор точки $\mathbb{C}v(\overline{\omega}_{l_0})$ -- стандартной
параболической подгруппой \cite[стр. 291]{Hum}.

Нами получен квантовый аналог вложения $\mathfrak{p}^-$ в это однородное
пространство, называемое обобщенным пространством флагов. В классическом
случае $q=1$ такое вложение хорошо известно и широко используется
\cite[стр. 34, 35]{Springer}. В частном случае
$\mathfrak{g}=\mathfrak{sl}_N$ близкий результат получен Фиорези в статье
\cite{Fioresi_cell}.


\subsubsection{$*$-Алгебра \boldmath
$(\mathbb{C}[\mathbb{X}_\mathbb{S}]_{q},*)$.}\label{star_emb}

Рассмотрим наименьшую унитальную $U_q\mathfrak{g}$-модульную подалгебру
$\mathbb{C}[\mathbb{X}_\mathbb{S}]_q\subset\mathbb{C}[G]_q$, содержащую
элементы $t$, $t^\star$.

Известно \cite[стр. 480]{Stokman_list}, что $
\mathbb{C}[\mathbb{X}_\mathbb{S}]_q$ является множеством тех
$f\in\mathbb{C}[G]_q$, для которых|
$$
L_\mathrm{reg}(E_j)f=L_\mathrm{reg}(F_j)f=L_\mathrm{reg}(K_j^{\pm 1}-1)f=0,
\qquad j\neq l_1\}.
$$

\begin{example} Если $\mathfrak{g}=\mathfrak{sl}_2$, то
$$\mathbb{C}[\mathbb{X}_\mathbb{S}]_q=\mathbb{C}[SL_2]_q.$$
\end{example}

\begin{proposition} \label{star_X_new}
Существует и единственна такая антилинейная инволюция $*$ в
$\mathbb{C}[\mathbb{X}_\mathbb{S}]_q$, что
$(\mathbb{C}[\mathbb{X}_\mathbb{S}]_q,*)$ является
$(U_q\mathfrak{g},*)$-модульной \hbox{$*$-алгеброй} и $t^*=t^\star$.
\end{proposition}

{\bf Доказательство.} Единственность требуемой инволюции очевидна. Ее
существование будет доказано, если удастся наделить $\mathbb{C}[G]_q$ такой
структурой $(U_q\mathfrak{g},*)$-модульной алгебры, что $t^*=t^\star$.
    Воспользуемся естественными вложениями
     $U_q\mathfrak{g}\hookrightarrow\mathbb{C}[G]_q^*$,
$\widetilde{w}_0U_q\mathfrak{g}\hookrightarrow\mathbb{C}[G]_q^*$. Из
обратимости элемента $\widetilde{w}_0$ следует, что образ
$\widetilde{w}_0U_q\mathfrak{g}$ при этом вложении плотен в
\hbox{$\mathbb{C}[G]_q^*\cong\mathop{\times}\limits_{\lambda \in P_+}
\operatorname{End}L(\lambda)$} в топологии прямого произведения.
Следовательно, каноническое спаривание $\langle\cdot,\cdot\rangle$
$$\mathbb{C}[G]_q\times\widetilde{w}_0U_q\mathfrak{g}\to\mathbb{C}$$
невырождено. Прервем доказательство предложения \ref{star_X_new} для того,
чтобы получить ряд вспомогательных результатов.

\begin{lemma}\label{def_star}
Существует и единственно такое антилинейное отображение $*$ в
$\mathbb{C}[G]_q$, что
\begin{equation}\label{dual_star}
\langle f^*,\widetilde{w}_0\xi\rangle=\overline{\langle
f,\widetilde{w}_0(S(\xi))^*\rangle}
\end{equation}
при всех $f\in\mathbb{C}[G]_q$, $\xi\in U_q\mathfrak{g}$.
\end{lemma}

{\bf Доказательство.} Единственность отображения $*$ очевидна. Докажем его
существование с помощью инволюции $\sharp$ в алгебре Хопфа
$\mathbb{C}[G_0]_q$ (см. \itemiii\ \ref{canon_example}). Пусть, как и
прежде, $L_\mathrm{reg}$ -- это представление алгебры
$U_q\mathfrak{g}^\mathrm{op}$ в пространстве $\mathbb{C}[G]_q$, определяемое
равенством
$$
(L_\mathrm{reg}(\xi)f)(\eta)=f(\xi\eta),\qquad
f\in\mathbb{C}[G]_q,\;\xi,\eta\in U_q\mathfrak{g}.
$$
Это представление допускает продолжение по непрерывности до представления
$\overline{L}_\mathrm{reg}$ алгебры $\mathbb{C}[G]_q^{*\mathrm{op}}$. Из
сравнения \eqref{dual_star} с \eqref{g_0} следует, что при всех
$f\in\mathbb{C}[G]_q$, $\xi\in U_q\mathfrak{g}$
$$
\left\langle\overline{L}_\mathrm{reg}(\widetilde{w}_0)f^*,\xi\right\rangle=
\left\langle\overline{L}_\mathrm{reg}(\widetilde{w}_0)f,
(S(\xi))^*\right\rangle=
\left\langle\left(\overline{L}_\mathrm{reg}(\widetilde{w}_0)f\right)^\sharp,
\xi\right\rangle.
$$
Значит,
\begin{equation}\label{intertw_stars}
\overline{L}_\mathrm{reg}(\widetilde{w}_0)\cdot
*=\sharp\cdot\overline{L}_\mathrm{reg}(\widetilde{w}_0).
\end{equation}
Наконец,
\begin{equation}\label{explicit_star_reg}
f^*=L_\mathrm{reg}(\widetilde{w}_0^{-1})
\left(L_\mathrm{reg}(\widetilde{w}_0)f\right)^\sharp.
\end{equation}
Последнее равенство доставляет антилинейный оператор с нужными свойствами.
\hfill $\square$

\bigskip Докажем, что $*$ является инволюцией алгебры $\mathbb{C}[G]_q$.
\begin{lemma}\label{prop_6.1}
1. $f^{**}=f$ для всех $f\in\mathbb{C}[G]_q$.

2. $(f_1f_2)^*=f_2^*f_1^*$ для всех $f_1,f_2\in\mathbb{C}[G]_q$.
\end{lemma}

{\bf Доказательство.} Квадрат антилинейного оператора $*S:U_q\mathfrak{g}\to
U_q\mathfrak{g}$ является тождественным отображением. Отсюда следует первое
из доказываемых утверждений:
$$
\langle f^{**},\widetilde{w}_0\xi\rangle=\langle
f,\widetilde{w}_0(S((S(\xi))^*))^*\rangle=\langle f,\xi\rangle
$$
при всех $f\in\mathbb{C}[G]_q$, $\xi\in U_q\mathfrak{g}$.

Докажем второе утверждение. Строго говоря, в приведенном ниже
доказательстве следовало бы заменить универсальную R-матрицу определяемым
ею оператором $R_{\mathbb{C}[G]_q\mathbb{C}[G]_q}$. Но это привело бы к
потере наглядности. Воспользуемся свойствами отображения $*S$, вытекающими
из равенств
\begin{gather*}
\triangle(S(\xi))=(S\otimes
S)\triangle^\mathrm{cop}(\xi),\qquad\triangle(\xi^*)=
(\triangle(\xi))^{*\otimes *},
\\ S\otimes S(R)=R,\qquad R^{*\otimes *}=R_{21},
\end{gather*}
где $R_{21}$ отличается от $R$ перестановкой тензорных сомножителей
(последнее равенство следует из хорошо известных свойств универсальной
R-матрицы; см., например, предложение 2.3.1 в \cite{KorSoib}).

Для любых $f_1,f_2\in\mathbb{C}[G]_q$, $\xi\in U_q\mathfrak{g}$
$$
\langle(f_1f_2)^*,\widetilde{w}_0\xi\rangle=\overline{\langle
f_1f_2,\widetilde{w}_0(S(\xi))^*\rangle}=\overline{\langle f_1\otimes
f_2,(\widetilde{w}_0\otimes\widetilde{w}_0)R\triangle(S(\xi))^*\rangle};
$$
\begin{multline*}
\langle f_2^*f_1^*,\widetilde{w}_0\xi\rangle=\langle f_2^*\otimes
f_1^*,(\widetilde{w}_0\otimes\widetilde{w}_0)R\triangle(\xi)\rangle=
\\ =\overline{\langle f_2\otimes
f_1,(\widetilde{w}_0\otimes\widetilde{w}_0)(*S\otimes
*S)(R\triangle(\xi))\rangle}=
\\ =\overline{\langle f_2\otimes
f_1,(\widetilde{w}_0\otimes\widetilde{w}_0))
R_{21}\triangle^\mathrm{cop}((S(\xi))^*)\rangle}=
\\ =\overline{\langle f_1\otimes
f_2,(\widetilde{w}_0\otimes\widetilde{w}_0))R\triangle((S(\xi))^*)\rangle}.
\end{multline*}
Значит, $\langle(f_1f_2)^*,\widetilde{w}_0\xi\rangle=\langle
f_2^*f_1^*,\widetilde{w}_0\xi\rangle$ при всех $f_1,f_2\in\mathbb{C}[G]_q$,
$\xi\in U_q\mathfrak{g}$. \hfill $\square$

\bigskip

Рассмотрим $*$-алгебру
$\mathbb{C}[w_0G_0]_q\overset{\mathrm{def}}{=}(\mathbb{C}[G]_q,*)$. Она
является $q$-аналогом алгебры $\mathbb{C}[w_0G_0]$ регулярных функций на
вещественном аффинном алгебраическом многообразии $w_0G_0$, где $G_0$ --
некомпактная вещественная форма группы $G$.

\begin{example}
Если $\mathfrak{g}=\mathfrak{sl}_2$, то $(\alpha,\alpha)=2$,
$(\overline{\omega},\overline{\omega})=\frac12$, $(H,H)=2$. С помощью
\eqref{s_bar_eq} находим явный вид элемента $\widetilde{w}_0$ в
фундаментальном представлении $\overline{\pi}$ алгебры
$\mathbb{C}[SL_2]_q^*$:
$$
\overline{\pi}(\widetilde{w}_0)=q^{-1/4}\begin{pmatrix}0 & -1\\ q^{-1} &
0\end{pmatrix}.
$$
Значит,
\begin{equation}\label{L_action}
\begin{pmatrix}\overline{L}_\mathrm{reg}(\widetilde{w}_0)t_{11} &\;
\overline{L}_\mathrm{reg}(\widetilde{w}_0)t_{12}\\
\overline{L}_\mathrm{reg}(\widetilde{w}_0)t_{21}&\;
\overline{L}_\mathrm{reg}(\widetilde{w}_0)t_{22}\end{pmatrix}=
q^{-1/4}\begin{pmatrix}0 &\; -1\\ q^{-1} &\;
0\end{pmatrix}\begin{pmatrix}t_{11} &\; t_{12}\\ t_{21} &\;
t_{22}\end{pmatrix}.
\end{equation}

Используя \eqref{L_action}, нетрудно убедиться в том, что инволюции
\eqref{inv}, \eqref{inv_su11} связаны равенством \eqref{intertw_stars}.
Например,
 $$
 \overline{L}_\mathrm{reg}(\widetilde{w}_0)t_{11}^*=
\overline{L}_\mathrm{reg}(\widetilde{w}_0)(-t_{22})=-q^{-5/4}t_{12}
=(-t_{21})^\sharp=
\left(\overline{L}_\mathrm{reg}(\widetilde{w}_0)t_{11}\right)^\sharp.
$$
Значит, $\overline{L}_\mathrm{reg}(\widetilde{w}_0)f^*=
\left(\overline{L}_\mathrm{reg}(\widetilde{w}_0)f\right)^\sharp$ при
$f=t_{11}$.

Таким образом, новое определение инволюции $*$ обобщает прежнее ее
определение \eqref{inv}, применимое только в частном случае
$\mathfrak{g}=\mathfrak{sl}_2$.
\end{example}

\begin{lemma}\label{mod_alg}
$*$-Алгебра $\mathbb{C}[w_0G_0]_q$ является $(U_q\mathfrak{g},*)$-модульной.
\end{lemma}

{\bf Доказательство.} При всех $f\in\mathbb{C}[G]_q$, $\xi,\eta\in
U_q\mathfrak{g}$
$$
\langle(\xi f)^*,\widetilde{w}_0\eta\rangle=\overline{\langle\xi
f,\widetilde{w}_0(S(\eta))^*\rangle}=\overline{\langle
f,\widetilde{w}_0(S(\eta))^*\xi\rangle};
$$
\begin{multline*}
\langle (S(\xi))^*f^*,\widetilde{w}_0\eta\rangle=\langle
f^*,\widetilde{w}_0\eta(S(\xi))^*\rangle=\overline{\langle
f,\widetilde{w}_0(S(\eta(S(\xi))^*))^*\rangle}=
\\ =\overline{\langle
f,\widetilde{w}_0(S(\eta))^*(S((S(\xi))^*))^*\rangle}=\overline{\langle
f,\widetilde{w}_0(S(\eta))^*\xi\rangle}.
\end{multline*}
Значит, при всех $f\in\mathbb{C}[G]_q$, $\xi\in U_q\mathfrak{g}$
$$(\xi f)^*=(S(\xi))^*f^*.\eqno\square$$

\bigskip
Перейдем к доказательству того, что $t^*=t^\star$.

\begin{lemma}\label{up_to_const}
\begin{equation}\label{t_star}
t^*=\mathrm{const}\cdot c_{-\overline{\omega}_{l_1},\overline{\omega}_{l_0}}
^{\overline{\omega}_{l_0}}.
\end{equation}
\end{lemma}

{\bf Доказательство.} Рассмотрим фундаментальное представление $\pi_{l_0}$
алгебры $U_q\mathfrak{g}$, отвечающее $U_q\mathfrak{g}$-модулю
$L(\overline{\omega}_{l_0})$, и двойственное представление $\pi_{l_0}^*$.
Элементы $ t^*$, $c_{-\overline{\omega}_{l_1},\overline{\omega}_{l_0}}
^{\overline{\omega}_{l_0}}$ принадлежат одному и тому же весовому
подпространству простого \hbox{$U_q\mathfrak{g}^\mathrm{op}\otimes
U_q\mathfrak{g}$}-модуля
$\left\{f\in(U_q\mathfrak{g})^*\left|\:\operatorname{Ker}f\supset
\operatorname{Ker}\pi_{l_0}^*\right.\right\}\subset\mathbb{C}[G]_q$. Но это
подпространство одномерно. \hfill $\square$

\medskip

\begin{lemma}\label{star_aster}
$t^*=t^\star$.
\end{lemma}

{\bf Доказательство.} Согласно \eqref{aster_standard}, достаточно доказать,
что константа в \eqref{t_star} равна
$(-q)^{(\overline{\omega}_{l_1}+\overline{\omega}_{l_0},\rho)}. $ Но
$\langle f^*,\widetilde{w}_0\rangle=\overline{\langle
f,\widetilde{w}_0\rangle}$ при всех \hbox{$f\in\mathbb{C}[G]_q$}. Значит,
достаточно доказать равенство
\begin{equation}\label{star_aster_1}
\overline{\left\langle c_{\overline{\omega}_{l_1},-\overline{\omega}_{l_0}}
^{\overline{\omega}_{l_1}},\widetilde{w}_0\right\rangle}=
(-q)^{(\overline{\omega}_{l_1}-w_0\overline{\omega}_{l_1},\rho)}\left\langle
c_{-\overline{\omega}_{l_1},\overline{\omega}_{l_0}}
^{\overline{\omega}_{l_0}},\widetilde{w}_0\right\rangle.
\end{equation}
Но, как следует из \eqref{w0_tilda}, равенство \eqref{star_aster_1}
равносильно равенству
\begin{equation}\label{star_aster_2}
\overline{\left\langle c_{\overline{\omega}_{l_1},-\overline{\omega}_{l_0}}
^{\overline{\omega}_{l_1}},\overline{w}_0^{-1}\right\rangle}=
(-q)^{(\overline{\omega}_{l_1}- u_1\,\overline{\omega}_{l_1},\rho)}
\left\langle c_{-\overline{\omega}_{l_1},\overline{\omega}_{l_0}}
^{\overline{\omega}_{l_0}},\overline{w}_0^{-1}\right\rangle,
\end{equation}
где $u_1$ -- элемент наименьшей длины в классе смежности
$w_0W_{\mathbb{S}_1}$ по подгруппе $W_{\mathbb{S}_1}\subset W$, порожденной
множеством $\{s_j\;|\; j\neq l_1\}$.

Если $\overline{\omega}_{l_1}- u_1\overline{\omega}_{l_1}=
\sum\limits_kn_k\alpha_k$, то
$(\overline{\omega}_{l_1}-u_1\overline{\omega}_{l_1},\rho)=
\sum\limits_kn_kd_k$, поскольку
$(\alpha_k,\overline{\omega}_j)=d_j\delta_{jk}$. Значит, если
$u_1=s_{i_1}s_{i_2}\ldots s_{i_M}$ -- приведенное разложение элемента $u_1$,
то
$$
(-q)^{(\overline{\omega}_{l_1}-u_1 \overline{\omega}_{l_1},\rho)}=
\prod_{j=1}^M(-q_{i_j}).
$$
Из последнего равенства и из определения $\overline{w}_0^{\,-1}$ следует,
что \eqref{star_aster_2} достаточно доказать в простейшем частном случае
$\mathfrak{g}=\mathfrak{sl}_2$. Остается воспользоваться равенством
\eqref{s_bar}. \hfill $\square$

\medskip
Предложение \ref{star_X_new} вытекает из доказанных лемм. \hfill $\square$

\bigskip
В дальнейшем будет использоваться встретившееся в доказательстве леммы
\ref{star_aster} обозначение $\mathbb{S}_1=\{1,2, \cdots, l\} \setminus
\{l_1\}$.

\begin{remark}\label{explicit_*} Для матричных элементов в  специальном
 базисе
векторного пространства $L(\overline{\omega}_{l_1})$ имеет место равенство
$$
\left(c_{\overline{\omega}_{l_1},\lambda}^{\overline{\omega}_{l_1}}\right)^*=
(-1)^{\lambda(H_{\mathbb{S}_1})+\overline{\omega}_{l_0}(H_{\mathbb{S}_1})}
\left(c_{\overline{\omega}_{l_1},\lambda}^{\overline{\omega}_{l_1}}\right)
^\star,\qquad\lambda\in W\overline{\omega}_{l_0}.
$$
В частности, полагая $l=m+n-1$, $l_0=n$,
$\mathfrak{g}=\mathfrak{sl}_{m+n}$, находим явный вид инволюции $*$
\cite[стр. 859]{SSV4}:
\begin{equation}\label{SL-*}
t_{ij}^*=
\mathrm{sign}\left(\left(i-m-\frac12\right)\left(j-n-\frac12\right)\right)
t_{ij}^\star.
\end{equation}
\end{remark}


\subsubsection{Каноническое вложение $*$-алгебры \boldmath
$\operatorname{Pol}(\mathfrak{p}^-)_q$.}\label{canon_embed}

Рассмотрим подалгебру
$\mathbb{C}[\mathbb{X}_\mathbb{S}]_{q,x}\subset\mathbb{C}[G]_{q,x}$,
порожденную $\mathbb{C}[\mathbb{X}_\mathbb{S}]_q$ и элементом $x^{-1}$. Она
является локализацией $\mathbb{C}[\mathbb{X}_\mathbb{S}]_q$ по
мультипликативному множеству $x^{\mathbb{Z}_+}$.

Инволюция $*$ единственным образом продолжается с
$\mathbb{C}[\mathbb{X}_\mathbb{S}]_q$ на
$\mathbb{C}[\mathbb{X}_\mathbb{S}]_{q,x}$. Именно, $(x^{-1})^*=x^{-1}$.
Возникающая при этом $U_q\mathfrak{g}$-модульная алгебра с инволюцией
 $
 (\mathbb{C}[\mathbb{X}_\mathbb{S}]_{q,x},*)$
является $(U_q\mathfrak{g},*)$-модульной, см. предложение
\ref{*_local_gen}.

 Аналогично определяется $(U_q\mathfrak{g},*)$-модульная
алгебра
$$\mathbb{C}[w_0G_0]_{q,x}\overset{\mathrm{def}}{=}(\mathbb{C}[G]_{q,x},*).$$
\begin{theorem}\label{c_e}\begin{itemize}
\item[1.] Существует и единствен такой гомоморфизм
$(U_q\mathfrak{g},*)$-модульных алгебр
$\mathcal{I}:\operatorname{Pol}(\mathfrak{p}^-)_q
 \to \mathbb{C}[\mathbb{X}_\mathbb{S}]_{q,x}$, при котором
$\mathcal{I}: \,z_{\mathrm{low}}\mapsto t^{-1}t'$.
 \item[2.] Гомоморфизм $\mathcal{I}$ инъективен.
\end{itemize}
\end{theorem}

{\bf Доказательство.} Пусть $\mathbb{C}[\mathbb{X}_\mathbb{S}]_{q,t}$ --
подалгебра $\mathbb{C}[\mathbb{X}_\mathbb{S}]_{q,x}$, порожденная
$\mathbb{C}[\mathbb{X}_\mathbb{S}]_q$ и элементом $t^{-1}$.

Из предложения \ref{hol_embed} следует, что отображение
$i:z_\mathrm{low}\mapsto t^{-1}t'$ единственным образом продолжается до
вложения $U_q\mathfrak{g}$-модульных алгебр $
i:\mathbb{C}[\mathfrak{p}^-]_q\hookrightarrow\mathbb{C}[\mathbb{X}_\mathbb{S}]_{q,t}
$.

Кроме того, $\operatorname{Pol}(\mathfrak{p}^-)_q=$
\hbox{$(\mathbb{C}[\mathfrak{p}^-\oplus\mathfrak{p}^+],*)$}, и умножение
$$m:\mathbb{C}[\mathfrak{p}^-]_q\otimes\mathbb{C}[\mathfrak{p}^+]_q\to
\mathbb{C}[\mathfrak{p}^-\oplus\mathfrak{p}^+]_q, \qquad
 m:f_-\otimes f_+\mapsto f_-f_+,$$
 биективно.

 Значит, гомоморфизм $\mathcal{I}$ единствен,
и, если он существует, то
\begin{equation}\label{def_I}
\mathcal{I}(f_-f_+)=i(f_-)(i(f_+))^*,\qquad
f_\pm\in\mathbb{C}[\mathfrak{p}^\pm]_q.
\end{equation}

\medskip Перейдем к доказательству существования. Рассмотрим линейный
оператор $\mathcal{I}:\operatorname{Pol}(\mathfrak{p}^-)_q\to
 \mathbb{C}[\mathbb{X}_\mathbb{S}]_{q,x}$, определенный равенством \eqref{def_I}. Так как
$\operatorname{Pol}(\mathfrak{p}^-)_q$ и
 $\mathbb{C}[\mathbb{X}_\mathbb{S}]_{q,x}$ являются
  $(U_q\mathfrak{g},*)$-модульными алгебрами, то
$\mathcal{I}$ -- морфизм $U_q\mathfrak{g}$-модулей. Остается доказать, что
этот линейный оператор является гомоморфизмом алгебр.

Воспользуемся вложением
$\mathbb{C}[\mathbb{X}_\mathbb{S}]_{q,x}\subset\mathbb{C}[G]_{q,x}$. Как
следует из определения умножения в
$\mathbb{C}[\mathfrak{p}^-\oplus\mathfrak{p}^+]_q$, достаточно получить
равенство
\begin{equation}\label{eq_m}
(\mathcal{I}f_+)(\mathcal{I}f_-)=
m\check{R}_{\mathcal{I}\mathbb{C}[\mathfrak{p}^+]_q,
\mathcal{I}\mathbb{C}[\mathfrak{p}^-]_q}
(\mathcal{I}f_+\otimes\mathcal{I}f_-),
\end{equation}
где $f_\mp\in\mathcal{I}\mathbb{C}[\mathfrak{p}^\pm]_{q,\pm 1}$, а $m$ --
умножение в $\mathbb{C}[G]_{q,x}$. Вместо \eqref{eq_m} будет получено более
общее равенство
\begin{equation}\label{eq_gen}
q^{ck^2}\mathcal{I}(f_+)t^{*k}\cdot t^k\cdot\mathcal{I}(f_-)=
m\check{R}_{\mathcal{I}\mathbb{C}[G]_{q,t^*},\mathcal{I}\mathbb{C}[G]_{q,t}}
(\mathcal{I}(f_+)t^{*k}\otimes t^k\mathcal{I}(f_-)),
\end{equation}
где $k\in\mathbb{Z}$,\
$\mathbb{C}[G]_{q,t^*}=\{f^*|\:f\in\mathbb{C}[G]_{q,t}\}$, а $c$ --
рациональное число, зависящее от $f_-,f_+$. Так же, как в \itemiiiе
\ref{w0_G_x}, доказывается, что последнее равенство достаточно получить для
больших $k$. Точнее, его достаточно доказать в предположении
\begin{equation}\label{special_case}
\mathcal{I}(f_+)t^{*k}\in\mathbb{C}[X_\mathbb{S}^+]_q,\qquad
\mathcal{I}(f_-)t^k\in\mathbb{C}[X^-_\mathbb{S}]_q,
\end{equation}
где $\mathbb{C}[X^-_\mathbb{S}]_q$, $\mathbb{C}[X_\mathbb{S}^+]_q\subset
\mathbb{C}[\mathbb{X}_\mathbb{S}]_q $ -- это $U_q\mathfrak{g}$-модульные
подалгебры, порожденные элементами $t$, $t^*$ соответственно, см. \itemiii\
\ref{w0_G_x}.

 Отметим, что
 $\mathbb{C}[X^+_\mathbb{S}]_q=(\mathbb{C}[X_\mathbb{S}^-]_q)^*$.

 Наделим $\mathbb{C}[\mathbb{X}_\mathbb{S}]_q$ градуировкой, полагая
$$
\deg\left(c_{\overline{\omega}_{l_1};\lambda,j}^{\overline{\omega}_{l_1}}
\right)=1,\qquad
\deg\left(c_{-\overline{\omega}_{l_1};\mu,k}^{\overline{\omega}_{l_0}}\right)
=-1.
$$
 Корректность этого определения доказывается так же, как аналогичное
утверждение в \itemiiiе \ref{w0_G_x}. В предположении \eqref{special_case}
равенство \eqref{eq_gen} следует из равенства
\begin{equation}\label{c_r_new}
q^{-\deg(\psi_+)\deg(\psi_-)/(H_\mathbb{S},H_\mathbb{S})}\psi_+\psi_-=
m\check{R}_{\mathbb{C}[G]_q,\mathbb{C}[G]_q} (\psi_+\otimes\psi_-),
\end{equation}
где $\psi_+\in \mathbb{C}[X_\mathbb{S}^+]_q]_q$, $\psi_-\in
\mathbb{C}[X^-_\mathbb{S}]_q$. В свою очередь, \eqref{c_r_new} следует из
коммутационных соотношений \eqref{QISM}, как нетрудно показать, используя
\eqref{t0_new}, равенства
$$
L_\mathrm{reg}(F_j)\psi_+=L_\mathrm{reg}(E_j)\psi_-=0,\qquad j=1,2,\ldots,l,
$$
$$ L_\mathrm{reg}(H_j)\psi_\pm=0,\quad j\ne l_1, \qquad
(H_{\mathbb{S}},H_{\mathbb{S}})=(H_{\mathbb{S}_1},H_{\mathbb{S}_1})
$$
 и определение градуировки в $\mathbb{C}[\mathbb{X}_\mathbb{S}]_q$.

Существование и единственность гомоморфизма $\mathcal{I}$ доказаны. Осталось
убедиться в том, что он инъективен. Без предположения о трансцендентности
$q$ доказать инъективность гомоморфизма $\mathcal{I}$ непросто. Для этого
потребуются вспомогательные результаты.

\bigskip
Наш план состоит в построении такого представления $\mathscr{T}$ алгебры
$\mathbb{C}[G]_{q,x}$, что отвечающее ему представление
$\mathscr{T}\circ\mathcal{I}$ алгебры
$\mathbb{C}[\mathfrak{p}^-\oplus\mathfrak{p}^+]_q$ является точным. Отсюда
будет следовать инъективность $\mathcal{I}$.

При построении $\mathscr{T}$ мы будем следовать описанному в \itemiiiе
\ref{IrrepK} алгоритму построения $*$-представлений алгебры
$\mathbb{C}[K]_q=(\mathbb{C}[G]_q,\star)$. Именно, пусть $w_0^\mathbb{S}$
-- элемент наименьшей длины в классе смежности $W_\mathbb{S}\,w_0\subset W$
по подгруппе $W_\mathbb{S}$, порожденной простыми отражениями $s_i$, $i\ne
l_0$, см. лемму \ref{w_s}. Выберем приведенное разложение этого элемента
$$w_0^\mathbb{S}=s_{i_1}s_{i_2}\ldots s_{i_N}$$
и рассмотрим неприводимое представление $*$-алгебры Хопфа $\mathbb{C}[K]_q$
в $l^2(\mathbb{Z}_+)^{\otimes N}$
$$\Pi=\Pi_{i_1}\otimes\Pi_{i_2}\otimes\cdots\otimes\Pi_{i_N},$$
где $\Pi_{i_j}$ -- $*$-представления алгебры $(\mathbb{C}[G]_q,\star)$ в
$l^2(\mathbb{Z}_+)$, описанные в \itemiiiе \ref{IrrepK}. Из определения
представлений $\Pi_{i_j}$ и представления $\Pi$ следует, что операторы
$\Pi(f)$, $f\in\mathbb{C}[G]_q$, допускают сужение на линейную оболочку
$\mathscr{L}$ элементов стандартного базиса $\{e_{k_1}\otimes
e_{k_2}\otimes\cdots\otimes e_{k_N}\}$ гильбертова пространства
$l^2(\mathbb{Z}_+)^{\otimes N}$. Как известно (\cite[стр. 314,
315]{DijStok1}), элементы этого базиса являются собственными векторами
ограниченного самосопряженного компактного линейного оператора $\Pi(x)$, и
отвечающие им собственные значения положительны. Кроме того, старшее
собственное значение однократно и отвечает собственному вектору
$\mathbf{e}=e_0\otimes e_0\otimes\cdots\otimes e_0$. В частности, линейный
оператор $\Pi(x)|_\mathscr{L}$ в предгильбертовом пространстве $\mathscr{L}$
обратим. Значит, представление $\Pi|_\mathscr{L}$ алгебры $\mathbb{C}[G]_q$
в этом предгильбертовом пространстве канонически продолжается до
представления $\mathscr{T}$ ее локализации $\mathbb{C}[G]_{q,x}$
$$
\mathscr{T}(f)v=\Pi(f)v,\qquad f\in\mathbb{C}[G]_q,\quad v\in\mathscr{L}.
$$
Остается доказать точность представления
$\mathscr{T}_\mathcal{I}=\mathscr{T}\cdot\mathcal{I}$ алгебры
$\mathbb{C}[\mathfrak{p}^-\oplus\mathfrak{p}^+]_q$. Прежде всего, докажем
существование такой инволюции $\star$ в
$\mathbb{C}[\mathfrak{p}^-\oplus\mathfrak{p}^+]_q$, что
\begin{equation}\label{new_*}
\mathcal{I}(f^\star)=\mathcal{I}(f)^\star,\qquad
f\in\mathbb{C}[\mathfrak{p}^-\oplus\mathfrak{p}^+]_q.
\end{equation}
Для этого воспользуемся стандартной градуировкой
$$
\mathbb{C}[\mathfrak{p}^-\oplus\mathfrak{p}^+]_{q,j}=
\left\{\left.f\in\mathbb{C}[\mathfrak{p}^-\oplus\mathfrak{p}^+]_q\right|\:
H_\mathbb{S}f=2jf\right\},\qquad j\in\mathbb{Z},
$$
и связанным с ней автоморфизмом алгебры
$\mathbb{C}[\mathfrak{p}^-\oplus\mathfrak{p}^+]_q$: $f\mapsto(-1)^{\deg
f}f$. Именно, положим
\begin{equation}\label{starnew_starold}
f^\star\stackrel{\operatorname{def}}{=}(-1)^{\deg f}f^*
\end{equation}
и введем обозначение $\operatorname{Pol}(U)_q=
(\mathbb{C}[\mathfrak{p}^-\oplus\mathfrak{p}^+]_q,\star)$. Отметим, что
алгебра $\operatorname{Pol}(U)_q$ является
$(U_q\mathfrak{g},\star)$-модульной, поскольку алгебра
$\operatorname{Pol}(\mathfrak{p}^-)_q$ является
$(U_q\mathfrak{g},*)$-модульной.

Равенство \eqref{new_*} вытекает из того, что
$\mathcal{I}(f^*)=\mathcal{I}(f)^*$ при всех
$f\in\operatorname{Pol}(\mathfrak{p}^-)_q$ и из определений инволюций.
Следовательно, $\mathscr{T}_\mathcal{I}$ является \hbox{$*$-представлением}
алгебры $\operatorname{Pol}(U)_q$ в предгильбертовом пространстве
$\mathscr{L}$. Так как собственное значение, отвечающее собственному вектору
$\mathbf{e}$, является наибольшим, то из коммутационных соотношений
\eqref{quasicom} получаем:
\begin{equation}\label{main_prop}
\mathscr{T}
\left(c_{\overline{\omega}_{l_1},\lambda}^{\overline{\omega}_{l_1}}\right)
\mathbf{e}=0,\quad\lambda\ne
w_0\overline{\omega}_{l_1};\qquad\mathscr{T}_\mathcal{I}(z)\mathbf{e}=0,
\quad z\in\mathbb{C}[\mathfrak{p}^-]_{q,1}.
\end{equation}

Второе из этих равенств указывает план дальнейших действий. Будет введено в
рассмотрение фоковское представление $\mathscr{T}_F$ $*$-алгебры
$\operatorname{Pol}(U)_q$ и ненулевой морфизм представления $\mathscr{T}_F$
в $\mathscr{T}_\mathcal{I}$. Точность $\mathscr{T}_\mathcal{I}$ будет
следовать из точности и неприводимости
представления $\mathscr{T}_F$.\\

 Из хорошо известного равенства $(S\otimes\mathrm{id})R=(\mathrm{id}\otimes
S^{-1})R=R^{-1}$ (см. \cite[стр. 34]{Drinf2}) вытекает обратимость линейного
оператора
$\check{R}_{\mathbb{C}[\mathfrak{p}^+]_q,\mathbb{C}[\mathfrak{p}^-]_q}$ (см.
\itemiii\ \ref{RM}). Это
 позволяет перейти от разложения \eqref{expan1} к разложению
$$
\operatorname{Pol}(U)_q= \newoplus_{i,j=0}^\infty
\mathbb{C}[\mathfrak{p}^+]_{q,-j}\mathbb{C}[\mathfrak{p}^-]_{q,i}.
$$
Очевидно, $\mathbb{C}[\mathfrak{p}^+]_{q,0}=
\mathbb{C}[\mathfrak{p}^-]_{q,0}=\mathbb{C}\cdot 1$. Значит, для любого
элемента $f\in\operatorname{Pol}(U)_q$ существует и единственно разложение,
аналогичное \eqref{expan2}:
$$
f=\langle f\rangle\cdot 1+\sum_{(i,j)\ne(0,0)}f_{ij},\qquad
f_{ij}\in\mathbb{C}[\mathfrak{p}^+]_{q,-j}\mathbb{C}[\mathfrak{p}^-]_{q,i},
$$
где $\langle f\rangle$ -- линейный функционал на $\operatorname{Pol}(U)_q$.

Следуя \itemiiiу \ref{finite}, введем $\operatorname{Pol}(U)_q$-модуль
$\mathscr{H}$ с образующей $v_0$ и определяющими соотношениями $fv_0=0$ при
$f\in\mathbb{C}[\mathfrak{p}^-]_{q,1}$. То есть, $fv_0=\langle f\rangle v_0$
при $f\in\mathbb{C}[\mathfrak{p}^-]_q$. Очевидно,
$$
\mathbb{C}[\mathfrak{p}^+]_q\overset{\approx}{\to}\mathscr{H},\qquad
f\mapsto fv_0,
$$
что позволяет наделить $\mathscr{H}$ полуторалинейной формой
\begin{equation}\label{h_f}
(f_1v_0,f_2v_0)=\langle f_2^\star f_1\rangle,\qquad
f_1,f_2\in\mathbb{C}[\mathfrak{p}^+]_q.
\end{equation}
Как следует из равенства \eqref{main_prop},
$$
\langle f\rangle=(\mathscr{T}_\mathcal{I}(f)\mathbf{e},\mathbf{e}),\qquad
f\in\operatorname{Pol}(U)_q.
$$
Значит, \eqref{h_f} -- неотрицательная эрмитова форма. Ее положительность
нетрудно доказать, используя методы теории функций на компактных квантовых
группах. Действительно, если $\langle f^\star f\rangle=0$ и
$f\in\mathbb{C}[\mathfrak{p}^+]_q$, то при больших $k\in\mathbb{N}$ элемент
$\mathcal{I}(f)t^{\star k}$ принадлежит $*$-алгебре
$(\mathbb{C}[G]_q,\star)$ и аннулируется всеми неприводимыми ее
представлениями, см. предложение \ref{Soib_list}. Значит, $L^2$-норма
элемента $\mathcal{I}(f)t^{\star k}$ равна нулю. Из соотношений
ортогональности следует равенство $\mathcal{I}(f)t^{\star k}=0$, а из
предложения \ref{hol_embed} -- равенство $f=0$.

Таким образом, $\mathscr{H}$ является предгильбертовым пространством и
$\operatorname{Pol}(U)_q$-модулем. Ему отвечает $*$-представление алгебры
$\operatorname{Pol}(U)_q$ в $\mathscr{H}$. Это и есть упоминавшееся выше
фоковское представление $\mathscr{T}_F$.

Из определения $\mathscr{T}_F$ вытекает существование и единственность
такого линейного оператора $\mathscr{I}:\mathscr{H}\to\mathscr{L}$,
сплетающего представления $\mathscr{T}_F$ и $\mathscr{T}_\mathcal{I}$, что
$\mathscr{I}:v_0\mapsto\mathbf{e}$.

Неприводимость представления $\mathscr{T}_F$ доказывается так же, как
неприводимость $T_F$ в \itemiiiе \ref{finite}. Именно, если $\mathscr{H}'\ne
0$ -- общее инвариантное подпространство операторов $\mathscr{T}_F(f)$,
$f\in\operatorname{Pol}(U)_q$, то $\mathscr{H}'$ содержит ненулевой вектор
$v$, для которого
$$
\mathscr{T}_F(f)v=0,\qquad
f\in\newoplus_{j=1}^\infty\mathbb{C}[\mathfrak{p}^-]_{q,j}.
$$
Значит, вектор $v$ ортогонален подпространству
$\left(\newoplus\limits_{j=1}^\infty
\mathbb{C}[\mathfrak{p}^+]_{q,-j}\right)v_0$ и, следовательно,
$\mathbb{C}v=\mathbb{C}v_0$, $\mathscr{H}'=\mathscr{H}$. Неприводимость
$\mathscr{T}_F$ доказана.

Наделим $\mathscr{H}$ градуировкой
$$
\mathscr{H}=\newoplus\limits_{j=0}^\infty\mathscr{H}_j,\qquad
\mathscr{H}_j=\mathbb{C}[\mathfrak{p}^+]_{q,-j}v_0.
$$
Подпространства $\mathscr{H}_j$ попарно ортогональны и конечномерны. Так же,
как в \itemiiiе \ref{finite}, доказывается, что при всех
$j,k\in\mathbb{Z}_+$ отображение
$$
\mathbb{C}[\mathfrak{p}^+]_{q,-j}\otimes\mathbb{C}[\mathfrak{p}^-]_{q,k}
\overset{m}{\to}
\mathbb{C}[\mathfrak{p}^+]_{q,-j}\mathbb{C}[\mathfrak{p}^-]_{q,k}\to
\operatorname{Hom}(\mathscr{H}_k,\mathscr{H}_j),\qquad
f\mapsto\mathscr{T}_F(f)|_{\mathscr{H}_k},
$$
биективно (см. лемму \ref{ball_I_8.7}). Отсюда, используя положительность
эрмитовой формы $(\cdot,\cdot)$ в $\mathscr{H}$, нетрудно получить точность
представления $\mathscr{T}_F$ (см. лемму \ref{ball_I_8.8}).

Неприводимость и точность представления $\mathscr{T}_F$, а, следовательно,
точность $\mathscr{T}_\mathcal{I}$ и инъективность гомоморфизма
$\mathcal{I}$ доказаны. \hfill $\square$

\bigskip
  Следующее предложение вытекает из того, что алгебра $\mathbb{C}[G]_q$
  является областью целостности, см. предложение \ref{dom_noether}.
\begin{corollary}\label{Pol-integrity} Алгебра
$\operatorname{Pol}(\mathfrak{p}^-)_q$ является областью целостности.
\end{corollary}


\begin{remark} Многие утверждения \itemiiiов \ref{finite}, \ref{gen},
\ref{Uq-mod}, \ref{InvInt} без труда переносятся на случай компактных
квантовых групп. Ограничимся некоторыми определениями и формулировками
результатов. Они обычно используются при сравнении квантовых аналогов
компактного и некомпактного симметрических пространств, двойственных по
Картану. Рассмотрим алгебру с единицей
$\operatorname{Fun}(U)_q\supset\operatorname{Pol}(U)_q$, порожденную своей
подалгеброй $\operatorname{Pol}(U)_q$ и элементом $f_0$, удовлетворяющим
определяющим соотношениям
\begin{equation}\label{rel_1_2_new}
f_0^2=f_0, \quad \psi f_0=\langle\psi\rangle f_0,\;
\psi\in\mathbb{C}[\mathfrak{p}^-]_q,\quad f_0\psi=\langle\psi\rangle f_0, \;
\psi\in\mathbb{C}[\mathfrak{p}^+]_q.
\end{equation}
 Наделим ее инволюцией, полагая $f_0^\star=f_0$. Существует и единственно
такое продолжение $*$-представления $\mathscr{T}_F$ с подалгебры
$\operatorname{Pol}(U)_q$ на алгебру $\operatorname{Fun}(U)_q$, что
$\mathscr{T}_F(f_0)$ -- ортогональный проектор на одномерное подпространство
$\mathbb{C}v_0$. Полученное представление алгебры $\operatorname{Fun}(U)_q$
будем также называть фоковским. Двусторонний идеал $\mathscr{D}(U)_q$
алгебры $\operatorname{Fun}(U)_q$, порожденный элементом $f_0$, будем
называть алгеброй финитных функций. Обобщенными функциями будем называть
формальные ряды
$$
f=\sum_{k,j=0}^\infty f_{k,j},\qquad
f_{k,j}\in\mathbb{C}[\mathfrak{p}^+]_{q,-k}\mathbb{C}[\mathfrak{p}^-]_{q,j}.
$$
Векторное пространство $\mathscr{D}(U)_q'$ всех обобщенных функций будем
наделять топологией почленной сходимости (напомним, что
$\dim(\mathbb{C}[\mathfrak{p}^+]_{q,-k}\mathbb{C}[\mathfrak{p}^-]_{q,j})
<\infty$). Инволюция $\star$ и действие алгебры $U_q\mathfrak{g}$ допускают
продолжение по непрерывности на $\mathscr{D}(U)_q'$ и последующее сужение на
$\mathscr{D}(U)_q$. $*$-Алгебра $\mathscr{D}(U)_q$ является
$(U_q\mathfrak{g},\star)$-модульной. Наделим $\mathscr{H}$ структурой
$U_q\mathfrak{q}^-$-модуля и введем отвечающее этому
$U_q\mathfrak{q}^-$-модулю представление $\Gamma$ с помощью биекции
$$
\mathbb{C}[\mathfrak{p}^+]_qf_0\overset{\approx}{\to}\mathscr{H},\qquad\psi
f_0\mapsto\mathscr{T}_F(\psi)v_0.
$$
Положительный $U_q\mathfrak{g}$-инвариантный интеграл на $\mathscr{D}(U)_q$
существует, единствен с точностью до числового множителя и равен
\begin{equation}\label{inv_int_U}
\int\limits_{U_q}fd\nu=
\mathrm{const}\cdot\operatorname{tr}_q\mathscr{T}_F(f),
\end{equation}
где использованы конечномерность линейного оператора $\mathscr{T}_F(f)$,
$f\in \mathscr{D}(U)_q$, и обозначение $
\operatorname{tr}_q\mathscr{T}_F(f)=
\operatorname{tr}(\mathscr{T}_F(f)\Gamma(K_{-2\rho})), $ введенное в
\itemiiiе \ref{finite}.

 Разумеется, инвариантный интеграл на
$\mathscr{D}(U)_q$ можно получить по-другому, воспользовавшись естественным
вложением $*$-алгебры $\mathscr{D}(U)_q$ в алгебру непрерывных функций на
компактной квантовой группе и результатами \itemiiiа \ref{IrrepK}.
\end{remark}


\subsubsection{\boldmath $U_q\mathfrak{k}$-инвариантные
полиномы.}\label{inv_polynomials}

Введем в рассмотрение неоднократно используемые в дальнейшем
 попарно коммутирующие $U_q\mathfrak{k}$-инвариантные элементы
 $y_1,\ldots,y_r$ алгебры $\mathbb{C}[\mathfrak{p}^-]_q$.
Здесь
 $r$ -- ранг неприводимой ограниченной области $\mathbb{D}$, см. \itemiii\
\ref{Hua-Schmid}.

\begin{lemma}\label{spec_H_S} Множество
$$
 \{\,\overline{\omega}_{l_1}(H_\mathbb{S})\,,\,
\overline{\omega}_{l_1}(H_\mathbb{S})-2\,,\,\ldots\,,\,
-\overline{\omega}_{l_1}(H_\mathbb{S})+2\,,\,
-\overline{\omega}_{l_0}(H_\mathbb{S})\,\}
$$
 является множеством всех
собственных значений линейного оператора $H_{\mathbb{S}}$ в векторном
пространстве $L(\overline{\omega}_{l_1})$. Число собственных значений равно
$r+1$, где $r$ -- ранг ограниченной симметрической области $\mathbb{D}$.
\end{lemma}

{\bf Доказательство.} Старший и младший веса $U_q\mathfrak{g}$-модуля
$L(\overline{\omega}_{l_1})$ равны $\overline{\omega}_{l_1}$ и
$-\overline{\omega}_{l_0}$ соответственно. Следовательно, первое утверждение
вытекает из простоты $U_q\mathfrak{g}$-модуля $L(\overline{\omega}_{l_1})$ и
из определения элемента $H_\mathbb{S}\in\mathfrak{h}$. Второе утверждение
равносильно равенству
\begin{equation}\label{2_r}
\overline{\omega}_{l_1}(H_\mathbb{S})+\overline{\omega}_{l_0}(H_\mathbb{S})
=2r.
\end{equation}
Это равенство легко получить, используя картановский список неприводимых
ограниченных областей \cite{Helg}. Действительно, ранги областей этого
списка хорошо известны \cite[стр. 387]{Helg}, а слагаемые в левой части
равенства \eqref{2_r} равны коэффициентам при $\alpha_{l_0}$ в разложении
фундаментальных весов $\overline{\omega}_{l_1}$, $\overline{\omega}_{l_0}$
по простым корням алгебры Ли $\mathfrak{g}$. Коэффициенты таких разложений
приведены в \cite{Bou4-6}. \hfill $\square$

\medskip
Будем использовать введенные в \itemiiiе \ref{regular_K} обозначения
$$
v_{\lambda,j}^{\overline{\omega}_{l_1}},\qquad
c_{\overline{\omega}_{l_1};\lambda,j}^{\overline{\omega}_{l_1}}
$$
для элементов специального базиса векторного пространства
$L(\overline{\omega}_{l_1})$ и матричных элементов соответствующего
фундаментального представления в этом базисе.

\begin{lemma}\label{z_lambda_j}
Для любого матричного элемента
$c_{\overline{\omega}_{l_1};\lambda,j}^{\overline{\omega}_{l_1}}$
 существует и единствен такой
элемент $z_{\lambda,j}$ алгебры $\mathbb{C}[\mathfrak{p}^-]_q$, что
$$
\mathcal{I}(z_{\lambda,j})=t^{-1}\cdot
c_{\overline{\omega}_{l_1};\lambda,j}^{\overline{\omega}_{l_1}}.
$$
\end{lemma}

{\bf Доказательство.} Единственность следует из инъективности гомоморфизма
$\mathcal{I}$. Докажем существование элемента $z_{\lambda,j}$. Рассмотрим
подпространство $t \cdot
\mathcal{I}\mathbb{C}[\mathfrak{p}^-]_q\subset\mathbb{C}[G]_{q,x}$. Оно
содержит элемент $t$ и является $U_q\mathfrak{b}^+$-подмодулем
$U_q\mathfrak{b}^+$-модуля $\mathbb{C}[G]_{q,x}$, поскольку при всех
$j=1,2,\ldots,l$
\begin{gather*}
K_j^{\pm 1}(t\mathcal{I}(f))=(K_j^{\pm 1}t)\mathcal{I}(K_j^{\pm 1}f),
\\ E_j(t\mathcal{I}(f))=(E_jt)\mathcal{I}f+(K_jt)\mathcal{I}(E_jf),\qquad
f\in\mathbb{C}[\mathfrak{p}^-]_q,
\end{gather*}
и элементы $K_j^{\pm 1}t$, $E_jt$ принадлежат подпространству $t\cdot
\mathcal{I}\mathbb{C}[\mathfrak{p}^-]_q$. Значит,
$$
t \cdot \mathcal{I}\mathbb{C}[\mathfrak{p}^-]_q\supset
U_q\mathfrak{b}^+\cdot t,
$$
и остается воспользоваться тем, что все матричные элементы
$c_{\overline{\omega}_{l_1};\lambda,j}^{\overline{\omega}_{l_1}}$ содержатся
в подпространстве $U_q\mathfrak{b}^+\cdot t$. \hfill $\square$

\medskip Будем использовать сокращенное обозначение $z_\lambda$ вместо
$z_{\lambda,1}$, если весовое подпространство
$L(\overline{\omega}_{l_1})_\lambda$ одномерно.

\medskip
Доказанные леммы \ref{spec_H_S}, \ref{z_lambda_j} позволяют ввести в
рассмотрение элементы
\begin{equation}\label{def_y_i}
y_i=\sum_{\{(\lambda,j)|\:
\lambda(H_\mathbb{S})=-\overline{\omega}_{l_0}(H_\mathbb{S})+2i\}}
z_{\lambda,j}\cdot z_{\lambda,j}^*,\qquad i=1,\ldots,r.
\end{equation}
Очевидно, $y_i=y_i^*$ при всех $i=1,2,\ldots,r$. Введем обозначение
$y_0=1$.

\begin{proposition}\label{commut_y_i}
1. Элементы $y_1,y_2,\ldots,y_r$ алгебры
$\operatorname{Pol}(\mathfrak{p}^-)_q$ являются
$U_q\mathfrak{k}$-инвариантными и не зависят от выбора ортономированного
базиса $\left\{v_{\lambda,j}^{\overline{\omega}_{l_1}}\right\}$ в
$L(\overline{\omega}_{l_1})$.
\begin{flalign}\label{commutative}
\it{2}. && y_iy_j=y_jy_i,\qquad i,j=1,2,\ldots,r.&&
\end{flalign}
\end{proposition}

{\bf Доказательство.} Первое утверждение следует из того, что
 алгебра $\operatorname{Pol}(\mathfrak{p}^-)_q$ является
$(U_q\mathfrak{k},*)$-модульной, эрмитова форма в
$L(\overline{\omega}_{l_1})$ -- $U_q\mathfrak{k}$-инвариантной, а базис
$\left\{v_{\lambda,j}^{\overline{\omega}_{l_1}}\right\}$ --
ортономированным.

Второе утверждение вытекает из $U_q\mathfrak{k}$-инвариантности элементов
$y_i$ и из следующего результата, полученного Шкляровым.

\begin{lemma}\label{k-invariants}
Подалгебра $U_q\mathfrak{k}$-инвариантных элементов алгебры
$\operatorname{Pol}(\mathfrak{p}^-)_q$ коммутативна.
\end{lemma}

{\bf Доказательство.} Рассмотрим $U_q\mathfrak{k}$-подмодуль
$\mathbb{C}[\mathfrak{p}^-]_qf_0\subset D(\mathbb{D})_q$. Используя теорему
Хуа-Шмида \ref{k-types}, нетрудно доказать, что в категории
$U_q\mathfrak{k}$-модулей
$$
\mathbb{C}[\mathfrak{p}^-]_qf_0\approx
\newoplus_{\lambda\in\mathcal{A}_+}L(\mathfrak{k},\lambda),
$$
где $\mathcal{A}_+\subset P_+^\mathbb{S}$ -- введенное в \itemiiiе
\ref{Hua-Schmid} подмножество. Значит, каждый эндоморфизм
$U_q\mathfrak{k}$-модуля $\mathbb{C}[\mathfrak{p}^-]_q$ скалярен на каждой
из его $U_q\mathfrak{k}$-изотипических компонент. Следовательно, алгебра
эндоморфизмов $U_q\mathfrak{k}$-модуля $\mathbb{C}[\mathfrak{p}^-]_q$
коммутативна.

В частности, коммутируют операторы умножения на $\varphi$, $\psi$
$$f\mapsto \varphi f,\qquad f\mapsto\psi f,$$
 если элементы
$\varphi,\psi\in\operatorname{Pol}(\mathfrak{p}^-)_q$ являются
$U_q\mathfrak{k}$-инвариантными. Другими словами, коммутируют операторы
фоковского представления $T_F(\varphi)$, $T_F(\psi)$ алгебры
$\operatorname{Pol}(\mathfrak{p}^-)_q$ в предгильбертовом пространстве
$\mathcal{H}\cong\mathbb{C}[\mathfrak{p}^-]_qf_0$, см. \itemiii \
\ref{Uq-mod}. Остается воспользоваться леммой \ref{ball_I_8.8} о точности
фоковского представления. \hfill $\square$

\begin{remark}\label{Fourier}
Как видно из доказательства леммы \ref{k-invariants}, в категории
$U_q\mathfrak{k}$-модулей
$$
\mathcal{H}=\bigoplus_{\lambda\in\mathcal{A}_+}\mathcal{H}_\lambda,\qquad
\mathcal{H}_\lambda \approx L(\mathfrak{k},\lambda),
$$
и каждый $U_q\mathfrak{k}$-инвариантный элемент
$\psi\in\operatorname{Pol}(\mathfrak{p}^-)_q$ однозначно определяется
своими ''коэффициентами Фурье''
$\widehat{\psi}(\lambda)=T_F(\psi)|_{\mathcal{H}_\lambda}$.
\end{remark}


\subsubsection{Представления $*$-алгебры \boldmath
$\operatorname{Pol}(\mathfrak{p}^-)_q$.} \label{reps_pol}

В этом \itemiiiе будет завершено доказательство существования и
единственности точного неприводимого $*$-представления алгебры
$\operatorname{Pol}(\mathfrak{p}^-)_q$ ограниченными операторами в
гильбертовом пространстве.

Воспользуемся элементами $y_1,y_2,\ldots, y_r$, введенными в предыдущем
\itemiiiе.

\begin{proposition}\label{def_y}
 Существует и единствен такой элемент $y$ алгебры
$\operatorname{Pol}(\mathfrak{p}^-)_q$, что $\mathcal{I}y=x^{-1}$. Этот
элемент $U_q\mathfrak{k}$-инвариантен и равен
\begin{equation}\label{explicit_y} y=1+\sum_{i=1}^r(-1)^iy_i.
\end{equation}
\end{proposition}

{\bf Доказательство.} Единственность элемента $y$ очевидна. Покажем, что
его существование и равенство \eqref{explicit_y} вытекают из
\eqref{unitary_1}. Действительно, согласно \eqref{unitary_1},
\begin{equation}\label{unitary_new}
\sum_{(\lambda,j)}
c_{\overline{\omega}_{l_1};\lambda,j}^{\overline{\omega}_{l_1}}\cdot
\left(c_{\overline{\omega}_{l_1};\lambda,j}^{\overline{\omega}_{l_1}}\right)
^\star=1.
\end{equation}
Но
\begin{equation}\label{star_aster_new}
\left(c_{\overline{\omega}_{l_1};\lambda,j}^{\overline{\omega}_{l_1}}\right)
^\star=(-1)^{\frac{\lambda(H_\mathcal{S})+\lambda(\overline{\omega}_{l_0})}2}
\left(c_{\overline{\omega}_{l_1};\lambda,j}^{\overline{\omega}_{l_1}}\right)
^*,
\end{equation}
поскольку
$$
t^\star=t^*,\qquad F_j^\star=
\begin{cases}
F_j^*, & j\ne l_0,
\\-F_j^*, & j=l_0,
\end{cases}
$$
и
$$
(\xi f)^\star=(S(\xi))^\star f^\star,\qquad(\xi f)^*=(S(\xi))^*f^*
$$
при всех $\xi\in U_q\mathfrak{g}$, $f\in\mathbb{C}[G]_q$. Из определения
элементов $y_i$ и из равенств \eqref{unitary_new}, \eqref{star_aster_new}
следует, что
$$
t\mathcal{I}\left(\sum_{i=0}^r(-1)^iy_i\right)t^*=1.
$$
Остается умножить обе части последнего равенства на $t^{-1}$ слева и на
$(t^*)^{-1}$ справа.

Из \eqref{explicit_y} и из предложения \ref{commut_y_i} вытекает
$U_q\mathfrak{k}$-инвариантность $y$. \hfill $\square$

\medskip

\begin{corollary}\label{quasicom_y}
Для всех $z\in\mathbb{C}[\mathfrak{p}^-]_{q,1}$
\begin{equation}\label{q_y}
z\, y=q_{l_0}^{-2}\,y\, z,\qquad z^*\, y=q_{l_0}^2\, y\, z^*.
\end{equation}
\end{corollary}

{\bf Доказательство.} В частном случае $z=z_{\rm low}$ коммутационные
соотношения \eqref{quasicom_y} легко получить с помощью равенства
$\mathcal{I}y=x^{-1}$ и следствия \ref{subalgebra}. Общий случай сводится к
этому частному случаю действием $U_q\mathfrak{k}$. \hfill $\square$

\bigskip Рассмотрим фоковское представление $T_F$ алгебры
$\operatorname{Pol}(\mathfrak{p}^-)_q$ в предгильбертовом пространстве
$\mathcal{H}$ и пополнение $\overline{\mathcal{H}}$ этого предгильбертового
пространства. Как показано в
\itemiiiе \ref{boundness}, операторы представления $T_F$ ограничены и,
следовательно, допускают продолжение по непрерывности до операторов
представления $\overline{T}_F$ в гильбертовом пространстве
$\overline{\mathcal{H}}$.

\begin{proposition}\label{irreps_Pol}
1. Представление $\overline{T}_F$ является точным неприводимым
$*$-представлением алгебры $\operatorname{Pol}(\mathfrak{p}^-)_q$
ограниченными операторами в гильбертовом пространстве.

2. Представление, обладающее всеми этими свойствами, единственно с
точностью до унитарной эквивалентности.
\end{proposition}

{\bf Доказательство.} 1. Точность представления $\overline{T}_F$ вытекает из
доказанной в \itemiiiе \ref{finite} точности представления $T_F$. Перейдем к
доказательству неприводимости.

Из \eqref{def_y_i} следует, что самосопряженный линейный оператор
$\overline{T}_F(y)$ компактен и его спектр является замыканием
геометрической прогрессии. Точнее,
\begin{equation}\label{spec_y}
\overline{T}_F(y)|_{\mathcal{H}_k}=q_{l_0}^{2k}\cdot 1,\qquad
k\in\mathbb{Z}_+,
\end{equation}
где $\mathcal{H}_k=\mathbb{C}[\mathfrak{p}^-]_{q,k}\cdot v_0$. Пусть
$\mathfrak{A}$ -- замыкание алгебры операторов
$\left\{\left.\overline{T}_F(f)\right|\:
f\in\operatorname{Pol}(\mathfrak{p}^-)_q\right\}$ по операторной норме.
Алгебра $\mathfrak{A}$ содержит ортогональные проекторы на все
подпространства $\mathcal{H}_k$, как следует из \eqref{spec_y}. Эти
подпространства попарно ортогональны, конечномерны, и их сумма плотна в
$\overline{\mathcal{H}}$. Значит, из леммы \ref{ball_I_8.7} следует, что
алгебра $\mathfrak{A}$ содержит все компактные линейные операторы в
$\overline{\mathcal{H}}$. Это влечет неприводимость представления
$\overline{T}_F$.

Очевидно, $\overline{T}_F$ является $*$-представлением.

\medskip 2. Пусть $T'$ -- точное неприводимое $*$-представление алгебры
$\operatorname{Pol}(\mathfrak{p}^-)_q$ ограниченными линейными операторами
в гильбертовом пространстве. Покажем, что представления $T'$ и
$\overline{T}_F$ унитарно эквивалентны.

Из точности $T'$ следует, что $T'(y)\neq 0$. Это позволяет так же, как в
\itemiiiе \ref{Fock_l}, использовать
 коммутационные соотношения \eqref{q_y} для доказательства того, что
 спектр самосопряженного линейного оператора $T'(y)$
является замыканием геометрической прогрессии.

Пусть $v'$ -- нормированный собственный вектор, отвечающий наибольшему по
модулю собственному значению. Тогда, как следует из \eqref{q_y},
$T'(z)^*v'=0$, для всех $z\in\mathbb{C}[\mathfrak{p}^-]_{q,1}$. Значит,
$(T'(f)v',v')=(\overline{T}_F(f)v_0,v_0)$, для всех
$f\in\operatorname{Pol}(\mathfrak{p}^-)_q$. Таким образом, отображение
$v_0\mapsto v'$ допускает продолжение до изометрического линейного
оператора, сплетающего представления $\overline{T}_F$ и $T'$. Этот оператор
сюръективен, поскольку представление $T'$ неприводимо.\hfill $\square$

\begin{remark} Во второй части доказательства, то есть при доказательстве
единственности, используется не точность представления $T'$, а лишь то, что
$T'(y)\ne 0$. Это наблюдение находит неожиданное приложение: если $J=J^*$--
ненулевой двусторонний идеал алгебры $\mathrm{Pol}(\mathfrak{p}^-)_q$ и
если фактор-алгебра $\mathrm{Pol}(\mathfrak{p}^-)_q/J$ обладает
неприводимыми $*$-представлениями, разделяющими ее точки, то $y\in J$. В
частном случае квантового матричного шара и определяющего идеала $J$ его
границы Шилова \cite[стр. 381-383]{Vak01} так доказывается включение $(y)
\subset J$, означающее, что граница Шилова принадлежит топологической
границе. Здесь $(y)$-- двусторонний идеал, порожденный элементом $y$.
\end{remark}

\bigskip

  Мультипликативное подмножество
$y^{\mathbb{Z}_+}$ является подмножеством Орэ области целостности
$\operatorname{Pol}(\mathfrak{p}^-)_q$, см. \eqref{q_y} и следствие
\ref{Pol-integrity}. Рассмотрим локализацию
$\operatorname{Pol}(\mathfrak{p}^-)_{q,y}$ алгебры
$\operatorname{Pol}(\mathfrak{p}^-)_q$ по подмножеству $y^{\mathbb{Z}_+}$.
Разумеется, $\operatorname{Pol}(\mathfrak{p}^-)_q\hookrightarrow
\operatorname{Pol}(\mathfrak{p}^-)_{q,y}$. Более того, естественное
продолжение представления $T_F$ на локализацию
$\operatorname{Pol}(\mathfrak{p}^-)_{q,y}$ алгебры
$\operatorname{Pol}(\mathfrak{p}^-)_q$ является точным представлением
$\operatorname{Pol}(\mathfrak{p}^-)_{q,y}$ в $\mathcal{H}$. Рассуждая так
же, как при доказательстве теоремы \ref{c_e}, получаем

\begin{proposition}\label{inj}
Естественное продолжение канонического вложения
  $\operatorname{Pol}(\mathfrak{p}^-)_q\hookrightarrow
  \mathbb{C}[\mathbb{X}_\mathbb{S}]_{q,x}$
  на локализацию
$\operatorname{Pol}(\mathfrak{p}^-)_{q,y}$ является вложением алгебр
\begin{equation}\label{ext_inj}
\operatorname{Pol}(\mathfrak{p}^-)_{q,y}\hookrightarrow
 \mathbb{C}[\mathbb{X}_\mathbb{S}]_{q,x}.
\end{equation}
\end{proposition}

\begin{remark}
Так же, как в предыдущих \itemiiiах, доказывается, что структура
$(U_q\mathfrak{g},*)$-модульной алгебры канонически продолжаются с
$\operatorname{Pol}(\mathfrak{p}^-)_q$ на
$\operatorname{Pol}(\mathfrak{p}^-)_{q,y}$ и что \eqref{ext_inj} -- морфизм
$(U_q\mathfrak{g},*)$-модульных алгебр (используется то, что элементы
$\xi(fy^n)y^{-n}$, где $\xi\in U_q\mathfrak{g}$,
$f\in\operatorname{Pol}(\mathfrak{p}^-)_q$, являются значениями
соответствующих полиномов Лорана на геометрической прогрессии).
\end{remark}

\bigskip

Сохраним обозначение $\mathcal{I}$ для {\it канонического вложения}
\eqref{ext_inj} и опишем
 образ алгебры $\operatorname{Pol}(\mathfrak{p}^-)_{q,y}$ при этом
 вложении.

Пусть $K$- связная аффинная алгебраическая подгруппа с алгеброй Ли
$\mathfrak{k}$, см. \itemiii\ \ref{prehomogeneous}, $w_0$-- элемент
максимальной длины группы Вейля $W$ и
\begin{equation}\label{K_1}
K_1=w_0 K w_0^{-1}.
\end{equation}
Очевидно, для элементов $a,b \in G$ равенство $w_0\cdot Ka=w_0\cdot Kb$
равносильно равенству $K_1\,(w_0a)=K_1\,(w_0b)$.

 Решения системы уравнений
$$ L_{\rm reg}(E_i)f=L_{\rm reg}(F_i)f=0,\qquad i\neq l_1,$$
$$ L_{\rm reg}(K_j^{\pm 1})f=f,\qquad j=1,2,\ldots,l$$
образуют подалгебру $U_q\mathfrak{g}$-модульной алгебры $\mathbb{C}[G]_q$.
Это квантовый аналог алгебры $\mathbb{C}[K_1\backslash G]$ регулярных
функций на аффинном алгебраическом многообразии $K_1\backslash G$, см.
\itemiii\ \ref{IrrepK} и предложение \ref{from_Jo-R}. Введем обозначение
$\mathbb{C}[K_1\backslash G]_q$ для рассматриваемой
$U_q\mathfrak{g}$-модульной алгебры.

 Алгебра
$\mathbb{C}[K_1\backslash G]_q$ порождена элементами
 $c_{\overline{\omega}_{l_1};\lambda,j}^{\overline{\omega}_{l_1}}
 (c_{\overline{\omega}_{l_1};\mu,k} ^{\overline{\omega}_{l_1}})^\star$, см.
 \itemiii\  \ref{IrrepK}.
 Значит, ее можно определить как  наименьшую
 $U_q\mathfrak{g}$-модульную подалгебру, содержащую элемент $x$.
 Именно это определение  используется в дальнейшем.

 \medskip В \itemiiiе \ref{canon_embed} алгебра
 $\mathbb{C}[\mathbb{X}_\mathbb{S}]_{q}$
 была наделена градуировкой.  Продолжим ее на локализацию
 $\mathbb{C}[\mathbb{X}_\mathbb{S}]_{q,x}$.

\begin{lemma}\label{gradation}
Существует и единственна такая $\mathbb{Z}$-градуировка алгебры
 $\mathbb{C}[\mathbb{X}_\mathbb{S}]_{q,x}$,
 что
$ \deg\left(c_{\overline{\omega}_{l_1};\lambda,j}
^{\overline{\omega}_{l_1}}\right)=1$,\;
$\deg\left(c_{-\overline{\omega}_{l_1};\mu,k}
^{\overline{\omega}_{l_0}}\right)=-1.
$
\end{lemma}

{\bf Доказательство.} Единственность очевидна, а построение градуировки с
нужными свойствами осуществляется так же, как в \itemiiiе \ref{w0_G_x}.
Именно,
$$
L_\mathrm{reg}(H_{l_1})c_{\overline{\omega}_{l_1};\lambda,j}
^{\overline{\omega}_{l_1}}=c_{\overline{\omega}_{l_1};\lambda,j}
^{\overline{\omega}_{l_1}},\qquad
L_\mathrm{reg}(H_{l_1})c_{-\overline{\omega}_{l_1};\mu,k}
^{\overline{\omega}_{l_0}}=-c_{-\overline{\omega}_{l_1};\mu,k}
^{\overline{\omega}_{l_0}},
$$
и дифференцирование $L_\mathrm{reg}(H_{l_1})$ алгебры
 $\mathbb{C}[\mathbb{X}_\mathbb{S}]_{q}$
  допускает продолжение до дифференцирования
ее локализации $\mathbb{C}[\mathbb{X}_\mathbb{S}]_{q,x}$. Остается положить
$\deg f=j$ для всех таких элементов $f
\in\mathbb{C}[\mathbb{X}_\mathbb{S}]_{q,x}$, что
$L_\mathrm{reg}(H_{l_1})f=jf$.\hfill $\square$

\medskip Отметим, что локализация $\mathbb{C}[K_1\backslash G]_{q,x}$
алгебры
  $\mathbb{C}[K_1\backslash G]_q$ по
  мультипликативному подмножеству
  $x^{\mathbb{Z}_+}$ естественно вложена в
  $\mathbb{C}[\mathbb{X}_\mathbb{S}]_{q,x}$.

\begin{proposition}\label{image_emb}\ \ \
$\mathcal{I}\operatorname{Pol}(\mathfrak{p}^-)_{q,y}=
\mathbb{C}[K_1\backslash G]_{q,x}$.
\end{proposition}

{\bf Доказательство.} Включение $U_q\mathfrak{g}$-модульных алгебр
$\mathcal{I}\operatorname{Pol}(\mathfrak{p}^-)_{q,y}\subset
\mathbb{C}[K_1\backslash G]_{q,x}$ очевидно. Для доказательства обратного
включения рассмотрим подпространство
$$
\left(\newoplus_{j=1}^\infty\,t^j\cdot\mathcal{I}
\operatorname{Pol}(\mathfrak{p}^-)_{q,y}\right)\;\newoplus\;
\mathcal{I}\operatorname{Pol}(\mathfrak{p}^-)_{q,y}\;\newoplus\;
\left(\newoplus_{j=1}^\infty\,
\mathcal{I}\operatorname{Pol}(\mathfrak{p}^-)_{q,y}\cdot t^{*j}\right).
$$
Оно является $U_q\mathfrak{g}$-модульной подалгеброй и содержит элементы
$t$, $t^*$, $x^{-1}$. Следовательно, оно совпадает с
$\mathbb{C}[\mathbb{X}_\mathbb{S}]_{q,x}$. Значит,
$\mathcal{I}\operatorname{Pol}(\mathfrak{p}^-)_{q,y}$ содержит все элементы
нулевой степени однородности из $\mathbb{C}[\mathbb{X}_\mathbb{S}]_{q,x}$.
Остается воспользоваться тем, что ненулевые элементы алгебры
$\mathbb{C}[K_1\backslash G]_{q,x}$ являются однородными степени 0. \hfill
$\square$


\subsubsection{ Дополнение о сферических
$(\mathfrak{g},\mathfrak{k})$-модулях и о теореме
Хуа-Шмида.}\label{Hua-Schmid}

В этом \itemiiiе приведены хорошо известные и часто используемые результаты
теории эрмитовых симметрических пространств и ограниченных симметрических
областей \cite{Shimeno}, \cite{Takeuchi}.

Рассмотрим эрмитову симметрическую пару $(\mathfrak{g},\mathfrak{k})$,
отвечающую подмножеству $\mathbb{S}=\{1,2,\cdots,l\}\setminus\{l_0\}$, и
элемент $H_{\mathbb{S}}\in\mathfrak{h}$, см. \itemiiiы \ref{Hermit},
\ref{GVerma}.

Множество $\Phi_c=\{\beta\in\Phi|\:\beta(H_\mathbb{S})=0\}$ называют
множеством компактных корней, \IND{корни ! компактные} а множество
$\Phi_n=\Phi\setminus\Phi_c$ -- множеством некомпактных корней. \IND{корни
! некомпактные} Компактные корни образуют систему корней полупростой
комплексной алгебры Ли $\mathfrak{k}_{ss}\subset\mathfrak{k}$, порожденной
множеством $\{E_i,F_i,H_i\}_{i\ne l_0}$. Группа Вейля $W_\mathbb{S}$ этой
системы корней вложена в $W$ и действует на множестве положительных
некомпактных корней.

Два линейно независимых корня $\gamma'$, $\gamma''$ называют строго
ортогональными, \IND{корни ! строго ортогональные} если ни
$\gamma'+\gamma''$, ни $\gamma'-\gamma''$ не являются корнями. В этом
случае $\gamma'$ и $\gamma''$ ортогональны, как легко доказать, рассмотрев
множество корней $\{\gamma''+j\gamma'\in\Phi\,|\,j\in\mathbb{Z}\}$
\cite[стр. 185--186]{Bou4-6}.

Следуя \cite[стр. 572]{EnrJos}, введем в рассмотрение семейство простых
комплексных алгебр Ли
$\mathfrak{g}_1,\,\mathfrak{g}_2,\,\ldots\,,\mathfrak{g}_r$,
рассматриваемых с точностью до изоморфизма, с системами корней, вложенными
в $\Phi$. Их старшие корни $\{\gamma_1,\gamma_2,\ldots,\gamma_r\}$ являются
строго ортогональными положительными некомпактными корнями алгебры Ли
$\mathfrak{g}$ и играют важную роль в теории эрмитовых симметрических
пространств.

Выделим на диаграмме Дынкина алгебры Ли $\mathfrak{g}$ вершину, отвечающую
простому корню $\alpha_{l_0}$. Алгебры Ли
$\mathfrak{g}_1,\mathfrak{g}_2,\ldots,\mathfrak{g}_r$ будут определяться
своими системами простых корней, вложенными в систему простых корней
алгебры Ли $\mathfrak{g}$ и содержащими $\alpha_{l_0}$.

Опишем индуктивное построение этих алгебр Ли. Первой из них будет
$\mathfrak{g}_1=\mathfrak{g}$. Для перехода от $\mathfrak{g}_n$ к
$\mathfrak{g}_{n+1}$ рассмотрим подграф диаграммы Дынкина алгебры Ли
$\mathfrak{g}_n$, отвечающий ее корням, ортогональным ее старшему корню
$\gamma_n$. Если этот подграф не содержит выделенную вершину, то построение
завершено. В противном случае имеется компонента связности подграфа,
содержащая выделенную вершину. Она является диаграммой Дынкина и определяет
простую комплексную алгебру Ли $\mathfrak{g}_{n+1}\subset\mathfrak{g}_n$ с
системой корней, вложенной в систему корней алгебры Ли $\mathfrak{g}_n$. В
результате получаем линейно упорядоченное множество положительных
некомпактных строго ортогональных корней
$\{\gamma_1,\gamma_2,\ldots,\gamma_r\}$, которое использовалось в работах
\cite[стр. 82]{Inoue}, \cite[стр. 3776-3777]{OO}, \cite[стр.
443]{Takeuchi}.

Число $r$ называется рангом неприводимой ограниченной симметрической
области, связанной с эрмитовой симметрической парой $(\mathfrak{g},
\mathfrak{k})$.

\begin{example}
Алгебры Ли $\mathfrak{g}_1,\mathfrak{g}_2,\ldots,\mathfrak{g}_r$ легко
найти, используя при вычислении скалярных произведений явный вид
коэффициентов разложения максимального корня $\delta$ по фундаментальным
весам $\overline{\omega}_1,\overline{\omega}_2,\ldots\overline{\omega}_l$,
см. таблицы в \cite{Bou4-6}. Например, для исключительных алгебр Ли
$E_6,\;E_7$ получаем:
$$E_6,\quad A_5;\qquad\qquad E_7,\quad D_6,\quad A_1\;.$$
В первом случае ранг равен двум, а во втором -- трем.
\end{example}

\begin{remark}\label{HCh-gammas}
В \cite[стр. 334]{Shimeno} корни $\{\gamma_1,\gamma_2,\ldots,\gamma_r\}$
нумеруются в обратном порядке. Многие авторы
\cite{HCh6,Kor,RubShiff,John,Heckman_Sch} отдают предпочтение системе
строго ортогональных некомпактных положительных корней
$\{w_{0,\mathbb{S}}\gamma_1,\;w_{0,\mathbb{S}}\gamma_2,\;\ldots,\;
w_{0,\mathbb{S}}\gamma_r\}$, где $w_{0,\mathbb{S}}$ -- элемент максимальной
длины группы Вейля $W_\mathbb{S}$.
\end{remark}

\medskip

 Будем использовать выбранное в \itemiiiе \ref{PBW} невырожденное
 скалярное произведение  в $\mathfrak{h}^*$ для отождествления
 $\mathfrak{h}^*$ с $\mathfrak{h}$.
 Сопоставим каждому корню $\beta\in\Phi$ элемент
$H_\beta$ картановской подалгебры $\mathfrak{h}$, отвечающий линейному
функционалу $\beta \in \mathfrak{h}^*$ при изоморфизме $\mathfrak{h} \cong
\mathfrak{h}^*$. Как отмечалось в \itemiiiе \ref{U_ext},
$$H_{\alpha_i}=d_iH_i, \qquad \quad i=1,2,\ldots,l.$$

Следуя \cite{Shimeno}, введем в рассмотрение линейную оболочку
$\mathfrak{h}^-$ элементов $H_{\gamma_1},H_{\gamma_2},\ldots,H_{\gamma_r}$
и подпространство
$$\mathfrak{h}^+=\bigcap_{j=1}\operatorname{Ker}\gamma_j.$$
Нетрудно доказать, что $\mathfrak{h}^+$ является ортогональным дополнением к
$\mathfrak{h}^-$ в $\mathfrak{h}$ по отношению к инвариантному скалярному
произведению в $\mathfrak{g}$:
\begin{equation}\label{h_pm}
\mathfrak{h}=\mathfrak{h}^+\oplus\mathfrak{h}^-.
\end{equation}

\begin{example}
Пусть $N=m+n$ и $n\ge m$. Рассмотрим эрмитово симметрическую пару
$(\mathfrak{g},\mathfrak{k})$:
$$
\mathfrak{g}=\mathfrak{sl}_N,\qquad
\mathfrak{k}=\mathfrak{s}(\mathfrak{gl}_n\times\mathfrak{gl}_m).
$$
Положительные некомпактные корни имеют вид
$$
\varepsilon_i-\varepsilon_j,\qquad
i\in\{1,2,\ldots,n\},\;j\in\{n+1,n+2,\ldots,N\},
$$
где
\begin{equation}\label{epsilon}
\varepsilon_i(H_k)=
\begin{cases}
1, & i=k
\\ -1, & i=k+1
\\ 0, & i\notin\{k,k+1\}.
\end{cases}
\end{equation}
Группа $W$ действует в $\mathfrak{h}^*$ перестановками элементов
$\varepsilon_i$. Из определений следует, что $r=m$ и
\begin{equation}\label{roots_gamma}
\gamma_j=\varepsilon_j-\varepsilon_{N+1-j}=\sum_{i=j}^{N-j}\alpha_i,\qquad
j=1,2,\ldots,m,
\end{equation}
$$H_{\alpha_i}=H_i, \qquad H_{\gamma_j}=\sum_{i=j}^{N-j}H_i.$$
Подпространство $\mathfrak{h}^-$ порождается элементами
$$H_j+H_{N-j},\quad j=1,2,\ldots,m-1;\qquad \sum_{i=m}^{N-m}H_i.$$
\end{example}

\bigskip

Рассмотрим ограничение линейных функционалов на $\mathfrak{h}^-$:
\begin{equation}\label{res_roots}
\mathfrak{h}^*\mapsto(\mathfrak{h}^-)^*,\qquad f\mapsto
f|_{\mathfrak{h}^-}.
\end{equation}
Множество всех ненулевых элементов образа $\Phi$ при отображении
\eqref{res_roots} обозначим $\Phi_\mathrm{res}$, и кратностью элементов
этого множества будем называть число элементов его полного прообраза.
Известно, что подмножество $\Phi_\mathrm{res}\subset(\mathfrak{h}^-)^*$
является системой корней. Ее элементы называют ограниченными корнями,
\IND{корни ! ограниченные} а ее группу Вейля $W_\mathrm{res}$ --
ограниченной группой Вейля. \IND{группа Вейля ! ограниченная} Ненулевые
ограничения положительных корней на $\mathfrak{h}^-$ образуют множество
$\Phi_\mathrm{res}^+$, называемое множеством положительных ограниченных
корней.

\begin{remark} \label{Caley_remark}  Широко используемая в теории
представлений максимальная абелева подалгебры Ли
 $\mathfrak{a}_\mathbb{C}\subset \mathfrak{p}^+ \oplus \mathfrak{p}^-$
 может быть получена из  $\mathfrak{h}^-$
с помощью преобразования Кэли \cite[стр. 212]{Kor}, \cite[стр.
334]{Shimeno}.
\end{remark}

 Будем отождествлять  линейные функционалы
  $\{\gamma_1,\gamma_2,\ldots,\gamma_r\}$
и их сужения на $\mathfrak{h}^-$, что допустимо, поскольку эти сужения
линейно независимы \cite{HCh6}.

\begin{proposition}(\cite[стр. 212]{Kor})\label{Koranyi}
Группа $W_\mathrm{res}$ порождена перестановками элементов базиса
$\{\gamma_1,\gamma_2,\ldots,\gamma_r\}$ и умножениями его элементов на $\pm
1$. Множество $\Phi_\mathrm{res}$ состоит из корней $\pm\gamma_j$ кратности
1, корней $\pm\frac12\gamma_i\pm\frac12\gamma_j$ кратности $a$ и, возможно,
корней $\pm\frac12\gamma_j$ четной кратности $2b$. Если
$\frac12\gamma_i\in\Phi_\mathrm{res}$, то $\Phi_\mathrm{res}$ является
неприведенной системой корней типа $BC$, а в противном случае --
приведенной системой корней типа $C$.
\end{proposition}

Нетрудно описать сужения простых компактных корней на подпространство
$\mathfrak{h}^-$ \cite[стр.83]{Inoue}:
\begin{itemize}
\item[]
$\{\frac{1}{2}(\gamma_1-\gamma_2),\frac{1}{2}(\gamma_2-\gamma_3),\ldots,
\frac{1}{2}(\gamma_{r-1}-\gamma_r), \frac{1}{2}\gamma_r\}, \qquad \text{в
случае\ \ } BC, $
 \item[] $
\{\frac{1}{2}(\gamma_1-\gamma_2),\frac{1}{2}(\gamma_2-\gamma_3),\ldots,
\frac{1}{2}(\gamma_{r-1}-\gamma_r)\}, \qquad \text{в случае\ \ } C.
$\end{itemize}



\medskip
\begin{remark}
Используя классификацию неприводимых ограниченных симметрических областей
\cite{Helg}, нетрудно показать, что такая область с точностью до изоморфизма
определяется тройкой чисел $r$, $a$, $b$ \cite[стр. 213]{Kor}, \cite[стр.
16]{Arazy}. Число $p=(r-1)a+b+2$ называют родом неприводимой ограниченной
симметрической области $\mathbb{D}$ \cite[стр.16]{Arazy}, \cite[стр.
239]{Kor}.
\end{remark}

\bigskip
В дальнейшем для $\mathfrak{g}$-модулей и $\mathfrak{k}$-модулей мы
используем такие же обозначения, как для $U_q\mathfrak{g}$-модулей и
$U_q\mathfrak{k}$-модулей, опуская индекс $q$.

Вес $\lambda\in P_+$ называют сферическим, \IND{вес ! сферический} если
простой конечномерный $\mathfrak{g}$-модуль $L(\lambda)$ со старшим весом
$\lambda$ является сферическим, \IND{модуль ! сферический} то есть обладает
ненулевым \hbox{$\mathfrak{k}$-инвариантным} вектором. В этом случае
подпространство $\mathfrak{k}$-инвариантных векторов одномерно \cite[стр.
452]{Helg}, и, следовательно, корректно определены отвечающие им матричные
элементы $\varphi_\lambda(g)$ -- зональные сферические функции.
\IND{зональные сферические функции}

Приведем хорошо известное описание множества $P_+^{\,\rm{spher}}$
сферических весов и порожденной им абелевой группы $P^{\,\rm spher}\subset
P$.

\begin{proposition}(\cite[стр. 593]{Helg1})\label{s_w}
 Простой конечномерный $\mathfrak{g}$-модуль
 является сферическим, если и только если  его старший вес $\lambda$
 удовлетворяет следующим условиям:
\begin{itemize}
\item[1.] $\lambda|_{\mathfrak{h}^+}=0$,
 \item[2.]
$\dfrac{(\lambda,\gamma)}{(\gamma,\gamma)}\in\mathbb{Z}$ для любого
ограниченного корня $\gamma$.
\end{itemize}
\end{proposition}

\begin{example} Пусть $n>1$.
Рассмотрим эрмитово симметрическую пару\\
$(\mathfrak{sl}_{n+1},\mathfrak{s}(\mathfrak{gl}_n\times\mathfrak{gl}_1))$.
В этом случае $r=1$,
$$
\gamma_1=\overline{\omega}_1+\overline{\omega}_n=\sum_{i=1}^n\alpha_i,\qquad
H_{\gamma_1}=\sum_{i=1}^n H_i=e_{1,1}-e_{n+1,n+1},
$$
где $\{e_{i,j}\}$ -- стандартный базис векторного пространства матриц.
Значит, система ограниченных корней имеет вид $ \Phi_\mathrm{res}=
\left\{-\gamma_1,-\frac12\gamma_1,\frac12\gamma_1,\gamma_1\right\}$, а
подпространство $\mathfrak{h}^+$ порождается элементами
$$H_1-H_n,\qquad H_j,\quad j=2,3,\ldots,n-1.$$
Условия $\lambda\in P_+$, $\lambda|_{\mathfrak{h}^+}=0$ означают, что
$\lambda\in\mathbb{Z}_+(\overline{\omega}_1+\overline{\omega}_n)$. Условие
целочисленности выполнено для всех таких весов $\lambda$, поскольку
$$
(\gamma_1,\gamma_1)=2,\qquad
(\overline{\omega}_1+\overline{\omega}_n,\gamma_1)=
\left(\overline{\omega}_1+\overline{\omega}_n,\sum_{i=1}^n\alpha_i\right)=2,
$$
как следует из \eqref{inv_inner_prod}, \eqref{inner_omega}.
\end{example}

Известно \cite[стр. 6]{Hoog}, что
$$
P_+^{\,\rm{spher}}\;=\;\newoplus_{i=1}^r\mathbb{Z}_+\,\mu_i,\qquad
P^{\,\rm{spher}}\;=\;\newoplus_{i=1}^r\mathbb{Z}\,\mu_i,
$$
где $\mu_1,\mu_2,\ldots,\mu_r$ -- фундаментальные сферические веса. Их
явный вид приведен, например, в \cite[стр. 478]{Lep}.


\medskip

Наделим $P_+^{\rm spher}$ отношением частичного порядка, полагая
$\sum_{j=1}^r m'_j\mu_j\preceq\sum_{j=1}^rm''_j\mu_j$ если $m'_j\le m''_j$
при всех $j=1,2,\ldots,r$. Нетрудно показать \cite{Vretare},\cite[стр.
297]{VilenkinKlimyk3}, что
\begin{equation}\label{spherical_filtration}
\varphi_{\lambda'}(g)\varphi_{\lambda''}(g)\,=\,
\sum_{\nu\preceq\lambda'+\lambda''}\,c_\nu(\lambda',\lambda'')\,
\varphi_\nu(g),
\end{equation}
где $c_\nu(\lambda',\lambda'')$ -- неотрицательные числа и $
c_{\lambda'+\lambda''}(\lambda',\lambda'')\ne 0$.

\bigskip\medskip

Рассмотрим алгебру $\operatorname{Pol}(\mathfrak{p}^-)$ полиномов на
овеществленном векторном пространстве $\mathfrak{p}^-$ и подалгебру ее
$\mathfrak{k}$-инвариантов
$\operatorname{Pol}(\mathfrak{p}^-)^\mathfrak{k}$. Используя матричные
элементы конечномерных сферических $\mathfrak{g}$-модулей $L(\lambda)$,
отвечающие их $\mathfrak{k}$-инвариантным векторам, нетрудно получить базис
векторного пространства $\operatorname{Pol}(\mathfrak{p}^-)^\mathfrak{k}$.
При этом переход к полиномам на $\mathfrak{p}^-$ осуществляется с помощью
канонического вложения
$\mathbb{C}[\mathfrak{p}^-]\hookrightarrow\mathbb{C}[G]_t$). Следующий
результат доставляет описание алгебры
$\operatorname{Pol}(\mathfrak{p}^-)^\mathfrak{k}$.

\begin{proposition}(\cite[стр. 74]{John})\label{free_k}
Существуют $r$ однородных полиномов $P_1,P_2,\ldots,P_r$ в
$\operatorname{Pol}(\mathfrak{p}^-)^\mathfrak{k}$, для которых

1. $\deg P_i=2i$, \ \ $i=1,2,\ldots,r,$.

2. $\operatorname{Pol}(\mathfrak{p}^-)^\mathfrak{k}$ -- свободная
коммутативная алгебра с образующими $P_1,P_2,\ldots,P_r$.
\end{proposition}

Аналогичный результат имеет место для подалгебры Ли
$\mathfrak{n}^+_\mathfrak{k}\subset\mathfrak{k}$, порожденной элементами
$E_i$, $i\in\mathbb{S}$. Пусть
$\mathbb{C}[\mathfrak{p}^-]^{\operatorname{inv}}$ -- подалгебра
\hbox{$\mathfrak{n}^+_\mathfrak{k}$-инвариантов} алгебры
$\mathbb{C}[\mathfrak{p}^-]$ голоморфных полиномов на $\mathfrak{p}^-$.

\begin{proposition}(\cite[стр. 73]{John})\label{free_n}
Существуют $r$ однородных полиномов $p_1,p_2,\ldots,p_r$ в
$\mathbb{C}[\mathfrak{p}^-]^{\operatorname{inv}}$, для которых

1. $\deg p_i=i$,

2. вес $p_i$ равен $\sum\limits_{k=1}^i \gamma_k$,

3. $\mathbb{C}[\mathfrak{p}^-]^{\operatorname{inv}}$ -- свободная
коммутативная алгебра с образующими $p_1,p_2,\ldots,p_r$.

4. Образующие $p_1,p_2,\ldots,p_r$ алгебры
$\mathbb{C}[\mathfrak{p}^-]^{\operatorname{inv}}$ определяются свойствами 1
-- 3 с точностью до числовых множителей.
\end{proposition}

В формулировке предложения \ref{free_n} подразумевается, что $r$ -- ранг
$\mathbb{D}$.

\bigskip

Следующий результат называется теоремой Хуа-Шмида. Он вытекает из
предложения \ref{free_n}.

\begin{proposition}

(\cite[стр. 73]{John}, \cite[стр. 443]{Takeuchi})\label{k-types}
 В категории $\mathfrak{k}$-модулей
$$\mathbb{C}[\mathfrak{p}^-]\approx
\newoplus_{\lambda\in\mathcal{A}_+}L(\mathfrak{k},\lambda),$$
где
\begin{equation}\label{HuaSchmid_lattice} \mathcal{A}_+=
\newoplus_{j=1}^r\;\, \mathbb{Z}_+\left(\sum\limits_{i=1}^j\gamma_i\right).
\end{equation}
\end{proposition}

\begin{corollary}\label{HuaSchmid_mult_free}
Кратности вхождений простых $\mathfrak{k}$-модулей в
$\mathbb{C}[\mathfrak{p}^-]$ не превосходят единицы.
\end{corollary}

\subsection{Ковариантные дифференциальные исчисления и инвариантные
дифференциальные операторы}\label{quantum_diff_calculus}

\subsubsection{Дифференциальные исчисления первого порядка.}\label{FODC}

Введем необходимые определения, следуя \cite{KlSch}. Рассмотрим алгебру $F$,
которая будет предполагаться унитальной, если не оговорено противное.

Дифференциальным исчислением первого порядка над $F$ \IND{дифференциальное
исчисление первого порядка} называют $F$-бимодуль $M$ и линейное
отображение $d:F\to M$ такие, что
\begin{itemize}
\item[1.] для всех $f_1,f_2\in F$
\begin{equation}\label{Leibnitz_rule}
d(f_1\cdot f_2)=df_1\cdot f_2+f_1\cdot df_2,
\end{equation}
\item[2.] $M$ является линейной оболочкой векторов $f_1\cdot df_2\cdot
    f_3$, где $f_1,f_2,f_3\in F$.
\end{itemize}
Из (\ref{Leibnitz_rule}) следует, что $d 1 = 0$. Второе условие равносильно
тому, что $M$ является линейной оболочкой элементов $f_1\cdot df_2$, где
$f_1,f_2\in F$.

Пусть $A$ -- алгебра Хопфа, и $F$ является $A$-модульной алгеброй.
Дифференциальное исчисление первого порядка $(M, d)$ над $F$ называют
ковариантным, \IND{дифференциальное исчисление первого порядка !
ковариантное} если $M$ является $A$-модульным $F$-бимодулем и $d$ --
морфизм $A$-модулей.

Будем использовать аббревиатуру CFODC для ковариантных дифференциальных
исчислений первого порядка. Два CFODC $(M',d')$, $(M'',d'')$ называют
изоморфными, \IND{дифференциальное исчисление первого порядка !
ковариантное ! изоморфизм} если существует такой изоморфизм $A$-модульных
$F$-бимодулей $i:M'\to M''$, что $d''=i\cdot d'$.

 Ограничимся  эрмитово симметрическим случаем и построим
 CFODC над
$U_q\mathfrak{g}$-модульной алгеброй $\mathbb{C}[\mathfrak{p}^-]_q$.

\medskip Как известно \cite[стр. 82]{BaEast}, в классическом случае $q=1$
дифференциал сопряжен к морфизму обобщенных модулей Верма
\begin{equation}\label{d_conjugate}
 d^*:\;N(\mathfrak{q}^+,-\alpha_{l_0}) \to N(\mathfrak{q}^+,0),\qquad
 d^*:\,v(\mathfrak{q}^+,-\alpha_{l_0})\mapsto F_{l_0}\,v(\mathfrak{q}^+,0)
\end{equation}
 в категории градуированных $U\mathfrak{g}$-модулей.
Перейдем к квантовой универсальной обертывающей алгебре.

\medskip

Пусть $d^*$ -- морфизм $U_q\mathfrak{g}$-модулей, определяемый той же
формулой \eqref{d_conjugate}, что и в классическом случае. Существование и
единственность такого морфизма очевидны, см. \itemiii\ \ref{GVerma}.
Напомним, что векторное пространство $N(\mathfrak{q}^+,-\alpha_{l_0})$
является градуированным $N(\mathfrak{q}^+,-\alpha_{l_0})=
\newoplus\limits_{j=1}^\infty N(\mathfrak{q}^+,-\alpha_{l_0})_{-j}$.
Введем в рассмотрение сопряженное градуированное векторное пространство
$$
\Lambda^1(\mathfrak{p}^-)_q=
\newoplus\limits_{j=1}^\infty\Lambda^1(\mathfrak{p}^-)_{q,j},\qquad
\Lambda^1(\mathfrak{p}^-)_{q,j}=N(\mathfrak{q}^+, -\alpha_{l_0})_{-j}^*
$$
и сопряженный линейный оператор
\begin{equation}\label{d_minus}
d:\mathbb{C}[\mathfrak{p}^-]_q\to\Lambda^1(\mathfrak{p}^-)_q.
\end{equation}

 Покажем,
что $\Lambda^1(\mathfrak{p}^-)_q$ является $U_q\mathfrak{g}$-модульным
$\mathbb{C}[\mathfrak{p}^-]_q$-бимодулем.

Пусть $\lambda \in P_+^\mathbb{S}$. Рассмотрим следующие морфизмы
$$
\Delta_{\mathrm{left},\lambda}^+: N(\mathfrak{q}^+, \lambda) \rightarrow
N(\mathfrak{q}^+, 0) \otimes N(\mathfrak{q}^+, \lambda),$$
$$
\Delta_{\mathrm{right},\lambda}^+: N(\mathfrak{q}^+, \lambda) \rightarrow
N(\mathfrak{q}^+, \lambda) \otimes N(\mathfrak{q}^+, 0)$$
 в категории $U_q\mathfrak{g}^{\mathrm{cop}}$-модулей, определяемые их
действием на образующие
$$
\Delta_{\mathrm{left},\lambda}^+: v(\mathfrak{q}^+, \lambda) \mapsto
v(\mathfrak{q}^+, 0) \otimes v(\mathfrak{q}^+, \lambda),
$$
$$
\Delta_{\mathrm{right},\lambda}^+: v(\mathfrak{q}^+, \lambda) \mapsto
v(\mathfrak{q}^+, \lambda) \otimes v(\mathfrak{q}^+, 0).
$$
Следующие равенства вытекают из определений
\begin{equation}\label{com1}
  (\mathrm{id} \otimes \Delta_{\mathrm{left}, \lambda}^+)
\Delta_{\mathrm{left}, \lambda}^+ = (\Delta^+ \otimes \mathrm{id})
\Delta_{\mathrm{left}, \lambda}^+,
\end{equation}
\begin{equation}\label{com2}
  (\Delta_{\mathrm{right}, \lambda}^+ \otimes \mathrm{id})
\Delta_{\mathrm{right}, \lambda}^+ = (\mathrm{id} \otimes \Delta^+)
\Delta_{\mathrm{right}, \lambda}^+,
\end{equation}
\begin{equation}\label{com3}
  (\varepsilon^+ \otimes \mathrm{id}) \Delta_{\mathrm{left}, \lambda}^+
= (\mathrm{id} \otimes \varepsilon^+) \Delta_{\mathrm{right}, \lambda}^+ =
\mathrm{id},
\end{equation}
\begin{equation}\label{com4}
  (\mathrm{id} \otimes \Delta_{\mathrm{right}, \lambda}^+)
\Delta_{\mathrm{left}, \lambda}^+ = (\Delta_{\mathrm{left}, \lambda}^+
\otimes \mathrm{id}) \Delta_{\mathrm{right}, \lambda}^+
\end{equation}
(коумножение $\Delta^+$ и коединица $\varepsilon^+$ в $N(\mathfrak{q}^+,0)$
 были введены в \itemiiiе \ref{vect}.)
 Например, \eqref{com4} следует из
$$
(\mathrm{id} \otimes \Delta_{\mathrm{right}, \lambda}^+)
\Delta_{\mathrm{left}, \lambda}^+ v(\mathfrak{q}^+, \lambda) =
v(\mathfrak{q}^+, 0) \otimes v(\mathfrak{q}^+, \lambda) \otimes
v(\mathfrak{q}^+, 0),
$$
$$
(\Delta_{\mathrm{left}, \lambda}^+ \otimes \mathrm{id})
\Delta_{\mathrm{right}, \lambda}^+ v(\mathfrak{q}^+, \lambda) =
v(\mathfrak{q}^+, 0) \otimes v(\mathfrak{q}^+, \lambda) \otimes
v(\mathfrak{q}^+, 0).
$$

Рассмотрим категорию градуированных векторных пространств и сопряженное к
$N(\mathfrak{q}^+,\lambda)$ пространство $\Gamma(\mathfrak{p}^-,\lambda)_q$
в этой категории. Наделим его структурой $U_q\mathfrak{g}$-модуля, а также
левым и правым действиями $\mathbb{C}[\mathfrak{p}^-]_q$ с помощью
операторов
$$
m_{\operatorname{left},\lambda}\stackrel{\operatorname{def}}{=}
(\Delta_{\mathrm{left},\lambda}^+)^*,\qquad
m_{\operatorname{right},\lambda}\stackrel{\operatorname{def}}{=}
(\Delta_{\mathrm{right},\lambda}^+)^*.
$$
Из равенств (\ref{com1}) -- (\ref{com4}) вытекает следующее утверждение.
\begin{proposition}\label{homog_bundle}
Для любого $\lambda\in P_+^\mathbb{S}$ градуированное векторное
пространство $\Gamma(\mathfrak{p}^-,\lambda)_q$ является
$U_q\mathfrak{g}$-модульным $\mathbb{C}[\mathfrak{p}^-]_q$-бимодулем.
\end{proposition}

\begin{remark}\label{the_same_V} Очевидно,
$$
\Gamma(\mathfrak{p}^-,\lambda)_q\cong
\left\{\left.f\in\mathrm{Hom}_{_{U_q\mathfrak{q}^+}}(U_q\mathfrak{g},
L(\mathfrak{q}^+,-w_0\lambda))\;\right|\;\dim(U_q\mathfrak{h}\cdot f)<\infty
\right\}.
$$
\end{remark}

Пусть $N(\mathfrak{q}^+,\lambda)_{\mathrm{highest}}$ -- старшая однородная
компонента градуированного векторного пространства
$N(\mathfrak{q}^+,\lambda)$ и $P_{\mathrm{highest}}$ -- проектор в
$N(\mathfrak{q}^+,\lambda)$ на подпространство
$N(\mathfrak{q}^+,\lambda)_{\mathrm{highest}}$ параллельно сумме остальных
однородных компонент.

\begin{lemma}\label{for_freedom}
Линейные операторы
$$
(\mathrm{id} \otimes P_{\mathrm{highest}}) \Delta_{\mathrm{left},
\lambda}^+:\; N(\mathfrak{q}^+, \lambda) \rightarrow N(\mathfrak{q}^+,0)
\otimes N(\mathfrak{q}^+, \lambda)_{\mathrm{highest}},
$$
$$
(P_{\mathrm{highest}} \otimes \mathrm{id}) \Delta_{\mathrm{right},
\lambda}^+: \; N(\mathfrak{q}^+, \lambda) \rightarrow N(\mathfrak{q}^+,
\lambda)_{\mathrm{highest}} \otimes N(\mathfrak{q}^+,0)
$$
инъективны.
\end{lemma}

{\bf Доказательство.}
 Воспользуемся теми же обозначениями, что в формулировке предложения
\ref{N_basis}. Инъективность рассматриваемых линейных операторов вытекает из
этого предложения, изоморфизма $N(\mathfrak{q}^+,
\lambda)_{\mathrm{highest}} \cong L(\mathfrak{q}^+, \lambda)$ и из того, что
для любого весового вектора ${v \in N(\mathfrak{q}^+,
\lambda)_{\mathrm{highest}}}$ имеют место равенства
$$
(\mathrm{id} \otimes P_{\mathrm{highest}}) \Delta_{\mathrm{left}, \lambda}^+
(F_{\beta_M}^{j_M} F_{\beta_{M-1}}^{j_{M-1}} \ldots
F_{\beta_{M'+1}}^{j_{M'+1}} v)=$$ $$= \mathrm{const}' \cdot
F_{\beta_M}^{j_M} F_{\beta_{M-1}}^{j_{M-1}} \ldots
F_{\beta_{M'+1}}^{j_{M'+1}} v(\mathfrak{q}^+,0) \otimes v,
$$
$$
(P_{\mathrm{highest}} \otimes \mathrm{id}) \Delta_{\mathrm{right},
\lambda}^+ (F_{\beta_M}^{j_M} F_{\beta_{M-1}}^{j_{M-1}} \ldots
F_{\beta_{M'+1}}^{j_{M'+1}} v)=$$ $$= \mathrm{const}'' \cdot v \otimes
F_{\beta_M}^{j_M} F_{\beta_{M-1}}^{j_{M-1}} \ldots
F_{\beta_{M'+1}}^{j_{M'+1}} v(\mathfrak{q}^+,0),
$$
где числовые множители $\mathrm{const}'$, $\mathrm{const}''$ отличны от
нуля.
\hfill $\square$

\medskip
\begin{proposition}\label{freedom}
1. Бимодуль $\Gamma(\mathfrak{p}^-, \lambda)_q$ над алгеброй
$\mathbb{C}[\mathfrak{p}^-]_q$ является как свободным левым, так и свободным
правым $\mathbb{C}[\mathfrak{p}^-]_q$-модулем.
\\
2. Если $\Gamma(\mathfrak{p}^-, \lambda)_{q,{\mathrm{lowest}}}$ -- младшая
однородная компонента $\Gamma(\mathfrak{p}^-, \lambda)_q$, то имеют место
изоморфизмы $U_q\mathfrak{k}$-модулей
\begin{equation}\label{morph1}
  \mathbb{C}[\mathfrak{p}^-]_q \otimes
\Gamma(\mathfrak{p}^-, \lambda)_{q,{\mathrm{lowest}}}
\stackrel{\approx}{_\rightarrow} \Gamma(\mathfrak{p}^-, \lambda)_q, \quad f
\otimes v \mapsto f v,
\end{equation}
\begin{equation}\label{morph2}
  \Gamma(\mathfrak{p}^-, \lambda)_{q,{\mathrm{lowest}}} \otimes
\mathbb{C}[\mathfrak{p}^-]_q \stackrel{\approx}{_\rightarrow}
\Gamma(\mathfrak{p}^-, \lambda)_q, \quad v \otimes f \mapsto v f.
\end{equation}
\end{proposition}

{\bf Доказательство.}
 Первое утверждение следует из второго. Морфизмы $U_q\mathfrak{t}$-модулей
(\ref{morph1}), (\ref{morph2}) являются морфизмами градуированных векторных
пространств. Из предложения \ref{N_basis} (по двойственности) следует, что
размерности соответствующих их однородных компонент конечны и равны. Значит,
биективность морфизмов (\ref{morph1}), (\ref{morph2}) следует из
сюръективности, а их сюръективность вытекает из леммы \ref{for_freedom}.
\hfill $\square$

\bigskip Разумеется, $\Gamma(\mathfrak{p}^-, \lambda)_q$ не является
свободным $\mathbb{C}[\mathfrak{p}^-]_q$-{\sl бимодулем}. Найдем соотношения
между образующими из $\Gamma(\mathfrak{p}^-,
\lambda)_{q,{\mathrm{lowest}}}$. Так как $\Gamma(\mathfrak{p}^-, \lambda)_q$
является $U_q\mathfrak{g}$-модулем с младшим весом, то корректно определен
морфизм $U_q\mathfrak{g}$-модулей
$$
\check{R}_{\mu, \lambda} \stackrel{\mathrm{def}}{=}
\check{R}_{\Gamma(\mathfrak{p}^-, \mu)_q, \Gamma(\mathfrak{p}^-, \lambda)_q}
: \Gamma(\mathfrak{p}^-, \mu)_q \otimes \Gamma(\mathfrak{p}^-, \lambda)_q
\rightarrow \Gamma(\mathfrak{p}^-, \lambda)_q \otimes \Gamma(\mathfrak{p}^-,
\mu)_q,
$$
введенный в \itemiiiе $\ref{RM}$ с помощью универсальной $R$-матрицы. В
частности,
$$
\check{R}_{0, \lambda}: \mathbb{C}[\mathfrak{p}^-]_q \otimes
\Gamma(\mathfrak{p}^-, \lambda)_q \rightarrow \Gamma(\mathfrak{p}^-,
\lambda)_q \otimes \mathbb{C}[\mathfrak{p}^-]_q,
$$
поскольку в категории $U_q\mathfrak{g}$-модульных
$\mathbb{C}[\mathfrak{p}^-]_q$-бимодулей $$\mathbb{C}[\mathfrak{p}^-]_q
\cong \Gamma(\mathfrak{p}^-, 0)_q.$$

Первое утверждение следующего предложения означает, что рассматриваемые
бимодули над коммутативной алгеброй $\mathbb{C}[\mathfrak{p}^-]_q$ в
категории $\mathcal{C}^-$ являются локальными, см. \itemiii\
\ref{braided_categories_sl_2}.

\begin{proposition}\label{rel_1}
\begin{itemize}
\item[1.]\ $m_{\mathrm{left, \lambda}} = m_{\mathrm{right, \lambda}} \cdot
\check{R}_{0, \lambda}$,\ \ \ $ m_{\mathrm{right, \mu}} = m_{\mathrm{left,
\mu}} \cdot \check{R}_{\mu, 0}$.
 \item[2.]\ \ $\check{R}_{0, \lambda}:\;
\mathbb{C}[\mathfrak{p}^-]_{q,1} \otimes \Gamma(\mathfrak{p}^-,
\lambda)_{q,{\mathrm{lowest}}} \rightarrow \Gamma(\mathfrak{p}^-,
\lambda)_{q,{\mathrm{lowest}}} \otimes \mathbb{C}[\mathfrak{p}^-]_{q,1}$.
\end{itemize}
\end{proposition}

{\bf Доказательство.}
 Первое равенство доказывается переходом к
 сопряженным операторам. Действительно,
оно равносильно равенству морфизмов $U_q\mathfrak{g}^{\mathrm{cop}}$-модулей
\begin{equation}\label{adjoint_eq}
  (\check{R}_{0, \lambda})^* \Delta_{\mathrm{right, \lambda}}^+ =
\Delta_{\mathrm{left}, \lambda}^+.
\end{equation}
В свою очередь, (\ref{adjoint_eq}) следует из того, что вектор
$v(\mathfrak{q}^+, \lambda)$ порождает
$U_q\mathfrak{g}^{\mathrm{cop}}$-модуль $N(\mathfrak{q}^+, \lambda)$ и из
равенств
$$ \Delta_{\mathrm{right}, \lambda}^+ v(\mathfrak{q}^+, \lambda)
= v(\mathfrak{q}^+, \lambda) \otimes v(\mathfrak{q}^+,0), \quad
\Delta_{\mathrm{left}, \lambda}^+ v(\mathfrak{q}^+, \lambda) =
v(\mathfrak{q}^+,0) \otimes v(\mathfrak{q}^+, \lambda),
$$
$$
\check{R}_{0, \lambda}^* \; v(\mathfrak{q}^+, \lambda) \otimes
v(\mathfrak{q}^+,0) = v(\mathfrak{q}^+,0) \otimes v(\mathfrak{q}^+,
\lambda).
$$
Второе равенство доказывается аналогично первому.

Последнее из утверждений передложения \ref{rel_1} является следствием того,
что $\check{R}_{0, \lambda}$ является морфизмом $U_q\mathfrak{g}$-модулей и,
следовательно, этот линейный оператор перестановочен с действием
$H_\mathbb{S} \otimes 1 +1\otimes H_\mathbb{S}$. \hfill $\square$

\begin{corollary}
Пусть $\{ z_i\}$ -- базис конечномерного векторного пространства
$\mathbb{C}[\mathfrak{p}^-]_{q,1}$, а $\{\gamma_j\}$ -- базис
конечномерного векторного пространства
$\Gamma(\mathfrak{p}^-,\lambda)_{q,{\mathrm{lowest}}}$. Существует и
единственна такая матрица $(\check{R}_{ij}^{km}(\lambda))$, что
\begin{equation}\label{finite_R}
\check{R}_{0,\lambda}(z_i\otimes\gamma_j)=
\sum\limits_{k,m}\check{R}_{ij}^{km}(\lambda)\,\gamma_k\otimes z_m,
\end{equation}
где \ $i,m\in 1,2,\ldots,\dim\mathbb{C}[\mathfrak{p}^-]_{q,1}$, \ $j,k\in
1,2,\ldots,\dim\Gamma(\mathfrak{p}^-, \lambda)_{q,\mathrm{lowest}}$.
\end{corollary}

\medskip

\begin{remark}\label{comm_mu_lambda}
Отыскание функций $\check{R}_{ij}^{km}(\lambda)$ кажется непростой задачей,
поскольку в их определении участвует действие универсальной
\hbox{$R$-матрицы} в тензорных произведениях бесконечномерных
$U_q\mathfrak{g}$-модулей. Покажем, что на самом деле проблемы нет. Так же,
как в \itemiiiе \ref{quadratic_first}, рассмотрим подалгебру Хопфа
$U_q\mathfrak{k}_\mathrm{ss}\subset U_q\mathfrak{g}$, порожденную
элементами $K^{\pm 1}_i$, $E_i$, $F_i$ с $i\ne l_0$. Действие универсальной
$R$-матрицы алгебры Хопфа $U_q\mathfrak{g}$ на векторы из
$\mathbb{C}[\mathfrak{p}^-]_{q,1}\otimes
\Gamma(\mathfrak{p}^-,\lambda)_{q,{\mathrm{lowest}}}$ отличается от
действия универсальной \hbox{$R$-матрицы} алгебры Хопфа
$U_q\mathfrak{k}_\mathrm{ss}$ лишь числовым множителем
\begin{equation}\label{const_R_new}
\operatorname{const}=q_{l_0}^{\frac{(\lambda,\overline{\omega}_{l_0})}
{(\overline{\omega}_{l_0},\overline{\omega}_{l_0})}}.
\end{equation}
Для доказательства \eqref{const_R_new} достаточно воспользоваться формулой
(\ref{Rmatrix}), приведенным разложением элемента $w_0$, описанным в
\itemiiiе \ref{GVerma}, и равенствами \eqref{t_canon},
$$
(\alpha_{l_0},\lambda)=(\alpha_{l_0\,|\mathfrak{h}\cap\mathfrak{k}},
\lambda_{|\mathfrak{h} \cap \mathfrak{k}})+
\frac{(\alpha_{l_0},\overline{\omega}_{l_0})
(\overline{\omega}_{l_0},\lambda)}
{(\overline{\omega}_{l_0},\overline{\omega}_{l_0})},\qquad
(\alpha_{l_0},\overline{\omega}_{l_0})=d_{l_0},
$$
позволяющими сравнить картановские множители (см. (\ref{Rmatrix}) и
замечание \ref{long_root}). Так как $\check{R}_{\mu,\lambda}$ является
морфизмом $U_q\mathfrak{g}$-модулей, то
\begin{equation*}
\check{R}_{\mu,\lambda}:\;\Gamma(\mathfrak{p}^-,\mu)_{q,{\mathrm{lowest}}}
\otimes\Gamma(\mathfrak{p}^-,\lambda)_{q,{\mathrm{lowest}}}\to
\Gamma(\mathfrak{p}^-,\lambda)_{q,{\mathrm{lowest}}}\otimes
\Gamma(\mathfrak{p}^-,\mu)_{q,{\mathrm{lowest}}} .
\end{equation*}
Так же, как прежде, находится числовой множитель, связывающий действия
$R$-матриц алгебр Хопфа $U_q\mathfrak{g}$ и $U_q\mathfrak{k}_\mathrm{ss}$:
\begin{equation}\label{const_R}
\operatorname{const}=q^{-\frac{(\mu,\overline{\omega}_{l_0})
(\lambda,\overline{\omega}_{l_0})}
{(\overline{\omega}_{l_0},\overline{\omega}_{l_0})}}.
\end{equation}
\end{remark}

\medskip

Следующее утверждение вытекает из предложений \ref{freedom}, \ref{rel_1}.
\begin{proposition}\label{rel_G}
Множество $\{\gamma_j\}_{j=1,2,\ldots,
\dim\Gamma(\mathfrak{p}^-,\lambda)_{q,\mathrm{lowest}}} $ порождает
$\mathbb{C}[\mathfrak{p}^-]_q$-бимодуль $\Gamma(\mathfrak{p}^-,\lambda)_q$,
и соотношения
\begin{equation}\label{rel_Gamma}
z_i\,\gamma_j=\sum\limits_{k,m}\check{R}_{ij}^{km}(\lambda)\gamma_k\,z_m
\end{equation}
являются определяющими.
\end{proposition}

 По построению, $\Lambda^1(\mathfrak{p}^-)_q
= \Gamma(\mathfrak{p}^-, -\alpha_{l_0})_q$. Значит,
$\Lambda^1(\mathfrak{p}^-)_q$ является $U_q\mathfrak{g}$-модульным
$\mathbb{C}[\mathfrak{p}^-]_q$-бимодулем.

\begin{corollary}\label{first_order_rel}
При всех $i,j=1,2,\ldots,\dim\,\mathfrak{p}^-$
\begin{equation}\label{rel_z_dz}
z_i\,dz_j=\sum\limits_{k,m=1}^{\dim\,\mathfrak{p}^-}\check{R}_{ij}^{km}\,
dz_k\,z_m,
\end{equation}
где $\check{R}_{ij}^{km}=\check{R}_{ij}^{km}(-\alpha_{l_0})$, и
\eqref{rel_z_dz} -- определяющий список соотношений между образующими
$\mathbb{C}[\mathfrak{p}^-]_q$-бимодуля $\Lambda^1(\mathfrak{p}^-)_q$.
\end{corollary}

\begin{remark}\label{z_gamma_low}
Пусть $\gamma_{\mathrm{low}}$ -- младший вектор $U_q\mathfrak{g}$-модуля
$\Gamma(\mathfrak{p}^-,\lambda)_q$ и $\{z_i\}$ -- базис весовых векторов
$U_q\mathfrak{k}$-модуля $\mathbb{C}[\mathfrak{p}^-]_{q,1}$:
$H_jz_i=\mu_i(H_j)z_i$. Тогда, как следует из (\ref{rel_Gamma}),
(\ref{Rmatrix}) и (\ref{fraction}),
$z_i\gamma_{\mathrm{low}}=q^{-(\lambda,\mu_i)}\gamma_{\mathrm{low}}z_i$. В
частности,
\begin{equation}\label{zl0_dzl0}
z_{\mathrm{low}}\;dz_{\mathrm{low}}=q_{l_0}^{-2}\,dz_{\mathrm{low}}\;
z_{\mathrm{low}}.
\end{equation}
Это -- обобщение соотношения \eqref{z_dz_sl_2}.
\end{remark}

\bigskip

Перейдем к доказательству того, что $(\Lambda^1(\mathfrak{p}^-)_q, d)$
является дифференциальным исчислением первого порядка над
$\mathbb{C}[\mathfrak{p}^-]_q$. Правило Лейбница \eqref{Leibnitz_rule}
нетрудно получить, перейдя к сопряженным линейным операторам в следующем
равенстве морфизмов $U_q\mathfrak{g}^{\mathrm{cop}}$-модулей.

\begin{lemma}
$ \Delta^+ d^* = (d^* \otimes \mathrm{id}) \Delta_{\mathrm{right},
-\alpha_{l_0}}^+ + (\mathrm{id} \otimes d^*) \Delta_{\mathrm{left},
-\alpha_{l_0}}^+$.
\end{lemma}
{\bf Доказательство.} Достаточно воспользоваться равенствами
$$
\Delta^+ d^* v(\mathfrak{q}^+, -\alpha_{l_0}) = \Delta^+(F_{l_0}
v(\mathfrak{q}^+,0)) = F_{l_0}(\Delta^+ v(\mathfrak{q}^+,0)) =
$$
$$
= F_{l_0}(v(\mathfrak{q}^+,0) \otimes v(\mathfrak{q}^+,0)) = F_{l_0}
v(\mathfrak{q}^+,0) \otimes v(\mathfrak{q}^+,0) + v(\mathfrak{q}^+,0)
\otimes F_{l_0} v(\mathfrak{q}^+,0),
$$
$$
((d^* \otimes \mathrm{id}) \Delta_{\mathrm{right}, -\alpha_{l_0}}^+ +
(\mathrm{id} \otimes d^*) \Delta_{\mathrm{left}, -\alpha_{l_0}}^+)
v(\mathfrak{q}^+,-\alpha_{l_0}) =$$ $$=d^* v(\mathfrak{q}^+, -\alpha_{l_0})
\otimes v(\mathfrak{q}^+,0) + v(\mathfrak{q}^+,0) \otimes d^*
v(\mathfrak{q}^+, -\alpha_{l_0}) =
$$
$$
= F_{l_0} v(\mathfrak{q}^+,0) \otimes v(\mathfrak{q}^+,0) +
v(\mathfrak{q}^+,0) \otimes F_{l_0} v(\mathfrak{q}^+,0). \eqno \square
$$

\bigskip
Остается показать, что элементы $df$, где $f \in
\mathbb{C}[\mathfrak{p}^-]_q$, порождают
$\mathbb{C}[\mathfrak{p}^-]_q$-бимодуль $\Lambda^1(\mathfrak{p}^-)_q$. Мы
получим более сильное утверждение.

\begin{proposition}\label{end_of_story}
Если $\{ z_j \}$ -- базис векторного пространства
$\mathbb{C}[\mathfrak{p}^-]_{q,1}$, то
$$
\Lambda^1(\mathfrak{p}^-)_q = \newoplus_{j=1}^{\dim\mathfrak{p}^-}
\mathbb{C}[\mathfrak{p}^-]_q\, d z_j =
\newoplus_{j=1}^{\dim\mathfrak{p}^-}\,
 dz_j \,\mathbb{C}[\mathfrak{p}^-]_q.
$$
\end{proposition}

{\bf Доказательство.}
 Как следует из предложения \ref{freedom}, достаточно показать, что линейное
отображение
$$d_{|\mathbb{C}[\mathfrak{p}^-]_{q,1}}: \mathbb{C}[\mathfrak{p}^-]_{q,1} \rightarrow
\Lambda^1(\mathfrak{p}^-)_{q,1}$$ биективно. Но сопряженное линейное
отображение является морфизмом простых весовых $U_q\mathfrak{k}$-модулей и
отлично от нуля, поскольку $$d^* v(\mathfrak{q}^+, -\alpha_{l_0})=
v(\mathfrak{q}^+, 0) \neq 0. \eqno \square$$

\subsubsection{Универсальные обертывающие дифференциальные
исчисления.}\label{diff_univ}

Будем использовать понятия, введенные в
\itemiiiе \ref{pol-forms_sl_2}. Напомним, что дифференциальным исчислением
    \IND{дифференциальное исчисление} над алгеброй $F$ называют такую
    дифференциальную градуированную алгебру $(\Omega,d)$, для которой
    $\Omega_0=F$ и которая порождается элементами из $F\oplus dF$.

Как известно \cite[стр. 463]{KlSch}, каждому дифференциальному исчислению
первого порядка $(M,d)$ над алгеброй $F$ можно сопоставить дифференциальное
исчисление $(\Omega^{\operatorname{univ}},d^{\operatorname{univ}})$ над этой
алгеброй, обладающее следующими свойствами:
\begin{itemize}
\item[1.] $\Omega^{\operatorname{univ}}_{ 1}=M$;

\item[2.] $d^{\operatorname{univ}}_{\quad|F}=d$;

\item[3.] для любого дифференциального исчисления $(\Omega',d')$ над
    алгеброй $F$, удовлетворяющего двум предыдущим требованиям
    ($\Omega'_1=M$, $d'|_F=d$), существует гомоморфизм дифференциальных
    градуированных алгебр $\Omega^{\operatorname{univ}}\to\Omega'$,
    тождественный на $F\oplus M$.
\end{itemize}
Такое дифференциальное исчисление единственно с точностью до изоморфизма,
тождественного на $F\oplus M$, и называется универсальным обертывающим
дифференциальным исчислением. \IND{универсальное обертывающее
дифференциальное исчисление} Оно существует и является ковариантным, если
ковариантным было исходное дифференциальным исчисление первого порядка
\cite[стр. 463, 464]{KlSch}. \IND{универсальное обертывающее
дифференциальное исчисление ! ковариантное}

\begin{remark}\label{our_covariance} В \cite{KlSch} используется
другая терминология и неявно предполагается, что все рассматриваемые
$U_q\mathfrak{g}$-модули являются локально конечномерными. Это не
сказывается на приведенном в \cite[стр. 464]{KlSch} доказательстве
ковариантности универсального обертывающего дифференциального исчисления.
\end{remark}

\bigskip
Приведем три примера построения универсального дифференциального исчисления.
Каждое следующее дифференциальное исчисление первого порядка будет
получаться из предыдущего ''увеличением списка определяющих соотношений''.

\begin{example}\label{tensors}
(Нет соотношений.) Рассмотрим векторное пространство $V$ и построим
дифференциальное исчисление над его тензорной алгеброй $T(V)$. Выберем
векторное пространство $V'$, изоморфное $V$, и изоморфизм $d: V\rightarrow
V'$. Наделим тензорную алгебру $\Omega'=T(V\oplus V')$ следующей
градуировкой
$$\deg v=0,\quad v\in V,\qquad\qquad\deg v'=1,\quad v'\in V'.$$
Очевидно, \ $\Omega'=\newoplus\limits_{j\in\mathbb{Z}_+}\Omega'_j$, где
$\Omega'_j=\{t\in\Omega'\,|\,\deg\,t=j\}$. Определим линейный оператор $d'$
в $\Omega'$ рекурсивно: $d' 1=0$,
$$d'v=dv,\;v\in V,\qquad d'v'=0,\;v'\in V'$$
$$d'(v\otimes t)=dv\otimes t+v\otimes d't,\qquad v\in V,\;t\in\Omega',$$
$$d'(v'\otimes t)=-v'\otimes d't,\qquad v'\in V',\;t\in\Omega'.$$
Равенства $(d')^2=0$ и \eqref{Leibnitz_form} легко следуют из определений.
Пара $(\Omega'_1,d'_{|\Omega'_0})$ является дифференциальным исчислением
первого порядка над тензорной алгеброй $T(V)$, а пара $(\Omega',d')$ -- его
универсальным обертывающим дифференциальным исчислением.
\end{example}

\begin{example}\label{J_0}
(вводятся соотношения между координатами.) Пусть $J_0$ -- двусторонний
идеал тензорной алгебры $T(V)$ и $J_0=\newoplus\limits_{j\ge 2}(J_0\cap
V^{\otimes\,j})$. Построим дифференциальное исчисление над алгеброй
$F=T(V)/J_0$.

Введем двусторонний идеал $J_F$ алгебры $\Omega'$, порожденный $J_0$ и
$d'J_0$. Алгебра $\Omega^F=\Omega'/J_F$ наследует градуировку
$\Omega^F=\newoplus\limits_{j\in\mathbb{Z}_+}\Omega^F_j$.

Очевидно, $V\hookrightarrow F$. Так как $d'J_F\subset J_F$, то действие
дифференциала $d'$ переносится в $\Omega^F=\Omega'/J_F$. Получаем линейный
оператор $d_F$ и дифференциальное исчисление $(\Omega^F,d^F)$ над алгеброй
$F=\Omega^F_0$. Оно является универсальным обертывающим для
дифференциального исчисления первого порядка $(\Omega^F_1,d^F|_F)$ и хорошо
известно в теории квантовых групп \cite[стр. 462]{Sch}.
\end{example}

\begin{example}\label{J_1}
(вводятся коммутационные соотношения между координатами и дифференциалами.)
Рассмотрим обратимый линейный оператор $\check{\mathcal{R}}:V\otimes
V'\rightarrow V'\otimes V$. Пусть $J_1$ -- подбимодуль $F$-бимодуля
$\Omega^F_1$, порожденный подпространством
\begin{equation}\label{rel_space}
\{vv'-v'v\,|\,v\otimes v'=\check{\mathcal{R}}(v\otimes v'),\quad v\in
V,\,v'\in V'\}\subset\Omega^F_1.
\end{equation}
Факторизация $(\Omega^F_1,d^F|_F)$ по $J_1$ приводит к дифференциальному
исчислению первого порядка над алгеброй $F$. Для того, чтобы получить
соответствующее универсальное обертывающее дифференциальное исчисление,
рассмотрим двусторонний идеал $J_M$ алгебры $\Omega_F$, порожденный
множествами $J_1$, $d^F J_1$, и профакторизуем $\Omega_F$ по $J_M$. На
$\Omega_F/J_M$ переносятся градуировка и действие дифференциала, поскольку
идеал $J_M$ однороден и $d^F J_M \subset J_M$.
\end{example}

\bigskip

Рассмотрим введенное в \itemiiiе \ref{FODC} ковариантное дифференциальное
исчисление первого порядка $(\Lambda^1(\mathfrak{p}^-)_q,d)$ над алгеброй
$\mathbb{C}[\mathfrak{p}^-]_q$. Пусть $(\Lambda(\mathfrak{p}^-)_q,d)$ --
отвечающее ему универсальное дифференциальное исчисление над
$\mathbb{C}[\mathfrak{p}^-]_q$ и $\wedge$ -- умножение в
$\Lambda(\mathfrak{p}^-)_q=\newoplus\limits_{j\in\mathbb{Z}_+}
\Lambda^j(\mathfrak{p}^-)_q$.

 Нетрудно описать дифференциальную алгебру
$\Lambda(\mathfrak{p}^-)_q$ в терминах образующих и соотношений, используя
результаты \itemiiiа \ref{FODC} и пример \ref{J_1}. Для этого достаточно
пополнить список определяющих соотношений квадратичной алгебры
$\mathfrak{p}^-$ соотношениями \eqref{rel_z_dz} и равенствами, получаемыми
  дифференцированием \eqref{rel_z_dz}. Именно,
\begin{equation}\label{rel_dz_dz}
dz_i\wedge
dz_j=-\sum\limits_{k,m=1}^{\dim\,\mathfrak{p}^-}\check{R}_{ij}^{km}dz_k
\wedge dz_m.
\end{equation}

\bigskip

\begin{remark} \label{second_def_diff}  Получим с помощью результатов
\itemiiiа \ref{braided_categories_sl_2} другое описание дифференциального
    исчисления $(\Lambda(\mathfrak{p}^-)_q,d)$. Из него будет очевидна
    ковариантность рассматриваемого дифференциального исчисления.
    Во-первых, алгебра $F=\mathbb{C}[\mathfrak{p}^-]_q$ коммутативна в
    категории $\mathcal{C}^-$. Это означает, что $F$ является объектом
    $\mathcal{C}^-$, умножение $m:F\otimes F\to F$ -- морфизмом в
    $\mathcal{C}^-$ и $m\;=\;m\,\check{R}_{FF}$. Во-вторых,
    $M=\Lambda^1(\mathfrak{p}^-)_q$ является локальным бимодулем над $F$,
    то есть $M$ является объектом $\mathcal{C}^-$, действия
    $m_\mathrm{left}:F\otimes M\to M$, $m_\mathrm{right}:M\otimes F\to M$
    -- морфизмами в $\mathcal{C}^-$ и
\begin{equation}\label{FODC_in_C}
\qquad m_\mathrm{right}=m_\mathrm{left}\,\check{R}_{MF}\qquad
m_\mathrm{left}=m_\mathrm{right}\,\check{R}_{FM}.
\end{equation}
В-третьих, $d:F \to M$ -- морфизм в $\mathcal{C}^-$ и $(M,d)$ является
дифференциальным исчислением первого порядка над $F$.

Согласно предложению \ref{dislectic}, автоморфизм $\check{R}_{MM}$ в
$M\otimes M$ порождает автоморфизм $\bar{R}_{MM}$ в $M\otimes_FM$. Именно,
$\bar{R}_{MM}\cdot j\,=\,j\cdot\check{R}_{MM}$, где $j$ -- канонический
эпиморфизм $M\otimes M\to M\otimes_FM$.

Рассмотрим тензорную алгебру
$$
T_F(M)=\newoplus\limits_{j\in \mathbb{Z}_+}T_F(M)_j,\qquad
T_F(M)_j=\underbrace{M\otimes_FM\otimes_F\cdots\otimes_F M}_j\;.
$$
Здесь $T_F(M)_0=F$, $T_F(M)_1=M$. Пусть $J$ -- двусторонний идеал этой
алгебры, порожденный множеством
$$
\{\omega_1\otimes_F\omega_2\,+\,\bar{R}_{MM}(\omega_1\otimes_F\omega_2)\;|\;
\omega_1,\omega_2\in M\}.
$$

 Из определений следует, что $\Lambda(\mathfrak{p}^-)_q \cong T_F(M)/J$,
и, следовательно, $\Lambda(\mathfrak{p}^-)_q$ является
$U_q\mathfrak{g}$-модульной алгеброй. Остается заметить, что, во-первых, $d$
является морфизмом $U_q\mathfrak{g}$-модулей, поскольку
$$
d\,(f_0\wedge df_1 \wedge df_2 \wedge \cdots \wedge df_n)\;=\; df_0\wedge
df_1 \wedge df_2 \wedge \cdots \wedge df_n
 $$
при всех $n \in \mathbb{Z}_+$, $f_0,\ldots,f_n \in F$, и, во-вторых,
дифференциальное исчисления первого порядка $(M,d)$ ковариантно.
\end{remark}

\bigskip
Докажем, что $\Lambda^j(\mathfrak{p}^-)_q$ является свободным левым
$\mathbb{C}[\mathfrak{p}^-]_q$-модулем ранга $\binom{\dim\mathfrak{p}^-}{j}$
и свободным правым $\mathbb{C}[\mathfrak{p}^-]_q$-модулем ранга
$\binom{\dim\mathfrak{p}^-}{j}$, как в классическом случае $q=1$.

Отметим, что линейная оболочка $j$-форм с постоянными коэффициентами
$$
dz_{i_1}\wedge dz_{i_2}\wedge\cdots\wedge dz_{i_j},\qquad i_1,i_2,\cdots,i_j
\in\{1,2,\cdots,\dim\mathfrak{p}^-\},
$$
является $U_q\mathfrak{k}$-модулем. Введем для нее обозначение
$\Lambda^j(\mathfrak{p}^-)^\mathrm{const}_q$.

\begin{proposition}\label{2_forms}
Имеют место изоморфизмы $U_q\mathfrak{k}$-модулей
\begin{equation}\label{morph1-2forms}
\mathbb{C}[\mathfrak{p}^-]_q\otimes
\Lambda^j(\mathfrak{p}^-)^\mathrm{const}_q
\stackrel{\sim}{_\rightarrow}\Lambda^j(\mathfrak{p}^-)_q,\quad
f\otimes\omega\mapsto f\omega,
\end{equation}
\begin{equation}\label{morph2-2forms}
\Lambda^j(\mathfrak{p}^-)^\mathrm{const}_q\otimes
\mathbb{C}[\mathfrak{p}^-]_q\stackrel{\sim}{_\rightarrow}
\Lambda^j(\mathfrak{p}^-)_q,\quad\omega\otimes f\mapsto\omega\,f.
\end{equation}
\end{proposition}

{\bf Доказательство.} Достаточно воспользоваться предложением \ref{freedom}
и определением алгебры $\Lambda(\mathfrak{p}^-)_q$. \hfill $\square$

\bigskip
Перейдем к вычислению $\dim \Lambda^j(\mathfrak{p}^-)^\mathrm{const}_q$.

\begin{lemma}\label{dim_2-forms}
Размерность векторного пространства
$\Lambda^2(\mathfrak{p}^-)^\mathrm{const}_q$ не превосходит своего значения
в классическом случае $q=1$:
\begin{equation}\label{ineq_2-forms}
\dim\Lambda^2(\mathfrak{p}^-)^\mathrm{const}_q\le
\frac{\dim\mathfrak{p}^-(\dim\mathfrak{p}^--1)}{2}.
\end{equation}
\end{lemma}

{\bf Доказательство.} Рассмотрим линейный оператор
\begin{equation*}
\widetilde{R}:dz_i\otimes dz_j\to\sum\limits_{k,m=1}^{\dim\,\mathfrak{p}^-}
\check{R}_{ij}^{km}d z_k\otimes dz_m
\end{equation*}
в пространстве $\Lambda^1(\mathfrak{p}^-)_q^{\mathrm{const}}\otimes
\Lambda^1(\mathfrak{p}^-)_q^{\mathrm{const}}$. Имеет место естественный
изоморфизм
$$
\Lambda^2(\mathfrak{p}^-)_q^{\mathrm{const}}\cong
\left\{v\in\Lambda^1(\mathfrak{p}^-)_q^{\mathrm{const}}
\otimes\Lambda^1(\mathfrak{p}^-)_q^{\mathrm{const}}\left|\:
\widetilde{R}v=-v\right.\right\}.
$$
Используя замечание \ref{comm_mu_lambda} и равенство (\ref{Resh_Dr}),
нетрудно показать, что все собственные значения этого линейного оператора
вещественны и отличны от нуля. Рассуждая так же, как при доказательстве
леммы \ref{eigenvectors_in_L}, считаем функции $\check{R}_{ij}^{km}$
переменной $q$ аналитическими и, следовательно, непрерывными на $(0,1]$.
При всех $q\in(0,1]$ количество отрицательных собственных значений
оператора $\widetilde{R}$, рассматриваемых с учетом кратностей, равно
$\frac{\dim{p}^-(\dim{p}^--1)}2$. Следовательно, размерность собственного
подпространства, отвечающего собственному значению $-1$, не превосходит
$\frac{\dim{p}^-(\dim{p}^--1)}{2}$. \hfill $\square$

\bigskip

Неравенство, противоположное \eqref{ineq_2-forms}, будет получено, если
удастся доказать, что -1 является собственным значением оператора
$\widetilde{R}$ кратности не менее $\frac{\dim{p}^-(\dim{p}^--1)}2$.
Достаточно получить такую оценку для линейного оператора
$(\widetilde{R}^{*})^{-1}$. Его можно отождествить со сужением линейного
оператора
$\check{R}_{N(\mathfrak{q}^+,-\alpha_{l_0})N(\mathfrak{q}^+,-\alpha_{l_0})}$
на тензорное произведение старших однородных компонент
$N(\mathfrak{q}^+,-\alpha_{l_0})_\mathrm{highest}\otimes
N(\mathfrak{q}^+,-\alpha_{l_0})_\mathrm{highest}$.

\medskip

Обратимся к классическому случаю $q=1$. Следующий результат принадлежит
Костанту \cite[стр. 359-360]{Kostant_1961}, \cite[стр. 348]{Rocha}.

\begin{lemma}\label{W-orbit}
Рассмотрим $U\mathfrak{k}$-модуль $(\mathfrak{p}^-)^{\wedge r}$, где
$r=1,2,\ldots,\dim\mathfrak{p}^-$. Его изотипические компоненты являются
простыми $U\mathfrak{k}$-модулями, весовые подпространства одномерны, и
множество весов имеет вид
$$\{w\rho-\rho\;|\;w\in W^\mathbb{S}\;\&\;l(w)=r\},$$
где $\mathbb{S}=\{1,2,\ldots,l\}\setminus l_0$, $W^\mathbb{S}\,=\,
(^\mathbb{S}W)^{-1})$, а $^\mathbb{S}W$ было введено в \itemiiiе
\ref{IrrepK}.
\end{lemma}

\begin{lemma}\label{dim_2-forms_2}
Размерность векторного пространства
$\Lambda^2(\mathfrak{p}^-)^\mathrm{const}_q$ больше или равна своему
значению при $q=1$, то есть
\begin{equation}\label{eq_2-forms}
\dim\Lambda^2(\mathfrak{p}^-)^\mathrm{const}_q\ge
\frac{\dim\mathfrak{p}^-(\dim\mathfrak{p}^--1)}{2}.
\end{equation}
\end{lemma}

{\bf Доказательство.} При $q=1$ неравенство \eqref{eq_2-forms} очевидно.
Согласно лемме \ref{W-orbit}, из него следует, что старшие веса простых
$U\mathfrak{k}$-модулей, допускающих вложение в
$\mathfrak{p}^-\wedge\mathfrak{p}^-$, принадлежат $W\rho-\rho$, и сумма их
размерностей равна $\frac{\dim{p}^-(\dim{p}^--1)}{2}$.

При переходе от классического случая к квантовому не изменяются ни старшие
веса, ни размерности простых $U_q\mathfrak{k}$-подмодулей конечномерного
$U_q\mathfrak{k}$-модуля
$N(\mathfrak{q}^+,-\alpha_{l_0})_{\mathrm{highest}}\otimes
N(\mathfrak{q}^+,-\alpha_{l_0})_{\mathrm{highest}}$. Не меняются и кратности
в разложениях тензорных произведений конечномерных простых
\hbox{$U_q\mathfrak{k}$-модулей}.

Остается убедиться в том, что на образе любого морфизма
{$U_q\mathfrak{g}$-модулей}
$$
N(\mathfrak{q}^+,w\rho-\rho)\to N(\mathfrak{q}^+,-\alpha_{l_0})\otimes
N(\mathfrak{q}^+,-\alpha_{l_0}),\qquad w\in W^{\mathbb{S}}\,\&\,l(w)=2
$$
линейный оператор $\check{R}_{N(\mathfrak{q}^+,-\alpha_{l_0}),
N(\mathfrak{q}^+,-\alpha_{l_0})}$ равен минус единице. Достаточно доказать,
что он равен $\pm 1$, поскольку из \eqref{Resh_Dr} и непрерывности по
параметру $q$ следует, что мы не покидаем спектральное подпространство,
отвечающее отрицательной части спектра.

Воспользуемся равенством \eqref{Resh_Dr}. Из его доказательства,
приведенного в \cite[стр 39]{Drinf2}, видно, что на образе морфизма $
N(\mathfrak{q}^+,\nu)\rightarrow N(\mathfrak{q}^+,\lambda)\otimes
N(\mathfrak{q}^+,\mu)$ линейный оператор
$\check{R}_{N(\mathfrak{q}^+,\mu),N(\mathfrak{q}^+,\lambda)}
\check{R}_{N(\mathfrak{q}^+,\lambda),N(\mathfrak{q}^+,\mu)}$ скалярен и
равен
$$
q^{-(\mu,\mu+2\rho)-(\lambda,\lambda+2\rho)+(\nu,\nu+2\rho)}=
q^{-(\mu+\rho,\mu+\rho)-(\lambda+\rho,\lambda+\rho)+(\nu+\rho,\nu+\rho)+
(\rho,\rho)}.
$$
Подставляя $\lambda=\mu=-\alpha_{l_0}$, $\nu=w\rho-\rho$, получаем единицу,
поскольку веса $-\alpha_{l_0}+\rho$, $\nu+\rho$, $\rho$ принадлежат одной
$W$-орбите и, следовательно, имеют равные длины. \hfill $\square$

\bigskip Перейдем к дифференциальным формам более высокого порядка.

\begin{lemma}\label{dim_j-forms}
При всех $j\ge 3$
\begin{equation}\label{ineq_2-forms-j}
\dim\Lambda^j(\mathfrak{p}^-)^\mathrm{const}_q\le
\binom{\dim\mathfrak{p}^-}{j}.
\end{equation}
\end{lemma}

{\bf Доказательство.} Выберем базис весовых векторов
\begin{equation}\label{distinguished_basis}
\{dz_1,dz_2,\ldots,dz_{\dim\mathfrak{p}^-}\}\subset
\Lambda^1(\mathfrak{p}^-)_q^{\mathrm{const}},
\end{equation}
следуя заключительному замечанию \itemiiiа \ref{quadratic_first}. Введем
лексикографическое отношение порядка на множестве $\{dz_i\otimes
dz_k\}_{i,k=1,2,\ldots,\dim\mathfrak{p}^-}$. В этом базисе действие
универсальной \hbox{$R$-матрицы} алгебры Хопфа $U_q\mathfrak{k}$
описывается треугольной матрицей с положительными элементами на главной
диагонали, см. \eqref{Rmatrix}. Значит, ненулевой вклад в правую часть
равенства
\begin{equation}\label{subst_rules}
dz_i\otimes dz_j=
-\sum\limits_{k,m=1}^{\dim\,\mathfrak{p}^-}\check{R}_{ij}^{km}dz_k\otimes
dz_m,\qquad i\ge j
\end{equation}
вносят только слагаемые $dz_k\otimes dz_m$ с $k\le m$. Следовательно,
каждый элемент пространства $\Lambda^j(\mathfrak{p}^-)^\mathrm{const}_q$
принадлежит линейной оболочке мономов
$$
dz_{i_1}\wedge dz_{i_2}\wedge\cdots\wedge dz_{i_j},\qquad 1\le
i_1<i_2<\ldots<i_j\le\dim\mathfrak{p}^-.\eqno\square
$$

\begin{remark}\label{Grobner_const_forms}
Выделенный базис \eqref{distinguished_basis} можно использовать для
отождествления тензорной алгебры $T((\mathfrak{p}^-)^*)$ со свободной
некоммутативной алгеброй $\mathbb{C}\langle dz_1, dz_2,\ldots,dz_{\dim
\mathfrak{p}^-}\rangle$. Рассмотрим двусторонний идеал $I$ свободной
алгебры, порожденный множеством $G$ разностей левых и правых частей равенств
\eqref{subst_rules}. В заключительной части доказательства леммы
\ref{dim_j-forms} эти равенства играют роль правил подстановки: левая часть
равенства, где бы она ни встретилась, заменяется его правой частью.
\end{remark}
\medskip

Доказательство следующей леммы повторяет доказательство аналогичного
результата И.~Геккенберже и С.~Колба \cite{HeckenbergerKolb03}.

\begin{lemma}\label{dim_j-forms_2}
При всех $j\ge 3$
\begin{equation}\label{eq_j-forms}
\dim\Lambda^j(\mathfrak{p}^-)^\mathrm{const}_q\ge
\binom{\dim\mathfrak{p}^-}{j}.
\end{equation}
\end{lemma}

{\bf Доказательство.} Если неравенство \eqref{eq_j-forms} выполняется при
$j=3$, то, как следует из даймонд-леммы \ref{diamond_compos}, множество $G$
является базисом Гребнера двустороннего идеала $I$, см. замечание
\ref{Grobner_const_forms}, и отсюда следует неравенство
\eqref{dim_j-forms_2}.

Значит, при доказательстве предложения можно ограничиться частным случаем
$j=3$. Отождествим $(\mathfrak{p}^-)^*$ с
$\Lambda^1(\mathfrak{p}^-)_q^{\mathrm{const}}$ и введем обозначение
$(\mathfrak{p}^-)^{*\wedge 2}$ для подпространства
$\left\{v\in(\mathfrak{p}^-)^{*\otimes
2}\,\left|\,\widetilde{R}\,v=-v\right.\right\}$, ср. с \eqref{subst_rules}.
Как следует из лемм \ref{dim_2-forms}, \ref{dim_2-forms_2},
\begin{equation}\label{final_2_forms}
\dim((\mathfrak{p}^-)^{*\wedge
2})=\frac{\dim\mathfrak{p}^-(\dim\mathfrak{p}^--1)}2.
\end{equation}

Рассмотрим подпространства $L_1=(\mathfrak{p}^-)^{*\wedge
2}\otimes(\mathfrak{p}^-)^*$,
$L_2=(\mathfrak{p}^-)^*\otimes(\mathfrak{p}^-)^{*\wedge 2}$ векторного
пространства $(\mathfrak{p}^-)^{*\otimes 3}$ и комплекс линейных отображений
$$
0\rightarrow L_1\cap L_2\rightarrow L_1\newoplus
L_2\stackrel{j}{\rightarrow}(\mathfrak{p}^-)^{*\otimes
3}\rightarrow\mathbb{C}[\mathfrak{p}^-]_{q,3}\rightarrow 0,
$$
$$j:v_1\oplus v_2\mapsto v_1-v_2,\qquad v_j\in L_j
\subset(\mathfrak{p}^-)^{*\otimes 3},
$$
точный во всех членах, кроме, быть может, члена $(\mathfrak{p}^-)^{*\otimes
3}$. Вычисляя эйлерову характеристику этого комплекса \cite[стр.
91]{KostrikinManin}, \cite[стр. 119]{Lang}, приходим к следующему
неравенству:
$$
-\dim(L_1\cap L_2)+\dim(L_1\newoplus L_2)-\dim((\mathfrak{p}^-)^{*\otimes
3})+\dim\mathbb{C}[\mathfrak{p}^-]_{q,3}\le 0.
$$
Используя \eqref{final_2_forms}, \eqref{dim_j}, получаем
$$
\dim(L_1\cap L_2)\ge
n^2(n-1)-n^3+\frac{n(n+1)(n+2)}{6}=\frac{n(n-1)(n-2)}{6},
$$
где $n=\dim\mathfrak{p}^-=\dim(\mathfrak{p}^-)^*$. Из приведенного в примере
\ref{J_0} описания универсального дифференциального исчисления вытекает
неравенство
$$
\dim\Lambda^3(\mathfrak{p}^-)^\mathrm{const}_q\ge\dim(L_1\cap
L_2)\ge\binom{\dim\mathfrak{p}^-}{3}.\eqno\square
$$

\bigskip

Следующее утверждение вытекает из доказанных в этом \itemiiiе лемм.
\begin{proposition}\label{higher_order}
Однородные компоненты $\Lambda^j(\mathfrak{p}^-)_q$ дифференциальной
градуированной алгебры $\Lambda(\mathfrak{p}^-)_q$ отличны от нуля при
$j\le\dim\mathfrak{p}^-$. Каждая из них является как свободным левым, так и
свободным правым $\mathbb{C}[\mathfrak{p}^-]_q$-модулем ранга
$\binom{\dim\mathfrak{p}^-}{j}$.
\end{proposition}

В важном частном случае $\mathfrak{g}=\mathfrak{sl}_N$ полученный результат
хорошо известен \cite[стр. 296]{Maltsiniotis}.

\bigskip
\begin{remark}\label{math.QA}
В предположении о трансцендентности $q$ нетрудно доказать
 точность комплекса де
Рама и то, что получаемая из него по двойственности резольвента
тривиального $U_q\mathfrak{g}$-модуля изоморфна обобщенной резольвенте
Бернштейна-Гельфанда-Гельфанда, см.
 \cite{SSV_diff_calc}. Недавние результаты Геккенберже и Колба
указывают на то, что от предположения о трансцендентности $q$ можно
освободиться, см. \cite{HeckenbergerKolb0605, HeckenbergerKolb0611}.
\end{remark}

\hyphenation{диф-фе-рен-ци-аль-ны-ми анти-изо-мор-физ-мом}

\subsubsection{Дифференциальные формы с полиномиальными
коэффициентами.}\label{pol_diff}

В этом и в следующем \itemiiiе будут получены обобщения результатов
\itemiiiов \ref{pol-forms_sl_2} и \ref{finite-forms_sl_2} соответственно.

Начнем с построения ковариантных дифференциальных исчислений первого порядка
над алгеброй $\operatorname{Pol}(\mathfrak{p}^-)_q$. Напомним, что в
категории $U_q\mathfrak{g}$-модулей
$\operatorname{Pol}(\mathfrak{p}^-)_q=\mathbb{C}[\mathfrak{p}^-]_q \otimes
\mathbb{C}[\mathfrak{p}^+]_q$.

Рассмотрим $U_q\mathfrak{g}$-модуль
$\Omega^{1,0}(\mathfrak{p}^-)_q\stackrel{\operatorname{def}}{=}
\Lambda^1(\mathfrak{p}^-)_q\otimes\mathbb{C}[\mathfrak{p}^+]_q$. Используя
универсальную $R$-матрицу так же, как в \itemiiiе \ref{braid_algebra},
наделим $\Omega^{1,0}(\mathfrak{p}^-)_q$ структурой
$U_q\mathfrak{g}$-модульного
$\operatorname{Pol}(\mathfrak{p}^-)_q$-бимодуля:
$$
\xymatrix{\operatorname{Pol}(\mathfrak{p}^-)_q
\otimes\Omega^{1,0}(\mathfrak{p}^-)_q\ar[d]^\cong
\\ \mathbb{C}[\mathfrak{p}^-]_q\otimes(\mathbb{C}[\mathfrak{p}^+]_q\otimes
\Lambda^1(\mathfrak{p}^-)_q)\otimes\mathbb{C}[\mathfrak{p}^+]_q
\ar[d]^{\mathrm{id}\otimes\check{R}\otimes\mathrm{id}}
\\ \mathbb{C}[\mathfrak{p}^-]_q\otimes(\Lambda^1(\mathfrak{p}^-)_q\otimes
\mathbb{C}[\mathfrak{p}^+]_q)\otimes\mathbb{C}[\mathfrak{p}^+]_q
\ar[d]^\cong
\\ (\mathbb{C}[\mathfrak{p}^-]_q\otimes\Lambda^1(\mathfrak{p}^-)_q)\otimes
(\mathbb{C}[\mathfrak{p}^+]_q\otimes\mathbb{C}[\mathfrak{p}^+]_q)\ar[d]
\\ \Omega^{1,0}(\mathfrak{p}^-)_q}
$$

\medskip

Введем дифференциал
$$
\partial:\mathbb{C}[\mathfrak{p}^-]_q\otimes\mathbb{C}[\mathfrak{p}^+]_q\to
\Lambda^1(\mathfrak{p}^-)_q\otimes\mathbb{C}[\mathfrak{p}^+]_q,\qquad
\partial=d\otimes\operatorname{id}.
$$

  Используя результаты \itemiiiа \ref{FODC} и предложение \ref{braiding1},
  нетрудно доказать следующее утверждение.
\begin{proposition}\label{CFODC_Pol_1}
1. $\Omega^{1,0}(\mathfrak{p}^-)_q$ является как свободным левым
$\mathrm{Pol}(\mathfrak{p}^-)_q$-модулем, так и свободным правым
$\mathrm{Pol}(\mathfrak{p}^-)_q$-модулем.

2. $(\Omega^{1,0}(\mathfrak{p}^-)_q,\partial)$ -- ковариантное
дифференциальное исчисление первого порядка над
$\mathrm{Pol}(\mathfrak{p}^-)_q$.
\end{proposition}

Аналогично вводится ковариантное дифференциальное исчисление первого порядка
$(\Omega^{0,1}(\mathfrak{p}^-)_q,\overline{\partial})$ над
$\mathrm{Pol}(\mathfrak{p}^-)_q$. Опишем требуемые построения. Начнем с
алгебры $\mathbb{C}[\mathfrak{p}^+]_q$.
 Воспользуемся антилинейным
отображением
$$
*: \mathbb{C}[\mathfrak{p}^-]_q \rightarrow \mathbb{C}[\mathfrak{p}^+]_q,
\qquad *: f \mapsto f^*
$$
и сопряженным к нему антилинейным отображением
$$
*: N(\mathfrak{p}^+, 0) \rightarrow N(\mathfrak{p}^-, 0), \qquad *: \xi
v(\mathfrak{p}^+, 0) \mapsto (S^{-1}(\xi))^* v(\mathfrak{p}^-,0).
$$
 Отметим, что
$$
(F_{l_0} v(\mathfrak{p}^+,0))^* = (S^{-1}(F_{l_0}))^* v(\mathfrak{p}^-,0) =
(-K_{l_0} F_{l_0})^* v(\mathfrak{p}^-,0) = E_{l_0} v(\mathfrak{p}^-,0).
$$
Заменяя в построениях раздела \ref{FODC} морфизм $U_q\mathfrak{g}$-модулей
$$
 N(\mathfrak{p}^+, -\alpha_{l_0}) \rightarrow N(\mathfrak{p}^+,0), \quad
v(\mathfrak{p}^+, -\alpha_{l_0}) \mapsto F_{l_0} v(\mathfrak{p}^+,0)
$$
морфизмом $U_q\mathfrak{g}$-модулей
$$
 N(\mathfrak{p}^-, \alpha_{l_0}) \rightarrow N(\mathfrak{p}^-,0), \quad
v(\mathfrak{p}^-, \alpha_{l_0}) \mapsto E_{l_0} v(\mathfrak{p}^-,0),
$$
получаем ковариантное дифференциальное исчисление первого порядка
$(\Lambda^1(\mathfrak{p}^+)_q, d)$ над $U_q\mathfrak{g}$-модульной алгеброй
$\mathbb{C}[\mathfrak{p}^+]_q$. Так же, как предложение \ref{star_Pol},
доказывается

\begin{proposition}\label{star_rel}
Пусть $*: \Lambda^1(\mathfrak{p}^-)_q \rightarrow
\Lambda^1(\mathfrak{p}^+)_q$ -- антилинейное отображение, сопряженное к
антилинейному отображению
$$
*: N(\mathfrak{p}^-, \alpha_{l_0}) \rightarrow N(\mathfrak{p}^+,
-\alpha_{l_0}), \quad *: \xi v(\mathfrak{p}^-, \alpha_{l_0}) \mapsto
(S^{-1}(\xi))^* v(\mathfrak{p}^+, -\alpha_{l_0}).
$$
Тогда при всех $f, f_1, f_2, f_3 \in \mathbb{C}[\mathfrak{p}^-]_q$
$$
(f_1 \; d f_2 \; f_3)^* = f_3^* (d f_2)^* f_1^*, \qquad (d f)^* = d f^*.
$$
\end{proposition}

Если $\{z_i\}$ -- базис векторного пространства
$\mathbb{C}[\mathfrak{p}^-]_{q,1}$, то, как следует из равенства
\eqref{rel_z_dz} и из предложения \ref{star_rel},
\begin{equation}\label{dzs_zs}
  d z_j^*\; z_i^* = \sum\limits_{k,m=1}^{\dim\,\mathfrak{p}^-}
\check{R}_{ij}^{km}\, z_m^* \; d z_k^*.
\end{equation}

Рассмотрим $U_q\mathfrak{g}$-модуль $
\Omega^{0,1}(\mathfrak{p}^-)_q\stackrel{\operatorname{def}}{=}
\mathbb{C}[\mathfrak{p}^-]_q \otimes \Lambda^1(\mathfrak{p}^+)_q$. Как и
прежде, с помощью
 универсальной $R$-матрицы  наделим $\Omega^{0,1}(\mathfrak{p}^-)_q$ структурой
$U_q\mathfrak{g}$-модульного
$\operatorname{Pol}(\mathfrak{p}^-)_q$-бимодуля. Введем дифференциал
$$ \overline{\partial}: \mathbb{C}[\mathfrak{p}^-]_q \otimes
\mathbb{C}[\mathfrak{p}^+]_q \rightarrow \mathbb{C}[\mathfrak{p}^+]_q
\otimes \Lambda^1(\mathfrak{p}^-)_q, \qquad
\overline{\partial}=\operatorname{id} \otimes d.
 $$
Следующее утверждение доказывается так же, как предложение
\ref{CFODC_Pol_1}.
\begin{proposition}\label{CFODC_Pol_2}
1. $\Omega^{0,1}(\mathfrak{p}^-)_q$ является как свободным левым
$\mathrm{Pol}(\mathfrak{p}^-)_q$-модулем, так и свободным правым
$\mathrm{Pol}(\mathfrak{p}^-)_q$-модулем.

2. $(\Omega^{0,1}(\mathfrak{p}^-)_q,\overline{\partial})$ -- ковариантное
дифференциальное исчисление первого порядка над
$\mathrm{Pol}(\mathfrak{p}^-)_q$.
\end{proposition}

  Рассмотрим $U_q\mathfrak{g}$-модульный
$\mathrm{Pol}(\mathfrak{p}^-)_q$-бимодуль
    $$\Omega^1(\mathfrak{p}^-)_q=
    \Omega^{1,0}(\mathfrak{p}^-)_q \newoplus
    \Omega^{0,1}(\mathfrak{p}^-)_q$$
и дифференциал $d=\partial + \overline{\partial}$.

\begin{corollary} 1. $\Omega^1(\mathfrak{p}^-)_q$ является как свободным левым
$\mathrm{Pol}(\mathfrak{p}^-)_q$-модулем, так и свободным правым
$\mathrm{Pol}(\mathfrak{p}^-)_q$-модулем.

2. $(\Omega^1(\mathfrak{p}^-)_q,d)$ -- ковариантное дифференциальное
исчисление первого порядка над $\mathrm{Pol}(\mathfrak{p}^-)_q$.
\end{corollary}

 Нетрудно показать, что
$y\,\omega= \omega\, y$ при всех $\omega \in \Lambda^1(\mathfrak{p}^-)_{q,1}
\oplus \Lambda^1(\mathfrak{p}^+)_{q,-1}$, см. доказательство приведенной
ниже леммы \ref{f0_omega}.
 C помощью локализации
$\mathrm{Pol}(\mathfrak{p}^-)_q$-бимодуля
$\Omega^1(\mathfrak{p}^-)_q$ по мультипликативному подмножеству
$y^{\mathbb{Z}_+}$ и рассуждений, описанных в \itemiiiе
\ref{w0_G_x}, получаем ковариантное дифференциальное исчисление
первого порядка над $U_q\mathfrak{g}$-модульной алгеброй
$\mathrm{Pol}(\mathfrak{p}^-)_{q,y}$.

   Перейдем от дифференциальных исчислений первого порядка к полным
дифференциальным исчислениям над алгеброй $\mathrm{Pol}(\mathfrak{p}^-)_q$.
Рассмотрим дифференциальное исчисление первого порядка
$(\Lambda^1(\mathfrak{p}^+)_q,d)$ над алгеброй
$\mathbb{C}[\mathfrak{p}^+]_q$ и отвечающее ему универсальное
дифференциальное исчисление $(\Lambda(\mathfrak{p}^+)_q,d)$. Из результатов
\itemiiiа \ref{diff_univ} следует, что $\Lambda^j(\mathfrak{p}^+)_q$
является как свободным левым, так и свободным правым
$\mathbb{C}[\mathfrak{p}^+]_q$-модулем ранга $\binom{\dim
\mathfrak{p}^+}{j}$.

Пусть $(\Omega^{*,0}(\mathfrak{p}^-)_q,\partial)$,
$(\Omega^{0,*}(\mathfrak{p}^-)_q,\overline{\partial})$ и
$(\Omega(\mathfrak{p}^-)_q,d)$ -- универсальные дифференциальные
исчисления, отвечающие дифференциальным исчислениям первого порядка
$(\Omega^{1,0}(\mathfrak{p}^-)_q,\partial)$,
$(\Omega^{0,1}(\mathfrak{p}^-)_q,\overline{\partial})$ и
$(\Omega^1(\mathfrak{p}^-)_q,d)$ соответственно. В категории
$U_q\mathfrak{g}$-модулей
$$
\Omega^{*,0}(\mathfrak{p}^-)_q\cong\Lambda(\mathfrak{p}^-)_q\otimes
\mathbb{C}[\mathfrak{p}^+]_q,\qquad\partial=d\otimes\operatorname{id},
$$
$$
\Omega^{0,*}(\mathfrak{p}^-)_q\cong\mathbb{C}[\mathfrak{p}^-]_q\otimes
\Lambda(\mathfrak{p}^+)_q,\qquad\overline{\partial}=
\operatorname{id}\otimes d,
$$
$$
\Omega(\mathfrak{p}^-)_q\cong\Omega^{*,0}(\mathfrak{p}^-)_q
\otimes_{\operatorname{Pol}(\mathfrak{p}^-)_q}
\Omega^{0,*}(\mathfrak{p}^-)_q\cong\Lambda(\mathfrak{p}^-)_q\otimes
\Lambda(\mathfrak{p}^+)_q,
$$
причем в последнем случае дифференциал равен
$d\otimes\operatorname{id}\oplus\operatorname{id}\otimes d$.

 Сохраним обозначения $\partial$, $\overline{\partial}$ для операторов $d
\otimes \operatorname{id}$, $\operatorname{id} \otimes d$ в
$\Omega(\mathfrak{p}^-)_q$. Алгебра $\Omega(\mathfrak{p}^-)_q$ является
градуированной: $\Omega(\mathfrak{p}^-)_q=\newoplus\limits_{k=1}^{\dim
\mathfrak{p}^-}\Omega^k(\mathfrak{p}^-)_q$, и
\begin{equation*}\label{bi_gr}
\Omega^k(\mathfrak{p}^-)_q=\newoplus_{i+j=k}
   \Omega^{i,j}(\mathfrak{p}^-)_q, \qquad
   \Omega^{i,j}(\mathfrak{p}^-)_q \cong \Lambda^i(\mathfrak{p}^-)_q \otimes
\Lambda^j(\mathfrak{p}^+)_q,
\end{equation*}
\begin{equation*}\label{bi_degree}
\partial: \Omega^{i,j}(\mathfrak{p}^-)_q \rightarrow
\Omega^{i+1,j}(\mathfrak{p}^-)_q, \qquad \overline{\partial}:
\Omega^{i,j}(\mathfrak{p}^-)_q \rightarrow \Omega^{i,j+1}(\mathfrak{p}^-)_q.
\end{equation*}
Равенство $d^2=0$ означает, что имеют место следующие равенства:
\begin{equation}\label{commute}
  \partial^2 = \overline{\partial}^{\,2} = 0, \quad \partial\, \overline{\partial}
= -\overline{\partial} \,\partial.
\end{equation}

Рассуждая так же, как в \itemiiiах \ref{braid_algebra}, \ref{Pol}, нетрудно
наделить $\Omega(\mathfrak{p}^-)_q$ структурой
$(U_q\mathfrak{g},*)$-модульной алгебры с помощью инволюции
$$
 *: \omega \otimes 1 \mapsto 1 \otimes \omega^*, \qquad
\omega \in \Lambda^1(\mathfrak{p}^-)_q.
$$
Из определений вытекает равенство
\begin{equation}\label{d_star}
  \overline{\partial}\, \omega =
  (-1)^k\,( \partial \,  \omega^*)^*, \qquad \omega \in \Omega^k(\mathfrak{p}^-)_q.
 \end{equation}

\medskip Приведем некоторые коммутационные соотношения. Пусть $\{z_i\}$ --
базис векторного пространства
 $\mathbb{C}[\mathfrak{p}^-]_{q,1}$. Используя универсальную
 \hbox{$R$-матрицу}
алгебры Хопфа $U_q\mathfrak{g}$, введем оператор перестановки тензорных
сомножителей
$$
\mathbb{C}[\mathfrak{p}^+]_{q,-1} \otimes \Lambda^1(\mathfrak{p}^-)_{q,1}
\rightarrow \Lambda^1(\mathfrak{p}^-)_{q,1} \otimes
\mathbb{C}[\mathfrak{p}^+]_{q,-1},
$$
$$
z_i^* \otimes d z_j \mapsto \sum\limits_{k,m=1}^{\dim\,\mathfrak{p}^-}
\widehat{R}_{ij}^{km} d z_k \otimes z_m^*.
$$
Из определений следует равенство
\begin{equation}\label{rel_zs_dz}
  z_i^* \, d z_j =  \sum\limits_{k,m=1}^{\dim\,\mathfrak{p}^-}
\widehat{R}_{ij}^{km} d z_k \, z_m^*.
\end{equation}
В частности,
\begin{equation}\label{zl_dzsl}
  z_{\operatorname{low}}^* \, d z_{\operatorname{low}} = q_{l_0}^2 d z_{\operatorname{low}} \, z_{\operatorname{low}}^*.
\end{equation}
Применяя инволюцию $*$, получаем
\begin{equation}\label{dzs_z}
  d z_j^* \, z_i =  \sum\limits_{k,m=1}^{\dim\,\mathfrak{p}^-}
\widehat{R}_{ij}^{km} z_m \, d z_k^*,
\end{equation}
\begin{equation}\label{dzls_zl}
  d z_{\operatorname{low}}^* \,
  z_{\operatorname{low}} =
  q_{l_0}^2 z_{\operatorname{low}} \, dz_{\operatorname{low}}^*.
\end{equation}

Из определения умножения в $\operatorname{Pol}(\mathfrak{p}^-)_q$ и из
\eqref{Rmatrix} следует равенство
\begin{equation}\label{rel_zs_z}
z_i^*\,z_j=\sum\limits_{k,m=1}^{\dim\,\mathfrak{p}^-}\widehat{R}_{ij}^{km}
z_k\,z_m^*+(z_j,z_i),
\end{equation}
где $(\cdot,\cdot)$ -- стандартная эрмитова форма в
$\mathbb{C}[\mathfrak{p}^-]_q$, см. \itemiii\ \ref{Pol}.

\begin{note}\label{FF}
Пусть $(z_i,z_j)=\delta_{i,j}$. Коммутационные соотношения \eqref{rel_zs_z}
являются \hbox{$q$-аналогами} канонических коммутационных соотношений и
определяют виковскую алгебру \cite{Jorg, Jorg_Pacific}. Кроме того, в
формальном пределе $q\to 1$ элементы
$\{z_i\}_{i=1,2,\cdots,\dim\mathfrak{p}^-}$ попарно коммутируют. Значит,
алгебра $\operatorname{Pol}(\mathfrak{p}^-)_q$ является \hbox{$q$-аналогом}
алгебры Вейля дифференциальных операторов с полиномиальными коэффициентами
на $\mathfrak{p}^-$, а скалярное произведение $(\cdot,\cdot)$ --
$q$-аналогом скалярного произведения Фишера-Фока в векторном пространстве
голоморфных полиномов на $\mathfrak{p}^-$ \cite[стр.~24]{Arazy}.
\end{note}

\subsubsection{Дифференциальные формы с обобщенными
коэффициентами.}\label{finite_diff}

Алгебры $\Lambda(\mathfrak{p}^-)_q$, $\Lambda(\mathfrak{p}^+)_q$ являются
весовыми $U_q\mathfrak{g}$-модулями и, следовательно, градуированными
векторными пространствами
$$
\Lambda(\mathfrak{p}^-)_q = \newoplus\limits_{i=0}^\infty
\Lambda(\mathfrak{p}^-)_{q,i}, \qquad \Lambda(\mathfrak{p}^+)_q =
\newoplus\limits_{j=0}^\infty \Lambda(\mathfrak{p}^+)_{q,-j}.
$$
Значит, введенное в предыдущем \itemiiiе векторное пространство
$\Omega(\mathfrak{p}^-)_q$ наделено биградуировкой
\begin{equation}\label{summ}
 \Omega(\mathfrak{p}^-)_q =
\newoplus\limits_{(i,j) \in \mathbb{Z}_+^2}
\Omega(\mathfrak{p}^-)_{q,i,-j}.
\end{equation}
Как следует из результатов \itemiiiа \ref{pol_diff},
$\dim(\Omega(\mathfrak{p}^-)_{q,i,-j}) < \infty$ при всех $(i,j) \in
\mathbb{Z}_+^2$.

Напомним, что топологическое векторное пространство
$\mathscr{D}(\mathbb{D})_q'$ обобщенных функций в квантовой ограниченной
симметрической области было получено в \itemiiiе \ref{gen} переходом от
конечных сумм (\ref{series_f}) к формальным рядам. Поступим аналогично.
Пусть
\begin{equation}\label{series_omega}
\Omega(\mathbb{D})_q'=
\prod_{(i,j)\in\mathbb{Z}_+^2}\Omega(\mathfrak{p}^-)_{q,i,-j}
\end{equation}
-- векторное пространство формальных рядов с членами из
$\Omega(\mathfrak{p}^-)_{q,i,-j}$ и топологией покомпонентной сходимости.
Очевидно, $\Omega(\mathfrak{p}^-)_q$ -- плотное линейное подмногообразие
$\Omega(\mathbb{D})_q'$. Другими словами, $\Omega(\mathbb{D})_q'$ является
пополнением векторного пространства $\Omega(\mathfrak{p}^-)_q$, наделенного
топологией покомпонентной сходимости. Разумеется, $\Omega(\mathbb{D})_q'=
\newoplus\limits_{i,j=0}^{\dim\mathfrak{p}^-}\Omega^{i,j}(\mathbb{D})_q'$,
где $\Omega^{i,j}(\mathbb{D})_q'$ -- замыкание линейного многообразия
$\Omega^{i,j}(\mathfrak{p}^-)_q\subset\Omega(\mathbb{D})_q'$. Элементы
$\omega\in\Omega(\mathbb{D})_q'$ будем называть дифференциальными формами.
\IND{дифференциальная форма ! с обобщенными коэффициентами}

Линейные операторы $\partial$, $\overline{\partial}$ и антилинейный оператор
$*$ допускают продолжение по непрерывности с $\Omega(\mathfrak{p}^-)_q$ на
$\Omega(\mathbb{D})_q'$. При этом остаются справедливыми равенства
(\ref{commute}), (\ref{d_star}). Алгебра $\Omega(\mathfrak{p}^-)_q$ является
$U_q\mathfrak{g}$-модульным $\Omega(\mathfrak{p}^-)_q$-бимодулем. Структура
$U_q\mathfrak{g}$-модульного $\Omega(\mathfrak{p}^-)_q$-бимодуля
распространяется с $\Omega(\mathfrak{p}^-)_q$ на $\Omega(\mathbb{D})_q'$ по
непрерывности. При этом равенство (\ref{Leibnitz_form}) остается
справедливым, если хотя бы одна из дифференциальных форм $\omega', \omega''$
принадлежит $\Omega(\mathfrak{p}^-)_q$.

Дифференциальную форму $\omega\in\Omega(\mathbb{D})_q'$ назовем финитной,
\IND{дифференциальная форма ! финитная} если при некотором
$N(\omega)\in\mathbb{N}$
$$
\newoplus\limits_{j=N(\omega)}^\infty\mathbb{C}[\mathfrak{p}^+]_{q,-j}\cdot
\omega=0,\qquad\omega\cdot\newoplus\limits_{j=N(\omega)}^\infty
\mathbb{C}[\mathfrak{p}^-]_{q,j}=0.
$$
(ср. с предложением \ref{image}). Подпространство финитных дифференциальных
форм ${\omega\in\Omega(\mathbb{D})_q'}$ обозначим $\Omega(\mathbb{D})_q$.

Очевидно,
 $ \Omega(\mathbb{D})_q =
\newoplus\limits_{i,j=0}^{\dim\mathfrak{p}^-} \Omega^{i,j}(\mathbb{D})_q
$, где $\Omega^{i,j}(\mathbb{D})_q = \Omega^{i,j}(\mathbb{D})_q' \cap
\Omega(\mathbb{D})_q$.

\begin{proposition}\label{inv_subspaces}
$ \partial: \Omega(\mathbb{D})_q \rightarrow \Omega(\mathbb{D})_q$,
$\overline{\partial}: \Omega(\mathbb{D})_q \rightarrow
\Omega(\mathbb{D})_q$.
\end{proposition}

{\bf Доказательство.} Линейные операторы
$$
\partial: \Lambda(\mathfrak{p}^-)_q \rightarrow \Lambda(\mathfrak{p}^-)_q,
\qquad \overline{\partial}: \Lambda(\mathfrak{p}^+)_q \rightarrow
\Lambda(\mathfrak{p}^+)_q
$$
являются морфизмами $U_q\mathfrak{g}$-модулей и, следовательно, сохраняют
степень однородности. Значит,
$$
\partial \mathbb{C}[\mathfrak{p}^-]_{q,j} \subset
\mathbb{C}[\mathfrak{p}^-]_{q,j-1} \cdot \Lambda^1(\mathfrak{p}^-)_{q,1},
\quad \overline{\partial} \mathbb{C}[\mathfrak{p}^+]_{q,-j} \subset
\Lambda^1(\mathfrak{p}^+)_{q,-1}\cdot \mathbb{C}[\mathfrak{p}^+]_{q,-(j-1)}
$$
и равенства
$$
\mathbb{C}[\mathfrak{p}^+]_{q,-j} \cdot \omega = \omega \cdot
\mathbb{C}[\mathfrak{p}^-]_{q,j} = 0
$$
влекут
$$
\mathbb{C}[\mathfrak{p}^+]_{q,-(j+1)} \cdot \partial \omega = \partial
\omega \cdot \mathbb{C}[\mathfrak{p}^-]_{q,j+1} = 0,
$$
$$
\mathbb{C}[\mathfrak{p}^+]_{q,-(j+1)} \cdot \overline{\partial} \omega =
\overline{\partial} \omega \cdot \mathbb{C}[\mathfrak{p}^-]_{q,j+1} = 0.
\eqno \square$$

\bigskip Покажем, что структура $U_q\mathfrak{g}$-модульного
$\mathrm{Pol}(\mathfrak{p}^-)_q$-бимодуля в $\Omega(\mathbb{D})_q$ и
соображения непрерывности позволяют наделить $\Omega(\mathbb{D})_q$
структурой \hbox{$U_q\mathfrak{g}$-модульного}
$\mathscr{D}(\mathbb{D})_q$-бимодуля.

Ограничимся на время обобщенными функциями конечного типа, \IND{обобщенная
функция ! конечного типа} то есть формальными рядами вида
$$
f=\sum\limits_{\{i,j|\;|i-j|\le N(f)\}}\;f_{ij}\in
\mathbb{C}[\mathfrak{p}^-]_{q,i}\cdot\mathbb{C}[\mathfrak{p}^+]_{q,-j}.
$$
Будем говорить, что последовательность $f_n$, $n\in\mathbb{N}$ обобщенных
функций конечного типа сходится к нулю, если она сходится к нулю в
$\mathscr{D}(\mathbb{D})_q'$ и числа $N(f_n)$ можно выбрать так, что
$\sup\limits_n\,N(f_n)<\infty$. Описанную в предыдущих разделах структуру
$U_q\mathfrak{g}$-модульной алгебры в $\mathrm{Fun}(\mathbb{D})_q$ можно
получить, исходя из $\mathrm{Pol}(\mathbb{D})_q$, продолжением по
непрерывности. (Действительно, если $f_n\to f$, $g_n\to g$, то $f_ng_n\to f
g$.) Точно так же вводятся элементы
$\omega\in\mathscr{T}^\mathrm{red}(\mathbb{D})_q'$ конечного типа,
сходимость $\omega_n\to\omega$ и продолжение по непрерывности структуры
$U_q\mathfrak{g}$-модульной алгебры с $\Omega(\mathfrak{p}^-)_q$ на
$\Omega(\mathfrak{p}^-)_q\oplus\Omega(\mathbb{D})_q$. При этом сохраняется
равенство (\ref{Leibnitz_form}).

В дальнейшем будет показано (см. замечание \ref{df0_dbarf0}), что
 $\partial \mathscr{D}(\mathbb{D})_q$ порождает
$\mathscr{D}(\mathbb{D})_q$-бимодуль $\Omega^{1,0}(\mathbb{D})_q$, а
$\overline{\partial} \mathscr{D}(\mathbb{D})_q$ порождает
$\mathscr{D}(\mathbb{D})_q$-бимодуль $\Omega^{0,1}(\mathbb{D})_q$. Значит,
имеет место

\begin{proposition}\label{d_dbar}
Пары $(\Omega^{1,0}(\mathbb{D})_q, \partial)$, $(\Omega^{0,1}(\mathbb{D})_q,
\overline{\partial})$ являются ковариантными дифференциальными исчислениями
первого порядка над \hbox{$U_q\mathfrak{g}$-модульной} алгеброй
$\mathscr{D}(\mathbb{D})_q$.
\end{proposition}

\medskip
Найдем явный вид элемента $\overline{\partial} f_0$, где $f_0$ --
стандартная образующая $\mathrm{Pol}(\mathfrak{p}^-)_q$-бимодуля
$\mathscr{D}(\mathbb{D})_q$ (см. \itemiii\ \ref{finite}).

\begin{lemma}\label{f0_omega}
Для любого $\omega \in \Lambda^1(\mathfrak{p}^-)_{q,1}$
$$
f_0 \cdot \omega= \omega \cdot f_0.
$$
\end{lemma}

{\bf Доказательство.}
 Пусть
$$
\widetilde{R}: \mathscr{D}(\mathbb{D})_q' \otimes
\Lambda^1(\mathfrak{p}^-)_q \rightarrow \Lambda^1(\mathfrak{p}^-)_q \otimes
\mathscr{D}(\mathbb{D})_q'
$$
морфизм $U_q\mathfrak{g}$-модулей, определяемый действием универсальной
$R$-матрицы алгебры Хопфа $U_q\mathfrak{g}$. Тогда $\widetilde{R}: f_0
\otimes \omega \mapsto \omega \otimes f_0$, как следует из явной формулы
(\ref{Rmatrix}), равенства $ F_{l_0} \omega = 0$ и из того, что
$$E_jf_0=0,\quad j\neq l_0,\qquad H_if_0=0, \quad i=1,2,\ldots,l,$$
в силу $U_q\mathfrak{k}$-инвариантности элемента $f_0$.

Остается отметить, что структуру $\Lambda^1(\mathfrak{p}^-)_q$-бимодуля в
$\Omega(\mathbb{D})_q'$, как и в $\Omega(\mathfrak{p}^-)_q$, можно описать с
помощью оператора $\widetilde{R}$. \hfill $\square$

\medskip Применяя инволюцию $*$, получаем равенство
$\omega \cdot f_0 = f_0 \cdot \omega$ для всех $\omega \in
\Lambda^1(\mathfrak{p}^+)_{q,-1}$. Будем использовать стандартные
обозначения $$ dz=\partial z,\quad dz^*=\overline{\partial} z^*, \qquad z\in
\mathbb{C}[\mathfrak{p}^-]_{q,1}$$ и эрмитову форму в
$\mathbb{C}[\mathfrak{p}^-]_q$, введенную в \itemiiiе \ref{Pol}.

\begin{proposition}\label{df0}
Пусть $\{z_j\}$ -- ортогональный базис эрмитова векторного пространства
$\mathbb{C}[\mathfrak{p}^-]_{q,1}$. Тогда
\begin{equation}\label{dbar_f0}
  \overline{\partial} f_0 = -\sum\limits_{j=1}^{\dim\, \mathfrak{p}^-}
\frac{1}{(z_j, z_j)} z_j \cdot f_0 \cdot d z_j^*.
\end{equation}
\end{proposition}

{\bf Доказательство.}
 Получим систему уравнений, однозначно определяющую дифференциальную форму
$\omega = \overline{\partial} f_0$, после чего нам останется проверить, что
правая часть равенства (\ref{dbar_f0}) является решением этой системы
уравнений. Пусть $$\Psi = \sum\limits_{i=1}^{\dim\mathfrak{p}^-} \frac{z_i
z_i^*}{(z_i, z_i)}.$$ Тогда $\overline{\partial} \psi = \sum\limits_i (z_i,
z_i)^{-1} z_i d z_i^*$. Применение оператора $\overline{\partial}$ к обеим
частям равенства $\Psi f_0 = 0$ дает
$$
\sum\limits_i (z_i, z_i)^{-1} \,z_i \, d z_i^* \, f_0 + \sum\limits_i (z_i,
z_i)^{-1}\, z_i\, z_i^*\, \overline{\partial} f_0 = 0.
$$
Переставляя $d z_i^*$ с $f_0$, получаем первое уравнение для
$\omega=\overline{\partial} f_0$:
\begin{equation}\label{eq_omega_1}
  \Psi\, \omega = - \sum\limits_i (z_i, z_i)^{-1} z_i f_0 d z_i^*.
\end{equation}

Второе уравнение вытекает из того, что оператор $\overline{\partial}$
коммутирует с оператором $H_\mathbb{S}$, продолженным по непрерывности на
$\mathscr{T}^{0,1}(\mathbb{D})_q'$:
\begin{equation}\label{eq_omega_2}
  H_\mathbb{S} \, \omega = 0.
\end{equation}

 Дифференциальная форма
 $-\sum\limits_i (z_i, z_i)^{-1} z_i f_0 d z_i^* $ является решением
системы уравнений \eqref{eq_omega_1}, \eqref{eq_omega_2} в
$\Omega^{0,1}(\mathbb{D})_q'$. Остается доказать единственность решения этой
системы уравнений, то есть доказать, что соответствующая однородная система
уравнений
\begin{equation}\label{sys_omega}
  \Psi \omega = H_\mathbb{S} \, \omega = 0, \qquad
  \omega \in \Omega^{0,1}(\mathbb{D})_q'
\end{equation}
имеет лишь нулевое решение. Согласно следствию \ref{embed_finite}
$$
\Omega^{0,1}(\mathbb{D})_q' \cong \prod_{_{j,k}}
\mathbb{C}[\mathfrak{p}^-]_{q,j} \; f_0 \; \Lambda^1(\mathfrak{p}^+)_{q,-k}.
$$
Если существует ненулевое решение системы уравнений (\ref{sys_omega}), то
существует ненулевой элемент
\begin{equation}\label{j_neq_0}
  f \in \mathbb{C}[\mathfrak{p}^-]_{q,j} f_0, \qquad j \neq 0,
\end{equation}
для которого $\Psi f = 0$ и, следовательно, $f^* \Psi f = 0$. Это
противоречит положительности стандартной эрмитовой формы в
$\mathbb{C}[\mathfrak{p}^-]_q$ (см. предложение \ref{positive_herm}). \hfill
$\square$

\begin{corollary}\label{(1,0)-form}
$
\partial f_0 = -\sum\limits_{j=1}^{\dim\, \mathfrak{p}^-} \frac{1}{(z_j,
z_j)}{\; d z_j} \cdot f_0 \cdot z_j^*.
$
\end{corollary}

\begin{note}\label{df0_dbarf0} Используя \eqref{rel_zs_z}, получаем
$$
\mathscr{D}(\mathbb{D})_q \cdot \partial f_0 \cdot \mathscr{D}(\mathbb{D})_q
= \Omega^{1,0}(\mathbb{D})_q, \quad \mathscr{D}(\mathbb{D})_q \cdot
\overline{\partial} f_0 \cdot \mathscr{D}(\mathbb{D})_q =
\Omega^{0,1}(\mathbb{D})_q.
$$
\end{note}

\bigskip Полученные результаты позволяют описать универсальные дифференциальные
исчисления над алгеброй \hbox{$\operatorname{Pol}(\mathfrak{p}^-)_q \oplus
\mathscr{D}(\mathbb{D})_q$}, отвечающие ковариантным дифференциальным
исчислениям первого порядка
$$ (\Omega^{1,0}(\mathfrak{p}^-)_q \oplus \Omega^{1,0}(\mathbb{D})_q, \partial),\quad
(\Omega^{0,1}(\mathfrak{p}^-)_q \oplus \Omega^{0,1}(\mathbb{D})_q,
\overline{\partial}),\quad (\Omega^1(\mathfrak{p}^-)_q \oplus
\Omega^1(\mathbb{D})_q, d).$$


\subsubsection{Инвариантные дифференциальные операторы
на $K\backslash G$.}\label{affine_diff}

Рассмотрим эрмитово-симметрическую пару $(\mathfrak{g},\mathfrak{k})$ и
$U_q\mathfrak{g}$-модульную алгебру
$$
 \mathbb{C}[K\backslash G]_q\,=\,
 \{f \in \mathbb{C}[G]_q\quad|\quad L_{\rm reg}(\xi)f =
 \varepsilon(\xi)f, \quad \xi \in U_q\mathfrak{k}\}.
$$
Как отмечалось в \itemiiiах \ref{IrrepK}, \ref{reps_pol}, она является
квантовым аналогом алгебры регулярных функций на аффинном алгебраическом
многообразии $K\backslash G$.

Введем в рассмотрение алгебру линейных операторов в $\mathbb{C}[K\backslash
G]_q$, являющуюся квантовым аналогом алгебры $G$-инвариантных
дифференциальных операторов на $K\backslash G$. \IND{квантовый аналог !
алгебры $G$-инвариантных дифференциальных операторов на $K\backslash G$} В
классическом случае $q=1$ используемый нами подход описан в
\cite{EastwoodRice, KoranyiReimann}.

Следуя Дринфельду \cite{DrinfEng}, будем называть алгебру
$(U_q\mathfrak{g})^*$ всех линейных функционалов на $U_q\mathfrak{g}$
алгеброй функций на формальной квантовой группе $G$. \IND{алгебра функций !
на формальной квантовой группе} Это $U_q\mathfrak{g}$-модульная алгебра:
$$
R_{\mathrm{reg}}(\xi)l\,(\eta)\;=\;l(\eta\xi),\qquad\xi,\eta\in
U_q\mathfrak{g},\;l\in(U_q\mathfrak{g})^*,
$$
и ее $U_q\mathfrak{g}$-модульную подалгебру
\begin{equation*}\label{pairing_e}
{\rm Hom}_{_{U_q\mathfrak{k}}}(U_q\mathfrak{g},\mathbb{C})\cong
(\mathbb{C} \otimes_{_{U_q\mathfrak{k}}}U_q\mathfrak{g})^*
\end{equation*}
мы будем называть алгеброй функций на формальном квантовом однородном
пространстве $K\backslash G$. \IND{алгебра функций ! на формальном
квантовом однородном пространстве}

По определению, $\mathbb{C}[G]_q$ является наибольшим весовым локально
конечномерным подмодулем $U_q\mathfrak{g}$-модуля $(U_q\mathfrak{g})^*$.
Значит, $\mathbb{C}[K\backslash G]_q$ -- наибольший весовой локально
конечномерный подмодуль $U_q\mathfrak{g}$-модуля ${\rm
Hom}_{_{U_q\mathfrak{k}}}(U_q\mathfrak{g},\mathbb{C})$.

Эндоморфизмы правого $U_q\mathfrak{g}$-модуля $U_q\mathfrak{g}$ нумеруются
элементами $a\in U_q\mathfrak{g}$ и имеют вид $\eta\mapsto a\eta$.
Сопряженные линейные операторы
$$
L_{\rm reg}(a):(U_q\mathfrak{g})^*\to(U_q\mathfrak{g})^*,\qquad L_{\rm
reg}(a):f(\eta)\mapsto f(a\eta),
$$
играют роль инвариантных дифференциальных операторов на формальной
квантовой группе $G$. \IND{инвариантный дифференциальный оператор ! на
формальной квантовой группе}

Аналогично, эндоморфизмы правого $U_q\mathfrak{g}$-модуля $\mathbb{C}
\otimes_{_{U_q\mathfrak{k}}} U_q\mathfrak{g}$ нумеруются его
$U_q\mathfrak{k}$ инвариантными элементами
$$
{\rm End}_{_{U_q\mathfrak{g}}}(\mathbb{C}\otimes_{_{U_q\mathfrak{k}}}
U_q\mathfrak{g})\cong{\rm
Hom}_{_{U_q\mathfrak{k}}}(\mathbb{C},\mathbb{C}\otimes_{_{U_q\mathfrak{k}}}
U_q\mathfrak{g})\cong (\mathbb{C}\otimes_{_{U_q\mathfrak{k}}}
U_q\mathfrak{g}))^{U_q\mathfrak{k}}
$$
и имеют вид $1\otimes\eta\mapsto(1\otimes a)(1\otimes\eta)$, где $1\otimes
a\in(\mathbb{C}\otimes_{_{U_q\mathfrak{l}}}
U_q\mathfrak{g})^{U_q\mathfrak{k}}$. Сопряженные линейные операторы будем
называть инвариантными дифференциальными операторами на формальном
квантовом однородном пространстве $K\backslash G$. \IND{инвариантный
дифференциальный оператор ! на формальном квантовом однородном
пространстве} Они образуют алгебру ${\mathbf D}(K\backslash G)_q\,\cong\,
{\rm End}_{_{U_q\mathfrak{g}}}(\mathbb{C}\otimes _{_{U_q\mathfrak{k}}}
U_q\mathfrak{g})^{\rm op} $.

 Операторы из ${\mathbf D}(K\backslash G)$ являются эндоморфизмами
$U_q\mathfrak{g}$-модуля ${\rm
Hom}_{_{U_q\mathfrak{l}}}(U_q\mathfrak{g},\mathbb{C})$, и, следовательно,
отображают подпространство $\mathbb{C}[K\backslash G]_q$ в себя. Как
вытекает из определений, ``инвариантный дифференциальный оператор``,
отвечающий элементу
 $$
 1\otimes a \in (\mathbb{C}\otimes
_{_{U_q\mathfrak{k}}} U_q\mathfrak{g})^{U_q\mathfrak{k}} \cong {\rm
End}_{_{U_q\mathfrak{g}}}(\mathbb{C}\otimes _{_{U_q\mathfrak{k}}}
U_q\mathfrak{g}),
 $$
 является сужением линейного оператора $L_{\rm reg}(a)$ на подпространство
$\mathbb{C}[K\backslash G]_q$.

\begin{remark} Используя антипод $S$, нетрудно исключить из рассмотрения правые
$U_q\mathfrak{g}$-модули, заменив их левыми $U_q\mathfrak{g}$-модулями.
Именно,
$$
{\mathbf D}(K\backslash G)_q\cong
 {\rm End}_{_{U_q\mathfrak{g}}} (U_q\mathfrak{g}
 \otimes _{_{U_q\mathfrak{k}}} \mathbb{C})^{\rm op},
$$
и ``инвариантный дифференциальный оператор``, отвечающий элементу
$$
1\otimes a \in (U_q\mathfrak{g}
 \otimes _{_{U_q\mathfrak{k}}} \mathbb{C})^{U_q\mathfrak{k}}\cong
{\rm Hom}_{_{U_q\mathfrak{k}}} (\mathbb{C},U_q\mathfrak{g}
 \otimes _{_{U_q\mathfrak{k}}} \mathbb{C})\cong
 {\rm End}_{_{U_q\mathfrak{g}}} (U_q\mathfrak{g}
 \otimes _{_{U_q\mathfrak{k}}} \mathbb{C}),
$$
 является сужением линейного оператора $L_{\rm reg}(S(a))$ на
подпространство $\mathbb{C}[K\backslash G]_q$.
\end{remark}

\bigskip
 Перейдем от пространств функций  к пространствам
сечений однородных векторных расслоений на $K\backslash G$.


Рассмотрим конечномерный весовой $U_q\mathfrak{k}$-модуль $V$. Векторное
пространство
\begin{equation}\label{space_of_sections}
{\rm Hom}_{_{U_q\mathfrak{k}}}(U_q\mathfrak{g},V)\cong
(V^*\otimes_{_{U_q\mathfrak{k}}}\,U_q\mathfrak{g})^*
\end{equation}
назовем пространством сечений однородного векторного расслоения на
``формальном квантовом однородном пространстве $K\backslash G$``.
\IND{пространство сечений однородного векторного расслоения на ``формальном
квантовом однородном пространстве''} Оно естественным образом наделяется
структурой $U_q\mathfrak{g}$-модульного ${\rm
Hom}_{_{U_q\mathfrak{k}}}(U_q\mathfrak{g},\mathbb{C})$-бимодуля: достаточно
воспользоваться отображением
\begin{equation}\label{fiberwise_otimes}
{\rm Hom}_{_{U_q\mathfrak{k}}}(U_q\mathfrak{g},V_1)\times
{\rm Hom}_{_{U_q\mathfrak{k}}}(U_q\mathfrak{g},V_2)\to
{\rm Hom}_{_{U_q\mathfrak{k}}}(U_q\mathfrak{g},V_1\otimes V_2),
\end{equation}
где $(f_1(\xi),f_2(\xi))\mapsto(f_1\otimes f_2)\Delta(\xi)$, полагая
сначала $V_1=V$, $V_2=\mathbb{C}$, а потом $V_1=\mathbb{C}$, $V_2=V$.

\medskip

Введем инвариантные дифференциальные операторы так же, как в частном случае
$V_1=V_2=\mathbb{C}$. Если линейный оператор
$A:V_2^*\otimes_{_{U_q\mathfrak{k}}} U_q\mathfrak{g} \to
V_1^*\otimes_{_{U_q\mathfrak{k}}}U_q\mathfrak{g}$ является морфизмом правых
$U_q\mathfrak{g}$-модулей, то сопряженный линейный оператор
$$
A^*:{\rm Hom}_{_{U_q\mathfrak{k}}}(U_q\mathfrak{g},V_1)\to{\rm
Hom}_{_{U_q\mathfrak{k}}}(U_q\mathfrak{g},V_2)
$$
также является морфизмом $U_q\mathfrak{g}$-модулей. Такие линейные
операторы служат квантовыми аналогами инвариантных дифференциальных
операторов на $K\backslash G$. \IND{квантовый аналог ! инвариантного
дифференциального оператора на $K\backslash G$}

Значит, инвариантные дифференциальные операторы нумеруются элементами
пространства $ {\rm Hom}_{_{U_q\mathfrak{g}}}(V_2^*\otimes
_{_{U_q\mathfrak{k}}} U_q\mathfrak{g},V_1^*\otimes _{_{U_q\mathfrak{k}}}
U_q\mathfrak{g})$, то есть элементами пространства
$$
 {\rm
Hom}_{_{U_q\mathfrak{k}}}(V_2^*,V_1^*\otimes _{_{U_q\mathfrak{k}}}
U_q\mathfrak{g})\cong ((V_1^*\otimes _{_{U_q\mathfrak{k}}}
U_q\mathfrak{g})\otimes V_2)^{U_q\mathfrak{k}}.
$$

Опишем действие этих операторов, следуя \cite{KoranyiReimann}. Каждому
элементу $\sum\limits_jl_j\otimes\xi_j\otimes v_j\in V_1^*\otimes
U_q\mathfrak{g}\otimes V_2$ сопоставим линейный оператор из ${\rm
Hom}(U_q\mathfrak{g},V_1)$ в ${\rm Hom}(U_q\mathfrak{g},V_2)$:
\begin{equation}\label{simple_formulae}
\qquad f(\eta)\mapsto\sum_jl_j(f(\xi_j\,\eta))\,v_j,\qquad\eta\in
U_q\mathfrak{g}.
\end{equation}
Этот линейный оператор, очевидно, является квантовым аналогом инвариантного
дифференциального оператора на группе $G$. \IND{квантовый аналог !
инвариантного дифференциального оператора на группе $G$}

Элементам
$$
A\in\left(V_1^*\otimes_{_{U_q{\mathfrak{k}}}}U_q\mathfrak{g}\right)\otimes
V_2\cong{\rm Hom}(V_2^*,V_1^*\otimes_{_{U_q{\mathfrak{k}}}}U_q\mathfrak{g})
$$
отвечают линейные операторы $D_A:{\rm
Hom}_{_{U_q\mathfrak{k}}}(U_q\mathfrak{g},V_1)\to{\rm
Hom}(U_q\mathfrak{g},V_2)$, являющиеся морфизмами
$U_q\mathfrak{g}$-модулей. Предположим, что
$$
A\in{\rm Hom}_{_{U_q\mathfrak{k}}}(V_2^*,V_1^*\otimes_{_{U_q{\mathfrak{k}}}}
U_q\mathfrak{g}).
$$
Тогда $D_A\cdot{\rm Hom}_{_{U_q\mathfrak{k}}}(U_q\mathfrak{g},V_1)\subset
{\rm Hom}_{_{U_q\mathfrak{k}}}(U_q\mathfrak{g},V_2)$. Действительно, если
$\pi_j$ -- представления, отвечающие $U_q\mathfrak{l}$-модулям $V_j$, и
$A=\sum\limits_il_i\otimes\xi_i\otimes v_i$, то для всех $\zeta\in
U_q\mathfrak{k}$, $\eta\in U_q\mathfrak{g}$ имеют место равенства
$$
D_Af(\zeta\eta)=\sum_il_i(\pi_1(\zeta)f(\xi_i\eta))v_i=
\sum_il_i(f(\xi_i\eta))\pi_2(\zeta)v_i=\pi_2(\zeta)\,D_Af(\eta),
$$
где $f\in{\rm Hom}_{_{U_q\mathfrak{k}}}(U_q\mathfrak{g},V_1)$. Получаем
морфизм $U_q\mathfrak{g}$-модулей
\begin{equation}\label{diff_op_formal}
D_A:{\rm Hom}_{_{U_q\mathfrak{k}}}(U_q\mathfrak{g},V_1)\to{\rm
Hom}_{_{U_q\mathfrak{k}}}(U_q\mathfrak{g},V_2).
\end{equation}

\medskip

Перейдем от формальных сечений к регулярным сечениям рассматриваемых
однородных векторных расслоений и их квантовых аналогов. Пусть
$\mathbb{V}_1$ и $\mathbb{V}_2$ -- наибольшие весовые локально
конечномерные подмодули $U_q\mathfrak{g}$-модулей ${\rm
Hom}_{_{U_q\mathfrak{k}}}(U_q\mathfrak{g},V_1)$, ${\rm
Hom}_{_{U_q\mathfrak{k}}}(U_q\mathfrak{g},V_2)$ соответственно. Это
$U_q\mathfrak{g}$-модульные $\mathbb{C}[K\backslash G]_q$-бимодули. Они
являются квантовыми аналогами пространств (регулярных) сечений однородных
векторных расслоений на $K\backslash G$. \IND{квантовый аналог !
пространств (регулярных) сечений однородных векторных расслоений на
$K\backslash G$} Очевидно, $D_A:\mathbb{V}_1\to\mathbb{V}_2$, и мы получаем
квантовые аналоги инвариантных дифференциальных операторов в пространствах
сечений однородных векторных расслоений. \IND{квантовый аналог !
инвариантных дифференциальных операторов в пространствах сечений однородных
векторных расслоений}

 \begin{remark} Так же, как в частном случае $V_1=V_2=\mathbb{C}$, с помощью антипода $S$
 можно исключить из рассмотрения правые $U_q\mathfrak{g}$-модули,
 заменив их левыми $U_q\mathfrak{g}$-модулями.
  \end{remark}


 \subsubsection{Инвариантные дифференциальные операторы
 на $K\backslash G$: продолжение.}\label{vect_diff}
 Хорошо известно \cite[стр. 329]{Helg1}, что алгебра
 ${\mathbf D}(K\backslash G)$  инвариантных дифференциальных операторов
 на $K\backslash G$ коммутативна.
 Получим $q$-аналог этого результата, доказав коммутативность алгебры
 ${\mathbf D}(K\backslash G)_q$.

  Естественное представление этой алгебры в векторном пространстве
   $\mathbb{C}[K\backslash G]_q$ является точным, поскольку спаривание
 $$
\mathbb{C}[K\backslash G]_q\times
 \left(\mathbb{C} \otimes_{_{U_q\mathfrak{k}}} U_q\mathfrak{g}\right)
 \to \mathbb{C},  \qquad f \times (1\otimes \xi) \mapsto f(\xi),
 $$
 невырождено (доказательство см., например, в \cite{HeckenbergerKolb04}).
 Значит, коммутативность  алгебры ${\mathbf D}(K\backslash G)_q$ следует
 из коммутативности алгебры эндоморфизмов $U_q\mathfrak{g}$-модуля
 $\mathbb{C}[K\backslash G]_q$, которая, в свою очередь, вытекает из
 отсутствия кратностей в разложении $U_q\mathfrak{g}$-модуля
 $\mathbb{C}[K\backslash G]_q$ в сумму простых $U_q\mathfrak{g}$-модулей
\begin{equation}\label{spher_submodules}
\mathbb{C}[K\backslash G]_q \approx \newoplus\limits_{\lambda \in P_+^{\rm
spher}} L(\lambda),
\end{equation}
 см. предложение \ref{PeterWeyl} и \itemiii\ \ref{loc_finite}.

\medskip
Приведем описание алгебры ${\mathbf D}(K\backslash G)_q$, используемое
Балдони и Фраджриа в статье \cite[стр. 1187]{Baldoni_Frajria}, и обсудим
основной результат этой работы.

Напомним, что $U_q\mathfrak{g}$ является $U_q\mathfrak{k}$-модульной
алгеброй по отношению к присоединенному действию и что морфизм
$U_q\mathfrak{k}$-модулей
$$
m:\,U_q\mathfrak{k}\otimes U_q\mathfrak{p}^-\otimes U_q\mathfrak{p}^+,
\qquad m:\,\xi\otimes \eta \otimes \zeta \mapsto \xi \eta \zeta
$$
биективен, см. \itemiii\ \ref{decomp_Dima}. Значит,
$$
 {\mathbf D}(K\backslash G)_q\cong
 (\mathbb{C}\otimes_{_{U_q\mathfrak{k}}}U_q\mathfrak{g})^{U_q\mathfrak{k}}
 \cong (U_q\mathfrak{p}^-\, U_q\mathfrak{p}^+)^{U_q\mathfrak{k}},
$$
причем инвариантный дифференциальный оператор, отвечающий элементу $a\in
(U_q\mathfrak{p}^-\, U_q\mathfrak{p}^+)^{U_q\mathfrak{k}}$, является
 сужением линейного оператора
 $L_{\rm reg}(a)$ на  инвариантное подпространство
$\mathbb{C}[K\backslash G]_q$.

\begin{example}\label{Laplace_L_reg} Однородные компоненты
$(U_q\mathfrak{p}^+)_1$, $(U_q\mathfrak{p}^-)_{-1}$ являются простыми и
 конечномерными весовыми $U_q\mathfrak{k}$-модулями, и они
сопряжены. Значит, $\dim\,
 \left(
 (U_q\mathfrak{p}^-)_{-1} (U_q\mathfrak{p}^+)_1\right)^{U_q\mathfrak{k}}\;=
 \dim\,
 \left(
 (U_q\mathfrak{p}^-)_{-1}\otimes (U_q\mathfrak{p}^+)_1\right)^{U_q\mathfrak{k}}\; =\; 1
$. Следовательно, существует и единствен с точностью до числового множителя
элемент алгебры $D(K\backslash G)_q$, отвечающий ненулевому
$U_q\mathfrak{k}$-инварианту из $
(U_q\mathfrak{p}^-)_{-1}(U_q\mathfrak{p}^+)_1$. Это квантовый аналог
инвариантного лапласиана.\IND{инвариантный ! лапласиан} В простейшем случае
$\mathfrak{g}=\mathfrak{sl}_2$ элемент $FE$ принадлежит
$\left((U_q\mathfrak{p}^-)_{-1}
(U_q\mathfrak{p}^+)_1\right)^{U_q\mathfrak{k}}$ и для всех $f\in
\mathbb{C}[K\backslash \mathrm{SL}_2]_q $ имеем:$
L_{\mathrm{reg}}(FE)f=L_{\mathrm{reg}}(C_q)f=R_{\mathrm{reg}}(C_q)f$,
  где $C_q$-- квантовый аналог элемента Казимира, см. \eqref{C_q}. О его
связи с инвариантным лапласианом в квантовом круге см.
 \itemiiiы  \ref{radial_part}, \ref{canon_example_sl_2}.
\end{example}

Аналогичные рассуждения, приведенные в \cite{Baldoni_Frajria}, позволяют
ввести в рассмотрение выделенный базис векторного пространства
$D(K\backslash G)_q$. Роль простых $U_q\mathfrak{k}$-модулей
$(U_q\mathfrak{p}^+)_1$, $(U_q\mathfrak{p}^-)_{-1}$ играют компоненты
 разложений
\begin{equation}\label{HuaSchmid_components}
U_q\mathfrak{p}^+=\newoplus\limits_{\lambda \in \mathcal{A}_+}
(U_q\mathfrak{p}^+)_\lambda,\qquad
U_q\mathfrak{p}^-=\newoplus\limits_{\lambda \in \mathcal{A}_+}
(U_q\mathfrak{p}^-)_\lambda,
\end{equation}
 в категории $U_q\mathfrak{k}$-модулей, где
$ (U_q\mathfrak{p}^+)_\lambda\approx L(\mathfrak{k},\lambda)$,
 $(U_q\mathfrak{p}^-)_\lambda\approx L(\mathfrak{k},\lambda)^*$,
а множество $\mathcal{A}_+$ вводится равенством \eqref{HuaSchmid_lattice}.
Разложения \eqref{HuaSchmid_components} легко получить с помощью
изоморфизма
 $U_q\mathfrak{k}$-модульных алгебр $U_q\mathfrak{p}^+ \cong
 \mathbb{C}[\mathfrak{p}^-]_q$ ( см. \itemiii\ \ref{Add_Uq}),
 разложения Хуа-Шмида (см. предложение  \ref{k-types})
 и результатов \itemiiiа \ref{loc_finite}.

\begin{remark}\label{shimura_operators} Замена в предыдущих рассуждениях
тривиального $U_q\mathfrak{k}$-модуля $\mathbb{C}$ одномерным
$U_q\mathfrak{k}$-модулем общего вида приводит к квантовым аналогам
операторов Шимуры \cite{Shimura}.
\end{remark}

 В работе \cite{Baldoni_Frajria}  доказано, что для неисключительных
  алгебр Ли $\mathfrak{g}$ все инвариантные дифференциальные операторы в
$\mathbb{C}[K\backslash G]_q$ можно получить, действием элементов центра
$Z(U_q\mathfrak{g})$ алгебры $U_q\mathfrak{g}$. Как известно
\cite{Helg_red_groups}, ограничение на $\mathfrak{g}$ нельзя отбросить даже
в случае $q=1$. В предположении о неисключительности алгебры Ли
$\mathfrak{g}$ получен и основной результат работы \cite{Baldoni_Frajria}--
описан образ алгебры ${\bf D}(K\backslash G)_q$ при гомоморфизме
Хариш-Чандры, который несущественным сдвигом отличается от вложения ${\bf
D}(K\backslash G)_q$ в алгебру функций на $P_+^{\rm spher}$, отвечающего
разложению \eqref{spher_submodules}.

\bigskip
 Построения \itemiiiа \ref{pol_diff} доставляют ковариантные
дифференциальные исчисления первого порядка
$(\Omega^{1,0}(\mathfrak{p}^-)_q,\partial)$,
$(\Omega^{0,1}(\mathfrak{p}^-)_q,\bar{\partial})$ над
$U_q\mathfrak{g}$-модульной алгеброй $\mathrm{Pol}(\mathfrak{p}^-)_q$.
Опишем в самых общих чертах другой подход к этим дифференциальным
исчислениям. Для упрощения изложения ограничимся первым из них и будем
предполагать, что $\mathbb{D}\subset \mathfrak{p}^-$ является ограниченной
симметрической областью трубчатого типа, см.
\itemiii\ \ref{prehomogeneous}.

Рассмотрим инвариантный дифференциальный оператор $\partial'$, отвечающий
морфизму $U_q\mathfrak{k}$-модулей
$$
(U_q\mathfrak{p}^+)_1\to \mathbb{C}\otimes_{_{U_q\mathfrak{k}}}
U_q\mathfrak{g}, \qquad E_{_{l_0}}\mapsto 1\otimes E_{_{l_0}}.
$$
Как показано в \cite{HeckenbergerKolb04},
$$
 \partial'(f_1 f_2)\,=\, \partial'(f_1)\, f_2 + f_1\, \partial'(f_2), \qquad
f_1, f_1 \in \mathbb{C}[K\backslash G]_q,
$$
то есть линейный оператор $\partial'$ доставляет ковариантное
дифференциальное исчисление первого порядка над алгеброй
$\mathbb{C}[K\backslash G]_q$.
 Продолжим его на локализацию
$\mathbb{C}[K\backslash G]_{q,x}$ алгебры $\mathbb{C}[K\backslash G]_q$ по
мультипликативному множеству $x^{\mathbb{Z}_+}$ (существование продолжения
доказывается так же, как предложение \ref{*_local_gen}, с использованием
коммутационных соотношений между $x$ и матричными элементами
$c^\Lambda_{\lambda,i;\mu,j}$, см.
\itemiii\ \ref{CommRelations}).

Остается перенести полученное ковариантное дифференциальное исчисление
первого порядка с $\mathbb{C}[K\backslash G]_{q,x}$ на локализацию
$\mathrm{Pol}(\mathfrak{p}^-)_{q,y}$ алгебры
$\mathrm{Pol}(\mathfrak{p}^-)_q$ по мультипликативному множеству
$y^{\mathbb{Z}_+}$ с помощью канонического изоморфизма \hbox{$\mathcal{I}:
 \mathrm{Pol}(\mathfrak{p}^-)_{q,y}\stackrel{\approx}{\to}
\mathbb{C}[K\backslash G]_{q,x}$,} см. \itemiii\ \ref{reps_pol}.

То, что возникающее ковариантное дифференциальное исчисление первого
 порядка над подалгеброй
 $\mathrm{Pol}(\mathfrak{p}^-)_q\subset \mathrm{Pol}(\mathfrak{p}^-)_q$
изоморфно $(\Omega^{1,0}(\mathfrak{p}^-)_q,\partial)$, нетрудно доказать
используя предложение \ref{from_Jo-R} и соображения двойственности, см.
\itemiiiы \ref{rel-for-f_0}, \ref{affine_diff}, \ref{FODC}, \ref{pol_diff}.

\bigskip Напомним обозначение $K_1=w_0Kw_0^{-1}$, см. \eqref{K_1}. Предложение \ref{image_emb} наводит на
мысль о том, что каждому инвариантному дифференциальному оператору в
$\mathbb{C}[K_1\backslash G]_q$ отвечает инвариантный дифференциальный
оператор в квантовой ограниченной симметрической области $\mathbb{D}$.
Опишем этот линейный оператор.

Рассмотрим подалгебру $U_q\mathfrak{k}_1$, порожденную множеством $\{E_j,
F_j\}_{j\neq l_1}\cup \{K_i^{\pm 1}\}_{i=1,2,\ldots,l}$, где $l_1$
определяется равенством $\bar{\omega}_{l_1}=-w_0\bar{\omega}_{l_0}$, см.
\itemiii\ \ref{w0_G_x}. Разумеется, $U_q\mathfrak{k}_1$-- квантовая
универсальная обертывающая алгебра, отвечающая группе $K_1$.

Рассмотрим алгебру $\mathbb{C}[G]_q^*$. Так как $\widetilde{w}_0\in
\mathbb{C}[G]_q^*$ и $U_q\mathfrak{g}$ вложено в $\mathbb{C}[G]_q^*$, то
$\widetilde{w}_0U_q\mathfrak{g}\subset \mathbb{C}[G]_q^*$. Напомним, что
$\mathbb{C}[G]_q^*$ является $U_q\mathfrak{g}$-бимодулем, а
$\widetilde{w}_0 U_q\mathfrak{g}$-- его подбимодулем
\begin{equation}\label{bi_mod_w_0}
\xi\times \eta: \, \widetilde{w}_0\zeta \mapsto \widetilde{w}_0
(\widetilde{w}_0^{-1}\xi\widetilde{w}_0)\zeta\eta,
\end{equation}
см. предложение \ref{from_Jo-R}.

\begin{proposition}\label{nonstandard_loc}
Если $A$-- инвариантный дифференциальный оператор в
$\mathbb{C}[K_1\backslash G]_q$ и $A^*$-- сопряженный линейный оператор в
$\mathbb{C}[K_1\backslash G]_q^*$, то
\begin{equation}\label{support_in_w_0}
A^*\, (\mathbb{C}\otimes_{_{U_q\mathfrak{k}_1}}\widetilde{w}_0
U_q\mathfrak{g}) \subset
\mathbb{C}\otimes_{_{U_q\mathfrak{k}_1}}\widetilde{w}_0 U_q\mathfrak{g}.
\end{equation}
\end{proposition}

{\bf Доказательство.} Как следует из предложения \ref{from_Jo-R}, включение
\eqref{support_in_w_0} равносильно включению
$$
A^*\, (\mathbb{C}\otimes_{_{U_q\mathfrak{k}_1}}
U_q\mathfrak{g}\,\widetilde{w}_0) \subset
\mathbb{C}\otimes_{_{U_q\mathfrak{k}_1}} U_q\mathfrak{g}\,\widetilde{w}_0,
$$
а оно, в свою очередь, вытекает из очевидного включения
$$
A^*\, (\mathbb{C}\otimes_{_{U_q\mathfrak{k}_1}} U_q\mathfrak{g}) \subset
\mathbb{C}\otimes_{_{U_q\mathfrak{k}_1}} U_q\mathfrak{g}
$$
и из того, что линейный оператор $A$ является морфизмом
$U_q\mathfrak{g}$-модулей. \hfill $\square$

\medskip Используя антипод $S$, наделим $\mathscr{D}(\mathbb{D})_q$
структурой правого $U_q\mathfrak{g}$-модуля, изоморфного
 $\mathbb{C}\otimes_{_{U_q\mathfrak{k}}} U_q\mathfrak{g}$:
 $$
\mathbb{C}\otimes_{_{U_q\mathfrak{k}}} U_q\mathfrak{g}\to
\mathscr{D}(\mathbb{D})_q,\qquad 1\otimes 1 \mapsto f_0,
 $$
см. \itemiii\ \ref{rel-for-f_0}.
 Нетрудно перенести описанный в предложении
\ref{nonstandard_loc} линейный оператор с векторного пространства
 $\mathbb{C}\otimes_{_{U_q\mathfrak{k}_1}}
 \left(\widetilde{w}_0 U_q\mathfrak{g}\right)$ на векторное пространство
  $\mathbb{C}\otimes_{_{U_q\mathfrak{k}}} U_q\mathfrak{g}$ и, наконец, на
$\mathscr{D}(\mathbb{D})_q$. Для этого достаточно воспользоваться
изоморфизмом правых $U_q\mathfrak{g}$-модулей
\begin{equation}\label{delete_w_0}
\mathbb{C}\otimes_{_{U_q\mathfrak{k}}} U_q\mathfrak{g}
\to
 \mathbb{C}\otimes_{_{U_q\mathfrak{k}_1}}
 \left(\widetilde{w}_0 U_q\mathfrak{g}\right),\qquad
  1\otimes 1 \mapsto 1\otimes \widetilde{w}_0.
\end{equation}
(Существование гомоморфизма \eqref{delete_w_0} вытекает из
$U_q\mathfrak{k}$-инвариантности элемента $1\otimes \widetilde{w}_0\in
 \mathbb{C}\otimes_{_{U_q\mathfrak{k}_1}} U_q\mathfrak{g}$,
 см. \eqref{bi_mod_w_0}. Единственность такого гомоморфизма очевидна, а его
 биективность следует из результатов
 \itemiiiов \ref{decomp_Dima}, \ref{w_0}.)

Мы получили инвариантный дифференциальный оператор в
$\mathscr{D}(\mathbb{D})_q$, отвечающий инвариантному дифференциальному
оператору в $\mathbb{C}[K_1\backslash G]_q$.
 Отметим, что нетрудно перейти от инвариантного линейного оператора в пространстве
$\mathscr{D}(\mathbb{D})_q$
 финитных функций к инвариантному линейному оператору в пространстве
 $\mathscr{D}(\mathbb{D})'_q$ обобщенных функций в квантовой ограниченной
 симметрической области.

\subsubsection{Инвариантные дифференциальные операторы на большой клетке
обобщенного пространства флагов.}\label{projective_diff}

Рассмотрим $\mathfrak{g}$-модуль $L(\overline{\omega}_{l_0})$ со старшим
весом $\overline{\omega}_{l_0}$ и старшим вектором
$v(\overline{\omega}_{l_0})$. Как известно, $G$-орбита точки
$\mathbb{C}v(\overline{\omega}_{l_0})$ проективного пространства
$P(L(\overline{\omega}_{l_0}))$ является обобщенным пространством флагов и
обладает стандартным разбиением на клетки Шуберта $X^w$ коразмерности
$l(w)$, нумеруемые элементами $w\in W^\mathbb{S}$. Точка
$\mathbb{C}v(\overline{\omega}_{l_0})$ принадлежит большой клетке $X^e$
этого разбиения.

\medskip

Повторим построения \itemiiiа \ref{affine_diff}, заменяя $\mathfrak{k}$
параболической подалгеброй Ли $\mathfrak{q}^+\supset\mathfrak{k}$,
порожденной множеством
$$\{E_i,H_i\;|\;i=1,2,\ldots,l\}\;\cup\;\{F_j\;|\;j\ne l_0\}.$$
Квантовым аналогом $U\mathfrak{g}$-модульной алгебры $\mathbb{C}[X^e]$
регулярных функций на $X^e$ будет служить $U_q\mathfrak{g}$-модульная
алгебра
$$
\mathbb{C}[U]_q\,=\,\{f\in{\rm
Hom}_{_{U_q\mathfrak{q}^+}}(U_q\mathfrak{g},\mathbb{C})\;|\;
\dim(U_q\mathfrak{h}\cdot f)<\infty\},
$$
а квантовым аналогом алгебры формальных рядов --
$U_q\mathfrak{g}$-модульная алгебра ${\rm
Hom}_{_{U_q\mathfrak{q}^+}}(U_q\mathfrak{g},\mathbb{C})\cong
(\mathbb{C}\otimes_{_{U_q\mathfrak{q}^+}}U_q\mathfrak{g})^*$.

\medskip

 Перейдем от функций на $X^e$ к сечениям однородных векторных расслоений.

Каждый конечномерный весовой $U_q\mathfrak{k}$-модуль $V$ наделим
структурой $U_q\mathfrak{q}^+$-модуля с помощью естественного эпиморфизма
алгебр $U_q\mathfrak{q}^+\to U_q\mathfrak{k}$. Квантовым аналогом
пространства сечений однородного векторного расслоения со слоем $V$ над
$X^e$ будет служить весовой $U_q\mathfrak{g}$-модуль
\begin{equation}\label{fiber_bundle}
\mathbb{V}\,=\,\{f\in{\rm
Hom}_{_{U_q\mathfrak{q}^+}}(U_q\mathfrak{g},V)\;|\;
\dim(U_q\mathfrak{h}\cdot f)<\infty\},
\end{equation}
а квантовым аналогом соответствующего пространства ``формальных рядов'' --
векторное пространство ${\rm
Hom}_{_{U_q\mathfrak{q}^+}}(U_q\mathfrak{g},V)$.

Описанный подход к построению пространств сечений однородных векторных
расслоений хорошо известен в классическом случае $q=1$ \cite[стр.
52]{Vogan_CMP}.

В дальнейшем существенно используется то, что естественное спаривание
$$
{\rm Hom}_{_{U_q\mathfrak{q}^+}}(U_q\mathfrak{g},V)\times
 (V^*\otimes_{_{U_q{\mathfrak{q}^+}}} U_q\mathfrak{g})\to \mathbb{C},
 \qquad
  f\times (l\otimes \xi) \mapsto l(f(\xi))
$$
доставляет канонический изоморфизм $U_q\mathfrak{g}$-модулей
\begin{equation*}\label{iso_Jakobsen}
 {\rm Hom}_{_{U_q\mathfrak{q}^+}}(U_q\mathfrak{g},V) \cong
 (V^*\otimes_{_{U_q\mathfrak{q}^+}} U_q\mathfrak{g})^*.
\end{equation*}
Здесь $V^*\otimes_{_{U_q\mathfrak{q}^+}} U_q\mathfrak{g}$ является
 правым $U_q\mathfrak{g}$-модулем, поскольку, как и в \itemiiiе
\ref{affine_diff}, мы не используем антипод $S$ при переходе к сопряженным
пространствам.

\medskip Введем квантовые аналоги инвариантных дифференциальных операторов
на $X^e$. Рассмотрим весовые конечномерные $U_q\mathfrak{k}$-модули $V_1$,
$V_2$, отвечающие им $U_q\mathfrak{q}^+$-модули $V_1$, $V_2$ и пространства
``сечений однородных векторных расслоений`` $\mathbb{V}_1 \subset
 {\rm Hom}_{_{U_q\mathfrak{q}^+}}(U_q\mathfrak{g},V_1)$,
 $\mathbb{V}_2 \subset
 {\rm Hom}_{_{U_q\mathfrak{q}^+}}(U_q\mathfrak{g},V_2)$.

Если линейный оператор
$A:V_2^*\otimes_{_{U_q\mathfrak{q}^+}}U_q\mathfrak{g}\to
V_1^*\otimes_{_{U_q\mathfrak{q}^+}}U_q\mathfrak{g}$ является морфизм правых
$U_q\mathfrak{g}$-модулей, то сопряженный линейный оператор
$$
A^*:\;{\rm Hom}_{_{U_q\mathfrak{q}^+}}(U_q\mathfrak{g},V_1)\to{\rm
Hom}_{_{U_q\mathfrak{q}^+}}(U_q\mathfrak{g},V_2)
$$
является морфизмом левых $U_q\mathfrak{g}$-модулей и отображает
$\mathbb{V}_1$ в $\mathbb{V}_2$. Такие операторы $A^*$ будем называть
квантовыми инвариантными дифференциальными операторами на большой клетке
$X^e$. \IND{квантовый инвариантный дифференциальный оператор на большой
клетке $X^e$}

 Они  нумеруются элементами
 $
{\rm Hom}_{_{U_q\mathfrak{g}}}(V_2^*\otimes_{_{U_q\mathfrak{q}^+}}
U_q\mathfrak{g},V_1^*\otimes_{_{U_q\mathfrak{q}^+}} U_q\mathfrak{g})
 $, то есть элементами
 $
{\rm Hom}_{_{U_q\mathfrak{q}^+}}(V_2,V_1^*\otimes_{_{U_q\mathfrak{q}^+}}
 U_q\mathfrak{g})\cong
((V_1^*\otimes_{_{U_q\mathfrak{q}^+}}
 U_q\mathfrak{g})\otimes V_2)^{U_q\mathfrak{q}^+}.
 $


\begin{remark} 1. В классическом пределе $q=1$ наш подход  к построению
инвариантных дифференциальных операторов согласуется с их описанием в
работах Харрисона и Якобсена \cite{Har-Jakob, Jakob_basic}.

2. Инвариантные дифференциальные операторы на большой клетке $X^e$
естественно возникают в задачах теории представлений
\cite{DavEnrSt1,DavEnrSt2}, математической физики \cite{BaEast} и
алгебраической геометрии \cite{Enr_class, Enr_except, EnrHunz}.
\end{remark}

\newpage

\section{ ЗАКЛЮЧЕНИЕ}\label{conclusion}

\subsection{Границы и сферические основные серии}

Подведем итоги:
\begin{itemize}
\item[-] в главе \ref{basic_theory} вводятся в рассмотрение квантовые
ограниченные симметрические области общего вида и, тем самым, предъявляется
список объектов исследования,

\item[-]
 результаты главы \ref{quantum_disc_part} о квантовом круге
позволяют сформировать список интересующих нас вопросов некоммутативного
комплексного анализа и некоммутативного гармонического анализа.
\end{itemize}

Типичная задача состоит в ответе на один из этих вопросов применительно к
некоторому множеству квантовых ограниченных симметрических областей.

Такой подход к постановке задач весьма наивен, поскольку многие
содержательные результаты об ограниченных симметрических областях
тривиализуются в случае вещественного ранга 1. Это означает, что список
интересующих нас вопросов должен быть существенно дополнен. Приведем
примеры вопросов, не возникающих в случае квантового круга.

 \bigskip
Начнем с предварительных сведений о простейшей серии квантовых ограниченных
симметрических областей трубчатого типа.

Пусть $ n \in \mathbb{N}$ и $N=2n$. Наделим векторное пространство линейных
операторов в $\mathbb{C}^n$ и канонически изоморфное ему векторное
пространство $\mathrm{Mat}_{n}$ комплексных $n\times n$-матриц операторными
нормами. Как известно \cite[стр. 387]{Helg}, единичный шар $ \mathbb{D} =
\{\mathbf{z} \in \mathrm{Mat}_{n}\ |\ \mathbf{z} \mathbf{z}^* < 1\} $
является ограниченной симметрической областью вещественного ранга $n$.
Отвечающую ей квантовую ограниченную симметрическую область называют
квантовым матричным шаром. В этом случае
$$l=N-1,\quad l_0=n,\quad U_q\mathfrak{g}=U_q\mathfrak{sl}_N,\quad
U_q\mathfrak{k}=U_q\mathfrak{s}(\mathfrak{gl}_n\times
\mathfrak{gl}_n),\quad \mathfrak{p}^-=\mathrm{Mat}_{n}$$
 и используется обозначение
$U_q\mathfrak{su}_{n,n}=(U_q\mathfrak{sl}_N,*)$, см. \eqref{star_1}.

Известны описания $U_q\mathfrak{sl}_N$-модульной алгебры
$\mathbb{C}[\mathfrak{p}^-]_q=\mathbb{C}[\mathrm{Mat}_{n}]_q$ и
$U_q\mathfrak{su}_{n,n}$-модульной $*$-алгебры
$\mathrm{Pol}(\mathfrak{p}^-)_q=\mathrm{Pol}(\mathrm{Mat}_{n})_q$ с помощью
множеств образующих
$$
\{z^\alpha_a\}_{\alpha,\, a=1,2,\ldots, n}\subset
\mathbb{C}[\mathfrak{p}^-]_{q,1},\qquad \{z^\alpha_a,
(z^\alpha_a)^*\}_{\alpha,\, a=1,2,\ldots, n}\subset
\mathbb{C}[\mathfrak{p}^-]_{q,1} \newoplus
\mathbb{C}[\mathfrak{p}^+]_{q,-1}
$$
 и коммутационных соотношений  \cite{SSV_CJP,SSV4}. Образующие
$z^\alpha_a$ являются квантовыми аналогами матричных элементов матрицы
$\mathbf{z}=(z^\alpha_a)$, рассматриваемых как линейные функции на
$\mathrm{Mat}_{n}$. Коммутационные соотношения между ними хорошо известны:
$$
 z_a^\alpha z_b^\beta=\left
\{\begin{array}{ccl}qz_b^\beta z_a^\alpha &,&
a=b \;\&\;\alpha<\beta \quad{\rm or}\quad a<b \;\&\;\alpha=\beta \\
z_b^\beta z_a^\alpha &,& a<b \;\&\;\alpha>\beta \\ z_b^\beta
z_a^\alpha+(q-q^{-1})z_a^\beta z_b^\alpha &,& a<b \;\&\;\alpha<\beta
\end{array} \right..
$$

\bigskip
\bigskip  При $n=1$, то есть в случае квантового круга $\mathbb{D}$,
имеется одна сферическая основная унитарная серия модулей Хариш-Чандры, см.
 \itemiii\ \ref{boundness_Laplace}. Она допускает геометрическую реализацию в
пространстве полиномов Лорана, то есть в пространстве функций на единичной
окружности $\overline{\mathbb{D}}\setminus \mathbb{D}$, а унитаризующее
скалярное произведение вводится с помощью положительного инвариантного
интеграла на $\overline{\mathbb{D}}\setminus \mathbb{D}$, см.
\itemiii\ \ref{q-cone}.

При переходе от $n=1$ к $n>1$ ситуация кардинально меняется уже в
классическом случае $q=1$. Во-первых, растет число основных сферических
унитарных серий. Во-вторых, действие группы $\mathrm{Aut}(\mathbb{D})$ на
$\overline{\mathbb{D}}\setminus \mathbb{D}$ перестает быть транзитивным,
что препятствует построению геометрических реализаций основных сферических
унитарных серий в пространствах функций на топологической границе области
$\mathbb{D}$. Нужны другие границы, например, граница Шилова
$S(\mathbb{D})$, либо максимальная граница Сатаке-Ферстенберга
\cite{Koranyi_Inv, Koranyi_harmonic}.

Квантовый аналог последней из них доставляет геометрическую реализацию {\sl
невырожденной} основной сферической серии унитаризумых
$U_q\mathfrak{su}_{n,n}$-модулей. Эта серия участвует в разложении
квазирегулярного представления $*$-алгебры $U_q\mathfrak{su}_{n,n}$ в
предгильбертовом пространстве $\mathscr{D}(\mathbb{D})_q$ в прямой интеграл
неприводимых $*$-представлений. Известен явный вид меры Планшереля.

 Квантовый аналог границы
Шилова $S(\mathbb{D})$ матричного шара и соответствующая {\sl вырожденная}
 основная сферическая серия  $U_q\mathfrak{su}_{n,n}$-модулей
изучались в работах \cite{Vak01} и \cite{Bershtein1}. В классическом случае
$q=1$ граница Шилова $S(\mathbb{D})$ является множеством унитарных матриц
$U_n$ и, следовательно, вещественным аффинным алгебраическим многообразием,
описываемым матричным уравнением $\mathbf{z}\mathbf{z}^* = I$. Приведем
аналогичное описание границы Шилова квантового матричного шара.

Рассмотрим двусторонний идеал $J$ алгебры $\mathrm{Pol}(\mathrm{Mat}_n)_q$,
порожденный элементами
$$
\sum\limits_{j=1}^n q^{2n-\alpha-\beta} z^\alpha_j(z^\beta_j)^*
-\delta^{\alpha \beta},\qquad \alpha,\beta=1,2,\ldots,n.
$$
Очевидно, $J=J^*$ и, как можно показать, этот двусторонний идеал
$U_q\mathfrak{sl}_N$-инвариантен. Следовательно, фактор-алгебра
$\mathbb{C}[S(\mathbb{D})]_q=\mathrm{Pol}(\mathrm{Mat}_n)_q / J$ является
$U_q\mathfrak{su}_{n,n}$-модульной $*$-алгеброй. Ее называют алгеброй
регулярных функций на границе Шилова квантового матричного шара.
Канонический гомоморфизм $\mathrm{Pol}(\mathrm{Mat}_n)_q\to
\mathbb{C}[S(\mathbb{D})]_q$ является квантовым аналогом гомоморфизма,
сопоставляющего полиномиальным функциям их сужения на границу Шилова.

\begin{remark}
 Геометрические методы исследования представлений вещественных
полупростых групп Ли имеют богатую историю. Ограничимся ссылками на работы
 М.~Граева \cite{Graev_MMO},  Хариш-Чандры \cite{HCh6},
 М.~Нарасимхана и К. Окамото \cite{NarasOk}, А.~Кнаппа и Н.~Уоллака
 \cite{KW},
 Дж.~Вольфа \cite{Wolf,Wolf_MAMS, Wolf_Contemp},
 У.~Шмида \cite{Schmid_thesis, Schmid}, Л.~Барчини \cite{Barch1,Barch2}, Р.~Жиро
 \cite{Zierau, Zierau_IAS}, Дж. Арази и Г.~Упмайера \cite{ArUp}, Г.~Вонга \cite{Wong},
 Дж.-Т. Чанга  \cite{Chang, Chang_Proc}. Р.~Бэстона и  М.~Иствуда
 \cite{BaEast}, Л.~Барчини, М.~Иствуда и А.~Говера \cite{BarEastGov},
Л.~Барчини, А.~Кнаппа и Р.~Жиро \cite{BKZ, Zierau_IAS},
 Дж.~Вольфа и Р.~Жиро \cite{Wolf_Zierau}.
    Как правило, используемые  в этих работах  результаты о вещественных и
 комплексных многообразиях были известны задолго до того, как нашли
 приложения к построению геометрических реализаций представлений
вещественных полупростых групп Ли.
 В теории квантовых ограниченных симметрических областей
 сложилась принципиально иная ситуация. В
литературе о некоммутативной геометрии практически отсутствуют результаты,
необходимые для построения
 геометрических реализаций  модулей Хариш-Чандры над алгеброй $U_q\mathfrak{g}$,
 см. \itemii\
\ref{Vect_and_HCh}.
 Это вынуждает  вводить необходимые геометрические понятия параллельно с
построением геометрических реализаций модулей Хариш-Чандры.
\end{remark}

\medskip
Используя изоморфизм $U_q \mathfrak{s}(\mathfrak{gl}_n \times
\mathfrak{gl}_n)$-модульных алгебр
$$
\mathbb{C}[S(\mathbb{D})]_q \stackrel{\approx}{\to} \mathbb{C}[U_n]_q,
\qquad
  z^\alpha_a \mapsto q^{\alpha-n} z^\alpha_a,\quad
\alpha,a=1,2,\ldots,n,
$$
легко доказать существование и единственность $U_q
\mathfrak{s}(\mathfrak{gl}_n \times
\mathfrak{gl}_n)$-ин\-ва\-ри\-ант\-но\-го интеграла
$$ \nu: \mathbb{C}[S(\mathbb{D}]_q \to \mathbb{C}, \qquad \nu:f
\mapsto \int\limits_{S(\mathbb{D})_q} f d\nu,$$ для которого
$\int\limits_{S(\mathbb{D})_q} 1 d\nu=1$. Этот интеграл положителен и
используется при построении интегральных операторов из
$\mathbb{C}[S(\mathbb{D})]_q$ в $\mathscr{D}(\mathbb{D})_q'$.

\bigskip В классическом
случае $q=1$ ядром Пуассона называют функцию
\begin{equation*}\label{classical_Poisson}
 P(\mathbf{z},\mathbf{\zeta})=\frac{\mathrm{det}(1-\mathbf{zz}^*)^n}
{|\mathrm{det}(1-{\mathbf{\zeta}^*\mathbf{z}})|^{2n}}, \qquad \mathbf{z}
\in \mathbb{D},\; \mathbf{\zeta}\in S(\mathbb{D}).
\end{equation*}
на декартовом произведении $\mathbb{D}\times S(\mathbb{D})$, см. \cite[стр.
98]{Hua}.

Определяющее свойство ядра Пуассона состоит в том, что интегральный
 оператор с этим ядром  отображает единицу в единицу и является
морфизмом $SU_{n,n}$-модуля функций на
  $S(\mathbb{D})$ в $SU_{n,n}$-модуль функций в $\mathbb{D}$.
Такое определение ядра Пуассона без изменений переносится на квантовый
случай. Например, для квантового круга ядро Пуассона несущественно
отличается от инвариантного ядра $k_{22}^{-1}\,k_{11}^{-1}$, см. \itemiii\
\ref{proof_inv_G_m}, и равно
\begin{equation*}\label{Poisson_old}
P= (1-z^*\zeta)^{-1}\;(1-z^*z)\;(1-z\zeta^*)^{-1}.
\end{equation*}
Аналогичные рассуждения позволяют ввести в рассмотрение ядро Пуассона для
квантового матричного шара. Из определений следует, что, как в
классическом, так и в квантовом случае, образ интегрального оператора с
ядром Пуассона состоит из решений уравнения $\widetilde{\square}u=0$, где
$\widetilde{\square}$ -- инвариантный Лапласиан.

В классическом случае $q=1$ при переходе от круга к $n\times n$-матричному
шару с $n>1$ возникает неожиданный эффект, обнаруженный Хуа. Именно,
интегралы Пуассона удовлетворяют инвариантной системе дифференциальных
уравнений-- системе Хуа
\begin{equation*}\label{classical_Hua}
\sum_{a,b=1}^n \left(\delta_{a,b} - \sum_{c=1}^n
\overline{z_c^\alpha}z_c^\beta\right) \frac{\partial^2\, u}
 { \partial z_c^\beta\, \partial \overline{z_c^\alpha}}=0,\qquad
 \alpha,\beta=1,2,\ldots,n,
\end{equation*}
см. \cite{JohnsonKoranyi}. Уравнение $\widetilde{\square}u=0$ является
следствием системы Хуа. Аналогичные результаты имеют место для квантового
$n\times n$-матричного шара.

\bigskip
Весьма общий подход к построению инвариантных систем дифференциальных
уравнений, решениями которых являются все интегралы Пуассона, предложен в
\cite[предложение 4.2]{BenBurDamRat}. Этот подход переносится на квантовый
случай, как следует из результатов \itemiiiа \ref{vect_diff}, где
обсуждается понятие инвариантной системы дифференциальных уравнений в
$\mathscr{D}(\mathbb{D})_q'$.

\subsection{Аналитическое продолжение голоморфной дискретной серии и
преобразование Пенроуза}

Пусть $\lambda \in \mathbb{C}$. Каноническое вложение, описанное в \itemiiе
\ref{canonical_embed}, позволяет ввести в рассмотрение представление
$\pi_\lambda$ алгебры $U_q\mathfrak{su}_{n,n}$ в векторном пространстве
$\mathbb{C}[\mathrm{Mat}_{n}]_q$. Приведем явный вид операторов
представления образующих $E_j,F_j,K^\pm_j$ :
\begin{equation*}\label{ball2_6.2}
  \pi_\lambda(E_j) f \;=\;
\left\{
\begin{array}{ll}
    E_j f, &\; j \neq  n, \\
    E_n f - q^{1/2} \frac{1 - q^{2\lambda}}{1 - q^2}
(K_n f) z^n_n, &\; j = n,
  \end{array}
\right.
\end{equation*}
\begin{equation*}\label{ball2_6.3}
  \pi_\lambda(F_j) f \;=\;
\left\{
\begin{array}{ll}
    F_j f, &\; j \neq  n, \\
    q^{-\lambda} F_n f, &\; j = n,
  \end{array}
\right.
\end{equation*}
\begin{equation*}\label{ball2_6.4}
  \pi_\lambda(K_j^{\pm 1}) f \;=\;
\left\{
\begin{array}{ll}
    K_j^{\pm 1} f, &\; j \neq  n, \\
    q^{\pm\lambda} K_n^{\pm 1} f, &\; j = n.
  \end{array}
\right.
\end{equation*}

Покажем, что этот $U_q\mathfrak{su}_{n,n}$-модуль унитаризуем при $\lambda
> N-1$. В дальнейшем будет сформулирован результат, позволяющий ослабить
требование на $\lambda$, см. предложение \ref{last_proposition}.

В
\itemiiiе \ref{InvInt} найден явный вид положительного $U_q\mathfrak{g}$-инвариантного
интеграла на $\mathscr{D}(\mathbb{D})_q$. В случае квантового $n\times
n$-матричного шара
\begin{equation*}\label{inv_int_mat}
  \int\limits_{\mathbb{D}_q} f d\nu = (1-q^2)^{n^2}\,
\mathrm{tr}(T_F(f) \Gamma(K_{-2\rho})).
\end{equation*}
Здесь
\begin{equation*}\label{K_ro}
  K_{-2\rho} = \prod\limits_{j=1}^{N-1} K_j^{j(N-j)}.
\end{equation*}
 Как объяснялось в \itemiiiе \ref{InvInt}, след в правой части
равенства (\ref{inv_int_mat}) корректно определен. Нетрудно доказать

\begin{proposition}\label{weighted_ball} Пусть $\lambda
> N-1$. Линейный функционал
$$
\nu_\lambda: \mathscr{D}(\mathbb{D})_q \rightarrow \mathbb{C}, \qquad
\nu_\lambda: f \mapsto C(\lambda)\mathrm{tr}(T_F(f) T_F(y)^\lambda
\Gamma(K_{-2\rho})),
$$
где $ C(\lambda) = \prod\limits_{j=0}^{n-1} \prod\limits_{k=0}^{n-1} (1 -
q^{2(\lambda+1-N)} \; q^{2(j+k)})$, допускает продолжение по непрерывности
на $\mathrm{Pol}(\mathrm{Mat}_{n})_q$, и
\begin{equation*}\label{explicit_nu}
  \int\limits_{\mathbb{D}_q} f d\nu_\lambda = C(\lambda) \mathrm{tr}(T_F(f)
T_F(y)^\lambda \, \Gamma(K_{-2\rho})), \quad f \in
\mathrm{Pol}(\mathrm{Mat}_{n})_q.\footnote{Используется топология
пространства $C(\overline{\mathbb{D}})_q$, см. \itemiii\ \ref{boundness}.}
\end{equation*}
\end{proposition}

При $\lambda > N-1$ представление $\pi_\lambda$ алгебры
$U_q\mathfrak{su}_{n,m}$ в предгильбертовом пространстве
$\mathbb{C}[\mathrm{Mat}_{n}]_q$ со скалярным произведением
\begin{equation}\label{inner_lambda}
(f_1,f_2)_\lambda=\int\limits_{\mathbb{D}_q} f_2^* f_1 d\nu_\lambda,\qquad
f_1,f_2 \in \mathbb{C}[\mathrm{Mat}_{n}]_q
\end{equation}
является $*$-представлением. Это так называемая голоморфная дискретная
серия.

В теории ограниченных симметрических областей хорошо известен результат
Орстеда \cite{Orsted} об аналитическом продолжении интеграла
(\ref{inner_lambda}) по параметру $\lambda$. Приведем $q$-аналог этого
результата. Используя разложение Хуа-Шмида, см. \itemiii\ \ref{Hua-Schmid},
нетрудно показать, что $U_q \mathfrak{s}(\mathfrak{gl}_n \times
\mathfrak{gl}_n)$-модуль $\mathbb{C}[\mathrm{Mat}_{n}]_q$ является суммой
попарно неизоморфных простых весовых $U_q \mathfrak{s}(\mathfrak{gl}_n
\times \mathfrak{gl}_n)$-подмодулей
$\mathbb{C}[\mathrm{Mat}_{n}]_q^{(\mathbf{k})}$, где
$\mathbf{k}=(k_1,k_2,\ldots,k_n) \in \mathbb{Z}_{+}^n$ и $k_1 \geq k_2 \geq
\ldots \geq k_n$.

Эрмитовы формы
$$
(f_1,f_2)_\lambda= \int \limits_{\mathbb{D}_q} f_2^* f_1 d\nu_\lambda,
\qquad (f_1,f_2)_{S(\mathbb{D})_q}=\int \limits_{S(\mathbb{D})_q} f_2^* f_1
d\nu
$$
в $\mathbb{C}[\mathrm{Mat}_{n}]_q$ положительно определены и $U_q
\mathfrak{s}(\mathfrak{u}_n \times \mathfrak{u}_n)$-инвариантны. Значит, на
каждом из подпространств $\mathbb{C}[\mathrm{Mat}_{n}]_q^{(\mathbf{k})}$
они отличаются лишь числовым множителем.\footnote{Это следствие леммы Шура,
а точнее, теоремы Бернсайда \cite[стр. 177]{CurtisReiner}.} Другими
словами, при $\lambda
> N-1$
\begin{equation}\label{unique_up_const}
 \int \limits_{\mathbb{D}_q} f_2^* f_1 d\nu_\lambda = c (\mathbf{k},
\lambda) \int \limits_{S(\mathbb{D})_q} f_2^* f_1 d\nu, \qquad f_1,f_2 \in
\mathbb{C}[\mathrm{Mat}_{m,n}]_q^{(\mathbf{k})}.
\end{equation}

 Квантовый аналог результата Орстеда состоит в том, что каждая из функций
$c(\mathbf{k},\lambda)$ переменной $\lambda$ допускает аналитическое
продолжение до мероморфной функции. Точнее, в \cite{ShZhang_diff} доказано

\begin{proposition}\label{last_proposition} Для всех  $\mathbf{k}=(k_1,k_2,\ldots,k_n)$
\begin{equation*}\label{c_k_lambda}
c(\mathbf{k},\lambda) = \frac{\prod\limits_{j=1}^n
(q^{2(n+1-j)};q^2)_{k_j}}{\prod\limits_{j=1}^n
(q^{2(\lambda+1-j)};q^2)_{k_j}}.
\end{equation*}
\end{proposition}

Доказательство существенно использует результат Милна \cite{Milne},
являющийся частным случаем биномиальной формулы Окунькова для полиномов
Макдональда.

\bigskip Перейдем к преобразованию Пенроуза.
Введем обозначение $(f_1,f_2)_{\lambda}$ для аналитического продолжения
интеграла $\int \limits_{\mathbb{D}_q} f_2^* f_1 d\nu_\lambda$ по параметру
$\lambda$ и обратимся к частному случаю $2\times 2$-матричного шара. В этом
случае
\begin{equation}\label{last_equation}
 c(\mathbf{k},\lambda) =
 \frac{(q^4;q^2)_{k_1}(q^2;q^2)_{k_2}}
 {(q^{2\lambda};q^2)_{k_1}(q^{2(\lambda-1)};q^2)_{k_2}},\qquad k_1 \ge k_2
\ge 0,
\end{equation}
 и функция $(f_1,f_2)_{\lambda}$ переменной $\lambda$ аналитична  при $\lambda > 1$. Как следует
из \eqref{last_equation}, ядро $L$ полуторалинейной формы
 $\lim_{\lambda\to 1+0} (1-q^{2(\lambda-1)})(f_1,f_2)_\lambda$ равно
$$
L=\newoplus_{k=0}^\infty \mathbb{C}[\mathrm{Mat}_2]_q^{(k,0)}.
$$
Разумеется, $L$ является подпредставлением представления $\pi_1$, и на $L$
корректно определена $U_q\mathfrak{su}_{2,2}$-инвариантная эрмитова форма
$$(f_1,f_2)=\lim_{\lambda\to 1+0}(f_1,f_2)_\lambda,$$
наделяющая $L$ структурой предгильбертова пространства.

Получено $*$-представление алгебры $U_q\mathfrak{su}_{2,2}$ в
 $L$-- так называемое {\sl лестничное}
представление. Нетрудно показать, что $L$ является ядром инвариантного
дифференциального оператора, точнее, дифференциального оператора,
сплетающего представления $\pi_1$ и $\pi_3$ алгебры
$U_q\mathfrak{su}_{2,2}$.

В классическом пределе $q\to 1$ получаем известную реализацию лестничного
представления алгебры $U\mathfrak{sl}_4$ в пространстве
$$
\{\psi \in \mathbb{C}[\mathbf{z}]\;\;|\;\; \frac{\partial^2\,\psi}{\partial
z^1_1\,\partial z^2_2}\;-\; \frac{\partial^2\,\psi}{\partial
z^1_2\,\partial z^2_1}\;=\;0\}.
$$
В классическом случае известна еще одна реализация этого представления.
Именно, рассмотрим когеррентный пучок $\mathcal{O}(-2)$ на $\mathbb{CP}^3$
и открытое по Зарисскому множество $$U=\{(u_1:u_2:u_3:u_4)\in
\mathbb{CP}^3\;|\; u_3\ne 0 \;\text{или}\; u_4\ne 0\}.$$ Действие алгебры
$U\mathfrak{sl}_4$ в когомологиях $H^1(U,\mathcal{O}(-2))$ определяет ее
представление, эквивалентное лестничному представлению.

Изоморфизм $U\mathfrak{sl}_4$-модулей
$$
H^1(U,\mathcal{O}(-2))\;\to\; \{\psi \in \mathbb{C}[\mathbf{z}]\;\;|\;\;
\frac{\partial^2\,\psi}{\partial z^1_1\,\partial z^2_2}\;-\;
\frac{\partial^2\,\psi}{\partial z^1_2\,\partial z^2_1}\;=\;0\}
$$
 единствен с точностью до числового
множителя и с точностью до числового множителя равен преобразованию
Пенроуза \cite{BaEast}. Описанные построения переносятся на квантовый
случай и доставляют квантовый аналог преобразования Пенроуза \cite{SSSV2}.

В книге Бэстона и Иствуда \cite{BaEast} в связи с преобразованием Пенроуза
рассматривается ряд более сложных примеров, чем описанный выше пример
лестничного представления алгебры $U\mathfrak{sl}_4$. К сожалению, в этих
более сложных ситуациях методы статьи \cite{SSSV2} применить не удается.
Дело в том, что в \cite{SSSV2} имитируются когомологии Чеха, то есть неявно
используются недостаточно развитые к настоящему времени методы
некоммутативной алгебраической геометрии.

  Альтернативный подход состоит в замене когомологий Чеха когомологиями  Дольбо,
точнее, в использовании квантового аналога метода когомологической индукции
Вогана-Цакермана \cite{KV}. По-видимому, такой квантовый аналог можно
развить, если не стремиться к чрезмерной общности построений, ограничившись
только тем, что необходимо для получения квантовых аналогов результатов
книги \cite{BaEast}. Первый шаг сделан в \cite{SStV_KV1}.

 \begin{remark}  Дополнительным источником  задач некоммутативного анализа в
квантовых ограниченных симметрических областях могут служить статьи и
монографии, посвященные обычным ограниченным симметрическим областям,
например, работы Арази, Жу, Кораньи и Упмайера
 \cite{Arazy,ArazyFisher,AFP,Arazy_Ark,ArUp,ArazyUpmeier, Zhu,
 Zhu_ball, Kor, Koranyi_Inv, Koranyi_harmonic,
 Upm_Proc_Symp,UU,Upm_Trans,Upm_Ann,Upm_Int_Eq,Upm_Indices}.
\end{remark}

\newpage
\renewcommand{\refname}{\centerline{\large\rm \bf ЛИТЕРАТУРА}}
\addcontentsline{toc}{section}{\large ЛИТЕРАТУРА}
\fontsize {14}{15}\selectfont

\end{document}